%% file: PKGGC.tex
\documentclass[10pt,oneside]{amsbook}

\usepackage{amssymb}
\usepackage{tikz-cd}
\usepackage[all,cmtip]{xy}
\usepackage{scalerel}
\usepackage[Euler]{upgreek}
\usepackage{mathtools}
\usepackage{upref}
\usepackage{url}

\usepackage[scr=stixtwofancy]{mathalpha}
\newcommand{\kos}[1]{\mathscr{#1}}

\newcommand{\mon}[1]{\mathcal{\MakeUppercase{#1}}}
\newcommand{\monbase}[1]{\mathfrak{#1}}

\newcommand{\cat}[1]{\mathrm{#1}}

\newcommand{\kosComm}{\kos{C\mkern-2.5mu o\mkern-2.2mu m\mkern-2.5mu m}\mkern-1.5mu}
\newcommand{\kosEns}{\kos{E\!n\!s}}
\newcommand{\kosAff}{\kos{A\mkern-5mu f\mkern-6.5mu f}\mkern-2mu}
\newcommand{\kosRep}{\kos{R\!e\!p}\mkern-1.5mu}
\newcommand{\kosEnd}{\kos{E\mkern-2.5mu n\mkern-3.5mu d}\mkern-1.5mu}
\newcommand{\kosSet}{\kos{S\mkern-2.7mu e\mkern-1.5mu t}\mkern-1.5mu}
\newcommand{\kosCov}{\kos{C\mkern-2.5mu o\mkern-2.5mu v}\mkern-1.5mu}
\newcommand{\kosMod}{\kos{M\mkern-3.5mu o\mkern-4mu d}\mkern-3.5mu}

\newcommand{\prpty}[1]{(\mathrm{#1})}

\newcommand{\CC}{\mathcal{C}}
\newcommand{\CF}{\mathcal{F}}
\newcommand{\CK}{\mathcal{K}}
\newcommand{\CS}{\mathcal{S}}
\newcommand{\CT}{\mathcal{T}}

\newcommand{\tensor}{\circledast}
\newcommand{\unit}{\mathit{1}}

\newcommand{\acts}{\diamond}
\newcommand{\bcts}{\mathbin{\bar\diamond}}

\newcommand{\id}{\mathit{id}}
\newcommand{\pr}{\prime}
\newcommand{\ppr}{\prime\mkern-2.1mu\prime}

\newcommand{\cp}{{\scaleobj{.9}{\vartriangle}}}
\newcommand{\pc}{{\scaleobj{.9}{\triangledown}}}

\newcommand{\dvs}{d}
\newcommand{\varg}{r}
\newcommand{\tw}[2]{\prescript{#1}{}{#2}}

\DeclareMathOperator{\Hom}{Hom}
\DeclareMathOperator{\Isom}{Isom}
\DeclareMathOperator{\End}{End}
\DeclareMathOperator{\Aut}{Aut}

\newcommand{\uHom}{\underline{\Hom}}
\newcommand{\uIsom}{\underline{\Isom}}
\newcommand{\uEnd}{\underline{\End}}
\newcommand{\uAut}{\underline{\Aut}}

\newcommand{\kpr}{\kos{p\!r}}
\newcommand{\slot}{{-}}

\newcommand{\etsf}{\textup{\textsf{\'{e}t}}}

\newcommand{\ev}{\mathrm{ev}}
\newcommand{\cv}{\mathrm{cv}}

\newcommand{\pt}{\{\cdot\}}

\DeclareMathOperator{\Spec}{Spec}
\DeclareMathOperator{\op}{op}

\usetikzlibrary{arrows.meta}
\tikzcdset{
  filled/.style={-{Latex[round, length=1mm, width=2.2mm]}}
}
\tikzcdset{
  2cell/.style={draw=none,
    "{\textstyle\Rightarrow\mkern-#1mu}"{anchor=center, sloped, allow upside down}},
  2cell/.default=4.5
}
\newdimen\morrule \newdimen\moraxis \newdimen\morwd
\newdimen\morht  \newdimen\morbear \newdimen\morlent \newdimen\morlenb
\newcommand{\mor}{\mathrel{
  \morrule=\fontdimen8\textfont3
  \moraxis=\fontdimen22\textfont2
  \sbox0{$\to$}\morwd=\wd0 \morht=\ht0 \advance\morht by -\moraxis
  \morbear=0.054\morwd
  \tikz[baseline=-\moraxis]{%
    \useasboundingbox (0,-\moraxis) rectangle (\morwd,\morht);
    \draw[-{Latex[round, length=0.9mm, width=1.6mm]}, line width=\morrule]
      (\morbear,0) -- (\morwd-\morbear,0);
    }
  }
}
\NewDocumentCommand{\xmor}{o m}{
  \mathrel{
    \morrule=\fontdimen8\textfont3
    \moraxis=\fontdimen22\textfont2
    \sbox0{$\to$}\morwd=\wd0 \morbear=0.054\morwd
    \sbox0{\scriptsize$\mathstrut#2$}\morlent=\wd0 \advance\morlent by 1.35mm
    \IfNoValueTF{#1}{\morlenb=0pt}{
      \sbox2{\scriptsize$\mathstrut#1$}\morlenb=\wd2 \advance\morlenb by 1.35mm}
    \tikz[baseline=-\moraxis]{
      \pgfmathsetlengthmacro{\morL}{max(\morlent,\morlenb,\morwd)}%
      \node[inner sep=0pt] at (-0.225mm, 0.6em) {\scriptsize$\mathstrut #2$};
      \IfNoValueF{#1}{\node[inner sep=0pt] at (-0.225mm, -0.7em)
        {\scriptsize$\mathstrut #1$};}
      \draw[-{Latex[round, length=0.9mm, width=1.6mm]}, line width=\morrule]
        (-0.5*\morL+\morbear, 0) -- (0.5*\morL-\morbear, 0);
    }
  }
}

\newdimen\ovw \newdimen\ovapex \newdimen\ovrise \newdimen\ovgap
\newdimen\ovmin \newdimen\ovrule

\newcommand*{\ovf}{}
\newcommand*{\ovtmp}{}

\newcommand{\ovcore}[4]{%
  \begingroup
  \ifx#3\scriptstyle
    \ovrule=\fontdimen8\scriptfont3 \renewcommand*{\ovf}{0.7}%
  \else\ifx#3\scriptscriptstyle
    \ovrule=\fontdimen8\scriptscriptfont3 \renewcommand*{\ovf}{0.5}%
  \else
    \ovrule=\fontdimen8\textfont3 \renewcommand*{\ovf}{1}%
  \fi\fi
  \ovrise=0.53ex \ovrise=\ovf\ovrise
  \ovgap=0.15ex  \ovgap=\ovf\ovgap
  \ovmin=0.6em  \ovmin=\ovf\ovmin
  \sbox0{$#3#4$}\ovw=\wd0
  \ifdim\ovw<\ovmin \ovw=\ovmin \fi
  \renewcommand*{\ovtmp}{#2}%
  \mathord{\vbox{\offinterlineskip
    \hbox to\ovw{\hss
      \tikz{%
        \ifx\ovtmp\empty
          \draw[#1, line width=\ovrule] (0,0) -- (\ovw,0);
        \else
          \ovapex=\ovw \ovapex=#2\ovapex
          \draw[#1, line width=\ovrule] (0,0) -- (\ovapex,\ovrise) -- (\ovw,0);
        \fi}%
      \hss}%
    \kern\ovgap
    \hbox to\ovw{\hss\box0\hss}}}%
  \endgroup
}
  
\newcommand{\vecar}[1]{\mathpalette{\ovcore{-{Computer Modern Rightarrow[scale=0.6]}}{}}{#1}}
\newcommand{\cevar}[1]{\mathpalette{\ovcore{{Computer Modern Rightarrow[scale=0.6]}-}{}}{#1}}
\newcommand{\hatar}[1]{\mathpalette{\ovcore{-{Computer Modern Rightarrow[scale=0.6]}}{0.4625}}{#1}}
\newcommand{\tahar}[1]{\mathpalette{\ovcore{{Computer Modern Rightarrow[scale=0.6]}-}{0.5375}}{#1}}
\newcommand{\what}[1]{\mathpalette{\ovcore{-}{0.5}}{#1}}

\newtheorem{theorem}{Theorem}[chapter]
\newtheorem{proposition}[theorem]{Proposition}
\newtheorem{lemma}[theorem]{Lemma}
\newtheorem{corollary}[theorem]{Corollary}

\theoremstyle{definition}
\newtheorem{definition}[theorem]{Definition}

\theoremstyle{remark}
\newtheorem{remark}[theorem]{Remark}

\numberwithin{section}{chapter}
\numberwithin{equation}{chapter}

\begin{document}

\frontmatter

\title{Prekosmic Grothendieck/Galois Categories}

\author{Jaehyeok Lee}
\address{Research Institute of Basic Sciences, Ajou University,
Suwon, Republic of Korea}
\email{jhlee1228lee@ajou.ac.kr}
\thanks{This work was supported by the National Research Foundation of Korea(NRF) grant funded by the Korea government(MSIT)(RS-2024-00334558).}

\author{Jae-Suk Park}
\address{Department of Mathematics,
Pohang University of Science and Technology,
Pohang, Republic of Korea}
\email{jaesuk@postech.ac.kr}
\thanks{}

\date{}



\dedicatory{}

\begin{abstract}
  We develop a theory
  on
  duality between groups
  and the categories of their representations,
  by providing an abstract characterization of such categories.
  The theory
  has two parallel formulations,
  each presenting a duality
  analogous to one of the two
  that were 
  originally proposed by Grothendieck:
  the one between
  profinite groups and the abstract categories of their
  continuous representations on finite sets,
  and
  the other between affine group schemes over a field and
  the abstract categories of their finite-dimensional linear representations.
  In the present theory,
  groups and their representations 
  are
  internal to a symmetric monoidal base category,
  with no finiteness condition
  imposed.
  The role of a fiber functor is played by a fiber point,
  an abstraction of
  the forgetful functor with its chosen adjoint.
  Over each base,
  a group determines the category of its representations
  equipped with the canonical fiber point,
  and a pointed abstract category
  determines the group that represents the presheaf of automorphisms of its chosen fiber point.
  These constructions yield mutually quasi-inverse equivalences of $(2,1)$-categories.
  By introducing change of base,
  we gather the dualities over individual bases into a single duality
  in each formulation.
  For a pointed abstract category,
  fiber points of the category are in perfect correspondence with
  torsors under the group associated to the chosen point.
  As an illustration of the theory,
  we reformulate the Galois theory of covering spaces
  for nonempty, path-connected, locally path-connected, semilocally simply connected spaces.
\end{abstract}

\maketitle

\tableofcontents
\mainmatter
\include{PKGGCintro}
\include{Groprekos}

\include{Grocatoverpk}

\include{PKGroCat}

\include{Galprekos}

\include{Galcatoverpk}

\include{PKGalCat}

\appendix
\include{Galcovsp}
\include{ltcatoversmc}
\include{cltcatoversmc}

\backmatter
\providecommand{\bysame}{\leavevmode\hbox to3em{\hrulefill}\thinspace}
\providecommand{\MR}{\relax\ifhmode\unskip\space\fi MR }
\providecommand{\MRhref}[2]{%
  \href{http://www.ams.org/mathscinet-getitem?mr=#1}{#2}
}
\providecommand{\href}[2]{#2}

\end{document}

%% file: PKGGCintro.tex
\chapter{Introduction}
\label{chap-Intro}

Following the idea and vision of Grothendieck,
we establish
dualities between groups
and the categories of their representations,
by providing an abstract characterization
of such categories.
We begin by
recalling
two dualities
that were originally proposed by Grothendieck.

A Galoisian category
is a category $\CC$ that is equivalent to
the category $\uppi\text{-}\cat{FSet}$
of finite sets with a continuous left $\uppi$-action
for some profinite group $\uppi$.\footnote{
  We use the term Galoisian category
  for what is called a Galois category
  in the literature.
  The term translates the \emph{cat\'{e}gorie galoisienne}
  of \cite[Expos\'{e}~V]{Grothendieck1971}
  and avoids confusion with concepts introduced in this paper.
}
Grothendieck
introduced the concept
in \cite[Expos\'{e}~V]{Grothendieck1971}
to describe the fundamental group of a scheme
in parallel with that of a topological space
by finding a common categorical principle
underlying 
the Galois theories of field extensions
and of covering spaces.
He 
gave an equivalent characterization of
a Galoisian category $\CC$
by a list of categorical axioms,
which do not presuppose that such $\uppi$ exists.
The list consists of properties of $\CC$
and the existence of a fiber functor $\CF$,
a functor from $\CC$ to 
the category $\cat{FSet}$ of finite sets,
where $\CF$ is required to have properties
of the forgetful functor
$\CF_{\uppi}\colon \uppi\text{-}\cat{FSet}\to \cat{FSet}$.

Let $\uppi$ be a profinite group.
Then $\uppi\text{-}\cat{FSet}$
is a Galoisian category 
with the fiber functor $\CF_{\uppi}$,
and $\uppi$ is canonically isomorphic to the group
$\Aut(\CF_{\uppi})$ of natural automorphisms of $\CF_{\uppi}$.
Conversely,
let $\CC$ be a Galoisian category.
Each choice of a fiber functor $\CF$ on $\CC$
determines the profinite group $\Aut(\CF)$,
and the chosen $\CF$ lifts to an equivalence
from $\CC$ to $\Aut(\CF)\text{-}\cat{FSet}$.
These correspondences form 
a duality between profinite groups
and Galoisian categories with a choice of a fiber functor.
With this duality at his disposal,
Grothendieck showed that
finite \'{e}tale covers of a connected, locally Noetherian scheme $X$
form a Galoisian category,
and that
taking the fiber along a geometric point $\bar{x}$ of $X$
determines a fiber functor $\CF_{\bar{x}}$,
and defined the \'{e}tale fundamental group
of $X$ based at $\bar{x}$ as the profinite group
$\Aut(\CF_{\bar{x}})$.

It was the idea of Grothendieck
that an analogous duality should also exist
for affine group schemes over a field $k$,
and more generally for gerbes in place of group schemes.
The idea originates in his search for the theory of motives,
the common quintessence of
different Weil cohomology theories for algebraic varieties.
Tannakian categories over $k$
are abstract categories that are in duality with
affine gerbes over the category $\cat{Aff}_k$ of affine $k$-schemes
for the fpqc topology.
The concept was first developed
by Saavedra Rivano \cite{SaavedraRivano1972}
and further developed by Deligne \cite{Deligne1990}.
We refer to \cite{Milne2025}
for an overview of the origins and the theory of Tannakian categories.
Among Tannakian categories over $k$,
we recall the neutral ones, following \cite{Deligne1990}.

A tensor category over $k$
\cite[2.1 and~2.5]{Deligne1990}
is a $k$-linear symmetric monoidal category
such that every object
admits a dual,
the underlying category is abelian,
and the canonical morphism from $k$
to the endomorphism $k$-algebra of the unit object
is an isomorphism.\footnote{
  A symmetric monoidal category is
  said to be $k$-linear
  if the underlying category is equipped with a $k$-linear structure
  in such a way that the monoidal product is $k$-bilinear.
}
We write $\cat{FVec}_k$
for the tensor category over $k$
of finite-dimensional $k$-vector spaces,
and $\cat{FRep}_k(G)$ for the tensor category over $k$
of finite-dimensional linear representations of
an affine group $k$-scheme $G$.
A neutral fiber functor $\upomega$ on a tensor category over $k$
is an exact $k$-linear strong symmetric monoidal functor to $\cat{FVec}_k$;
we follow \cite[2.8]{Deligne1990}.
An example is the forgetful functor
$\upomega_G\colon  \cat{FRep}_k(G)\to \cat{FVec}_k$.
For a neutral fiber functor $\upomega$,
$\uAut_k(\upomega)$ denotes the presheaf of groups on $\cat{Aff}_k$
of monoidal natural automorphisms of $\upomega$.

A neutral Tannakian category over $k$
is a tensor category over $k$
that admits a neutral fiber functor.
A tensor category over $k$
is a neutral Tannakian category over $k$
if and only if
it is equivalent
to $\cat{FRep}_k(G)$
as tensor categories over $k$
for some affine group $k$-scheme $G$.

Let $G$ be an affine group $k$-scheme.
Then $\cat{FRep}_k(G)$
is a neutral Tannakian category over $k$
with the neutral fiber functor $\upomega_G$,
and $G$ canonically represents
$\uAut_k(\upomega_G)$.
Conversely, let $\mathfrak{T}$ be a neutral Tannakian category over $k$.
Each choice of a neutral fiber functor $\upomega$ on $\mathfrak{T}$
determines an affine group $k$-scheme $G(\upomega)$
that represents $\uAut_k(\upomega)$,
and the chosen $\upomega$ lifts to an equivalence 
from $\mathfrak{T}$ to $\cat{FRep}_k(G(\upomega))$
as tensor categories over $k$.
These correspondences form a duality
between affine group $k$-schemes
and neutral Tannakian categories over $k$
with a choice of a neutral fiber functor.

In this paper we establish two types of dualities
in which groups and 
their representations
are internal to a symmetric monoidal category.
Each of them is analogous to one of the two dualities recalled above,
but differs in that no finiteness condition is imposed
on the groups or the representations.
The body of the paper is divided into Part~\ref{part-GroContext},
on the duality analogous to the one above for affine group $k$-schemes,
and Part~\ref{part-GalContext},
on the duality analogous to the other one for profinite groups.
Part~\ref{part-GroContext} and Part~\ref{part-GalContext}
are parallel yet categorically dual to each other.

We start with Part~\ref{part-GroContext}.
A \emph{Grothendieck prekosmos}
is a symmetric monoidal category
whose underlying category has coreflexive equalizers.
For a Grothendieck prekosmos $\kos{K}$,
$\kosAff(\kos{K})$
denotes the cartesian monoidal category 
of \emph{affine schemes in $\kos{K}$},
the opposite of the cocartesian monoidal category of 
commutative monoids in $\kos{K}$.
A \emph{Grothendieck object in $\kos{K}$}
is a group object $\Spec(\pi)$ in $\kosAff(\kos{K})$
such that tensoring with the coordinate monoid $\pi$
preserves coreflexive equalizers.
An \emph{internal representation} of $\Spec(\pi)$
is an object in the underlying category of $\kos{K}$
equipped with a functorial action of $\Spec(\pi)$.
The internal representations
of $\Spec(\pi)$ form a Grothendieck prekosmos,
denoted by $\kosRep_{\kos{K}}(\pi)$.
In the motivating case,
where $\kos{K}$ is the Grothendieck prekosmos
$\kosMod_R$
of modules over a commutative ring $R$,
Grothendieck objects in $\kos{K}$
are the flat affine group $R$-schemes,
and their internal representations
are the linear representations over $R$.

The primary objective of Part~\ref{part-GroContext}
is to find an abstract characterization of
$\kosRep_{\kos{K}}(\pi)$,
consisting only of categorical axioms.
For this purpose, we introduce the $2$-category
of Grothendieck prekosmoi.
A $1$-cell,
called a \emph{Grothendieck morphism},
is an adjunction $\kos{f}=(\kos{f}^*\dashv\kos{f}_*)$
between Grothendieck prekosmoi
such that the left adjoint part $\kos{f}^*$ is strong symmetric monoidal,
and the direction of $\kos{f}$
follows the right adjoint part $\kos{f}_*$.
A $2$-cell from $\kos{f}$ to $\kos{g}$,
called a \emph{Grothendieck transformation},
is a monoidal natural transformation
$\kos{f}^*\Rightarrow\kos{g}^*$
between the left adjoint parts.
Grothendieck morphisms have their own lexicon,
each of whose entries reflects
properties of the adjoint parts individually or in their interplay;
see Definition~\ref{def-Groprekos-Gromorprop}.
For instance, a Grothendieck morphism $\kos{f}$ is said to be
\emph{\'{e}tale}
if it satisfies the projection formula
and $\kos{f}_*$ is conservative,
and to be
\emph{surjective}
if $\kos{f}^*$ is conservative and preserves coreflexive equalizers.
The names of both
the $1$-cells and the $2$-cells,
the entries of the lexicon,
and the corresponding names in Part~\ref{part-GalContext}
are adapted from the vocabulary of topos theory
\cite{Johnstone2002a,Johnstone2002b,MacLaneMoerdijk1994}.
For the relation between the lexicons of
Grothendieck morphisms, their counterparts in Part~\ref{part-GalContext},
and 
essential geometric morphisms between topoi,
see Remark~\ref{rem-Galprekos-Galmornames}.

A \emph{Grothendieck $\kos{K}$-prekosmos}
is a Grothendieck prekosmos $\kos{T}$
equipped with a Grothendieck morphism $\kos{t}$ to $\kos{K}$.
An example is $\kos{K}$ itself
with its identity adjunction $\vphantom{\id}\smash[t]{\id}_{\kos{K}}^{\vphantom{*}}$.
By taking the pseudoslice over $\kos{K}$,
we obtain the $2$-category of Grothendieck $\kos{K}$-prekosmoi.
A $1$-cell from $(\kos{T}^{\pr},\kos{t}^{\pr})$ to $(\kos{T},\kos{t})$
is a Grothendieck morphism $\kos{f}$ from $\kos{T}^{\pr}$ to $\kos{T}$
with an invertible Grothendieck transformation
$\kos{t}^{\pr}\cong \kos{t}\kos{f}$;
we write the $1$-cell as $\kos{f}$.
The $1$-cells and the $2$-cells are called
\emph{Grothendieck $\kos{K}$-morphisms}
and \emph{Grothendieck $\kos{K}$-transformations}, respectively.
A Grothendieck $\kos{K}$-morphism
inherits each entry of the lexicon
from the underlying Grothendieck morphism.

A \emph{fiber point} of a Grothendieck $\kos{K}$-prekosmos
$(\kos{T},\kos{t})$
is an \'{e}tale Grothendieck $\kos{K}$-morphism $\omega$
from $(\kos{K},\vphantom{\id}\smash[t]{\id}_{\kos{K}}^{\vphantom{*}})$ to $(\kos{T},\kos{t})$.
For a Grothendieck object $\Spec(\pi)$ in $\kos{K}$,
$\kosRep_{\kos{K}}(\pi)$
is equipped with a Grothendieck morphism $\kos{t}_{\pi}^{\vphantom{*}}$ to $\kos{K}$
whose left adjoint part
endows an object with the trivial action of $\Spec(\pi)$,
and $(\kosRep_{\kos{K}}(\pi),\kos{t}_{\pi}^{\vphantom{*}})$
is equipped with a surjective fiber point $\kos{e}_{\pi}^{\vphantom{*}}$
whose left adjoint part
is the forgetful functor.
For a fiber point $\omega$,
$\uAut_{\kos{K}}(\omega)$ denotes the presheaf of groups on
$\cat{Aff}(\kos{K})$ of invertible Grothendieck $\kos{K}$-transformations
from $\omega$ to itself,
where $\cat{Aff}(\kos{K})$ is the category of affine schemes in $\kos{K}$.

A \emph{Grothendieck category over $\kos{K}$}
is a Grothendieck $\kos{K}$-prekosmos
that admits a surjective fiber point.\footnote{
  A Grothendieck category over a Grothendieck prekosmos
  is to be distinguished from
  what is called a Grothendieck category in the literature,
  an AB5 abelian category admitting a generator.
}
A Grothendieck $\kos{K}$-prekosmos
is a Grothendieck category over $\kos{K}$
if and only if it is equivalent to
$(\kosRep_{\kos{K}}(\pi),\kos{t}_{\pi}^{\vphantom{*}})$
as Grothendieck $\kos{K}$-prekosmoi
for some Grothendieck object $\Spec(\pi)$ in $\kos{K}$;
see Theorem~\ref{thm-Grocatoverpk-PrimaryObj}.

Let $\Spec(\pi)$ be a Grothendieck object in $\kos{K}$.
Then the pair $(\kosRep_{\kos{K}}(\pi),\kos{t}_{\pi}^{\vphantom{*}})$
is a Grothendieck category over $\kos{K}$
with the surjective fiber point $\kos{e}_{\pi}^{\vphantom{*}}$,
and $\Spec(\pi)$ canonically represents
$\uAut_{\kos{K}}(\kos{e}_{\pi}^{\vphantom{*}})$.
Conversely, let $(\kos{T},\kos{t})$ be a Grothendieck category over $\kos{K}$.
Each choice of a surjective fiber point $\omega$ of $(\kos{T},\kos{t})$
determines the Grothendieck object
$\Spec(\pi(\omega))$ in $\kos{K}$
that represents $\uAut_{\kos{K}}(\omega)$,
and the chosen $\omega$ extends to an equivalence
$\breve{\omega}$
from $(\kosRep_{\kos{K}}(\pi(\omega)),\kos{t}_{\pi(\omega)})$
to $(\kos{T},\kos{t})$
as Grothendieck $\kos{K}$-prekosmoi,
where $\pi(\omega)\coloneqq \omega^*\omega_*(\kappa)$
and $\kappa$ is the unit object of $\kos{K}$;
see Theorem~\ref{thm-Grocatoverpk-mainthm}.
These correspondences form a duality
between Grothendieck objects in $\kos{K}$
and \emph{pointed Grothendieck categories over $\kos{K}$},
those with a choice of a surjective fiber point.
By this we mean that we have mutually quasi-inverse
equivalences of $(2,1)$-categories
\begin{equation}
  \label{eq-Intro-Groequiv21cat}
  \xymatrix@C=35pt{
    \mathsf{pGroObj}(\kos{K})
    \ar@<1ex>[r]^-{\mathsf{Rep}_{\kos{K}}}
    \ar@{}[r]|-{\simeq}
    &\mathsf{pGroCat}_{\mathsf{pt}}(\kos{K})
    \ar@<1ex>[l]^-{\mathsf{Aut}_{\kos{K}}}
  }
\end{equation}
between the $(2,1)$-category 
of Grothendieck objects in $\kos{K}$
and the $(2,1)$-category
of pointed Grothendieck categories over $\kos{K}$;
see Theorem~\ref{thm-dualityoverGropk-equiv21cat}.

An affine scheme in $\kos{K}$
is said to be an \emph{fpqc covering of $\Spec(\kappa)$}
if tensoring with its coordinate monoid is conservative
and preserves coreflexive equalizers.
Given a Grothendieck object $\Spec(\pi)$ in $\kos{K}$,
an \emph{fpqc left $\Spec(\pi)$-torsor over $\Spec(\kappa)$}
is an affine scheme in $\kos{K}$
that is an fpqc covering of $\Spec(\kappa)$
and is equipped with a left action of $\Spec(\pi)$
such that the shear morphism is invertible.
In the motivating case $\kos{K}=\kosMod_R$,
such torsors are the left torsors over $\Spec(R)$
under a flat affine group $R$-scheme for the fpqc topology.

Let $(\kos{T},\kos{t})$ be a Grothendieck category over $\kos{K}$
with a choice of a surjective fiber point $\omega$.
Then fiber points of $(\kos{T},\kos{t})$
are in perfect correspondence with
fpqc $\Spec(\pi)$-torsors over $\Spec(\kappa)$,
where  $\Spec(\pi)$ is the Grothendieck object in $\kos{K}$
that represents $\uAut_{\kos{K}}(\omega)$.
By this we mean that we have an adjoint equivalence
\begin{equation}
  \label{eq-Intro-GroTorsperfcorres}
  \cat{pFib}_{\kos{K}}(\kos{t})
  \simeq
  (\Spec(\pi)\text{-}\cat{Tors})_{\kappa}^{\op}
\end{equation}
between the groupoid of fiber points of $(\kos{T},\kos{t})$
and the opposite of the groupoid of fpqc left $\Spec(\pi)$-torsors over $\Spec(\kappa)$;
see Theorem~\ref{thm-pFibGroTors-perfcorres}.
This is analogous to the case of 
\cite[Theorem~3.2]{DeligneMilne1982}
in which the fiber functors are neutral.
Each fiber point $\omega^{\pr}$ 
is sent to the representing object of the presheaf
$\uIsom_{\kos{K}}(\omega^{\pr},\omega)$ on $\cat{Aff}(\kos{K})$
of invertible Grothendieck $\kos{K}$-transformations
from $\omega^{\pr}$ to $\omega$.
Each fpqc left $\Spec(\pi)$-torsor $\Spec(p)$ over $\Spec(\kappa)$
is sent to the \emph{$\Spec(p)$-twisted fiber point}
of $(\kosRep_{\kos{K}}(\pi),\kos{t}_{\pi}^{\vphantom{*}})$
followed by the equivalence $\breve{\omega}$ extending $\omega$.

We gather the dualities \eqref{eq-Intro-Groequiv21cat},
one for each Grothendieck prekosmos $\kos{K}$,
into a single duality.
The change of base along an \'{e}tale Grothendieck morphism
from $\kos{K}^{\pr}$ to $\kos{K}$
sends each Grothendieck object in $\kos{K}$
to one in $\kos{K}^{\pr}$,
generalizing the change of base of a flat affine group $R$-scheme
along a morphism $\Spec(R^{\pr})\to\Spec(R)$ of affine schemes.
We write
$\mathsf{pGroKos}_{\etsf}$
for the $(2,1)$-category of
Grothendieck prekosmoi,
\'{e}tale Grothendieck morphisms,
and invertible Grothendieck transformations.
A \emph{prekosmic Grothendieck object}
consists of a Grothendieck prekosmos $\kos{K}$
and a Grothendieck object in $\kos{K}$.
Prekosmic Grothendieck objects form a $(2,1)$-category
$\mathsf{pGroObj}$ over $\mathsf{pGroKos}_{\etsf}$,
and the strict fiber over each $\kos{K}$
is $2$-isomorphic to $\mathsf{pGroObj}(\kos{K})$.

A \emph{prekosmic Grothendieck category}
consists of a Grothendieck prekosmos $\kos{K}$
and a Grothendieck category over $\kos{K}$.
A \emph{pointed prekosmic Grothendieck category}
is a prekosmic Grothendieck category
with a choice of a surjective fiber point.
Pointed prekosmic Grothendieck categories likewise form a $(2,1)$-category
$\mathsf{pGroCat}_{\mathsf{pt}}$ over $\mathsf{pGroKos}_{\etsf}$,
and the strict fiber over each $\kos{K}$ is $2$-isomorphic to
$\mathsf{pGroCat}_{\mathsf{pt}}(\kos{K})$.
We have mutually quasi-inverse equivalences
\begin{equation}
  \label{eq-Intro-Groequiv21catoveret}
  \vcenter{\hbox{
    \xymatrix@R=25pt@C=0pt{
      \mathsf{pGroObj}
      \ar@<1ex>[rr]^-{\mathsf{Rep}}
      \ar[dr]_-{\mathsf{p}_{\mathsf{obj}}}
      \ar@{}[rr]|-{\simeq}
      &\text{ }
      &\mathsf{pGroCat}_{\mathsf{pt}}
      \ar@<1ex>[ll]^-{\mathsf{Aut}}
      \ar[dl]^-{\mathsf{p}_{\mathsf{ptcat}}}
      \\%1
      \text{ }
      &\mathsf{pGroKos}_{\etsf}
      &\text{ }
    }
  }}
\end{equation}
of $(2,1)$-categories over $\mathsf{pGroKos}_{\etsf}$;
see Theorem~\ref{thm-PKGroCat-equiv21catoveret}.

We turn to Part~\ref{part-GalContext}.
A \emph{Galois prekosmos}
is a symmetric monoidal category
whose underlying category has reflexive coequalizers.
For a Galois prekosmos $\kos{S}$,
$\kosEns(\kos{S})$
denotes the cartesian monoidal category
of cocommutative comonoids in $\kos{S}$.
A \emph{Galois object in $\kos{S}$}
is a group object $\pi$ in $\kosEns(\kos{S})$
such that tensoring with $\pi$
preserves reflexive coequalizers.
An \emph{internal representation} of $\pi$
is an object in the underlying category of $\kos{S}$
equipped with a functorial action of $\pi$.
The internal representations
of $\pi$ form a Galois prekosmos,
denoted by $\kosRep_{\kos{S}}(\pi)$.
In the motivating case,
where $\kos{S}$ is the Galois prekosmos $\kosSet$ of sets,
Galois objects in $\kos{S}$ are the groups,
and their internal representations
are the sets with a left action of the group.

The primary objective of Part~\ref{part-GalContext}
is to find an abstract characterization of
$\kosRep_{\kos{S}}(\pi)$,
consisting only of categorical axioms.
For this purpose, we introduce the $2$-category
of Galois prekosmoi.
A $1$-cell,
called a \emph{Galois morphism},
is an adjunction $\kos{f}=(\kos{f}_!\dashv\kos{f}^!)$
between Galois prekosmoi
such that the right adjoint part $\kos{f}^!$ is strong symmetric monoidal,
and the direction of $\kos{f}$
follows the left adjoint part $\kos{f}_!$.
A $2$-cell from $\kos{f}$ to $\kos{g}$,
called a \emph{Galois transformation},
is a comonoidal natural transformation
$\kos{f}^!\Rightarrow\kos{g}^!$
between the right adjoint parts.
Galois morphisms have their own lexicon;
see Definition~\ref{def-Galprekos-Galmorprop}
and Remark~\ref{rem-Galprekos-Galmornames}.
For instance, a Galois morphism $\kos{f}$ is said to be
\emph{\'{e}tale}
if it satisfies the projection formula
and $\kos{f}_!$ is conservative,
and to be
\emph{surjective}
if $\kos{f}^!$ is conservative and preserves reflexive coequalizers.

A \emph{Galois $\kos{S}$-prekosmos}
is a Galois prekosmos $\kos{T}$
equipped with a Galois morphism $\kos{t}$ to $\kos{S}$.
An example is $\kos{S}$ itself
with its identity adjunction $\vphantom{\id}\smash[t]{\id}_{\kos{S}\vphantom{!}}^{\vphantom{!}}$.
By taking the pseudoslice over $\kos{S}$,
we obtain the $2$-category of Galois $\kos{S}$-prekosmoi.
A $1$-cell from $(\kos{T}^{\pr},\kos{t}^{\pr})$ to $(\kos{T},\kos{t})$
is a Galois morphism $\kos{f}$ from $\kos{T}^{\pr}$ to $\kos{T}$
with an invertible Galois transformation
$\kos{t}\kos{f}\cong \kos{t}^{\pr}$;
we write the $1$-cell as $\kos{f}$.
The $1$-cells and the $2$-cells are called
\emph{Galois $\kos{S}$-morphisms}
and \emph{Galois $\kos{S}$-transformations}, respectively.
A Galois $\kos{S}$-morphism
inherits each entry of the lexicon from the underlying Galois morphism.

A \emph{fiber point} of a Galois $\kos{S}$-prekosmos
$(\kos{T},\kos{t})$
is an \'{e}tale Galois $\kos{S}$-morphism $\omega$
from $(\kos{S},\vphantom{\id}\smash[t]{\id}_{\kos{S}\vphantom{!}}^{\vphantom{!}})$ to $(\kos{T},\kos{t})$.
For a Galois object $\pi$ in $\kos{S}$,
$\kosRep_{\kos{S}}(\pi)$
is equipped with a Galois morphism $\kos{t}_{\pi\vphantom{!}}^{\vphantom{!}}$ to $\kos{S}$
whose right adjoint part
endows an object with the trivial action of $\pi$,
and $(\kosRep_{\kos{S}}(\pi),\kos{t}_{\pi\vphantom{!}}^{\vphantom{!}})$
is equipped with a surjective fiber point $\kos{u}_{\pi\vphantom{!}}^{\vphantom{!}}$
whose right adjoint part is the forgetful functor.
For a fiber point $\omega$,
$\uAut_{\kos{S}}(\omega)$ denotes the presheaf of groups on
$\cat{Ens}(\kos{S})$ of invertible Galois $\kos{S}$-transformations
from $\omega$ to itself,
where $\cat{Ens}(\kos{S})$ is the category of cocommutative comonoids in $\kos{S}$.

A \emph{Galois category over $\kos{S}$}
is a Galois $\kos{S}$-prekosmos
that admits a surjective fiber point.
A Galois $\kos{S}$-prekosmos
is a Galois category over $\kos{S}$
if and only if it is equivalent to
$(\kosRep_{\kos{S}}(\pi),\kos{t}_{\pi\vphantom{!}}^{\vphantom{!}})$
as Galois $\kos{S}$-prekosmoi
for some Galois object $\pi$ in $\kos{S}$;
see Theorem~\ref{thm-Galcatoverpk-PrimaryObj}.

Let $\pi$ be a Galois object in $\kos{S}$.
Then $(\kosRep_{\kos{S}}(\pi),\kos{t}_{\pi\vphantom{!}}^{\vphantom{!}})$
is a Galois category over $\kos{S}$
with the surjective fiber point $\kos{u}_{\pi\vphantom{!}}^{\vphantom{!}}$,
and $\pi$ canonically represents
$\uAut_{\kos{S}}(\kos{u}_{\pi\vphantom{!}}^{\vphantom{!}})$.
Conversely, let $(\kos{T},\kos{t})$ be a Galois category over $\kos{S}$.
Each choice of a surjective fiber point $\omega$ of $(\kos{T},\kos{t})$
determines the Galois object
$\pi(\omega)$ in $\kos{S}$ that represents $\uAut_{\kos{S}}(\omega)$,
and the chosen $\omega$ extends to an equivalence $\breve{\omega}$
from $(\kosRep_{\kos{S}}(\pi(\omega)),\kos{t}_{\pi(\omega)})$
to $(\kos{T},\kos{t})$
as Galois $\kos{S}$-prekosmoi,
where $\pi(\omega)\coloneqq \omega^!\omega_!(s)$
and $s$ is the unit object of $\kos{S}$;
see Theorem~\ref{thm-Galcatoverpk-mainthm}.
These correspondences form a duality
between Galois objects in $\kos{S}$
and \emph{pointed Galois categories over $\kos{S}$},
those with a choice of a surjective fiber point.
By this we mean that we have mutually quasi-inverse
equivalences of $(2,1)$-categories
\begin{equation}
  \label{eq-Intro-Galequiv21cat}
  \xymatrix@C=35pt{
    \mathsf{pGalObj}(\kos{S})
    \ar@<1ex>[r]^-{\mathsf{Rep}_{\kos{S}}}
    \ar@{}[r]|-{\simeq}
    &\mathsf{pGalCat}_{\mathsf{pt}}(\kos{S})
    \ar@<1ex>[l]^-{\mathsf{Aut}_{\kos{S}}}
  }
\end{equation}
between the $(2,1)$-category 
of Galois objects in $\kos{S}$
and the $(2,1)$-category
of pointed Galois categories over $\kos{S}$;
see Theorem~\ref{thm-dualityoverGalpk-equiv21cat}.

A cocommutative comonoid in $\kos{S}$
is said to be a \emph{co-fpqc covering of $s$}
if tensoring with it is conservative
and preserves reflexive coequalizers.
Given a Galois object $\pi$ in $\kos{S}$,
a \emph{co-fpqc right $\pi$-torsor over $s$}
is a cocommutative comonoid in $\kos{S}$
that is a co-fpqc covering of $s$
and is equipped with a right action of $\pi$
such that the shear morphism is invertible.
In the motivating case $\kos{S}=\kosSet$,
such torsors are the right torsors under a group.

Let $(\kos{T},\kos{t})$ be a Galois category over $\kos{S}$
with a choice of a surjective fiber point $\omega$.
Then fiber points of $(\kos{T},\kos{t})$
are in perfect correspondence
with co-fpqc $\pi$-torsors over $s$,
where $\pi$ is the Galois object in $\kos{S}$ that represents $\uAut_{\kos{S}}(\omega)$.
By this we mean that we have an adjoint equivalence
\begin{equation}
  \label{eq-Intro-GalTorsperfcorres}
  \cat{pFib}_{\kos{S}}(\kos{t})
  \simeq
  (\cat{Tors}\text{-}\pi)_{s}
\end{equation}
between the groupoid of fiber points of $(\kos{T},\kos{t})$
and the groupoid of co-fpqc right $\pi$-torsors over $s$;
see Theorem~\ref{thm-pFibGalTors-perfcorres}.
Each fiber point $\omega^{\pr}$
is sent to the representing object of
 the presheaf
$\uIsom_{\kos{S}}(\omega,\omega^{\pr})$ on $\cat{Ens}(\kos{S})$
of invertible Galois $\kos{S}$-transformations
from $\omega$ to $\omega^{\pr}$.
Each co-fpqc right $\pi$-torsor $p$ over $s$
is sent to the \emph{$p$-twisted fiber point}
of $(\kosRep_{\kos{S}}(\pi),\kos{t}_{\pi\vphantom{!}}^{\vphantom{!}})$
followed by the equivalence $\breve{\omega}$ extending $\omega$.

We gather the dualities \eqref{eq-Intro-Galequiv21cat},
one for each Galois prekosmos $\kos{S}$,
into a single duality.
The change of base along an \'{e}tale Galois morphism
from $\kos{S}^{\pr}$ to $\kos{S}$
sends each Galois object in $\kos{S}$
to one in $\kos{S}^{\pr}$.
We write
$\mathsf{pGalKos}_{\etsf}$
for the $(2,1)$-category of
Galois prekosmoi,
\'{e}tale Galois morphisms,
and invertible Galois transformations.
A \emph{prekosmic Galois object}
consists of a Galois prekosmos $\kos{S}$
and a Galois object in $\kos{S}$.
Prekosmic Galois objects form a $(2,1)$-category
$\mathsf{pGalObj}$ over $\mathsf{pGalKos}_{\etsf}$,
and the strict fiber over each $\kos{S}$
is $2$-isomorphic to $\mathsf{pGalObj}(\kos{S})$.

A \emph{prekosmic Galois category}
consists of a Galois prekosmos $\kos{S}$
and a Galois category over $\kos{S}$.
A \emph{pointed prekosmic Galois category}
is a prekosmic Galois category
with a choice of a surjective fiber point.
Pointed prekosmic Galois categories likewise form a $(2,1)$-category
$\mathsf{pGalCat}_{\mathsf{pt}}$ over $\mathsf{pGalKos}_{\etsf}$,
and the strict fiber over each $\kos{S}$ is $2$-isomorphic to
$\mathsf{pGalCat}_{\mathsf{pt}}(\kos{S})$.
We have mutually quasi-inverse equivalences
\begin{equation}
  \label{eq-Intro-Galequiv21catoveret}
  \vcenter{\hbox{
    \xymatrix@R=25pt@C=0pt{
      \mathsf{pGalObj}
      \ar@<1ex>[rr]^-{\mathsf{Rep}}
      \ar[dr]_-{\mathsf{p}_{\mathsf{obj}}}
      \ar@{}[rr]|-{\simeq}
      &\text{ }
      &\mathsf{pGalCat}_{\mathsf{pt}}
      \ar@<1ex>[ll]^-{\mathsf{Aut}}
      \ar[dl]^-{\mathsf{p}_{\mathsf{ptcat}}}
      \\%1
      \text{ }
      &\mathsf{pGalKos}_{\etsf}
      &\text{ }
    }
  }}
\end{equation}
of $(2,1)$-categories over $\mathsf{pGalKos}_{\etsf}$;
see Theorem~\ref{thm-PKGalCat-equiv21catoveret}.

As an illustration,
we reformulate the Galois theory of covering spaces
in the language of Galois categories
over the Galois prekosmos $\kosSet$ of sets.
By a \emph{uc-space}
we mean a nonempty, path-connected, locally path-connected,
semilocally simply connected space.
Covering spaces of a uc-space $X$ form
a Galois category over $\kosSet$,
denoted by $(\kosCov(X),\kos{t}_{X\vphantom{!}}^{\vphantom{!}})$,
and each point $x\in X$
determines the surjective fiber point of $(\kosCov(X),\kos{t}_{X\vphantom{!}}^{\vphantom{!}})$ at $x$
whose right adjoint part
sends a covering space of $X$ to its fiber over $x$.
The reformulation states that 
the assignment $X\mapsto (\kosCov(X),\kos{t}_{X\vphantom{!}}^{\vphantom{!}})$
defines a pseudofunctor
\begin{equation}
  \label{eq-Intro-Cov21funct}
  \mathsf{Cov}
  \colon
  \mathsf{Top}^{\mathsf{uc}}
  \to 
  \mathsf{pGalCat}(\kosSet)
\end{equation}
from the $(2,1)$-category of uc-spaces 
to the $(2,1)$-category of Galois categories over $\kosSet$;
see Theorem~\ref{thm-Galcovsp-Cov}.
Under the equivalences
\eqref{eq-Intro-Galequiv21cat} for $\kos{S}=\kosSet$,
the pointed variant of \eqref{eq-Intro-Cov21funct}
recovers the $2$-functor
that sends each pointed uc-space
to its fundamental group;
see Theorem~\ref{thm-Galcovsp-Covpt}.

In Part~\ref{part-GroContext}
we utilize the action of the base Grothendieck prekosmos $\kos{K}$
on a Grothendieck $\kos{K}$-prekosmos
to recognize endofunctors on $\kos{K}$
as tensoring with a commutative monoid in $\kos{K}$.
We study this action in Appendix~\ref{appdix-ltcatoversmc}
for a symmetric monoidal category $\kos{K}$ in general,
by introducing the $2$-category of
symmetric monoidal categories equipped with
a lax monoidal functor from $\kos{K}$,
and provide an abstract characterization of
the endofunctor on $\kos{K}$ given by tensoring with
a commutative monoid, a commutative bimonoid,
or a commutative Hopf monoid in $\kos{K}$.
Likewise, in Part~\ref{part-GalContext}
we utilize the action of the base Galois prekosmos $\kos{S}$
on a Galois $\kos{S}$-prekosmos
to recognize endofunctors on $\kos{S}$
as tensoring with a cocommutative comonoid in $\kos{S}$.
We study this action in Appendix~\ref{appdix-cltcatoversmc}
for a symmetric monoidal category $\kos{S}$ in general,
by introducing the $2$-category of
symmetric monoidal categories equipped with
a colax monoidal functor from $\kos{S}$,
and provide an abstract characterization of
the endofunctor on $\kos{S}$ given by tensoring with
a cocommutative comonoid, a cocommutative bimonoid,
or a cocommutative Hopf monoid in $\kos{S}$.

\input{PKGGCai}

\section{Notation}

We fix Grothendieck universes as needed
and tacitly assume that the categories under consideration
are small relative to them.
We use primes, double primes, and tildes
to distinguish mathematical entities of the same kind,
and the notation attached to a decorated symbol inherits the decoration.

Abstract categories are denoted by calligraphic letters such as $\CT$, $\CK$, and $\CS$.
Categories are typically named after their objects in roman type,
and such names may carry arguments or subscripts.
For instance, we write $\cat{Set}$, $\cat{Grp}$, and $\cat{Mon}$
for the categories of sets, groups, and monoids.
The opposite category of $\CT$ is denoted by $\CT^{\op}$.
We write $\mon{f}$, $\mon{g}\colon\CT^{\pr}\to\CT$ for functors
and $\vartheta\colon\mon{f}\Rightarrow\mon{g}$ for a natural transformation.
The symbol $\slot$ marks a functor argument, as in $\mon{f}(\slot)$.
The identity morphism of an object $X$ is denoted by $I_X$,
and the identity natural transformation of $\mon{f}$ by $I_{\mon{f}}$.
The identity functor of $\CT$ is denoted by $\id_{\CT}$
and written as an equality $\id_{\CT}\colon\CT=\CT$.
We write $g\circ f$ for composition of morphisms in a category
and $\mon{f}\mon{f}^{\pr}$ for composition of functors.
For natural transformations,
we write $\tilde{\vartheta}\circ\vartheta$ for vertical composition,
$\vartheta\vartheta^{\pr}$ for horizontal composition,
and $\mon{f}\vartheta^{\pr}$ and $\vartheta\mon{f}^{\pr}$ for whiskering.
Isomorphisms are denoted by $\cong$ and equivalences by $\simeq$.

We refer to \cite{Johnson2021} for $2$-categorical notions.
We denote $2$-categories by typewriter letters such as $\mathtt{B}$,
and $(2,1)$-categories by sans serif letters such as $\mathsf{B}$
to emphasize that their $2$-cells are invertible.
We also use sans serif type for pseudofunctors and $2$-functors,
as in $\mathsf{Rep}$ and $\mathsf{Aut}$,
unless otherwise specified.
The conventions for identities, isomorphisms, equivalences,
and composition apply to every $2$-category,
unless otherwise specified.
For $0$-cells $X$, $Y$ in a $2$-category $\mathtt{B}$,
we write $\mathtt{B}[X,Y]$ for the hom-category
whose objects are $1$-cells from $X$ to $Y$
and whose morphisms are $2$-cells between them.
The pseudoslice $2$-category of $\mathtt{B}$ over $X$ is denoted by $\mathtt{B}/X$,
and the coslice $2$-category of $\mathtt{B}$ under $X$ by $X/\mathtt{B}$;
we indicate in words whether the coslice is lax or colax.
An equivalence of pseudofunctors is an invertible strong transformation.
Since the $2$-cell components of a strong transformation are invertible,
we write them in whichever direction is convenient,
taking their inverses when necessary.
The strict fiber of a $2$-functor over a $0$-cell $X$
consists of the cells lying over $X$,
over $I_X$,
and over the identity $2$-cell of $I_X$.

We refer to \cite{Johnson2021} for monoidal categorical notions.
We use script type to distinguish a symmetric monoidal category
from its underlying category in notation.
We write $\kos{T}=(\CT,\tensor,\unit)$
for a symmetric monoidal category with underlying category $\CT$,
monoidal product $\tensor$, and unit object $\unit$.
Objects in $\CT$ are denoted by uppercase letters such as $X$, $Y$, and $Z$.
The coherence isomorphisms of $\kos{T}$ are denoted by
$a_{X,Y,Z}\colon X\tensor (Y\tensor Z)\cong (X\tensor Y)\tensor Z$,
$\imath_X\colon X\cong \unit\tensor X$,
$\jmath_X\colon X\cong X\tensor \unit$,
and $\sigma_{X,Y}\colon X\tensor Y\cong Y\tensor X$.
We omit $a$, $\imath$, and $\jmath$ from the notation
except when they need to be indicated,
as justified by the Mac Lane coherence theorem.
The symbol $\kos{K}=(\CK,\otimes,\kappa)$
is reserved for the base Grothendieck prekosmos in Part~\ref{part-GroContext}
and for the base symmetric monoidal category in Appendix~\ref{appdix-ltcatoversmc}.
The symbol $\kos{S}=(\CS,\otimes,s)$
is reserved for the base Galois prekosmos in Part~\ref{part-GalContext}
and for the base symmetric monoidal category in Appendix~\ref{appdix-cltcatoversmc}.
Objects in $\CK$ and in $\CS$
are denoted by lowercase letters such as $x$, $y$, and $z$.

Monoidal functors between symmetric monoidal categories
are understood to be symmetric.
We write $\mon{f}\colon\kos{T}^{\pr}\to\kos{T}$
for a monoidal functor,
specifying whether it is lax, colax, or strong,
and the same symbol denotes its underlying functor
$\mon{f}\colon\CT^{\pr}\to\CT$.
The tensor and unit coherence morphisms of $\mon{f}$
are omitted from the notation except when they need to be indicated,
in which case they are denoted by
$(\mon{f}_{\tensor})_{\slot,\slot}$ and $\mon{f}_{\unit}$,
unless otherwise specified.
The coherence morphisms point in the lax direction
throughout Part~\ref{part-GroContext} and Appendix~\ref{appdix-ltcatoversmc},
and in the colax direction
throughout Part~\ref{part-GalContext} and Appendix~\ref{appdix-cltcatoversmc}.
Natural transformations compatible with lax coherence morphisms
are called monoidal,
and those compatible with colax coherence morphisms comonoidal.

We use the standard notions of monoids, comonoids, bimonoids,
and Hopf monoids in a symmetric monoidal category,
and of their actions, coactions, modules, and comodules,
with their usual morphisms.
Monoids are understood to be unital and comonoids counital.
We also use the notions of monoid objects and group objects
in cartesian monoidal categories,
and of comonoid objects in cocartesian monoidal categories.

Let $\kos{T}=(\CT,\tensor,\unit)$ be a symmetric monoidal category.
A dual of an object $X\in\CT$ in $\kos{T}$
is a triple $(X^{\vee},\ev,\cv)$
consisting of an object $X^{\vee}\in\CT$
and morphisms
$\ev\colon X^{\vee}\tensor X\to\unit$,
$\cv\colon\unit\to X\tensor X^{\vee}$
satisfying the triangle identities.
We say that $X$ is dualizable if it admits a dual in $\kos{T}$.

Homotopies between continuous maps are denoted by $\simeq$,
and they appear only in Appendix~\ref{appdix-Galcovsp}.
 

%% file: PKGGCai.tex
\subsection*{Acknowledgements}

The authors developed and verified the mathematical content of
Chapters~\ref{chap-Groprekos}--\ref{chap-PKGalCat}
and Appendices~\ref{appdix-ltcatoversmc} and~\ref{appdix-cltcatoversmc},
and formulated the main statements in
Appendix~\ref{appdix-Galcovsp}, without AI assistance,
apart from a few supplementary remarks
and assistance with formulas in two proofs in
Chapter~\ref{chap-PKGalCat}, as specified below.
Fable (Anthropic, Claude Fable~5 and~5.1) was otherwise used only
to assist with the English and mathematical exposition.
An initial draft of Chapter~\ref{chap-Intro} was also completed
with AI assistance limited to this editorial role.
At a late stage of preparation,
Fable and ChatGPT (OpenAI, GPT-6 Astra)
assisted in drafting and revising the exposition.
ChatGPT was also used to proofread the manuscript
and to help complete its remaining parts.

ChatGPT suggested Remarks~\ref{rem-Groprekos-omegalowerstarconservative}
and~\ref{rem-Galprekos-omegalowershriekconservative}.
It also helped clarify Remark~\ref{rem-Galprekos-Galmornames}~(2), (3)
and verify Remarks~\ref{rem-PKGroCat-cleavagetransport}
and~\ref{rem-PKGalCat-cleavagetransport}.
In Appendix~\ref{appdix-Galcovsp},
it assisted in drafting and refining the proofs.
Fable helped in finding and verifying
the formula for the component $\vartheta_X$
in the proof of Proposition~\ref{prop-dualityPKGal-Rep},
and the formulas for 
$\grave{\vartheta}^{\sharp}$ and $\vartheta^{\sharp}$
in the proof of Proposition~\ref{prop-dualityPKGal-Aut}.

The authors take full responsibility for the mathematical content,
exposition, and references, including all material prepared with AI assistance.

%% file: Groprekos.tex
\part{On Prekosmic Grothendieck Categories}
\label{part-GroContext}

\chapter{Grothendieck Prekosmoi}
\label{chap-Groprekos}

\begin{definition}
  \label{def-Groprekos-Groprekos}
  A \emph{Grothendieck prekosmos} is a symmetric monoidal category
  whose underlying category has coreflexive equalizers.
\end{definition}

Throughout this chapter,
we fix a Grothendieck prekosmos 
$\kos{K}=(\CK,\otimes,\kappa)$,
with underlying category $\CK$,
monoidal product $\otimes\colon\CK\times\CK\to \CK$,
and unit object $\kappa$.
Objects in $\CK$ are denoted by lowercase letters
such as $x$, $y$, and $z$.
We write $\kos{T}=(\CT,\tensor,\unit)$
for an arbitrary Grothendieck prekosmos,
with underlying category $\CT$,
monoidal product $\tensor\colon\CT\times\CT\to \CT$,
and unit object $\unit$.
Objects in $\CT$ are denoted by uppercase letters
such as $X$, $Y$, and $Z$.

\begin{definition}
  \label{def-Groprekos-Gromor}
  Let $\kos{T}=(\CT,\tensor,\unit)$,
  $\kos{T}^{\pr}=(\CT^{\pr},\tensor^{\pr},\unit^{\pr})$
  be Grothendieck prekosmoi.
  A \emph{Grothendieck morphism}
  $\kos{f}\colon\kos{T}^{\pr}\mor \kos{T}$,
  denoted by an arrow with a filled tip,
  is a tuple $\kos{f}=(\kos{f}^*,\kos{f}_*,\eta,\epsilon)$
  consisting of a strong monoidal functor
  $\kos{f}^*\colon\kos{T}\to \kos{T}^{\pr}$,
  a functor $\kos{f}_*\colon\CT^{\pr}\to \CT$,
  and natural transformations
  $\eta\colon\vphantom{\id}\smash[t]{\id}_{\CT}^{\vphantom{*}}\Rightarrow \kos{f}_*\kos{f}^*$,
  $\epsilon\colon\kos{f}^*\kos{f}_*\Rightarrow \vphantom{\id}\smash[t]{\id}_{\CT^{\pr}}^{\vphantom{*}}$
  satisfying the triangle identities.
  Thus $\kos{f}^*$ is left adjoint to $\kos{f}_*$,
  and we refer to $\kos{f}^*$ (resp.\ $\kos{f}_*$)
  as the \emph{left adjoint part} (resp.\ \emph{right adjoint part})
  of $\kos{f}$.
  We omit $\eta$ and $\epsilon$ from the notation
  and write $\kos{f}=(\kos{f}^*\dashv \kos{f}_*)$.
\end{definition}

\begin{definition}
  \label{def-Groprekos-Grotrans}
  Let $\kos{T}$, $\kos{T}^{\pr}$ be Grothendieck prekosmoi
  and $\kos{f}$, $\kos{g}\colon\kos{T}^{\pr}\mor \kos{T}$
  be Grothendieck morphisms.
  A \emph{Grothendieck transformation}
  $\vartheta\colon\kos{f}\Rightarrow \kos{g}$
  is a monoidal natural transformation
  $\vartheta\colon\kos{f}^*\Rightarrow \kos{g}^*$
  between the left adjoint parts.
\end{definition}

We denote by $\mathtt{pGroKos}$ the $2$-category
(resp.\ by $\mathsf{pGroKos}$ the $(2,1)$-category)
of Grothendieck prekosmoi,
Grothendieck morphisms,
and Grothendieck transformations
(resp.\ invertible Grothendieck transformations).
Composition of Grothendieck morphisms 
is given by that of adjunctions.
The identity Grothendieck morphism
of a Grothendieck prekosmos $\kos{T}$
is the tuple consisting of
identity functors and identity natural transformations;
we denote it by $\vphantom{\id}\smash[t]{\id}_{\kos{T}}^{\vphantom{*}}$
and write it as an equality
$\vphantom{\id}\smash[t]{\id}_{\kos{T}}^{\vphantom{*}}\colon\kos{T}=\kos{T}$.
Horizontal and vertical compositions
of Grothendieck transformations
are given by those of monoidal natural transformations
between the left adjoint parts.
Identity Grothendieck transformations
are given by the identity natural transformations
of the left adjoint parts.

Let $\kos{T}$, $\kos{T}^{\pr}$ be Grothendieck prekosmoi.
For each Grothendieck morphism
$\kos{f}\colon\kos{T}^{\pr}\mor \kos{T}$,
the right adjoint part $\kos{f}_*$
has a unique lax monoidal functor structure
such that $\eta$ and $\epsilon$ are monoidal natural transformations.
Conversely, if we are given lax monoidal functors
$\kos{f}^*\colon\kos{T}\to \kos{T}^{\pr}$,
$\kos{f}_*\colon\kos{T}^{\pr}\to \kos{T}$
and monoidal natural transformations
$\eta\colon\vphantom{\id}\smash[t]{\id}_{\kos{T}}^{\vphantom{*}}\Rightarrow \kos{f}_*\kos{f}^*$,
$\epsilon\colon\kos{f}^*\kos{f}_*\Rightarrow \vphantom{\id}\smash[t]{\id}_{\kos{T}^{\pr}}^{\vphantom{*}}$
satisfying the triangle identities,
then $\kos{f}^*$ is necessarily strong monoidal,
and $\kos{f}=(\kos{f}^*\dashv \kos{f}_*)\colon\kos{T}^{\pr}\mor \kos{T}$
is a Grothendieck morphism.
This is an instance of doctrinal adjunction \cite[Theorem~1.5]{Kelly1974}
applied to the $2$-monad on the $2$-category of categories
whose algebras are symmetric monoidal categories.

To each Grothendieck morphism
$\kos{f}\colon\kos{T}^{\pr}\mor \kos{T}$
between Grothendieck prekosmoi
$\kos{T}=(\CT,\tensor,\unit)$,
$\kos{T}^{\pr}=(\CT^{\pr},\tensor^{\pr},\unit^{\pr})$,
we associate a natural transformation
$\varphi\colon
\slot\tensor \kos{f}_*(\slot)\Rightarrow \kos{f}_*(\kos{f}^*(\slot)\tensor^{\pr}\slot)$
between functors from $\CT\times\CT^{\pr}$ to $\CT$
whose component at $Z\in \CT$, $X^{\pr}\in \CT^{\pr}$ is
\begin{equation}
  \label{eq-Groprekos-projformula}
  \xymatrix@C=35pt{
    \varphi_{Z, X^{\pr}}\colon
    Z\tensor \kos{f}_*(X^{\pr})
    \ar[r]^-{\eta_Z\tensor I_{\kos{f}_*(X^{\pr})}}
    &\kos{f}_*\kos{f}^*(Z)\tensor \kos{f}_*(X^{\pr})
    \ar[r]^-{(\kos{f}_{*\tensor})_{\kos{f}^*(Z), X^{\pr}}}
    &\kos{f}_*(\kos{f}^*(Z)\tensor^{\pr} X^{\pr})
    .
  }
\end{equation}
We say that $\kos{f}$ \emph{satisfies the projection formula},
or that \emph{the projection formula holds for} $\kos{f}$,
if the associated $\varphi$ is a natural isomorphism.

Grothendieck morphisms have their own lexicon,
each of whose entries reflects
properties of the adjoint parts individually or in their interplay.
The terminology is adapted from that of
essential geometric morphisms between topoi;
see Remark~\ref{rem-Galprekos-Galmornames} for an explanation.

\begin{definition}
  \label{def-Groprekos-Gromorprop}
  A Grothendieck morphism $\kos{f}$ between Grothendieck prekosmoi
  is said to be
  \begin{enumerate}
    \item 
    \emph{connected}
    if the left adjoint part $\kos{f}^*$ is fully faithful.

    \item 
    \emph{locally connected}
    if it satisfies the projection formula.
    
    \item
    \emph{\'{e}tale}
    if it is locally connected
    and the right adjoint part $\kos{f}_*$ is conservative.
    
    \item 
    \emph{surjective}
    if the left adjoint part $\kos{f}^*$
    is conservative and preserves coreflexive equalizers.
    
    \item
    an \emph{embedding}
    if the right adjoint part $\kos{f}_*$
    is fully faithful.

  \end{enumerate}  
\end{definition}

We write $\mathtt{pGroKos}/\kos{K}$
for the pseudoslice $2$-category of $\mathtt{pGroKos}$ over $\kos{K}$.
\begin{itemize}
  \item 
  A $0$-cell in $\mathtt{pGroKos}/\kos{K}$
  is a \emph{Grothendieck $\kos{K}$-prekosmos},
  or a \emph{Grothendieck prekosmos over $\kos{K}$}.
  It is a pair $(\kos{T},\kos{t})$ consisting of a Grothendieck prekosmos $\kos{T}$
  and a Grothendieck morphism $\kos{t}\colon\kos{T}\mor \kos{K}$.
  We write the pair as $\kos{T}\xmor{\kos{t}} \kos{K}$,
  or simply as $\kos{t}$.
  
  \item
  A $1$-cell in $\mathtt{pGroKos}/\kos{K}$
  is a \emph{Grothendieck $\kos{K}$-morphism},
  or a \emph{Grothendieck morphism over $\kos{K}$}.
  Let $\kos{T}\xmor{\kos{t}} \kos{K}$,
  $\kos{T}^{\pr}\xmor{\kos{t}^{\pr}} \kos{K}$
  be Grothendieck $\kos{K}$-prekosmoi.
  A Grothendieck $\kos{K}$-morphism
  $(\kos{f},\what{\kos{f}})\colon\kos{t}^{\pr}\mor \kos{t}$,
  denoted by an arrow with a filled tip,
  is a pair consisting of a Grothendieck morphism
  $\kos{f}\colon\kos{T}^{\pr}\mor \kos{T}$
  and an invertible Grothendieck transformation
  $\what{\kos{f}}\colon\kos{t}^{\pr}\cong \kos{t}\kos{f}$.
  We omit $\what{\kos{f}}$ from the notation and
  write $\kos{f}\colon\kos{t}^{\pr}\mor \kos{t}$.
  
  \item
  A $2$-cell in $\mathtt{pGroKos}/\kos{K}$
  is a \emph{Grothendieck $\kos{K}$-transformation},
  also called a \emph{Grothendieck transformation over $\kos{K}$}.
  Let 
  $\kos{f}$, $\kos{g}\colon\kos{t}^{\pr}\mor \kos{t}$
  be Grothendieck $\kos{K}$-morphisms.
  By a Grothendieck $\kos{K}$-transformation
  $\vartheta\colon\kos{f}\Rightarrow \kos{g}$
  we mean a Grothendieck transformation
  between the underlying Grothendieck morphisms
  $\kos{f}$, $\kos{g}\colon\kos{T}^{\pr}\mor \kos{T}$
  that satisfies
  $\what{\kos{g}}=
  (\kos{t}\vartheta)\circ \what{\kos{f}}
  \colon\kos{t}^{\pr}\Rightarrow \kos{t}\kos{g}$.
\end{itemize}
We illustrate Grothendieck $\kos{K}$-morphisms
$\kos{f}$, $\kos{g}\colon\kos{t}^{\pr}\mor\kos{t}$
and a Grothendieck $\kos{K}$-transformation
$\vartheta\colon\kos{f}\Rightarrow\kos{g}$
as follows.
\begin{equation*}
  \begin{tikzcd}[row sep=5pt,column sep=5pt]
    \kos{T}^{\pr}
    \ar[rr,filled,"\kos{f}"]
    \ar[ddr,filled,"\kos{t}^{\pr}"']
    &\text{ }
    &\kos{T}
    \ar[ddl,filled,"\kos{t}"]
    &\text{ }
    &\kos{T}^{\pr}
    \ar[rr,filled,"\kos{g}"]
    \ar[ddr,filled,"\kos{t}^{\pr}"']
    &\text{ }
    &\kos{T}
    \ar[ddl,filled,"\kos{t}"]
    &\text{ }
    &\kos{T}^{\pr}
    \ar[rr,filled,shift left=0.5ex,bend left=22pt,"\kos{g}" description]
    \ar[rr,filled,shift right=0.5ex,bend right=22pt,"\kos{f}" description]
    \ar[ddr,filled,"\kos{t}^{\pr}"']
    &\text{ }
    &\kos{T}
    \ar[ddl,filled,"\kos{t}"]
    \\%1
    \text{ }
    \ar[rr,2cell,shift right=0.5ex,"{\what{\kos{f}}}"]
    &\text{ }
    &\text{ }
    &\text{ }
    &\text{ }
    \ar[rr,2cell,shift right=0.5ex,"{\!\what{\kos{g}}}"]
    &\text{ }
    &\text{ }
    &\text{ }
    &\text{ }
    &\text{ }
    \ar[u,2cell,pos=2.2,"{\vartheta}"']
    &\text{ }
    \\%2
    \text{ }
    &\kos{K}
    &\text{ }
    &\text{ }
    &\text{ }
    &\kos{K}
    &\text{ }
    &\text{ }
    &\text{ }
    &\kos{K}
    &\text{ }
  \end{tikzcd}
\end{equation*}

An \emph{equivalence} of Grothendieck prekosmoi
(resp.\ of Grothendieck $\kos{K}$-prekosmoi)
is an equivalence internal to $\mathtt{pGroKos}$
(resp.\ $\mathtt{pGroKos}/\kos{K}$).
A Grothendieck morphism $\kos{f}$
between Grothendieck prekosmoi
is an equivalence if and only if
$\eta$ and $\epsilon$ are natural isomorphisms;
the if part follows from doctrinal adjunction.
A Grothendieck $\kos{K}$-morphism
between Grothendieck $\kos{K}$-prekosmoi
is an equivalence
if and only if the underlying Grothendieck morphism
is an equivalence.

Let $\prpty{P}$ be one of the following:
connected,
locally connected,
\'{e}tale,
surjective,
an embedding.
We say that a Grothendieck $\kos{K}$-prekosmos
is $\prpty{P}$ if it is $\prpty{P}$
as a Grothendieck morphism to $\kos{K}$.
We say that a Grothendieck $\kos{K}$-morphism
between Grothendieck $\kos{K}$-prekosmoi
is $\prpty{P}$
if the underlying Grothendieck morphism
is $\prpty{P}$.
The property of being $\prpty{P}$
is preserved under composition with equivalences,
both for Grothendieck morphisms
and for Grothendieck $\kos{K}$-morphisms.

\begin{definition}
  \label{def-Groprekos-PKfiberpt}
  Let $\kos{T}\xmor{\kos{t}} \kos{K}$
  be a Grothendieck $\kos{K}$-prekosmos.
  A \emph{fiber point} of $\kos{t}$
  is an \'{e}tale Grothendieck $\kos{K}$-morphism
  $\omega\colon\vphantom{\id}\smash[t]{\id}_{\kos{K}}^{\vphantom{*}}\mor \kos{t}$.
  We illustrate $\omega$ as follows.
  \begin{equation*}
    \begin{tikzcd}[row sep=5pt,column sep=5pt]
      \kos{K}
      \ar[rr,filled,"\omega"]
      \ar[ddr,equal]
      &\text{ }
      &\kos{T}
      \ar[ddl,filled,"\kos{t}"]
      \\%1
      \text{ }
      \ar[rr,2cell,shift right=0.5ex,"{\!\what{\omega}}"]
      &\text{ }
      &\text{ }
      \\%2
      \text{ }
      &\kos{K}
      &\text{ }
    \end{tikzcd}
  \end{equation*}
\end{definition}

\begin{definition}
  \label{def-Groprekos-Grocatoverpk}
  A \emph{Grothendieck $\kos{K}$-category},
  or a \emph{Grothendieck category over $\kos{K}$},
  is a Grothendieck $\kos{K}$-prekosmos
  that admits a surjective fiber point.
\end{definition}

Let 
$\kos{T}\xmor{\kos{t}}\kos{K}$
be a Grothendieck $\kos{K}$-prekosmos.
Note that a fiber point
$\omega\colon\vphantom{\id}\smash[t]{\id}_{\kos{K}}^{\vphantom{*}}\mor\kos{t}$
of $\kos{t}$ is surjective
if and only if the underlying Grothendieck morphism
$\omega\colon\kos{K}\mor\kos{T}$ is surjective.
A \emph{morphism} of fiber points of $\kos{t}$
is a Grothendieck $\kos{K}$-transformation between them.
We denote by $\cat{pFib}_{\kos{K}}(\kos{t})$
the category of fiber points of $\kos{t}$.

Let $\kos{f}\colon\kos{t}^{\pr}\mor\kos{t}$
be a Grothendieck $\kos{K}$-morphism
between Grothendieck $\kos{K}$-prekosmoi.
When describing 
$\what{\kos{f}}\colon\kos{t}^{\pr}\cong \kos{t}\kos{f}$
as the monoidal natural isomorphism between the left adjoint parts,
we write
$\what{\kos{f}^*}\colon\kos{t}^{\pr*}\cong \kos{f}^*\kos{t}^*$.
Let $\kos{f}^{\pr}\colon\kos{t}^{\ppr}\mor\kos{t}^{\pr}$
be a Grothendieck $\kos{K}$-morphism,
and consider the composite
$\kos{f}\kos{f}^{\pr}\colon\kos{t}^{\ppr}\mor\kos{t}$.
When describing
$\what{\kos{f}\kos{f}^{\pr}}=(\what{\kos{f}}\kos{f}^{\pr})\circ \what{\kos{f}^{\pr}}
\colon\kos{t}^{\ppr}\cong\kos{t}\kos{f}\kos{f}^{\pr}$
as the monoidal natural isomorphism between the left adjoint parts,
we write
$\what{\kos{f}^{\pr*}\kos{f}^*}
=(\kos{f}^{\pr*}\what{\kos{f}^*})\circ \what{\kos{f}^{\pr*}}\colon\kos{t}^{\ppr*}\cong \kos{f}^{\pr*}\kos{f}^*\kos{t}^*$.

\begin{remark}
  \label{rem-Groprekos-omegalowerstarconservative}
  Let $\kos{T}\xmor{\kos{t}}\kos{K}$
  be a Grothendieck $\kos{K}$-prekosmos
  and
  $\omega\colon \vphantom{\id}\smash[t]{\id}_{\kos{K}}^{\vphantom{*}}\mor\kos{t}$ be a Grothendieck $\kos{K}$-morphism.
  Then the right adjoint part $\omega_*$ is conservative,
  and thus $\omega$ is 
  a fiber point of $\kos{t}$
  if and only if it is locally connected.
  Indeed, by uniqueness of right adjoints,
  the natural isomorphism
  $\what{\omega^*}\colon \vphantom{\id}\smash[t]{\id}_{\kos{K}}^{\vphantom{*}}\cong \omega^*\kos{t}^*$
  gives $\vphantom{\id}\smash[t]{\id}_{\kos{K}}^{\vphantom{*}}\cong \kos{t}_*\omega_*$.
  Hence the composite $\kos{t}_*\omega_*$ is conservative,
  and so is $\omega_*$.
\end{remark}

\section{Grothendieck objects in a Grothendieck prekosmos}
\label{sec-Groobjinpk}

We denote by $\kosComm(\kos{K})=(\cat{Comm}(\kos{K}),\otimes,\kappa)$
the cocartesian monoidal category
of commutative monoids in $\kos{K}$.
The product and unit morphisms
of a commutative monoid $b$ in $\kos{K}$
are denoted by $\pc_b\colon b\otimes b\to b$
and $u_b\colon \kappa\to b$.
Note that $\kosComm(\kos{K})$
is a Grothendieck prekosmos,
since the forgetful functor
$\cat{Comm}(\kos{K})\to \CK$
creates coreflexive equalizers.

We write $\cat{Aff}(\kos{K})$
for the opposite category of $\cat{Comm}(\kos{K})$,
and by an \emph{affine scheme in $\kos{K}$}
we mean an object in $\cat{Aff}(\kos{K})$.
The affine scheme in $\kos{K}$
corresponding to a commutative monoid $b$ in $\kos{K}$
is denoted by $\Spec(b)$.
For a morphism $f\colon b\to b^{\pr}$ in $\cat{Comm}(\kos{K})$,
we use the same symbol $f\colon \Spec(b^{\pr})\to\Spec(b)$
for the corresponding morphism in $\cat{Aff}(\kos{K})$;
whenever confusion may arise,
we indicate the category in which the morphism lies,
in words or by specifying the source and target. 
We write
$\kosAff(\kos{K})=(\cat{Aff}(\kos{K}),\times,\Spec(\kappa))$
for the cartesian monoidal category
of affine schemes in $\kos{K}$,
the opposite monoidal category of $\kosComm(\kos{K})$,
with $\Spec(b)\times \Spec(b^{\pr})=\Spec(b\otimes b^{\pr})$.

An \emph{affine group scheme in $\kos{K}$}
(resp.\ \emph{affine monoid scheme in $\kos{K}$})
is a group object in $\kosAff(\kos{K})$
(resp.\ monoid object in $\kosAff(\kos{K})$).
An affine scheme $\Spec(\pi)$ in $\kos{K}$
is an affine group scheme in $\kos{K}$
if and only if $\pi$ is a commutative Hopf monoid in $\kos{K}$,
in which case we write
$\cp_{\pi}\colon \Spec(\pi)\times\Spec(\pi)\to \Spec(\pi)$,
$e_{\pi}\colon \Spec(\kappa)\to \Spec(\pi)$,
and
$\varsigma_{\pi}\colon \Spec(\pi)\to \Spec(\pi)$
for the product, unit, and inverse morphisms
of $\Spec(\pi)$ in $\cat{Aff}(\kos{K})$;
the corresponding morphisms
$\cp_{\pi}\colon \pi\to \pi\otimes\pi$,
$e_{\pi}\colon \pi\to \kappa$,
and
$\varsigma_{\pi}\colon \pi\to \pi$
in $\cat{Comm}(\kos{K})$
are the coproduct, counit, and antipode of $\pi$.
For such $\pi$,
the antipode $\varsigma_{\pi}$
is involutive and anticommutes with $\cp_{\pi}$,
and the functor $\slot\otimes\pi\colon \CK\to \CK$ is conservative.
A \emph{morphism} of affine group schemes in $\kos{K}$
is a morphism of group objects in $\kosAff(\kos{K})$.

Let $\Spec(\pi)$ be an affine group scheme in $\kos{K}$.
For each affine scheme $\Spec(b)$ in $\kos{K}$, the set
$\Hom_{\cat{Aff}(\kos{K})}(\Spec(b),\Spec(\pi))
=\Hom_{\cat{Comm}(\kos{K})}(\pi,b)$
is a group under the convolution product $\star$,
which sends each pair
$\tilde{g}$, $g\colon \pi\to b$ to
$\tilde{g}\star g\coloneqq\pc_b\circ (\tilde{g}\otimes g)\circ \cp_{\pi}\colon \pi\to b$.
The identity element is
$u_b\circ e_{\pi}\colon \pi\to b$,
and the inverse of $g\colon \pi\to b$
is $g\circ \varsigma_{\pi}\colon \pi\to b$.
It follows from the Yoneda lemma that the assignment
$\Spec(\pi)\mapsto \Hom_{\cat{Aff}(\kos{K})}(\slot,\Spec(\pi))$
is an equivalence
from the category of affine group schemes in $\kos{K}$
to the category of presheaves of groups on $\cat{Aff}(\kos{K})$
whose underlying presheaves of sets are representable.

\begin{definition}
  \label{def-Groobjinpk-Groobj}
  A \emph{Grothendieck object in $\kos{K}$}
  is an affine group scheme $\Spec(\pi)$ in $\kos{K}$
  such that the functor $\slot\otimes\pi\colon \CK\to\CK$
  preserves coreflexive equalizers.
\end{definition}

A \emph{morphism}
of Grothendieck objects in $\kos{K}$
is a morphism of affine group schemes in $\kos{K}$.
Let $f$, $g\colon\Spec(\pi^{\pr})\to\Spec(\pi)$
be morphisms of Grothendieck objects in $\kos{K}$.
A \emph{$2$-isomorphism}
$\theta\colon f\cong g$ is a morphism
$\theta\colon \Spec(\kappa)\to \Spec(\pi)$ in $\cat{Aff}(\kos{K})$
that satisfies
$\cp_{\pi}\circ (\theta\times f)=\cp_{\pi}\circ (g\times \theta)
\colon \Spec(\pi^{\pr})\to\Spec(\pi)$.
Grothendieck objects in $\kos{K}$,
morphisms between them,
and $2$-isomorphisms between these morphisms
form a $(2,1)$-category,
which we denote by
\begin{equation}
  \label{eq-Groobjinpk-pGroObjK21cat}
  \mathsf{pGroObj}(\kos{K}).
\end{equation}
Let $\tilde{g}\colon \Spec(\pi^{\pr})\to\Spec(\pi)$
be a morphism
and $\tilde{\theta}\colon g\cong \tilde{g}$
be a $2$-isomorphism.
The vertical composite of
$\theta$ and $\tilde{\theta}$
is the convolution product
$\tilde{\theta}\star\theta\colon f\cong \tilde{g}$.
The identity $2$-isomorphism of $f$
is $e_{\pi}\colon \Spec(\kappa)\to \Spec(\pi)$.
Each $2$-isomorphism $\theta\colon f\cong g$ is indeed invertible
with respect to $\star$;
its inverse $\theta^{-1}\colon g\cong f$ is
$\varsigma_{\pi}\circ \theta\colon \Spec(\kappa)\to \Spec(\pi)$.
Let $f^{\pr}$, $g^{\pr}\colon \Spec(\pi^{\ppr})\to\Spec(\pi^{\pr})$
be morphisms,
and 
$\theta^{\pr}\colon f^{\pr}\cong g^{\pr}$
be a $2$-isomorphism.
We write
$ff^{\pr}$, $gg^{\pr}\colon \Spec(\pi^{\ppr})\to\Spec(\pi)$
for the composites.
The horizontal composite 
$\theta\theta^{\pr}\colon ff^{\pr}\cong gg^{\pr}$
of $\theta$ and $\theta^{\pr}$ is
$\theta\star (f\circ \theta^{\pr})=(g\circ \theta^{\pr})\star\theta
\colon \Spec(\kappa)\to \Spec(\pi)$.

\section{Internal representations of a Grothendieck object}
\label{sec-RepGroobj}

Let $\Spec(b)$ be an affine scheme in $\kos{K}$.
The endofunctor $\slot\otimes b$ on $\CK$ is a monad,
with product and unit natural transformations
induced by the product and unit morphisms of $b$.
We denote by $\CK_b^{\vphantom{*}}$ the Kleisli category for this monad.
Objects in $\CK_b^{\vphantom{*}}$ are written as $\kos{u}_b^*(x)$,
one for each object $x$ in $\CK$.
A morphism $\kos{u}_b^*(x)\to \kos{u}_b^*(y)$ in $\CK_b^{\vphantom{*}}$
is a morphism $x\to y\otimes b$ in $\CK$.
There is a canonical adjunction between $\CK$ and $\CK_b^{\vphantom{*}}$;
the left adjoint $\kos{u}_b^*\colon \CK\to \CK_b^{\vphantom{*}}$
sends $x$ to $\kos{u}_b^*(x)$,
and the right adjoint $\kos{u}_{b*}^{\vphantom{*}}\colon \CK_b^{\vphantom{*}}\to \CK$
sends $\kos{u}_b^*(x)$ to $x\otimes b$.
The Kleisli category $\CK_b^{\vphantom{*}}$ has a symmetric monoidal category structure,
with monoidal product
$\kos{u}_b^*(x)\otimes_b^{\vphantom{*}}\kos{u}_b^*(y)\coloneqq\kos{u}_b^*(x\otimes y)$
and unit object $\kos{u}_b^*(\kappa)$.
We write 
\begin{equation*}
  \kos{K}_b^{\vphantom{*}}=(\CK_b^{\vphantom{*}},\otimes_b^{\vphantom{*}},\kos{u}_b^*(\kappa))
\end{equation*}
for this symmetric monoidal category.
The left adjoint $\kos{u}_b^*$,
equipped with identity morphisms as its coherence isomorphisms,
is a strong monoidal functor from $\kos{K}$ to $\kos{K}_b^{\vphantom{*}}$.
The adjunction $\kos{u}_b^{\vphantom{*}}\coloneqq(\kos{u}_b^*\dashv \kos{u}_{b*}^{\vphantom{*}})$
thus has the structure of a Grothendieck morphism,
although $\kos{K}_b^{\vphantom{*}}$ is not, in general, a Grothendieck prekosmos.
Nonetheless, we shall denote the adjunction by an arrow
with a filled tip
\begin{equation*}
  \kos{u}_b^{\vphantom{*}}\colon \kos{K}_b^{\vphantom{*}}\mor\kos{K}.
\end{equation*}
The right adjoint part $\kos{u}_{b*}^{\vphantom{*}}$ of $\kos{u}_b^{\vphantom{*}}$ is conservative,
and, again by doctrinal adjunction,
has a unique lax monoidal functor structure
such that the unit and counit
of $\kos{u}_b^{\vphantom{*}}$ are monoidal natural transformations.
We have $\kos{u}_{b*}^{\vphantom{*}}\kos{u}_b^*=\slot\otimes b$
as lax monoidal endofunctors on $\kos{K}$,
where $\slot\otimes b$ is equipped with coherence morphisms
written as $(\otimes b)_{\tensor}$ and $(\otimes b)_{\unit}$,
explicitly given by
$(\otimes b)_{\tensor,x,y}
\coloneqq(I_{x\otimes y}\otimes \pc_b)\circ (I_x\otimes \sigma_{b,y}\otimes I_b)
\colon x\otimes b\otimes y\otimes b\to x\otimes y\otimes b$
and
$(\otimes b)_{\unit}\coloneqq I_{\kappa}\otimes u_b\colon \kappa\to \kappa\otimes b$.
When $b=\kappa$, we identify $\kos{u}_{\kappa}^{\vphantom{*}}\colon \kos{K}_{\kappa}^{\vphantom{*}}\mor\kos{K}$
with $\vphantom{\id}\smash[t]{\id}_{\kos{K}}^{\vphantom{*}}\colon \kos{K}=\kos{K}$.

Let
$f\colon \Spec(b^{\pr})\to \Spec(b)$
be a morphism of affine schemes in $\kos{K}$.
The \emph{base change} along $f$ is the pair
$(\kos{f}^*,\what{\kos{f}^*})$
consisting of the strong monoidal functor
$\kos{f}^*\colon \kos{K}_b^{\vphantom{*}}\to\kos{K}_{b^{\pr}}^{\vphantom{*}}$
that sends $\kos{u}_b^*(x)$ to $\kos{u}_{b^{\pr}}^*(x)$,
and the equality 
$\what{\kos{f}^*}\colon \kos{u}_{b^{\pr}}^*=\kos{f}^*\kos{u}_b^*$
of strong monoidal functors.
We omit $\what{\kos{f}^*}$ from the notation
and denote the pair by $\kos{f}^*\colon \kos{u}_b^*\to \kos{u}_{b^{\pr}}^*$.
The underlying functor 
$\kos{f}^*\colon\CK_b^{\vphantom{*}}\to \CK_{b^{\pr}}^{\vphantom{*}}$
is also referred to as the base change along $f$.
By the base change of a morphism in $\CK_b^{\vphantom{*}}$ along $f$
we mean its image under $\kos{f}^*$,
and we use the same terminology for
natural transformations whose components lie in $\CK_b^{\vphantom{*}}$.
The base change along $u_b\colon \Spec(b)\to \Spec(\kappa)$
is the pair $(\kos{u}_b^*,\what{\kos{u}_b^*})$ consisting of
$\kos{u}_b^*\colon \kos{K}\to \kos{K}_b^{\vphantom{*}}$
and the equality
$\what{\kos{u}_b^*}\colon \kos{u}_b^*=\kos{u}_b^*$,
and we denote the pair by
$\kos{u}_b^*\colon \vphantom{\id}\smash[t]{\id}_{\kos{K}}^{\vphantom{*}}\to \kos{u}_b^*$.
The base change along $u_b$ is also described
in terms of the adjunction $\kos{u}_b^{\vphantom{*}}\colon \kos{K}_b^{\vphantom{*}}\mor\kos{K}$;
it is the pair 
$(\kos{u}_b^{\vphantom{*}},\what{\vphantom{\kos{u}_b}\smash[t]{\kos{u}_b^{\vphantom{*}}}})$
consisting of the adjunction $\kos{u}_b^{\vphantom{*}}$ and the equality 
$\what{\vphantom{\kos{u}_b}\smash[t]{\kos{u}_b^{\vphantom{*}}}}\colon \kos{u}_b^{\vphantom{*}}=\kos{u}_b^{\vphantom{*}}$,
and we denote the pair by an arrow with a filled tip
$\kos{u}_b^{\vphantom{*}}\colon \kos{u}_b^{\vphantom{*}}\mor\vphantom{\id}\smash[t]{\id}_{\kos{K}}^{\vphantom{*}}$.
We illustrate these pairs as follows.
\begin{equation}
  \label{eq-RepGroobj-basechange}
  \begin{tikzcd}[row sep=5pt,column sep=5pt]
    \kos{K}_{b^{\pr}}
    &\text{ }
    &\kos{K}_b
    \ar[ll,"\kos{f}^*"']
    &\text{ }
    &\text{ }
    &\kos{K}_b
    &\text{ }
    &\kos{K}
    \ar[ll,"\kos{u}_b^*"']
    &\text{ }
    &\text{ }
    &\kos{K}_b
    \ar[rr,filled,"\kos{u}_b"]
    \ar[ddr,filled,"\kos{u}_b"']
    &\text{ }
    &\kos{K}
    \ar[ddl,equal]
    \\%1
    \text{ }
    \ar[rr,2cell,shift right=0.5ex,"{\,\,\what{\kos{f}^*}}"]
    &\text{ }
    &\text{ }
    &\text{ }
    &\text{ }
    &\text{ }
    \ar[rr,2cell,shift right=0.5ex,"{\!\!\!\what{\kos{u}_b^*}}"]
    &\text{ }
    &\text{ }
    &\text{ }
    &\text{ }
    &\text{ }
    \ar[rr,2cell,shift right=0.5ex,"{\!\!\!\what{\kos{u}_b}}"]
    &\text{ }
    &\text{ }
    \\%2
    \text{ }
    &\kos{K}
    \ar[uul,"\kos{u}_{b^{\pr}}^*"]
    \ar[uur,"\kos{u}_b^*"']
    &\text{ }
    &\text{ }
    &\text{ }
    &\text{ }
    &\kos{K}
    \ar[uul,"\kos{u}_b^*"]
    \ar[uur,equal]
    &\text{ }
    &\text{ }
    &\text{ }
    &\text{ }
    &\kos{K}
    &\text{ }
  \end{tikzcd}
\end{equation}

Let $x$ be an object in $\CK$.
The assignment 
$\Spec(b)\mapsto \Aut_{\CK_b^{\vphantom{*}}}(\kos{u}_b^*(x))$,
sending each affine scheme $\Spec(b)$ in $\kos{K}$
to the group of automorphisms of $\kos{u}_b^*(x)$ in $\CK_b^{\vphantom{*}}$,
defines a presheaf of groups on $\cat{Aff}(\kos{K})$,
denoted by $\uAut(x)$,
with restriction morphisms given by base change.
Likewise, the assignment
$\Spec(b)\mapsto \End_{\CK_b^{\vphantom{*}}}(\kos{u}_b^*(x))$
defines a presheaf of monoids on $\cat{Aff}(\kos{K})$,
denoted by $\uEnd(x)$.

Let $\Spec(\pi)$ be a Grothendieck object in $\kos{K}$
and $x\in \CK$.
A \emph{functorial action} of $\Spec(\pi)$ on $x$
is a datum $\uprho_x$,
for each affine scheme $\Spec(b)$ in $\kos{K}$
and each morphism
$g\colon \Spec(b)\to \Spec(\pi)$ in $\cat{Aff}(\kos{K})$,
of a morphism
\begin{equation}
  \label{eq-RepGroobj-actionofg}
  \uprho_x(g)\colon \kos{u}_b^*(x)\to \kos{u}_b^*(x)
\end{equation}
in $\CK_b^{\vphantom{*}}$,
and we require that
(i) $\uprho_x(g)$ is compatible with base change
(that is, $\kos{f}^*(\uprho_x(g))=\uprho_x(g\circ f)\colon \kos{u}_{b^{\pr}}^*(x)\to\kos{u}_{b^{\pr}}^*(x)$
for every morphism $f\colon \Spec(b^{\pr})\to \Spec(b)$ in $\cat{Aff}(\kos{K})$);
(ii) $\uprho_x(\tilde{g}\star g)=\uprho_x(\tilde{g})\circ \uprho_x(g)$
for all $\tilde{g}$, $g\colon \Spec(b)\to \Spec(\pi)$;
and
(iii) $\uprho_x(g)$ is the identity morphism of $\kos{u}_b^*(x)$
when $g=e_{\pi}\circ u_b\colon \Spec(b)\to \Spec(\pi)$.
Since every $g$ is invertible with respect to $\star$,
the morphisms $\uprho_x(g)$ are invertible.
We say that $\uprho_x$ is \emph{trivial}
if $\uprho_x(g)$ is the identity morphism of $\kos{u}_b^*(x)$
for every $g$.

\begin{definition}
  \label{def-RepGroobj-Internalrep}
  Let $\Spec(\pi)$ be a Grothendieck object in $\kos{K}$.
  A \emph{representation} of $\Spec(\pi)$ \emph{in $\kos{K}$},
  or an \emph{internal representation} of $\Spec(\pi)$,
  is a pair consisting of an object $x$ in $\CK$
  and a functorial action of $\Spec(\pi)$ on $x$.
  We say that a representation of $\Spec(\pi)$ in $\kos{K}$
  is \emph{trivial} if its functorial action is trivial.
\end{definition}

A functorial action $\uprho_x$ of $\Spec(\pi)$
on $x\in \CK$ is the same as a morphism
$\uprho_x\colon \Hom_{\cat{Aff}(\kos{K})}(\slot,\Spec(\pi))\Rightarrow \uAut(x)$
of presheaves of groups on $\cat{Aff}(\kos{K})$,
and is, by the Yoneda lemma,
uniquely determined by $\uprho_x(g)$ in the universal case
where $g$ is the identity morphism of $\Spec(\pi)$.
We write
\begin{equation*}
  \xi_x\colon \kos{u}_{\pi}^*(x)\cong \kos{u}_{\pi}^*(x)
\end{equation*}
for the universal element,
which is an isomorphism in $\CK_{\pi}^{\vphantom{*}}$.
Since $\uprho_x$ is compatible with the group structures on both sides,
$\xi_x$ must satisfy the following:
(i)
in $\CK_{\pi\otimes\pi}^{\vphantom{*}}$,
the base change of $\xi_x$ along
$\cp_{\pi}\colon \Spec(\pi)\times \Spec(\pi)\to \Spec(\pi)$
is equal to $\kpr_1^*(\xi_x)\circ \kpr_2^*(\xi_x)$,
where each $\kpr_i^*$ is the base change along the $i$-th projection
$\mathrm{pr}_i\colon \Spec(\pi)\times \Spec(\pi)\to\Spec(\pi)$;
(ii)
in $\CK$,
the base change of $\xi_x$ along
$e_{\pi}\colon \Spec(\kappa)\to \Spec(\pi)$
is equal to the identity morphism of $x$.
When we regard the isomorphism $\xi_x$ in $\CK_{\pi}^{\vphantom{*}}$
as a morphism in $\CK$, we denote it by
\begin{equation*}
  \rho_x\colon x\to x\otimes \pi.
\end{equation*}
Then $\rho_x$ must satisfy the right $\pi$-coaction relations,
and the assignment $\uprho_x\mapsto \rho_x$ is a one-to-one correspondence
between functorial actions of $\Spec(\pi)$ on $x$
and right $\pi$-coactions on $x$.
We henceforth identify a representation
$(x,\uprho_x)$ of $\Spec(\pi)$ in $\kos{K}$
with a right $\pi$-comodule $(x,\rho_x)$.
A \emph{morphism} of representations of $\Spec(\pi)$ in $\kos{K}$
is a morphism between the underlying objects in $\CK$
that is equivariant under the functorial actions
of $\Spec(\pi)$, or, equivalently, one that is 
compatible with the corresponding right $\pi$-coactions.
We denote by $\cat{Rep}_{\kos{K}}(\pi)$
the category of representations of $\Spec(\pi)$ in $\kos{K}$.

\begin{lemma}
  \label{lem-RepGroobj-Reppi}
  Let $\Spec(\pi)$ be a Grothendieck object in $\kos{K}$.
  \begin{enumerate}
    \item 
    The representations of $\Spec(\pi)$ in $\kos{K}$
    form a Grothendieck prekosmos, denoted by
    $\kosRep_{\kos{K}}(\pi)=(\cat{Rep}_{\kos{K}}(\pi),\tensor_{\pi},\unit_{\pi})$.
    
    \item
    The Grothendieck prekosmos 
    $\kosRep_{\kos{K}}(\pi)$
    is equipped with a connected and surjective Grothendieck morphism
    to $\kos{K}$,
    denoted by $\kos{t}_{\pi}^{\vphantom{*}}\colon \kosRep_{\kos{K}}(\pi)\mor\kos{K}$,
    whose left adjoint part $\kos{t}_{\pi}^*$
    sends each $z\in \CK$
    to the trivial representation of $\Spec(\pi)$ in $\kos{K}$
    with underlying object $z$.
    
    \item 
    The Grothendieck $\kos{K}$-prekosmos
    $\kosRep_{\kos{K}}(\pi)\xmor{\kos{t}_{\pi}^{\vphantom{*}}}\kos{K}$
    is equipped with a surjective fiber point,
    denoted by
    $\kos{e}_{\pi}^{\vphantom{*}}\colon \vphantom{\id}\smash[t]{\id}_{\kos{K}}^{\vphantom{*}}\mor\kos{t}_{\pi}^{\vphantom{*}}$,
    whose left adjoint part $\kos{e}_{\pi}^*$
    sends each $(x,\rho_x)\in\cat{Rep}_{\kos{K}}(\pi)$
    to its underlying object $x$ in $\CK$.
  \end{enumerate}
  \begin{equation*}
    \begin{tikzcd}[row sep=5pt,column sep=4pt]
      \kos{K}
      \ar[rr,filled,"\kos{e}_{\pi}"]
      \ar[ddr,equal]
      &\text{ }
      &\kosRep_{\kos{K}}(\pi)
      \ar[ddl,filled,"\kos{t}_{\pi}"]
      \\%1
      \text{ }
      \ar[rr,2cell,pos=0.38,shift right=0.3ex,"{\!\!\what{\kos{e}_{\pi}}}"]
      &\text{ }
      &\text{ }
      \\%2
      \text{ }
      &\kos{K}
      &\text{ }
    \end{tikzcd}
  \end{equation*}
\end{lemma}
\begin{proof}
  We describe the structures in detail.
  The monoidal product $\tensor_{\pi}$
  sends each pair $(x,\rho_x)$, $(y,\rho_y)\in\cat{Rep}_{\kos{K}}(\pi)$
  to $(x\otimes y,\rho_{x\otimes y})$,
  where 
  $\rho_{x\otimes y}\coloneqq(I_{x\otimes y}\otimes \pc_{\pi})
  \circ (I_x\otimes \sigma_{\pi,y}\otimes I_{\pi})
  \circ (\rho_x\otimes \rho_y)
  \colon x\otimes y\to x\otimes y\otimes \pi$.
  The unit object is
  $\unit_{\pi}\coloneqq(\kappa,I_{\kappa}\otimes u_{\pi})$.
  The left adjoint part $\kos{t}_{\pi}^*$
  sends each $z\in \CK$ to $(z,I_z\otimes u_{\pi})$,
  and the right adjoint part $\kos{t}_{\pi*}^{\vphantom{*}}$
  sends each $(x,\rho_x)\in\cat{Rep}_{\kos{K}}(\pi)$
  to the equalizer in $\CK$ of the coreflexive pair
  $\rho_x$, $I_x\otimes u_{\pi}\colon x\to x\otimes \pi$,
  with common retraction $I_x\otimes e_{\pi}$.
  The left adjoint part $\kos{e}_{\pi}^*$
  sends each $(x,\rho_x)\in\cat{Rep}_{\kos{K}}(\pi)$ to $x$,
  and the right adjoint part $\kos{e}_{\pi*}^{\vphantom{*}}$
  sends each $z\in \CK$ to $(z\otimes \pi,I_z\otimes \cp_{\pi})$.
  The unit and counit of $\kos{e}_{\pi}^{\vphantom{*}}$ are given by
  $\eta_{(x,\rho_x)}=\rho_x\colon (x,\rho_x)\to \kos{e}_{\pi*}^{\vphantom{*}}\kos{e}_{\pi}^*(x,\rho_x)$
  and $\epsilon_z=I_z\otimes e_{\pi}\colon \kos{e}_{\pi}^*\kos{e}_{\pi*}^{\vphantom{*}}(z)\to z$.
  The left adjoint parts
  $\kos{t}_{\pi}^*$ and $\kos{e}_{\pi}^*$
  are strong monoidal functors,
  whose coherence isomorphisms are identity morphisms.
  The invertible Grothendieck transformation
  $\what{\vphantom{\kos{e}_{\pi}}\smash[t]{\kos{e}_{\pi}^{\vphantom{*}}}}\colon \vphantom{\id}\smash[t]{\id}_{\kos{K}}^{\vphantom{*}}\cong \kos{t}_{\pi}^{\vphantom{*}}\kos{e}_{\pi}^{\vphantom{*}}$
  is the equality
  $\what{\kos{e}_{\pi}^*}\colon \vphantom{\id}\smash[t]{\id}_{\kos{K}}^{\vphantom{*}}=\kos{e}_{\pi}^*\kos{t}_{\pi}^*$
  between the left adjoint parts.
  It is straightforward to check that
  $\kosRep_{\kos{K}}(\pi)$ is a Grothendieck prekosmos,
  that $\kos{t}_{\pi}^{\vphantom{*}}\colon \kosRep_{\kos{K}}(\pi)\mor\kos{K}$
  is a connected Grothendieck morphism,
  and that $\kos{e}_{\pi}^{\vphantom{*}}\colon \vphantom{\id}\smash[t]{\id}_{\kos{K}}^{\vphantom{*}}\mor\kos{t}_{\pi}^{\vphantom{*}}$
  is a surjective fiber point of $\kos{t}_{\pi}^{\vphantom{*}}$.
  The composite $\kos{t}_{\pi}^{\vphantom{*}}\kos{e}_{\pi}^{\vphantom{*}}$ is surjective,
  as it is isomorphic to $\vphantom{\id}\smash[t]{\id}_{\kos{K}}^{\vphantom{*}}$.
  Since $\kos{e}_{\pi}^{\vphantom{*}}$ is surjective,
  we deduce that $\kos{t}_{\pi}^{\vphantom{*}}$ is surjective.
\end{proof}

%% file: Grocatoverpk.tex
\chapter{Grothendieck Categories over a Grothendieck Prekosmos}
\label{chap-Grocatoverpk}

Throughout this chapter,
we fix a Grothendieck prekosmos $\kos{K}=(\CK,\otimes,\kappa)$.
Recall from Definition~\ref{def-Groprekos-Grocatoverpk}
the definition of a Grothendieck $\kos{K}$-category,
and from Lemma~\ref{lem-RepGroobj-Reppi}
the Grothendieck $\kos{K}$-prekosmos
$\kosRep_{\kos{K}}(\pi)\xmor{\kos{t}_{\pi}^{\vphantom{*}}}\kos{K}$
of internal representations of
a Grothendieck object $\Spec(\pi)$ in $\kos{K}$.
The primary objective of Part~\ref{part-GroContext}
is to prove the following theorem.

\begin{theorem}
  \label{thm-Grocatoverpk-PrimaryObj}
  A Grothendieck $\kos{K}$-prekosmos
  is a Grothendieck $\kos{K}$-category
  if and only if it is equivalent
  to $\kosRep_{\kos{K}}(\pi)\xmor{\kos{t}_{\pi}^{\vphantom{*}}}\kos{K}$
  as Grothendieck $\kos{K}$-prekosmoi
  for some Grothendieck object $\Spec(\pi)$ in $\kos{K}$.
\end{theorem}

For each Grothendieck object $\Spec(\pi)$ in $\kos{K}$,
$\kos{t}_{\pi}^{\vphantom{*}}$ admits a surjective fiber point
$\kos{e}_{\pi}^{\vphantom{*}}\colon\vphantom{\id}\smash[t]{\id}_{\kos{K}}^{\vphantom{*}}\mor\kos{t}_{\pi}^{\vphantom{*}}$
by Lemma~\ref{lem-RepGroobj-Reppi}.
Thus if a Grothendieck $\kos{K}$-prekosmos
is equivalent to $\kos{t}_{\pi}^{\vphantom{*}}$
as Grothendieck $\kos{K}$-prekosmoi
for some $\Spec(\pi)$,
then it is a Grothendieck $\kos{K}$-category.
To prove the converse,
we consider the automorphisms of
a chosen surjective fiber point.

Let $\kos{T}\xmor{\kos{t}}\kos{K}$
be a Grothendieck prekosmos over $\kos{K}$
and
$\omega\colon\vphantom{\id}\smash[t]{\id}_{\kos{K}}^{\vphantom{*}}\mor\kos{t}$
be a Grothendieck $\kos{K}$-morphism.
For each $\Spec(b)\in\cat{Aff}(\kos{K})$,
recall the base change
$\kos{u}_b^{\vphantom{*}}\colon\kos{u}_b^{\vphantom{*}}\mor\vphantom{\id}\smash[t]{\id}_{\kos{K}}^{\vphantom{*}}$
along $u_b\colon\Spec(b)\to\Spec(\kappa)$
in \eqref{eq-RepGroobj-basechange}.
We denote by
$\omega\kos{u}_b^{\vphantom{*}}\colon\kos{u}_b^{\vphantom{*}}\mor\kos{t}$
the composite of $\omega$ with $\kos{u}_b^{\vphantom{*}}$,
and write
$\what{\kos{u}_b^*\omega^*}
=\kos{u}_b^*\what{\omega^*}
\colon\kos{u}_b^*\cong \kos{u}_b^*\omega^*\kos{t}^*$.
Let $\omega^{\pr}\colon\vphantom{\id}\smash[t]{\id}_{\kos{K}}^{\vphantom{*}}\mor\kos{t}$
also be a Grothendieck $\kos{K}$-morphism.
Although $\kos{K}_b^{\vphantom{*}}$ need not be a Grothendieck prekosmos,
by a \emph{Grothendieck $\kos{K}$-transformation}
$\vartheta\colon\omega^{\pr}\!\kos{u}_b^{\vphantom{*}}\Rightarrow\omega\kos{u}_b^{\vphantom{*}}$
we mean a monoidal natural transformation
$\vartheta\colon\kos{u}_b^*\omega^{\pr*}\Rightarrow \kos{u}_b^*\omega^*$
that satisfies
$\what{\kos{u}_b^*\omega^*}=(\vartheta\kos{t}^*)\circ \what{\kos{u}_b^*\omega^{\pr*}}$.
We say that $\vartheta$ is \emph{invertible}
if every component of $\vartheta$ is an isomorphism.
We illustrate
$\omega\kos{u}_b^{\vphantom{*}}$, $\omega^{\pr}\!\kos{u}_b^{\vphantom{*}}$,
and $\vartheta$ as follows.
\begin{equation*}
  \begin{tikzcd}[row sep=5pt,column sep=15pt]
    \kos{K}_b
    \ar[r,filled,"\kos{u}_b"]
    \ar[ddr,filled,"\kos{u}_b"']
    &\kos{K}
    \ar[dd,equal]
    \ar[r,filled,"\omega"]
    &\kos{T}
    \ar[ddl,filled,"\kos{t}"]
    &\kos{K}_b
    \ar[r,filled,"\kos{u}_b"]
    \ar[ddr,filled,"\kos{u}_b"']
    &\kos{K}
    \ar[dd,equal]
    \ar[r,filled,"\omega^{\pr}"]
    &\kos{T}
    \ar[ddl,filled,"\kos{t}"]
    &\kos{K}_b
    \ar[ddr,filled,"\kos{u}_b"']
    \ar[rr,filled,bend left=23pt,"\omega\kos{u}_b" description]
    \ar[rr,filled,bend right=23pt,"\omega^{\pr}\!\kos{u}_b" description]
    &\text{ }
    &\kos{T}
    \ar[ddl,filled,"\kos{t}"]
    \\%1
    \text{ }
    \ar[r,2cell,pos=0.9,"{\!\!\!\what{\kos{u}_b}}"]
    &\text{ }
    \ar[r,2cell,pos=0.1,"{\!\what{\omega}}"]
    &\text{ }
    &\text{ }
    \ar[r,2cell,pos=0.9,"{\!\!\!\what{\kos{u}_b}}"]
    &\text{ }
    \ar[r,2cell,pos=0.1,"{\,\what{\omega^{\pr}}}"]
    &\text{ }
    &\text{ }
    &\text{ }
    \ar[u,2cell,pos=1.9,"{\vartheta}"']
    &\text{ }
    \\%2
    \text{ }
    &\kos{K}
    &\text{ }
    &\text{ }
    &\kos{K}
    &\text{ }
    &\text{ }
    &\kos{K}
    &\text{ }
  \end{tikzcd}
\end{equation*}
The assignment sending each $\Spec(b)$
to the set of Grothendieck $\kos{K}$-transformations
$\omega^{\pr}\!\kos{u}_b^{\vphantom{*}}\Rightarrow \omega\kos{u}_b^{\vphantom{*}}$
defines a presheaf on $\cat{Aff}(\kos{K})$
with restriction morphisms given by base change;
we call it the 
\emph{presheaf on $\cat{Aff}(\kos{K})$
of Grothendieck $\kos{K}$-transformations}
from $\omega^{\pr}$ to $\omega$,
and denote it by $\uHom_{\kos{K}}(\omega^{\pr},\omega)$.
Likewise, we define the 
\emph{presheaf on $\cat{Aff}(\kos{K})$
of invertible Grothendieck $\kos{K}$-transformations}
from $\omega^{\pr}$ to $\omega$,
denoted by $\uIsom_{\kos{K}}(\omega^{\pr},\omega)$,
as the presheaf sending each $\Spec(b)$
to the set of invertible Grothendieck $\kos{K}$-transformations
$\omega^{\pr}\!\kos{u}_b^{\vphantom{*}}\cong\omega\kos{u}_b^{\vphantom{*}}$.
When $\omega^{\pr}=\omega$,
we define the 
\emph{presheaf of groups on $\cat{Aff}(\kos{K})$
of invertible Grothendieck $\kos{K}$-transformations}
from $\omega$ to itself,
denoted by $\uAut_{\kos{K}}(\omega)$,
as the presheaf sending each $\Spec(b)$
to the group of invertible Grothendieck $\kos{K}$-transformations
$\omega\kos{u}_b^{\vphantom{*}}\cong\omega\kos{u}_b^{\vphantom{*}}$.
We also define the
\emph{presheaf of monoids on $\cat{Aff}(\kos{K})$
of Grothendieck $\kos{K}$-transformations}
from $\omega$ to itself,
denoted by $\uEnd_{\kos{K}}(\omega)$,
as the presheaf sending each $\Spec(b)$
to the monoid of Grothendieck $\kos{K}$-transformations
$\omega\kos{u}_b^{\vphantom{*}}\Rightarrow\omega\kos{u}_b^{\vphantom{*}}$.

The following is the main theorem of Chapter~\ref{chap-Grocatoverpk}.

\begin{theorem}
  \label{thm-Grocatoverpk-mainthm}
  Let $\kos{T}\xmor{\kos{t}}\kos{K}$
  be a Grothendieck $\kos{K}$-prekosmos
  and $\omega\colon\vphantom{\id}\smash[t]{\id}_{\kos{K}}^{\vphantom{*}}\mor\kos{t}$
  be a surjective fiber point of $\kos{t}$.
  \begin{enumerate}
    \item 
    The presheaf of groups $\uAut_{\kos{K}}(\omega)$ on $\cat{Aff}(\kos{K})$
    of invertible Grothendieck $\kos{K}$-transformations
    from $\omega$ to itself
    is represented by a Grothendieck object
    $\Spec(\pi)$ in $\kos{K}$,
    with $\pi\in\cat{Comm}(\kos{K})$
    given by $\omega^*\omega_*(\kappa)$.
    
    \item
    The fiber point $\omega$ extends to an equivalence
    $\breve{\omega}\colon\kos{t}_{\pi}^{\vphantom{*}}\simeq\kos{t}$
    of Grothendieck $\kos{K}$-prekosmoi
    up to an invertible Grothendieck $\kos{K}$-transformation
    $\omega\cong \breve{\omega}\kos{e}_{\pi}^{\vphantom{*}}$
    given by the equality
    $\omega^*=\kos{e}_{\pi}^*\breve{\omega}^*$
    between the left adjoint parts.
  \end{enumerate}
  \begin{equation*}
    \Hom_{\cat{Aff}(\kos{K})}(\slot,\Spec(\pi))
    \cong
    \uAut_{\kos{K}}(\omega)
    \qquad\quad
    \begin{tikzcd}[row sep=5pt,column sep=6pt]
      \text{ }
      &\id_{\kos{K}}
      \ar[ddr,filled,"{\omega}"]
      \ar[ddl,filled,"{\kos{e}_{\pi}}"']
      &\text{ }
      \\%1
      \text{ }
      &\text{ }
      &\text{ }
      \ar[ll,2cell,shift left=0.5ex,pos=0.475]
      \\%2
      \kos{t}_{\pi}
      \ar[rr,filled,"{\breve{\omega}}","\simeq"']
      &\text{ }
      &\kos{t}
    \end{tikzcd}
  \end{equation*}
\end{theorem}

Theorem~\ref{thm-Grocatoverpk-PrimaryObj} follows immediately
once we prove Theorem~\ref{thm-Grocatoverpk-mainthm}.

\begin{proof}[Proof of Theorem~\ref{thm-Grocatoverpk-PrimaryObj}]
  Let $\kos{T}\xmor{\kos{t}}\kos{K}$
  be a Grothendieck $\kos{K}$-prekosmos.
  If an equivalence
  $\breve{\omega}\colon\kos{t}_{\pi}^{\vphantom{*}}\simeq\kos{t}$
  is given for some
  Grothendieck object $\Spec(\pi)$ in $\kos{K}$,
  then $\kos{t}$ is a Grothendieck $\kos{K}$-category,
  since it admits a surjective fiber point
  $\breve{\omega}\kos{e}_{\pi}^{\vphantom{*}}\colon\vphantom{\id}\smash[t]{\id}_{\kos{K}}^{\vphantom{*}}\mor\kos{t}$
  by Lemma~\ref{lem-RepGroobj-Reppi}.
  Conversely, if $\kos{t}$ is a Grothendieck $\kos{K}$-category,
  we may choose a surjective fiber point $\omega$ of $\kos{t}$,
  and by Theorem~\ref{thm-Grocatoverpk-mainthm},
  we obtain the Grothendieck object $\Spec(\pi)$
  in $\kos{K}$ that represents $\uAut_{\kos{K}}(\omega)$,
  as well as the equivalence $\breve{\omega}\colon\kos{t}_{\pi}^{\vphantom{*}}\simeq\kos{t}$.
\end{proof}

We shall prove Theorem~\ref{thm-Grocatoverpk-mainthm}
in Section~\ref{sec-GroMainthm}.
The proof is based on the following proposition.

\begin{proposition}\label{prop-Grocatoverpk-HomisIsom}
  Let $\kos{T}\xmor{\kos{t}}\kos{K}$ be a Grothendieck $\kos{K}$-prekosmos
  and let $\omega$, $\omega^{\pr}\colon\vphantom{\id}\smash[t]{\id}_{\kos{K}}^{\vphantom{*}}\mor\kos{t}$
  be Grothendieck $\kos{K}$-morphisms.
  Put $p\coloneqq\omega^{\pr*}\omega_*(\kappa)\in \cat{Comm}(\kos{K})$.
  Assume that $\omega$ is a fiber point of $\kos{t}$.
  Then the following hold:
  \begin{enumerate}
    \item
    We have a monoidal natural isomorphism
    $\slot\otimes p\cong \omega^{\pr*}\omega_*$
    between lax monoidal endofunctors on $\kos{K}$.
    
    \item 
    The presheaf 
    $\uHom_{\kos{K}}(\omega^{\pr},\omega)$
    on $\cat{Aff}(\kos{K})$
    is represented by $\Spec(p)$.

    \item
    We have an equality
    $\uHom_{\kos{K}}(\omega^{\pr},\omega)
    =\uIsom_{\kos{K}}(\omega^{\pr},\omega)$
    between presheaves on $\cat{Aff}(\kos{K})$.
  \end{enumerate}
\end{proposition}

We will prove Proposition~\ref{prop-Grocatoverpk-HomisIsom}
in Section~\ref{sec-ActionbaseGro}.
We present two direct consequences of our results.

\begin{corollary}
  \label{cor-Grocatoverpk-connsurj}
  Every Grothendieck $\kos{K}$-category is connected and surjective.
\end{corollary}
\begin{proof}
  This is a direct consequence of
  Lemma~\ref{lem-RepGroobj-Reppi}
  and Theorem~\ref{thm-Grocatoverpk-PrimaryObj}.
\end{proof}

\begin{corollary}
  \label{cor-Grocatoverpk-pFibgroupoid}
  Let $\kos{T}\xmor{\kos{t}}\kos{K}$ be a Grothendieck $\kos{K}$-prekosmos.
  The category $\cat{pFib}_{\kos{K}}(\kos{t})$
  of fiber points of $\kos{t}$ is a groupoid.
\end{corollary}
\begin{proof}
  This is a direct consequence of
  Proposition~\ref{prop-Grocatoverpk-HomisIsom}.
\end{proof}

Let $\kos{T}\xmor{\kos{t}}\kos{K}$
be a Grothendieck $\kos{K}$-category.
In Section~\ref{sec-pFibGroTors},
we show that once a surjective fiber point
$\omega$ of $\kos{t}$ is chosen,
there is a perfect correspondence
between fiber points of $\kos{t}$
and fpqc $\Spec(\pi)$-torsors over $\Spec(\kappa)$,
where $\Spec(\pi)$ is the Grothendieck object in $\kos{K}$
determined by $\omega$ in Theorem~\ref{thm-Grocatoverpk-mainthm}.
Here we introduce properties of fiber points
that are relevant to this topic.

\begin{definition}
  \label{def-Grocatoverpk-fpqccovering}
  We say that an affine scheme $\Spec(b)$ in $\kos{K}$
  is an \emph{fpqc covering of $\Spec(\kappa)$}
  if $\slot\otimes b\colon\CK\to \CK$
  is conservative and preserves coreflexive equalizers.
\end{definition}

Note that an affine scheme $\Spec(b)$ in $\kos{K}$
is an fpqc covering of $\Spec(\kappa)$
if and only if the base change $\kos{u}_b^*\colon\CK\to \CK_b^{\vphantom{*}}$
along $u_b\colon\Spec(b)\to\Spec(\kappa)$
is conservative and preserves coreflexive equalizers.

\begin{definition}
  \label{def-Grocatoverpk-locisom}
  Let $\kos{T}\xmor{\kos{t}}\kos{K}$
  be a Grothendieck $\kos{K}$-prekosmos
  and let
  $\omega$, $\omega^{\pr}\colon\vphantom{\id}\smash[t]{\id}_{\kos{K}}^{\vphantom{*}}\mor\kos{t}$
  be Grothendieck $\kos{K}$-morphisms.
  We say that $\omega$ is \emph{locally isomorphic} to $\omega^{\pr}$
  if the presheaf
  $\uIsom_{\kos{K}}(\omega,\omega^{\pr})$
  evaluates to a nonempty set
  at some affine scheme $\Spec(b)$ in $\kos{K}$
  that is an fpqc covering of $\Spec(\kappa)$.
  In this case, $\omega^{\pr}$ is locally isomorphic to $\omega$,
  and we also say that
  $\omega$ and $\omega^{\pr}$ are locally isomorphic.
\end{definition}

\begin{lemma}
  \label{lem-Grocatoverpk-pFiblocallyisom}
  Let $\kos{T}\xmor{\kos{t}}\kos{K}$
  be a Grothendieck $\kos{K}$-prekosmos
  and $\omega$, $\omega^{\pr}\colon\vphantom{\id}\smash[t]{\id}_{\kos{K}}^{\vphantom{*}}\mor\kos{t}$
  be fiber points of $\kos{t}$.
  Assume that $\omega$ is surjective.
  Then the following hold:
  \begin{enumerate}
    \item 
    The fiber points $\omega$ and $\omega^{\pr}$ are locally isomorphic.
    
    \item
    The fiber point $\omega^{\pr}$ is also surjective.
  \end{enumerate}
\end{lemma}
\begin{proof}
  Put $p^{\pr}\coloneqq\omega^*\omega^{\pr}_*(\kappa)\in\cat{Comm}(\kos{K})$.
  Since $\omega^{\pr}$ is a fiber point of $\kos{t}$,
  we exchange the roles of $\omega$ and $\omega^{\pr}$ in
  Proposition~\ref{prop-Grocatoverpk-HomisIsom}
  and deduce that $\slot\otimes p^{\pr}\cong \omega^*\omega^{\pr}_*$
  and that the presheaf
  $\uHom_{\kos{K}}(\omega,\omega^{\pr})
  =\uIsom_{\kos{K}}(\omega,\omega^{\pr})$
  is represented by $\Spec(p^{\pr})$.
  The evaluation of the presheaf at $\Spec(p^{\pr})$ is nonempty,
  as it contains the element
  $\omega\kos{u}_{p^{\pr}}^{\vphantom{*}}\cong \omega^{\pr}\!\kos{u}_{p^{\pr}}^{\vphantom{*}}$
  corresponding to the identity morphism of $\Spec(p^{\pr})$.
  For the moment, let us say that a functor has property $\prpty{Q}$
  if it is conservative and preserves coreflexive equalizers.
  Then $\omega^{\pr}_*$ has property $\prpty{Q}$.
  Since $\omega$ is surjective,
  $\omega^*$ also has property $\prpty{Q}$.
  Hence the composite $\omega^*\omega^{\pr}_*$ has property $\prpty{Q}$,
  and therefore $\Spec(p^{\pr})$ is an fpqc covering of $\Spec(\kappa)$.
  This proves (1).
  Moreover, the functor $\kos{u}_{p^{\pr}}^*$ has property $\prpty{Q}$,
  and hence so does the composite $\kos{u}_{p^{\pr}}^*\omega^*$.
  From $\omega\kos{u}_{p^{\pr}}^{\vphantom{*}}\cong \omega^{\pr}\!\kos{u}_{p^{\pr}}^{\vphantom{*}}$
  we deduce that
  $\kos{u}_{p^{\pr}}^*\omega^{\pr*}$ has property $\prpty{Q}$.
  Since $\CK$ has coreflexive equalizers
  and $\kos{u}_{p^{\pr}}^*$ has property $\prpty{Q}$,
  we conclude that $\omega^{\pr*}$ has property $\prpty{Q}$.
  This proves (2).
\end{proof}

\begin{corollary}
  \label{cor-Grocatoverpk-pFibsurjective}
  Let $\kos{T}\xmor{\kos{t}}\kos{K}$
  be a Grothendieck $\kos{K}$-category.
  Every fiber point of $\kos{t}$
  is surjective.
\end{corollary}
\begin{proof}
  This is a direct consequence of 
  Lemma~\ref{lem-Grocatoverpk-pFiblocallyisom}.
\end{proof}

\section{Action of the base Grothendieck prekosmos}
\label{sec-ActionbaseGro}

In this section, we prove Proposition~\ref{prop-Grocatoverpk-HomisIsom}.
We utilize the facts that
(i) a Grothendieck $\kos{K}$-prekosmos
is equipped with the induced action of $\kos{K}$;
(ii) the left and right adjoint parts
of a Grothendieck $\kos{K}$-morphism
are equipped with equivariances
under the induced $\kos{K}$-actions;
and
(iii)
a Grothendieck $\kos{K}$-transformation
between Grothendieck $\kos{K}$-morphisms
is precisely a Grothendieck transformation
between the underlying Grothendieck morphisms
that is compatible with the $\kos{K}$-equivariances
of the left adjoint parts.
These facts are systematically
studied in a broader context
throughout Appendix~\ref{appdix-ltcatoversmc}.
In what follows, we present only what we need.

Let $\kos{T}\xmor{\kos{t}}\kos{K}$
be a Grothendieck $\kos{K}$-prekosmos.
By the \emph{induced $\kos{K}$-action} on $\kos{T}$
we mean the functor $\acts\colon\CK\times\CT\to\CT$
sending each pair $(z,X)$ to $z\acts X\coloneqq\kos{t}^*(z)\tensor X$.
Let $\kos{T}^{\pr}\xmor{\kos{t}^{\pr}}\kos{K}$
be a Grothendieck $\kos{K}$-prekosmos,
with the induced $\kos{K}$-action $\acts^{\pr}$ on $\kos{T}^{\pr}$,
and let $\kos{f}\colon\kos{t}^{\pr}\mor\kos{t}$
be a Grothendieck $\kos{K}$-morphism.
We denote by 
$\what{\kos{f}_*}\colon\kos{t}^*\Rightarrow \kos{f}_*\kos{t}^{\pr*}$
the mate of $(\what{\kos{f}^*})^{-1}\colon\kos{f}^*\kos{t}^*\cong \kos{t}^{\pr*}$
with respect to the adjunctions
$\kos{f}$ and $\vphantom{\id}\smash[t]{\id}_{\kos{K}}^{\vphantom{*}}$.
We illustrate $\kos{f}\colon\kos{t}^{\pr}\mor\kos{t}$,
$\what{\kos{f}^*}$, and $\what{\kos{f}_*}$ as follows.
\begin{equation*}
  \begin{tikzcd}[row sep=5pt,column sep=5pt]
    \kos{T}^{\pr}
    \ar[rr,filled,"\kos{f}"]
    \ar[ddr,filled,"\kos{t}^{\pr}"']
    &\text{ }
    &\kos{T}
    \ar[ddl,filled,"\kos{t}"]
    &\text{ }
    &\text{ }
    &\kos{T}^{\pr}
    &\text{ }
    &\kos{T}
    \ar[ll,"\kos{f}^*"']
    &\text{ }
    &\text{ }
    &\kos{T}^{\pr}
    \ar[rr,"\kos{f}_*"]
    &\text{ }
    &\kos{T}
    \\%1
    \text{ }
    \ar[rr,2cell,shift right=0.5ex,"{\what{\kos{f}}}"]
    &\text{ }
    &\text{ }
    &\text{ }
    &\text{ }
    &\text{ }
    \ar[rr,2cell,shift right=0.5ex,"{\,\,\what{\kos{f}^*}}"]
    &\text{ }
    &\text{ }
    &\text{ }
    &\text{ }
    &\text{ }
    &\text{ }
    &\text{ }
    \ar[ll,2cell,shift left=0.5ex,pos=0.46,"{\what{\kos{f}_*}}"']
    \\%2
    \text{ }
    &\kos{K}
    &\text{ }
    &\text{ }
    &\text{ }
    &\text{ }
    &\kos{K}
    \ar[uul,"\kos{t}^{\pr*}"]
    \ar[uur,"\kos{t}^*"']
    &\text{ }
    &\text{ }
    &\text{ }
    &\text{ }
    &\kos{K}
    \ar[uul,"\kos{t}^{\pr*}"]
    \ar[uur,"\kos{t}^*"']
    &\text{ }
  \end{tikzcd}
\end{equation*}
The \emph{$\kos{K}$-equivariance} of $\kos{f}^*$
is the natural isomorphism
$\cevar{\kos{f}^*}\colon\slot\acts^{\pr}\kos{f}^*(\slot)\cong \kos{f}^*(\slot\acts\slot)$
between functors from $\CK\times \CT$ to $\CT^{\pr}$
whose component at $z\in \CK$, $X\in \CT$ is
\begin{equation*}
  \xymatrix@C=50pt{
    \kos{t}^{\pr*}(z)\tensor^{\pr}\kos{f}^*(X)
    \ar[r]^-{(\what{\kos{f}^*})_z\tensor^{\pr}I_{\kos{f}^*(X)}}_-{\cong}
    &\kos{f}^*\kos{t}^*(z)\tensor^{\pr}\kos{f}^*(X)
    \ar[r]^-{(\kos{f}^*_{\tensor})_{\kos{t}^*(z),X}}_-{\cong}
    &\kos{f}^*(\kos{t}^*(z)\tensor X)
    .
  }
\end{equation*}
The \emph{$\kos{K}$-equivariance} of $\kos{f}_*$
is the natural transformation
$\cevar{\kos{f}_*}\colon\slot\acts\kos{f}_*(\slot)\Rightarrow \kos{f}_*(\slot\acts^{\pr}\slot)$
between functors from $\CK\times \CT^{\pr}$ to $\CT$
whose component at $z\in \CK$, $X^{\pr}\in \CT^{\pr}$ is
\begin{equation*}
  \xymatrix@C=50pt{
    \kos{t}^*(z)\tensor \kos{f}_*(X^{\pr})
    \ar[r]^-{(\what{\kos{f}_*})_z\tensor I_{\kos{f}_*(X^{\pr})}}
    &\kos{f}_*\kos{t}^{\pr*}(z)\tensor \kos{f}_*(X^{\pr})
    \ar[r]^-{(\kos{f}_{*\tensor})_{\kos{t}^{\pr*}(z),X^{\pr}}}
    &\kos{f}_*(\kos{t}^{\pr*}(z)\tensor^{\pr} X^{\pr})
    .
  }
\end{equation*}
  
\begin{lemma}
  \label{lem-ActionbaseGro-projimplystrongequiv}
  Let $\kos{f}\colon\kos{t}^{\pr}\mor\kos{t}$
  be a Grothendieck $\kos{K}$-morphism.
  Recall the natural transformation
  $\varphi$ associated to $\kos{f}$ in \eqref{eq-Groprekos-projformula}.
  \begin{enumerate}
    \item 
    The component of $\cevar{\kos{f}_*}$
    at $z\in \CK$, $X^{\pr}\in\CT^{\pr}$
    can be redescribed as
    \begin{equation*}
      \xymatrix@C=15pt{
        \kos{t}^*(z)\tensor\kos{f}_*(X^{\pr})
        \ar[rr]^-{\varphi_{\kos{t}^*(z),X^{\pr}}}
        &\text{ }
        &\kos{f}_*(\kos{f}^*\kos{t}^*(z)\tensor^{\pr}X^{\pr})
        \ar[rrr]^-{\kos{f}_*((\what{\kos{f}^*})_z^{-1}\tensor^{\pr} I_{X^{\pr}})}_-{\cong}
        &\text{ }
        &\text{ }
        &\kos{f}_*(\kos{t}^{\pr*}(z)\tensor^{\pr}X^{\pr})
        .
      }
    \end{equation*}

    \item
    If $\kos{f}$ is locally connected, then $\cevar{\kos{f}_*}$ is a natural isomorphism.
  \end{enumerate}
\end{lemma}
\begin{proof}
  (1) follows by substituting 
  $\what{\kos{f}_*}=(\kos{f}_*(\what{\kos{f}^*})^{-1})\circ (\eta\kos{t}^*)
  \colon\kos{t}^*\Rightarrow \kos{f}_*\kos{t}^{\pr*}$,
  and (2) is a direct consequence of (1).
\end{proof}

We now show that 
a Grothendieck $\kos{K}$-transformation
between Grothendieck $\kos{K}$-morphisms
is precisely a Grothendieck transformation
between the underlying Grothendieck morphisms
that is compatible with the $\kos{K}$-equivariances
of the left adjoint parts.
In fact, we show a slightly more general statement.

Let $\kos{T}=(\CT,\tensor,\unit)$,
$\kos{T}^{\pr}=(\CT^{\pr},\tensor^{\pr},\unit^{\pr})$
be symmetric monoidal categories
and $\monbase{t}\colon\kos{K}\to \kos{T}$,
$\monbase{t}^{\pr}\colon\kos{K}\to \kos{T}^{\pr}$
be strong monoidal functors.
By the \emph{induced $\kos{K}$-action} $\acts$ on $\kos{T}$
we mean the functor given by 
$z\acts X\coloneqq\monbase{t}(z)\tensor X$,
and we write $\acts^{\pr}$
for the induced $\kos{K}$-action on $\kos{T}^{\pr}$.
We introduce the following terminology.
\begin{itemize}
  \item 
  A \emph{lax $\kos{K}$-tensor functor}
  from $\monbase{t}^{\pr}$ to $\monbase{t}$ is a pair
  $(\mon{f},\what{\mon{f}})$
  consisting of a lax monoidal functor
  $\mon{f}\colon\kos{T}^{\pr}\to \kos{T}$
  and a monoidal natural transformation
  $\what{\mon{f}}\colon\monbase{t}\Rightarrow\mon{f}\monbase{t}^{\pr}$.
  We omit $\what{\mon{f}}$ from the notation
  and write $\mon{f}\colon\monbase{t}^{\pr}\to \monbase{t}$.
  
  \item
  The \emph{$\kos{K}$-equivariance}
  of a lax $\kos{K}$-tensor functor
  $\mon{f}\colon \monbase{t}^{\pr}\to \monbase{t}$
  is the natural transformation
  $\cevar{\mon{f}}\colon\slot\acts \mon{f}(\slot)
  \Rightarrow\mon{f}(\slot \acts^{\pr} \slot)$
  between functors from $\CK\times \CT^{\pr}$ to $\CT$
  whose component at $z\in \CK$, $X^{\pr}\in \CT^{\pr}$ is
  \begin{equation*}
    \xymatrix@C=50pt{
      \monbase{t}(z)\tensor \mon{f}(X^{\pr})
      \ar[r]^-{(\what{\mon{f}})_z\tensor I_{\mon{f}(X^{\pr})}}
      &\mon{f}\monbase{t}^{\pr}(z)\tensor \mon{f}(X^{\pr})
      \ar[r]^-{(\mon{f}_{\tensor})_{\monbase{t}^{\pr}(z),X^{\pr}}}
      &\mon{f}(\monbase{t}^{\pr}(z)\tensor^{\pr} X^{\pr})
      .
    }
  \end{equation*}

  \item
  Let $\mon{f}$, $\mon{g}\colon\monbase{t}^{\pr}\to \monbase{t}$
  be lax $\kos{K}$-tensor functors.
  A \emph{lax $\kos{K}$-tensor natural transformation}
  $\vartheta\colon\mon{f}\Rightarrow\mon{g}$
  is a monoidal natural transformation
  between the underlying lax monoidal functors
  $\mon{f}$, $\mon{g}\colon\kos{T}^{\pr}\to\kos{T}$
  that satisfies
  $\what{\mon{g}}=(\vartheta\monbase{t}^{\pr})\circ \what{\mon{f}}\colon
  \monbase{t}\Rightarrow\mon{g}\monbase{t}^{\pr}$.
  We say that $\vartheta$ is a
  \emph{lax $\kos{K}$-tensor natural isomorphism}
  if every component of $\vartheta$ is invertible.
\end{itemize}
Let $\kos{T}^{\ppr}$ be a symmetric monoidal category
and $\monbase{t}^{\ppr}\colon\kos{K}\to \kos{T}^{\ppr}$
be a strong monoidal functor,
with the induced $\kos{K}$-action $\acts^{\ppr}$ on $\kos{T}^{\ppr}$.
The composite 
$\mon{f}\mon{f}^{\pr}\colon\monbase{t}^{\ppr}\to \monbase{t}$
of lax $\kos{K}$-tensor functors
$\mon{f}\colon\monbase{t}^{\pr}\to\monbase{t}$,
$\mon{f}^{\pr}\colon\monbase{t}^{\ppr}\to \monbase{t}^{\pr}$
consists of
$\mon{f}\mon{f}^{\pr}\colon\kos{T}^{\ppr}\to\kos{T}$
and $\what{\mon{f}\mon{f}^{\pr}}
=(\mon{f}\what{\mon{f}^{\pr}})\circ \what{\mon{f}}
\colon\monbase{t}
\Rightarrow \mon{f}\mon{f}^{\pr}\monbase{t}^{\ppr}$,
and its $\kos{K}$-equivariance is
$(\cevar{\mon{f}\mon{f}^{\pr}})_{z,X^{\ppr}}
=\mon{f}((\cevar{\mon{f}^{\pr}})_{z,X^{\ppr}})
\circ (\cevar{\mon{f}})_{z,\mon{f}^{\pr}(X^{\ppr})}
\colon z\acts\mon{f}\mon{f}^{\pr}(X^{\ppr})
\to \mon{f}\mon{f}^{\pr}(z\acts^{\ppr} X^{\ppr})$.
The identity lax $\kos{K}$-tensor functor of
$\monbase{t}\colon\kos{K}\to\kos{T}$,
denoted by $\vphantom{\id}\smash[t]{\id}_{\kos{T}}^{\vphantom{*}}$,
is the pair consisting of the identity functor of $\kos{T}$
and the identity natural transformation of $\monbase{t}$.

Let
$\kos{f}\colon\kos{t}^{\pr}\mor\kos{t}$
be a Grothendieck $\kos{K}$-morphism
between Grothendieck $\kos{K}$-prekosmoi
$\kos{T}\xmor{\kos{t}}\kos{K}$,
$\kos{T}^{\pr}\xmor{\kos{t}^{\pr}}\kos{K}$.
The left adjoint part $\kos{f}^*$ is a lax $\kos{K}$-tensor functor
from $\kos{t}^*$ to $\kos{t}^{\pr*}$
with $\what{\kos{f}^*}\colon\kos{t}^{\pr*}\cong \kos{f}^*\kos{t}^*$,
and the right adjoint part $\kos{f}_*$ is a lax $\kos{K}$-tensor functor
from $\kos{t}^{\pr*}$ to $\kos{t}^*$
with $\what{\kos{f}_*}\colon\kos{t}^*\Rightarrow\kos{f}_*\kos{t}^{\pr*}$.
One can check that 
the unit $\eta\colon\vphantom{\id}\smash[t]{\id}_{\kos{T}}^{\vphantom{*}}\Rightarrow\kos{f}_*\kos{f}^*$
and the counit $\epsilon\colon\kos{f}^*\kos{f}_*\Rightarrow\vphantom{\id}\smash[t]{\id}_{\kos{T}^{\pr}}^{\vphantom{*}}$
are lax $\kos{K}$-tensor natural transformations.
Let $\kos{g}\colon\kos{t}^{\pr}\mor\kos{t}$
also be a Grothendieck $\kos{K}$-morphism.
A Grothendieck $\kos{K}$-transformation
$\vartheta\colon\kos{f}\Rightarrow\kos{g}$
is precisely a lax $\kos{K}$-tensor natural transformation
$\vartheta\colon\kos{f}^*\Rightarrow\kos{g}^*$
between the left adjoint parts
$\kos{f}^*$, $\kos{g}^*\colon\kos{t}^*\to \kos{t}^{\pr*}$.

\begin{lemma}
  \label{lem-ActionbaseGro-laxKtnatequiv}
  Let $\kos{T}$, $\kos{T}^{\pr}$
  be symmetric monoidal categories
  and
  $\monbase{t}\colon\kos{K}\to \kos{T}$,
  $\monbase{t}^{\pr}\colon\kos{K}\to \kos{T}^{\pr}$
  be strong monoidal functors.
  Let $\acts$ and $\acts^{\pr}$ be the induced $\kos{K}$-actions
  on $\kos{T}$ and $\kos{T}^{\pr}$, respectively.
  Let $\mon{f}$, $\mon{g}\colon\monbase{t}^{\pr}\to \monbase{t}$
  be lax $\kos{K}$-tensor functors
  and $\vartheta\colon\mon{f}\Rightarrow\mon{g}$
  be a monoidal natural transformation
  between the underlying lax monoidal functors
  $\mon{f}$, $\mon{g}\colon\kos{T}^{\pr}\to \kos{T}$.
  Then $\vartheta$ is a lax $\kos{K}$-tensor natural transformation
  if and only if
  the following diagram commutes for every pair
  $z\in \CK$, $X^{\pr}\in \CT^{\pr}$.
  \begin{equation}
    \label{eq-ActionbaseGro-laxKtnatequiv}
    \vcenter{\hbox{
      \xymatrix@R=15pt@C=35pt{
        z\acts \mon{f}(X^{\pr})
        \ar[r]^-{(\cevar{\mon{f}})_{z,X^{\pr}}}
        \ar[d]_-{I_z\acts \vartheta_{X^{\pr}}}
        &\mon{f}(z\acts^{\pr} X^{\pr})
        \ar[d]^-{\vartheta_{z\acts^{\pr} X^{\pr}}}
        \\
        z\acts \mon{g}(X^{\pr})
        \ar[r]^-{(\cevar{\mon{g}})_{z,X^{\pr}}}
        &\mon{g}(z\acts^{\pr} X^{\pr})
      }
    }}
  \end{equation}
\end{lemma}
\begin{proof}
  If $\vartheta$ is a lax $\kos{K}$-tensor natural transformation,
  the diagram \eqref{eq-ActionbaseGro-laxKtnatequiv}
  commutes by the definitions.
  Conversely, assume that
  \eqref{eq-ActionbaseGro-laxKtnatequiv} commutes
  for every pair $z$, $X^{\pr}$.
  One can check that
  $(\what{\mon{f}})_z\colon\monbase{t}(z)\to \mon{f}\monbase{t}^{\pr}(z)$
  is equal to
  $(\cevar{\mon{f}})_{z,\unit^{\pr}}\circ
  (I_{\monbase{t}(z)}\tensor\mon{f}_{\unit})
  \colon\monbase{t}(z)\tensor \unit\to
  \mon{f}(\monbase{t}^{\pr}(z)\tensor^{\pr} \unit^{\pr})$
  up to coherence isomorphisms.
  Using this description for both $\what{\mon{f}}$ and $\what{\mon{g}}$
  and applying \eqref{eq-ActionbaseGro-laxKtnatequiv} at $X^{\pr}=\unit^{\pr}$,
  we deduce that
  $(\what{\mon{g}})_z=
  \vartheta_{\monbase{t}^{\pr}(z)}\circ (\what{\mon{f}})_z
  \colon\monbase{t}(z)\to\mon{g}\monbase{t}^{\pr}(z)$
  for every $z$.
  Thus $\vartheta$ is a lax $\kos{K}$-tensor natural transformation.
\end{proof}

Let $\Spec(b)$ be an affine scheme in $\kos{K}$.
Then $\slot\otimes b$ is a lax $\kos{K}$-tensor functor
from $\vphantom{\id}\smash[t]{\id}_{\kos{K}}^{\vphantom{*}}$ to $\vphantom{\id}\smash[t]{\id}_{\kos{K}}^{\vphantom{*}}$;
it is equipped with the monoidal natural transformation
$\what{\otimes b}\colon\vphantom{\id}\smash[t]{\id}_{\kos{K}}^{\vphantom{*}}\Rightarrow \slot\otimes b$
whose component at $z\in \CK$ is 
$(\what{\otimes b})_z=I_z\otimes u_b$,
and its $\kos{K}$-equivariance is
$a_{z,x,b}\colon z\otimes (x\otimes b)\cong (z\otimes x)\otimes b$
for $z$, $x\in \CK$.
Recall $\kos{u}_b^{\vphantom{*}}\colon\kos{u}_b^{\vphantom{*}}\mor\vphantom{\id}\smash[t]{\id}_{\kos{K}}^{\vphantom{*}}$
in \eqref{eq-RepGroobj-basechange}.
The left adjoint part $\kos{u}_b^*$
is a lax $\kos{K}$-tensor functor
from $\vphantom{\id}\smash[t]{\id}_{\kos{K}}^{\vphantom{*}}$ to $\kos{u}_b^*$
with $\what{\kos{u}_b^*}\colon\kos{u}_b^*=\kos{u}_b^*$,
and the right adjoint part $\kos{u}_{b*}^{\vphantom{*}}$
is a lax $\kos{K}$-tensor functor
from $\kos{u}_b^*$ to $\vphantom{\id}\smash[t]{\id}_{\kos{K}}^{\vphantom{*}}$
with $\what{\vphantom{\kos{u}_{b*}}\smash[t]{\kos{u}_{b*}^{\vphantom{*}}}}\colon\vphantom{\id}\smash[t]{\id}_{\kos{K}}^{\vphantom{*}}\Rightarrow\kos{u}_{b*}^{\vphantom{*}}\kos{u}_b^*$
given by the adjunction unit.
We have $\kos{u}_{b*}^{\vphantom{*}}\kos{u}_b^*=\slot\otimes b\colon\vphantom{\id}\smash[t]{\id}_{\kos{K}}^{\vphantom{*}}\to \vphantom{\id}\smash[t]{\id}_{\kos{K}}^{\vphantom{*}}$
as lax $\kos{K}$-tensor functors.

Let $f\colon\Spec(b^{\pr})\to \Spec(b)$
be a morphism in $\cat{Aff}(\kos{K})$.
The base change along $f$
is a lax $\kos{K}$-tensor functor
$\kos{f}^*\colon\kos{u}_b^*\to \kos{u}_{b^{\pr}}^*$
with $\what{\kos{f}^*}\colon\kos{u}_{b^{\pr}}^*=\kos{f}^*\kos{u}_b^*$.
We have a lax $\kos{K}$-tensor natural transformation
$\slot\otimes b\Rightarrow \slot\otimes b^{\pr}$
whose component at $z\in \CK$ is $I_z\otimes f$.

We give special attention to lax $\kos{K}$-tensor functors
from $\vphantom{\id}\smash[t]{\id}_{\kos{K}}^{\vphantom{*}}$ to $\vphantom{\id}\smash[t]{\id}_{\kos{K}}^{\vphantom{*}}$.
Let $\mon{f}\colon\vphantom{\id}\smash[t]{\id}_{\kos{K}}^{\vphantom{*}}\to \vphantom{\id}\smash[t]{\id}_{\kos{K}}^{\vphantom{*}}$
be a lax $\kos{K}$-tensor functor.
We denote by
\begin{equation}
  \label{eq-ActionbaseGro-tahardef}
  \tahar{\mon{f}}\colon\slot\otimes \mon{f}(\kappa)\Rightarrow\mon{f}
\end{equation}
the natural transformation whose component at $z\in \CK$ is
$(\tahar{\mon{f}})_z
=\mon{f}(\jmath^{-1}_z)\circ (\cevar{\mon{f}})_{z,\kappa}\colon
z\otimes \mon{f}(\kappa)\to \mon{f}(z)$.
One can check that
$\tahar{\mon{f}}$ is a lax $\kos{K}$-tensor natural transformation,
and that
$(\tahar{\mon{f}})_{\kappa}
=\imath^{-1}_{\mon{f}(\kappa)}\colon\kappa\otimes\mon{f}(\kappa)\cong \mon{f}(\kappa)$.
Let $\mon{g}\colon\vphantom{\id}\smash[t]{\id}_{\kos{K}}^{\vphantom{*}}\to \vphantom{\id}\smash[t]{\id}_{\kos{K}}^{\vphantom{*}}$
also be a lax $\kos{K}$-tensor functor.
We have
$(\tahar{\mon{g}\mon{f}})_z=
\mon{g}((\tahar{\mon{f}})_z)\circ (\cevar{\mon{g}})_{z,\mon{f}(\kappa)}
\colon z\otimes \mon{g}\mon{f}(\kappa)\to \mon{g}\mon{f}(z)$
for $z\in \CK$.
Let $\vartheta\colon\mon{f}\Rightarrow\mon{g}$
be a lax $\kos{K}$-tensor natural transformation.
By taking $\monbase{t}=\monbase{t}^{\pr}=\vphantom{\id}\smash[t]{\id}_{\kos{K}}^{\vphantom{*}}$,
$X^{\pr}=\kappa$ in Lemma~\ref{lem-ActionbaseGro-laxKtnatequiv},
we obtain the following commutative diagram for every $z\in \CK$.
\begin{equation}\label{eq-ActionbaseGro-laxKtnattahar}
  \vcenter{\hbox{
    \xymatrix@R=15pt@C=35pt{
      z\otimes \mon{f}(\kappa)
      \ar[r]^-{(\tahar{\mon{f}})_z}
      \ar[d]_-{I_z\otimes \vartheta_{\kappa}}
      &\mon{f}(z)
      \ar[d]^-{\vartheta_z}
      \\%1
      z\otimes \mon{g}(\kappa)
      \ar[r]^-{(\tahar{\mon{g}})_z}
      &\mon{g}(z)
    }
  }}
\end{equation}

We are now ready to prove Proposition~\ref{prop-Grocatoverpk-HomisIsom}.

\begin{proof}[Proof of Proposition~\ref{prop-Grocatoverpk-HomisIsom}]
  Let $\kos{T}\xmor{\kos{t}}\kos{K}$ be a Grothendieck $\kos{K}$-prekosmos
  and let $\omega$, $\omega^{\pr}\colon\vphantom{\id}\smash[t]{\id}_{\kos{K}}^{\vphantom{*}}\mor\kos{t}$
  be Grothendieck $\kos{K}$-morphisms.
  Put $p\coloneqq\omega^{\pr*}\omega_*(\kappa)\in\cat{Comm}(\kos{K})$.
  Assume that $\omega$ is a fiber point of $\kos{t}$.
  The left adjoint part of $\omega^{\pr}$
  is a lax $\kos{K}$-tensor functor
  $\omega^{\pr*}\colon\kos{t}^*\to \vphantom{\id}\smash[t]{\id}_{\kos{K}}^{\vphantom{*}}$
  whose $\kos{K}$-equivariance
  is a natural isomorphism.
  The right adjoint part of $\omega$
  is a lax $\kos{K}$-tensor functor
  $\omega_*\colon\vphantom{\id}\smash[t]{\id}_{\kos{K}}^{\vphantom{*}}\to \kos{t}^*$
  whose $\kos{K}$-equivariance
  is a natural isomorphism by
  Lemma~\ref{lem-ActionbaseGro-projimplystrongequiv}~(2).
  Hence the composite
  $\omega^{\pr*}\omega_*\colon\vphantom{\id}\smash[t]{\id}_{\kos{K}}^{\vphantom{*}}\to \vphantom{\id}\smash[t]{\id}_{\kos{K}}^{\vphantom{*}}$
  is a lax $\kos{K}$-tensor functor whose $\kos{K}$-equivariance
  is a natural isomorphism.
  Thus, following the notation in \eqref{eq-ActionbaseGro-tahardef},
  we have a lax $\kos{K}$-tensor natural isomorphism
  \begin{equation}\label{eq-ActionbaseGro-taharupperlowerstar}
    \tahar{\omega^{\pr*}\omega_*}
    \colon\slot\otimes p\cong \omega^{\pr*}\omega_*
  \end{equation}
  whose component at $z\in\CK$ is
  $(\tahar{\omega^{\pr*}\omega_*})_z=
  \omega^{\pr*}\omega_*(\jmath^{-1}_z)\circ 
  (\cevar{\omega^{\pr*}\omega_*})_{z,\kappa}
  \colon z\otimes \omega^{\pr*}\omega_*(\kappa)\cong \omega^{\pr*}\omega_*(z)$.
  This in particular proves (1).

  Let $\Spec(b)$ be an affine scheme in $\kos{K}$.
  A Grothendieck $\kos{K}$-transformation
  $\vartheta\colon\omega^{\pr}\!\kos{u}_b^{\vphantom{*}}\Rightarrow\omega\kos{u}_b^{\vphantom{*}}$
  is precisely a lax $\kos{K}$-tensor natural transformation
  $\vartheta\colon\kos{u}_b^*\omega^{\pr*}\Rightarrow \kos{u}_b^*\omega^*$
  between the left adjoint parts.
  By taking mates, we see that such $\vartheta$ are in bijection with
  lax $\kos{K}$-tensor natural transformations
  $\tilde{\vartheta}\colon\omega^{\pr*}\omega_*\Rightarrow\kos{u}_{b*}^{\vphantom{*}}\kos{u}_b^*
  =\slot\otimes b$.
  For each $\vartheta$,
  the component of the corresponding $\tilde{\vartheta}$ at $z\in \CK$ is 
  $(\epsilon_z\otimes I_b)\circ \vartheta_{\omega_*(z)}
  \colon\omega^{\pr*}\omega_*(z)\to z\otimes b$,
  where $\vartheta_{\omega_*(z)}\colon\omega^{\pr*}\omega_*(z)\to \omega^*\omega_*(z)\otimes b$
  is the component of $\vartheta$ at $\omega_*(z)\in\CT$,
  regarded as a morphism in $\CK$.
  For each $\tilde{\vartheta}$,
  the component of the corresponding $\vartheta$ at $X\in \CT$,
  regarded as a morphism in $\CK$,
  is $\tilde{\vartheta}_{\omega^*(X)}\circ \omega^{\pr*}(\eta_X)
  \colon\omega^{\pr*}(X)\to \omega^*(X)\otimes b$.
  We illustrate $\vartheta$ and $\tilde{\vartheta}$ as follows.
  \begin{equation*}
    \begin{tikzcd}[row sep=0pt, column sep=15pt]
      \text{ }
      &\id_{\kos{K}}
      \ar[dl,"\kos{u}_b^*"']
      &\text{ }
      &\text{ }
      &\text{ }
      &\id_{\kos{K}}
      \ar[dl,"\kos{u}_b^*"']
      \ar[dr,"\omega_*"]
      &\text{ }
      \\%1
      \kos{u}_b^*
      &\text{ }
      &\kos{t}^*
      \ar[ul,"\omega^*"']
      \ar[dl,"\omega^{\pr*}"]
      &\text{ }
      &\kos{u}_b^*
      \ar[dr,"\kos{u}_{b*}"']
      &\text{ }
      &\kos{t}^*
      \ar[dl,"\omega^{\pr*}"]
      \ar[ll,2cell,shift left=0.5ex,"{\tilde{\vartheta}}"']
      \\%2
      \text{ }
      &\id_{\kos{K}}
      \ar[ul,"\kos{u}_b^*"]
      \ar[uu,2cell,"{\vartheta}"']
      &\text{ }
      &\text{ }
      &\text{ }
      &\id_{\kos{K}}
      &\text{ }
    \end{tikzcd}
  \end{equation*}
  On the other hand, lax $\kos{K}$-tensor natural transformations
  $\tilde{\vartheta}\colon\omega^{\pr*}\omega_*\Rightarrow \slot\otimes b$
  correspond bijectively to morphisms
  $g\colon p\to b$ in $\cat{Comm}(\kos{K})$;
  each $\tilde{\vartheta}$ corresponds to
  $\imath_b^{-1}\circ \tilde{\vartheta}_{\kappa}\colon p\to b$,
  and each $g$ corresponds to the $\tilde{\vartheta}$
  whose component at $z\in\CK$ is
  $(I_z\otimes g)\circ (\tahar{\omega^{\pr*}\omega_*})_z^{-1}
  \colon\omega^{\pr*}\omega_*(z)\to z\otimes b$.
  One can check this by observing that each $\tilde{\vartheta}$
  satisfies the relation \eqref{eq-ActionbaseGro-laxKtnattahar}
  for every $z\in \CK$,
  which takes the following form.
  \begin{equation*}
    \vcenter{\hbox{
      \xymatrix@R=12pt@C=35pt{
        z\otimes \omega^{\pr*}\omega_*(\kappa)
        \ar[d]_-{I_z\otimes \tilde{\vartheta}_{\kappa}}
        \ar[r]^-{(\tahar{\omega^{\pr*}\omega_*})_z}_-{\cong}
        &\omega^{\pr*}\omega_*(z)
        \ar[d]^-{\tilde{\vartheta}_z}
        \\%1
        z\otimes (\kappa\otimes b)
        \ar[r]^-{I_z\otimes \imath^{-1}_b}_-{\cong}
        &z\otimes b
      }
    }}
  \end{equation*}
  Thus we obtain a bijection between
  Grothendieck $\kos{K}$-transformations
  $\vartheta\colon\omega^{\pr}\!\kos{u}_b^{\vphantom{*}}\Rightarrow\omega\kos{u}_b^{\vphantom{*}}$
  and morphisms $g\colon p\to b$ in $\cat{Comm}(\kos{K})$.
  For each $X\in\CT$, we define a morphism
  $\rho_X\colon \omega^{\pr*}(X)\to \omega^*(X)\otimes p$
  in $\CK$ as
  \begin{equation}\label{eq-ActionbaseGro-rhoX}
    \xymatrix@C=12pt{
      \rho_X\colon
      \omega^{\pr*}(X)
      \ar[rr]^-{\omega^{\pr*}(\eta_X)}
      &\text{ }
      &\omega^{\pr*}\omega_*\omega^*(X)
      \ar[rrr]^-{(\tahar{\omega^{\pr*}\omega_*})_{\omega^*(X)}^{-1}}_-{\cong}
      &\text{ }
      &\text{ }
      &\omega^*(X)\otimes \omega^{\pr*}\omega_*(\kappa).
    }
  \end{equation}
  Each $\vartheta$ corresponds to 
  $g=\imath_b^{-1}\circ (\epsilon_{\kappa}\otimes I_b)\circ \vartheta_{\omega_*(\kappa)}\colon p\to b$,
  where $\vartheta_{\omega_*(\kappa)}\colon\omega^{\pr*}\omega_*(\kappa)\to \omega^*\omega_*(\kappa)\otimes b$
  is the component of $\vartheta$ at $\omega_*(\kappa)\in \CT$,
  regarded as a morphism in $\CK$.
  Each $g\colon p\to b$ corresponds to the
  $\vartheta\colon\omega^{\pr}\!\kos{u}_b^{\vphantom{*}}\Rightarrow\omega\kos{u}_b^{\vphantom{*}}$
  whose component at $X\in \CT$,
  regarded as a morphism in $\CK$,
  is $(I_{\omega^*(X)}\otimes g)\circ \rho_X
  \colon\omega^{\pr*}(X)\to \omega^*(X)\otimes b$.
  The correspondence between
  $\vartheta$ and $g$ is natural in $\Spec(b)$
  and induces an isomorphism
  $\Hom_{\cat{Aff}(\kos{K})}(\slot,\Spec(p))
  \cong\uHom_{\kos{K}}(\omega^{\pr},\omega)$
  of presheaves on $\cat{Aff}(\kos{K})$.
  We write $\xi\colon\omega^{\pr}\!\kos{u}_p^{\vphantom{*}}\Rightarrow \omega\kos{u}_p^{\vphantom{*}}$
  for the universal element
  corresponding to the identity morphism of $\Spec(p)$.
  The component of $\xi$ at $X\in \CT$,
  regarded as a morphism in $\CK$,
  is $\rho_X$ in \eqref{eq-ActionbaseGro-rhoX}.
  This proves (2).
  
  The equality
  $\uHom_{\kos{K}}(\omega^{\pr},\omega)
  =\uIsom_{\kos{K}}(\omega^{\pr},\omega)$
  holds precisely when the universal element
  $\xi\colon\omega^{\pr}\!\kos{u}_p^{\vphantom{*}}\Rightarrow \omega\kos{u}_p^{\vphantom{*}}$
  is invertible;
  for if $\xi$ is invertible,
  since every
  $\vartheta\colon\omega^{\pr}\!\kos{u}_b^{\vphantom{*}}\Rightarrow\omega\kos{u}_b^{\vphantom{*}}$
  is the base change of $\xi$ along a unique morphism 
  $g\colon\Spec(b)\to\Spec(p)$ in $\cat{Aff}(\kos{K})$,
  the base change of $\xi^{-1}$ along the same $g$
  is the inverse of $\vartheta$.
  On the other hand,
  $\xi\colon\kos{u}_p^*\omega^{\pr*}\Rightarrow \kos{u}_p^*\omega^*$
  is invertible if and only if
  $\kos{u}_{p*}^{\vphantom{*}}\xi\colon\omega^{\pr*}(\slot)\otimes p
  \Rightarrow \omega^*(\slot)\otimes p$
  is invertible,
  since $\kos{u}_{p*}^{\vphantom{*}}\colon\CK_p^{\vphantom{*}}\to \CK$ is conservative.
  For each $X\in \CT$, we have the
  following explicit description of $\kos{u}_{p*}^{\vphantom{*}}(\xi_X)$.
  \begin{equation}\label{eq-ActionbaseGro-uplowerstarxi}
    \vcenter{\hbox{
      \xymatrix@R=15pt@C=4.8pt{
        \omega^{\pr*}(X)\otimes\omega^{\pr*}\omega_*(\kappa)
        \ar[d]_-{(\omega^{\pr*}_{\tensor})_{X,\omega_*(\kappa)}}^-{\cong}
        \ar[rrrrr]^-{\kos{u}_{p*}(\xi_X)}
        &\text{ }
        &\text{ }
        &\text{ }
        &\text{ }
        &\omega^*(X)\otimes \omega^{\pr*}\omega_*(\kappa)
        \\%1
        \omega^{\pr*}(X\tensor \omega_*(\kappa))
        \ar[rr]^-{\omega^{\pr*}(\varphi_{X,\kappa})}_-{\cong}
        &\text{ }
        &\omega^{\pr*}\omega_*(\omega^*(X)\otimes \kappa)
        \ar[rrr]^-{\omega^{\pr*}\omega_*(\jmath_{\omega^*(X)}^{-1})}_-{\cong}
        &\text{ }
        &\text{ }
        &\omega^{\pr*}\omega_*\omega^*(X)
        \ar[u]_-{(\tahar{\omega^{\pr*}\omega_*})_{\omega^*(X)}^{-1}}^-{\cong}
      }
    }}
  \end{equation}
  Here $\varphi$ is the natural transformation
  associated to $\omega$ in \eqref{eq-Groprekos-projformula}.
  Since $\omega$ is locally connected,
  the morphism $\omega^{\pr*}(\varphi_{X,\kappa})$
  in \eqref{eq-ActionbaseGro-uplowerstarxi}
  is invertible for every $X\in\CT$.
  Thus $\kos{u}_{p*}^{\vphantom{*}}\xi$ is invertible, and so is $\xi$.
  This proves (3).
  Throughout the proof, we have used only that
  $\omega$ is locally connected,
  which is equivalent to our assumption by
  Remark~\ref{rem-Groprekos-omegalowerstarconservative}.
\end{proof}

Let $\kos{t}$, $\omega$, $\omega^{\pr}$, and $p$
be as in Proposition~\ref{prop-Grocatoverpk-HomisIsom}
and recall the universal element
$\xi$ of the representation
$\Hom_{\cat{Aff}(\kos{K})}(\slot,\Spec(p))
\cong\uHom_{\kos{K}}(\omega^{\pr},\omega)$.
We denote by
$\mathring{\rho}_X\colon\omega^*(X)\to \omega^{\pr*}(X)\otimes p$
the component of the inverse $\xi^{-1}$ at $X\in \CT$,
regarded as a morphism in $\CK$;
it is explicitly given by
\begin{equation}
  \label{eq-ActionbaseGro-xiinverse}
  \xymatrix@C=45pt{
    \mathring{\rho}_X\colon
    \omega^*(X)
    \ar[r]^-{I_{\omega^*(X)}\otimes u_p}
    &\omega^*(X)\otimes \omega^{\pr*}\omega_*(\kappa)
    \ar[r]^-{\kos{u}_{p*}(\xi_X^{-1})}_-{\cong}
    &\omega^{\pr*}(X)\otimes \omega^{\pr*}\omega_*(\kappa)
    ,
  }
\end{equation}
where $\kos{u}_{p*}^{\vphantom{*}}(\xi_X^{-1})=\kos{u}_{p*}^{\vphantom{*}}(\xi_X)^{-1}$
is the composite, in reverse order, of the inverses of the isomorphisms in 
\eqref{eq-ActionbaseGro-uplowerstarxi}.
Then for each $\Spec(b)\in\cat{Aff}(\kos{K})$
and each 
$\vartheta\colon\omega^{\pr}\!\kos{u}_b^{\vphantom{*}}\Rightarrow\omega\kos{u}_b^{\vphantom{*}}$,
the component of $\vartheta^{-1}$ at $X\in \CT$,
regarded as a morphism in $\CK$,
is given by
$(I_{\omega^{\pr*}(X)}\otimes g)\circ \mathring{\rho}_X
\colon\omega^*(X)\to \omega^{\pr*}(X)\otimes b$,
where $g=\imath_b^{-1}\circ (\epsilon_{\kappa}\otimes I_b)\circ \vartheta_{\omega_*(\kappa)}\colon p\to b$
is the morphism in $\cat{Comm}(\kos{K})$ corresponding to $\vartheta$.

The following lemma, along with Proposition~\ref{prop-Grocatoverpk-HomisIsom},
will be used in the proof of Theorem~\ref{thm-Grocatoverpk-mainthm}.

\begin{lemma}
  \label{lem-ActionbaseGro-inverseandcomposition}
  Let $\kos{T}\xmor{\kos{t}}\kos{K}$ be a Grothendieck $\kos{K}$-prekosmos
  and $\omega$, $\omega^{\pr}\colon\vphantom{\id}\smash[t]{\id}_{\kos{K}}^{\vphantom{*}}\mor\kos{t}$
  be 
  fiber points of $\kos{t}$.
  \begin{enumerate}
    \item 
    The isomorphism
    \begin{equation*}
      \uIsom_{\kos{K}}(\omega^{\pr},\omega)
      \cong \uIsom_{\kos{K}}(\omega,\omega^{\pr})
    \end{equation*}
    of presheaves
    on $\cat{Aff}(\kos{K})$,
    obtained by taking inverses,
    corresponds to the isomorphism
    $\Spec(\omega^{\pr*}\omega_*(\kappa))
    \cong \Spec(\omega^*\omega^{\pr}_*(\kappa))$
    between the representing objects,
    given by
    \begin{equation*}
      \xymatrix@C=11pt{
        \omega^*\omega^{\pr}_*(\kappa)
        \ar[rr]^-{\mathring{\rho}_{\omega^{\pr}_*(\kappa)}}
        &\text{ }
        &\omega^{\pr*}\omega^{\pr}_*(\kappa)\otimes \omega^{\pr*}\omega_*(\kappa)
        \ar[rrr]^-{\epsilon^{\pr}_{\kappa}\otimes I_{\omega^{\pr*}\omega_*(\kappa)}}
        &\text{ }
        &\text{ }
        &\omega^{\pr*}\omega_*(\kappa)
        .
      }
    \end{equation*}

    \item
    Let $\omega^{\ppr}\colon\vphantom{\id}\smash[t]{\id}_{\kos{K}}^{\vphantom{*}}\mor\kos{t}$
    be a Grothendieck $\kos{K}$-morphism.
    The morphism 
    \begin{equation*}
      \uIsom_{\kos{K}}(\omega^{\pr},\omega)\times \uIsom_{\kos{K}}(\omega^{\ppr},\omega^{\pr})
      \to \uIsom_{\kos{K}}(\omega^{\ppr},\omega)
    \end{equation*}
    of presheaves
    on $\cat{Aff}(\kos{K})$,
    obtained by the vertical composition,
    corresponds to the morphism
    $\Spec(\omega^{\pr*}\omega_*(\kappa))
    \times
    \Spec(\omega^{\ppr*}\omega^{\pr}_*(\kappa))
    \to \Spec(\omega^{\ppr*}\omega_*(\kappa))$
    between the representing objects,
    given by
    \begin{equation*}
      \xymatrix@C=16pt{
        \omega^{\ppr*}\omega_*(\kappa)
        \ar[rr]^-{\omega^{\ppr*}(\eta^{\pr}_{\omega_*(\kappa)})}
        &\text{ }
        &\omega^{\ppr*}\omega^{\pr}_*\omega^{\pr*}\omega_*(\kappa)
        \ar[rrr]^-{(\tahar{\omega^{\ppr*}\omega^{\pr}_*})_{\omega^{\pr*}\omega_*(\kappa)}^{-1}}_-{\cong}
        &\text{ }
        &\text{ }
        &\omega^{\pr*}\omega_*(\kappa)
        \otimes
        \omega^{\ppr*}\omega^{\pr}_*(\kappa)
        .
      }
    \end{equation*}
  \end{enumerate}
\end{lemma}
\begin{proof}
  Put $p\coloneqq\omega^{\pr*}\omega_*(\kappa)\in\cat{Comm}(\kos{K})$.
  Since $\omega$ is a fiber point of $\kos{t}$,
  Proposition~\ref{prop-Grocatoverpk-HomisIsom}
  implies that 
  $\Spec(p)$
  represents
  $\uHom_{\kos{K}}(\omega^{\pr},\omega)
  =\uIsom_{\kos{K}}(\omega^{\pr},\omega)$.
  Let $\xi\colon\omega^{\pr}\!\kos{u}_p^{\vphantom{*}}\cong \omega\kos{u}_p^{\vphantom{*}}$
  be the universal element of this representation.
  Since $\omega^{\pr}$ is a fiber point of $\kos{t}$ as well,
  by exchanging the roles of $\omega$ and $\omega^{\pr}$
  in Proposition~\ref{prop-Grocatoverpk-HomisIsom},
  we deduce that 
  $\Spec(\omega^*\omega^{\pr}_*(\kappa))$
  represents
  $\uHom_{\kos{K}}(\omega,\omega^{\pr})
  =\uIsom_{\kos{K}}(\omega,\omega^{\pr})$.
  The isomorphism $\Spec(p)\cong \Spec(\omega^*\omega^{\pr}_*(\kappa))$
  in (1) is the unique morphism in $\cat{Aff}(\kos{K})$
  that corresponds to
  $\xi^{-1}\colon\omega\kos{u}_p^{\vphantom{*}}\cong \omega^{\pr}\!\kos{u}_p^{\vphantom{*}}$,
  and is explicitly given by
  $\imath_p^{-1}\circ (\epsilon^{\pr}_{\kappa}\otimes I_p)\circ \mathring{\rho}_{\omega^{\pr}_*(\kappa)}
  \colon\omega^*\omega^{\pr}_*(\kappa)\cong p$.
  Here $\mathring{\rho}_{\omega^{\pr}_*(\kappa)}$
  is the component of $\xi^{-1}$ at $\omega^{\pr}_*(\kappa)\in\CT$,
  regarded as a morphism in $\CK$,
  given by substituting $X=\omega^{\pr}_*(\kappa)$
  in \eqref{eq-ActionbaseGro-xiinverse}.
  This proves (1).
  
  Let $\Spec(b)$ be an affine scheme in $\kos{K}$
  and $\vartheta\colon\omega^{\pr}\!\kos{u}_b^{\vphantom{*}}\Rightarrow \omega\kos{u}_b^{\vphantom{*}}$,
  $\vartheta^{\pr}\colon\omega^{\ppr}\!\kos{u}_b^{\vphantom{*}}\Rightarrow \omega^{\pr}\!\kos{u}_b^{\vphantom{*}}$
  be Grothendieck $\kos{K}$-transformations,
  with their vertical composite 
  $\vartheta^{\ppr}\coloneqq\vartheta\circ\vartheta^{\pr}\colon
  \omega^{\ppr}\!\kos{u}_b^{\vphantom{*}}\Rightarrow \omega\kos{u}_b^{\vphantom{*}}$.
  Note that $\vartheta$, $\vartheta^{\pr}$, and $\vartheta^{\ppr}$
  are all invertible by Proposition~\ref{prop-Grocatoverpk-HomisIsom}~(3).
  We denote by 
  $\tilde{\vartheta}\colon\omega^{\pr*}\omega_*\Rightarrow \slot\otimes b$,
  $\tilde{\vartheta}^{\pr}\colon\omega^{\ppr*}\omega^{\pr}_*\Rightarrow \slot\otimes b$,
  and
  $\tilde{\vartheta}^{\ppr}\colon\omega^{\ppr*}\omega_*\Rightarrow \slot\otimes b$
  the mates of $\vartheta$, $\vartheta^{\pr}$, and $\vartheta^{\ppr}$,
  respectively.
  One can check that the component of $\tilde{\vartheta}^{\ppr}$
  at $z\in\CK$ is given by
  \begin{equation}\label{eq-ActionbaseGro-inverseandcomposition}
    \xymatrix@C=35pt{
      \tilde{\vartheta}^{\ppr}_z\colon
      \omega^{\ppr*}\omega_*(z)
      \ar[r]^-{\omega^{\ppr*}(\eta^{\pr}_{\omega_*(z)})}
      &\omega^{\ppr*}\omega^{\pr}_*\omega^{\pr*}\omega_*(z)
      \ar[r]^-{(\tilde{\vartheta}^{\pr}\tilde{\vartheta})_z}
      &z\otimes b\otimes b
      \ar[r]^-{I_z\otimes \pc_b}
      &z\otimes b
      .
    }
  \end{equation}
  Let $g\colon\omega^{\pr*}\omega_*(\kappa)\to b$,
  $g^{\pr}\colon\omega^{\ppr*}\omega^{\pr}_*(\kappa)\to b$,
  and $g^{\ppr}\colon\omega^{\ppr*}\omega_*(\kappa)\to b$
  be the morphisms in $\cat{Comm}(\kos{K})$
  corresponding to $\tilde{\vartheta}$, $\tilde{\vartheta}^{\pr}$, and $\tilde{\vartheta}^{\ppr}$,
  respectively.
  Using the relation \eqref{eq-ActionbaseGro-inverseandcomposition}
  for $z=\kappa$,
  and the relation
  $\tilde{\vartheta}^{\pr}_z=
  (I_z\otimes g^{\pr})\circ (\tahar{\omega^{\ppr*}\omega^{\pr}_*})_z^{-1}
  \colon\omega^{\ppr*}\omega^{\pr}_*(z)\to z\otimes b$
  for $z=\omega^{\pr*}\omega_*(\kappa)$,
  we deduce that
  $g^{\ppr}=\pc_b\circ (g\otimes g^{\pr})\circ (\tahar{\omega^{\ppr*}\omega^{\pr}_*})_{\omega^{\pr*}\omega_*(\kappa)}^{-1}\circ \omega^{\ppr*}(\eta^{\pr}_{\omega_*(\kappa)})$.
  This proves (2).
\end{proof}

\section{The main theorem}
\label{sec-GroMainthm}

We now prove Theorem~\ref{thm-Grocatoverpk-mainthm},
thereby completing the proof of Theorem~\ref{thm-Grocatoverpk-PrimaryObj}.

\begin{proof}[Proof of Theorem~\ref{thm-Grocatoverpk-mainthm}]
  Let $\kos{T}\xmor{\kos{t}}\kos{K}$
  be a Grothendieck $\kos{K}$-prekosmos
  and $\omega\colon\vphantom{\id}\smash[t]{\id}_{\kos{K}}^{\vphantom{*}}\mor\kos{t}$
  be a surjective fiber point of $\kos{t}$.
  Put $\pi\coloneqq\omega^*\omega_*(\kappa)\in\cat{Comm}(\kos{K})$.
  Proposition~\ref{prop-Grocatoverpk-HomisIsom}
  implies that the affine scheme $\Spec(\pi)$ in $\kos{K}$
  represents the presheaf
  $\uHom_{\kos{K}}(\omega,\omega)
  =\uIsom_{\kos{K}}(\omega,\omega)$
  on $\cat{Aff}(\kos{K})$.
  Let $\xi\colon\omega\kos{u}_{\pi}^{\vphantom{*}}\cong \omega\kos{u}_{\pi}^{\vphantom{*}}$
  be the universal element of this representation.
  By \eqref{eq-ActionbaseGro-rhoX} for $\omega^{\pr}=\omega$,
  the component of $\xi$ at $X\in \CT$,
  regarded as a morphism in $\CK$,
  is given by
  \begin{equation}\label{eq-GroMainthm-rho}
    \xymatrix@C=12pt{
      \rho_X\colon
      \omega^*(X)
      \ar[rr]^-{\omega^*(\eta_X)}
      &\text{ }
      &\omega^*\omega_*\omega^*(X)
      \ar[rrr]^-{(\tahar{\omega^*\omega_*})_{\omega^*(X)}^{-1}}_-{\cong}
      &\text{ }
      &\text{ }
      &\omega^*(X)\otimes \omega^*\omega_*(\kappa)
      .
    }
  \end{equation}
  By \eqref{eq-ActionbaseGro-xiinverse} for $\omega^{\pr}=\omega$,
  the component of the inverse $\xi^{-1}$ at $X\in \CT$,
  regarded as a morphism in $\CK$,
  is given by
  \begin{equation}
    \label{eq-GroMainthm-mathringrho}
    \xymatrix@C=45pt{
      \mathring{\rho}_X\colon
      \omega^*(X)
      \ar[r]^-{I_{\omega^*(X)}\otimes u_{\pi}}
      &\omega^*(X)\otimes \omega^*\omega_*(\kappa)
      \ar[r]^-{\kos{u}_{\pi*}(\xi_X^{-1})}_-{\cong}
      &\omega^*(X)\otimes \omega^*\omega_*(\kappa)
      .
    }
  \end{equation}
  One can check that the identity Grothendieck $\kos{K}$-transformation
  of $\omega$ corresponds to the morphism
  $\epsilon_{\kappa}\colon\omega^*\omega_*(\kappa)\to \kappa$
  in $\cat{Comm}(\kos{K})$.
  Combining this with Lemma~\ref{lem-ActionbaseGro-inverseandcomposition}
  for $\omega^{\pr}=\omega^{\ppr}=\omega$,
  we deduce that $\Spec(\pi)$ equipped with
  \begin{equation}\label{eq-GroMainthm-Specpigroupstr}
    \begin{aligned}
      \cp_{\omega^*\omega_*(\kappa)}
      \colon&\!\!
      \xymatrix@C=12pt{
        \omega^*\omega_*(\kappa)
        \ar[rr]^-{\omega^*(\eta_{\omega_*(\kappa)})}
        &\text{ }
        &\omega^*\omega_*\omega^*\omega_*(\kappa)
        \ar[rrr]^-{(\tahar{\omega^*\omega_*})_{\omega^*\omega_*(\kappa)}^{-1}}_-{\cong}
        &\text{ }
        &\text{ }
        &\omega^*\omega_*(\kappa)\otimes\omega^*\omega_*(\kappa)
      }
      \\%1
      e_{\omega^*\omega_*(\kappa)}
      \colon&\!\!
      \xymatrix@C=25pt{
        \omega^*\omega_*(\kappa)
        \ar[r]^-{\epsilon_{\kappa}}
        &\kappa
      }
      \\%2
      \varsigma_{\omega^*\omega_*(\kappa)}
      \colon&\!\!
      \xymatrix@C=10pt{
        \omega^*\omega_*(\kappa)
        \ar[rr]^-{\mathring{\rho}_{\omega_*(\kappa)}}
        &\text{ }
        &\omega^*\omega_*(\kappa)\otimes \omega^*\omega_*(\kappa)
        \ar[rrr]^-{\epsilon_{\kappa}\otimes I_{\omega^*\omega_*(\kappa)}}
        &\text{ }
        &\text{ }
        &\omega^*\omega_*(\kappa)
      }
    \end{aligned}
  \end{equation}
  is an affine group scheme in $\kos{K}$
  that represents the presheaf of groups
  $\uAut_{\kos{K}}(\omega)$ on $\cat{Aff}(\kos{K})$,
  where the morphism $\mathring{\rho}_{\omega_*(\kappa)}$
  in the description of the antipode
  is given by substituting $X=\omega_*(\kappa)$ in \eqref{eq-GroMainthm-mathringrho}.
  The right adjoint part $\omega_*$ preserves coreflexive equalizers,
  and the left adjoint part $\omega^*$
  also preserves coreflexive equalizers,
  since $\omega$ is surjective.
  From $\tahar{\omega^*\omega_*}\colon\slot\otimes \pi\cong \omega^*\omega_*$
  we deduce that
  $\Spec(\pi)$ is a Grothendieck object in $\kos{K}$.
  This proves (1).

  It remains to show (2).
  Put $\varpi\coloneqq\omega^*\omega_*$.
  We begin by showing that
  $\tahar{\varpi}\colon\slot\otimes\pi\cong \varpi$
  is an isomorphism of comonads on $\CK$.
  Let $z\in\CK$.
  One can check that
  $\omega^*(\eta_{\omega_*(z)})\circ (\tahar{\varpi})_z
  =(\tahar{\varpi}\tahar{\varpi})_z\circ (I_z\otimes \cp_{\pi})
  \colon z\otimes \pi\to \varpi\varpi(z)$,
  using \eqref{eq-ActionbaseGro-laxKtnattahar}
  for $\vartheta=\omega^*\eta\omega_*\colon\varpi\Rightarrow\varpi\varpi$,
  the relation
  $(\tahar{\varpi\varpi})_z
  =\varpi((\tahar{\varpi})_z)\circ (\cevar{\varpi})_{z,\pi}
  \colon z\otimes \varpi(\pi)\cong \varpi\varpi(z)$,
  and \eqref{eq-ActionbaseGro-laxKtnatequiv}
  for $\vartheta=\tahar{\varpi}\colon\slot\otimes \pi\cong \varpi$,
  $X^{\pr}=\pi$.
  We also have
  $I_z\otimes e_{\pi}=\epsilon_z\circ (\tahar{\varpi})_z
  \colon z\otimes \pi\to z$,
  as it is the relation \eqref{eq-ActionbaseGro-laxKtnattahar}
  for $\vartheta=\epsilon\colon\varpi\Rightarrow\vphantom{\id}\smash[t]{\id}_{\kos{K}}^{\vphantom{*}}$.
  Hence $\tahar{\varpi}\colon\slot\otimes \pi\cong \varpi$
  is an isomorphism of comonads on $\CK$.

  The assignment
  $X\mapsto (\omega^*(X),\rho_X)$
  defines a functor
  $\breve{\omega}^*\colon\CT\to\cat{Rep}_{\kos{K}}(\pi)$,
  since $\tahar{\varpi}\colon\slot\otimes\pi\cong \varpi$
  is an isomorphism of comonads on $\CK$.
  The coherence isomorphisms
  $\omega^*_{\tensor}$, $\omega^*_{\unit}$ of $\omega^*$
  lift to isomorphisms
  $\breve{\omega}^*_{\tensor}$, $\breve{\omega}^*_{\unit}$
  in $\cat{Rep}_{\kos{K}}(\pi)$,
  since $\rho\colon\omega^*(\slot)\Rightarrow\omega^*(\slot)\otimes \pi$
  is a monoidal natural transformation;
  they are the coherence isomorphisms of the strong monoidal functor
  $\breve{\omega}^*\colon\kos{T}\to\kosRep_{\kos{K}}(\pi)$.
  The monoidal natural isomorphism
  $\what{\omega^*}\colon\vphantom{\id}\smash[t]{\id}_{\kos{K}}^{\vphantom{*}}\cong \omega^*\kos{t}^*$
  lifts to a monoidal natural isomorphism
  $\what{\breve{\omega}^*}\colon\kos{t}_{\pi}^*\cong \breve{\omega}^*\kos{t}^*$,
  since $\rho\colon\omega^*(\slot)\Rightarrow\omega^*(\slot)\otimes \pi$
  is moreover a lax $\kos{K}$-tensor natural transformation;
  it is the structure of the lax $\kos{K}$-tensor functor
  $\breve{\omega}^*\colon\kos{t}^*\to \kos{t}_{\pi}^*$,
  omitted from the notation.
  We have an equality
  $\omega^*=\kos{e}_{\pi}^*\breve{\omega}^*\colon\kos{t}^*\to \vphantom{\id}\smash[t]{\id}_{\kos{K}}^{\vphantom{*}}$
  of lax $\kos{K}$-tensor functors.
  
  We present the rest of the proof in three steps.
  The arguments are essentially the same as
  those in the proof of the crude comonadicity theorem,
  the dual of \cite[Theorem~2.1]{Beck1968}.
  Nonetheless, we carry them out in full,
  as they involve the terminology
  specific to the present setting.
  Given $(x,\rho_x)\in\cat{Rep}_{\kos{K}}(\pi)$,
  we write
  $\dot{\rho}_x\coloneqq(\tahar{\varpi})_x\circ\rho_x
  \colon x\to \omega^*\omega_*(x)$
  for short.

  First, we show that
  $\omega\colon\vphantom{\id}\smash[t]{\id}_{\kos{K}}^{\vphantom{*}}\mor\kos{t}$
  extends to a Grothendieck $\kos{K}$-morphism
  $\breve{\omega}\colon\kos{t}_{\pi}^{\vphantom{*}}\mor\kos{t}$
  up to an invertible Grothendieck $\kos{K}$-transformation
  $\omega\cong \breve{\omega}\kos{e}_{\pi}^{\vphantom{*}}$.
  For each $(x,\rho_x)\in\cat{Rep}_{\kos{K}}(\pi)$,
  we fix an equalizer
  $\breve{\omega}_*(x,\rho_x)\hookrightarrow \omega_*(x)$
  in $\CT$ of the coreflexive pair
  $\omega_*(\dot{\rho}_x)$,
  $\eta_{\omega_*(x)}
  \colon\omega_*(x)\to \omega_*\omega^*\omega_*(x)$,
  with common retraction
  $\omega_*(\epsilon_x)$.
  Then $\breve{\omega}_*$ is right adjoint to $\breve{\omega}^*$.
  The component of the adjunction unit
  $\breve{\eta}\colon\vphantom{\id}\smash[t]{\id}_{\CT}^{\vphantom{*}}\Rightarrow\breve{\omega}_*\breve{\omega}^*$
  at $X\in \CT$ is the unique morphism
  $X\to \breve{\omega}_*\breve{\omega}^*(X)$
  in $\CT$ whose composite with
  $\breve{\omega}_*\breve{\omega}^*(X)
  \hookrightarrow\omega_*\omega^*(X)$ is $\eta_X$.
  The component of the adjunction counit
  $\breve{\epsilon}\colon\breve{\omega}^*\breve{\omega}_*\Rightarrow\id_{\cat{Rep}_{\kos{K}}(\pi)}$
  at $(x,\rho_x)$,
  seen as a morphism in $\CK$,
  is the composite of
  $\omega^*\breve{\omega}_*(x,\rho_x)
  \hookrightarrow \omega^*\omega_*(x)$
  followed by $\epsilon_x\colon \omega^*\omega_*(x)\to x$.
  Thus 
  $\breve{\omega}=(\breve{\omega}^*\dashv \breve{\omega}_*)\colon
  \kosRep_{\kos{K}}(\pi)\mor\kos{T}$
  is a Grothendieck morphism.
  It is moreover a Grothendieck $\kos{K}$-morphism
  $\breve{\omega}\colon\kos{t}_{\pi}^{\vphantom{*}}\mor\kos{t}$,
  where $\what{\breve{\omega}}\colon\kos{t}_{\pi}^{\vphantom{*}}\cong \kos{t}\breve{\omega}$
  is the monoidal natural isomorphism
  $\what{\breve{\omega}^*}\colon\kos{t}_{\pi}^*\cong \breve{\omega}^*\kos{t}^*$
  between the left adjoint parts.
  The invertible Grothendieck $\kos{K}$-transformation
  $\omega\cong \breve{\omega}\kos{e}_{\pi}^{\vphantom{*}}$
  is the equality
  $\omega^*=\kos{e}_{\pi}^*\breve{\omega}^*$
  between the left adjoint parts.

  Next, we show that
  $\breve{\omega}\colon\kos{t}_{\pi}^{\vphantom{*}}\mor\kos{t}$ is an embedding.
  Let $(x,\rho_x)\in\cat{Rep}_{\kos{K}}(\pi)$.
  Since $\omega^*$ preserves coreflexive equalizers,
  $\omega^*\breve{\omega}_*(x,\rho_x)\hookrightarrow \omega^*\omega_*(x)$
  is an equalizer of the coreflexive pair
  $\omega^*\omega_*(\dot{\rho}_x)$,
  $\omega^*(\eta_{\omega_*(x)})
  \colon\omega^*\omega_*(x)\to \omega^*\omega_*\omega^*\omega_*(x)$.
  On the other hand, $\dot{\rho}_x\colon
  x\to \omega^*\omega_*(x)$
  is a split equalizer of the same pair
  with respect to the splitting
  $\epsilon_x$, $\epsilon_{\omega^*\omega_*(x)}$.
  Thus the underlying morphism in $\CK$
  of $\breve{\epsilon}_{(x,\rho_x)}$ is an isomorphism,
  as it is the unique morphism between the two equalizers.
  Hence $\breve{\epsilon}$ is a natural isomorphism,
  and $\breve{\omega}\colon\kos{t}_{\pi}^{\vphantom{*}}\mor\kos{t}$ is an embedding.

  Finally, we show that
  $\breve{\omega}\colon\kos{t}_{\pi}^{\vphantom{*}}\mor\kos{t}$
  is an equivalence.
  It suffices to show that
  $\breve{\eta}$ is a natural isomorphism as well.
  Let $X\in\CT$.
  Observe that $\dot{\rho}_{\omega^*(X)}=\omega^*(\eta_X)$.
  Since $\omega^*$ preserves coreflexive equalizers,
  $\omega^*\breve{\omega}_*\breve{\omega}^*(X)
  \hookrightarrow \omega^*\omega_*\omega^*(X)$
  is an equalizer of the coreflexive pair
  $\omega^*\omega_*\omega^*(\eta_X)$,
  $\omega^*(\eta_{\omega_*\omega^*(X)})
  \colon\omega^*\omega_*\omega^*(X)\to \omega^*\omega_*\omega^*\omega_*\omega^*(X)$.
  On the other hand, $\omega^*(\eta_X)\colon\omega^*(X)\to \omega^*\omega_*\omega^*(X)$
  is a split equalizer of the same pair
  with respect to the splitting
  $\epsilon_{\omega^*(X)}$,
  $\epsilon_{\omega^*\omega_*\omega^*(X)}$.
  Thus $\omega^*(\breve{\eta}_X)$ is an isomorphism,
  as it is the unique morphism between the two equalizers.
  Since $\omega^*$ is conservative,
  $\breve{\eta}_X$ is an isomorphism.
  Hence $\breve{\eta}$ is a natural isomorphism,
  and $\breve{\omega}\colon\kos{t}_{\pi}^{\vphantom{*}}\mor\kos{t}$
  is an equivalence.
  This completes the proof of Theorem~\ref{thm-Grocatoverpk-mainthm}.
\end{proof}

\begin{remark}
  \label{rem-GroMainthm-reconstruction}
  Let $\Spec(\pi)$ be a Grothendieck object in $\kos{K}$
  and $\kos{e}_{\pi}^{\vphantom{*}}\colon\vphantom{\id}\smash[t]{\id}_{\kos{K}}^{\vphantom{*}}\mor\kos{t}_{\pi}^{\vphantom{*}}$
  be the surjective fiber point of $\kos{t}_{\pi}^{\vphantom{*}}$
  in Lemma~\ref{lem-RepGroobj-Reppi}.
  We recover $\Spec(\pi)$ from $\kos{e}_{\pi}^{\vphantom{*}}$ as
  the representing object of $\uAut_{\kos{K}}(\kos{e}_{\pi}^{\vphantom{*}})$.
  Indeed, by Theorem~\ref{thm-Grocatoverpk-mainthm}~(1),
  $\uAut_{\kos{K}}(\kos{e}_{\pi}^{\vphantom{*}})$
  is represented by
  $\Spec(\kos{e}_{\pi}^*\kos{e}_{\pi*}^{\vphantom{*}}(\kappa))
  =\Spec(\kappa\otimes \pi)$
  equipped with
  $\cp_{\kappa\otimes \pi}$,
  $e_{\kappa\otimes \pi}$,
  and $\varsigma_{\kappa\otimes \pi}$
  in \eqref{eq-GroMainthm-Specpigroupstr}.
  One can check that 
  $\imath_{\pi}\colon\Spec(\kappa\otimes \pi)\cong \Spec(\pi)$
  is an isomorphism of Grothendieck objects in $\kos{K}$.
  Thus we have
  $\Hom_{\cat{Aff}(\kos{K})}(\slot,\Spec(\pi))
  \cong \uAut_{\kos{K}}(\kos{e}_{\pi}^{\vphantom{*}})$;
  each morphism $g\colon\Spec(b)\to\Spec(\pi)$ in $\cat{Aff}(\kos{K})$
  is sent to the $\kos{e}_{\pi}^{\vphantom{*}}\kos{u}_b^{\vphantom{*}}\cong \kos{e}_{\pi}^{\vphantom{*}}\kos{u}_b^{\vphantom{*}}$
  whose component at $(x,\rho_x)\in\cat{Rep}_{\kos{K}}(\pi)$ is
  $\uprho_x(g)\colon\kos{u}_b^*(x)\cong \kos{u}_b^*(x)$
  in \eqref{eq-RepGroobj-actionofg}.
\end{remark}

\begin{remark}
  \label{rem-GroMainthm-flatfiberpt}
  Let $\kos{T}\xmor{\kos{t}}\kos{K}$ be a Grothendieck $\kos{K}$-prekosmos
  and $\omega\colon\vphantom{\id}\smash[t]{\id}_{\kos{K}}^{\vphantom{*}}\mor\kos{t}$
  be a fiber point of $\kos{t}$.
  Assume that $\omega^*$ preserves coreflexive equalizers.
  Following the proof of Theorem~\ref{thm-Grocatoverpk-mainthm},
  we see that
  $\uAut_{\kos{K}}(\omega)$
  is represented by the Grothendieck object
  $\Spec(\pi)$ in $\kos{K}$ with $\pi\coloneqq\omega^*\omega_*(\kappa)$,
  that $\omega$ extends to a Grothendieck $\kos{K}$-morphism
  $\breve{\omega}\colon\kos{t}_{\pi}^{\vphantom{*}}\mor \kos{t}$
  up to an invertible Grothendieck $\kos{K}$-transformation
  $\omega\cong\breve{\omega}\kos{e}_{\pi}^{\vphantom{*}}$,
  and that $\breve{\omega}$ is an embedding.
\end{remark}

We end this section with a remark
that revisits the proof of Theorem~\ref{thm-Grocatoverpk-mainthm}
under a weaker assumption
and relates, for each $X\in\CT$,
the invertibility of the component $\xi_X$ of the universal element $\xi$
to the dualizability of $\breve{\omega}^*(X)$ in $\kosRep_{\kos{K}}(\pi)$.

\begin{remark}
  \label{rem-GroMainthm-dualizableobj}
  Let $\kos{T}\xmor{\kos{t}}\kos{K}$
  be a Grothendieck $\kos{K}$-prekosmos
  and $\omega\colon\vphantom{\id}\smash[t]{\id}_{\kos{K}}^{\vphantom{*}}\mor\kos{t}$
  be a Grothendieck $\kos{K}$-morphism.
  Put $\pi\coloneqq\omega^*\omega_*(\kappa)$.
  Assume that
  \begin{equation*}
    \omega^*(\varphi_{\kos{t}^*(\slot),\slot})
    \colon\omega^*(\kos{t}^*(\slot)\tensor \omega_*(\slot))
    \Rightarrow \omega^*\omega_*(\omega^*\kos{t}^*(\slot)\otimes \slot)
  \end{equation*}
  is a natural isomorphism,
  where $\varphi$ is the natural transformation
  associated to $\omega$ in \eqref{eq-Groprekos-projformula}.
  By Lemma~\ref{lem-ActionbaseGro-projimplystrongequiv}~(1)
  for $\kos{f}=\omega$,
  our assumption is equivalent to the condition that
  the $\kos{K}$-equivariance of the lax $\kos{K}$-tensor functor
  $\omega^*\omega_*\colon\vphantom{\id}\smash[t]{\id}_{\kos{K}}^{\vphantom{*}}\to \vphantom{\id}\smash[t]{\id}_{\kos{K}}^{\vphantom{*}}$
  is a natural isomorphism.
  Then Proposition~\ref{prop-Grocatoverpk-HomisIsom}~(1) and~(2)
  for $\omega^{\pr}=\omega$ still hold,
  and Lemma~\ref{lem-ActionbaseGro-inverseandcomposition}~(2)
  for $\omega^{\pr}=\omega^{\ppr}=\omega$
  also holds with $\uIsom$ replaced by $\uHom$.
  Hence $\Spec(\pi)$,
  equipped with $\cp_{\pi}$, $e_{\pi}$ in \eqref{eq-GroMainthm-Specpigroupstr},
  is an affine monoid scheme in $\kos{K}$
  that represents $\uEnd_{\kos{K}}(\omega)$.
  Let $\xi$ be the universal element of this representation,
  whose component at $X$, regarded as a morphism in $\CK$,
  is $\rho_X$ in \eqref{eq-GroMainthm-rho}.
  Using \eqref{eq-ActionbaseGro-uplowerstarxi} for $\omega^{\pr}=\omega$,
  we deduce that for each $X\in \CT$,
  the following are equivalent:
  \begin{itemize}
    \item
    $\omega^*(\varphi_{X,\kappa})\colon
    \omega^*(X\tensor \omega_*(\kappa))
    \to \omega^*\omega_*(\omega^*(X)\otimes \kappa)$
    is an isomorphism in $\CK$.

    \item
    The component of the universal element $\xi$ at $X$
    is an isomorphism in $\CK_{\pi}^{\vphantom{*}}$.

    \item
    For each $\Spec(b)\in \cat{Aff}(\kos{K})$
    and each Grothendieck $\kos{K}$-transformation
    $\vartheta\colon\omega\kos{u}_b^{\vphantom{*}}\Rightarrow \omega\kos{u}_b^{\vphantom{*}}$,
    the component of $\vartheta$ at $X$
    is an isomorphism in $\CK_b^{\vphantom{*}}$.
  \end{itemize}
  We similarly define a representation of $\Spec(\pi)$ in $\kos{K}$
  as an object in $\CK$ equipped with a functorial action of $\Spec(\pi)$,
  and identify it with a right $\pi$-comodule.
  We still have the symmetric monoidal category
  $\kosRep_{\kos{K}}(\pi)=(\cat{Rep}_{\kos{K}}(\pi),\tensor_{\pi},\unit_{\pi})$,
  the strong monoidal functor
  $\kos{e}_{\pi}^*\colon\kosRep_{\kos{K}}(\pi)\to \kos{K}$,
  the strong monoidal functor
  $\breve{\omega}^*\colon\kos{T}\to \kosRep_{\kos{K}}(\pi)$
  that sends each $X\in \CT$ to $(\omega^*(X),\rho_X)$,
  and the equality
  $\omega^*=\kos{e}_{\pi}^*\breve{\omega}^*$.
  For each object $X$ in $\CT$, the following are equivalent:
  \begin{itemize}
    \item
    $\breve{\omega}^*(X)$ admits a dual in $\kosRep_{\kos{K}}(\pi)$.

    \item
    $\omega^*(X)$ admits a dual in $\kos{K}$
    and the component of the universal element $\xi$ at $X$
    is an isomorphism in $\CK_{\pi}^{\vphantom{*}}$.
  \end{itemize}
  Indeed, we write $x\coloneqq\omega^*(X)$ and $\rho_x\coloneqq\rho_X\colon x\to x\otimes \pi$
  so that $\breve{\omega}^*(X)=(x,\rho_x)$.
  Whenever $\xi_X$ is an isomorphism in $\CK_{\pi}^{\vphantom{*}}$,
  we write $\mathring{\rho}_x\colon x\to x\otimes \pi$
  for the inverse $\xi_X^{-1}$ regarded as a morphism in $\CK$.
  If $(x,\rho_x)$ admits a dual $((x^{\vee},\rho^{\vee}),\ev,\cv)$
  in $\kosRep_{\kos{K}}(\pi)$,
  then $(x^{\vee},\ev,\cv)$ is a dual of $x$ in $\kos{K}$,
  and one can check that the inverse $\xi_X^{-1}$ exists
  with $\mathring{\rho}_x\colon x\to x\otimes \pi$ given as follows.
  \begin{equation*}
    \xymatrix@R=12pt@C=12pt{
      x
      \ar[d]_-{\mathring{\rho}_x}
      \ar[rr]^-{\cv\otimes I_x}
      &\text{ }
      &x\otimes x^{\vee}\otimes x
      \ar[rrr]^-{I_x\otimes \rho^{\vee}\otimes I_x}
      &\text{ }
      &\text{ }
      &x\otimes x^{\vee}\otimes \pi\otimes x
      \ar[d]^-{I_{x\otimes x^{\vee}}\otimes \sigma_{\pi,x}}_-{\cong}
      \\%1
      x\otimes \pi
      &\text{ }
      &\text{ }
      &\text{ }
      &\text{ }
      &x\otimes x^{\vee}\otimes x\otimes \pi
      \ar[lllll]_-{I_x\otimes \ev\otimes I_{\pi}}
      }
  \end{equation*}
  Conversely, assume that $x$ admits a dual $(x^{\vee},\ev,\cv)$ in $\kos{K}$
  and that $\xi_X$ is an isomorphism.
  Then one can check that
  $((x^{\vee},\rho^{\vee}),\ev,\cv)$ is a dual of $(x,\rho_x)$
  in $\kosRep_{\kos{K}}(\pi)$,
  with the right $\pi$-coaction $\rho^{\vee}$ of $x^{\vee}$ given as follows.
  \begin{equation*}
    \xymatrix@R=12pt@C=12pt{
      x^{\vee}
      \ar[d]_-{\rho^{\vee}}
      \ar[rr]^-{I_{x^{\vee}}\otimes \cv}
      &\text{ }
      &x^{\vee}\otimes x\otimes x^{\vee}
      \ar[rrr]^-{I_{x^{\vee}}\otimes \mathring{\rho}_x\otimes I_{x^{\vee}}}
      &\text{ }
      &\text{ }
      &x^{\vee}\otimes x\otimes \pi\otimes x^{\vee}
      \ar[d]^-{I_{x^{\vee}\otimes x}\otimes \sigma_{\pi,x^{\vee}}}_-{\cong}
      \\%1
      x^{\vee}\otimes \pi
      &\text{ }
      &\text{ }
      &\text{ }
      &\text{ }
      &x^{\vee}\otimes x\otimes x^{\vee}\otimes \pi
      \ar[lllll]_-{\ev\otimes I_{x^{\vee}\otimes \pi}}
    }
  \end{equation*}
\end{remark}

\section{Fpqc torsors under a Grothendieck object}
\label{sec-GroTors}

Throughout this section,
$\Spec(\pi)$ is a Grothendieck object in $\kos{K}$.
An \emph{fpqc left} $\Spec(\pi)$\emph{-torsor over $\Spec(\kappa)$}
is a pair $(\Spec(p),\lambda_p)$
consisting of an affine scheme $\Spec(p)$ in $\kos{K}$
that is an fpqc covering of $\Spec(\kappa)$,
and a morphism
$\lambda_p\colon\Spec(\pi)\times \Spec(p)\to \Spec(p)$
in $\cat{Aff}(\kos{K})$
that is a left $\Spec(\pi)$-action,
such that the \emph{shear morphism}
$\tau_p\coloneqq(\lambda_p\times I_p)\circ (I_{\pi}\times \pc_p)
\colon\Spec(\pi)\times \Spec(p)\to \Spec(p)\times \Spec(p)$
is an isomorphism.
We omit $\lambda_p$ from the notation and write $\Spec(p)$.
The \emph{division morphism} of $\Spec(p)$
is $\dvs_p\coloneqq\mathrm{pr}_1\circ\tau_p^{-1}\colon\Spec(p)\times\Spec(p)\to\Spec(\pi)$.
As morphisms in $\cat{Comm}(\kos{K})$,
we have
$\tau_p=(I_{\pi}\otimes \pc_p)\circ (\lambda_p\otimes I_p)
\colon p\otimes p\cong \pi\otimes p$
and 
$\dvs_p=\tau_p^{-1}\circ (I_{\pi}\otimes u_p)
\colon\pi\to p\otimes p$.

The pair consisting of $\Spec(\pi)$
and its product morphism $\cp_{\pi}$ in $\cat{Aff}(\kos{K})$
is an fpqc left $\Spec(\pi)$-torsor over $\Spec(\kappa)$.
As morphisms in $\cat{Comm}(\kos{K})$,
the shear morphism is
$\tau_{\pi}=(I_{\pi}\otimes \pc_{\pi})\circ (\cp_{\pi}\otimes I_{\pi})
\colon\pi\otimes \pi\cong \pi\otimes \pi$,
with its inverse
$\tau_{\pi}^{-1}=(I_{\pi}\otimes \pc_{\pi})\circ (I_{\pi}\otimes \varsigma_{\pi}\otimes I_{\pi})\circ (\cp_{\pi}\otimes I_{\pi})
\colon\pi\otimes \pi\cong \pi\otimes \pi$,
and the division morphism is
$\dvs_{\pi}=(I_{\pi}\otimes \varsigma_{\pi})\circ \cp_{\pi}
\colon\pi\to \pi\otimes \pi$.

An \emph{fpqc right} $\Spec(\pi)$\emph{-torsor over $\Spec(\kappa)$}
is a pair consisting of an affine scheme $\Spec(p)$ in $\kos{K}$
that is an fpqc covering of $\Spec(\kappa)$,
and a morphism 
$\varrho_p\colon\Spec(p)\times \Spec(\pi)\to \Spec(p)$
in $\cat{Aff}(\kos{K})$ that is a right $\Spec(\pi)$-action,
such that 
$(I_p\times \varrho_p)\circ (\pc_p\times I_{\pi})
\colon\Spec(p)\times \Spec(\pi)\to \Spec(p)\times \Spec(p)$
is an isomorphism.
Every fpqc left $\Spec(\pi)$-torsor $\Spec(p)$ over $\Spec(\kappa)$
can be regarded as an fpqc right $\Spec(\pi)$-torsor over $\Spec(\kappa)$,
and vice versa.
As morphisms in $\cat{Comm}(\kos{K})$,
the two $\Spec(\pi)$-actions are related by
$\varrho_p=(I_p\otimes \varsigma_{\pi})\circ \sigma_{\pi,p}\circ \lambda_p\colon p\to p\otimes \pi$
and 
$\lambda_p=(\varsigma_{\pi}\otimes I_p)\circ \sigma_{p,\pi}\circ \varrho_p\colon p\to \pi\otimes p$.

The following lemma collects basic relations
among the morphisms introduced above.
These relations will be used freely in what follows;
we cite individual relations
only when we find it helpful to the reader.

\begin{lemma}
  \label{lem-GroTors-technicalcalculation}
  Let $\Spec(p)$ be an fpqc left $\Spec(\pi)$-torsor over $\Spec(\kappa)$.
  The shear morphism $\tau_p\colon p\otimes p\cong \pi\otimes p$ in $\cat{Comm}(\kos{K})$
  satisfies the following relations:
  \begin{enumerate}
    \item
    $\tau_p\circ (I_p\otimes u_p)=\lambda_p\colon p\to \pi\otimes p$.
    
    \item
    $\tau_p\circ (u_p\otimes I_p)=u_{\pi}\otimes I_p\colon p\to \pi\otimes p$.

    \item
    $\tau_p\circ (I_p\otimes \pc_p)
    =(I_{\pi}\otimes \pc_p)\circ (\tau_p\otimes I_p)
    \colon p\otimes p\otimes p\to \pi\otimes p$.

    \item
    $(\cp_{\pi}\otimes I_p)\circ \tau_p
    =(I_{\pi}\otimes \tau_p)\circ (\lambda_p\otimes I_p)
    \colon p\otimes p\to \pi\otimes \pi\otimes p$.

    \item 
    $(e_{\pi}\otimes I_p)\circ \tau_p=\pc_p\colon p\otimes p\to p$.
    
    \item 
    $(I_{\pi}\otimes \tau_p)\circ (\tau_p\otimes I_p)
    =(\tau_{\pi}\otimes I_p)\circ(I_{\pi}\otimes \tau_p)\circ(I_{\pi\otimes p}\otimes \pc_p)\circ (I_{\pi}\otimes \sigma_{p,p}\otimes I_p)\circ (\lambda_p\otimes I_{p\otimes p})
    \colon p\otimes p\otimes p\cong \pi\otimes \pi\otimes p$.
  \end{enumerate}
  The division morphism $\dvs_p\colon \pi\to p\otimes p$ in $\cat{Comm}(\kos{K})$
  satisfies the following relations:
  \begin{enumerate}
    \item[(7)]
    $(I_p\otimes \pc_p)\circ (\dvs_p\otimes I_p)=\tau_p^{-1}
    \colon \pi\otimes p\cong p\otimes p$.
    
    \item[(8)]
    $(\lambda_p\otimes I_p)\circ \dvs_p=(I_{\pi}\otimes\dvs_p)\circ \cp_{\pi}
    \colon \pi\to \pi\otimes p\otimes p$.
    
    \item[(9)]
    $\pc_p\circ \dvs_p=u_p\circ e_{\pi}\colon \pi\to p$.
    
    \item[(10)]
    $\dvs_p\circ \varsigma_{\pi}=\sigma_{p,p}\circ \dvs_p\colon \pi\to p\otimes p$.

    \item[(11)]
    $(\pc_p\otimes I_p)\circ (I_p\otimes \dvs_p)
    =\sigma_{p,p}
    \circ \tau_p^{-1}
    \circ \sigma_{p,\pi}
    \circ (I_p\otimes \varsigma_{\pi})
    \colon p\otimes \pi\cong p\otimes p$.
  \end{enumerate}
  The induced right $\Spec(\pi)$-action
  $\varrho_p=(I_p\otimes \varsigma_{\pi})\circ \sigma_{\pi,p}\circ \lambda_p\colon p\to p\otimes \pi$,
  as a morphism in $\cat{Comm}(\kos{K})$,
  satisfies the following relations:
  \begin{enumerate}
    \item[(12)]
    $(\pc_p\otimes I_{\pi})\circ (I_p\otimes \varrho_p)
    =(I_p\otimes \varsigma_{\pi})
    \circ \sigma_{\pi,p}\circ \tau_p \circ \sigma_{p,p}
    \colon p\otimes p\cong p\otimes \pi$.

    \item[(13)]
    $(\pc_p\otimes I_p)\circ (I_p\otimes \dvs_p)\circ \varrho_p
    =u_p\otimes I_p\colon p\to p\otimes p$.

    \item[(14)]
    $(I_p\otimes \varrho_p)\circ \dvs_p
    =(\dvs_p\otimes I_{\pi})\circ \cp_{\pi}\colon \pi\to p\otimes p\otimes \pi$.
  \end{enumerate}
\end{lemma}
\begin{proof}
  (1)--(5) are immediate.
  We obtain (6)
  using $\lambda_p\circ \pc_p
  =\pc_{\pi\otimes p}\circ (\lambda_p\otimes \lambda_p)
  \colon p\otimes p\to \pi\otimes p$.
  (7), (8), and (9) follow from 
  (3), (4), and (5), respectively.
  We obtain (10) by substituting
  $\varsigma_{\pi}=(e_{\pi}\otimes I_{\pi})\circ \tau_{\pi}^{-1}\circ (I_{\pi}\otimes u_{\pi})
  \colon \pi\cong \pi$
  and then using (6) in its equivalent form
  $(I_{\pi}\otimes \tau_p^{-1})\circ (\tau_{\pi}^{-1}\otimes I_p)
  =(I_{\pi\otimes p}\otimes \pc_p)\circ (I_{\pi}\otimes \sigma_{p,p}\otimes I_p)\circ (\lambda_p\otimes I_{p\otimes p})\circ (\tau_p^{-1}\otimes I_p)\circ (I_{\pi}\otimes \tau_p^{-1})
  \colon \pi\otimes \pi\otimes p\cong \pi\otimes p\otimes p$.
  We obtain (11) by rewriting $\tau_p^{-1}$ using (7) and then using (10).
  (12) is immediate.
  We see directly from (11) and (12) that
  $(\pc_p\otimes I_p)\circ (I_p\otimes \dvs_p)\colon p\otimes \pi\cong p\otimes p$
  and
  $(\pc_p\otimes I_{\pi})\circ (I_p\otimes \varrho_p)\colon p\otimes p\cong p\otimes \pi$
  are inverse to each other.
  Thus we obtain (13).
  (14) is equivalent to
  $(\tau_p\otimes I_{\pi})\circ (I_p\otimes \varrho_p)\circ \dvs_p
  =(I_{\pi}\otimes u_p\otimes I_{\pi})\circ \cp_{\pi}
  \colon \pi\to \pi\otimes p\otimes \pi$,
  and one can check this using (8), (10), and (12).
\end{proof}

\begin{remark}
  \label{rem-GroTors-fpqcexactseq}
  Let $\Spec(p)$ be an fpqc left $\Spec(\pi)$-torsor over $\Spec(\kappa)$
  and $z\in\CK$.
  Then $I_z\otimes u_p\colon z\to z\otimes p$
  is an equalizer in $\CK$ of the coreflexive pair
  $I_z\otimes u_p\otimes I_p$, $I_{z\otimes p}\otimes u_p\colon 
  z\otimes p\to z\otimes p\otimes p$,
  since the diagram becomes a split equalizer diagram
  after applying $\slot\otimes p$;
  the argument is the same as that used
  in proving the exactness of the Amitsur complex
  of a faithfully flat morphism of commutative rings.
  By composing the pair with
  $I_z\otimes \tau_p\colon z\otimes p\otimes p\cong z\otimes \pi\otimes p$,
  we deduce that $I_z\otimes u_p$
  is an equalizer in $\CK$ of the coreflexive pair
  $I_z\otimes \lambda_p$, $I_z\otimes u_{\pi}\otimes I_p
  \colon z\otimes p\to z\otimes \pi\otimes p$.
  In particular, $u_p\colon \Spec(p)\to\Spec(\kappa)$
  is compatible with the left $\Spec(\pi)$-actions
  when $\Spec(\kappa)$ is equipped with the trivial action
  $\mathrm{pr}_2\colon \Spec(\pi)\times \Spec(\kappa)\to\Spec(\kappa)$.
\end{remark}

We denote by $(\Spec(\pi)\text{-}\cat{Tors})_{\kappa}$
the category of fpqc left $\Spec(\pi)$-torsors over $\Spec(\kappa)$,
whose morphisms are $\Spec(\pi)$-equivariant morphisms
between the underlying affine schemes in $\kos{K}$.
Note that $(\Spec(\pi)\text{-}\cat{Tors})_{\kappa}$ is a groupoid:
for each morphism
$f\colon \Spec(q)\to \Spec(p)$ in $(\Spec(\pi)\text{-}\cat{Tors})_{\kappa}$,
one can check that $f\otimes I_q\colon p\otimes q\to q\otimes q$
is inverse to $(I_p\otimes \pc_q)\circ (I_p\otimes f\otimes I_q)\circ (\dvs_p\otimes I_q)\circ \tau_q\colon q\otimes q\to p\otimes q$,
and since $\slot\otimes q$ is conservative, $f$ is an isomorphism.
The category of fpqc right $\Spec(\pi)$-torsors over $\Spec(\kappa)$
is also a groupoid and is isomorphic to $(\Spec(\pi)\text{-}\cat{Tors})_{\kappa}$.

Let $\Spec(\pi^{\pr})$ also be a Grothendieck object in $\kos{K}$.
An \emph{fpqc $(\Spec(\pi),\Spec(\pi^{\pr}))$-bitorsor over $\Spec(\kappa)$}
is a tuple $(\Spec(p),\lambda_p,\varrho^{\pr}_p)$
consisting of an affine scheme $\Spec(p)$ in $\kos{K}$
that is an fpqc covering of $\Spec(\kappa)$,
a left $\Spec(\pi)$-action $\lambda_p$,
and a right $\Spec(\pi^{\pr})$-action $\varrho^{\pr}_p$, such that
(i) $\Spec(p)$ with $\lambda_p$
is an fpqc left $\Spec(\pi)$-torsor over $\Spec(\kappa)$;
(ii) $\Spec(p)$ with $\varrho^{\pr}_p$
is an fpqc right $\Spec(\pi^{\pr})$-torsor over $\Spec(\kappa)$;
and
(iii) the actions $\lambda_p$ and $\varrho^{\pr}_p$ are compatible,
in the sense that
$(\lambda_p\otimes I_{\pi^{\pr}})\circ \varrho^{\pr}_p=(I_{\pi}\otimes \varrho^{\pr}_p)\circ \lambda_p\colon p\to \pi\otimes p\otimes \pi^{\pr}$.
We omit $\lambda_p$, $\varrho^{\pr}_p$ from the notation
and write $\Spec(p)$.

\begin{proposition}\label{prop-GroTors-gaugegroup}
  Let $\Spec(p)$ be an fpqc left $\Spec(\pi)$-torsor over $\Spec(\kappa)$.
  We have a Grothendieck object $\Spec(\pi^{p})$ in $\kos{K}$
  that acts on $\Spec(p)$ from the right
  and makes $\Spec(p)$ an fpqc $(\Spec(\pi),\Spec(\pi^p))$-bitorsor over $\Spec(\kappa)$.
\end{proposition}

We refer to $\Spec(\pi^p)$ as the \emph{gauge group} of $\Spec(p)$.

\begin{proof}[Proof of Proposition~\ref{prop-GroTors-gaugegroup}]
  Let $\varrho_p$ be the induced right $\Spec(\pi)$-action on $\Spec(p)$.
  We fix an equalizer $\mathrm{eq}_p\colon \pi^p\hookrightarrow p\otimes p$
  in $\cat{Comm}(\kos{K})$ of the coreflexive pair
  $\varrho_p\otimes I_p$, $I_p\otimes \lambda_p\colon p\otimes p\to p\otimes \pi\otimes p$,
  with common retraction $I_p\otimes e_{\pi}\otimes I_p$.
  We write $\varrho^p_p\colon p\to p\otimes \pi^p$
  for the unique morphism in $\cat{Comm}(\kos{K})$ that satisfies 
  $(I_p\otimes \mathrm{eq}_p)\circ \varrho^p_p
  =(\dvs_p\otimes I_p)\circ \lambda_p\colon p\to p\otimes p\otimes p$;
  one can check that $\varrho^p_p$ exists
  using Lemma~\ref{lem-GroTors-technicalcalculation}~(14).
  Moreover, $\lambda_p$ and $\varrho^p_p$ are compatible;
  that is, we have 
  $(\lambda_p\otimes I_{\pi^p})\circ \varrho^p_p
  =(I_{\pi}\otimes \varrho^p_p)\circ \lambda_p\colon 
  p\to \pi\otimes p\otimes \pi^p$,
  as this can be checked after postcomposing with
  the monomorphism
  $I_{\pi \otimes p}\otimes \mathrm{eq}_p$.
  
  Put
  $\Xi_p\coloneqq
  (I_p\otimes \pc_p)\circ (\sigma_{p,p}\otimes I_p)\circ (\mathrm{eq}_p\otimes I_p)
  \colon \pi^p\otimes p\to p\otimes p$
  and
  $\Xi_p^{-1}\coloneqq
  (I_{\pi^p}\otimes \pc_p)\circ (\sigma_{p,\pi^p}\otimes I_p)\circ (\varrho^p_p\otimes I_p)
  \colon p\otimes p\to \pi^p\otimes p$.
  One can check that
  $(I_p\otimes \pc_p)\circ (\sigma_{p,p}\otimes I_p)\circ (\mathrm{eq}_p\otimes I_p)\circ \sigma_{p,\pi^p}\circ \varrho^p_p
  =I_p\otimes u_p\colon p\to p\otimes p$
  and that
  $(I_{\pi^p}\otimes \pc_p)\circ(\sigma_{p,\pi^p}\otimes I_p)\circ (\varrho^p_p\otimes I_p)\circ \sigma_{p,p}\circ \mathrm{eq}_p
  =I_{\pi^p}\otimes u_p\colon \pi^p\to \pi^p\otimes p$;
  the latter can be checked after postcomposing
  with $\mathrm{eq}_p\otimes I_p$
  using Lemma~\ref{lem-GroTors-technicalcalculation}~(13).
  It thus follows that $\Xi_p$ and $\Xi_p^{-1}$ are inverse to each other.
  From this and the fact that $\Spec(p)$ is an fpqc covering of $\Spec(\kappa)$,
  we deduce that $\Spec(\pi^p)$ is an fpqc covering of $\Spec(\kappa)$.

  The coproduct $\cp_{\pi^p}$ is the unique morphism in $\cat{Comm}(\kos{K})$
  that satisfies $(\mathrm{eq}_p\otimes I_{\pi^p})\circ \cp_{\pi^p}=(I_p\otimes \varrho^p_p)\circ \mathrm{eq}_p
  \colon \pi^p\to p\otimes p\otimes \pi^p$;
  one can check that $\cp_{\pi^p}$ exists
  using the fact that $\lambda_p$ and $\varrho^p_p$ are compatible.
  We have
  $(I_p\otimes \cp_{\pi^p})\circ \varrho^p_p
  =(\varrho^p_p\otimes I_{\pi^p})\circ \varrho^p_p\colon 
  p\to p\otimes \pi^p\otimes \pi^p$,
  as this can be checked after postcomposing with
  $I_p\otimes \mathrm{eq}_p\otimes I_{\pi^p}$.
  From this we deduce that $\cp_{\pi^p}$ is coassociative,
  as the relation can be checked after postcomposing with
  $\mathrm{eq}_p\otimes I_{\pi^p\otimes \pi^p}$.

  The counit $e_{\pi^p}$ is the unique morphism in $\cat{Comm}(\kos{K})$
  that satisfies $u_p\circ e_{\pi^p}=\pc_p\circ \mathrm{eq}_p\colon \pi^p\to p$;
  one can check that $e_{\pi^p}$ exists
  using Lemma~\ref{lem-GroTors-technicalcalculation}~(1), (7), and~(13).
  We have
  $(I_p\otimes e_{\pi^p})\circ \varrho^p_p=I_p\colon p\to p$,
  as this can be checked after postcomposing with $I_p\otimes u_p$.
  From this we deduce that $e_{\pi^p}$ is right counital,
  as the relation can be checked
  after postcomposing with $\mathrm{eq}_p$.
  We also have the left counital relation of $e_{\pi^p}$,
  as the relation can be checked after postcomposing with
  $(I_p\otimes \mathrm{eq}_p)\circ (u_p\otimes I_{\pi^p})$
  using Lemma~\ref{lem-GroTors-technicalcalculation}~(13).

  We have shown that $\Spec(\pi^p)$ is an affine monoid scheme in $\kos{K}$,
  that it is an fpqc covering of $\Spec(\kappa)$,
  that $\varrho^p_p$ is a right $\Spec(\pi^p)$-action on $\Spec(p)$,
  and that the two actions $\lambda_p$ and $\varrho^p_p$ 
  on $\Spec(p)$ are compatible.
  One can check that
  $(\pc_p\otimes I_{\pi^p})\circ (I_p\otimes \varrho^p_p)
  =\sigma_{\pi^p,p}\circ \Xi_p^{-1}\circ \sigma_{p,p}\colon p\otimes p\cong p\otimes \pi^p$,
  from which we deduce that
  $(I_p\times \varrho^p_p)\circ (\pc_p\times I_{\pi^p})
  \colon \Spec(p)\times \Spec(\pi^p)\to \Spec(p)\times \Spec(p)$
  is an isomorphism.
  
  It remains to show that $\pi^p$ has an antipode.
  The antipode of $\pi^p$ exists if and only if 
  $\nu\coloneqq (I_{\pi^p}\otimes \pc_{\pi^p})
  \circ(\cp_{\pi^p}\otimes I_{\pi^p})
  \colon \pi^p\otimes \pi^p\to \pi^p\otimes \pi^p$
  is an isomorphism,
  in which case the antipode is given by
  $\varsigma_{\pi^p}
  \coloneqq(e_{\pi^p}\otimes I_{\pi^p})\circ \nu^{-1}\circ (I_{\pi^p}\otimes u_{\pi^p})
  \colon \pi^p\to \pi^p$.
  Since $\slot\otimes p$ is conservative,
  the existence of $\varsigma_{\pi^p}$
  is also equivalent to the condition that
  $\nu\otimes I_p\colon 
  \pi^p\otimes \pi^p\otimes p\to \pi^p\otimes \pi^p\otimes p$
  is an isomorphism.

  Put
  $\mathring{\Upsilon}_p\coloneqq
  (I_{\pi^p}\otimes \sigma_{p,p})
  \circ (I_{\pi^p}\otimes \mathrm{eq}_p)
  \circ \cp_{\pi^p}
  \colon \pi^p\to \pi^p\otimes p\otimes p$.
  One can check that
  $\Xi_p\circ (\pc_{\pi^p}\otimes I_p)
  =\pc_{p\otimes p}\circ (\sigma_{p,p}\otimes I_{p\otimes p})
  \circ (\mathrm{eq}_p\otimes I_{p\otimes p})
  \circ (I_{\pi^p}\otimes \Xi_p)
  \colon \pi^p\otimes \pi^p\otimes p\to p\otimes p$,
  and from this we deduce that
  $(I_{\pi^p}\otimes \pc_{p\otimes p})
  \circ (\mathring{\Upsilon}_p\otimes I_{p\otimes p})
  \circ (I_{\pi^p}\otimes \Xi_p)
  =(I_{\pi^p}\otimes \Xi_p)\circ (\nu\otimes I_p)
  \colon \pi^p\otimes \pi^p\otimes p
  \to \pi^p\otimes p\otimes p$.
  Hence $\varsigma_{\pi^p}$ exists if and only if
  $(I_{\pi^p}\otimes \pc_{p\otimes p})\circ (\mathring{\Upsilon}_p\otimes I_{p\otimes p})
  \colon \pi^p\otimes p\otimes p\to \pi^p\otimes p\otimes p$
  is an isomorphism.

  Put
  $\Upsilon_p\coloneqq
  (I_{\pi^p}\otimes \mathrm{eq}_p)\circ \cp_{\pi^p}
  \colon \pi^p\to \pi^p\otimes p\otimes p$.
  One can check that
  $(I_p\otimes \pc_p\otimes I_p)\circ (\mathrm{eq}_p\otimes \mathrm{eq}_p)\circ \cp_{\pi^p}
  =(I_p\otimes u_p\otimes I_p)\circ \mathrm{eq}_p\colon \pi^p\to p\otimes p\otimes p$.
  From this, along with the coassociativity of $\cp_{\pi^p}$
  and the right counital relation of $e_{\pi^p}$,
  we deduce that
  $(I_{\pi^p}\otimes \pc_{p\otimes p})\circ (\mathring{\Upsilon}_p\otimes I_{p\otimes p})\circ \Upsilon_p
  =(I_{\pi^p}\otimes \pc_{p\otimes p})\circ (\Upsilon_p\otimes I_{p\otimes p})\circ \mathring{\Upsilon}_p
  =I_{\pi^p}\otimes u_{p\otimes p}\colon 
  \pi^p\to \pi^p\otimes p\otimes p$.
  It thus follows that
  $(I_{\pi^p}\otimes \pc_{p\otimes p})\circ (\Upsilon_p\otimes I_{p\otimes p})$
  and 
  $(I_{\pi^p}\otimes \pc_{p\otimes p})\circ (\mathring{\Upsilon}_p\otimes I_{p\otimes p})$
  are inverse to each other.
  This shows that $\pi^p$ has an antipode.
  We conclude that
  $\Spec(\pi^p)$ is a Grothendieck object in $\kos{K}$
  and that 
  $\Spec(p)$ is an fpqc $(\Spec(\pi),\Spec(\pi^p))$-bitorsor over $\Spec(\kappa)$.
\end{proof}

\section[Fiber points and fpqc torsors]{A perfect correspondence between fiber points and fpqc torsors}
\label{sec-pFibGroTors}

Recall from Corollary~\ref{cor-Grocatoverpk-pFibgroupoid}
that given any Grothendieck $\kos{K}$-prekosmos
$\kos{T}\xmor{\kos{t}}\kos{K}$,
the category $\cat{pFib}_{\kos{K}}(\kos{t})$ of fiber points of $\kos{t}$
is a groupoid.
The goal of this section is to prove the following theorem.

\begin{theorem}
  \label{thm-pFibGroTors-perfcorres}
  Let $\kos{T}\xmor{\kos{t}}\kos{K}$ be a Grothendieck $\kos{K}$-category.
  For each surjective fiber point $\omega$ of $\kos{t}$,
  we have an adjoint equivalence of groupoids
  $\cat{pFib}_{\kos{K}}(\kos{t})\simeq (\Spec(\pi)\text{-}\cat{Tors})_{\kappa}^{\op}$,
  where $\Spec(\pi)$ with $\pi\coloneqq\omega^*\omega_*(\kappa)$
  is the Grothendieck object in $\kos{K}$
  that represents $\uAut_{\kos{K}}(\omega)$
  in Theorem~\ref{thm-Grocatoverpk-mainthm}.
\end{theorem}

The proof of Theorem~\ref{thm-pFibGroTors-perfcorres}
is based on the following proposition.

\begin{proposition}
  \label{prop-pFibGroTors-twistedfib}
  Let $\Spec(\pi)$ be a Grothendieck object in $\kos{K}$
  and $\Spec(p)$ be an fpqc left $\Spec(\pi)$-torsor over $\Spec(\kappa)$.
  Recall from Lemma~\ref{lem-RepGroobj-Reppi}
  the surjective fiber point
  $\kos{e}_{\pi}^{\vphantom{*}}\colon \vphantom{\id}\smash[t]{\id}_{\kos{K}}^{\vphantom{*}}\mor\kos{t}_{\pi}^{\vphantom{*}}$
  of $\kos{t}_{\pi}^{\vphantom{*}}$.
  \begin{enumerate}
    \item 
    We have a fiber point
    $\omega^p\colon \vphantom{\id}\smash[t]{\id}_{\kos{K}}^{\vphantom{*}}\mor\kos{t}_{\pi}^{\vphantom{*}}$
    of $\kos{t}_{\pi}^{\vphantom{*}}$,
    whose left adjoint part $\omega^{p*}$
    sends each $(x,\rho_x)\in \cat{Rep}_{\kos{K}}(\pi)$
    to its $\Spec(p)$-twisted form in $\CK$.
    
    \item
    The presheaf
    $\uIsom_{\kos{K}}(\omega^p,\kos{e}_{\pi}^{\vphantom{*}})$
    on $\cat{Aff}(\kos{K})$
    is represented by
    $\Spec(p)$.
  \end{enumerate}
  \begin{equation*}
    \begin{tikzcd}[row sep=5pt,column sep=4pt]
      \kos{K}
      \ar[rr,filled,"\omega^p"]
      \ar[ddr,equal]
      &\text{ }
      &\kosRep_{\kos{K}}(\pi)
      \ar[ddl,filled,"\kos{t}_{\pi}"]
      \\%1
      \text{ }
      \ar[rr,2cell,pos=0.38,shift right=0.3ex,"{\,\,\what{\omega^p}}"]
      &\text{ }
      &\text{ }
      \\%2
      \text{ }
      &\kos{K}
      &\text{ }
    \end{tikzcd}
    \qquad
    \Hom_{\cat{Aff}(\kos{K})}(\slot,\Spec(p))
    \cong
    \uIsom_{\kos{K}}(\omega^p,\kos{e}_{\pi}^{\vphantom{*}})
  \end{equation*}
\end{proposition}

We refer to $\omega^p\colon \vphantom{\id}\smash[t]{\id}_{\kos{K}}^{\vphantom{*}}\mor\kos{t}_{\pi}^{\vphantom{*}}$
as the \emph{$\Spec(p)$-twisted fiber point} of $\kos{t}_{\pi}^{\vphantom{*}}$.

\begin{proof}[Proof of Proposition~\ref{prop-pFibGroTors-twistedfib}]
  We write
  $X=(x,\rho_x)$, $Y=(y,\rho_y)$
  for objects in $\cat{Rep}_{\kos{K}}(\pi)$.
  For each $X$, we fix an equalizer
  $\mathrm{eq}^p_X\colon \omega^{p*}(X)\hookrightarrow x\otimes p$
  in $\CK$ of the coreflexive pair
  $\rho_x\otimes I_p$, $I_x\otimes \lambda_p\colon x\otimes p\to x\otimes \pi\otimes p$,
  with common retraction 
  $I_x\otimes e_{\pi}\otimes I_p$.
  The object $\omega^{p*}(X)$ in $\CK$
  is what we mean by the $\Spec(p)$-twisted form of $X$ in (1).
  Then $\omega^{p*}$ is a functor,
  sending each morphism $l\colon X\to Y$ in $\cat{Rep}_{\kos{K}}(\pi)$
  to the unique morphism
  $\omega^{p*}(l)\colon \omega^{p*}(X)\to \omega^{p*}(Y)$
  in $\CK$ that satisfies
  $\mathrm{eq}^p_Y\circ \omega^{p*}(l)=(l\otimes I_p)\circ \mathrm{eq}^p_X
  \colon \omega^{p*}(X)\to y\otimes p$.
  When we regard
  $\mathrm{eq}^p_X\colon \omega^{p*}(X)\hookrightarrow x\otimes p=\kos{e}_{\pi}^*(X)\otimes p$
  as a morphism in the Kleisli category $\CK_p^{\vphantom{*}}$,
  we denote it by
  $\xi^p_X\colon \kos{u}_p^*\omega^{p*}(X)\to \kos{u}_p^*\kos{e}_{\pi}^*(X)$.

  We first show that $\xi^p_X$ is an isomorphism in $\CK_p^{\vphantom{*}}$ for every $X$.
  It suffices to show that
  $\kos{u}_{p*}^{\vphantom{*}}(\xi^p_X)=(I_x\otimes \pc_p)\circ (\mathrm{eq}^p_X\otimes I_p)
  \colon \omega^{p*}(X)\otimes p\to x\otimes p$
  is an isomorphism in $\CK$,
  since $\kos{u}_{p*}^{\vphantom{*}}\colon \CK_p^{\vphantom{*}}\to \CK$ is conservative.
  In the diagram
  \begin{equation*}
    \vcenter{\hbox{
      \xymatrix@R=20pt@C=50pt{
        \omega^{p*}(X)\otimes p
        \ar@{^{(}->}[r]^-{\mathrm{eq}^p_X\otimes I_p}
        \ar@{.>}[d]_-{\exists!\, \kos{u}_{p*}(\xi^p_X)}^-{\cong}
        &x\otimes p\otimes p
        \ar@<0.5ex>[r]^-{\rho_x\otimes I_{p\otimes p}}
        \ar@<-0.5ex>[r]_-{I_x\otimes \lambda_p\otimes I_p}
        \ar[d]_-{I_x\otimes \tau_p}^-{\cong}
        &x\otimes \pi\otimes p\otimes p
        \ar[d]^-{I_{x\otimes \pi}\otimes \tau_p}_-{\cong}
        \\%1
        x\otimes p
        \ar[r]^-{\rho_x\otimes I_p}
        &x\otimes \pi\otimes p
        \ar@<0.5ex>[r]^-{\rho_x\otimes I_{\pi\otimes p}}
        \ar@<-0.5ex>[r]_-{I_x\otimes \cp_{\pi}\otimes I_p}
        &x\otimes \pi\otimes \pi\otimes p
      }
    }}
  \end{equation*}
  we see that
  (i) the top horizontal morphisms form an equalizer diagram,
  since $\slot\otimes p$ preserves coreflexive equalizers;
  (ii) the bottom horizontal morphisms
  form a split equalizer diagram
  with respect to the splitting
  $I_x\otimes e_{\pi}\otimes I_p$,
  $I_{x\otimes \pi}\otimes e_{\pi}\otimes I_p$;
  (iii) the middle and right vertical isomorphisms
  commute with the horizontal parallel pairs;
  and (iv)
  the unique isomorphism between the two equalizers
  that makes the left square commute
  is given by
  $(I_x\otimes e_{\pi}\otimes I_p)\circ (I_x\otimes \tau_p)\circ (\mathrm{eq}^p_X\otimes I_p)
  =\kos{u}_{p*}^{\vphantom{*}}(\xi^p_X)$.
  Thus $\kos{u}_{p*}^{\vphantom{*}}(\xi^p_X)$ is an isomorphism in $\CK$,
  and hence $\xi^p_X$ is an isomorphism in $\CK_p^{\vphantom{*}}$.
  
  We next show that $\omega^{p*}$ is a strong monoidal functor.
  Given $X$ and $Y$, recall the coherence morphism
  $(\otimes p)_{\tensor,x,y}\colon x\otimes p\otimes y\otimes p\to x\otimes y\otimes p$
  of the lax monoidal endofunctor
  $\slot\otimes p$ on $\kos{K}$.
  \begin{itemize}
    \item 
    Let $(\omega^{p*}_{\tensor})_{X,Y}
    \colon \omega^{p*}(X)\otimes \omega^{p*}(Y)\to
    \omega^{p*}(X\tensor_{\pi} Y)$
    be the unique morphism in $\CK$ that satisfies 
    $\mathrm{eq}^p_{X\tensor_{\pi}Y}\circ (\omega^{p*}_{\tensor})_{X,Y}
    =(\otimes p)_{\tensor,x,y}\circ (\mathrm{eq}^p_X\otimes \mathrm{eq}^p_Y)$.
    One can check that $(\omega^{p*}_{\tensor})_{X,Y}$ exists
    and that
    $(\sigma_{y,x}\otimes I_p)
    \circ (I_y\otimes \kos{u}_{p*}^{\vphantom{*}}(\xi^p_X))
    \circ (\sigma_{\omega^{p*}(X),y}\otimes I_p)
    \circ (I_{\omega^{p*}(X)}\otimes \kos{u}_{p*}^{\vphantom{*}}(\xi^p_Y))
    = \kos{u}_{p*}^{\vphantom{*}}(\xi^p_{X\tensor_{\pi}Y})
    \circ ((\omega^{p*}_{\tensor})_{X,Y}\otimes I_p)
    \colon \omega^{p*}(X)\otimes \omega^{p*}(Y)\otimes p
    \cong x\otimes y\otimes p$.
    Hence $(\omega^{p*}_{\tensor})_{X,Y}\otimes I_p$ is an isomorphism,
    and since $\slot\otimes p$ is conservative,
    we deduce that $(\omega^{p*}_{\tensor})_{X,Y}$ is an isomorphism.
    
    \item
    Let 
    $\omega^{p*}_{\unit}\colon \kappa\to \omega^{p*}(\unit_{\pi})$
    be the unique morphism in $\CK$ that satisfies
    $\mathrm{eq}^p_{\unit_{\pi}}\circ \omega^{p*}_{\unit}=I_{\kappa}\otimes u_p$.
    The morphism $\omega^{p*}_{\unit}$ exists and is an isomorphism,
    since $u_p\colon \kappa \to p$
    is an equalizer of the coreflexive pair
    $\lambda_p$, $u_{\pi}\otimes I_p\colon p\to \pi\otimes p$
    by Remark~\ref{rem-GroTors-fpqcexactseq}.
  \end{itemize}
  One can check that
  $\omega^{p*}_{\tensor}$, $\omega^{p*}_{\unit}$
  are compatible with the coherence isomorphisms of
  $\kosRep_{\kos{K}}(\pi)$ and $\kos{K}$,
  using the compatibility relations that
  $(\otimes p)_{\tensor}$, $(\otimes p)_{\unit}$ satisfy.
  Thus $\omega^{p*}$ equipped with $\omega^{p*}_{\tensor}$, $\omega^{p*}_{\unit}$
  is a strong monoidal functor.

  We then construct the Grothendieck morphism
  $\omega^p=(\omega^{p*},\omega^p_*,\eta^p,\epsilon^p)
  \colon \kos{K}\mor \kosRep_{\kos{K}}(\pi)$.
  Let $\varrho_p\colon p\to p\otimes \pi$
  be the induced right $\Spec(\pi)$-action on $\Spec(p)$.
  \begin{itemize}
    \item
    The right adjoint part $\omega^p_*$ sends each $z\in\CK$
    to the pair $(z\otimes p,I_z\otimes \varrho_p)$.

    \item 
    The component $\eta^p_X\colon X\to \omega^p_*\omega^{p*}(X)$
    of $\eta^p$ at each $X$
    is the unique morphism
    in $\cat{Rep}_{\kos{K}}(\pi)$
    whose underlying morphism
    $\eta^p_X\colon x\to \omega^{p*}(X)\otimes p$
    in $\CK$ satisfies
    $(\mathrm{eq}^p_X\otimes I_p)\circ \eta^p_X
    =(I_x\otimes \dvs_p)\circ \rho_x\colon x\to x\otimes p\otimes p$;
    the existence of the morphism $\eta^p_X$ in $\CK$
    follows from Lemma~\ref{lem-GroTors-technicalcalculation}~(8);
    $\eta^p_X$ is a morphism in $\cat{Rep}_{\kos{K}}(\pi)$,
    since the compatibility relation can be checked
    after postcomposing with
    $\mathrm{eq}^p_X\otimes I_{p\otimes \pi}$
    using Lemma~\ref{lem-GroTors-technicalcalculation}~(14).

    \item
    The component $\epsilon^p_z\colon \omega^{p*}\omega^p_*(z)\to z$
    of $\epsilon^p$ at each $z\in \CK$ is the unique morphism in $\CK$
    that satisfies 
    $(I_z\otimes u_p)\circ \epsilon^p_z
    =(I_z\otimes \pc_p)\circ \mathrm{eq}^p_{\omega^p_*(z)}
    \colon \omega^{p*}\omega^p_*(z)\to z\otimes p$;
    one can check that $\epsilon^p_z$ exists
    using Lemma~\ref{lem-GroTors-technicalcalculation}~(1), (7), and~(13).

    \item
    We have one of the triangle identities
    $\epsilon^p_{\omega^{p*}(X)}\circ \omega^{p*}(\eta^p_X)=I_{\omega^{p*}(X)}$,
    as this can be checked after postcomposing with
    $(\mathrm{eq}^p_X\otimes I_p)\circ (I_{\omega^{p*}(X)}\otimes u_p)$.
    We also have the other triangle identity
    $\omega^p_*(\epsilon^p_z)\circ \eta^p_{\omega^p_*(z)}
    =I_{\omega^p_*(z)}$,
    as this can be checked, as an equality in $\CK$,
    after postcomposing with
    $I_z\otimes u_p\otimes I_p$.
  \end{itemize}
  Thus $\omega^p\colon \kos{K}\mor\kosRep_{\kos{K}}(\pi)$
  is a Grothendieck morphism.
  The right adjoint part $\omega^p_*$ is lax monoidal by doctrinal adjunction;
  we have $\kos{e}_{\pi}^*\omega^p_*=\slot\otimes p$
  as lax monoidal endofunctors on $\kos{K}$.
  
  We now lift $\omega^p$ to a Grothendieck $\kos{K}$-morphism
  $\vphantom{\id}\smash[t]{\id}_{\kos{K}}^{\vphantom{*}}\mor\kos{t}_{\pi}^{\vphantom{*}}$.
  For $z\in \CK$, 
  let $(\what{\omega^{p*}})_z\colon z\to \omega^{p*}\kos{t}_{\pi}^*(z)$
  be the unique morphism in $\CK$
  that satisfies $\mathrm{eq}^p_{\kos{t}_{\pi}^*(z)}\circ (\what{\omega^{p*}})_z=I_z\otimes u_p$;
  the morphism $(\what{\omega^{p*}})_z$ exists and is an isomorphism,
  since $I_z\otimes u_p$ is an equalizer of
  the coreflexive pair $I_z\otimes \lambda_p$, $I_z\otimes u_{\pi}\otimes I_p$
  by Remark~\ref{rem-GroTors-fpqcexactseq}.
  One can check that $\what{\omega^{p*}}\colon \vphantom{\id}\smash[t]{\id}_{\kos{K}}^{\vphantom{*}}\cong \omega^{p*}\kos{t}_{\pi}^*$
  is a monoidal natural isomorphism,
  using the fact that
  $\what{\otimes p}\colon \vphantom{\id}\smash[t]{\id}_{\kos{K}}^{\vphantom{*}}\Rightarrow \slot\otimes p$
  is a monoidal natural transformation.
  Thus $\omega^p$
  is a Grothendieck $\kos{K}$-morphism $\vphantom{\id}\smash[t]{\id}_{\kos{K}}^{\vphantom{*}}\mor\kos{t}_{\pi}^{\vphantom{*}}$,
  with $\what{\omega^p}\colon \vphantom{\id}\smash[t]{\id}_{\kos{K}}^{\vphantom{*}}\cong \kos{t}_{\pi}^{\vphantom{*}}\omega^p$
  given by $\what{\omega^{p*}}\colon \vphantom{\id}\smash[t]{\id}_{\kos{K}}^{\vphantom{*}}\cong \omega^{p*}\kos{t}_{\pi}^*$.
  We have $z\acts X=\kos{t}_{\pi}^*(z)\tensor_{\pi}X=(z\otimes x,I_z\otimes \rho_x)$.
  One can check that the $\kos{K}$-equivariance
  $(\cevar{\omega^{p*}})_{z,X}\colon z\otimes \omega^{p*}(X)\cong \omega^{p*}(z\acts X)$
  of $\omega^{p*}$ is the unique morphism that satisfies
  $\mathrm{eq}^p_{z\acts X}\circ (\cevar{\omega^{p*}})_{z,X}
  =a_{z,x,p}\circ (I_z\otimes \mathrm{eq}^p_X)
  \colon z\otimes \omega^{p*}(X)\to (z\otimes x)\otimes p$.

  Finally, we show that
  $\omega^p\colon \vphantom{\id}\smash[t]{\id}_{\kos{K}}^{\vphantom{*}}\mor\kos{t}_{\pi}^{\vphantom{*}}$
  is a fiber point of $\kos{t}_{\pi}^{\vphantom{*}}$.
  The right adjoint part $\omega^p_*$ is conservative by Remark~\ref{rem-Groprekos-omegalowerstarconservative}.
  It remains to show that 
  $\omega^p$ satisfies the projection formula.
  We write $\varphi^p$ for the natural transformation
  associated to $\omega^p$ in \eqref{eq-Groprekos-projformula}.
  Let $x\in \CK$ and $Z=(z,\rho_z)\in\cat{Rep}_{\kos{K}}(\pi)$.
  The underlying morphism in $\CK$ of the component
  $\varphi^p_{Z,x}\colon Z\tensor_{\pi}\omega^p_*(x)\to \omega^p_*(\omega^{p*}(Z)\otimes x)$
  of $\varphi^p$ at $Z$, $x$ is given by
  $\kos{e}_{\pi}^*(\varphi^p_{Z,x})
  =(\otimes p)_{\tensor,\omega^{p*}(Z),x}\circ (\eta^p_Z\otimes I_{x\otimes p})
  \colon z\otimes x\otimes p\to \omega^{p*}(Z)\otimes x\otimes p$.
  One can check that 
  $(I_z\otimes \pc_p)\circ (\mathrm{eq}^p_Z\otimes I_p)\circ \eta^p_Z
  =I_z\otimes u_p\colon z\to z\otimes p$.
  We also have
  $(I_{\omega^{p*}(Z)}\otimes \pc_p)\circ (\eta^p_Z\otimes I_p)\circ \mathrm{eq}^p_Z
  =I_{\omega^{p*}(Z)}\otimes u_p
  \colon \omega^{p*}(Z)\to \omega^{p*}(Z)\otimes p$,
  as this can be checked after postcomposing with 
  $\mathrm{eq}^p_Z\otimes I_p$.
  From these two equalities we deduce that
  $\kos{e}_{\pi}^*(\varphi^p_{Z,x})$ is inverse to
  $(\otimes p)_{\tensor,z,x}\circ (\mathrm{eq}^p_Z\otimes I_{x\otimes p})
  \colon \omega^{p*}(Z)\otimes x\otimes p\to z\otimes x\otimes p$.
  Hence $\kos{e}_{\pi}^*\varphi^p$ is a natural isomorphism,
  and so is $\varphi^p$, since $\kos{e}_{\pi}^*$ is conservative.
  This proves (1).

  By Proposition~\ref{prop-Grocatoverpk-HomisIsom}
  for $\omega=\kos{e}_{\pi}^{\vphantom{*}}$, $\omega^{\pr}=\omega^p$,
  we deduce that
  $\Spec(\omega^{p*}\kos{e}_{\pi*}^{\vphantom{*}}(\kappa))$
  represents
  $\uHom_{\kos{K}}(\omega^p,\kos{e}_{\pi}^{\vphantom{*}})
  =\uIsom_{\kos{K}}(\omega^p,\kos{e}_{\pi}^{\vphantom{*}})$.
  We have an isomorphism 
  $i_p\colon p\cong \omega^{p*}\kos{e}_{\pi*}^{\vphantom{*}}(\kappa)$
  in $\cat{Comm}(\kos{K})$,
  which is the unique one that satisfies
  $\mathrm{eq}^p_{\kos{e}_{\pi*}^{\vphantom{*}}(\kappa)}\circ i_p
  =(\imath_{\pi}\otimes I_p)\circ \lambda_p
  \colon p\to \kappa\otimes \pi\otimes p$.
  Thus $\Spec(p)$ represents
  $\uIsom_{\kos{K}}(\omega^p,\kos{e}_{\pi}^{\vphantom{*}})$.
  This proves (2).
  For later use,
  we show that the component of the universal element
  $\omega^p\kos{u}_p^{\vphantom{*}}\cong \kos{e}_{\pi}^{\vphantom{*}}\kos{u}_p^{\vphantom{*}}$
  of this representation
  at $X$ is given by
  $\xi^p_X\colon \kos{u}_p^*\omega^{p*}(X)\cong \kos{u}_p^*\kos{e}_{\pi}^*(X)$.
  Indeed, one can check that for $z\in \CK$,
  the isomorphism $(\tahar{\omega^{p*}\kos{e}_{\pi*}^{\vphantom{*}}})_z\colon 
  z\otimes \omega^{p*}\kos{e}_{\pi*}^{\vphantom{*}}(\kappa)
  \cong \omega^{p*}\kos{e}_{\pi*}^{\vphantom{*}}(z)$
  in \eqref{eq-ActionbaseGro-taharupperlowerstar} satisfies
  $\mathrm{eq}^p_{\kos{e}_{\pi*}^{\vphantom{*}}(z)}\circ (\tahar{\omega^{p*}\kos{e}_{\pi*}^{\vphantom{*}}})_z
  =(\jmath_z^{-1}\otimes I_{\pi\otimes p})
  \circ (I_z\otimes \mathrm{eq}^p_{\kos{e}_{\pi*}^{\vphantom{*}}(\kappa)})
  \colon z\otimes \omega^{p*}\kos{e}_{\pi*}^{\vphantom{*}}(\kappa)
  \to z\otimes \pi\otimes p$.
  From this we deduce that the morphism
  $\rho_X\colon \omega^{p*}(X)\to x\otimes \omega^{p*}\kos{e}_{\pi*}^{\vphantom{*}}(\kappa)$
  in \eqref{eq-ActionbaseGro-rhoX} is equal to 
  $(I_x\otimes i_p)\circ \mathrm{eq}^p_X$,
  as this can be checked after postcomposing with
  $I_x\otimes \mathrm{eq}^p_{\kos{e}_{\pi*}^{\vphantom{*}}(\kappa)}$.
  Hence the isomorphisms
  $\xi^p_X\colon \kos{u}_p^*\omega^{p*}(X)\cong \kos{u}_p^*\kos{e}_{\pi}^*(X)$
  in $\CK_p^{\vphantom{*}}$ constitute the universal element
  $\xi^p\colon \omega^p\kos{u}_p^{\vphantom{*}}\cong \kos{e}_{\pi}^{\vphantom{*}}\kos{u}_p^{\vphantom{*}}$.
\end{proof}

The following lemma,
along with Proposition~\ref{prop-pFibGroTors-twistedfib},
will be used in the proof of Theorem~\ref{thm-pFibGroTors-perfcorres}.

\begin{lemma}
  \label{lem-pFibGroTors-torsfibmor}
  Let $\Spec(\pi)$ be a Grothendieck object in $\kos{K}$
  and $f\colon\Spec(q)\cong \Spec(p)$
  be an isomorphism of fpqc left $\Spec(\pi)$-torsors over $\Spec(\kappa)$.
  We have an isomorphism
  $\vartheta^f\colon \omega^p\cong \omega^q$
  between the twisted fiber points
  $\omega^p$, $\omega^q\colon \vphantom{\id}\smash[t]{\id}_{\kos{K}}^{\vphantom{*}}\mor\kos{t}_{\pi}^{\vphantom{*}}$
  in Proposition~\ref{prop-pFibGroTors-twistedfib}.
\end{lemma}
\begin{proof}
  For each $X=(x,\rho_x)\in \cat{Rep}_{\kos{K}}(\pi)$,
  let $\vartheta^f_X\colon \omega^{p*}(X)\to \omega^{q*}(X)$
  be the unique morphism in $\CK$ that satisfies
  $\mathrm{eq}^q_X\circ \vartheta^f_X
  =(I_x \otimes f)\circ \mathrm{eq}^p_X
  \colon \omega^{p*}(X)\to x\otimes q$.
  One can check that the morphisms $\vartheta^f_X$ exist,
  are isomorphisms,
  and constitute an invertible Grothendieck $\kos{K}$-transformation
  $\vartheta^f$ from $\omega^p$ to $\omega^q$.
\end{proof}

\begin{proof}[Proof of Theorem~\ref{thm-pFibGroTors-perfcorres}]
  Let $\kos{T}\xmor{\kos{t}}\kos{K}$ be a Grothendieck $\kos{K}$-category
  and $\omega\colon \vphantom{\id}\smash[t]{\id}_{\kos{K}}^{\vphantom{*}}\mor\kos{t}$
  be a surjective fiber point of $\kos{t}$.
  Put $\pi\coloneqq\omega^*\omega_*(\kappa)$
  and $\varpi\coloneqq \omega^*\omega_*$.
  By Theorem~\ref{thm-Grocatoverpk-mainthm},
  $\uAut_{\kos{K}}(\omega)$
  is represented by the Grothendieck object
  $\Spec(\pi)$ in $\kos{K}$,
  and $\omega$ extends to $\breve{\omega}\colon \kos{t}_{\pi}^{\vphantom{*}}\simeq\kos{t}$
  up to an invertible Grothendieck $\kos{K}$-transformation
  $\omega\cong \breve{\omega}\kos{e}_{\pi}^{\vphantom{*}}$.

  We define the functor
  $\cat{pFib}_{\kos{K}}(\kos{t})\to (\Spec(\pi)\text{-}\cat{Tors})_{\kappa}^{\op}$
  as follows.
  Each fiber point $\omega^{\pr}$ of $\kos{t}$
  is sent to $\Spec(\omega^{\pr*}\omega_*(\kappa))$ equipped with
  $\lambda_{\omega^{\pr*}\omega_*(\kappa)}\coloneqq
  (\tahar{\omega^{\pr*}\omega_*})^{-1}_{\pi}\circ \omega^{\pr*}(\eta_{\omega_*(\kappa)})
  \colon \omega^{\pr*}\omega_*(\kappa)\to \pi\otimes \omega^{\pr*}\omega_*(\kappa)$.
  By Proposition~\ref{prop-Grocatoverpk-HomisIsom},
  $\Spec(\omega^{\pr*}\omega_*(\kappa))$
  represents 
  $\uHom_{\kos{K}}(\omega^{\pr},\omega)
  =\uIsom_{\kos{K}}(\omega^{\pr},\omega)$.
  Observe that $\lambda_{\omega^{\pr*}\omega_*(\kappa)}$
  is the morphism in Lemma~\ref{lem-ActionbaseGro-inverseandcomposition}~(2)
  with $(\omega,\omega^{\pr},\omega^{\ppr})$
  replaced by $(\omega,\omega,\omega^{\pr})$.
  It thus follows that
  $\lambda_{\omega^{\pr*}\omega_*(\kappa)}$
  is a left $\Spec(\pi)$-action,
  and that the shear morphism
  $\tau_{\omega^{\pr*}\omega_*(\kappa)}$
  is an isomorphism.
  Moreover, $\Spec(\omega^{\pr*}\omega_*(\kappa))$
  is an fpqc covering of $\Spec(\kappa)$;
  this follows from Lemma~\ref{lem-ActionbaseGro-inverseandcomposition}~(1)
  and from the conclusion in the proof
  of Lemma~\ref{lem-Grocatoverpk-pFiblocallyisom}
  that the representing object 
  $\Spec(\omega^*\omega^{\pr}_*(\kappa))$
  of $\uIsom_{\kos{K}}(\omega,\omega^{\pr})$
  is an fpqc covering of $\Spec(\kappa)$.
  Thus $\Spec(\omega^{\pr*}\omega_*(\kappa))$
  equipped with $\lambda_{\omega^{\pr*}\omega_*(\kappa)}$
  is indeed an fpqc left $\Spec(\pi)$-torsor over $\Spec(\kappa)$.
  Each isomorphism $\vartheta\colon \omega^{\pr}\cong \omega^{\ppr}$
  between fiber points of $\kos{t}$
  is sent to the $\Spec(\pi)$-equivariant isomorphism
  $\vartheta_{\omega_*(\kappa)}\colon 
  \Spec(\omega^{\ppr*}\omega_*(\kappa))
  \cong 
  \Spec(\omega^{\pr*}\omega_*(\kappa))$
  in $\cat{Aff}(\kos{K})$.
  
  We define the functor
  $(\Spec(\pi)\text{-}\cat{Tors})_{\kappa}^{\op}\to \cat{pFib}_{\kos{K}}(\kos{t})$
  as follows.
  Each fpqc left $\Spec(\pi)$-torsor $\Spec(p)$ over $\Spec(\kappa)$
  is sent to the fiber point
  $\breve{\omega}\omega^p$ of $\kos{t}$,
  where $\omega^p\colon \vphantom{\id}\smash[t]{\id}_{\kos{K}}^{\vphantom{*}}\mor\kos{t}_{\pi}^{\vphantom{*}}$
  is the $\Spec(p)$-twisted fiber point of $\kos{t}_{\pi}^{\vphantom{*}}$
  in Proposition~\ref{prop-pFibGroTors-twistedfib}.
  Each isomorphism $f\colon \Spec(q)\cong \Spec(p)$
  of fpqc left $\Spec(\pi)$-torsors over $\Spec(\kappa)$
  is sent to the isomorphism
  $\breve{\omega}\vartheta^f
  \colon \breve{\omega}\omega^p\cong \breve{\omega}\omega^q$
  of fiber points of $\kos{t}$,
  where $\vartheta^f\colon \omega^p\cong \omega^q$
  is the isomorphism in Lemma~\ref{lem-pFibGroTors-torsfibmor}.

  We show that the functors defined above
  constitute an adjoint equivalence.

  Let $\omega^{\pr}$ be a fiber point of $\kos{t}$.
  Put
  $\varpi^{\pr}\coloneqq\omega^{\pr*}\omega_*$
  and $p\coloneqq \omega^{\pr*}\omega_*(\kappa)$.
  We define an isomorphism 
  $\vartheta^{\omega^{\pr}}\colon \omega^{\pr}\cong \breve{\omega}\omega^p$
  of fiber points of $\kos{t}$ as follows.
  Recall that for $X\in \CT$,
  we have $\breve{\omega}^*(X)=(\omega^*(X),\rho_X)$
  with $\rho_X=(\tahar{\varpi})_{\omega^*(X)}^{-1}\circ \omega^*(\eta_X)
  \colon \omega^*(X)\to \omega^*(X)\otimes \pi$,
  and $\lambda_p=(\tahar{\varpi^{\pr}})^{-1}_{\pi}\circ \omega^{\pr*}(\eta_{\omega_*(\kappa)})
  \colon p\to \pi\otimes p$.
  One can check that in the diagram
  \begin{equation*}
    \vcenter{\hbox{
      \xymatrix@R=20pt@C=15pt{
        \omega^{\pr*}(X)
        \ar[rr]^-{\omega^{\pr*}(\eta_X)}
        \ar@{.>}[d]_-{\exists!\,\vartheta^{\omega^{\pr}}_X}
        &\text{ }
        &\omega^{\pr*}\omega_*\omega^*(X)
        \ar@<0.5ex>[rrr]^-{\omega^{\pr*}\omega_*\omega^*(\eta_X)}
        \ar@<-0.5ex>[rrr]_-{\omega^{\pr*}(\eta_{\omega_*\omega^*(X)})}
        \ar[d]_-{(\tahar{\varpi^{\pr}})_{\omega^*(X)}^{-1}}^-{\cong}
        &\text{ }
        &\text{ }
        &\omega^{\pr*}\omega_*\omega^*\omega_*\omega^*(X)
        \ar[d]^-{(\tahar{\varpi^{\pr}}\tahar{\varpi})_{\omega^*(X)}^{-1}}_-{\cong}
        \\%1
        \omega^{p*}\breve{\omega}^*(X)
        \ar@{^{(}->}[rr]^-{\mathrm{eq}^p_{\breve{\omega}^*(X)}}
        &\text{ }
        &\omega^*(X)\otimes p
        \ar@<0.5ex>[rrr]^-{\rho_X\otimes I_p}
        \ar@<-0.5ex>[rrr]_-{I_{\omega^*(X)}\otimes \lambda_p}
        &\text{ }
        &\text{ }
        &\omega^*(X)\otimes \pi\otimes p
      }
    }}
  \end{equation*}
  the vertical isomorphisms commute with the horizontal parallel pairs,
  using the relation
  $(\tahar{\varpi^{\pr}\varpi})_{\omega^*(X)}\circ (I_{\omega^*(X)}\otimes (\tahar{\varpi^{\pr}})_{\pi})
  =(\tahar{\varpi^{\pr}}\tahar{\varpi})_{\omega^*(X)}
  \colon \omega^*(X)\otimes \pi\otimes p\cong \varpi^{\pr}\varpi\omega^*(X)$
  and \eqref{eq-ActionbaseGro-laxKtnattahar}
  for $\omega^{\pr*}\eta\omega_*\colon \varpi^{\pr}\Rightarrow\varpi^{\pr}\varpi$ with $z=\omega^*(X)$;
  thus the left vertical morphism
  $\vartheta^{\omega^{\pr}}_X$
  exists.
  One can also check that the morphisms
  $\vartheta^{\omega^{\pr}}_X$
  constitute a morphism
  $\vartheta^{\omega^{\pr}}\colon \omega^{\pr}\Rightarrow \breve{\omega}\omega^p$
  of fiber points of $\kos{t}$,
  using the fact that $\tahar{\varpi^{\pr}}\colon \slot\otimes p\cong \varpi^{\pr}$
  is a lax $\kos{K}$-tensor natural isomorphism.
  By Corollary~\ref{cor-Grocatoverpk-pFibgroupoid},
  $\vartheta^{\omega^{\pr}}$ is an isomorphism in $\cat{pFib}_{\kos{K}}(\kos{t})$.
  
  Moreover,
  $\vartheta^{\omega^{\pr}}\colon \omega^{\pr}\cong \breve{\omega}\omega^p$
  is natural in $\omega^{\pr}$.
  Indeed, let $\omega^{\ppr}$ also be a fiber point of $\kos{t}$
  and $\vartheta\colon \omega^{\pr}\cong\omega^{\ppr}$
  be an isomorphism in $\cat{pFib}_{\kos{K}}(\kos{t})$.
  Put $q\coloneqq\omega^{\ppr*}\omega_*(\kappa)$
  and $f\coloneqq\vartheta_{\omega_*(\kappa)}\colon \Spec(q)\cong \Spec(p)$.
  Then we have
  $\vartheta^f_{\breve{\omega}^*(X)}\circ \vartheta^{\omega^{\pr}}_X
  =\vartheta^{\omega^{\ppr}}_X\circ \vartheta_X
  \colon \omega^{\pr*}(X)\cong \omega^{q*}\breve{\omega}^*(X)$,
  as this can be checked after postcomposing with
  $\mathrm{eq}^q_{\breve{\omega}^*(X)}$
  using \eqref{eq-ActionbaseGro-laxKtnattahar}
  for $\vartheta\omega_*\colon \omega^{\pr*}\omega_*\cong\omega^{\ppr*}\omega_*$
  with $z=\omega^*(X)$.

  Let $\Spec(p)$ be an fpqc left $\Spec(\pi)$-torsor over $\Spec(\kappa)$.
  We write $f^p\colon p\to \omega^{p*}\breve{\omega}^*\omega_*(\kappa)$
  for the unique morphism in $\cat{Comm}(\kos{K})$ that satisfies
  $\mathrm{eq}^p_{\breve{\omega}^*\omega_*(\kappa)}\circ f^p=
  \lambda_p\colon p\to \pi\otimes p$.
  One can check that $f^p$ exists,
  is an isomorphism,
  and is compatible with the left $\Spec(\pi)$-actions;
  the compatibility relation can be checked after postcomposing with
  $I_{\pi}\otimes \mathrm{eq}^p_{\breve{\omega}^*\omega_*(\kappa)}$
  using the relation
  $\mathrm{eq}^p_{\breve{\omega}^*\omega_*(\pi)}
  \circ (\tahar{\omega^{p*}\breve{\omega}^*\omega_*})_{\pi}
  =((\tahar{\varpi})_{\pi}\otimes I_p)
  \circ (I_{\pi}\otimes\mathrm{eq}^p_{\breve{\omega}^*\omega_*(\kappa)})
  \colon \pi\otimes \omega^{p*}\breve{\omega}^*\omega_*(\kappa)
  \to \varpi(\pi)\otimes p$.

  Moreover,
  $f^p\colon p\cong \omega^{p*}\breve{\omega}^*\omega_*(\kappa)$
  is natural in $\Spec(p)$.
  Indeed, let $\Spec(q)$ also be an fpqc left $\Spec(\pi)$-torsor over $\Spec(\kappa)$
  and $f\colon \Spec(q)\cong\Spec(p)$
  be an isomorphism in $(\Spec(\pi)\text{-}\cat{Tors})_{\kappa}$.
  We have
  $\vartheta^f_{\breve{\omega}^*\omega_*(\kappa)}\circ f^p
  =f^q\circ f
  \colon p\cong \omega^{q*}\breve{\omega}^*\omega_*(\kappa)$,
  as this can be checked after postcomposing with
  $\mathrm{eq}^q_{\breve{\omega}^*\omega_*(\kappa)}$.

  It remains to check the triangle identities.
  Let $\omega^{\pr}$ be a fiber point of $\kos{t}$
  and put $p\coloneqq\omega^{\pr*}\omega_*(\kappa)$.
  We have one of the triangle identities
  $(f^p)^{-1}\circ \vartheta^{\omega^{\pr}}_{\omega_*(\kappa)}
  =I_p\colon p\to p$,
  as this can be checked after postcomposing with $\lambda_p$.

  Let $\Spec(p)$ be an fpqc left $\Spec(\pi)$-torsor over $\Spec(\kappa)$.
  Put $\omega^{\pr}\coloneqq\breve{\omega}\omega^p$,
  $q\coloneqq \omega^{p*}\breve{\omega}^*\omega_*(\kappa)$,
  and $f\coloneqq f^p\colon \Spec(q)\cong \Spec(p)$.
  Then we have the other triangle identity 
  $(\breve{\omega}\vartheta^f)^{-1}\circ \vartheta^{\omega^{\pr}}
  =I_{\omega^{\pr}}$;
  for each $X\in \CT$ we have
  $(\vartheta^f_{\breve{\omega}^*(X)})^{-1}\circ \vartheta^{\omega^{\pr}}_X
  =I_{\omega^{\pr*}(X)}$,
  as this can be checked after postcomposing with
  $(I_{\omega^*(X)}\otimes \lambda_p)\circ \mathrm{eq}^p_{\breve{\omega}^*(X)}$
  using the relation
  $\mathrm{eq}^p_{\breve{\omega}^*\omega_*\omega^*(X)}
  \circ (\tahar{\omega^{p*}\breve{\omega}^*\omega_*})_{\omega^*(X)}
  =((\tahar{\varpi})_{\omega^*(X)}\otimes I_p)
  \circ (I_{\omega^*(X)}\otimes\mathrm{eq}^p_{\breve{\omega}^*\omega_*(\kappa)})
  \colon \omega^*(X)\otimes \omega^{p*}\breve{\omega}^*\omega_*(\kappa)
  \to \varpi\omega^*(X)\otimes p$.
  This completes the proof of Theorem~\ref{thm-pFibGroTors-perfcorres}.
\end{proof}

We end this section with a lemma
that exhibits the relation between
the gauge group and the twisted fiber point
associated to an fpqc left torsor over $\Spec(\kappa)$
under a Grothendieck object in $\kos{K}$.

\begin{lemma}
  \label{lem-pFibGroTors-pipomegap}
  Let $\Spec(\pi)$ be a Grothendieck object in $\kos{K}$
  and
  $\Spec(p)$ be an fpqc left $\Spec(\pi)$-torsor over $\Spec(\kappa)$.
  Recall from Proposition~\ref{prop-GroTors-gaugegroup}
  the gauge group $\Spec(\pi^p)$ of $\Spec(p)$,
  and from Proposition~\ref{prop-pFibGroTors-twistedfib}
  the $\Spec(p)$-twisted fiber point
  $\omega^p$ of $\kos{t}_{\pi}^{\vphantom{*}}$.
  \begin{enumerate}
    \item
    The presheaf of groups
    $\uAut_{\kos{K}}(\omega^p)$
    on $\cat{Aff}(\kos{K})$
    is represented by $\Spec(\pi^p)$.
    
    \item
    The fiber point $\omega^p$
    extends to an equivalence
    $\kos{t}_{\pi^p}^{\vphantom{*}}\simeq\kos{t}_{\pi}^{\vphantom{*}}$
    of Grothendieck $\kos{K}$-prekosmoi
    up to an invertible Grothendieck $\kos{K}$-transformation.
  \end{enumerate}
  \begin{equation*}
    \Hom_{\cat{Aff}(\kos{K})}(\slot,\Spec(\pi^p))
    \cong
    \uAut_{\kos{K}}(\omega^p)
    \qquad\quad
    \vcenter{\hbox{
      \begin{tikzcd}[row sep=5pt,column sep=6pt]
        \text{ }
        &\id_{\kos{K}}
        \ar[ddr,filled,"{\omega^p}"]
        \ar[ddl,filled,"{\kos{e}_{\pi^p}}"']
        &\text{ }
        \\%1
        \text{ }
        &\text{ }
        &\text{ }
        \ar[ll,2cell,shift left=0.5ex,pos=0.475]
        \\%2
        \kos{t}_{\pi^p}
        \ar[rr,filled,"\simeq"]
        &\text{ }
        &\kos{t}_{\pi}
      \end{tikzcd}
    }}
  \end{equation*}
\end{lemma}
\begin{proof}
  By Corollary~\ref{cor-Grocatoverpk-pFibsurjective},
  the twisted fiber point
  $\omega^p$ of $\kos{t}_{\pi}^{\vphantom{*}}$ is surjective.
  By Theorem~\ref{thm-Grocatoverpk-mainthm}
  for $\kos{t}=\kos{t}_{\pi}^{\vphantom{*}}$, $\omega=\omega^p$,
  we deduce that $\uAut_{\kos{K}}(\omega^p)$
  is represented by the Grothendieck object
  $\Spec(\omega^{p*}\omega^p_*(\kappa))$ in $\kos{K}$,
  and that $\omega^p$ extends to
  $\breve{\omega}^p\colon \kos{t}_{\omega^{p*}\omega^p_*(\kappa)}\simeq \kos{t}_{\pi}^{\vphantom{*}}$
  up to an invertible Grothendieck $\kos{K}$-transformation
  $\omega^p\cong \breve{\omega}^p\kos{e}_{\omega^{p*}\omega^p_*(\kappa)}$.
  The coproduct, counit, and antipode of $\omega^{p*}\omega^p_*(\kappa)$
  are given as in \eqref{eq-GroMainthm-Specpigroupstr}.
  Let
  $i_{\pi^p}\colon \pi^p\cong \omega^{p*}\omega^p_*(\kappa)$
  be the unique isomorphism in $\cat{Comm}(\kos{K})$ that satisfies
  $\mathrm{eq}^p_{\omega^p_*(\kappa)}\circ i_{\pi^p}
  =(\imath_p\otimes I_p)\circ \mathrm{eq}_p
  \colon \pi^p\to \kappa\otimes p\otimes p$.
  Then the isomorphism
  $(\tahar{\omega^{p*}\omega^p_*})_z\colon 
  z\otimes \omega^{p*}\omega^p_*(\kappa)
  \cong \omega^{p*}\omega^p_*(z)$
  in \eqref{eq-ActionbaseGro-taharupperlowerstar}
  satisfies
  $\mathrm{eq}^p_{\omega^p_*(z)}\circ (\tahar{\omega^{p*}\omega^p_*})_z
  =(I_z\otimes \mathrm{eq}_p)\circ (I_z\otimes i_{\pi^p}^{-1})
  \colon z\otimes \omega^{p*}\omega^p_*(\kappa)\to z\otimes p\otimes p$,
  and using this one can check that $i_{\pi^p}\colon \Spec(\omega^{p*}\omega^p_*(\kappa))\cong \Spec(\pi^p)$
  is an isomorphism of Grothendieck objects in $\kos{K}$.
  This proves (1) and (2).
\end{proof}

\section{Fundamental morphisms over a Grothendieck prekosmos}
\label{sec-fundmoroverGropk}

In this section,
we define the $(2,1)$-category
$\mathsf{pGroCat}(\kos{K})$
whose $0$-cells are Grothendieck $\kos{K}$-categories;
we begin by describing the $1$-cells.

\begin{definition}
  \label{def-fundmoroverGropk-fundmor}
  Let
  $\kos{T}\xmor{\kos{t}}\kos{K}$,
  $\kos{T}^{\pr}\xmor{\kos{t}^{\pr}}\kos{K}$
  be Grothendieck $\kos{K}$-categories.
  A \emph{fundamental $\kos{K}$-morphism},
  or a \emph{fundamental morphism over $\kos{K}$},
  from $\kos{t}^{\pr}$ to $\kos{t}$
  is a Grothendieck $\kos{K}$-morphism
  $\kos{f}\colon \kos{t}^{\pr}\mor\kos{t}$
  for which there exists a surjective fiber point
  $\omega^{\pr}\colon \vphantom{\id}\smash[t]{\id}_{\kos{K}}^{\vphantom{*}}\mor\kos{t}^{\pr}$ of $\kos{t}^{\pr}$
  such that the composite $\kos{f}\omega^{\pr}\colon \vphantom{\id}\smash[t]{\id}_{\kos{K}}^{\vphantom{*}}\mor\kos{t}$
  is a surjective fiber point of $\kos{t}$.
  \begin{equation*}
    \begin{tikzcd}[row sep=4pt,column sep=7pt]
      \text{ }
      &\id_{\kos{K}}
      \ar[ddr,filled,"{\kos{f}\omega^{\pr}}"]
      \ar[ddl,filled,"{\omega^{\pr}}"']
      &\text{ }
      \\%1
      \text{ }
      &\text{ }
      &\text{ }
      \\%2
      \kos{t}^{\pr}
      \ar[rr,filled,"{\kos{f}}"]
      &\text{ }
      &\kos{t}
    \end{tikzcd}
  \end{equation*}
\end{definition}

\begin{lemma}
  \label{lem-fundmoroverGropk-fundmor}
  Let 
  $f\colon \Spec(\pi^{\pr})\to\Spec(\pi)$
  be a morphism of Grothendieck objects in $\kos{K}$.
  We have a fundamental $\kos{K}$-morphism
  $\kos{f}\colon \kos{t}_{\pi^{\pr}}^{\vphantom{*}}\mor\kos{t}_{\pi}^{\vphantom{*}}$
  and an invertible Grothendieck $\kos{K}$-transformation
  $\kos{e}_{\pi}^{\vphantom{*}}\cong \kos{f}\kos{e}_{\pi^{\pr}}^{\vphantom{*}}$.
  \begin{equation*}
    \begin{tikzcd}[row sep=5pt,column sep=6pt]
      \text{ }
      &\id_{\kos{K}}
      \ar[ddr,filled,"{\kos{e}_{\pi}}"]
      \ar[ddl,filled,"{\kos{e}_{\pi^{\pr}}}"']
      &\text{ }
      \\%1
      \text{ }
      &\text{ }
      &\text{ }
      \ar[ll,2cell,shift left=0.5ex,pos=0.475]
      \\%2
      \kos{t}_{\pi^{\pr}}
      \ar[rr,filled,"{\kos{f}}"]
      &\text{ }
      &\kos{t}_{\pi}
    \end{tikzcd}
  \end{equation*}
\end{lemma}
\begin{proof}
  We write $X=(x,\rho_x)$
  for an object in $\cat{Rep}_{\kos{K}}(\pi)$
  and $Y=(y,\rho^{\pr}_y)$
  for an object in $\cat{Rep}_{\kos{K}}(\pi^{\pr})$.
  The left adjoint part $\kos{f}^*\colon \kosRep_{\kos{K}}(\pi)\to \kosRep_{\kos{K}}(\pi^{\pr})$
  is the strong monoidal functor
  that sends each $X$ to
  $(x,(I_x\otimes f)\circ \rho_x)\in \cat{Rep}_{\kos{K}}(\pi^{\pr})$
  and whose coherence isomorphisms are the identity morphisms.
  For each $Y$, we fix an equalizer
  $\kos{f}_*(Y)\hookrightarrow (y\otimes \pi,I_y\otimes \cp_{\pi})$
  in $\cat{Rep}_{\kos{K}}(\pi)$
  of the coreflexive pair
  $\rho^{\pr}_y\otimes I_{\pi}$,
  $(I_y\otimes f\otimes I_{\pi})\circ (I_y\otimes \cp_{\pi})
  \colon (y\otimes \pi,I_y\otimes \cp_{\pi})
  \to (y\otimes \pi^{\pr}\otimes \pi, I_{y\otimes \pi^{\pr}}\otimes \cp_{\pi})$,
  with common retraction
  $I_y\otimes e_{\pi^{\pr}}\otimes I_{\pi}$.
  Then $\kos{f}_*$ is the right adjoint part of $\kos{f}=(\kos{f}^*,\kos{f}_*,\eta^f,\epsilon^f)$.
  The component $\eta^f_X\colon X\to \kos{f}_*\kos{f}^*(X)$
  of $\eta^f$ at $X$ is the unique morphism in $\cat{Rep}_{\kos{K}}(\pi)$
  whose composite with 
  $\kos{f}_*\kos{f}^*(X)\hookrightarrow (x\otimes \pi,I_x\otimes \cp_{\pi})$
  is $\rho_x$.
  The component $\epsilon^f_Y\colon \kos{f}^*\kos{f}_*(Y)\to Y$
  of $\epsilon^f$ at $Y$,
  seen as a morphism in $\CK$,
  is the composite of
  $\kos{e}_{\pi}^*\kos{f}_*(Y)\hookrightarrow y\otimes \pi$
  followed by $I_y\otimes e_{\pi}\colon y\otimes \pi\to y$.
  The invertible Grothendieck transformation
  $\what{\kos{f}}\colon \kos{t}_{\pi^{\pr}}^{\vphantom{*}}\cong \kos{t}_{\pi}^{\vphantom{*}}\kos{f}$
  is the equality
  $\what{\kos{f}^*}\colon \kos{t}_{\pi^{\pr}}^*=\kos{f}^*\kos{t}_{\pi}^*$
  between the left adjoint parts,
  and the invertible Grothendieck $\kos{K}$-transformation
  $\kos{e}_{\pi}^{\vphantom{*}}\cong \kos{f}\kos{e}_{\pi^{\pr}}^{\vphantom{*}}$
  is the equality
  $\kos{e}_{\pi}^*=\kos{e}_{\pi^{\pr}}^*\kos{f}^*$
  between the left adjoint parts.
  Hence $\kos{f}\kos{e}_{\pi^{\pr}}^{\vphantom{*}}$
  is a surjective fiber point of $\kos{t}_{\pi}^{\vphantom{*}}$,
  and $\kos{f}\colon \kos{t}_{\pi^{\pr}}^{\vphantom{*}}\mor\kos{t}_{\pi}^{\vphantom{*}}$
  is a fundamental $\kos{K}$-morphism.
\end{proof}

\begin{proposition}
  \label{prop-fundmoroverGropk-inducedtor}
  Let $f\colon \Spec(\pi^{\pr})\to\Spec(\pi)$
  be a morphism of Grothendieck objects in $\kos{K}$
  and $\Spec(p^{\pr})$ be an fpqc left $\Spec(\pi^{\pr})$-torsor over $\Spec(\kappa)$.
  We have an fpqc left $\Spec(\pi)$-torsor
  $\Spec(p)$ over $\Spec(\kappa)$
  and an invertible Grothendieck $\kos{K}$-transformation
  $\omega^p\cong \kos{f}\omega^{p^{\pr}}$,
  where
  $\kos{f}\colon \kos{t}_{\pi^{\pr}}^{\vphantom{*}}\mor\kos{t}_{\pi}^{\vphantom{*}}$
  is the fundamental $\kos{K}$-morphism
  associated to $f$
  in Lemma~\ref{lem-fundmoroverGropk-fundmor},
  and 
  $\omega^p$, $\omega^{p^{\pr}}$
  are the twisted fiber points
  associated to $\Spec(p)$, $\Spec(p^{\pr})$, respectively,
  in Proposition~\ref{prop-pFibGroTors-twistedfib}.
  \begin{equation*}
    \begin{tikzcd}[row sep=5pt,column sep=6pt]
      \text{ }
      &\id_{\kos{K}}
      \ar[ddr,filled,"{\omega^p}"]
      \ar[ddl,filled,"{\omega^{p^{\pr}}}"']
      &\text{ }
      \\%1
      \text{ }
      &\text{ }
      &\text{ }
      \ar[ll,2cell,shift left=0.5ex,pos=0.475]
      \\%2
      \kos{t}_{\pi^{\pr}}
      \ar[rr,filled,"{\kos{f}}"]
      &\text{ }
      &\kos{t}_{\pi}
    \end{tikzcd}
  \end{equation*}
\end{proposition}

We refer to $\Spec(p)$ as the \emph{extension of structure group}
of $\Spec(p^{\pr})$ along $f$.

\begin{proof}[Proof of Proposition~\ref{prop-fundmoroverGropk-inducedtor}]
  Let 
  $\lambda^{\pr}_{p^{\pr}}\colon p^{\pr}\to \pi^{\pr}\otimes p^{\pr}$
  be the left $\Spec(\pi^{\pr})$-action
  on $\Spec(p^{\pr})$.
  We write
  $\cp^f\coloneqq (I_{\pi}\otimes f)\circ \cp_{\pi}
  \colon \pi\to \pi\otimes \pi^{\pr}$
  for the right $\Spec(\pi^{\pr})$-action on
  $\Spec(\pi)$ induced by $f$.
  We fix an equalizer
  $\mathrm{eq}\colon p\hookrightarrow \pi\otimes p^{\pr}$
  in $\cat{Comm}(\kos{K})$ of the coreflexive pair
  $I_{\pi}\otimes \lambda^{\pr}_{p^{\pr}}$,
  $\cp^f\otimes I_{p^{\pr}}\colon 
  \pi\otimes p^{\pr}\to \pi\otimes \pi^{\pr}\otimes p^{\pr}$,
  with common retraction
  $I_{\pi}\otimes e_{\pi^{\pr}}\otimes I_{p^{\pr}}$.
  Put $\tau^f\coloneqq (I_{\pi}\otimes \pc_{p^{\pr}})\circ (\mathrm{eq}\otimes I_{p^{\pr}})
  \colon p\otimes p^{\pr}\to \pi\otimes p^{\pr}$.
  In the diagram
  \begin{equation*}
    \vcenter{\hbox{
      \xymatrix@R=23pt@C=60pt{
        p\otimes p^{\pr}
        \ar@{^{(}->}[r]^-{\mathrm{eq}\otimes I_{p^{\pr}}}
        \ar@{.>}[d]_-{\exists!\,\tau^f}^-{\cong}
        &\pi\otimes p^{\pr}\otimes p^{\pr}
        \ar@<0.5ex>[r]^-{\cp^f\otimes I_{p^{\pr}\otimes p^{\pr}}}
        \ar@<-0.5ex>[r]_-{I_{\pi}\otimes \lambda^{\pr}_{p^{\pr}}\otimes I_{p^{\pr}}}
        \ar[d]_-{I_{\pi}\otimes \tau^{\pr}_{p^{\pr}}}^-{\cong}
        &\pi\otimes \pi^{\pr}\otimes p^{\pr}\otimes p^{\pr}
        \ar[d]^-{I_{\pi\otimes \pi^{\pr}}\otimes \tau^{\pr}_{p^{\pr}}}_-{\cong}
        \\%1
        \pi\otimes p^{\pr}
        \ar[r]^-{\cp^f\otimes I_{p^{\pr}}}
        &\pi\otimes \pi^{\pr}\otimes p^{\pr}
        \ar@<0.5ex>[r]^-{\cp^f\otimes I_{\pi^{\pr}\otimes p^{\pr}}}
        \ar@<-0.5ex>[r]_-{I_{\pi}\otimes \cp_{\pi^{\pr}}\otimes I_{p^{\pr}}}
        &\pi\otimes \pi^{\pr}\otimes \pi^{\pr}\otimes p^{\pr}
      }
    }}
  \end{equation*}
  we see that
  (i) the top horizontal morphisms form an equalizer diagram,
  since $\slot\otimes p^{\pr}$ preserves coreflexive equalizers;
  (ii) the bottom horizontal morphisms
  form a split equalizer diagram with respect to the splitting
  $I_{\pi}\otimes e_{\pi^{\pr}}\otimes I_{p^{\pr}}$,
  $I_{\pi\otimes \pi^{\pr}}\otimes e_{\pi^{\pr}}\otimes I_{p^{\pr}}$;
  (iii) the middle and right vertical isomorphisms
  commute with the horizontal parallel pairs;
  and
  (iv) the unique isomorphism between the two equalizers
  that makes the left square commute is given by
  $(I_{\pi}\otimes e_{\pi^{\pr}}\otimes I_{p^{\pr}})\circ (I_{\pi}\otimes \tau^{\pr}_{p^{\pr}})\circ (\mathrm{eq}\otimes I_{p^{\pr}})
  =\tau^f$.
  Since $\tau^f$ is an isomorphism in $\cat{Comm}(\kos{K})$,
  we deduce that $\Spec(p)$ is an fpqc covering of $\Spec(\kappa)$.

  Let $\lambda_p\colon p\to \pi\otimes p$
  be the unique morphism in $\cat{Comm}(\kos{K})$
  that satisfies
  $(I_{\pi}\otimes \mathrm{eq})\circ \lambda_p
  =(\cp_{\pi}\otimes I_{p^{\pr}})\circ \mathrm{eq}
  \colon p\to \pi\otimes \pi\otimes p^{\pr}$;
  one can check that $\lambda_p$
  exists and is a left $\Spec(\pi)$-action
  on $\Spec(p)$.
  Using $\mathrm{eq}\circ \pc_p=\pc_{\pi\otimes p^{\pr}}\circ (\mathrm{eq}\otimes\mathrm{eq})
  \colon p\otimes p\to \pi\otimes p^{\pr}$,
  we deduce that
  $(I_{\pi}\otimes \tau^f)\circ (\tau_p\otimes I_{p^{\pr}})
  =(\tau_{\pi}\otimes I_{p^{\pr}})
  \circ (I_{\pi}\otimes \tau^f)
  \circ (\sigma_{p,\pi}\otimes I_{p^{\pr}})
  \circ (I_p\otimes \tau^f)
  \circ (\sigma_{p,p}\otimes I_{p^{\pr}})
  \colon p\otimes p\otimes p^{\pr}\cong \pi\otimes \pi\otimes p^{\pr}$.
  Hence $\tau_p\otimes I_{p^{\pr}}$ is an isomorphism,
  and since $\slot\otimes p^{\pr}$ is conservative,
  $\tau_p$ is an isomorphism.
  Thus $\Spec(p)$ equipped with $\lambda_p$
  is an fpqc left $\Spec(\pi)$-torsor over $\Spec(\kappa)$.

  It remains to construct the invertible Grothendieck $\kos{K}$-transformation
  $\omega^p\cong \kos{f}\omega^{p^{\pr}}$.
  Put
  $g\coloneqq (e_{\pi}\otimes I_{p^{\pr}})\circ \mathrm{eq}
  \colon p\to p^{\pr}$.
  For each $X=(x,\rho_x)\in\cat{Rep}_{\kos{K}}(\pi)$,
  let $\vartheta_X\colon \omega^{p*}(X)\to \omega^{p^{\pr}*}\kos{f}^*(X)$
  be the unique morphism in $\CK$ that satisfies 
  $\mathrm{eq}^{p^{\pr}}_{\kos{f}^*(X)}\circ \vartheta_X
  =(I_x\otimes g)\circ \mathrm{eq}^p_X
  \colon \omega^{p*}(X)\to x\otimes p^{\pr}$.
  One can check that $\vartheta_X$ exists
  and that 
  $\xi^{p^{\pr}}_{\kos{f}^*(X)}\circ \kos{u}_{p^{\pr}}^*(\vartheta_X)
  \colon \kos{u}_{p^{\pr}}^*\omega^{p*}(X)\to \kos{u}_{p^{\pr}}^*\kos{e}_{\pi^{\pr}}^*\kos{f}^*(X)$
  is equal to
  $\kos{g}^*(\xi^p_X)\colon \kos{u}_{p^{\pr}}^*\omega^{p*}(X)
  \cong \kos{u}_{p^{\pr}}^*\kos{e}_{\pi}^*(X)$
  under the equality
  $\kos{u}_{p^{\pr}}^*\kos{e}_{\pi}^*
  =\kos{u}_{p^{\pr}}^*\kos{e}_{\pi^{\pr}}^*\kos{f}^*$.
  Thus $\kos{u}_{p^{\pr}}^*(\vartheta_X)$ is an isomorphism in $\CK_{p^{\pr}}^{\vphantom{*}}$,
  and since $\kos{u}_{p^{\pr}}^*$ is conservative,
  $\vartheta_X$ is an isomorphism in $\CK$.
  One can also check that the isomorphisms
  $\vartheta_X$ constitute
  an invertible Grothendieck $\kos{K}$-transformation
  $\vartheta\colon \omega^p\cong \kos{f}\omega^{p^{\pr}}$.
\end{proof}

\begin{proposition}
  \label{prop-fundmoroverGropk-realization}
  Let $\kos{T}\xmor{\kos{t}}\kos{K}$,
  $\kos{T}^{\pr}\xmor{\kos{t}^{\pr}}\kos{K}$
  be Grothendieck $\kos{K}$-categories
  and $\kos{f}\colon \kos{t}^{\pr}\mor\kos{t}$
  be a Grothendieck $\kos{K}$-morphism.
  The following are equivalent:
  \begin{enumerate}
    \item 
    $\kos{f}\colon \kos{t}^{\pr}\mor\kos{t}$
    is a fundamental $\kos{K}$-morphism.

    \item
    There exist a surjective fiber point
    $\omega$ of $\kos{t}$,
    a surjective fiber point $\omega^{\pr}$ of $\kos{t}^{\pr}$,
    and an invertible Grothendieck $\kos{K}$-transformation
    $\alpha\colon \omega\cong\kos{f}\omega^{\pr}$.
  \end{enumerate}
  In this case, 
  let $\omega$, $\omega^{\pr}$, and $\alpha$
  be as in \textup{(2)}.
  Put $\pi\coloneqq \omega^*\omega_*(\kappa)$
  and
  $\pi^{\pr}\coloneqq \omega^{\pr*}\omega^{\pr}_*(\kappa)$,
  and recall the equivalences
  $\breve{\omega}\colon \kos{t}_{\pi}^{\vphantom{*}}\simeq\kos{t}$,
  $\breve{\omega}^{\pr}\colon \kos{t}_{\pi^{\pr}}^{\vphantom{*}}\simeq\kos{t}^{\pr}$
  in Theorem~\ref{thm-Grocatoverpk-mainthm}.
  We have a morphism
  $\breve{f}\colon \Spec(\pi^{\pr})\to\Spec(\pi)$
  of Grothendieck objects in $\kos{K}$
  and an invertible Grothendieck $\kos{K}$-transformation
  $\breve{\alpha}\colon \breve{\omega}\breve{\kos{f}}\cong \kos{f}\breve{\omega}^{\pr}$
  that satisfies the following equality of pasting diagrams,
  where $\breve{\kos{f}}\colon \kos{t}_{\pi^{\pr}}^{\vphantom{*}}\mor\kos{t}_{\pi}^{\vphantom{*}}$
  is the fundamental $\kos{K}$-morphism
  associated to $\breve{f}$
  in Lemma~\ref{lem-fundmoroverGropk-fundmor}.
  \begin{equation*}
    \begin{tikzcd}[row sep=5pt,column sep=10pt]
      \text{ }
      &\text{ }
      &\id_{\kos{K}}
      \ar[ddll,bend right=1.8pc,filled,"{\kos{e}_{\pi^{\pr}}}"']
      \ar[dd,filled,"{\kos{e}_{\pi}}"]
      \ar[ddddrr,bend left=2.5pc,filled,"{\omega}"]
      &\text{ }
      &\text{ }
      &\text{ }
      &\text{ }
      &\text{ }
      &\id_{\kos{K}}
      \ar[ddll,bend right=1.8pc,filled,"{\kos{e}_{\pi^{\pr}}}"']
      \ar[dddd,filled,"{\omega^{\pr}}" description]
      \ar[ddddrr,bend left=2.5pc,filled,"{\omega}"]
      &\text{ }
      &\text{ }
      \\%1
      \text{ }
      &\text{ }
      &\text{ }
      \ar[ll,2cell,shift left=0.5ex,pos=0.25]
      &\text{ }
      &\text{ }
      &\text{ }
      &\text{ }
      &\text{ }
      &\text{ }
      &\text{ }
      &\text{ }
      \\%2
      \kos{t}_{\pi^{\pr}}
      \ar[rr,filled,"\breve{\kos{f}}"]
      \ar[ddrr,filled,"\breve{\omega}^{\pr}"',"\simeq"]
      &\text{ }
      &\kos{t}_{\pi}
      \ar[ddrr,filled,"\breve{\omega}","\simeq"']
      &\text{ }
      &\text{ }
      \ar[ll,2cell,shift right=0.5ex,pos=0.6]
      &=
      &\kos{t}_{\pi^{\pr}}
      \ar[ddrr,filled,"\breve{\omega}^{\pr}"',"\simeq"]
      &\text{ }
      &\text{ }
      \ar[ll,2cell,shift right=1ex,pos=0.4]
      &\text{ }
      &\text{ }
      \ar[ll,2cell,shift left=1.5ex,pos=0.55,"{\,\alpha}"']
      \\%3
      \text{ }
      &\text{ }
      &\text{ }
      &\text{ }
      \ar[ll,2cell,shift left=1ex,"{\,\breve{\alpha}}"']
      &\text{ }
      &\text{ }
      &\text{ }
      &\text{ }
      &\text{ }
      &\text{ }
      &\text{ }
      \\%4
      \text{ }
      &\text{ }
      &\kos{t}^{\pr}
      \ar[rr,filled,"\kos{f}"']
      &\text{ }
      &\kos{t}
      &\text{ }
      &\text{ }
      &\text{ }
      &\kos{t}^{\pr}
      \ar[rr,filled,"\kos{f}"']
      &\text{ }
      &\kos{t}
    \end{tikzcd}
  \end{equation*}
\end{proposition}
\begin{proof}
  It is immediate that (1) and (2) are equivalent.
  Let $\omega$, $\omega^{\pr}$, and $\alpha$ be as in (2).
  We write $\alpha_*\colon \kos{f}^*\omega_*\Rightarrow \omega^{\pr}_*$
  for the mate of $\alpha^{-1}\colon \omega^{\pr*}\kos{f}^*\cong \omega^*$
  with respect to the adjunctions
  $\omega$ and $\omega^{\pr}$;
  it is the unique lax $\kos{K}$-tensor natural transformation
  that satisfies the following relations.
  \begin{equation}\label{eq-fundmoroverGropk-matealpha}
    \vcenter{\hbox{
      \xymatrix@R=15pt@C=30pt{
        \kos{f}^*
        \ar@2{->}[r]^-{\eta^{\pr}\kos{f}^*}
        \ar@2{->}[d]_-{\kos{f}^*\eta}
        &\omega^{\pr}_*\omega^{\pr*}\kos{f}^*
        \ar@2{->}[d]^-{\omega^{\pr}_*\alpha^{-1}}_-{\cong}
        &\text{ }
        &\omega^{\pr*}\kos{f}^*\omega_*
        \ar@2{->}[r]^-{\omega^{\pr*}\alpha_*}
        \ar@2{->}[d]_-{\alpha^{-1}\omega_*}^-{\cong}
        &\omega^{\pr*}\omega^{\pr}_*
        \ar@2{->}[d]^-{\epsilon^{\pr}}
        \\%1
        \kos{f}^*\omega_*\omega^*
        \ar@2{->}[r]^-{\alpha_*\omega^*}
        &\omega^{\pr}_*\omega^*
        &\text{ }
        &\omega^*\omega_*
        \ar@2{->}[r]^-{\epsilon}
        &\id_{\kos{K}}
      }
    }}
  \end{equation}
  Put $\alpha^{\sharp}\coloneqq (\omega^{\pr*}\alpha_*)\circ (\alpha\omega_*)
  \colon \omega^*\omega_*\Rightarrow\omega^{\pr*}\omega^{\pr}_*$
  and define 
  $\breve{f}\colon \Spec(\pi^{\pr})\to \Spec(\pi)$ as 
  $\alpha^{\sharp}_{\kappa}\colon \omega^*\omega_*(\kappa)\to \omega^{\pr*}\omega^{\pr}_*(\kappa)$.
  By substituting $\vartheta=\alpha^{\sharp}$
  in \eqref{eq-ActionbaseGro-laxKtnattahar},
  we obtain the following commutative diagram for every $z\in \CK$.
  \begin{equation}
    \label{eq-fundmoroverGropk-alphasharpequiv}
    \vcenter{\hbox{
      \xymatrix@R=15pt@C=40pt{
        z\otimes \omega^*\omega_*(\kappa)
        \ar[r]^-{(\tahar{\omega^*\omega_*})_z}_-{\cong}
        \ar[d]_-{I_z\otimes \breve{f}}
        &\omega^*\omega_*(z)
        \ar[d]^-{\alpha^{\sharp}_z}
        \\%1
        z\otimes \omega^{\pr*}\omega^{\pr}_*(\kappa)
        \ar[r]^-{(\tahar{\omega^{\pr*}\omega^{\pr}_*})_z}_-{\cong}
        &\omega^{\pr*}\omega^{\pr}_*(z)
      }
    }}
  \end{equation}
  One can check that
  $\alpha^{\sharp}\colon \omega^*\omega_*\Rightarrow\omega^{\pr*}\omega^{\pr}_*$
  is a morphism of comonads,
  using \eqref{eq-fundmoroverGropk-matealpha}.
  From this,
  along with \eqref{eq-fundmoroverGropk-alphasharpequiv} for $z=\pi$,
  we deduce that $\breve{f}\colon \Spec(\pi^{\pr})\to\Spec(\pi)$
  is a morphism of Grothendieck objects in $\kos{K}$.
  It remains to show the existence of $\breve{\alpha}$.
  It suffices to show that for each $X\in \CT$,
  the isomorphism
  $\alpha_X\colon \omega^*(X)\cong \omega^{\pr*}\kos{f}^*(X)$
  in $\CK$ lifts to an isomorphism
  $\breve{\alpha}_X\colon \breve{\kos{f}}^*\breve{\omega}^*(X)\cong \breve{\omega}^{\pr*}\kos{f}^*(X)$
  in $\cat{Rep}_{\kos{K}}(\pi^{\pr})$.
  One can check this using
  \eqref{eq-fundmoroverGropk-alphasharpequiv}
  for $z=\omega^*(X)$ and
  \eqref{eq-fundmoroverGropk-matealpha}.
\end{proof}

\begin{lemma}
  \label{lem-fundmoroverGropk-sendingptpt}
  Let $\kos{t}$, $\kos{t}^{\pr}$
  be Grothendieck $\kos{K}$-categories
  and $\kos{f}\colon \kos{t}^{\pr}\mor\kos{t}$
  be a fundamental $\kos{K}$-morphism.
  The composite of $\kos{f}$
  with any fiber point of $\kos{t}^{\pr}$
  is a fiber point of $\kos{t}$.
\end{lemma}
\begin{proof}
  Let $\omega$, $\omega^{\pr}$, $\alpha$,
  $\pi$, $\pi^{\pr}$,
  $\breve{\omega}$,
  $\breve{\omega}^{\pr}$,
  $\breve{f}$, $\breve{\kos{f}}$,
  and $\breve{\alpha}$
  be as in Proposition~\ref{prop-fundmoroverGropk-realization},
  and $\omega^{\ppr}\colon \vphantom{\id}\smash[t]{\id}_{\kos{K}}^{\vphantom{*}}\mor\kos{t}^{\pr}$
  be a fiber point of $\kos{t}^{\pr}$.
  Recall from the proof of Theorem~\ref{thm-pFibGroTors-perfcorres} that
  $\Spec(p^{\pr})$ with $p^{\pr}\coloneqq \omega^{\ppr*}\omega^{\pr}_*(\kappa)$
  is an fpqc left $\Spec(\pi^{\pr})$-torsor over $\Spec(\kappa)$,
  and that $\omega^{\ppr}\cong \breve{\omega}^{\pr}\omega^{p^{\pr}}$.
  Let $\Spec(p)$
  be the extension of structure group
  of $\Spec(p^{\pr})$ along $\breve{f}$
  with 
  $\omega^p\cong \breve{\kos{f}}\omega^{p^{\pr}}$
  as in Proposition~\ref{prop-fundmoroverGropk-inducedtor}.
  From these isomorphisms, along with 
  $\breve{\alpha}\colon \breve{\omega}\breve{\kos{f}}\cong \kos{f}\breve{\omega}^{\pr}$,
  we deduce that
  $\kos{f}\omega^{\ppr}
  \cong \kos{f}\breve{\omega}^{\pr}\omega^{p^{\pr}}
  \cong \breve{\omega}\breve{\kos{f}}\omega^{p^{\pr}}
  \cong \breve{\omega}\omega^p$.
  Thus $\kos{f}\omega^{\ppr}$ is a fiber point of $\kos{t}$.
\end{proof}

\begin{lemma}
  \label{lem-fundmoroverGropk-closedundercomp}
  Fundamental $\kos{K}$-morphisms
  between Grothendieck $\kos{K}$-categories
  are closed under composition.
\end{lemma}
\begin{proof}
  Let
  $\kos{f}\colon \kos{t}^{\pr}\mor\kos{t}$,
  $\kos{f}^{\pr}\colon \kos{t}^{\ppr}\mor\kos{t}^{\pr}$
  be fundamental $\kos{K}$-morphisms
  between Grothendieck $\kos{K}$-categories
  $\kos{t}$, $\kos{t}^{\pr}$, $\kos{t}^{\ppr}$.
  Choose a surjective fiber point $\omega^{\ppr}$ of $\kos{t}^{\ppr}$
  such that $\kos{f}^{\pr}\omega^{\ppr}$ is a surjective fiber point of $\kos{t}^{\pr}$.
  By Lemma~\ref{lem-fundmoroverGropk-sendingptpt},
  the composite $\kos{f}\kos{f}^{\pr}\omega^{\ppr}$
  is a fiber point of $\kos{t}$,
  and it is surjective by Corollary~\ref{cor-Grocatoverpk-pFibsurjective}.
  Thus $\kos{f}\kos{f}^{\pr}\colon \kos{t}^{\ppr}\mor\kos{t}$
  is a fundamental $\kos{K}$-morphism.
\end{proof}

\begin{lemma}
  \label{lem-fundmoroverGropk-HomisIsom}
  Let $\kos{t}$, $\kos{t}^{\pr}$
  be Grothendieck $\kos{K}$-categories
  and
  $\kos{f}$, $\kos{g}\colon \kos{t}^{\pr}\mor\kos{t}$
  be fundamental $\kos{K}$-morphisms.
  Every Grothendieck $\kos{K}$-transformation
  from $\kos{f}$ to $\kos{g}$ is invertible.
\end{lemma}
\begin{proof}
  Let $\vartheta\colon \kos{f}\Rightarrow\kos{g}$
  be a Grothendieck $\kos{K}$-transformation
  and choose a surjective fiber point
  $\omega^{\pr}$ of $\kos{t}^{\pr}$.
  By Lemma~\ref{lem-fundmoroverGropk-sendingptpt},
  both $\kos{f}\omega^{\pr}$ and $\kos{g}\omega^{\pr}$
  are fiber points of $\kos{t}$.
  Thus 
  $\vartheta\omega^{\pr}\colon \kos{f}\omega^{\pr}\Rightarrow \kos{g}\omega^{\pr}$
  is invertible by Corollary~\ref{cor-Grocatoverpk-pFibgroupoid},
  and since $\omega^{\pr}$ is surjective, $\vartheta$ is invertible.
\end{proof}

Note that $\vphantom{\id}\smash[t]{\id}_{\kos{K}}^{\vphantom{*}}$ is a Grothendieck $\kos{K}$-category,
as the identity Grothendieck $\kos{K}$-morphism of $\vphantom{\id}\smash[t]{\id}_{\kos{K}}^{\vphantom{*}}$
is a surjective fiber point of $\vphantom{\id}\smash[t]{\id}_{\kos{K}}^{\vphantom{*}}$.
Conversely, every fiber point $\omega$ of $\vphantom{\id}\smash[t]{\id}_{\kos{K}}^{\vphantom{*}}$
is isomorphic to the identity Grothendieck $\kos{K}$-morphism
of $\vphantom{\id}\smash[t]{\id}_{\kos{K}}^{\vphantom{*}}$;
we have $\what{\omega}\colon \vphantom{\id}\smash[t]{\id}_{\kos{K}}^{\vphantom{*}}\cong\omega$.

\begin{lemma}
  \label{lem-fundmoroverGropk-pointismor}
  Let $\kos{t}$ be a Grothendieck $\kos{K}$-category.
  Fundamental $\kos{K}$-morphisms
  from $\vphantom{\id}\smash[t]{\id}_{\kos{K}}^{\vphantom{*}}$ to $\kos{t}$
  are precisely the surjective fiber points of $\kos{t}$.
\end{lemma}
\begin{proof}
  Let $\omega\colon \vphantom{\id}\smash[t]{\id}_{\kos{K}}^{\vphantom{*}}\mor\kos{t}$ be a Grothendieck $\kos{K}$-morphism.
  It is immediate that if $\omega$ is a surjective fiber point of $\kos{t}$,
  then $\omega$ is a fundamental $\kos{K}$-morphism.
  Conversely, assume that $\omega$ is a fundamental $\kos{K}$-morphism.
  Then $\omega$ is isomorphic to a surjective fiber point of $\kos{t}$,
  and is thus a surjective fiber point of $\kos{t}$.
\end{proof}

\begin{definition}
  \label{def-fundmoroverGropk-pGroCatK}
  We define the $(2,1)$-category
  $\mathsf{pGroCat}(\kos{K})$
  as follows.
  \begin{itemize}
    \item 
    A $0$-cell in $\mathsf{pGroCat}(\kos{K})$
    is a Grothendieck $\kos{K}$-category.
    
    \item
    A $1$-cell in $\mathsf{pGroCat}(\kos{K})$
    is a fundamental $\kos{K}$-morphism.

    \item
    A $2$-cell in $\mathsf{pGroCat}(\kos{K})$
    is an invertible Grothendieck $\kos{K}$-transformation.

    \item
    The rest of the $2$-categorical structure of
    $\mathsf{pGroCat}(\kos{K})$
    is given by that of $\mathtt{pGroKos}/\kos{K}$;
    the composition of $1$-cells is well defined by
    Lemma~\ref{lem-fundmoroverGropk-closedundercomp}.
  \end{itemize}
\end{definition}

\section{A duality over a Grothendieck prekosmos}
\label{sec-dualityoverGropk}

In this section, we establish a duality
between Grothendieck objects in $\kos{K}$
and Grothendieck $\kos{K}$-categories
with a choice of a surjective fiber point;
the latter are defined as follows.

\begin{definition}
  \label{def-dualityoverGropk-ptGrocat}
  A \emph{pointed Grothendieck $\kos{K}$-category},
  also called a \emph{pointed Grothendieck category over $\kos{K}$},
  is a pair $(\kos{t},\omega)$ 
  consisting of a Grothendieck $\kos{K}$-category $\kos{t}$
  and a surjective fiber point $\omega$ of $\kos{t}$.
\end{definition}

Recall from Section~\ref{sec-Groobjinpk}
the $(2,1)$-category $\mathsf{pGroObj}(\kos{K})$
of Grothendieck objects in $\kos{K}$,
defined in \eqref{eq-Groobjinpk-pGroObjK21cat}.
The goal of this section is to
define the $(2,1)$-category
$\mathsf{pGroCat}_{\mathsf{pt}}(\kos{K})$
whose $0$-cells are pointed Grothendieck $\kos{K}$-categories,
and prove the following theorem.

\begin{theorem}
  \label{thm-dualityoverGropk-equiv21cat}
  We have mutually quasi-inverse equivalences
  \begin{equation*}
    \xymatrix@C=35pt{
      \mathsf{pGroObj}(\kos{K})
      \ar@<1ex>[r]^-{\mathsf{Rep}_{\kos{K}}}
      \ar@{}[r]|-{\simeq}
      &\mathsf{pGroCat}_{\mathsf{pt}}(\kos{K})
      \ar@<1ex>[l]^-{\mathsf{Aut}_{\kos{K}}}
    }
  \end{equation*}
  of $(2,1)$-categories.
\end{theorem}

\begin{definition}
  \label{def-dualityoverGropk-pGroCatptK}
  We define the $(2,1)$-category
  $\mathsf{pGroCat}_{\mathsf{pt}}(\kos{K})$
  as follows.
  \begin{itemize}
    \item 
    A $0$-cell in $\mathsf{pGroCat}_{\mathsf{pt}}(\kos{K})$
    is a pointed Grothendieck $\kos{K}$-category.
 
    \item
    A $1$-cell in $\mathsf{pGroCat}_{\mathsf{pt}}(\kos{K})$
    is a \emph{pointed fundamental $\kos{K}$-morphism}.
    Let $(\kos{t},\omega)$, $(\kos{t}^{\pr},\omega^{\pr})$
    be pointed Grothendieck $\kos{K}$-categories.
    A pointed fundamental $\kos{K}$-morphism
    $(\kos{f},\alpha)\colon (\kos{t}^{\pr},\omega^{\pr})\mor (\kos{t},\omega)$,
    denoted by an arrow with a filled tip,
    is a pair consisting of a Grothendieck $\kos{K}$-morphism
    $\kos{f}\colon \kos{t}^{\pr}\mor\kos{t}$
    and an invertible Grothendieck $\kos{K}$-transformation
    $\alpha\colon \omega\cong \kos{f}\omega^{\pr}$.
    Such $\kos{f}$ is necessarily a fundamental $\kos{K}$-morphism
    by Proposition~\ref{prop-fundmoroverGropk-realization}.
    We illustrate
    $(\kos{f},\alpha)\colon (\kos{t}^{\pr},\omega^{\pr})\mor (\kos{t},\omega)$
    as follows.
    \begin{equation*}
      \begin{tikzcd}[row sep=5pt,column sep=7pt]
        \text{ }
        &\id_{\kos{K}}
        \ar[ddr,filled,"{\omega}"]
        \ar[ddl,filled,"{\omega^{\pr}}"']
        &\text{ }
        \\%1
        \text{ }
        &\text{ }
        &\text{ }
        \ar[ll,2cell,shift left=0.8ex,pos=0.475,"{\alpha}"']
        \\%2
        \kos{t}^{\pr}
        \ar[rr,filled,"{\kos{f}}"]
        &\text{ }
        &\kos{t}
      \end{tikzcd}
    \end{equation*}

    \item
    A $2$-cell 
    $\vartheta\colon (\kos{f},\alpha)\cong(\kos{g},\beta)$
    between $1$-cells
    $(\kos{f},\alpha)$, $(\kos{g},\beta)\colon (\kos{t}^{\pr},\omega^{\pr})\mor(\kos{t},\omega)$
    in $\mathsf{pGroCat}_{\mathsf{pt}}(\kos{K})$
    is an invertible Grothendieck $\kos{K}$-transformation
    $\vartheta\colon \kos{f}\cong \kos{g}$
    between the underlying fundamental $\kos{K}$-morphisms.

    \item
    The composite of $1$-cells
    $(\kos{f},\alpha)\colon (\kos{t}^{\pr},\omega^{\pr})\mor(\kos{t},\omega)$,
    $(\kos{f}^{\pr},\alpha^{\pr})\colon (\kos{t}^{\ppr},\omega^{\ppr})\mor(\kos{t}^{\pr},\omega^{\pr})$
    is the pair $(\kos{f}\kos{f}^{\pr},(\kos{f}\alpha^{\pr})\circ \alpha)
    \colon (\kos{t}^{\ppr},\omega^{\ppr})\mor(\kos{t},\omega)$.
    The identity $1$-cell of a $0$-cell 
    $(\kos{t},\omega)$
    is the pair consisting of the identity Grothendieck $\kos{K}$-morphism
    of $\kos{t}$
    and the identity Grothendieck $\kos{K}$-transformation of $\omega\colon \vphantom{\id}\smash[t]{\id}_{\kos{K}}^{\vphantom{*}}\mor\kos{t}$.
    Horizontal and vertical compositions of $2$-cells 
    are given by those of invertible Grothendieck $\kos{K}$-transformations
    between the underlying fundamental $\kos{K}$-morphisms.
    Identity $2$-cells are given by the identity Grothendieck $\kos{K}$-transformations
    of the underlying fundamental $\kos{K}$-morphisms.
  \end{itemize}
\end{definition}

\begin{remark}
  \label{rem-dualityoverGropk-notcoslice}
  The $0$-cells and the $1$-cells
  of $\mathsf{pGroCat}_{\mathsf{pt}}(\kos{K})$
  are those of the coslice $(2,1)$-category
  of $\mathsf{pGroCat}(\kos{K})$ under $\vphantom{\id}\smash[t]{\id}_{\kos{K}}^{\vphantom{*}}$,
  by Lemma~\ref{lem-fundmoroverGropk-pointismor},
  whereas the $2$-cells of $\mathsf{pGroCat}_{\mathsf{pt}}(\kos{K})$
  are inherited from $\mathsf{pGroCat}(\kos{K})$
  without any compatibility conditions.
\end{remark}

\begin{proposition}
  \label{prop-dualityoverGropk-RepKpseudofunct}
  We have a pseudofunctor
  \begin{equation*}
    \mathsf{Rep}_{\kos{K}}\colon \mathsf{pGroObj}(\kos{K})\to \mathsf{pGroCat}_{\mathsf{pt}}(\kos{K})
  \end{equation*}
  that sends each Grothendieck object $\Spec(\pi)$ in $\kos{K}$
  to the pointed Grothendieck $\kos{K}$-category
  $(\kos{t}_{\pi}^{\vphantom{*}},\kos{e}_{\pi}^{\vphantom{*}})$.
\end{proposition}
\begin{proof}
  We describe how $\mathsf{Rep}_{\kos{K}}$ sends $1$-cells and $2$-cells.
  \begin{itemize}
    \item
    A $1$-cell $f\colon \Spec(\pi^{\pr})\to\Spec(\pi)$
    in $\mathsf{pGroObj}(\kos{K})$
    is sent to the pair consisting of
    the fundamental $\kos{K}$-morphism
    $\kos{f}\colon \kos{t}_{\pi^{\pr}}^{\vphantom{*}}\mor \kos{t}_{\pi}^{\vphantom{*}}$
    and the invertible Grothendieck $\kos{K}$-transformation
    $\kos{e}_{\pi}^{\vphantom{*}}\cong \kos{f}\kos{e}_{\pi^{\pr}}^{\vphantom{*}}$
    in Lemma~\ref{lem-fundmoroverGropk-fundmor}.

    \item
    Let $f$, $g\colon \Spec(\pi^{\pr})\to\Spec(\pi)$ be $1$-cells 
    in $\mathsf{pGroObj}(\kos{K})$
    and write $\kos{f}$, $\kos{g}\colon \kos{t}_{\pi^{\pr}}^{\vphantom{*}}\mor\kos{t}_{\pi}^{\vphantom{*}}$
    for the associated fundamental $\kos{K}$-morphisms.
    A $2$-cell $\theta\colon f\cong g$
    in $\mathsf{pGroObj}(\kos{K})$
    is sent to the $2$-cell
    $\vartheta\colon \kos{f}\cong \kos{g}$
    in $\mathsf{pGroCat}_{\mathsf{pt}}(\kos{K})$,
    defined as follows.
    The component $\vartheta_X\colon \kos{f}^*(X)\cong \kos{g}^*(X)$
    of $\vartheta$ at $X=(x,\rho_x)\in\cat{Rep}_{\kos{K}}(\pi)$
    is the isomorphism in $\cat{Rep}_{\kos{K}}(\pi^{\pr})$
    whose underlying isomorphism in $\CK$
    is $(I_x\otimes \theta)\circ \rho_x\colon x\cong x$.
    It is straightforward to check that
    $\vartheta$ is indeed an invertible Grothendieck $\kos{K}$-transformation
    from $\kos{f}$ to $\kos{g}$.
  \end{itemize}
  It is straightforward to check that 
  $\mathsf{Rep}_{\kos{K}}$ preserves vertical composition of $2$-cells and identity $2$-cells;
  thus $\mathsf{Rep}_{\kos{K}}$ defines a functor
  on each hom-groupoid.
  The constraints of the pseudofunctor $\mathsf{Rep}_{\kos{K}}$
  are defined as follows.
  \begin{itemize}
    \item
    Let $f$, $f^{\pr}$
    be a composable pair of $1$-cells in $\mathsf{pGroObj}(\kos{K})$,
    with composite $f^{\ppr}\coloneqq ff^{\pr}$.
    We write
    $\kos{f}$, $\kos{f}^{\pr}$, and $\kos{f}^{\ppr}$
    for the fundamental $\kos{K}$-morphisms
    associated to $f$, $f^{\pr}$, and $f^{\ppr}$, respectively.
    The component of the functorial constraint of $\mathsf{Rep}_{\kos{K}}$
    at $f$, $f^{\pr}$
    is the invertible Grothendieck $\kos{K}$-transformation
    $\kos{f}\kos{f}^{\pr}\cong \kos{f}^{\ppr}$
    given by the equality
    $\kos{f}^{\pr*}\kos{f}^*=\kos{f}^{\ppr*}$
    between the left adjoint parts.

    \item
    Let $\Spec(\pi)$ be a $0$-cell in $\mathsf{pGroObj}(\kos{K})$.
    We write 
    $\kos{I}_{\pi}^{\vphantom{*}}\colon \kos{t}_{\pi}^{\vphantom{*}}\mor\kos{t}_{\pi}^{\vphantom{*}}$
    for the fundamental $\kos{K}$-morphism
    associated to the identity morphism of $\Spec(\pi)$
    in Lemma~\ref{lem-fundmoroverGropk-fundmor}.
    The component of the unity constraint of $\mathsf{Rep}_{\kos{K}}$
    at $\Spec(\pi)$
    is the invertible Grothendieck $\kos{K}$-transformation
    from the identity Grothendieck $\kos{K}$-morphism of $\kos{t}_{\pi}^{\vphantom{*}}$
    to $\kos{I}_{\pi}^{\vphantom{*}}$,
    given by the equality between the left adjoint parts.
  \end{itemize}
  One can check that $\mathsf{Rep}_{\kos{K}}$ preserves
  horizontal composition of $2$-cells
  up to the functorial constraint,
  or equivalently,
  that the functorial constraint
  $\kos{f}\kos{f}^{\pr}\cong \kos{f}^{\ppr}$
  is natural in $f$, $f^{\pr}$.
  It is immediate that the functorial and unity constraints of $\mathsf{Rep}_{\kos{K}}$
  satisfy the associativity and left/right unity axioms,
  as the axioms are just equalities between the left adjoint parts.
  This completes the construction of
  $\mathsf{Rep}_{\kos{K}}\colon \mathsf{pGroObj}(\kos{K})\to \mathsf{pGroCat}_{\mathsf{pt}}(\kos{K})$.
\end{proof}

\begin{proposition}
  \label{prop-dualityoverGropk-AutK2funct}
  We have a $2$-functor
  \begin{equation*}
    \mathsf{Aut}_{\kos{K}}\colon \mathsf{pGroCat}_{\mathsf{pt}}(\kos{K})\to \mathsf{pGroObj}(\kos{K})
  \end{equation*}
  that sends each pointed Grothendieck $\kos{K}$-category
  $(\kos{t},\omega)$
  to the Grothendieck object $\Spec(\omega^*\omega_*(\kappa))$ in $\kos{K}$
  that represents $\uAut_{\kos{K}}(\omega)$ in Theorem~\ref{thm-Grocatoverpk-mainthm}.
\end{proposition}
\begin{proof}
  We describe how $\mathsf{Aut}_{\kos{K}}$ sends $1$-cells and $2$-cells.
  \begin{itemize}
    \item
    A $1$-cell $(\kos{f},\alpha)\colon (\kos{t}^{\pr},\omega^{\pr})\mor(\kos{t},\omega)$
    in $\mathsf{pGroCat}_{\mathsf{pt}}(\kos{K})$ is sent to
    the $1$-cell
    $\breve{f}\colon \Spec(\omega^{\pr*}\omega^{\pr}_*(\kappa))\to \Spec(\omega^*\omega_*(\kappa))$
    in $\mathsf{pGroObj}(\kos{K})$
    described in the proof of Proposition~\ref{prop-fundmoroverGropk-realization}.

    \item
    A $2$-cell $\vartheta\colon (\kos{f},\alpha)\cong (\kos{g},\beta)$
    between $1$-cells
    $(\kos{f},\alpha)$, $(\kos{g},\beta)\colon (\kos{t}^{\pr},\omega^{\pr})\mor (\kos{t},\omega)$
    in $\mathsf{pGroCat}_{\mathsf{pt}}(\kos{K})$
    is sent to the $2$-cell
    $\breve{\theta}\colon \breve{f}\cong \breve{g}$
    in $\mathsf{pGroObj}(\kos{K})$,
    defined as follows.
    Put $\pi\coloneqq \omega^*\omega_*(\kappa)$
    and 
    $\vartheta^{\sharp}\coloneqq 
    \epsilon\circ(\beta^{-1}\omega_*)\circ(\omega^{\pr*}\vartheta\omega_*)\circ(\alpha\omega_*)
    \colon \omega^*\omega_*\Rightarrow\vphantom{\id}\smash[t]{\id}_{\kos{K}}^{\vphantom{*}}$.
    We define
    $\breve{\theta}\colon \Spec(\kappa)\to \Spec(\pi)$
    as $\vartheta^{\sharp}_{\kappa}\colon \pi\to \kappa$.
    By \eqref{eq-ActionbaseGro-laxKtnattahar}
    for the lax $\kos{K}$-tensor natural transformation $\vartheta^{\sharp}$,
    we deduce that the following diagram commutes
    for every $z\in \CK$.
    \begin{equation}
      \label{eq-dualityoverGropk-varthetasharpequiv}
      \vcenter{\hbox{
        \xymatrix@R=12pt@C=35pt{
          z\otimes \omega^*\omega_*(\kappa)
          \ar[r]^-{(\tahar{\omega^*\omega_*})_z}_-{\cong}
          \ar[d]_-{I_z\otimes \breve{\theta}}
          &\omega^*\omega_*(z)
          \ar[d]^-{\vartheta^{\sharp}_z}
          \\%1
          z\otimes \kappa
          \ar[r]^-{\jmath_z^{-1}}_-{\cong}
          &z
        }
      }}
    \end{equation}
    Using \eqref{eq-fundmoroverGropk-matealpha},
    one can check that 
    $\vartheta^{\sharp}=\epsilon^{\pr}\circ (\omega^{\pr*}\beta_*)\circ (\omega^{\pr*}\vartheta\omega_*) \circ(\alpha\omega_*)
    \colon \omega^*\omega_*\Rightarrow\vphantom{\id}\smash[t]{\id}_{\kos{K}}^{\vphantom{*}}$
    and that $(\alpha^{\sharp}\vartheta^{\sharp})
    \circ (\omega^*\eta\omega_*)
    =(\omega^{\pr*}\beta_*)
    \circ(\omega^{\pr*}\vartheta\omega_*)
    \circ(\alpha\omega_*)
    =(\vartheta^{\sharp}\beta^{\sharp})
    \circ (\omega^*\eta\omega_*)
    \colon \omega^*\omega_*\Rightarrow\omega^{\pr*}\omega^{\pr}_*$.
    From these equalities,
    along with \eqref{eq-fundmoroverGropk-alphasharpequiv} for $z=\pi$
    and \eqref{eq-dualityoverGropk-varthetasharpequiv} for $z=\pi$,
    we deduce that
    $\breve{\theta}$
    is indeed a $2$-isomorphism from $\breve{f}$ to $\breve{g}$.
  \end{itemize}
  It is straightforward to check that
  $\mathsf{Aut}_{\kos{K}}$
  preserves identity $1$-cells and identity $2$-cells.
  \begin{itemize}
    \item 
    Let $\vartheta\colon (\kos{f},\alpha)\cong (\kos{g},\beta)$
    and
    $\tilde{\vartheta}\colon (\kos{g},\beta)\cong (\tilde{\kos{g}},\tilde{\beta})$
    be $2$-cells in $\mathsf{pGroCat}_{\mathsf{pt}}(\kos{K})$,
    with vertical composite
    $\bar{\vartheta}\coloneqq \tilde{\vartheta}\circ \vartheta
    \colon (\kos{f},\alpha)\cong (\tilde{\kos{g}},\tilde{\beta})$.
    Put $\pi\coloneqq \omega^*\omega_*(\kappa)$.
    One can check that
    $(\vartheta^{\sharp}\tilde{\vartheta}^{\sharp})\circ (\omega^*\eta\omega_*)
    =\bar{\vartheta}^{\sharp}\colon \omega^*\omega_*\Rightarrow\vphantom{\id}\smash[t]{\id}_{\kos{K}}^{\vphantom{*}}$,
    using \eqref{eq-fundmoroverGropk-matealpha}.
    From this,
    along with \eqref{eq-dualityoverGropk-varthetasharpequiv}
    for $z=\pi$,
    we deduce that $\breve{\tilde{\theta}}\star\breve{\theta}=\breve{\bar{\theta}}$.
    This shows that $\mathsf{Aut}_{\kos{K}}$
    preserves vertical composition of $2$-cells.
    
    \item 
    Let $(\kos{f},\alpha)\colon (\kos{t}^{\pr},\omega^{\pr})\mor(\kos{t},\omega)$
    and
    $(\kos{f}^{\pr},\alpha^{\pr})\colon (\kos{t}^{\ppr},\omega^{\ppr})\mor(\kos{t}^{\pr},\omega^{\pr})$
    be $1$-cells in $\mathsf{pGroCat}_{\mathsf{pt}}(\kos{K})$,
    with composite denoted by
    $(\kos{f}^{\ppr},\alpha^{\ppr})\colon (\kos{t}^{\ppr},\omega^{\ppr})\mor(\kos{t},\omega)$,
    where $\kos{f}^{\ppr}\coloneqq \kos{f}\kos{f}^{\pr}$
    and $\alpha^{\ppr}\coloneqq (\kos{f}\alpha^{\pr})\circ \alpha\colon \omega\cong \kos{f}\kos{f}^{\pr}\omega^{\ppr}$.
    One can check that 
    $\alpha^{\ppr}_*=\alpha^{\pr}_*\circ (\kos{f}^{\pr*}\alpha_*)
    \colon \kos{f}^{\pr*}\kos{f}^*\omega_*\Rightarrow\omega^{\ppr}_*$
    and that
    $\alpha^{\ppr\sharp}=\alpha^{\pr\sharp}\circ \alpha^{\sharp}
    \colon \omega^*\omega_*\Rightarrow\omega^{\ppr*}\omega^{\ppr}_*$,
    using \eqref{eq-fundmoroverGropk-matealpha}.
    Thus we have $\breve{f}^{\ppr}=\breve{f}\circ \breve{f}^{\pr}
    \colon \Spec(\omega^{\ppr*}\omega^{\ppr}_*(\kappa))
    \to \Spec(\omega^*\omega_*(\kappa))$.
    This shows that $\mathsf{Aut}_{\kos{K}}$
    preserves composition of $1$-cells.
  
    \item
    Let $(\kos{f},\alpha)$,
    $(\kos{g},\beta)\colon (\kos{t}^{\pr},\omega^{\pr})\mor (\kos{t},\omega)$,
    $(\kos{f}^{\pr},\alpha^{\pr})$,
    $(\kos{g}^{\pr},\beta^{\pr})\colon (\kos{t}^{\ppr},\omega^{\ppr})\mor (\kos{t}^{\pr},\omega^{\pr})$
    be $1$-cells and
    $\vartheta\colon (\kos{f},\alpha)\cong (\kos{g},\beta)$,
    $\vartheta^{\pr}\colon (\kos{f}^{\pr},\alpha^{\pr})\cong (\kos{g}^{\pr},\beta^{\pr})$
    be $2$-cells in $\mathsf{pGroCat}_{\mathsf{pt}}(\kos{K})$.
    Write
    $(\kos{f}^{\ppr},\alpha^{\ppr})$,
    $(\kos{g}^{\ppr},\beta^{\ppr})\colon (\kos{t}^{\ppr},\omega^{\ppr})\mor(\kos{t},\omega)$
    for the composite $1$-cells
    and
    $\vartheta^{\ppr}\colon (\kos{f}^{\ppr},\alpha^{\ppr})\cong(\kos{g}^{\ppr},\beta^{\ppr})$
    for the horizontal composite of $\vartheta$ and $\vartheta^{\pr}$.
    One can check that
    $\vartheta^{\ppr\sharp}
    =(\vartheta^{\sharp}\vartheta^{\pr\sharp})
    \circ (\omega^*\omega_*\beta^{\sharp})
    \circ (\omega^*\eta\omega_*)
    \colon \omega^*\omega_*\Rightarrow\vphantom{\id}\smash[t]{\id}_{\kos{K}}^{\vphantom{*}}$,
    using \eqref{eq-fundmoroverGropk-matealpha}.
    From this,
    along with \eqref{eq-dualityoverGropk-varthetasharpequiv}
    for $z=\pi\coloneqq \omega^*\omega_*(\kappa)$,
    we deduce that
    $\breve{\theta}^{\ppr}
    =(\breve{g}\circ\breve{\theta}^{\pr})\star \breve{\theta}$.
    This shows that $\mathsf{Aut}_{\kos{K}}$
    preserves horizontal composition of $2$-cells.
  \end{itemize}
  This completes the construction of
  $\mathsf{Aut}_{\kos{K}}\colon \mathsf{pGroCat}_{\mathsf{pt}}(\kos{K})\to \mathsf{pGroObj}(\kos{K})$.
\end{proof}

\begin{proof}[Proof of Theorem~\ref{thm-dualityoverGropk-equiv21cat}]
  We show that the pseudofunctor
  $\mathsf{Rep}_{\kos{K}}$
  in Proposition~\ref{prop-dualityoverGropk-RepKpseudofunct}
  and the $2$-functor
  $\mathsf{Aut}_{\kos{K}}$
  in Proposition~\ref{prop-dualityoverGropk-AutK2funct}
  are mutually quasi-inverse equivalences of $(2,1)$-categories;
  that is, the composites
  $\mathsf{Aut}_{\kos{K}}\circ \mathsf{Rep}_{\kos{K}}$
  and $\mathsf{Rep}_{\kos{K}}\circ \mathsf{Aut}_{\kos{K}}$
  are equivalent to the identity $2$-functors
  of $\mathsf{pGroObj}(\kos{K})$
  and $\mathsf{pGroCat}_{\mathsf{pt}}(\kos{K})$,
  respectively.
  
  The composite
  $\mathsf{Aut}_{\kos{K}}\circ \mathsf{Rep}_{\kos{K}}\colon 
  \mathsf{pGroObj}(\kos{K})\to \mathsf{pGroObj}(\kos{K})$
  sends each $0$-cell $\Spec(\pi)$
  to the $0$-cell $\Spec(\kappa\otimes \pi)$,
  and we have an isomorphism
  $\imath_{\pi}\colon \Spec(\kappa\otimes \pi)\cong \Spec(\pi)$
  of Grothendieck objects in $\kos{K}$;
  see Remark~\ref{rem-GroMainthm-reconstruction}.
  It is straightforward to check that
  $\mathsf{Aut}_{\kos{K}}\circ \mathsf{Rep}_{\kos{K}}$
  is a $2$-functor, and that the isomorphisms $\imath_{\pi}$
  constitute a $2$-natural isomorphism
  from $\mathsf{Aut}_{\kos{K}}\circ \mathsf{Rep}_{\kos{K}}$
  to the identity $2$-functor of $\mathsf{pGroObj}(\kos{K})$.

  We now construct an equivalence
  from the composite
  $\mathsf{Rep}_{\kos{K}}\circ \mathsf{Aut}_{\kos{K}}$
  to the identity $2$-functor of $\mathsf{pGroCat}_{\mathsf{pt}}(\kos{K})$
  as follows.
  Such an equivalence is an invertible strong transformation
  in the sense of
  \cite[Definition~4.3.1, Lemma~4.3.9(2), Corollary~4.4.13, and Definition~5.1.18]{Johnson2021};
  it is an equivalence internal to the bicategory of
  pseudofunctors from $\mathsf{pGroCat}_{\mathsf{pt}}(\kos{K})$ to itself,
  strong transformations, and modifications.
  We take the component $1$-cell of this transformation
  at each $0$-cell $(\kos{t},\omega)$ to be the $1$-cell
  $(\kos{t}_{\omega^*\omega_*(\kappa)},\kos{e}_{\omega^*\omega_*(\kappa)})
  \simeq
  (\kos{t},\omega)$
  consisting of the equivalence
  $\breve{\omega}\colon \kos{t}_{\omega^*\omega_*(\kappa)}\simeq\kos{t}$
  and 
  the invertible Grothendieck $\kos{K}$-transformation
  $\omega\cong \breve{\omega}\kos{e}_{\omega^*\omega_*(\kappa)}$
  in Theorem~\ref{thm-Grocatoverpk-mainthm}.
  We take its component $2$-cell
  at each $1$-cell 
  $(\kos{f},\alpha)\colon (\kos{t}^{\pr},\omega^{\pr})\mor (\kos{t},\omega)$
  to be the invertible Grothendieck $\kos{K}$-transformation
  $\breve{\alpha}\colon \breve{\omega}\breve{\kos{f}}\cong \kos{f}\breve{\omega}^{\pr}$
  in Proposition~\ref{prop-fundmoroverGropk-realization}.
  It is straightforward to check,
  using \cite[Proposition~6.2.16]{Johnson2021}, that
  these data indeed define an invertible strong transformation.
  This completes the proof of Theorem~\ref{thm-dualityoverGropk-equiv21cat}.
\end{proof}

%% file: PKGroCat.tex
\chapter{Prekosmic Grothendieck Categories}
\label{chap-PKGroCat}

In this chapter, we gather the dualities
established in Section~\ref{sec-dualityoverGropk},
one for each Grothendieck prekosmos $\kos{K}$,
into a single duality.
We begin with the following definitions;
recall Definitions~\ref{def-Groprekos-Grocatoverpk},
\ref{def-Groobjinpk-Groobj},
and~\ref{def-dualityoverGropk-ptGrocat}.

\begin{definition}
  \label{def-PKGroCat-PKGroCat}
  A \emph{prekosmic Grothendieck category}
  is a pair $(\kos{K},\kos{t})$
  consisting of a Grothendieck prekosmos $\kos{K}$
  and a Grothendieck $\kos{K}$-category
  $\kos{T}\xmor{\kos{t}}\kos{K}$.
\end{definition}

\begin{definition}
  \label{def-PKGroCat-PKGroobj}
  A \emph{prekosmic Grothendieck object}
  is a pair $(\kos{K},\Spec(\pi))$
  consisting of a Grothendieck prekosmos $\kos{K}$
  and a Grothendieck object $\Spec(\pi)$ in $\kos{K}$.
\end{definition}

\begin{definition}
  \label{def-PKGroCat-ptPKGroCat}
  A \emph{pointed prekosmic Grothendieck category}
  is a tuple $(\kos{K},\kos{t},\omega)$
  consisting of a Grothendieck prekosmos $\kos{K}$,
  a Grothendieck $\kos{K}$-category
  $\kos{T}\xmor{\kos{t}}\kos{K}$,
  and a surjective fiber point $\omega\colon \vphantom{\id}\smash[t]{\id}_{\kos{K}}^{\vphantom{*}}\mor\kos{t}$ of $\kos{t}$.
\end{definition}

The change of base along an \'{e}tale Grothendieck morphism
$\kos{K}^{\pr}\mor\kos{K}$
between Grothendieck prekosmoi
sends each Grothendieck object in $\kos{K}$
to one in $\kos{K}^{\pr}$,
generalizing the change of base
of a flat affine group $R$-scheme
along a morphism $\Spec(R^{\pr})\to \Spec(R)$
of affine schemes.
We write $\mathsf{pGroKos}_{\etsf}$
for the $(2,1)$-category 
of Grothendieck prekosmoi,
\'{e}tale Grothendieck morphisms,
and invertible Grothendieck transformations.
Prekosmic Grothendieck objects
form a $(2,1)$-category $\mathsf{pGroObj}$
over $\mathsf{pGroKos}_{\etsf}$,
and the strict fiber 
over each Grothendieck prekosmos $\kos{K}$
is $2$-isomorphic to
$\mathsf{pGroObj}(\kos{K})$
in \eqref{eq-Groobjinpk-pGroObjK21cat}.
Pointed prekosmic Grothendieck categories
likewise form a $(2,1)$-category $\mathsf{pGroCat}_{\mathsf{pt}}$
over $\mathsf{pGroKos}_{\etsf}$,
and the strict fiber 
over each Grothendieck prekosmos $\kos{K}$
is $2$-isomorphic to
$\mathsf{pGroCat}_{\mathsf{pt}}(\kos{K})$
in Definition~\ref{def-dualityoverGropk-pGroCatptK}.
The following theorem gathers,
for each Grothendieck prekosmos $\kos{K}$,
the mutually quasi-inverse equivalences
$\mathsf{Rep}_{\kos{K}}$ and $\mathsf{Aut}_{\kos{K}}$
in Theorem~\ref{thm-dualityoverGropk-equiv21cat}.

\begin{theorem}
  \label{thm-PKGroCat-equiv21catoveret}
  We have mutually quasi-inverse equivalences
  \begin{equation*}
    \vcenter{\hbox{
      \xymatrix@R=25pt@C=0pt{
        \mathsf{pGroObj}
        \ar@<1ex>[rr]^-{\mathsf{Rep}}
        \ar[dr]_-{\mathsf{p}_{\mathsf{obj}}}
        \ar@{}[rr]|-{\simeq}
        &\text{ }
        &\mathsf{pGroCat}_{\mathsf{pt}}
        \ar@<1ex>[ll]^-{\mathsf{Aut}}
        \ar[dl]^-{\mathsf{p}_{\mathsf{ptcat}}}
        \\%1
        \text{ }
        &\mathsf{pGroKos}_{\etsf}
        &\text{ }
      }
    }}
  \end{equation*}
  of $(2,1)$-categories over $\mathsf{pGroKos}_{\etsf}$.
\end{theorem}

\begin{remark}
  \label{rem-PKGroCat-cleavagetransport}
  The $2$-functor $\mathsf{p}_{\mathsf{obj}}$
  is a split $2$-fibration,
  with the canonical cleavage described in Section~\ref{sec-PKGroobj}.
  The equivalences in Theorem~\ref{thm-PKGroCat-equiv21catoveret}
  transport this cleavage to $\mathsf{p}_{\mathsf{ptcat}}$,
  making it a cloven fibration of bicategories.
  We refer to \cite{Buckley2014} for the terminology.
\end{remark}

The following lemma ensures that
the $(2,1)$-category $\mathsf{pGroKos}_{\etsf}$ is well defined.

\begin{lemma}
  \label{lem-PKGroCat-etalecomposite}
  \'{E}tale Grothendieck morphisms between Grothendieck prekosmoi
  are closed under composition.
\end{lemma}
\begin{proof}
  Let $\phi\colon \kos{K}^{\pr}\mor\kos{K}$,
  $\phi^{\pr}\colon \kos{K}^{\ppr}\mor\kos{K}^{\pr}$
  be \'{e}tale Grothendieck morphisms between Grothendieck prekosmoi,
  and $\varphi$, $\varphi^{\pr}$, and $\varphi^{\ppr}$
  be the natural transformations in \eqref{eq-Groprekos-projformula}
  associated to $\phi$, $\phi^{\pr}$, and $\phi^{\ppr}\coloneqq\phi\phi^{\pr}$, respectively.
  One can check that
  $\varphi^{\ppr}_{z,x^{\ppr}}=
  \phi_*(\varphi^{\pr}_{\phi^*(z),x^{\ppr}})\circ \varphi_{z,\phi^{\pr}_*(x^{\ppr})}$
  for every $z\in \CK$, $x^{\ppr}\in\CK^{\ppr}$,
  and since $\varphi$ and $\varphi^{\pr}$ are invertible,
  we deduce that $\varphi^{\ppr}$ is invertible.
  The right adjoint part $\phi^{\ppr}_*=\phi_*\phi^{\pr}_*$
  is conservative, since $\phi_*$ and $\phi^{\pr}_*$ are.
  Hence $\phi^{\ppr}$ is \'{e}tale.
\end{proof}

\section{Prekosmic Grothendieck objects}
\label{sec-PKGroobj}

Recall from Definition~\ref{def-PKGroCat-PKGroobj}
the definition of a prekosmic Grothendieck object.
Let $(\kos{K},\Spec(\pi))$,
$(\kos{K}^{\pr},\Spec(\pi^{\pr}))$
be prekosmic Grothendieck objects.
A \emph{morphism}
$(\phi,f)\colon (\kos{K}^{\pr},\Spec(\pi^{\pr}))\to (\kos{K},\Spec(\pi))$
consists of an \'{e}tale Grothendieck morphism
$\phi\colon \kos{K}^{\pr}\mor\kos{K}$
and a morphism
$f\colon \Spec(\pi^{\pr})\to \Spec(\phi^*(\pi))$
of Grothendieck objects in $\kos{K}^{\pr}$;
see Remark~\ref{rem-PKGroobj-etalebasechange}.
Let 
$(\psi,g)\colon (\kos{K}^{\pr},\Spec(\pi^{\pr}))\to (\kos{K},\Spec(\pi))$
be a morphism.
A \emph{$2$-isomorphism}
$(\chi,\theta)\colon (\phi,f)\cong (\psi,g)$
is a pair consisting of 
an invertible Grothendieck transformation $\chi\colon \phi\cong \psi$
and a $2$-isomorphism
$\theta\colon f\cong \chi_{\pi}g$
between morphisms
$f$, $\chi_{\pi}g\colon \Spec(\pi^{\pr})\to \Spec(\phi^*(\pi))$.
Prekosmic Grothendieck objects,
morphisms between them,
and $2$-isomorphisms between these morphisms
form a $(2,1)$-category,
which we denote by
\begin{equation}
  \label{eq-PKGroobj-21cat}
  \mathsf{pGroObj}.
\end{equation}
We describe the rest of the structure of $\mathsf{pGroObj}$ as follows.

Let 
$(\tilde{\psi},\tilde{g})\colon (\kos{K}^{\pr},\Spec(\pi^{\pr}))\to (\kos{K},\Spec(\pi))$
be a morphism
and $(\tilde{\chi},\tilde{\theta})\colon (\psi,g)\cong (\tilde{\psi},\tilde{g})$
be a $2$-isomorphism.
The vertical composite of
$(\chi,\theta)$ and $(\tilde{\chi},\tilde{\theta})$
is the pair $(\tilde{\chi}\circ \chi,(\chi_{\pi}\tilde{\theta})\star\theta)
\colon (\phi,f)\cong (\tilde{\psi},\tilde{g})$
consisting of
$\tilde{\chi}\circ\chi\colon \phi\cong \tilde{\psi}$
and
$(\chi_{\pi}\tilde{\theta})\star\theta\colon f\cong \chi_{\pi}\tilde{\chi}_{\pi}\tilde{g}$.
The identity $2$-isomorphism of $(\phi,f)$
is the pair consisting of the identity Grothendieck transformation of $\phi$
and the identity $2$-isomorphism of $f$.
Each $2$-isomorphism $(\chi,\theta)\colon (\phi,f)\cong (\psi,g)$
is indeed invertible with respect to the vertical composition;
its inverse is the pair
$(\chi^{-1},\chi_{\pi}^{-1}\theta^{-1})\colon (\psi,g)\cong (\phi,f)$
consisting of
$\chi^{-1}\colon \psi\cong \phi$
and
$\chi_{\pi}^{-1}\theta^{-1}\colon g\cong \chi_{\pi}^{-1}f$,
where $\theta^{-1}\colon \chi_{\pi}g\cong f$
is the inverse of $\theta$ with respect to $\star$.

Let $(\phi^{\pr},f^{\pr})$, $(\psi^{\pr},g^{\pr})\colon (\kos{K}^{\ppr},\Spec(\pi^{\ppr}))\to (\kos{K}^{\pr},\Spec(\pi^{\pr}))$
be morphisms.
The composite of $(\phi,f)$ and $(\phi^{\pr},f^{\pr})$ is the pair
$(\phi\phi^{\pr},\phi^{\pr*}(f)f^{\pr})$
consisting of
$\phi\phi^{\pr}\colon \kos{K}^{\ppr}\mor\kos{K}$
and 
$\phi^{\pr*}(f)f^{\pr}\colon \Spec(\pi^{\ppr})\to 
\Spec(\phi^{\pr*}\phi^*(\pi))$.
The identity morphism of $(\kos{K},\Spec(\pi))$
is the pair consisting of 
the identity Grothendieck morphism $\vphantom{\id}\smash[t]{\id}_{\kos{K}}^{\vphantom{*}}$ of $\kos{K}$
and the identity morphism of $\Spec(\pi)$.
Let 
$(\chi^{\pr},\theta^{\pr})\colon (\phi^{\pr},f^{\pr})\cong (\psi^{\pr},g^{\pr})$
be a $2$-isomorphism.
The horizontal composite of 
$(\chi,\theta)$ and $(\chi^{\pr},\theta^{\pr})$
is the pair
$(\chi\chi^{\pr},\phi^{\pr*}(\theta)\theta^{\pr})\colon 
(\phi\phi^{\pr},\phi^{\pr*}(f)f^{\pr})\cong (\psi\psi^{\pr},\psi^{\pr*}(g)g^{\pr})$
consisting of
$\chi\chi^{\pr}\colon \phi\phi^{\pr}\cong \psi\psi^{\pr}$
and the horizontal composite $\phi^{\pr*}(\theta)\theta^{\pr}$ of
$\phi^{\pr*}(\theta)\colon \phi^{\pr*}(f)\cong \phi^{\pr*}(\chi_{\pi})\phi^{\pr*}(g)$
and
$\theta^{\pr}\colon f^{\pr}\cong \chi^{\pr}_{\pi^{\pr}}g^{\pr}$.

\begin{remark}
  \label{rem-PKGroobj-etalebasechange}
  When we define a morphism of prekosmic Grothendieck objects
  $(\phi,f)\colon (\kos{K}^{\pr},\Spec(\pi^{\pr}))\to (\kos{K},\Spec(\pi))$,
  we require $\phi$ to be \'{e}tale,
  for then $\Spec(\phi^*(\pi))$ is
  a Grothendieck object in $\kos{K}^{\pr}$.
  Indeed, the functor
  $\phi_*(\slot\otimes^{\pr}\phi^*(\pi))$
  preserves coreflexive equalizers,
  as it is isomorphic to $\phi_*(\slot)\otimes \pi$
  by \eqref{eq-Groprekos-projformula},
  and since $\phi_*$ is conservative and preserves coreflexive equalizers,
  we deduce that the functor $\slot\otimes^{\pr} \phi^*(\pi)$
  preserves coreflexive equalizers.
  The argument still applies
  when we are merely given 
  a strong monoidal functor $\phi^*\colon \kos{K}\to\kos{K}^{\pr}$,
  a functor $\phi_*\colon \CK^{\pr}\to \CK$ that is conservative
  and preserves coreflexive equalizers,
  and a natural isomorphism
  $\varphi\colon
  \slot\otimes \phi_*(\slot)
  \cong \phi_*(\phi^*(\slot)\otimes^{\pr}\slot)$.
  Thus such data $(\phi^*,\phi_*,\varphi)$,
  in which $\phi^*$ need not be left adjoint to $\phi_*$,
  could serve as $1$-cells of
  $\mathsf{pGroKos}_{\etsf}$,
  with $\mathsf{pGroObj}$ and $\mathsf{pGroCat}_{\mathsf{pt}}$
  enlarged accordingly.
  We do not pursue this extension in this paper.
\end{remark}

We denote by 
$\mathsf{p}_{\mathsf{obj}}\colon \mathsf{pGroObj}\to\mathsf{pGroKos}_{\etsf}$
the $2$-functor that sends each cell to its first component.
It is immediate that $\mathsf{pGroObj}$
is a split fibered $2$-category over $\mathsf{pGroKos}_{\etsf}$
with respect to $\mathsf{p}_{\mathsf{obj}}$,
whose cleavage is given as follows.
\begin{itemize}
  \item 
  For a $0$-cell $(\kos{K},\Spec(\pi))$ in $\mathsf{pGroObj}$
  and a $1$-cell $\phi\colon \kos{K}^{\pr}\mor\kos{K}$ in $\mathsf{pGroKos}_{\etsf}$,
  the cartesian $1$-cell
  $(\kos{K}^{\pr},\Spec(\phi^*(\pi)))\to (\kos{K},\Spec(\pi))$
  over $\phi$
  is the pair consisting of 
  $\phi$ and the identity morphism of $\Spec(\phi^*(\pi))$.

  \item
  For a $1$-cell $(\psi,g)\colon (\kos{K}^{\pr},\Spec(\pi^{\pr}))\to (\kos{K},\Spec(\pi))$
  in $\mathsf{pGroObj}$
  and a $2$-cell 
  $\chi\colon \phi\cong \psi$
  between $1$-cells $\phi$, $\psi\colon \kos{K}^{\pr}\mor\kos{K}$
  in $\mathsf{pGroKos}_{\etsf}$,
  the cartesian $2$-cell
  $(\phi,\chi_{\pi}g)\cong (\psi,g)$
  over $\chi$
  is the pair consisting of 
  $\chi$ and the identity $2$-isomorphism
  of $\chi_{\pi}g\colon \Spec(\pi^{\pr})\to \Spec(\phi^*(\pi))$.
\end{itemize}
The strict fiber of $\mathsf{p}_{\mathsf{obj}}$
over each Grothendieck prekosmos $\kos{K}$
is $2$-isomorphic to $\mathsf{pGroObj}(\kos{K})$,
defined in \eqref{eq-Groobjinpk-pGroObjK21cat}.
We illustrate $\mathsf{p}_{\mathsf{obj}}$ as follows.
\begin{equation}
  \label{eq-PKGroobj-pobj21funct}
  \begin{tikzcd}[row sep=22pt,column sep=20pt]
    \mathsf{pGroObj}
    \ar[d,"{\mathsf{p}_{\mathsf{obj}}}"']
    &(\kos{K}^{\pr},\Spec(\pi^{\pr}))
    \ar[rr,shift left=0ex,bend left=1.2pc,"{(\psi,g)}" description]
    \ar[rr,shift right=0ex,bend right=1.2pc,"{(\phi,f)}" description]
    &\text{ }
    &(\kos{K},\Spec(\pi))
    \\%1
    {\mathsf{pGroKos}_{\etsf}}
    &\kos{K}^{\pr}
    \ar[rr,filled,shift left=0.3ex,bend left=1.2pc,"\psi" description]
    \ar[rr,filled,shift right=0.3ex,bend right=1.2pc,"\phi" description]
    &\text{ }
    \ar[u,2cell,pos=1.17,"{\,{(\chi,\theta)}}"']
    \ar[u,2cell,pos=-0.05,"{\,\chi}"']
    &\kos{K}
  \end{tikzcd}
\end{equation}

\section{Change of base along a Grothendieck morphism}
\label{sec-changebaseGro}
In this section, we introduce the notation for
changing the base Grothendieck prekosmos,
which we use throughout the rest of this chapter.
Let $\phi\colon \kos{K}^{\pr}\mor\kos{K}$
be a Grothendieck morphism between Grothendieck prekosmoi
$\kos{K}$, $\kos{K}^{\pr}$.
We write
$(\slot)_{\phi}\colon \mathtt{pGroKos}/\kos{K}^{\pr}\to \mathtt{pGroKos}/\kos{K}$
for the change of base $2$-functor
induced by $\phi$
between the pseudoslice $2$-categories.
We explain the $2$-functor $(\slot)_{\phi}$ in detail.
\begin{itemize}
  \item 
  Let $\kos{T}^{\pr}\xmor{\kos{t}^{\pr}}\kos{K}^{\pr}$
  be a Grothendieck $\kos{K}^{\pr}$-prekosmos.
  When we regard the composite 
  $\phi\kos{t}^{\pr}\colon \kos{T}^{\pr}\mor \kos{K}$
  as a Grothendieck $\kos{K}$-prekosmos,
  we write it as $\kos{T}^{\pr}\xmor{\kos{t}^{\pr}_{\phi}}\kos{K}$,
  or simply as $\kos{t}^{\pr}_{\phi}$.

  \item
  In particular,
  when we regard $\phi$
  as a Grothendieck $\kos{K}$-prekosmos,
  we write it as $\kos{K}^{\pr}\xmor{\phi}\kos{K}$, or simply as $\phi$.
  
  \item
  Let $\kos{f}^{\pr}\colon \tilde{\kos{t}}^{\pr}\mor \kos{t}^{\pr}$
  be a Grothendieck $\kos{K}^{\pr}$-morphism
  between Grothendieck $\kos{K}^{\pr}$-prekosmoi
  $\kos{t}^{\pr}$, $\tilde{\kos{t}}^{\pr}$;
  it consists of a Grothendieck morphism
  $\kos{f}^{\pr}\colon \tilde{\kos{T}}^{\pr}\mor\kos{T}^{\pr}$
  and an invertible Grothendieck transformation
  $\what{\kos{f}^{\pr}}\colon \tilde{\kos{t}}^{\pr}\cong \kos{t}^{\pr}\kos{f}^{\pr}$.
  We write
  $\kos{f}^{\pr}_{\phi}\colon \tilde{\kos{t}}^{\pr}_{\phi}\mor \kos{t}^{\pr}_{\phi}$
  for the Grothendieck $\kos{K}$-morphism
  whose underlying Grothendieck morphism is that of
  $\kos{f}^{\pr}\colon \tilde{\kos{t}}^{\pr}\mor \kos{t}^{\pr}$
  and whose invertible Grothendieck transformation
  $\what{\kos{f}^{\pr}_{\phi}}\colon 
  \tilde{\kos{t}}^{\pr}_{\phi}\cong \kos{t}^{\pr}_{\phi}\kos{f}^{\pr}_{\phi}$
  is $\phi\what{\kos{f}^{\pr}}\colon 
  \phi\tilde{\kos{t}}^{\pr}\cong \phi\kos{t}^{\pr}\kos{f}^{\pr}$.
  We illustrate
  $\kos{f}^{\pr}_{\phi}\colon \tilde{\kos{t}}^{\pr}_{\phi}\mor\kos{t}^{\pr}_{\phi}$
  and its description as follows.
  \begin{equation*}
    \begin{tikzcd}[row sep=7pt,column sep=7pt]
      \tilde{\kos{T}}^{\pr}
      \ar[rr,filled,"\kos{f}^{\pr}_{\phi}"]
      \ar[ddr,filled,"\tilde{\kos{t}}^{\pr}_{\phi}"']
      &\text{ }
      &\kos{T}^{\pr}
      \ar[ddl,filled,"\kos{t}^{\pr}_{\phi}"]
      \\%1
      \text{ }
      \ar[rr,2cell,shift right=0.7ex,"{\!\!\what{\kos{f}^{\pr}_{\phi}}}"]
      &\text{ }
      &\text{ }
      \\%2
      \text{ }
      &\kos{K}
      &\text{ }
    \end{tikzcd}
    \,\,=\,\,
    \begin{tikzcd}[row sep=3pt,column sep=8pt]
      \tilde{\kos{T}}^{\pr}
      \ar[rr,filled,"\kos{f}^{\pr}"]
      \ar[ddr,filled,"\tilde{\kos{t}}^{\pr}"']
      &\text{ }
      &\kos{T}^{\pr}
      \ar[ddl,filled,"\kos{t}^{\pr}"]
      \\%1
      \text{ }
      \ar[rr,2cell,shift right=1ex,pos=0.475,"{\,\what{\kos{f}^{\pr}}}"]
      &\text{ }
      &\text{ }
      \\%2
      \text{ }
      &\kos{K}^{\pr}
      \ar[dd,filled,"\phi"]
      &\text{ }
      \\%3
      \text{ }
      &\text{ }
      &\text{ }
      \\%4
      \text{ }
      &\kos{K}
      &\text{ }
    \end{tikzcd}
  \end{equation*}

  \item
  Let $\kos{t}^{\pr}$,
  $\tilde{\kos{t}}^{\pr}$
  be Grothendieck $\kos{K}^{\pr}$-prekosmoi,
  $\kos{f}^{\pr}$, $\kos{g}^{\pr}\colon \tilde{\kos{t}}^{\pr}\mor \kos{t}^{\pr}$
  be Grothendieck $\kos{K}^{\pr}$-morphisms,
  and $\vartheta^{\pr}\colon \kos{f}^{\pr}\Rightarrow \kos{g}^{\pr}$
  be a Grothendieck $\kos{K}^{\pr}$-transformation.
  Then $\vartheta^{\pr}$ is a Grothendieck $\kos{K}$-transformation
  from $\kos{f}^{\pr}_{\phi}$ to $\kos{g}^{\pr}_{\phi}$,
  and we write it as
  $\vartheta^{\pr}_{\phi}\colon \kos{f}^{\pr}_{\phi}\Rightarrow\kos{g}^{\pr}_{\phi}$.
\end{itemize}
The pair $(\phi,\what{\phi})$
consisting of $\phi$
and the equality $\what{\phi}\colon \phi=\phi$
is a Grothendieck $\kos{K}$-morphism from $\phi$ to $\vphantom{\id}\smash[t]{\id}_{\kos{K}}^{\vphantom{*}}$;
we omit $\what{\phi}$ from the notation and write
$\phi\colon \phi\mor\vphantom{\id}\smash[t]{\id}_{\kos{K}}^{\vphantom{*}}$.

Let 
$\psi\colon \kos{K}^{\pr}\mor\kos{K}$
be a Grothendieck morphism
and $\chi\colon \phi\cong \psi$ be an invertible Grothendieck transformation.
We explain the $2$-natural isomorphism
$(\slot)_{\chi}$
induced by $\chi$
from $(\slot)_{\phi}$ to $(\slot)_{\psi}$.
\begin{itemize}
  \item 
  Let $\kos{T}^{\pr}\xmor{\kos{t}^{\pr}}\kos{K}^{\pr}$
  be a Grothendieck $\kos{K}^{\pr}$-prekosmos.
  We write $\kos{t}^{\pr}_{\chi}\colon \kos{t}^{\pr}_{\phi}\mor \kos{t}^{\pr}_{\psi}$
  for the Grothendieck $\kos{K}$-morphism
  whose underlying Grothendieck morphism is
  the identity Grothendieck morphism $\vphantom{\id}\smash[t]{\id}_{\kos{T}^{\pr}}^{\vphantom{*}}$ of $\kos{T}^{\pr}$
  and whose invertible Grothendieck transformation
  $\what{\kos{t}^{\pr}_{\chi}}$ is
  $\chi\kos{t}^{\pr}\colon \phi\kos{t}^{\pr}\cong \psi\kos{t}^{\pr}$.
  We illustrate $\kos{t}^{\pr}_{\chi}\colon \kos{t}^{\pr}_{\phi}\mor\kos{t}^{\pr}_{\psi}$
  and its left adjoint part
  $\kos{t}^{\pr*}_{\chi}\colon \kos{t}^{\pr*}_{\psi}\to \kos{t}^{\pr*}_{\phi}$,
  as well as their descriptions as follows.
  \begin{equation*}
    \begin{tikzcd}[row sep=5pt,column sep=5pt]
      \kos{T}^{\pr}
      \ar[rr,filled,"\kos{t}^{\pr}_{\chi}"]
      \ar[ddr,filled,"\kos{t}^{\pr}_{\phi}"']
      &\text{ }
      &\kos{T}^{\pr}
      \ar[ddl,filled,"\kos{t}^{\pr}_{\psi}"]
      \\%1
      \text{ }
      \ar[rr,2cell,shift right=1.2ex,"{\!\!\what{\kos{t}^{\pr}_{\chi}}}"]
      &\text{ }
      &\text{ }
      \\%2
      \text{ }
      &\kos{K}
      &\text{ }
    \end{tikzcd}
    =\,\,
    \begin{tikzcd}[row sep=2pt,column sep=1pt]
      \kos{T}^{\pr}
      \ar[dd,filled,"\kos{t}^{\pr}"']
      \ar[rr,equal]
      &\text{ }
      &\kos{T}^{\pr}
      \ar[dd,filled,"\kos{t}^{\pr}"]
      \\%1
      \text{ }
      &\text{ }
      &\text{ }
      \\%2
      \kos{K}^{\pr}
      \ar[ddr,filled,"\phi"']
      \ar[rr,equal]
      &\text{ }
      &\kos{K}^{\pr}
      \ar[ddl,filled,"\psi"]
      \\%3
      \text{ }
      \ar[rr,2cell,shift right=0.5ex,"{\!\chi}"]
      &\text{ }
      &\text{ }
      \\%4
      \text{ }
      &\kos{K}
      &\text{ }
    \end{tikzcd}
    \qquad
    \begin{tikzcd}[row sep=5pt,column sep=5pt]
      \kos{T}^{\pr}
      &\text{ }
      &\kos{T}^{\pr}
      \ar[ll,"\kos{t}^{\pr*}_{\chi}"']
      \\%1
      \text{ }
      \ar[rr,2cell,shift right=1.2ex,"{\!\what{\kos{t}^{\pr*}_{\chi}}}"]
      &\text{ }
      &\text{ }
      \\%2
      \text{ }
      &\kos{K}
      \ar[uul,"\kos{t}^{\pr*}_{\phi}"]
      \ar[uur,"\kos{t}^{\pr*}_{\psi}"']
      &\text{ }
    \end{tikzcd}
    =\,\,
    \begin{tikzcd}[row sep=2pt,column sep=1pt]
      \kos{T}^{\pr}
      \ar[rr,equal]
      &\text{ }
      &\kos{T}^{\pr}
      \\%1
      \text{ }
      &\text{ }
      &\text{ }
      \\%2
      \kos{K}^{\pr}
      \ar[uu,"\kos{t}^{\pr*}"]
      \ar[rr,equal]
      &\text{ }
      &\kos{K}^{\pr}
      \ar[uu,"\kos{t}^{\pr*}"']
      \\%3
      \text{ }
      \ar[rr,2cell,shift right=0.5ex,"{\!\chi}"]
      &\text{ }
      &\text{ }
      \\%4
      \text{ }
      &\kos{K}
      \ar[uul,"\phi^*"]
      \ar[uur,"\psi^*"']
      &\text{ }
    \end{tikzcd}
  \end{equation*}

  \item
  In particular, when $\kos{t}^{\pr}=\vphantom{\id}\smash[t]{\id}_{\kos{K}^{\pr}}^{\vphantom{*}}$,
  we write $\id_{\chi}^{\pr}\colon \phi\mor\psi$
  for the Grothendieck $\kos{K}$-morphism
  consisting of the identity Grothendieck morphism $\vphantom{\id}\smash[t]{\id}_{\kos{K}^{\pr}}^{\vphantom{*}}$ of $\kos{K}^{\pr}$
  and the invertible Grothendieck transformation
  $\chi\colon \phi\cong\psi$.
  We illustrate $\id_{\chi}^{\pr}\colon \phi\mor\psi$ and
  its left adjoint part
  $\id_{\chi}^{\pr*}\colon \psi^*\to \phi^*$,
  as well as their descriptions as follows.
  \begin{equation*}
    \begin{tikzcd}[row sep=4pt,column sep=4pt]
      \kos{K}^{\pr}
      \ar[rr,filled,"\id_{\chi}^{\pr}"]
      \ar[ddr,filled,"\phi"']
      &\text{ }
      &\kos{K}^{\pr}
      \ar[ddl,filled,"\psi"]
      &\text{ }
      &\kos{K}^{\pr}
      \ar[ddr,filled,"\phi"']
      \ar[rr,equal]
      &\text{ }
      &\kos{K}^{\pr}
      \ar[ddl,filled,"\psi"]
      \\%1
      \text{ }
      \ar[rr,2cell,shift right=0.5ex,"{\!\!\!\!\what{\id_{\chi}^{\pr}}}"]
      &\text{ }
      &\text{ }
      &=
      &\text{ }
      \ar[rr,2cell,shift right=0ex,"{\!\chi}"]
      &\text{ }
      &\text{ }
      \\%2
      \text{ }
      &\kos{K}
      &\text{ }
      &\text{ }
      &\text{ }
      &\kos{K}
      &\text{ }
    \end{tikzcd}
    \qquad
    \begin{tikzcd}[row sep=4pt,column sep=4pt]
      \kos{K}^{\pr}
      &\text{ }
      &\kos{K}^{\pr}
      \ar[ll,"\id_{\chi}^{\pr*}"']
      &\text{ }
      &\kos{K}^{\pr}
      \ar[rr,equal]
      &\text{ }
      &\kos{K}^{\pr}
      \\%1
      \text{ }
      \ar[rr,2cell,shift right=0.5ex,"{\!\!\!\!\what{\id_{\chi}^{\pr*}}}"]
      &\text{ }
      &\text{ }
      &=
      &\text{ }
      \ar[rr,2cell,shift right=0ex,"{\!\chi}"]
      &\text{ }
      &\text{ }
      \\%2
      \text{ }
      &\kos{K}
      \ar[uul,"\phi^*"]
      \ar[uur,"\psi^*"']
      &\text{ }
      &\text{ }
      &\text{ }
      &\kos{K}
      \ar[uul,"\phi^*"]
      \ar[uur,"\psi^*"']
      &\text{ }
    \end{tikzcd}
  \end{equation*}
\end{itemize}
For each Grothendieck $\kos{K}^{\pr}$-morphism
$\kos{f}^{\pr}\colon \tilde{\kos{t}}^{\pr}\mor\kos{t}^{\pr}$,
we have
$\kos{f}^{\pr}_{\psi}\tilde{\kos{t}}^{\pr}_{\chi}
=\kos{t}^{\pr}_{\chi}\kos{f}^{\pr}_{\phi}
\colon \tilde{\kos{t}}^{\pr}_{\phi}\mor\kos{t}^{\pr}_{\psi}$.

Recall from Section~\ref{sec-ActionbaseGro}
the definitions of a lax $\kos{K}$-tensor functor
and a lax $\kos{K}$-tensor natural transformation.
\begin{itemize}
  \item 
  The left adjoint part of $\phi\colon \phi\mor\vphantom{\id}\smash[t]{\id}_{\kos{K}}^{\vphantom{*}}$
  is the lax $\kos{K}$-tensor functor
  $\phi^*\colon \vphantom{\id}\smash[t]{\id}_{\kos{K}}^{\vphantom{*}}\to \phi^*$
  consisting of $\phi^*\colon \kos{K}\to \kos{K}^{\pr}$
  and $\what{\phi^*}\colon \phi^*=\phi^*$,
  and its $\kos{K}$-equivariance is
  $(\cevar{\phi^*})_{z,x}=(\phi^*_{\tensor})_{z,x}\colon \phi^*(z)\otimes^{\pr}\phi^*(x)\cong \phi^*(z\otimes x)$,
  where $z$, $x\in \CK$.

  \item
  Let $\kos{f}^{\pr}\colon \tilde{\kos{t}}^{\pr}\mor\kos{t}^{\pr}$
  be a Grothendieck $\kos{K}^{\pr}$-morphism.
  We write
  $(\kos{f}^{\pr*})_{\phi}$
  and $(\kos{f}^{\pr}_*)_{\phi}$
  for the left and right adjoint parts
  of the Grothendieck $\kos{K}$-morphism $\kos{f}^{\pr}_{\phi}$,
  both of which are lax $\kos{K}$-tensor functors;
  the notation indicates that they are obtained from 
  $\kos{f}^{\pr*}$ and $\kos{f}^{\pr}_*$
  by change of base along $\phi$.
  The relation
  $\kos{f}^{\pr}_{\psi}\tilde{\kos{t}}^{\pr}_{\chi}
  =\kos{t}^{\pr}_{\chi}\kos{f}^{\pr}_{\phi}$
  above is the equality
  $\tilde{\kos{t}}^{\pr*}_{\chi}(\kos{f}^{\pr*})_{\psi}
  =(\kos{f}^{\pr*})_{\phi}\kos{t}^{\pr*}_{\chi}
  \colon \kos{t}^{\pr*}_{\psi}\to \tilde{\kos{t}}^{\pr*}_{\phi}$
  between the left adjoint parts.

  \item
  Let $\mon{g}^{\pr}\colon \vphantom{\id}\smash[t]{\id}_{\kos{K}^{\pr}}^{\vphantom{*}}\to \vphantom{\id}\smash[t]{\id}_{\kos{K}^{\pr}}^{\vphantom{*}}$ be a lax $\kos{K}^{\pr}$-tensor functor
  consisting of a lax monoidal functor $\mon{g}^{\pr}\colon \kos{K}^{\pr}\to \kos{K}^{\pr}$
  and a monoidal natural transformation $\what{\mon{g}^{\pr}}\colon \vphantom{\id}\smash[t]{\id}_{\kos{K}^{\pr}}^{\vphantom{*}}\Rightarrow \mon{g}^{\pr}$.
  We write $\mon{g}^{\pr}_{\phi}\colon \phi^*\to \phi^*$
  for the lax $\kos{K}$-tensor functor
  whose underlying lax monoidal functor is that of 
  $\mon{g}^{\pr}\colon \vphantom{\id}\smash[t]{\id}_{\kos{K}^{\pr}}^{\vphantom{*}}\to \vphantom{\id}\smash[t]{\id}_{\kos{K}^{\pr}}^{\vphantom{*}}$
  and whose monoidal natural transformation
  $\what{\mon{g}^{\pr}_{\phi}}\colon \phi^*\Rightarrow \mon{g}^{\pr}_{\phi}\phi^*$
  is $\what{\mon{g}^{\pr}}\phi^*\colon \phi^*\Rightarrow\mon{g}^{\pr}\phi^*$.
\end{itemize}
Let $\mon{f}\colon \vphantom{\id}\smash[t]{\id}_{\kos{K}}^{\vphantom{*}}\to \vphantom{\id}\smash[t]{\id}_{\kos{K}}^{\vphantom{*}}$ be a lax $\kos{K}$-tensor functor,
$\mon{g}^{\pr}\colon \vphantom{\id}\smash[t]{\id}_{\kos{K}^{\pr}}^{\vphantom{*}}\to \vphantom{\id}\smash[t]{\id}_{\kos{K}^{\pr}}^{\vphantom{*}}$
be a lax $\kos{K}^{\pr}$-tensor functor,
and $\vartheta\colon \phi^*\mon{f}\Rightarrow\mon{g}^{\pr}_{\phi}\phi^*$
be a lax $\kos{K}$-tensor natural transformation.
By taking
$\monbase{t}=\phi^*$,
$\monbase{t}^{\pr}=\vphantom{\id}\smash[t]{\id}_{\kos{K}}^{\vphantom{*}}$,
$X^{\pr}=\kappa$ in Lemma~\ref{lem-ActionbaseGro-laxKtnatequiv},
we obtain the following commutative diagram for every $z\in \CK$.
\begin{equation}
  \label{eq-changebaseGro-laxKtnattahar}
  \vcenter{\hbox{
    \xymatrix@R=12pt@C=60pt{
      \phi^*(z)\otimes^{\pr}\phi^*\mon{f}(\kappa)
      \ar[r]^-{(\phi^*_{\tensor})_{z,\mon{f}(\kappa)}}_-{\cong}
      \ar[d]_-{I_{\phi^*(z)}\otimes^{\pr} \vartheta_{\kappa}}
      &\phi^*(z\otimes \mon{f}(\kappa))
      \ar[r]^-{\phi^*((\tahar{\mon{f}})_z)}
      &\phi^*\mon{f}(z)
      \ar[d]^-{\vartheta_z}
      \\%1
      \phi^*(z)\otimes^{\pr} \mon{g}^{\pr}\phi^*(\kappa)
      \ar[r]^-{I_{\phi^*(z)}\otimes^{\pr}\mon{g}^{\pr}(\phi^*_{\unit})^{-1}}_-{\cong}
      &\phi^*(z)\otimes^{\pr}\mon{g}^{\pr}(\kappa^{\pr})
      \ar[r]^-{(\tahar{\mon{g}^{\pr}})_{\phi^*(z)}}
      &\mon{g}^{\pr}\phi^*(z)
    }
  }}
\end{equation}
Observe that when $\phi=\vphantom{\id}\smash[t]{\id}_{\kos{K}}^{\vphantom{*}}$,
the commutative diagram \eqref{eq-changebaseGro-laxKtnattahar}
reduces to the commutative diagram \eqref{eq-ActionbaseGro-laxKtnattahar}.

\section{Prekosmic fundamental morphisms}
\label{sec-PKGrofund}

\begin{definition}
  \label{def-PKGrofund-pkfundmor}
  Let $(\kos{K},\kos{t})$, $(\kos{K}^{\pr},\kos{t}^{\pr})$
  be prekosmic Grothendieck categories.
  A \emph{prekosmic fundamental morphism}
  $(\phi,\kos{f})\colon(\kos{K}^{\pr},\kos{t}^{\pr})\mor (\kos{K},\kos{t})$,
  denoted by an arrow with a filled tip,
  is a pair consisting of an \'{e}tale Grothendieck morphism
  $\phi\colon \kos{K}^{\pr}\mor\kos{K}$
  and a Grothendieck $\kos{K}$-morphism
  $\kos{f}\colon \kos{t}^{\pr}_{\phi}\mor \kos{t}$
  such that there exist
  a surjective fiber point $\omega\colon\vphantom{\id}\smash[t]{\id}_{\kos{K}}^{\vphantom{*}}\mor\kos{t}$ of $\kos{t}$,
  a surjective fiber point $\omega^{\pr}\colon\vphantom{\id}\smash[t]{\id}_{\kos{K}^{\pr}}^{\vphantom{*}}\mor\kos{t}^{\pr}$ of $\kos{t}^{\pr}$,
  and an invertible Grothendieck $\kos{K}$-transformation
  $\omega\phi\cong \kos{f}\omega^{\pr}_{\phi}$.
  We illustrate 
  $(\phi,\kos{f})$
  and 
  $\omega\phi\cong \kos{f}\omega^{\pr}_{\phi}$
  as follows.
  \begin{equation*}
    \begin{tikzcd}[row sep=5pt,column sep=20pt]
      \kos{T}^{\pr}
      \ar[dd,filled,"\kos{t}^{\pr}"']
      \ar[r,filled,"\kos{f}"]
      &\kos{T}
      \ar[dd,filled,"\kos{t}"]
      &\text{ }
      &\phi
      \ar[dd,filled,"\omega^{\pr}_{\phi}"']
      \ar[r,filled,"\phi"]
      &\id_{\kos{K}}
      \ar[dd,filled,"\omega"]
      \\%1
      \text{ }
      \ar[r,2cell,shift right=1ex,"{\what{\kos{f}}}"]
      &\text{ }
      &\text{ }
      &\text{ }
      &\text{ }
      \ar[l,2cell,shift left=0.7ex,pos=0.475,"{\,\exists}"']
      \\%2
      \kos{K}^{\pr}
      \ar[r,filled,"{\phi}"]
      &\kos{K}
      &\text{ }
      &\kos{t}^{\pr}_{\phi}
      \ar[r,filled,"{\kos{f}}"]
      &\kos{t}
    \end{tikzcd}
  \end{equation*}
\end{definition}

\begin{remark}
  \label{rem-PKGrofund-fundmordef}
  By taking $\phi=\vphantom{\id}\smash[t]{\id}_{\kos{K}}^{\vphantom{*}}$ in Definition~\ref{def-PKGrofund-pkfundmor},
  we recover the definition of a fundamental $\kos{K}$-morphism,
  given in Definition~\ref{def-fundmoroverGropk-fundmor}.
  Indeed, let $\kos{K}$ be a Grothendieck prekosmos,
  $\kos{t}$, $\kos{t}^{\pr}$ be Grothendieck $\kos{K}$-categories,
  and
  $\kos{f}\colon \kos{t}^{\pr}\mor\kos{t}$
  be a Grothendieck $\kos{K}$-morphism.
  Then $(\vphantom{\id}\smash[t]{\id}_{\kos{K}}^{\vphantom{*}},\kos{f})\colon(\kos{K},\kos{t}^{\pr})\mor (\kos{K},\kos{t})$
  is a prekosmic fundamental morphism
  if and only if
  $\kos{f}\colon \kos{t}^{\pr}\mor\kos{t}$
  is a fundamental $\kos{K}$-morphism.
\end{remark}

\begin{lemma}
  \label{lem-PKGrofund-inducedPKfundmor}
  Let 
  $(\phi,f)\colon (\kos{K}^{\pr},\Spec(\pi^{\pr}))\to (\kos{K},\Spec(\pi))$
  be a morphism of prekosmic Grothendieck objects.
  We have a Grothendieck $\kos{K}$-morphism
  $\kos{f}\colon (\kos{t}^{\pr}_{\pi^{\pr}})_{\phi}\mor\kos{t}_{\pi}^{\vphantom{*}}$
  and an invertible Grothendieck $\kos{K}$-transformation
  $\kos{e}_{\pi}^{\vphantom{*}}\phi\cong \kos{f}(\kos{e}^{\pr}_{\pi^{\pr}})_{\phi}$,
  so that the pair $(\phi,\kos{f})\colon (\kos{K}^{\pr},\kos{t}^{\pr}_{\pi^{\pr}})\mor (\kos{K},\kos{t}_{\pi}^{\vphantom{*}})$
  is a prekosmic fundamental morphism.
  \begin{equation*}
    \begin{tikzcd}[row sep=5pt,column sep=15pt]
      \kosRep_{\kos{K}^{\pr}}(\pi^{\pr})
      \ar[dd,filled,"\kos{t}^{\pr}_{\pi^{\pr}}"']
      \ar[r,filled,"\kos{f}"]
      &\kosRep_{\kos{K}}(\pi)
      \ar[dd,filled,"\kos{t}_{\pi}"]
      &\text{ }
      &\phi
      \ar[dd,filled,"(\kos{e}^{\pr}_{\pi^{\pr}})_{\phi}"']
      \ar[r,filled,"\phi"]
      &\id_{\kos{K}}
      \ar[dd,filled,"\kos{e}_{\pi}"]
      \\%1
      \text{ }
      \ar[r,2cell,shift right=0.6ex,"{\what{\kos{f}}}"]
      &\text{ }
      &\text{ }
      &\text{ }
      &\text{ }
      \ar[l,2cell,shift left=0.3ex,pos=0.475]
      \\%2
      \kos{K}^{\pr}
      \ar[r,filled,"{\phi}"]
      &\kos{K}
      &\text{ }
      &(\kos{t}^{\pr}_{\pi^{\pr}})_{\phi}
      \ar[r,filled,"{\kos{f}}"]
      &\kos{t}_{\pi}
    \end{tikzcd}
  \end{equation*}
\end{lemma}
\begin{proof}
  We write 
  $X=(x,\rho_x)$ for an object in $\cat{Rep}_{\kos{K}}(\pi)$
  and
  $Y^{\pr}=(y^{\pr},\rho^{\pr}_{y^{\pr}})$
  for an object in $\cat{Rep}_{\kos{K}^{\pr}}(\pi^{\pr})$.
  Put
  $\kos{e}\coloneqq \kos{e}_{\pi}^{\vphantom{*}}$
  and
  $\kos{e}^{\pr}\coloneqq \kos{e}^{\pr}_{\pi^{\pr}}$,
  and write 
  $\kos{e}=(\kos{e}^*,\kos{e}_*,\eta,\epsilon)$,
  $\kos{e}^{\pr}=(\kos{e}^{\pr*},\kos{e}^{\pr}_*,\eta^{\pr},\epsilon^{\pr})$,
  and $\phi=(\phi^*,\phi_*,\eta^{\phi},\epsilon^{\phi})$.
  The left adjoint part
  $\kos{f}^*\colon \kosRep_{\kos{K}}(\pi)\to \kosRep_{\kos{K}^{\pr}}(\pi^{\pr})$
  is a strong monoidal functor,
  which sends each $X$ to the pair
  $\kos{f}^*(X)=(\phi^*(x),(I_{\phi^*(x)}\otimes^{\pr} f)\circ (\phi^*_{\tensor})^{-1}_{x,\pi}\circ \phi^*(\rho_x))$,
  and whose coherence isomorphisms are
  given by those of $\phi^*\colon \kos{K}\to \kos{K}^{\pr}$
  so that we have an equality
  $\phi^*\kos{e}^*=\kos{e}^{\pr*}\kos{f}^*$
  of strong monoidal functors.
  We write
  $\zeta\colon \phi_*(\slot)\otimes \pi\cong \phi_*(\slot\otimes^{\pr}\phi^*(\pi))$
  for the natural isomorphism
  whose component at $z^{\pr}\in \CK^{\pr}$ is
  $(\phi_{*\tensor})_{z^{\pr},\phi^*(\pi)}\circ (I_{\phi_*(z^{\pr})}\otimes \eta^{\phi}_{\pi})$;
  observe that $\zeta$ agrees with $\varphi^{\phi}_{\pi,\slot}
  \colon \pi\otimes \phi_*(\slot)\cong \phi_*(\phi^*(\pi)\otimes^{\pr}\slot)$
  up to symmetric coherence isomorphisms,
  where $\varphi^{\phi}$
  is the natural isomorphism associated to $\phi$
  in \eqref{eq-Groprekos-projformula}.
  For each $Y^{\pr}$,
  we fix an equalizer $\kos{f}_*(Y^{\pr})$ in $\cat{Rep}_{\kos{K}}(\pi)$
  of the coreflexive pair
  $\phi_*(\rho^{\pr}_{y^{\pr}})\otimes I_{\pi}$,
  $(\phi_*(I_{y^{\pr}}\otimes^{\pr}f)\otimes I_{\pi})\circ 
  (\zeta_{y^{\pr}}\otimes I_{\pi})\circ (I_{\phi_*(y^{\pr})}\otimes \cp_{\pi})
  \colon
  (\phi_*(y^{\pr})\otimes \pi,I_{\phi_*(y^{\pr})}\otimes \cp_{\pi})
  \to
  (\phi_*(y^{\pr}\otimes^{\pr}\pi^{\pr})\otimes \pi,I_{\phi_*(y^{\pr}\otimes^{\pr}\pi^{\pr})}\otimes \cp_{\pi})$,
  with common retraction
  $\phi_*(I_{y^{\pr}}\otimes^{\pr}e^{\pr}_{\pi^{\pr}})\otimes I_{\pi}$.
  Then $\kos{f}_*$ is right adjoint to $\kos{f}^*$
  by the dual statement of the adjoint lifting theorem
  \cite[Theorem~4.5.6 and Exercise~4.8.6]{Borceux1994}.
  Hence $\kos{f}=(\kos{f}^*\dashv\kos{f}_*)$
  is a Grothendieck morphism from
  $\kosRep_{\kos{K}^{\pr}}(\pi^{\pr})$ to $\kosRep_{\kos{K}}(\pi)$.
  The invertible Grothendieck transformation
  $\what{\kos{f}}\colon \phi\kos{t}^{\pr}_{\pi^{\pr}}\cong \kos{t}_{\pi}^{\vphantom{*}}\kos{f}$
  is the equality
  $\kos{t}^{\pr*}_{\pi^{\pr}}\phi^*=\kos{f}^*\kos{t}_{\pi}^*$
  between the left adjoint parts,
  and the invertible Grothendieck $\kos{K}$-transformation
  $\kos{e}\phi\cong \kos{f}\kos{e}^{\pr}_{\phi}$
  is the equality
  $\phi^*\kos{e}^*=(\kos{e}^{\pr*})_{\phi}\kos{f}^*$
  between the left adjoint parts
  as lax $\kos{K}$-tensor functors.
\end{proof}

\begin{proposition}
  \label{prop-PKGrofund-inducedPKfundmor}
  Let 
  $(\phi,\kos{f})\colon (\kos{K}^{\pr},\kos{t}^{\pr})\mor(\kos{K},\kos{t})$
  be a prekosmic fundamental morphism.
  Choose a surjective fiber point $\omega$ of $\kos{t}$,
  a surjective fiber point $\omega^{\pr}$ of $\kos{t}^{\pr}$,
  and an invertible Grothendieck $\kos{K}$-transformation
  $\alpha\colon \omega\phi\cong \kos{f}\omega^{\pr}_{\phi}$.
  Put $\pi\coloneqq \omega^*\omega_*(\kappa)$
  and
  $\pi^{\pr}\coloneqq \omega^{\pr*}\omega^{\pr}_*(\kappa^{\pr})$,
  and recall the equivalences
  $\breve{\omega}\colon \kos{t}_{\pi}^{\vphantom{*}}\simeq\kos{t}$,
  $\breve{\omega}^{\pr}\colon \kos{t}^{\pr}_{\pi^{\pr}}\simeq\kos{t}^{\pr}$
  in Theorem~\ref{thm-Grocatoverpk-mainthm}.
  Then we have a morphism
  $(\phi,\breve{f})\colon (\kos{K}^{\pr},\Spec(\pi^{\pr}))\to (\kos{K},\Spec(\pi))$
  of prekosmic Grothendieck objects
  and an invertible Grothendieck $\kos{K}$-transformation
  $\breve{\alpha}\colon \breve{\omega}\breve{\kos{f}}\cong \kos{f}\breve{\omega}^{\pr}_{\phi}$
  such that the following equality of pasting diagrams holds,
  where $\breve{\kos{f}}\colon (\kos{t}^{\pr}_{\pi^{\pr}})_{\phi}\mor\kos{t}_{\pi}^{\vphantom{*}}$
  is the Grothendieck $\kos{K}$-morphism
  in Lemma~\ref{lem-PKGrofund-inducedPKfundmor}
  associated to $(\phi,\breve{f})$.
  \begin{equation*}
    \begin{tikzcd}[row sep=4pt,column sep=8pt]
      \phi
      \ar[rr,filled,"\phi"]
      \ar[dd,filled,"{(\kos{e}^{\pr}_{\pi^{\pr}})_{\phi}}"']
      &\text{ }
      &\id_{\kos{K}}
      \ar[dd,filled,"{\kos{e}_{\pi}}"]
      \ar[ddddrr,bend left=2.5pc,filled,"{\omega}"]
      &\text{ }
      &\text{ }
      &\text{ }%center column
      &\phi
      \ar[rr,filled,"\phi"]
      \ar[dd,filled,"{(\kos{e}^{\pr}_{\pi^{\pr}})_{\phi}}"']
      \ar[ddddrr,filled,bend left=2.5pc,"{\omega^{\pr}_{\phi}}" description]
      &\text{ }
      &\id_{\kos{K}}
      \ar[ddddrr,bend left=2.5pc,filled,"{\omega}"]
      &\text{ }
      &\text{ }
      \\%1
      \text{ }
      &\text{ }
      &\text{ }
      \ar[ll,2cell,shift left=0.3ex,pos=0.45]
      &\text{ }
      &\text{ }
      &\text{ }%center column
      &\text{ }
      &\text{ }
      &\text{ }
      &\text{ }
      &\text{ }
      \\%2
      (\kos{t}^{\pr}_{\pi^{\pr}})_{\phi}
      \ar[rr,filled,"\breve{\kos{f}}"]
      \ar[ddrr,filled,"\breve{\omega}^{\pr}_{\phi}"',"\simeq"]
      &\text{ }
      &\kos{t}_{\pi}
      \ar[ddrr,filled,"\breve{\omega}","\simeq"']
      &\text{ }
      &\text{ }
      \ar[ll,2cell,shift right=0ex,pos=0.6]
      &=%center column
      &(\kos{t}^{\pr}_{\pi^{\pr}})_{\phi}
      \ar[ddrr,filled,"\breve{\omega}^{\pr}_{\phi}"',"\simeq"]
      &\text{ }
      &\text{ }
      \ar[ll,2cell,shift right=0ex,pos=0.65]
      &\text{ }
      &\text{ }
      \ar[ll,2cell,shift right=0ex,pos=0.75,"{\,\alpha}"']
      \\%3
      \text{ }
      &\text{ }
      &\text{ }
      &\text{ }
      \ar[ll,2cell,shift left=1ex,"{\,\breve{\alpha}}"']
      &\text{ }
      &\text{ }%center column
      &\text{ }
      &\text{ }
      &\text{ }
      &\text{ }
      &\text{ }
      \\%4
      \text{ }
      &\text{ }
      &\kos{t}^{\pr}_{\phi}
      \ar[rr,filled,"\kos{f}"']
      &\text{ }
      &\kos{t}
      &\text{ }%center column
      &\text{ }
      &\text{ }
      &\kos{t}^{\pr}_{\phi}
      \ar[rr,filled,"\kos{f}"']
      &\text{ }
      &\kos{t}
    \end{tikzcd}
  \end{equation*}
\end{proposition}
\begin{proof}
  We write $\alpha_*\colon \kos{f}^*\omega_*\Rightarrow (\omega^{\pr}_*)_{\phi}\phi^*$
  for the mate of $\alpha^{-1}\colon (\omega^{\pr*})_{\phi}\kos{f}^*\cong \phi^*\omega^*$
  with respect to the adjunctions
  $\omega=(\omega^*\dashv\omega_*)$ and
  $\omega^{\pr}_{\phi}=((\omega^{\pr*})_{\phi}\dashv(\omega^{\pr}_*)_{\phi})$;
  it is the unique lax $\kos{K}$-tensor natural transformation
  that satisfies the following relations.
  \begin{equation}\label{eq-PKGrofund-matealpha}
    \vcenter{\hbox{
      \xymatrix@R=12pt@C=30pt{
        \kos{f}^*
        \ar@2{->}[r]^-{\eta^{\pr}\kos{f}^*}
        \ar@2{->}[d]_-{\kos{f}^*\eta}
        &(\omega^{\pr}_*)_{\phi}(\omega^{\pr*})_{\phi}\kos{f}^*
        \ar@2{->}[d]^-{(\omega^{\pr}_*)_{\phi}\alpha^{-1}}_-{\cong}
        &(\omega^{\pr*})_{\phi}\kos{f}^*\omega_*
        \ar@2{->}[r]^-{(\omega^{\pr*})_{\phi}\alpha_*}
        \ar@2{->}[d]_-{\alpha^{-1}\omega_*}^-{\cong}
        &(\omega^{\pr*})_{\phi}(\omega^{\pr}_*)_{\phi}\phi^*
        \ar@2{->}[d]^-{\epsilon^{\pr}\phi^*}
        \\%1
        \kos{f}^*\omega_*\omega^*
        \ar@2{->}[r]^-{\alpha_*\omega^*}
        &(\omega^{\pr}_*)_{\phi}\phi^*\omega^*
        &\phi^*\omega^*\omega_*
        \ar@2{->}[r]^-{\phi^*\epsilon}
        &\phi^*
      }
    }}
  \end{equation}
  Put
  $\alpha^\sharp\coloneqq ((\omega^{\pr*})_{\phi}\alpha_*)\circ (\alpha\omega_*)
  \colon \phi^*\omega^*\omega_*\Rightarrow (\omega^{\pr*})_{\phi}(\omega^{\pr}_*)_{\phi}\phi^*$
  and define 
  $\breve{f}\colon \Spec(\pi^{\pr})\to \Spec(\phi^*(\pi))$
  as
  $\omega^{\pr*}\omega^{\pr}_*(\phi^*_{\unit})^{-1}\circ \alpha^{\sharp}_{\kappa}
  \colon \phi^*\omega^*\omega_*(\kappa)
  \to
  \omega^{\pr*}\omega^{\pr}_*(\kappa^{\pr})$.
  Put
  $\varpi\coloneqq\omega^*\omega_*$,
  $\varpi^{\pr}\coloneqq\omega^{\pr*}\omega^{\pr}_*$,
  $\delta\coloneqq\omega^*\eta\omega_*\colon \varpi\Rightarrow\varpi\varpi$,
  and
  $\delta^{\pr}\coloneqq\omega^{\pr*}\eta^{\pr}\omega^{\pr}_*
  \colon \varpi^{\pr}\Rightarrow\varpi^{\pr}\varpi^{\pr}$.
  Observe that
  $\varpi^{\pr}_{\phi}=(\omega^{\pr*})_{\phi}(\omega^{\pr}_*)_{\phi}$
  as lax $\kos{K}$-tensor functors.
  By substituting 
  $\mon{f}=\varpi$, $\mon{g}^{\pr}=\varpi^{\pr}$,
  $\vartheta=\alpha^{\sharp}$
  in \eqref{eq-changebaseGro-laxKtnattahar},
  we obtain the following commutative diagram for every $z\in \CK$.
  \begin{equation}
    \label{eq-PKGrofund-alphasharpequiv}
    \vcenter{\hbox{
      \xymatrix@R=15pt@C=50pt{
        \phi^*(z)\otimes^{\pr}\phi^*\omega^*\omega_*(\kappa)
        \ar[r]^-{(\phi^*_{\tensor})_{z,\omega^*\omega_*(\kappa)}}_-{\cong}
        \ar[d]_-{I_{\phi^*(z)}\otimes^{\pr} \breve{f}}
        &\phi^*(z\otimes \omega^*\omega_*(\kappa))
        \ar[r]^-{\phi^*((\tahar{\omega^*\omega_*})_z)}_-{\cong}
        &\phi^*\omega^*\omega_*(z)
        \ar[d]^-{\alpha^{\sharp}_z}
        \\%1
        \phi^*(z)\otimes^{\pr}\omega^{\pr*}\omega^{\pr}_*(\kappa^{\pr})
        \ar[rr]^-{(\tahar{\omega^{\pr*}\omega^{\pr}_*})_{\phi^*(z)}}_-{\cong}
        &\text{ }
        &\omega^{\pr*}\omega^{\pr}_*\phi^*(z)
      }
    }}
  \end{equation}
  One can check that
  $(\varpi^{\pr}\alpha^{\sharp})\circ
  (\alpha^{\sharp}\varpi)\circ (\phi^*\delta)
  =(\delta^{\pr}\phi^*)\circ \alpha^{\sharp}
  \colon \phi^*\varpi\Rightarrow \varpi^{\pr}\varpi^{\pr}\phi^*$
  and that $\phi^*\epsilon=(\epsilon^{\pr}\phi^*)\circ \alpha^{\sharp}\colon \phi^*\varpi\Rightarrow\phi^*$
  as monoidal natural transformations between lax monoidal functors,
  using \eqref{eq-PKGrofund-matealpha}.
  We also have
  $(I_{\phi^*(\pi)}\otimes^{\pr}\breve{f})\circ \cp^{\pr}_{\phi^*(\pi)}
  =(\tahar{\varpi^{\pr}})_{\phi^*(\pi)}^{-1}\circ \alpha^{\sharp}_{\pi}\circ \phi^*(\delta_{\kappa})
  \colon \phi^*(\pi)\to \phi^*(\pi)\otimes^{\pr}\pi^{\pr}$,
  which follows from the relation
  \eqref{eq-PKGrofund-alphasharpequiv} for $z=\pi$.
  From these equalities, we deduce that
  $\breve{f}\colon \Spec(\pi^{\pr})\to \Spec(\phi^*(\pi))$
  is a morphism of Grothendieck objects in $\kos{K}^{\pr}$.
  It remains to show the existence of $\breve{\alpha}$.
  It suffices to show that for each $X\in \CT$,
  the isomorphism
  $\alpha_X\colon \phi^*\omega^*(X)\cong \omega^{\pr*}\kos{f}^*(X)$
  in $\CK^{\pr}$ lifts to an isomorphism
  $\breve{\alpha}_X\colon \breve{\kos{f}}^*\breve{\omega}^*(X)\cong \breve{\omega}^{\pr*}\kos{f}^*(X)$
  in $\cat{Rep}_{\kos{K}^{\pr}}(\pi^{\pr})$.
  One can check this
  using \eqref{eq-PKGrofund-alphasharpequiv} for $z=\omega^*(X)$
  and \eqref{eq-PKGrofund-matealpha}.
\end{proof}

\section{The duality over Grothendieck prekosmoi}
\label{sec-dualityPKGro}

In this section, we construct the $(2,1)$-category
$\mathsf{pGroCat}_{\mathsf{pt}}$
and prove Theorem~\ref{thm-PKGroCat-equiv21catoveret}.

\begin{definition}
  \label{def-dualityPKGro-pGroCatpt21cat}
  We define 
  the $(2,1)$-category
  $\mathsf{pGroCat}_{\mathsf{pt}}$
  as follows.
  \begin{itemize}
    \item 
    A $0$-cell in $\mathsf{pGroCat}_{\mathsf{pt}}$
    is a pointed prekosmic Grothendieck category.

    \item
    A $1$-cell in $\mathsf{pGroCat}_{\mathsf{pt}}$
    is a \emph{pointed prekosmic fundamental morphism}.
    Let 
    $(\kos{K},\kos{t},\omega)$,
    $(\kos{K}^{\pr},\kos{t}^{\pr},\omega^{\pr})$
    be $0$-cells.
    A pointed prekosmic fundamental morphism
    $(\phi,\kos{f},\alpha)\colon (\kos{K}^{\pr},\kos{t}^{\pr},\omega^{\pr})\mor(\kos{K},\kos{t},\omega)$,
    denoted by an arrow with a filled tip,
    is a tuple consisting of
    an \'{e}tale Grothendieck morphism $\phi\colon \kos{K}^{\pr}\mor\kos{K}$,
    a Grothendieck $\kos{K}$-morphism
    $\kos{f}\colon \kos{t}^{\pr}_{\phi}\mor\kos{t}$,
    and an invertible Grothendieck $\kos{K}$-transformation
    $\alpha\colon \omega\phi\cong \kos{f}\omega^{\pr}_\phi$.
    The pair $(\phi,\kos{f})$
    is a prekosmic fundamental morphism
    from $(\kos{K}^{\pr},\kos{t}^{\pr})$ to $(\kos{K},\kos{t})$.
    We illustrate $(\phi,\kos{f},\alpha)$ as follows.
    \begin{equation*}
      \begin{tikzcd}[row sep=5pt,column sep=20pt]
        \kos{T}^{\pr}
        \ar[dd,filled,"\kos{t}^{\pr}"']
        \ar[r,filled,"\kos{f}"]
        &\kos{T}
        \ar[dd,filled,"\kos{t}"]
        &\text{ }%empty
        &\phi
        \ar[dd,filled,"\omega^{\pr}_{\phi}"']
        \ar[r,filled,"\phi"]
        &\id_{\kos{K}}
        \ar[dd,filled,"\omega"]
        \\%1
        \text{ }
        \ar[r,2cell,shift right=1ex,"{\what{\kos{f}}}"]
        &\text{ }
        &\text{ }%empty
        &\text{ }
        &\text{ }
        \ar[l,2cell,shift left=0.7ex,pos=0.475,"{\alpha}"']
        \\%2
        \kos{K}^{\pr}
        \ar[r,filled,"{\phi}"]
        &\kos{K}
        &\text{ }%empty
        &\kos{t}^{\pr}_{\phi}
        \ar[r,filled,"{\kos{f}}"]
        &\kos{t}
      \end{tikzcd}
    \end{equation*}
    
    \item
    A $2$-cell
    $(\chi,\vartheta)\colon (\phi,\kos{f},\alpha)\cong (\psi,\kos{g},\beta)$
    in $\mathsf{pGroCat}_{\mathsf{pt}}$
    between $1$-cells
    $(\phi,\kos{f},\alpha)$, $(\psi,\kos{g},\beta)
    \colon (\kos{K}^{\pr},\kos{t}^{\pr},\omega^{\pr})\mor (\kos{K},\kos{t},\omega)$
    is a pair consisting of an invertible Grothendieck transformation
    $\chi\colon \phi\cong\psi$
    and an invertible Grothendieck $\kos{K}$-transformation
    $\vartheta\colon \kos{f}\cong \kos{g}\kos{t}^{\pr}_{\chi}$
    between Grothendieck $\kos{K}$-morphisms
    $\kos{f}$, $\kos{g}\kos{t}^{\pr}_{\chi}\colon \kos{t}^{\pr}_{\phi}\mor\kos{t}$.
  \end{itemize}
  We define the rest of the structure of $\mathsf{pGroCat}_{\mathsf{pt}}$ as follows.
  Let 
  $(\chi,\vartheta)\colon (\phi,\kos{f},\alpha)\cong (\psi,\kos{g},\beta)$
  be a $2$-cell
  between $1$-cells
  $(\phi,\kos{f},\alpha)$, $(\psi,\kos{g},\beta)
  \colon (\kos{K}^{\pr},\kos{t}^{\pr},\omega^{\pr})\mor (\kos{K},\kos{t},\omega)$.

  Let
  $(\tilde{\chi},\tilde{\vartheta})\colon (\psi,\kos{g},\beta)\cong (\tilde{\psi},\tilde{\kos{g}},\tilde{\beta})$
  be a $2$-cell.
  The vertical composite of
  $(\chi,\vartheta)$ and $(\tilde{\chi},\tilde{\vartheta})$
  is the pair
  $(\tilde{\chi}\circ\chi,(\tilde{\vartheta}\kos{t}^{\pr}_{\chi})
  \circ \vartheta)
  \colon (\phi,\kos{f},\alpha)\cong (\tilde{\psi},\tilde{\kos{g}},\tilde{\beta})$
  consisting of 
  $\tilde{\chi}\circ\chi\colon \phi\cong \tilde{\psi}$
  and
  $(\tilde{\vartheta}\kos{t}^{\pr}_{\chi})
  \circ \vartheta
  \colon \kos{f}\cong \tilde{\kos{g}}\kos{t}^{\pr}_{\tilde{\chi}}\kos{t}^{\pr}_{\chi}$.
  The identity $2$-cell of $(\phi,\kos{f},\alpha)$
  is the pair consisting of the identity Grothendieck transformation of $\phi$
  and the identity Grothendieck $\kos{K}$-transformation of $\kos{f}$.
  Each $2$-cell
  $(\chi,\vartheta)\colon (\phi,\kos{f},\alpha)\cong (\psi,\kos{g},\beta)$
  is indeed invertible with respect to the vertical composition;
  its inverse $(\chi^{-1},\vartheta^{-1}\kos{t}^{\pr}_{\chi^{-1}})\colon(\psi,\kos{g},\beta)\cong (\phi,\kos{f},\alpha)$
  is the pair consisting of 
  $\chi^{-1}\colon \psi\cong \phi$
  and 
  $\vartheta^{-1}\kos{t}^{\pr}_{\chi^{-1}}\colon \kos{g}\cong \kos{f}\kos{t}^{\pr}_{\chi^{-1}}$.

  Let 
  $(\phi^{\pr},\kos{f}^{\pr},\alpha^{\pr})$,
  $(\psi^{\pr},\kos{g}^{\pr},\beta^{\pr})
  \colon (\kos{K}^{\ppr},\kos{t}^{\ppr},\omega^{\ppr})\mor(\kos{K}^{\pr},\kos{t}^{\pr},\omega^{\pr})$
  be $1$-cells.
  The composite of
  $(\phi,\kos{f},\alpha)$
  and
  $(\phi^{\pr},\kos{f}^{\pr},\alpha^{\pr})$
  is the tuple 
  $(\phi^{\ppr},\kos{f}^{\ppr},\alpha^{\ppr})\colon (\kos{K}^{\ppr},\kos{t}^{\ppr},\omega^{\ppr})\mor (\kos{K},\kos{t},\omega)$
  consisting of 
  $\phi^{\ppr}\coloneqq \phi\phi^{\pr}\colon \kos{K}^{\ppr}\mor\kos{K}$,
  $\kos{f}^{\ppr}\coloneqq \kos{f}\kos{f}^{\pr}_{\phi}\colon \kos{t}^{\ppr}_{\phi^{\ppr}}\mor\kos{t}$,
  and 
  $\alpha^{\ppr}\coloneqq (\kos{f}\alpha^{\pr}_{\phi})\circ (\alpha\phi^{\pr}_{\phi})
  \colon \omega\phi^{\ppr}\cong \kos{f}^{\ppr}\omega^{\ppr}_{\phi^{\ppr}}$.
  We illustrate the composite 
  and its description as follows.
  \begin{equation*}
    \begin{tikzcd}[row sep=5pt,column sep=20pt]
      \phi^{\ppr}
      \ar[dd,filled,"\omega^{\ppr}_{\phi^{\ppr}}"']
      \ar[r,filled,"\phi^{\ppr}"]
      &\id_{\kos{K}}
      \ar[dd,filled,"\omega"]
      &\text{ }
      &\phi^{\pr}_{\phi}
      \ar[dd,filled,"(\omega^{\ppr}_{\phi^{\pr}})_{\phi}"']
      \ar[r,filled,"\phi^{\pr}_{\phi}"]
      &\phi
      \ar[dd,filled,"\omega^{\pr}_{\phi}"description]
      \ar[r,filled,"\phi"]
      &\id_{\kos{K}}
      \ar[dd,filled,"\omega"]
      \\%1
      \text{ }
      &\text{ }
      \ar[l,2cell,shift left=0ex,pos=0.475,"{\,\,\alpha^{\ppr}}"']
      &=
      &\text{ }
      &\text{ }
      \ar[l,2cell,shift left=0ex,pos=0.475,"{\alpha^{\pr}_{\phi}}"']
      &\text{ }
      \ar[l,2cell,shift left=0ex,pos=0.475,"{\alpha}"']
      \\%2
      \kos{t}^{\ppr}_{\phi^{\ppr}}
      \ar[r,filled,"{\kos{f}^{\ppr}}"]
      &\kos{t}
      &\text{ }
      &(\kos{t}^{\ppr}_{\phi^{\pr}})_{\phi}
      \ar[r,filled,"{\kos{f}^{\pr}_{\phi}}"]
      &\kos{t}^{\pr}_{\phi}
      \ar[r,filled,"{\kos{f}}"]
      &\kos{t}
    \end{tikzcd}
  \end{equation*}
  The identity $1$-cell of $(\kos{K},\kos{t},\omega)$ is 
  the tuple consisting of the identity Grothendieck morphism $\vphantom{\id}\smash[t]{\id}_{\kos{K}}^{\vphantom{*}}$ of $\kos{K}$,
  the identity Grothendieck $\kos{K}$-morphism of $\kos{t}$,
  and the identity Grothendieck $\kos{K}$-transformation of $\omega\colon \vphantom{\id}\smash[t]{\id}_{\kos{K}}^{\vphantom{*}}\mor\kos{t}$.
  Let
  $(\chi^{\pr},\vartheta^{\pr})\colon (\phi^{\pr},\kos{f}^{\pr},\alpha^{\pr})\cong (\psi^{\pr},\kos{g}^{\pr},\beta^{\pr})$
  be a $2$-cell.
  We write $(\psi^{\ppr},\kos{g}^{\ppr},\beta^{\ppr})\colon (\kos{K}^{\ppr},\kos{t}^{\ppr},\omega^{\ppr})\mor (\kos{K},\kos{t},\omega)$
  for the composite of
  $(\psi,\kos{g},\beta)$
  and
  $(\psi^{\pr},\kos{g}^{\pr},\beta^{\pr})$.
  The horizontal composite of
  $(\chi,\vartheta)$ and $(\chi^{\pr},\vartheta^{\pr})$
  is the pair
  $(\chi^{\ppr},\vartheta^{\ppr})
  \colon (\phi^{\ppr},\kos{f}^{\ppr},\alpha^{\ppr})\cong (\psi^{\ppr},\kos{g}^{\ppr},\beta^{\ppr})$
  consisting of
  $\chi^{\ppr}\colon \phi^{\ppr}\cong \psi^{\ppr}$,
  given by
  $\chi\chi^{\pr}\colon \phi\phi^{\pr}\cong \psi\psi^{\pr}$,
  and 
  $\vartheta^{\ppr}\colon \kos{f}^{\ppr}\cong \kos{g}^{\ppr}\kos{t}^{\ppr}_{\chi^{\ppr}}$,
  given by
  $\vartheta\vartheta^{\pr}_{\phi}
  \colon \kos{f}\kos{f}^{\pr}_{\phi}
  \cong \kos{g}\kos{t}^{\pr}_{\chi}\kos{g}^{\pr}_{\phi}(\kos{t}^{\ppr}_{\chi^{\pr}})_{\phi}
  =\kos{g}\kos{g}^{\pr}_{\psi}(\kos{t}^{\ppr}_{\psi^{\pr}})_{\chi}(\kos{t}^{\ppr}_{\chi^{\pr}})_{\phi}
  =\kos{g}\kos{g}^{\pr}_{\psi}\kos{t}^{\ppr}_{\chi\chi^{\pr}}$.
\end{definition}

We denote by $\mathsf{p}_{\mathsf{ptcat}}
\colon \mathsf{pGroCat}_{\mathsf{pt}}\to\mathsf{pGroKos}_{\etsf}$
the $2$-functor that sends each cell to its first component.
The strict fiber of $\mathsf{p}_{\mathsf{ptcat}}$
over each Grothendieck prekosmos $\kos{K}$
is $2$-isomorphic to $\mathsf{pGroCat}_{\mathsf{pt}}(\kos{K})$,
defined in Definition~\ref{def-dualityoverGropk-pGroCatptK}.
We illustrate $\mathsf{p}_{\mathsf{ptcat}}$ as follows.
\begin{equation*}
  \begin{tikzcd}[row sep=22pt,column sep=20pt]
    \mathsf{pGroCat}_{\mathsf{pt}}
    \ar[d,"\mathsf{p}_{\mathsf{ptcat}}"']
    &(\kos{K}^{\pr},\kos{t}^{\pr},\omega^{\pr})
    \ar[rr,filled,shift left=0.3ex,bend left=1.2pc,"{(\psi,\kos{g},\beta)}" description]
    \ar[rr,filled,shift right=0.3ex,bend right=1.2pc,"{(\phi,\kos{f},\alpha)}" description]
    &\text{ }
    &(\kos{K},\kos{t},\omega)
    \\%1
    {\mathsf{pGroKos}_{\etsf}}
    &\kos{K}^{\pr}
    \ar[rr,filled,shift left=0.3ex,bend left=1.2pc,"\psi" description]
    \ar[rr,filled,shift right=0.3ex,bend right=1.2pc,"\phi" description]
    &\text{ }
    \ar[u,2cell,pos=1.17,"{{(\chi,\vartheta)}}"']
    \ar[u,2cell,pos=-0.05,"{\,\chi}"']
    &\kos{K}
  \end{tikzcd}
\end{equation*}

\begin{proposition}
  \label{prop-dualityPKGro-Rep}
  We have a pseudofunctor
  \begin{equation*}
    \mathsf{Rep}\colon \mathsf{pGroObj}\to \mathsf{pGroCat}_{\mathsf{pt}}
  \end{equation*}
  over $\mathsf{pGroKos}_{\etsf}$
  that sends each prekosmic Grothendieck object $(\kos{K},\Spec(\pi))$
  to the pointed prekosmic Grothendieck category
  $(\kos{K},\kos{t}_{\pi}^{\vphantom{*}},\kos{e}_{\pi}^{\vphantom{*}})$.
\end{proposition}
\begin{proof}
  We describe how $\mathsf{Rep}$ sends
  $1$-cells and $2$-cells.
  \begin{itemize}
    \item
    A $1$-cell $(\phi,f)\colon (\kos{K}^{\pr},\Spec(\pi^{\pr}))\to(\kos{K},\Spec(\pi))$
    in $\mathsf{pGroObj}$
    is sent to the $1$-cell 
    $(\kos{K}^{\pr},\kos{t}^{\pr}_{\pi^{\pr}},\kos{e}^{\pr}_{\pi^{\pr}})
    \mor(\kos{K},\kos{t}_{\pi}^{\vphantom{*}},\kos{e}_{\pi}^{\vphantom{*}})$
    in $\mathsf{pGroCat}_{\mathsf{pt}}$
    consisting of 
    $\phi\colon \kos{K}^{\pr}\mor\kos{K}$,
    the Grothendieck $\kos{K}$-morphism
    $\kos{f}\colon (\kos{t}^{\pr}_{\pi^{\pr}})_{\phi}\mor \kos{t}_{\pi}^{\vphantom{*}}$,
    and the invertible Grothendieck $\kos{K}$-transformation
    $\kos{e}_{\pi}^{\vphantom{*}}\phi\cong \kos{f}(\kos{e}^{\pr}_{\pi^{\pr}})_{\phi}$,
    described in the proof of Lemma~\ref{lem-PKGrofund-inducedPKfundmor}.

    \item
    Let
    $(\phi,f)$, $(\psi,g)\colon (\kos{K}^{\pr},\Spec(\pi^{\pr}))\to(\kos{K},\Spec(\pi))$
    be $1$-cells
    and write
    $\kos{f}\colon (\kos{t}^{\pr}_{\pi^{\pr}})_{\phi}\mor \kos{t}_{\pi}^{\vphantom{*}}$,
    $\kos{g}\colon (\kos{t}^{\pr}_{\pi^{\pr}})_{\psi}\mor \kos{t}_{\pi}^{\vphantom{*}}$
    for the associated Grothendieck $\kos{K}$-morphisms.
    A $2$-cell $(\chi,\theta)\colon (\phi,f)\cong (\psi,g)$
    in $\mathsf{pGroObj}$
    is sent to the $2$-cell
    in $\mathsf{pGroCat}_{\mathsf{pt}}$
    consisting of $\chi\colon \phi\cong \psi$
    and the invertible Grothendieck $\kos{K}$-transformation
    $\vartheta\colon \kos{f}\cong \kos{g}(\kos{t}^{\pr}_{\pi^{\pr}})_{\chi}$,
    defined as follows.
    The component $\vartheta_X\colon \kos{f}^*(X)\cong \kos{g}^*(X)$
    of $\vartheta$ at $X=(x,\rho_x)\in\cat{Rep}_{\kos{K}}(\pi)$
    is the isomorphism in $\cat{Rep}_{\kos{K}^{\pr}}(\pi^{\pr})$
    whose underlying isomorphism in $\CK^{\pr}$ is
    $\chi_x\circ (I_{\phi^*(x)}\otimes^{\pr} \theta)\circ (\phi^*_{\tensor})_{x,\pi}^{-1}\circ \phi^*(\rho_x)\colon \phi^*(x)\cong \psi^*(x)$.
    It is straightforward to check that
    $\vartheta$ is indeed an invertible Grothendieck $\kos{K}$-transformation
    from $\kos{f}$ to $\kos{g}(\kos{t}^{\pr}_{\pi^{\pr}})_{\chi}$.
  \end{itemize}
  It is straightforward to check that $\mathsf{Rep}$
  preserves vertical composition of $2$-cells and identity $2$-cells;
  thus $\mathsf{Rep}$ defines a functor on each hom-groupoid.
  The constraints of the pseudofunctor $\mathsf{Rep}$
  are defined as follows.
  \begin{itemize}
    \item
    Let 
    $(\phi,f)$,
    $(\phi^{\pr},f^{\pr})$
    be a composable pair of $1$-cells in $\mathsf{pGroObj}$,
    with composite $(\phi^{\ppr},f^{\ppr})\coloneqq (\phi\phi^{\pr},\phi^{\pr*}(f)f^{\pr})$,
    where $\phi\colon \kos{K}^{\pr}\mor\kos{K}$
    and
    $\phi^{\pr}\colon\kos{K}^{\ppr}\mor\kos{K}^{\pr}$.
    We write 
    $\kos{f}$, $\kos{f}^{\pr}$, and $\kos{f}^{\ppr}$
    for the associated Grothendieck morphisms
    over $\kos{K}$, $\kos{K}^{\pr}$, and $\kos{K}$, respectively.
    The component of the functorial constraint of $\mathsf{Rep}$
    at $(\phi,f)$, $(\phi^{\pr},f^{\pr})$
    is the invertible Grothendieck $\kos{K}$-transformation
    $\kos{f}\kos{f}^{\pr}_{\phi}\cong\kos{f}^{\ppr}$
    given by the equality
    $(\kos{f}^{\pr*})_{\phi}\kos{f}^*=\kos{f}^{\ppr*}$
    between the left adjoint parts.

    \item
    Let $(\kos{K},\Spec(\pi))$ be a $0$-cell in $\mathsf{pGroObj}$.
    We write 
    $\kos{I}_{\pi}^{\vphantom{*}}\colon \kos{t}_{\pi}^{\vphantom{*}}\mor\kos{t}_{\pi}^{\vphantom{*}}$
    for the Grothendieck $\kos{K}$-morphism
    associated to its identity $1$-cell
    in Lemma~\ref{lem-PKGrofund-inducedPKfundmor},
    which agrees with the $\kos{I}_{\pi}^{\vphantom{*}}$
    introduced in the proof of Proposition~\ref{prop-dualityoverGropk-RepKpseudofunct}.
    The component
    of the unity constraint of $\mathsf{Rep}$
    at $(\kos{K},\Spec(\pi))$
    is the $2$-cell
    consisting of the identity Grothendieck transformation of
    $\vphantom{\id}\smash[t]{\id}_{\kos{K}}^{\vphantom{*}}$
    and the invertible Grothendieck $\kos{K}$-transformation
    from the identity Grothendieck $\kos{K}$-morphism
    of $\kos{t}_{\pi}^{\vphantom{*}}$
    to $\kos{I}_{\pi}^{\vphantom{*}}$,
    given by the equality between the left adjoint parts.
  \end{itemize}
  One can check that $\mathsf{Rep}$ preserves
  horizontal composition of $2$-cells
  up to the functorial constraint,
  or equivalently,
  that the functorial constraint
  $\kos{f}\kos{f}^{\pr}_{\phi}\cong \kos{f}^{\ppr}$
  is natural in $(\phi,f)$, $(\phi^{\pr},f^{\pr})$.
  It is immediate that 
  the functorial and unity constraints of $\mathsf{Rep}$
  satisfy the associativity and left/right unity axioms,
  as the axioms are just equalities between the left adjoint parts.
  This completes the construction of the pseudofunctor
  $\mathsf{Rep}\colon \mathsf{pGroObj}\to \mathsf{pGroCat}_{\mathsf{pt}}$
  over $\mathsf{pGroKos}_{\etsf}$.
  The strict fiber of $\mathsf{Rep}$ over
  each Grothendieck prekosmos $\kos{K}$
  is $2$-isomorphic to the pseudofunctor $\mathsf{Rep}_{\kos{K}}$
  in Proposition~\ref{prop-dualityoverGropk-RepKpseudofunct}.
\end{proof}

\begin{proposition}
  \label{prop-dualityPKGro-Aut}
  We have a $2$-functor
  \begin{equation*}
    \mathsf{Aut}\colon \mathsf{pGroCat}_{\mathsf{pt}}\to \mathsf{pGroObj}
  \end{equation*}
  over $\mathsf{pGroKos}_{\etsf}$
  that sends each pointed prekosmic Grothendieck category
  $(\kos{K},\kos{t},\omega)$
  to the prekosmic Grothendieck object $(\kos{K},\Spec(\omega^*\omega_*(\kappa)))$,
  where $\Spec(\omega^*\omega_*(\kappa))$
  is the representing object of $\uAut_{\kos{K}}(\omega)$
  in Theorem~\ref{thm-Grocatoverpk-mainthm}.
\end{proposition}
\begin{proof}
  We describe how $\mathsf{Aut}$ sends $1$-cells and $2$-cells.
  \begin{itemize}
    \item
    A $1$-cell $(\phi,\kos{f},\alpha)\colon (\kos{K}^{\pr},\kos{t}^{\pr},\omega^{\pr})\mor(\kos{K},\kos{t},\omega)$
    in $\mathsf{pGroCat}_{\mathsf{pt}}$ is sent to
    the $1$-cell
    $(\phi,\breve{f})\colon (\kos{K}^{\pr},\Spec(\omega^{\pr*}\omega^{\pr}_*(\kappa^{\pr})))\to (\kos{K},\Spec(\omega^*\omega_*(\kappa)))$
    in $\mathsf{pGroObj}$
    described in the proof of Proposition~\ref{prop-PKGrofund-inducedPKfundmor}.

    \item
    A $2$-cell $(\chi,\vartheta)\colon (\phi,\kos{f},\alpha)\cong (\psi,\kos{g},\beta)$
    in $\mathsf{pGroCat}_{\mathsf{pt}}$
    between $1$-cells
    $(\phi,\kos{f},\alpha)$, $(\psi,\kos{g},\beta)
    \colon (\kos{K}^{\pr},\kos{t}^{\pr},\omega^{\pr})\mor (\kos{K},\kos{t},\omega)$
    is sent to the $2$-cell
    $(\chi,\breve{\theta})\colon (\phi,\breve{f})\cong (\psi,\breve{g})$
    in $\mathsf{pGroObj}$,
    defined as follows.
    Put 
    $\pi\coloneqq \omega^*\omega_*(\kappa)$
    and
    $\vartheta^{\sharp}\coloneqq 
    (\id_{\chi}^{\pr*}\psi^*\epsilon)
    \circ(\id_{\chi}^{\pr*}\beta^{-1}\omega_*)
    \circ((\omega^{\pr*})_{\phi}\vartheta\omega_*)
    \circ(\alpha\omega_*)
    \colon \phi^*\omega^*\omega_*\Rightarrow \id_{\chi}^{\pr*}\psi^*$.
    We define
    $\breve{\theta}\colon \Spec(\kappa^{\pr})\to \Spec(\phi^*(\pi))$
    as 
    $(\psi^*_{\unit})^{-1}\circ \vartheta^{\sharp}_{\kappa}
    \colon \phi^*(\pi)\to \kappa^{\pr}$.
    By \eqref{eq-ActionbaseGro-laxKtnatequiv}
    for the lax $\kos{K}$-tensor natural transformation $\vartheta^{\sharp}$
    with $\monbase{t}=\phi^*$, $\monbase{t}^{\pr}=\vphantom{\id}\smash[t]{\id}_{\kos{K}}^{\vphantom{*}}$, $X^{\pr}=\kappa$,
    we deduce that the following diagram commutes
    for every $z\in \CK$.
    \begin{equation}
      \label{eq-dualityPKGro-thetasharpequiv}
      \vcenter{\hbox{
        \xymatrix@R=12pt@C=45pt{
          \phi^*(z)\otimes^{\pr}\phi^*\omega^*\omega_*(\kappa)
          \ar[r]^-{(\phi^*_{\tensor})_{z,\omega^*\omega_*(\kappa)}}_-{\cong}
          \ar[d]_-{I_{\phi^*(z)}\otimes^{\pr}\breve{\theta}}
          &\phi^*(z\otimes \omega^*\omega_*(\kappa))
          \ar[r]^-{\phi^*((\tahar{\omega^*\omega_*})_z)}_-{\cong}
          &\phi^*\omega^*\omega_*(z)
          \ar[d]^-{\vartheta^{\sharp}_z}
          \\%1
          \phi^*(z)
          \ar[rr]^-{\chi_z}_-{\cong}
          &\text{ }
          &\psi^*(z)
        }
      }}
    \end{equation}
    Put 
    $\varpi\coloneqq \omega^*\omega_*$, $\varpi^{\pr}\coloneqq \omega^{\pr*}\omega^{\pr}_*$,
    and $\delta\coloneqq \omega^*\eta\omega_*\colon \varpi\Rightarrow\varpi\varpi$.
    One can check that 
    $\vartheta^{\sharp}
    =(\epsilon^{\pr}\psi^*)
    \circ (\omega^{\pr*}\beta_*)
    \circ (\omega^{\pr*}\vartheta\omega_*)
    \circ(\alpha\omega_*)
    \colon \phi^*\omega^*\omega_*\Rightarrow \psi^*$
    and that    
    $(\varpi^{\pr}\vartheta^{\sharp})
    \circ (\alpha^{\sharp}\varpi)
    \circ (\phi^*\delta)
    =(\omega^{\pr*}\beta_*)
    \circ(\omega^{\pr*}\vartheta\omega_*)
    \circ(\alpha\omega_*)
    =(\beta^{\sharp})
    \circ (\vartheta^{\sharp}\varpi)
    \circ (\phi^*\delta)
    \colon \phi^*\varpi\Rightarrow \varpi^{\pr}\psi^*$
    as monoidal natural transformations between lax monoidal functors,
    using \eqref{eq-PKGrofund-matealpha}.
    From these equalities,
    along with the relations
    \eqref{eq-PKGrofund-alphasharpequiv} for $z=\pi$
    and
    \eqref{eq-dualityPKGro-thetasharpequiv} for $z=\pi$,
    we deduce that $\breve{\theta}$ is indeed a $2$-isomorphism
    from $\breve{f}$ to $\chi_{\pi}\breve{g}$.
  \end{itemize}
  It is straightforward to check that
  $\mathsf{Aut}$ preserves identity $1$-cells and identity $2$-cells.
  \begin{itemize}
    \item 
    Let $(\chi,\vartheta)\colon (\phi,\kos{f},\alpha)\cong (\psi,\kos{g},\beta)$
    and $(\tilde{\chi},\tilde{\vartheta})\colon (\psi,\kos{g},\beta)\cong (\tilde{\psi},\tilde{\kos{g}},\tilde{\beta})$
    be $2$-cells in $\mathsf{pGroCat}_{\mathsf{pt}}$,
    with vertical composite written as
    $(\bar{\chi},\bar{\vartheta})\colon (\phi,\kos{f},\alpha)\cong (\tilde{\psi},\tilde{\kos{g}},\tilde{\beta})$.
    Put $\pi\coloneqq \omega^*\omega_*(\kappa)$.
    One can check that
    $\bar{\vartheta}^{\sharp}
    =\tilde{\vartheta}^{\sharp}
    \circ (\vartheta^{\sharp}\omega^*\omega_*)
    \circ (\phi^*\omega^*\eta\omega_*)
    \colon \phi^*\omega^*\omega_*\Rightarrow \tilde{\psi}^*$
    as monoidal natural transformations
    between lax monoidal functors,
    using \eqref{eq-PKGrofund-matealpha}.
    From this,
    along with \eqref{eq-dualityPKGro-thetasharpequiv}
    for $z=\pi$,
    we deduce that 
    $(\breve{\tilde{\theta}}\circ \chi_{\pi})\star\breve{\theta}=\breve{\bar{\theta}}
    \colon \phi^*(\pi)\to \kappa^{\pr}$
    as morphisms in $\cat{Comm}(\kos{K}^{\pr})$.
    Hence $(\bar{\chi},\breve{\bar{\theta}})$
    is the vertical composite of $(\chi,\breve{\theta})$ and $(\tilde{\chi},\breve{\tilde{\theta}})$;
    this shows that $\mathsf{Aut}$ preserves vertical composition of $2$-cells.
    
    \item 
    Let 
    $(\phi^{\ppr},\kos{f}^{\ppr},\alpha^{\ppr})$
    be the composite of
    $(\phi,\kos{f},\alpha)\colon (\kos{K}^{\pr},\kos{t}^{\pr},\omega^{\pr})\mor(\kos{K},\kos{t},\omega)$
    and
    $(\phi^{\pr},\kos{f}^{\pr},\alpha^{\pr})\colon (\kos{K}^{\ppr},\kos{t}^{\ppr},\omega^{\ppr})\mor(\kos{K}^{\pr},\kos{t}^{\pr},\omega^{\pr})$
    in $\mathsf{pGroCat}_{\mathsf{pt}}$,
    with
    $\phi^{\ppr}\coloneqq \phi\phi^{\pr}$,
    $\kos{f}^{\ppr}\coloneqq \kos{f}\kos{f}^{\pr}_{\phi}$,
    and
    $\alpha^{\ppr}\coloneqq (\kos{f}\alpha^{\pr}_{\phi})\circ (\alpha\phi^{\pr}_{\phi})$.
    One can check that
    $\alpha^{\ppr}_*=(\alpha^{\pr}_*\phi^*)\circ (\kos{f}^{\pr*}\alpha_*)
    \colon \kos{f}^{\ppr*}\omega_*\Rightarrow \omega^{\ppr}_*\phi^{\ppr*}$
    and that
    $\alpha^{\ppr\sharp}
    =(\alpha^{\pr\sharp}\phi^*)\circ (\phi^{\pr*}\alpha^{\sharp})
    \colon \phi^{\ppr*}\omega^*\omega_*\Rightarrow
    \omega^{\ppr*}\omega^{\ppr}_*\phi^{\ppr*}$
    as monoidal natural transformations between lax monoidal functors,
    using \eqref{eq-PKGrofund-matealpha}.
    Thus we have
    $\breve{f}^{\ppr}=\phi^{\pr*}(\breve{f})\circ \breve{f}^{\pr}
    \colon \Spec(\omega^{\ppr*}\omega^{\ppr}_*(\kappa^{\ppr}))
    \to \Spec(\phi^{\ppr*}\omega^*\omega_*(\kappa))$.
    This shows that $\mathsf{Aut}$ preserves composition of $1$-cells.
    
    \item
    Let
    $(\chi,\vartheta)\colon (\phi,\kos{f},\alpha)\cong (\psi,\kos{g},\beta)$
    be a $2$-cell
    between $1$-cells from $(\kos{K}^{\pr},\kos{t}^{\pr},\omega^{\pr})$ to $(\kos{K},\kos{t},\omega)$,
    and $(\chi^{\pr},\vartheta^{\pr})\colon (\phi^{\pr},\kos{f}^{\pr},\alpha^{\pr})\cong (\psi^{\pr},\kos{g}^{\pr},\beta^{\pr})$
    be a $2$-cell
    between $1$-cells from
    $(\kos{K}^{\ppr},\kos{t}^{\ppr},\omega^{\ppr})$
    to $(\kos{K}^{\pr},\kos{t}^{\pr},\omega^{\pr})$.
    Write
    $(\phi^{\ppr},\kos{f}^{\ppr},\alpha^{\ppr})$,
    $(\psi^{\ppr},\kos{g}^{\ppr},\beta^{\ppr})$
    for the composite $1$-cells,
    and
    $(\chi^{\ppr},\vartheta^{\ppr})\colon (\phi^{\ppr},\kos{f}^{\ppr},\alpha^{\ppr})\cong(\psi^{\ppr},\kos{g}^{\ppr},\beta^{\ppr})$
    for the horizontal composite of $(\chi,\vartheta)$ and $(\chi^{\pr},\vartheta^{\pr})$.
    Put $\pi\coloneqq \omega^*\omega_*(\kappa)$.
    One can check that
    $\vartheta^{\ppr\sharp}
    =(\vartheta^{\pr\sharp}\psi^*)
    \circ (\phi^{\pr*}\beta^{\sharp})
    \circ (\phi^{\pr*}\vartheta^{\sharp}\omega^*\omega_*)
    \circ (\phi^{\ppr*}\omega^*\eta\omega_*)
    \colon \phi^{\ppr*}\omega^*\omega_*\Rightarrow
    \psi^{\ppr*}$
    as monoidal natural transformations between lax monoidal functors,
    using \eqref{eq-PKGrofund-matealpha}.
    From this,
    along with \eqref{eq-dualityPKGro-thetasharpequiv} for $z=\pi$,
    we deduce that
    $\breve{\theta}^{\ppr}
    =(\breve{\theta}^{\pr}\circ \phi^{\pr*}(\breve{g}\circ \chi_{\pi}))
    \star \phi^{\pr*}(\breve{\theta})
    \colon \phi^{\ppr*}(\pi)\to \kappa^{\ppr}$
    as morphisms in $\cat{Comm}(\kos{K}^{\ppr})$.
    Hence $(\chi^{\ppr},\breve{\theta}^{\ppr})$
    is the horizontal composite of
    $(\chi,\breve{\theta})$ and $(\chi^{\pr},\breve{\theta}^{\pr})$;
    this shows that $\mathsf{Aut}$ preserves horizontal composition of $2$-cells.
  \end{itemize}
  This completes the construction of the $2$-functor
  $\mathsf{Aut}\colon \mathsf{pGroCat}_{\mathsf{pt}}\to \mathsf{pGroObj}$
  over $\mathsf{pGroKos}_{\etsf}$.
  The strict fiber of $\mathsf{Aut}$ over each Grothendieck prekosmos $\kos{K}$
  is $2$-isomorphic to the $2$-functor
  $\mathsf{Aut}_{\kos{K}}$ in Proposition~\ref{prop-dualityoverGropk-AutK2funct}.
\end{proof}

\begin{proof}[Proof of Theorem~\ref{thm-PKGroCat-equiv21catoveret}]
  We show that the pseudofunctor
  $\mathsf{Rep}$
  in Proposition~\ref{prop-dualityPKGro-Rep}
  and the $2$-functor
  $\mathsf{Aut}$
  in Proposition~\ref{prop-dualityPKGro-Aut}
  are mutually quasi-inverse equivalences of $(2,1)$-categories
  over $\mathsf{pGroKos}_{\etsf}$;
  that is, the composites
  $\mathsf{Aut}\circ \mathsf{Rep}$
  and $\mathsf{Rep}\circ \mathsf{Aut}$
  are equivalent, over $\mathsf{pGroKos}_{\etsf}$,
  to the identity $2$-functors
  of $\mathsf{pGroObj}$
  and $\mathsf{pGroCat}_{\mathsf{pt}}$,
  respectively.
  
  The composite
  $\mathsf{Aut}\circ \mathsf{Rep}\colon 
  \mathsf{pGroObj}\to \mathsf{pGroObj}$
  sends each $0$-cell $(\kos{K},\Spec(\pi))$
  to the $0$-cell $(\kos{K},\Spec(\kappa\otimes \pi))$,
  and $\imath_{\pi}\colon \Spec(\kappa\otimes \pi)\cong \Spec(\pi)$
  is an isomorphism of Grothendieck objects in $\kos{K}$;
  see Remark~\ref{rem-GroMainthm-reconstruction}.
  It is straightforward to check that
  $\mathsf{Aut}\circ \mathsf{Rep}$
  is a $2$-functor over $\mathsf{pGroKos}_{\etsf}$,
  and that the $1$-cells
  $(\vphantom{\id}\smash[t]{\id}_{\kos{K}}^{\vphantom{*}},\imath_{\pi})$
  constitute a $2$-natural isomorphism over $\mathsf{pGroKos}_{\etsf}$
  from $\mathsf{Aut}\circ \mathsf{Rep}$
  to the identity $2$-functor of $\mathsf{pGroObj}$.

  We now construct an equivalence
  over $\mathsf{pGroKos}_{\etsf}$
  from the composite $\mathsf{Rep}\circ \mathsf{Aut}$
  to the identity $2$-functor of $\mathsf{pGroCat}_{\mathsf{pt}}$
  as follows.
  Such an equivalence is an invertible strong transformation
  over $\mathsf{pGroKos}_{\etsf}$,
  with the same convention as in the proof of
  Theorem~\ref{thm-dualityoverGropk-equiv21cat};
  it is an equivalence internal to the $2$-category of
  pseudofunctors from $\mathsf{pGroCat}_{\mathsf{pt}}$ to itself,
  strong transformations, and modifications,
  all over $\mathsf{pGroKos}_{\etsf}$.
  We take the component $1$-cell of this transformation
  at each $0$-cell
  $(\kos{K},\kos{t},\omega)$ 
  to be the equivalence
  $(\kos{K},\kos{t}_{\omega^*\omega_*(\kappa)},\kos{e}_{\omega^*\omega_*(\kappa)})
  \simeq
  (\kos{K},\kos{t},\omega)$
  in the strict fiber over $\kos{K}$,
  consisting of $\vphantom{\id}\smash[t]{\id}_{\kos{K}}^{\vphantom{*}}$,
  the equivalence
  $\breve{\omega}\colon \kos{t}_{\omega^*\omega_*(\kappa)}\simeq\kos{t}$,
  and the invertible Grothendieck $\kos{K}$-transformation
  $\omega\cong \breve{\omega}\kos{e}_{\omega^*\omega_*(\kappa)}$;
  see Theorem~\ref{thm-Grocatoverpk-mainthm}.
  We take its component $2$-cell at each $1$-cell 
  $(\phi,\kos{f},\alpha)\colon (\kos{K}^{\pr},\kos{t}^{\pr},\omega^{\pr})\mor (\kos{K},\kos{t},\omega)$
  to be the $2$-cell
  consisting of the identity Grothendieck transformation
  of $\phi$ and the invertible Grothendieck $\kos{K}$-transformation
  $\breve{\alpha}\colon \breve{\omega}\breve{\kos{f}}\cong \kos{f}\breve{\omega}^{\pr}_{\phi}$
  in Proposition~\ref{prop-PKGrofund-inducedPKfundmor}.
  It is straightforward to check,
  using \cite[Proposition~6.2.16]{Johnson2021}, that
  these data indeed define an invertible strong transformation
  in the $2$-category of pseudofunctors
  over $\mathsf{pGroKos}_{\etsf}$.
  This completes the proof of Theorem~\ref{thm-PKGroCat-equiv21catoveret}.
\end{proof}

%% file: Galprekos.tex
\part{On Prekosmic Galois Categories}
\label{part-GalContext}

\chapter{Galois Prekosmoi}
\label{chap-Galprekos}

\begin{definition}
  \label{def-Galprekos-Galprekos}
  A \emph{Galois prekosmos} is a symmetric monoidal category
  whose underlying category has reflexive coequalizers.
\end{definition}

Throughout this chapter,
we fix a Galois prekosmos 
$\kos{S}=(\CS,\otimes,s)$,
with underlying category $\CS$,
monoidal product $\otimes\colon\CS\times\CS\to \CS$,
and unit object $s$.
Objects in $\CS$ are denoted by lowercase letters
such as $x$, $y$, and $z$.
We write $\kos{T}=(\CT,\tensor,\unit)$
for an arbitrary Galois prekosmos,
with underlying category $\CT$,
monoidal product $\tensor\colon\CT\times\CT\to \CT$,
and unit object $\unit$.
Objects in $\CT$ are denoted by uppercase letters
such as $X$, $Y$, and $Z$.

\begin{definition}
  \label{def-Galprekos-Galmor}
  Let $\kos{T}=(\CT,\tensor,\unit)$,
  $\kos{T}^{\pr}=(\CT^{\pr},\tensor^{\pr},\unit^{\pr})$
  be Galois prekosmoi.
  A \emph{Galois morphism}
  $\kos{f}\colon\kos{T}^{\pr}\mor \kos{T}$,
  denoted by an arrow with a filled tip,
  is a tuple $\kos{f}=(\kos{f}^!,\kos{f}_!,\eta,\epsilon)$
  consisting of a strong monoidal functor
  $\kos{f}^!\colon\kos{T}\to \kos{T}^{\pr}$,
  a functor $\kos{f}_!\colon\CT^{\pr}\to \CT$,
  and natural transformations
  $\eta\colon\vphantom{\id}\smash[t]{\id}_{\CT^{\pr}\vphantom{!}}^{\vphantom{!}}\Rightarrow \kos{f}^!\kos{f}_!$,
  $\epsilon\colon\kos{f}_!\kos{f}^!\Rightarrow \vphantom{\id}\smash[t]{\id}_{\CT\vphantom{!}}^{\vphantom{!}}$
  satisfying the triangle identities.
  Thus $\kos{f}^!$ is right adjoint to $\kos{f}_!$,
  and we refer to $\kos{f}^!$ (resp.\ $\kos{f}_!$)
  as the \emph{right adjoint part} (resp.\ \emph{left adjoint part})
  of $\kos{f}$.
  We omit $\eta$ and $\epsilon$ from the notation
  and write $\kos{f}=(\kos{f}_!\dashv \kos{f}^!)$.
\end{definition}

\begin{definition}
  \label{def-Galprekos-Galtrans}
  Let $\kos{T}$, $\kos{T}^{\pr}$ be Galois prekosmoi
  and $\kos{f}$, $\kos{g}\colon\kos{T}^{\pr}\mor \kos{T}$
  be Galois morphisms.
  A \emph{Galois transformation}
  $\vartheta\colon\kos{f}\Rightarrow \kos{g}$
  is a comonoidal natural transformation
  $\vartheta\colon\kos{f}^!\Rightarrow \kos{g}^!$
  between the right adjoint parts.
\end{definition}

We denote by $\mathtt{pGalKos}$ the $2$-category
(resp.\ by $\mathsf{pGalKos}$ the $(2,1)$-category)
of Galois prekosmoi,
Galois morphisms,
and Galois transformations
(resp.\ invertible Galois transformations).
Composition of Galois morphisms 
is given by that of adjunctions.
The identity Galois morphism
of a Galois prekosmos $\kos{T}$
is the tuple consisting of
identity functors and identity natural transformations;
we denote it by $\vphantom{\id}\smash[t]{\id}_{\kos{T}\vphantom{!}}^{\vphantom{!}}$
and write it as an equality
$\vphantom{\id}\smash[t]{\id}_{\kos{T}\vphantom{!}}^{\vphantom{!}}\colon\kos{T}=\kos{T}$.
Horizontal and vertical compositions
of Galois transformations
are given by those of comonoidal natural transformations
between the right adjoint parts.
Identity Galois transformations
are given by the identity natural transformations
of the right adjoint parts.

Let $\kos{T}$, $\kos{T}^{\pr}$ be Galois prekosmoi.
For each Galois morphism
$\kos{f}\colon\kos{T}^{\pr}\mor \kos{T}$,
the left adjoint part $\kos{f}_!$
has a unique colax monoidal functor structure
such that $\eta$ and $\epsilon$ are comonoidal natural transformations.
Conversely, if we are given colax monoidal functors
$\kos{f}^!\colon\kos{T}\to \kos{T}^{\pr}$,
$\kos{f}_!\colon\kos{T}^{\pr}\to \kos{T}$
and comonoidal natural transformations
$\eta\colon\vphantom{\id}\smash[t]{\id}_{\kos{T}^{\pr}\vphantom{!}}^{\vphantom{!}}\Rightarrow \kos{f}^!\kos{f}_!$,
$\epsilon\colon\kos{f}_!\kos{f}^!\Rightarrow \vphantom{\id}\smash[t]{\id}_{\kos{T}\vphantom{!}}^{\vphantom{!}}$
satisfying the triangle identities,
then $\kos{f}^!$ is necessarily strong monoidal,
and $\kos{f}=(\kos{f}_!\dashv \kos{f}^!)\colon\kos{T}^{\pr}\mor \kos{T}$
is a Galois morphism.
This is an instance of doctrinal adjunction \cite[Theorem~1.5]{Kelly1974}
applied to the opposite of the $2$-monad considered in Chapter~\ref{chap-Groprekos};
see \cite[section~1.3]{Kelly1974}.

To each Galois morphism
$\kos{f}\colon\kos{T}^{\pr}\mor \kos{T}$
between Galois prekosmoi
$\kos{T}=(\CT,\tensor,\unit)$,
$\kos{T}^{\pr}=(\CT^{\pr},\tensor^{\pr},\unit^{\pr})$,
we associate a natural transformation
$\varphi\colon
\kos{f}_!(\slot\tensor^{\pr}\kos{f}^!(\slot))
\Rightarrow \kos{f}_!(\slot)\tensor \slot$
between functors from $\CT^{\pr}\times\CT$ to $\CT$
whose component at $X^{\pr}\in \CT^{\pr}$, $Z\in \CT$ is
\begin{equation}
  \label{eq-Galprekos-projformula}
  \xymatrix@C=40pt{
    \varphi_{X^{\pr},Z}\colon
    \kos{f}_!(X^{\pr}\tensor^{\pr} \kos{f}^!(Z))
    \ar[r]^-{(\kos{f}_{!\tensor})_{X^{\pr},\kos{f}^!(Z)}}
    &\kos{f}_!(X^{\pr})\tensor \kos{f}_!\kos{f}^!(Z)
    \ar[r]^-{I_{\kos{f}_!(X^{\pr})}\tensor \epsilon_Z}
    &\kos{f}_!(X^{\pr})\tensor Z
    .
  }
\end{equation}
We say that $\kos{f}$ \emph{satisfies the projection formula},
or that \emph{the projection formula holds for} $\kos{f}$,
if the associated $\varphi$ is a natural isomorphism.

Galois morphisms have their own lexicon,
each of whose entries reflects
properties of the adjoint parts individually or in their interplay.
Each property listed below is parallel yet dual to
the property of Grothendieck morphisms
in Definition~\ref{def-Groprekos-Gromorprop}
that shares its name.
The terminology is adapted from that of
essential geometric morphisms between topoi;
see Remark~\ref{rem-Galprekos-Galmornames} for an explanation.

\begin{definition}
  \label{def-Galprekos-Galmorprop}
  A Galois morphism $\kos{f}$ between Galois prekosmoi is said to be
  \begin{enumerate}
    \item 
    \emph{connected}
    if the right adjoint part $\kos{f}^!$ is fully faithful.

    \item 
    \emph{locally connected}
    if it satisfies the projection formula.
    
    \item
    \emph{\'{e}tale}
    if it is locally connected
    and the left adjoint part $\kos{f}_!$ is conservative.
    
    \item 
    \emph{surjective}
    if the right adjoint part $\kos{f}^!$
    is conservative and preserves reflexive coequalizers.
    
    \item
    an \emph{embedding}
    if the left adjoint part $\kos{f}_!$
    is fully faithful.
  \end{enumerate}  
\end{definition}

\begin{remark}
  \label{rem-Galprekos-Galmornames}
  We explain the relation between each property of Galois morphisms
  in Definition~\ref{def-Galprekos-Galmorprop}
  and the property of essential geometric morphisms between topoi
  that shares its name.
  Every topos is a Galois prekosmos
  with respect to its cartesian monoidal structure.
  Let $f\colon \kos{T}^{\pr}\to\kos{T}$
  be an essential geometric morphism between topoi
  \cite[section~A4.1]{Johnstone2002a},
  with adjoint triple $f_!\dashv f^*\dashv f_*$.
  Then $f^*$ is strong monoidal, since it preserves finite products,
  and thus $\kos{f}=(f_!\dashv f^*)\colon\kos{T}^{\pr}\mor\kos{T}$
  is a Galois morphism,
  with $\kos{f}_!=f_!$ and $\kos{f}^!=f^*$.
  We write $f$ for the geometric morphism
  and $\kos{f}$ for the Galois morphism.
  Then the following hold:
  \begin{enumerate}
    \item
    $f$ is connected if and only if $\kos{f}$ is connected.
    Recall that $f$ is connected if $f^*$ is fully faithful
    \cite[section~C1.5]{Johnstone2002b}.

    \item
    If $f$ is locally connected, then $\kos{f}$ is locally connected.
    Recall that $f$ is locally connected if
    $f^*$ has a $\kos{T}$-indexed left adjoint
    \cite[section~C3.3]{Johnstone2002b}.
    Frobenius reciprocity for $f_!\dashv f^*$
    \cite[Lemma~A1.5.8]{Johnstone2002a}
    is the invertibility of
    the natural transformation $\varphi$
    in \eqref{eq-Galprekos-projformula},
    under the identifications above.
    Local connectedness of $f$ requires Frobenius reciprocity
    for every slice of $f$
    \cite[Proposition~C3.3.1]{Johnstone2002b},
    whereas local connectedness of $\kos{f}$ requires it
    only for $f$.
    
    \item
    If $f$ is \'{e}tale, then $\kos{f}$ is \'{e}tale.
    Recall that $f$ is \'{e}tale,
    or a local homeomorphism in \cite[section~C3.3]{Johnstone2002b},
    if $\kos{T}^{\pr}\simeq\kos{T}/A$ as topoi over $\kos{T}$
    for an object $A$ in $\kos{T}$.
    It follows that if
    $f$ is \'{e}tale,
    then $f$ is locally connected and $f_!$ is conservative,
    and by (2) $\kos{f}$ is then \'{e}tale.

    \item 
    $f$ is surjective if and only if $\kos{f}$ is surjective.
    Recall that $f$ is surjective if $f^*$ is conservative
    \cite[section~A4.2]{Johnstone2002a},
    and that $f^*$ preserves reflexive coequalizers,
    since it is a left adjoint.
    
    \item
    $f$ is an embedding if and only if $\kos{f}$ is an embedding.
    Recall that $f$ is an embedding if $f_*$ is fully faithful
    \cite[section~A4.2]{Johnstone2002a}, \cite[section~VII.4]{MacLaneMoerdijk1994},
    and that for any adjoint triple $f_!\dashv f^*\dashv f_*$,
    $f_!$ is fully faithful if and only if $f_*$ is.
  \end{enumerate}
  The terminology for Grothendieck morphisms
  in Definition~\ref{def-Groprekos-Gromorprop}
  is inherited from that for Galois morphisms
  in Definition~\ref{def-Galprekos-Galmorprop}.
\end{remark}

We write $\mathtt{pGalKos}/\kos{S}$
for the pseudoslice $2$-category of $\mathtt{pGalKos}$ over $\kos{S}$.
\begin{itemize}
  \item 
  A $0$-cell in $\mathtt{pGalKos}/\kos{S}$
  is a \emph{Galois $\kos{S}$-prekosmos},
  or a \emph{Galois prekosmos over $\kos{S}$}.
  It is a pair $(\kos{T},\kos{t})$ consisting of a Galois prekosmos $\kos{T}$
  and a Galois morphism $\kos{t}\colon\kos{T}\mor \kos{S}$.
  We write the pair as $\kos{T}\xmor{\kos{t}} \kos{S}$,
  or simply as $\kos{t}$.
  
  \item
  A $1$-cell in $\mathtt{pGalKos}/\kos{S}$
  is a \emph{Galois $\kos{S}$-morphism},
  or a \emph{Galois morphism over $\kos{S}$}.
  Let $\kos{T}\xmor{\kos{t}} \kos{S}$,
  $\kos{T}^{\pr}\xmor{\kos{t}^{\pr}} \kos{S}$
  be Galois $\kos{S}$-prekosmoi.
  A Galois $\kos{S}$-morphism
  $(\kos{f},\what{\kos{f}})\colon\kos{t}^{\pr}\mor \kos{t}$,
  denoted by an arrow with a filled tip,
  is a pair consisting of a Galois morphism
  $\kos{f}\colon\kos{T}^{\pr}\mor \kos{T}$
  and an invertible Galois transformation
  $\what{\kos{f}}\colon \kos{t}\kos{f}\cong \kos{t}^{\pr}$.
  We omit $\what{\kos{f}}$ from the notation and
  write $\kos{f}\colon\kos{t}^{\pr}\mor \kos{t}$.
  
  \item
  A $2$-cell in $\mathtt{pGalKos}/\kos{S}$
  is a \emph{Galois $\kos{S}$-transformation},
  or a \emph{Galois transformation over $\kos{S}$}.
  Let 
  $\kos{f}$, $\kos{g}\colon\kos{t}^{\pr}\mor \kos{t}$
  be Galois $\kos{S}$-morphisms.
  A Galois $\kos{S}$-transformation
  $\vartheta\colon\kos{f}\Rightarrow \kos{g}$
  is a Galois transformation
  between the underlying Galois morphisms
  $\kos{f}$, $\kos{g}\colon\kos{T}^{\pr}\mor \kos{T}$
  that satisfies
  $\what{\kos{f}}=
  \what{\kos{g}}\circ (\kos{t}\vartheta)
  \colon\kos{t}\kos{f}\Rightarrow \kos{t}^{\pr}$.
\end{itemize}
We illustrate Galois $\kos{S}$-morphisms
$\kos{f}$, $\kos{g}\colon\kos{t}^{\pr}\mor\kos{t}$
and a Galois $\kos{S}$-transformation
$\vartheta\colon\kos{f}\Rightarrow\kos{g}$
as follows.
\begin{equation*}
  \begin{tikzcd}[row sep=5pt,column sep=5pt]
    \kos{T}^{\pr}
    \ar[rr,filled,"\kos{f}"]
    \ar[ddr,filled,"\kos{t}^{\pr}"']
    &\text{ }
    &\kos{T}
    \ar[ddl,filled,"\kos{t}"]
    &\text{ }
    &\kos{T}^{\pr}
    \ar[rr,filled,"\kos{g}"]
    \ar[ddr,filled,"\kos{t}^{\pr}"']
    &\text{ }
    &\kos{T}
    \ar[ddl,filled,"\kos{t}"]
    &\text{ }
    &\kos{T}^{\pr}
    \ar[rr,filled,shift left=0.5ex,bend left=22pt,"\kos{f}" description]
    \ar[rr,filled,shift right=0.5ex,bend right=22pt,"\kos{g}" description]
    \ar[ddr,filled,"\kos{t}^{\pr}"']
    &\text{ }
    \ar[d,2cell,pos=-1.5,"{\,\vartheta}"]
    &\kos{T}
    \ar[ddl,filled,"\kos{t}"]
    \\%1
    \text{ }
    &\text{ }
    &\text{ }
    \ar[ll,2cell,shift left=0.3ex,"{\what{\kos{f}}}"']
    &\text{ }
    &\text{ }
    &\text{ }
    &\text{ }
    \ar[ll,2cell,shift left=0.3ex,"{\what{\kos{g}}}"']
    &\text{ }
    &\text{ }
    &\text{ }
    &\text{ }
    \\%2
    \text{ }
    &\kos{S}
    &\text{ }
    &\text{ }
    &\text{ }
    &\kos{S}
    &\text{ }
    &\text{ }
    &\text{ }
    &\kos{S}
    &\text{ }
  \end{tikzcd}
\end{equation*}

An \emph{equivalence} of Galois prekosmoi
(resp.\ of Galois $\kos{S}$-prekosmoi)
is an equivalence internal to $\mathtt{pGalKos}$
(resp.\ $\mathtt{pGalKos}/\kos{S}$).
A Galois morphism $\kos{f}$
between Galois prekosmoi
is an equivalence if and only if
$\eta$ and $\epsilon$ are natural isomorphisms;
the if part follows from doctrinal adjunction.
A Galois $\kos{S}$-morphism
between Galois $\kos{S}$-prekosmoi
is an equivalence
if and only if the underlying Galois morphism
is an equivalence.

Let $\prpty{P}$ be one of the following:
connected,
locally connected,
\'{e}tale,
surjective,
an embedding.
We say that a Galois $\kos{S}$-prekosmos
is $\prpty{P}$ if it is $\prpty{P}$
as a Galois morphism to $\kos{S}$.
We say that a Galois $\kos{S}$-morphism
between Galois $\kos{S}$-prekosmoi
is $\prpty{P}$
if the underlying Galois morphism
is $\prpty{P}$.
The property of being $\prpty{P}$
is preserved under composition with equivalences,
both for Galois morphisms
and for Galois $\kos{S}$-morphisms.

\begin{definition}
  \label{def-Galprekos-PKfiberpt}
  Let $\kos{T}\xmor{\kos{t}} \kos{S}$
  be a Galois $\kos{S}$-prekosmos.
  A \emph{fiber point} of $\kos{t}$
  is an \'{e}tale Galois $\kos{S}$-morphism
  $\omega\colon\vphantom{\id}\smash[t]{\id}_{\kos{S}\vphantom{!}}^{\vphantom{!}}\mor \kos{t}$.
  We illustrate $\omega$ as follows.
  \begin{equation*}
    \begin{tikzcd}[row sep=5pt,column sep=5pt]
      \kos{S}
      \ar[rr,filled,"\omega"]
      \ar[ddr,equal]
      &\text{ }
      &\kos{T}
      \ar[ddl,filled,"\kos{t}"]
      \\%1
      \text{ }
      &\text{ }
      &\text{ }
      \ar[ll,2cell,shift left=0.3ex,"{\,\what{\omega}}"']
      \\%2
      \text{ }
      &\kos{S}
      &\text{ }
    \end{tikzcd}
  \end{equation*}
\end{definition}

\begin{definition}
  \label{def-Galprekos-Galcatoverpk}
  A \emph{Galois $\kos{S}$-category},
  or a \emph{Galois category over $\kos{S}$},
  is a Galois $\kos{S}$-prekosmos
  that admits a surjective fiber point.
\end{definition}

Let $\kos{T}\xmor{\kos{t}}\kos{S}$
be a Galois $\kos{S}$-prekosmos.
Note that a fiber point
$\omega\colon\vphantom{\id}\smash[t]{\id}_{\kos{S}\vphantom{!}}^{\vphantom{!}}\mor\kos{t}$
of $\kos{t}$ is surjective
if and only if the underlying Galois morphism
$\omega\colon\kos{S}\mor\kos{T}$ is surjective.
A \emph{morphism} of fiber points of $\kos{t}$
is a Galois $\kos{S}$-transformation between them.
We denote by $\cat{pFib}_{\kos{S}}(\kos{t})$
the category of fiber points of $\kos{t}$.

Let $\kos{f}\colon\kos{t}^{\pr}\mor\kos{t}$
be a Galois $\kos{S}$-morphism
between Galois $\kos{S}$-prekosmoi.
When describing 
$\what{\kos{f}}\colon \kos{t}\kos{f}\cong \kos{t}^{\pr}$
as the comonoidal natural isomorphism between the right adjoint parts,
we write
$\what{\kos{f}^!}\colon \kos{f}^!\kos{t}^!\cong \kos{t}^{\pr!}$.
Let $\kos{f}^{\pr}\colon\kos{t}^{\ppr}\mor\kos{t}^{\pr}$
be a Galois $\kos{S}$-morphism,
and consider the composite
$\kos{f}\kos{f}^{\pr}\colon\kos{t}^{\ppr}\mor\kos{t}$.
When describing
$\what{\kos{f}\kos{f}^{\pr}}=\what{\kos{f}^{\pr}}\circ (\what{\kos{f}}\kos{f}^{\pr})
\colon \kos{t}\kos{f}\kos{f}^{\pr}\cong \kos{t}^{\ppr}$
as the comonoidal natural isomorphism between the right adjoint parts,
we write
$\what{\kos{f}^{\pr!}\kos{f}^!}
=\what{\kos{f}^{\pr!}}\circ (\kos{f}^{\pr!}\what{\kos{f}^!})
\colon \kos{f}^{\pr!}\kos{f}^!\kos{t}^!\cong \kos{t}^{\ppr!}$.

\begin{remark}
  \label{rem-Galprekos-omegalowershriekconservative}
  Let $\kos{T}\xmor{\kos{t}}\kos{S}$
  be a Galois $\kos{S}$-prekosmos
  and
  $\omega\colon \vphantom{\id}\smash[t]{\id}_{\kos{S}\vphantom{!}}^{\vphantom{!}}\mor\kos{t}$ be a Galois $\kos{S}$-morphism.
  Then the left adjoint part $\omega_!$ is conservative,
  and thus $\omega$ is 
  a fiber point of $\kos{t}$
  if and only if it is locally connected.
  Indeed, by uniqueness of left adjoints,
  the natural isomorphism
  $\what{\omega^!}\colon \omega^!\kos{t}^!\cong \vphantom{\id}\smash[t]{\id}_{\kos{S}\vphantom{!}}^{\vphantom{!}}$
  gives $\kos{t}_!\omega_!\cong \vphantom{\id}\smash[t]{\id}_{\kos{S}\vphantom{!}}^{\vphantom{!}}$.
  Hence the composite $\kos{t}_!\omega_!$ is conservative,
  and so is $\omega_!$.
\end{remark}

\section{Galois objects in a Galois prekosmos}
\label{sec-Galobjinpk}

We denote by $\kosEns(\kos{S})=(\cat{Ens}(\kos{S}),\otimes,s)$
the cartesian monoidal category
of cocommutative comonoids in $\kos{S}$.
The coproduct and counit morphisms
of a cocommutative comonoid $c$ in $\kos{S}$
are denoted by $\cp_c\colon c\to c\otimes c$
and $e_c\colon c\to s$.
Note that $\kosEns(\kos{S})$
is a Galois prekosmos,
since the forgetful functor
$\cat{Ens}(\kos{S})\to \CS$
creates reflexive coequalizers.

A cocommutative comonoid $\pi$ in $\kos{S}$
is a group object in $\kosEns(\kos{S})$ if and only if
$\pi$ is a cocommutative Hopf monoid in $\kos{S}$,
in which case we write
$\pc_{\pi}\colon \pi\otimes\pi\to \pi$,
$u_{\pi}\colon s\to \pi$,
and
$\varsigma_{\pi}\colon \pi\to \pi$
for the product, unit, and inverse morphisms
of $\pi$ in $\cat{Ens}(\kos{S})$;
they are the product, unit, and antipode of $\pi$ as a Hopf monoid in $\kos{S}$.
For such $\pi$, the antipode $\varsigma_{\pi}$
is involutive and anticommutes with $\pc_{\pi}$,
and the functor $\pi\otimes \slot\colon \CS\to \CS$ is conservative.

Let $\pi$ be a group object in $\kosEns(\kos{S})$.
For each cocommutative comonoid $c$ in $\kos{S}$,
the set $\Hom_{\cat{Ens}(\kos{S})}(c,\pi)$
is a group under the convolution product $\star$,
which sends each pair $\tilde{g}$, $g\colon c\to \pi$
to
$\tilde{g}\star g\coloneqq \pc_{\pi}\circ (\tilde{g}\otimes g)\circ \cp_c\colon c\to \pi$.
The identity element is
$u_{\pi}\circ e_c\colon c\to \pi$,
and the inverse of $g\colon c\to \pi$
is $\varsigma_{\pi}\circ g\colon c\to \pi$.
It follows from the Yoneda lemma that the assignment
$\pi\mapsto \Hom_{\cat{Ens}(\kos{S})}(\slot,\pi)$
is an equivalence
from the category of group objects in $\kosEns(\kos{S})$
to the category of presheaves of groups on $\cat{Ens}(\kos{S})$
whose underlying presheaves of sets are representable.

\begin{definition}
  \label{def-Galobjinpk-Galobj}
  A \emph{Galois object in $\kos{S}$}
  is a group object $\pi$ in $\kosEns(\kos{S})$
  such that the functor $\pi\otimes\slot\colon \CS\to \CS$
  preserves reflexive coequalizers.
\end{definition}

A \emph{morphism}
of Galois objects in $\kos{S}$
is a morphism of group objects in $\kosEns(\kos{S})$.
Let $f$, $g\colon \pi^{\pr}\to \pi$
be morphisms of Galois objects in $\kos{S}$.
A \emph{$2$-isomorphism}
$\theta\colon f\cong g$
is a morphism
$\theta\colon s\to \pi$ in $\cat{Ens}(\kos{S})$
that satisfies
$\pc_{\pi}\circ (\theta\otimes f)=\pc_{\pi}\circ (g\otimes \theta)
\colon \pi^{\pr}\to \pi$.
Galois objects in $\kos{S}$,
morphisms between them,
and $2$-isomorphisms between these morphisms
form a $(2,1)$-category,
which we denote by
\begin{equation}
  \label{eq-Galobjinpk-pGalObjS21cat}
  \mathsf{pGalObj}(\kos{S}).
\end{equation}
Let $\tilde{g}\colon \pi^{\pr}\to\pi$
be a morphism
and $\tilde{\theta}\colon g\cong \tilde{g}$
be a $2$-isomorphism.
The vertical composite of $\theta$ and $\tilde{\theta}$
is the convolution product $\tilde{\theta}\star\theta\colon f\cong \tilde{g}$.
The identity $2$-isomorphism of $f$
is $u_{\pi}\colon s\to \pi$.
Each $2$-isomorphism $\theta\colon f\cong g$
is indeed invertible with respect to $\star$;
its inverse $\theta^{-1}\colon g\cong f$ is
$\varsigma_{\pi}\circ \theta\colon s\to \pi$.
Let $f^{\pr}$, $g^{\pr}\colon \pi^{\ppr}\to \pi^{\pr}$
be morphisms,
and
$\theta^{\pr}\colon f^{\pr}\cong g^{\pr}$
be a $2$-isomorphism.
We write $ff^{\pr}$, $gg^{\pr}\colon \pi^{\ppr}\to \pi$
for the composites.
The horizontal composite
$\theta\theta^{\pr}\colon ff^{\pr}\cong gg^{\pr}$
of $\theta$ and $\theta^{\pr}$ is
$\theta\star (f\circ \theta^{\pr})=(g\circ \theta^{\pr})\star \theta
\colon s\to \pi$.

\section{Internal representations of a Galois object}
\label{sec-RepGalobj}

Let $c$ be a cocommutative comonoid in $\kos{S}$.
The endofunctor $c\otimes\slot$ on $\CS$ is a comonad,
with coproduct and counit natural transformations
induced by the coproduct and counit morphisms of $c$.
We denote by $\CS_{c\vphantom{!}}^{\vphantom{!}}$ the Kleisli category for this comonad.
Objects in $\CS_{c\vphantom{!}}^{\vphantom{!}}$ are written as $\kos{e}_{c\vphantom{!}}^!(x)$,
one for each object $x$ in $\CS$.
A morphism $\kos{e}_{c\vphantom{!}}^!(x)\to \kos{e}_{c\vphantom{!}}^!(y)$ in $\CS_{c\vphantom{!}}^{\vphantom{!}}$
is a morphism $c\otimes x\to y$ in $\CS$.
There is a canonical adjunction between $\CS$ and $\CS_{c\vphantom{!}}^{\vphantom{!}}$;
the right adjoint $\kos{e}_{c\vphantom{!}}^!\colon \CS\to \CS_{c\vphantom{!}}^{\vphantom{!}}$
sends $x$ to $\kos{e}_{c\vphantom{!}}^!(x)$,
and the left adjoint $\kos{e}_{c!}^{\vphantom{!}}\colon \CS_{c\vphantom{!}}^{\vphantom{!}}\to \CS$
sends $\kos{e}_{c\vphantom{!}}^!(x)$ to $c\otimes x$.
The Kleisli category $\CS_{c\vphantom{!}}^{\vphantom{!}}$ has a symmetric monoidal category structure,
with monoidal product
$\kos{e}_{c\vphantom{!}}^!(x)\otimes_{c\vphantom{!}}^{\vphantom{!}}\kos{e}_{c\vphantom{!}}^!(y)\coloneqq\kos{e}_{c\vphantom{!}}^!(x\otimes y)$
and unit object $\kos{e}_{c\vphantom{!}}^!(s)$.
We write 
\begin{equation*}
  \kos{S}_{c\vphantom{!}}^{\vphantom{!}}=(\CS_{c\vphantom{!}}^{\vphantom{!}},\otimes_{c\vphantom{!}}^{\vphantom{!}},\kos{e}_{c\vphantom{!}}^!(s))
\end{equation*}
for this symmetric monoidal category.
The right adjoint $\kos{e}_{c\vphantom{!}}^!$,
equipped with identity morphisms as its coherence isomorphisms,
is a strong monoidal functor from $\kos{S}$ to $\kos{S}_{c\vphantom{!}}^{\vphantom{!}}$.
The adjunction $\kos{e}_{c\vphantom{!}}^{\vphantom{!}}\coloneqq(\kos{e}_{c!}^{\vphantom{!}}\dashv \kos{e}_{c\vphantom{!}}^!)$
thus has the structure of a Galois morphism,
although $\kos{S}_{c\vphantom{!}}^{\vphantom{!}}$ is not, in general, a Galois prekosmos.
Nonetheless, we shall denote the adjunction by an arrow
with a filled tip
\begin{equation*}
  \kos{e}_{c\vphantom{!}}^{\vphantom{!}}\colon \kos{S}_{c\vphantom{!}}^{\vphantom{!}}\mor\kos{S}.
\end{equation*}
The left adjoint part $\kos{e}_{c!}^{\vphantom{!}}$ of $\kos{e}_{c\vphantom{!}}^{\vphantom{!}}$ is conservative,
and, again by doctrinal adjunction,
has a unique colax monoidal functor structure
such that the unit and counit
of $\kos{e}_{c\vphantom{!}}^{\vphantom{!}}$ are comonoidal natural transformations.
We have $\kos{e}_{c!}^{\vphantom{!}}\kos{e}_{c\vphantom{!}}^!=c\otimes\slot$
as colax monoidal endofunctors on $\kos{S}$,
where $c\otimes \slot$ is equipped with coherence morphisms
written as $(c\otimes)_{\tensor}$ and $(c\otimes)_{\unit}$,
explicitly given by
$(c\otimes)_{\tensor,x,y}
\coloneqq(I_c\otimes \sigma_{c,x}\otimes I_y)\circ (\cp_c\otimes I_{x\otimes y})
\colon c\otimes x\otimes y\to c\otimes x\otimes c\otimes y$
and
$(c\otimes)_{\unit}\coloneqq
e_c\otimes I_{s}\colon c\otimes s\to s$.
When $c=s$, we identify
$\kos{e}_{s\vphantom{!}}^{\vphantom{!}}\colon \kos{S}_{s\vphantom{!}}^{\vphantom{!}}\mor\kos{S}$
with $\vphantom{\id}\smash[t]{\id}_{\kos{S}\vphantom{!}}^{\vphantom{!}}\colon \kos{S}=\kos{S}$.

Let $f\colon c^{\pr}\to c$ be a morphism in $\cat{Ens}(\kos{S})$.
The \emph{base change} along $f$ is the pair
$(\kos{f}^!,\what{\kos{f}^!})$
consisting of the strong monoidal functor
$\kos{f}^!\colon \kos{S}_{c\vphantom{!}}^{\vphantom{!}}\to\kos{S}_{c^{\pr}\vphantom{!}}^{\vphantom{!}}$
that sends $\kos{e}_{c\vphantom{!}}^!(x)$ to $\kos{e}_{c^{\pr}}^!(x)$,
and the equality 
$\what{\kos{f}^!}\colon \kos{f}^!\kos{e}_{c\vphantom{!}}^!=\kos{e}_{c^{\pr}}^!$
of strong monoidal functors.
We omit $\what{\kos{f}^!}$ from the notation
and denote the pair by $\kos{f}^!\colon \kos{e}_{c\vphantom{!}}^!\to \kos{e}_{c^{\pr}}^!$.
The underlying functor 
$\kos{f}^!\colon\CS_{c\vphantom{!}}^{\vphantom{!}}\to \CS_{c^{\pr}\vphantom{!}}^{\vphantom{!}}$
is also referred to as the base change along $f$.
By the base change of a morphism in $\CS_{c\vphantom{!}}^{\vphantom{!}}$ along $f$
we mean its image under $\kos{f}^!$,
and we use the same terminology for
natural transformations whose components lie in $\CS_{c\vphantom{!}}^{\vphantom{!}}$.
The base change along $e_c\colon c\to s$
is the pair $(\kos{e}_{c\vphantom{!}}^!,\what{\kos{e}_{c\vphantom{!}}^!})$ consisting of
$\kos{e}_{c\vphantom{!}}^!\colon \kos{S}\to \kos{S}_{c\vphantom{!}}^{\vphantom{!}}$
and the equality
$\what{\kos{e}_{c\vphantom{!}}^!}\colon \kos{e}_{c\vphantom{!}}^!=\kos{e}_{c\vphantom{!}}^!$,
and we denote the pair by
$\kos{e}_{c\vphantom{!}}^!\colon \vphantom{\id}\smash[t]{\id}_{\kos{S}\vphantom{!}}^{\vphantom{!}}\to \kos{e}_{c\vphantom{!}}^!$.
The base change along $e_c$ is also described
in terms of the adjunction $\kos{e}_{c\vphantom{!}}^{\vphantom{!}}\colon \kos{S}_{c\vphantom{!}}^{\vphantom{!}}\mor\kos{S}$;
it is the pair 
$(\kos{e}_{c\vphantom{!}}^{\vphantom{!}},\what{\vphantom{\kos{e}_{c}}\smash[t]{\kos{e}_{c\vphantom{!}}^{\vphantom{!}}}})$
consisting of the adjunction $\kos{e}_{c\vphantom{!}}^{\vphantom{!}}$ and the equality 
$\what{\vphantom{\kos{e}_{c}}\smash[t]{\kos{e}_{c\vphantom{!}}^{\vphantom{!}}}}\colon \kos{e}_{c\vphantom{!}}^{\vphantom{!}}=\kos{e}_{c\vphantom{!}}^{\vphantom{!}}$,
and we denote the pair by an arrow with a filled tip
$\kos{e}_{c\vphantom{!}}^{\vphantom{!}}\colon \kos{e}_{c\vphantom{!}}^{\vphantom{!}}\mor\vphantom{\id}\smash[t]{\id}_{\kos{S}\vphantom{!}}^{\vphantom{!}}$.
We illustrate these pairs as follows.
\begin{equation}
  \label{eq-RepGalobj-basechange}
  \begin{tikzcd}[row sep=5pt,column sep=5pt]
    \kos{S}_{c^{\pr}}
    &\text{ }
    &\kos{S}_{c}
    \ar[ll,"\kos{f}^!"']
    &\text{ }
    &\text{ }
    &\kos{S}_{c}
    &\text{ }
    &\kos{S}
    \ar[ll,"\kos{e}_{c}^!"']
    &\text{ }
    &\text{ }
    &\kos{S}_{c}
    \ar[rr,filled,"\kos{e}_{c}"]
    \ar[ddr,filled,"\kos{e}_{c}"']
    &\text{ }
    &\kos{S}
    \ar[ddl,equal]
    \\%1
    \text{ }
    &\text{ }
    &\text{ }
    \ar[ll,2cell,shift left=0.3ex,"{\,\,\what{\kos{f}^!}}"']
    &\text{ }
    &\text{ }
    &\text{ }
    &\text{ }
    &\text{ }
    \ar[ll,2cell,shift left=0.3ex,"{\!\!\!\what{\kos{e}_{c}^!}}"']
    &\text{ }
    &\text{ }
    &\text{ }
    &\text{ }
    &\text{ }
    \ar[ll,2cell,shift left=0.3ex,"{\!\!\!\what{\kos{e}_{c}}}"']
    \\%2
    \text{ }
    &\kos{S}
    \ar[uul,"\kos{e}_{c^{\pr}}^!"]
    \ar[uur,"\kos{e}_{c}^!"']
    &\text{ }
    &\text{ }
    &\text{ }
    &\text{ }
    &\kos{S}
    \ar[uul,"\kos{e}_{c}^!"]
    \ar[uur,equal]
    &\text{ }
    &\text{ }
    &\text{ }
    &\text{ }
    &\kos{S}
    &\text{ }
  \end{tikzcd}
\end{equation}

Let $x$ be an object in $\CS$.
The assignment 
$c\mapsto \Aut_{\CS_{c\vphantom{!}}^{\vphantom{!}}}(\kos{e}_{c\vphantom{!}}^!(x))$,
sending each cocommutative comonoid $c$ in $\kos{S}$
to the group of automorphisms of $\kos{e}_{c\vphantom{!}}^!(x)$ in $\CS_{c\vphantom{!}}^{\vphantom{!}}$,
defines a presheaf of groups on $\cat{Ens}(\kos{S})$,
denoted by $\uAut(x)$,
with restriction morphisms given by base change.
Likewise, the assignment
$c\mapsto \End_{\CS_{c\vphantom{!}}^{\vphantom{!}}}(\kos{e}_{c\vphantom{!}}^!(x))$
defines a presheaf of monoids on $\cat{Ens}(\kos{S})$,
denoted by $\uEnd(x)$.

Let $\pi$ be a Galois object in $\kos{S}$
and $x\in \CS$.
A \emph{functorial action} of $\pi$ on $x$
is a datum $\upgamma_x$,
for each cocommutative comonoid $c$ in $\kos{S}$
and each morphism
$g\colon c\to \pi$ in $\cat{Ens}(\kos{S})$,
of a morphism
\begin{equation}
  \label{eq-RepGalobj-actionofg}
  \upgamma_x(g)\colon \kos{e}_{c\vphantom{!}}^!(x)\to \kos{e}_{c\vphantom{!}}^!(x)
\end{equation}
in $\CS_{c\vphantom{!}}^{\vphantom{!}}$,
and we require that
(i) $\upgamma_x(g)$ is compatible with base change
(that is, $\kos{f}^!(\upgamma_x(g))=\upgamma_x(g\circ f)\colon \kos{e}_{c^{\pr}}^!(x)\to\kos{e}_{c^{\pr}}^!(x)$
for every morphism $f\colon c^{\pr}\to c$ in $\cat{Ens}(\kos{S})$);
(ii) $\upgamma_x(\tilde{g}\star g)=\upgamma_x(\tilde{g})\circ \upgamma_x(g)$
for all $\tilde{g}$, $g\colon c\to \pi$;
and
(iii) $\upgamma_x(g)$ is the identity morphism of $\kos{e}_{c\vphantom{!}}^!(x)$
when $g=u_{\pi}\circ e_c\colon c\to \pi$.
Since every $g$ is invertible with respect to $\star$,
the morphisms $\upgamma_x(g)$ are invertible.
We say that $\upgamma_x$ is \emph{trivial}
if $\upgamma_x(g)$ is the identity morphism of $\kos{e}_{c\vphantom{!}}^!(x)$
for every $g$.

\begin{definition}
  \label{def-RepGalobj-Internalrep}
  Let $\pi$ be a Galois object in $\kos{S}$.
  A \emph{representation} of $\pi$ \emph{in $\kos{S}$},
  or an \emph{internal representation} of $\pi$,
  is a pair consisting of an object $x$ in $\CS$
  and a functorial action of $\pi$ on $x$.
  We say that a representation of $\pi$ in $\kos{S}$
  is \emph{trivial} if its functorial action is trivial.
\end{definition}

A functorial action $\upgamma_x$ of $\pi$
on an object $x$ in $\CS$ is the same as a morphism
$\upgamma_x\colon \Hom_{\cat{Ens}(\kos{S})}(\slot,\pi)\Rightarrow \uAut(x)$
of presheaves of groups on $\cat{Ens}(\kos{S})$,
and is, by the Yoneda lemma,
uniquely determined by $\upgamma_x(g)$ in the universal case
where $g$ is the identity morphism of $\pi$.
We write
\begin{equation*}
  \xi_x\colon \kos{e}_{\pi\vphantom{!}}^!(x)\cong \kos{e}_{\pi\vphantom{!}}^!(x)
\end{equation*}
for the universal element,
which is an isomorphism in $\CS_{\pi\vphantom{!}}^{\vphantom{!}}$.
Since $\upgamma_x$ is compatible with the group structures on both sides,
$\xi_x$ must satisfy the following:
(i)
in $\CS_{\pi\otimes\pi\vphantom{!}}^{\vphantom{!}}$,
the base change of $\xi_x$ along
$\pc_{\pi}\colon \pi\otimes \pi\to \pi$
is equal to $\kpr_1^!(\xi_x)\circ \kpr_2^!(\xi_x)$,
where each $\kpr_i^!$ is the base change along the $i$-th projection
$\mathrm{pr}_i\colon \pi\otimes \pi\to \pi$;
(ii)
in $\CS$,
the base change of $\xi_x$ along
$u_{\pi}\colon s\to \pi$
is equal to the identity morphism of $x$.
When we regard the isomorphism $\xi_x$ in $\CS_{\pi\vphantom{!}}^{\vphantom{!}}$
as a morphism in $\CS$, we denote it by
\begin{equation*}
  \gamma_x\colon \pi\otimes x\to x.
\end{equation*}
Then $\gamma_x$ must satisfy the left $\pi$-action relations,
and the assignment $\upgamma_x\mapsto \gamma_x$ is a one-to-one correspondence
between functorial actions of $\pi$ on $x$
and left $\pi$-actions on $x$.
We henceforth identify a representation
$(x,\upgamma_x)$ of $\pi$ in $\kos{S}$
with a left $\pi$-module $(x,\gamma_x)$.
A \emph{morphism} of representations of $\pi$ in $\kos{S}$
is a morphism between the underlying objects in $\CS$
that is equivariant under the functorial actions
of $\pi$, or, equivalently, one that is 
compatible with the corresponding left $\pi$-actions.
We denote by $\cat{Rep}_{\kos{S}}(\pi)$
the category of representations of $\pi$ in $\kos{S}$.

\begin{lemma}
  \label{lem-RepGalobj-Reppi}
  Let $\pi$ be a Galois object in $\kos{S}$.
  \begin{enumerate}
    \item 
    The representations of $\pi$ in $\kos{S}$
    form a Galois prekosmos, which we denote by
    $\kosRep_{\kos{S}}(\pi)=(\cat{Rep}_{\kos{S}}(\pi),\tensor_{\pi},\unit_{\pi})$.
    
    \item
    The Galois prekosmos 
    $\kosRep_{\kos{S}}(\pi)$
    is equipped with a connected and surjective Galois morphism
    to $\kos{S}$,
    denoted by $\kos{t}_{\pi\vphantom{!}}^{\vphantom{!}}\colon \kosRep_{\kos{S}}(\pi)\mor\kos{S}$,
    whose right adjoint part $\kos{t}_{\pi\vphantom{!}}^!$
    sends each $z\in \CS$
    to the trivial representation of $\pi$ in $\kos{S}$
    with underlying object $z$.
    
    \item 
    The Galois $\kos{S}$-prekosmos
    $\kosRep_{\kos{S}}(\pi)\xmor{\kos{t}_{\pi\vphantom{!}}^{\vphantom{!}}}\kos{S}$
    is equipped with a surjective fiber point,
    denoted by
    $\kos{u}_{\pi\vphantom{!}}^{\vphantom{!}}\colon \vphantom{\id}\smash[t]{\id}_{\kos{S}\vphantom{!}}^{\vphantom{!}}\mor\kos{t}_{\pi\vphantom{!}}^{\vphantom{!}}$,
    whose right adjoint part $\kos{u}_{\pi\vphantom{!}}^!$
    sends each $(x,\gamma_x)\in\cat{Rep}_{\kos{S}}(\pi)$
    to its underlying object $x$ in $\CS$.
  \end{enumerate}
  \begin{equation*}
    \begin{tikzcd}[row sep=5pt,column sep=4pt]
      \kos{S}
      \ar[rr,filled,"\kos{u}_{\pi}"]
      \ar[ddr,equal]
      &\text{ }
      &\kosRep_{\kos{S}}(\pi)
      \ar[ddl,filled,"\kos{t}_{\pi}"]
      \\%1
      \text{ }
      &\text{ }
      &\text{ }
      \ar[ll,2cell,pos=0.58,shift left=0.3ex,"{\!\!\what{\kos{u}_{\pi}}}"']
      \\%2
      \text{ }
      &\kos{S}
      &\text{ }
    \end{tikzcd}
  \end{equation*}
\end{lemma}
\begin{proof}
  We describe the structures in detail.
  The monoidal product $\tensor_{\pi}$
  sends each pair $(x,\gamma_x)$, $(y,\gamma_y)\in\cat{Rep}_{\kos{S}}(\pi)$
  to $(x\otimes y,\gamma_{x\otimes y})$,
  where 
  $\gamma_{x\otimes y}\coloneqq
  (\gamma_x\otimes \gamma_y)
  \circ (I_{\pi}\otimes \sigma_{\pi,x}\otimes I_y)
  \circ (\cp_{\pi}\otimes I_{x\otimes y})
  \colon \pi\otimes x\otimes y\to x\otimes y$.
  The unit object is
  $\unit_{\pi}\coloneqq(s,e_{\pi}\otimes I_s)$.
  The right adjoint part $\kos{t}_{\pi\vphantom{!}}^!$
  sends each $z\in \CS$ to $(z,e_{\pi}\otimes I_z)$,
  and the left adjoint part $\kos{t}_{\pi!}^{\vphantom{!}}$
  sends each $(x,\gamma_x)\in\cat{Rep}_{\kos{S}}(\pi)$
  to the coequalizer in $\CS$ of the reflexive pair
  $\gamma_x$, $e_{\pi}\otimes I_x\colon \pi\otimes x\to x$,
  with common section $u_{\pi}\otimes I_x$.
  The right adjoint part $\kos{u}_{\pi\vphantom{!}}^!$
  sends each $(x,\gamma_x)\in\cat{Rep}_{\kos{S}}(\pi)$ to $x$,
  and the left adjoint part $\kos{u}_{\pi!}^{\vphantom{!}}$
  sends each $z\in \CS$ to $(\pi\otimes z,\pc_{\pi}\otimes I_z)$.
  The unit and counit of $\kos{u}_{\pi\vphantom{!}}^{\vphantom{!}}$ are given by
  $\eta_z=u_{\pi}\otimes I_z\colon z\to \kos{u}_{\pi\vphantom{!}}^!\kos{u}_{\pi!}^{\vphantom{!}}(z)$
  and $\epsilon_{(x,\gamma_x)}=\gamma_x\colon \kos{u}_{\pi!}^{\vphantom{!}}\kos{u}_{\pi\vphantom{!}}^!(x,\gamma_x)\to (x,\gamma_x)$.
  The right adjoint parts
  $\kos{t}_{\pi\vphantom{!}}^!$ and $\kos{u}_{\pi\vphantom{!}}^!$
  are strong monoidal functors,
  whose coherence isomorphisms are identity morphisms.
  The invertible Galois transformation
  $\what{\vphantom{\kos{u}_{\pi}}\smash[t]{\kos{u}_{\pi\vphantom{!}}^{\vphantom{!}}}}\colon \kos{t}_{\pi\vphantom{!}}^{\vphantom{!}}\kos{u}_{\pi\vphantom{!}}^{\vphantom{!}}\cong \vphantom{\id}\smash[t]{\id}_{\kos{S}\vphantom{!}}^{\vphantom{!}}$
  is the equality
  $\what{\kos{u}_{\pi\vphantom{!}}^!}\colon \kos{u}_{\pi\vphantom{!}}^!\kos{t}_{\pi\vphantom{!}}^!=\vphantom{\id}\smash[t]{\id}_{\kos{S}\vphantom{!}}^{\vphantom{!}}$
  between the right adjoint parts.
  It is straightforward to check that
  $\kosRep_{\kos{S}}(\pi)$ is a Galois prekosmos,
  that $\kos{t}_{\pi\vphantom{!}}^{\vphantom{!}}\colon \kosRep_{\kos{S}}(\pi)\mor\kos{S}$
  is a connected Galois morphism,
  and that $\kos{u}_{\pi\vphantom{!}}^{\vphantom{!}}\colon \vphantom{\id}\smash[t]{\id}_{\kos{S}\vphantom{!}}^{\vphantom{!}}\mor\kos{t}_{\pi\vphantom{!}}^{\vphantom{!}}$
  is a surjective fiber point of $\kos{t}_{\pi\vphantom{!}}^{\vphantom{!}}$.
  The composite $\kos{t}_{\pi\vphantom{!}}^{\vphantom{!}}\kos{u}_{\pi\vphantom{!}}^{\vphantom{!}}$ is surjective,
  as it is isomorphic to $\vphantom{\id}\smash[t]{\id}_{\kos{S}\vphantom{!}}^{\vphantom{!}}$.
  Since $\kos{u}_{\pi\vphantom{!}}^{\vphantom{!}}$ is surjective,
  we deduce that $\kos{t}_{\pi\vphantom{!}}^{\vphantom{!}}$ is surjective.
\end{proof}

%% file: Galcatoverpk.tex
\chapter{Galois Categories over a Galois Prekosmos}
\label{chap-Galcatoverpk}

Throughout this chapter,
we fix a Galois prekosmos $\kos{S}=(\CS,\otimes,s)$.
Recall from Definition~\ref{def-Galprekos-Galcatoverpk}
the definition of a Galois $\kos{S}$-category,
and from Lemma~\ref{lem-RepGalobj-Reppi}
the Galois $\kos{S}$-prekosmos
$\kosRep_{\kos{S}}(\pi)\xmor{\kos{t}_{\pi\vphantom{!}}^{\vphantom{!}}}\kos{S}$
of internal representations of
a Galois object $\pi$ in $\kos{S}$.
The primary objective of Part~\ref{part-GalContext}
is to prove the following theorem.

\begin{theorem}
  \label{thm-Galcatoverpk-PrimaryObj}
  A Galois $\kos{S}$-prekosmos
  is a Galois $\kos{S}$-category
  if and only if it is equivalent
  to $\kosRep_{\kos{S}}(\pi)\xmor{\kos{t}_{\pi\vphantom{!}}^{\vphantom{!}}}\kos{S}$
  as Galois $\kos{S}$-prekosmoi
  for some Galois object $\pi$ in $\kos{S}$.
\end{theorem}

For each Galois object $\pi$ in $\kos{S}$,
$\kos{t}_{\pi\vphantom{!}}^{\vphantom{!}}$ admits a surjective fiber point
$\kos{u}_{\pi\vphantom{!}}^{\vphantom{!}}\colon\vphantom{\id}\smash[t]{\id}_{\kos{S}\vphantom{!}}^{\vphantom{!}}\mor\kos{t}_{\pi\vphantom{!}}^{\vphantom{!}}$
by Lemma~\ref{lem-RepGalobj-Reppi}.
Thus if a Galois $\kos{S}$-prekosmos
is equivalent to $\kos{t}_{\pi\vphantom{!}}^{\vphantom{!}}$
as Galois $\kos{S}$-prekosmoi
for some $\pi$,
then it is a Galois $\kos{S}$-category.
To prove the converse,
we consider the automorphisms of
a chosen surjective fiber point.

Let $\kos{T}\xmor{\kos{t}}\kos{S}$
be a Galois $\kos{S}$-prekosmos
and
$\omega\colon\vphantom{\id}\smash[t]{\id}_{\kos{S}\vphantom{!}}^{\vphantom{!}}\mor\kos{t}$
be a Galois $\kos{S}$-morphism.
For each $c\in\cat{Ens}(\kos{S})$,
recall the base change
$\kos{e}_{c\vphantom{!}}^{\vphantom{!}}\colon\kos{e}_{c\vphantom{!}}^{\vphantom{!}}\mor\vphantom{\id}\smash[t]{\id}_{\kos{S}\vphantom{!}}^{\vphantom{!}}$
along $e_c\colon c\to s$
in \eqref{eq-RepGalobj-basechange}.
We denote by
$\omega\kos{e}_{c\vphantom{!}}^{\vphantom{!}}\colon\kos{e}_{c\vphantom{!}}^{\vphantom{!}}\mor\kos{t}$
the composite of $\omega$ with $\kos{e}_{c\vphantom{!}}^{\vphantom{!}}$,
and write
$\what{\kos{e}_{c\vphantom{!}}^!\omega^!}
=\kos{e}_{c\vphantom{!}}^!\what{\omega^!}
\colon\kos{e}_{c\vphantom{!}}^!\omega^!\kos{t}^!\cong \kos{e}_{c\vphantom{!}}^!$.
Let $\omega^{\pr}\colon\vphantom{\id}\smash[t]{\id}_{\kos{S}\vphantom{!}}^{\vphantom{!}}\mor\kos{t}$
also be a Galois $\kos{S}$-morphism.
Although $\kos{S}_{c\vphantom{!}}^{\vphantom{!}}$ need not be a Galois prekosmos,
by a \emph{Galois $\kos{S}$-transformation}
$\vartheta\colon\omega\kos{e}_{c\vphantom{!}}^{\vphantom{!}}\Rightarrow\omega^{\pr}\!\kos{e}_{c\vphantom{!}}^{\vphantom{!}}$
we mean a comonoidal natural transformation
$\vartheta\colon\kos{e}_{c\vphantom{!}}^!\omega^!\Rightarrow \kos{e}_{c\vphantom{!}}^!\omega^{\pr!}$
that satisfies
$\what{\kos{e}_{c\vphantom{!}}^!\omega^!}=\what{\kos{e}_{c\vphantom{!}}^!\omega^{\pr!}}\circ (\vartheta\kos{t}^!)$.
We say that $\vartheta$ is \emph{invertible}
if every component of $\vartheta$ is an isomorphism.
We illustrate
$\omega\kos{e}_{c\vphantom{!}}^{\vphantom{!}}$, $\omega^{\pr}\!\kos{e}_{c\vphantom{!}}^{\vphantom{!}}$,
and $\vartheta$ as follows.
\begin{equation*}
  \begin{tikzcd}[row sep=5pt,column sep=15pt]
    \kos{S}_{c}
    \ar[r,filled,"\kos{e}_{c}"]
    \ar[ddr,filled,"\kos{e}_{c}"']
    &\kos{S}
    \ar[dd,equal]
    \ar[r,filled,"\omega"]
    &\kos{T}
    \ar[ddl,filled,"\kos{t}"]
    &\kos{S}_{c}
    \ar[r,filled,"\kos{e}_{c}"]
    \ar[ddr,filled,"\kos{e}_{c}"']
    &\kos{S}
    \ar[dd,equal]
    \ar[r,filled,"\omega^{\pr}"]
    &\kos{T}
    \ar[ddl,filled,"\kos{t}"]
    &\kos{S}_{c}
    \ar[ddr,filled,"\kos{e}_{c}"']
    \ar[rr,filled,bend left=25pt,"\omega\kos{e}_{c}" description]
    \ar[rr,filled,bend right=25pt,"\omega^{\pr}\!\kos{e}_{c}" description]
    &\text{ }
    \ar[d,2cell,pos=-1.2,"{\,\vartheta}"]
    &\kos{T}
    \ar[ddl,filled,"\kos{t}"]
    \\%1
    \text{ }
    &\text{ }
    \ar[l,2cell,pos=0.1,"{\!\!\!\what{\kos{e}_{c}}}"']
    &\text{ }
    \ar[l,2cell,pos=0.9,"{\!\what{\omega}}"']
    &\text{ }
    &\text{ }
    \ar[l,2cell,pos=0.1,"{\!\!\!\what{\kos{e}_{c}}}"']
    &\text{ }
    \ar[l,2cell,pos=0.9,"{\,\what{\omega^{\pr}}}"']
    &\text{ }
    &\text{ }
    &\text{ }
    \\%2
    \text{ }
    &\kos{S}
    &\text{ }
    &\text{ }
    &\kos{S}
    &\text{ }
    &\text{ }
    &\kos{S}
    &\text{ }
  \end{tikzcd}
\end{equation*}
The assignment sending each $c$
to the set of Galois $\kos{S}$-transformations
$\omega\kos{e}_{c\vphantom{!}}^{\vphantom{!}}\Rightarrow \omega^{\pr}\!\kos{e}_{c\vphantom{!}}^{\vphantom{!}}$
defines a presheaf on $\cat{Ens}(\kos{S})$
with restriction morphisms given by base change;
we call it the 
\emph{presheaf on $\cat{Ens}(\kos{S})$
of Galois $\kos{S}$-transformations}
from $\omega$ to $\omega^{\pr}$,
and denote it by $\uHom_{\kos{S}}(\omega,\omega^{\pr})$.
Likewise, we define the 
\emph{presheaf on $\cat{Ens}(\kos{S})$
of invertible Galois $\kos{S}$-transformations}
from $\omega$ to $\omega^{\pr}$,
denoted by $\uIsom_{\kos{S}}(\omega,\omega^{\pr})$,
as the presheaf sending each $c$
to the set of invertible Galois $\kos{S}$-transformations
$\omega\kos{e}_{c\vphantom{!}}^{\vphantom{!}}\cong\omega^{\pr}\!\kos{e}_{c\vphantom{!}}^{\vphantom{!}}$.
When $\omega^{\pr}=\omega$,
we define the 
\emph{presheaf of groups on $\cat{Ens}(\kos{S})$
of invertible Galois $\kos{S}$-transformations}
from $\omega$ to itself,
denoted by $\uAut_{\kos{S}}(\omega)$,
as the presheaf sending each $c$
to the group of invertible Galois $\kos{S}$-transformations
$\omega\kos{e}_{c\vphantom{!}}^{\vphantom{!}}\cong\omega\kos{e}_{c\vphantom{!}}^{\vphantom{!}}$.
We also define the
\emph{presheaf of monoids on $\cat{Ens}(\kos{S})$
of Galois $\kos{S}$-transformations}
from $\omega$ to itself,
denoted by $\uEnd_{\kos{S}}(\omega)$,
as the presheaf sending each $c$
to the monoid of Galois $\kos{S}$-transformations
$\omega\kos{e}_{c\vphantom{!}}^{\vphantom{!}}\Rightarrow\omega\kos{e}_{c\vphantom{!}}^{\vphantom{!}}$.

The following is the main theorem of Chapter~\ref{chap-Galcatoverpk}.

\begin{theorem}
  \label{thm-Galcatoverpk-mainthm}
  Let $\kos{T}\xmor{\kos{t}}\kos{S}$
  be a Galois $\kos{S}$-prekosmos
  and $\omega\colon\vphantom{\id}\smash[t]{\id}_{\kos{S}\vphantom{!}}^{\vphantom{!}}\mor\kos{t}$
  be a surjective fiber point of $\kos{t}$.
  \begin{enumerate}
    \item 
    The presheaf of groups 
    on $\cat{Ens}(\kos{S})$
    of invertible Galois $\kos{S}$-transformations
    from $\omega$ to itself
    is represented by a Galois object
    $\pi$ in $\kos{S}$,
    with $\pi\in\cat{Ens}(\kos{S})$
    given by $\omega^!\omega_!(s)$.
    
    \item
    The fiber point $\omega$ extends to an equivalence
    $\breve{\omega}\colon\kos{t}_{\pi\vphantom{!}}^{\vphantom{!}}\simeq\kos{t}$
    of Galois $\kos{S}$-prekosmoi
    up to an invertible Galois $\kos{S}$-transformation
    $\breve{\omega}\kos{u}_{\pi\vphantom{!}}^{\vphantom{!}}\cong \omega$
    given by the equality
    $\kos{u}_{\pi\vphantom{!}}^!\breve{\omega}^!=\omega^!$
    between the right adjoint parts.
  \end{enumerate}
  \begin{equation*}
    \Hom_{\cat{Ens}(\kos{S})}(\slot,\pi)
    \cong
    \uAut_{\kos{S}}(\omega)
    \qquad\quad
    \begin{tikzcd}[row sep=5pt,column sep=6pt]
      \text{ }
      &\id_{\kos{S}}
      \ar[ddr,filled,"{\omega}"]
      \ar[ddl,filled,"{\kos{u}_{\pi}}"']
      &\text{ }
      \\%1
      \text{ }
      \ar[rr,2cell,shift right=0.7ex,pos=0.5]
      &\text{ }
      &\text{ }
      \\%2
      \kos{t}_{\pi}
      \ar[rr,filled,"{\breve{\omega}}","\simeq"']
      &\text{ }
      &\kos{t}
    \end{tikzcd}
  \end{equation*}
\end{theorem}

Theorem~\ref{thm-Galcatoverpk-PrimaryObj} follows immediately
once we prove Theorem~\ref{thm-Galcatoverpk-mainthm}.

\begin{proof}[Proof of Theorem~\ref{thm-Galcatoverpk-PrimaryObj}]
  Let $\kos{T}\xmor{\kos{t}}\kos{S}$
  be a Galois $\kos{S}$-prekosmos.
  If an equivalence
  $\breve{\omega}\colon\kos{t}_{\pi\vphantom{!}}^{\vphantom{!}}\simeq\kos{t}$
  is given for some
  Galois object $\pi$ in $\kos{S}$,
  then $\kos{t}$ is a Galois $\kos{S}$-category,
  since it admits a surjective fiber point
  $\breve{\omega}\kos{u}_{\pi\vphantom{!}}^{\vphantom{!}}\colon\vphantom{\id}\smash[t]{\id}_{\kos{S}\vphantom{!}}^{\vphantom{!}}\mor\kos{t}$
  by Lemma~\ref{lem-RepGalobj-Reppi}.
  Conversely, if $\kos{t}$ is a Galois $\kos{S}$-category,
  we may choose a surjective fiber point $\omega$ of $\kos{t}$,
  and by Theorem~\ref{thm-Galcatoverpk-mainthm},
  we obtain the Galois object $\pi$
  in $\kos{S}$ that represents $\uAut_{\kos{S}}(\omega)$,
  as well as the equivalence $\breve{\omega}\colon\kos{t}_{\pi\vphantom{!}}^{\vphantom{!}}\simeq\kos{t}$.
\end{proof}

We shall prove Theorem~\ref{thm-Galcatoverpk-mainthm}
in Section~\ref{sec-GalMainthm}.
The proof is based on the following proposition.

\begin{proposition}
  \label{prop-Galcatoverpk-HomisIsom}
  Let $\kos{T}\xmor{\kos{t}}\kos{S}$ be a Galois $\kos{S}$-prekosmos
  and let $\omega$, $\omega^{\pr}\colon\vphantom{\id}\smash[t]{\id}_{\kos{S}\vphantom{!}}^{\vphantom{!}}\mor\kos{t}$
  be Galois $\kos{S}$-morphisms.
  Put $p\coloneqq\omega^{\pr!}\omega_!(s)\in \cat{Ens}(\kos{S})$.
  Assume that $\omega$ is a fiber point of $\kos{t}$.
  Then the following hold:
  \begin{enumerate}
    \item
    We have a comonoidal natural isomorphism
    $\omega^{\pr!}\omega_!\cong p\otimes \slot$
    between colax monoidal endofunctors on $\kos{S}$.
    
    \item 
    The presheaf 
    $\uHom_{\kos{S}}(\omega,\omega^{\pr})$
    on $\cat{Ens}(\kos{S})$
    is represented by $p$.

    \item
    We have an equality
    $\uHom_{\kos{S}}(\omega,\omega^{\pr})
    =\uIsom_{\kos{S}}(\omega,\omega^{\pr})$
    between presheaves on $\cat{Ens}(\kos{S})$.
  \end{enumerate}
\end{proposition}

We will prove Proposition~\ref{prop-Galcatoverpk-HomisIsom}
in Section~\ref{sec-ActionbaseGal}.
We present two direct consequences of our results.

\begin{corollary}
  \label{cor-Galcatoverpk-connsurj}
  Every Galois $\kos{S}$-category is connected and surjective.
\end{corollary}
\begin{proof}
  This is a direct consequence of
  Lemma~\ref{lem-RepGalobj-Reppi}
  and Theorem~\ref{thm-Galcatoverpk-PrimaryObj}.
\end{proof}

\begin{corollary}
  \label{cor-Galcatoverpk-pFibgroupoid}
  Let $\kos{T}\xmor{\kos{t}}\kos{S}$ be a Galois $\kos{S}$-prekosmos.
  The category $\cat{pFib}_{\kos{S}}(\kos{t})$
  of fiber points of $\kos{t}$ is a groupoid.
\end{corollary}
\begin{proof}
  This is a direct consequence of
  Proposition~\ref{prop-Galcatoverpk-HomisIsom}.
\end{proof}

Let $\kos{T}\xmor{\kos{t}}\kos{S}$
be a Galois $\kos{S}$-category.
In Section~\ref{sec-pFibGalTors},
we show that once a surjective fiber point
$\omega$ of $\kos{t}$ is chosen,
there is a perfect correspondence
between fiber points of $\kos{t}$
and co-fpqc $\pi$-torsors over $s$,
where $\pi$ is the Galois object in $\kos{S}$
determined by $\omega$ in Theorem~\ref{thm-Galcatoverpk-mainthm}.
Here we introduce properties of fiber points
that are relevant to this topic.

\begin{definition}
  \label{def-Galcatoverpk-cofpqccovering}
  We say that a cocommutative comonoid $c$ in $\kos{S}$
  is a \emph{co-fpqc covering of $s$}
  if $c\otimes\slot\colon\CS\to \CS$
  is conservative and preserves reflexive coequalizers.
\end{definition}

Note that a cocommutative comonoid $c$ in $\kos{S}$
is a co-fpqc covering of $s$
if and only if the base change $\kos{e}_{c\vphantom{!}}^!\colon\CS\to \CS_{c\vphantom{!}}^{\vphantom{!}}$
along $e_c\colon c\to s$
is conservative and preserves reflexive coequalizers.

\begin{definition}
  \label{def-Galcatoverpk-locisom}
  Let $\kos{T}\xmor{\kos{t}}\kos{S}$ be a Galois $\kos{S}$-prekosmos
  and $\omega$, $\omega^{\pr}\colon\vphantom{\id}\smash[t]{\id}_{\kos{S}\vphantom{!}}^{\vphantom{!}}\mor\kos{t}$
  be Galois $\kos{S}$-morphisms.
  We say that $\omega$ is \emph{locally isomorphic} to $\omega^{\pr}$
  if the presheaf
  $\uIsom_{\kos{S}}(\omega,\omega^{\pr})$
  evaluates to a nonempty set
  at some cocommutative comonoid $c$ in $\kos{S}$
  that is a co-fpqc covering of $s$.
  In this case, $\omega^{\pr}$ is locally isomorphic to $\omega$,
  and we also say that
  $\omega$ and $\omega^{\pr}$ are locally isomorphic.
\end{definition}

\begin{lemma}
  \label{lem-Galcatoverpk-pFiblocallyisom}
  Let $\kos{T}\xmor{\kos{t}}\kos{S}$
  be a Galois $\kos{S}$-prekosmos
  and $\omega$, $\omega^{\pr}\colon\vphantom{\id}\smash[t]{\id}_{\kos{S}\vphantom{!}}^{\vphantom{!}}\mor\kos{t}$
  be fiber points of $\kos{t}$.
  Assume that $\omega$ is surjective.
  Then the following hold:
  \begin{enumerate}
    \item 
    The fiber points $\omega$ and $\omega^{\pr}$ are locally isomorphic.
    
    \item
    The fiber point $\omega^{\pr}$ is also surjective.
  \end{enumerate}
\end{lemma}
\begin{proof}
  Put $p^{\pr}\coloneqq\omega^!\omega^{\pr}_!(s)\in\cat{Ens}(\kos{S})$.
  Since $\omega^{\pr}$ is a fiber point of $\kos{t}$,
  we exchange the roles of $\omega$ and $\omega^{\pr}$ in
  Proposition~\ref{prop-Galcatoverpk-HomisIsom}
  and deduce that $\omega^!\omega^{\pr}_!\cong p^{\pr}\otimes\slot$
  and that the presheaf
  $\uHom_{\kos{S}}(\omega^{\pr},\omega)
  =\uIsom_{\kos{S}}(\omega^{\pr},\omega)$
  is represented by $p^{\pr}$.
  The evaluation of the presheaf at $p^{\pr}$ is nonempty,
  as it contains the element
  $\omega^{\pr}\!\kos{e}_{p^{\pr}\vphantom{!}}^{\vphantom{!}}\cong \omega\kos{e}_{p^{\pr}\vphantom{!}}^{\vphantom{!}}$
  corresponding to the identity morphism of $p^{\pr}$.
  For the moment, let us say that a functor has property $\prpty{Q}$
  if it is conservative and preserves reflexive coequalizers.
  Then $\omega^{\pr}_!$ has property $\prpty{Q}$.
  Since $\omega$ is surjective, $\omega^!$ also has property $\prpty{Q}$.
  Hence the composite $\omega^!\omega^{\pr}_!$ has property $\prpty{Q}$,
  and therefore $p^{\pr}$ is a co-fpqc covering of $s$.
  This proves (1).
  Moreover, the functor $\kos{e}_{p^{\pr}}^!$ has property $\prpty{Q}$,
  and hence so does the composite $\kos{e}_{p^{\pr}}^!\omega^!$.
  From $\omega^{\pr}\!\kos{e}_{p^{\pr}\vphantom{!}}^{\vphantom{!}}\cong \omega\kos{e}_{p^{\pr}\vphantom{!}}^{\vphantom{!}}$
  we deduce that
  $\kos{e}_{p^{\pr}}^!\omega^{\pr!}$ has property $\prpty{Q}$.
  Since $\CS$ has reflexive coequalizers
  and $\kos{e}_{p^{\pr}}^!$ has property $\prpty{Q}$,
  we conclude that $\omega^{\pr!}$ has property $\prpty{Q}$.
  This proves (2).
\end{proof}

\begin{corollary}
  \label{cor-Galcatoverpk-pFibsurjective}
  Let $\kos{T}\xmor{\kos{t}}\kos{S}$
  be a Galois $\kos{S}$-category.
  Every fiber point of $\kos{t}$
  is surjective.
\end{corollary}
\begin{proof}
  This is a direct consequence of 
  Lemma~\ref{lem-Galcatoverpk-pFiblocallyisom}.
\end{proof}

\section{Action of the base Galois prekosmos}
\label{sec-ActionbaseGal}

In this section, we prove Proposition~\ref{prop-Galcatoverpk-HomisIsom}.
We utilize the facts that
(i) a Galois $\kos{S}$-prekosmos
is equipped with the induced action of $\kos{S}$;
(ii) the left and right adjoint parts
of a Galois $\kos{S}$-morphism
are equipped with equivariances
under the induced $\kos{S}$-actions;
and
(iii)
a Galois $\kos{S}$-transformation
between Galois $\kos{S}$-morphisms
is precisely a Galois transformation
between the underlying Galois morphisms
that is compatible with the $\kos{S}$-equivariances
of the right adjoint parts.
These facts are systematically
studied in a broader context
throughout Appendix~\ref{appdix-cltcatoversmc}.
In what follows, we present only what we need.

Let $\kos{T}\xmor{\kos{t}}\kos{S}$
be a Galois $\kos{S}$-prekosmos.
By the \emph{induced $\kos{S}$-action} on $\kos{T}$
we mean the functor $\bcts\colon\CT\times\CS\to\CT$
sending each pair $(X,z)$ to $X\bcts z\coloneqq X\tensor \kos{t}^!(z)$.
Let $\kos{T}^{\pr}\xmor{\kos{t}^{\pr}}\kos{S}$
be a Galois $\kos{S}$-prekosmos,
with the induced $\kos{S}$-action $\bcts^{\pr}$ on $\kos{T}^{\pr}$,
and let $\kos{f}\colon\kos{t}^{\pr}\mor\kos{t}$
be a Galois $\kos{S}$-morphism.
We denote by 
$\what{\kos{f}_!}\colon\kos{f}_!\kos{t}^{\pr!}\Rightarrow \kos{t}^!$
the mate of $(\what{\kos{f}^!})^{-1}\colon\kos{t}^{\pr!}\cong \kos{f}^!\kos{t}^!$
with respect to the adjunctions
$\kos{f}$ and $\vphantom{\id}\smash[t]{\id}_{\kos{S}\vphantom{!}}^{\vphantom{!}}$.
We illustrate $\kos{f}\colon\kos{t}^{\pr}\mor\kos{t}$,
$\what{\kos{f}^!}$, and $\what{\kos{f}_!}$ as follows.
\begin{equation*}
  \begin{tikzcd}[row sep=5pt,column sep=5pt]
    \kos{T}^{\pr}
    \ar[rr,filled,"\kos{f}"]
    \ar[ddr,filled,"\kos{t}^{\pr}"']
    &\text{ }
    &\kos{T}
    \ar[ddl,filled,"\kos{t}"]
    &\text{ }
    &\text{ }
    &\kos{T}^{\pr}
    &\text{ }
    &\kos{T}
    \ar[ll,"\kos{f}^!"']
    &\text{ }
    &\text{ }
    &\kos{T}^{\pr}
    \ar[rr,"\kos{f}_!"]
    &\text{ }
    &\kos{T}
    \\%1
    \text{ }
    &\text{ }
    &\text{ }
    \ar[ll,2cell,shift left=0.3ex,"{\what{\kos{f}}}"']
    &\text{ }
    &\text{ }
    &\text{ }
    &\text{ }
    &\text{ }
    \ar[ll,2cell,shift left=0.3ex,"{\,\,\what{\kos{f}^!}}"']
    &\text{ }
    &\text{ }
    &\text{ }
    \ar[rr,2cell,shift right=0.5ex,pos=0.46,"{\!\!\what{\kos{f}_!}}"]
    &\text{ }
    &\text{ }
    \\%2
    \text{ }
    &\kos{S}
    &\text{ }
    &\text{ }
    &\text{ }
    &\text{ }
    &\kos{S}
    \ar[uul,"\kos{t}^{\pr!}"]
    \ar[uur,"\kos{t}^!"']
    &\text{ }
    &\text{ }
    &\text{ }
    &\text{ }
    &\kos{S}
    \ar[uul,"\kos{t}^{\pr!}"]
    \ar[uur,"\kos{t}^!"']
    &\text{ }
  \end{tikzcd}
\end{equation*}
The \emph{$\kos{S}$-equivariance} of $\kos{f}^!$
is the natural isomorphism
$\vecar{\kos{f}^!}\colon \kos{f}^!(\slot\bcts\slot)\cong \kos{f}^!(\slot)\bcts^{\pr}\slot$
between functors from $\CT\times \CS$ to $\CT^{\pr}$
whose component at $X\in \CT$, $z\in \CS$ is
\begin{equation*}
  \xymatrix@C=50pt{
    \kos{f}^!(X\tensor \kos{t}^!(z))
    \ar[r]^-{(\kos{f}^!_{\tensor})_{X,\kos{t}^!(z)}}_-{\cong}
    &\kos{f}^!(X)\tensor^{\pr} \kos{f}^!\kos{t}^!(z)
    \ar[r]^-{I_{\kos{f}^!(X)}\tensor^{\pr}(\what{\kos{f}^!})_z}_-{\cong}
    &\kos{f}^!(X)\tensor^{\pr}\kos{t}^{\pr!}(z)
    .
  }
\end{equation*}
The \emph{$\kos{S}$-equivariance} of $\kos{f}_!$
is the natural transformation
$\vecar{\kos{f}_!}\colon\kos{f}_!(\slot\bcts^{\pr}\slot)\Rightarrow \kos{f}_!(\slot)\bcts\slot$
between functors from $\CT^{\pr}\times \CS$ to $\CT$
whose component at $X^{\pr}\in \CT^{\pr}$, $z\in \CS$ is
\begin{equation*}
  \xymatrix@C=50pt{
    \kos{f}_!(X^{\pr}\tensor^{\pr} \kos{t}^{\pr!}(z))
    \ar[r]^-{(\kos{f}_{!\tensor})_{X^{\pr},\kos{t}^{\pr!}(z)}}
    &\kos{f}_!(X^{\pr})\tensor \kos{f}_!\kos{t}^{\pr!}(z)
    \ar[r]^-{I_{\kos{f}_!(X^{\pr})}\tensor (\what{\kos{f}_!})_z}
    &\kos{f}_!(X^{\pr})\tensor \kos{t}^!(z)
    .
  }
\end{equation*}
  
\begin{lemma}
  \label{lem-ActionbaseGal-projimplystrongequiv}
  Let $\kos{f}\colon\kos{t}^{\pr}\mor\kos{t}$
  be a Galois $\kos{S}$-morphism.
  Recall the natural transformation
  $\varphi$ associated to $\kos{f}$ in \eqref{eq-Galprekos-projformula}.
  \begin{enumerate}
    \item 
    The component of $\vecar{\kos{f}_!}$
    at $X^{\pr}\in\CT^{\pr}$, $z\in \CS$
    can be redescribed as
    \begin{equation*}
      \xymatrix@C=15pt{
        \kos{f}_!(X^{\pr}\tensor^{\pr}\kos{t}^{\pr!}(z))
        \ar[rrr]^-{\kos{f}_!(I_{X^{\pr}}\tensor^{\pr} (\what{\kos{f}^!})_z^{-1})}_-{\cong}
        &\text{ }
        &\text{ }
        &\kos{f}_!(X^{\pr}\tensor^{\pr}\kos{f}^!\kos{t}^!(z))
        \ar[rr]^-{\varphi_{X^{\pr},\kos{t}^!(z)}}
        &\text{ }
        &\kos{f}_!(X^{\pr})\tensor\kos{t}^!(z)
        .
      }
    \end{equation*}

    \item
    If $\kos{f}$ is locally connected, then $\vecar{\kos{f}_!}$ is a natural isomorphism.
  \end{enumerate}
\end{lemma}
\begin{proof}
  (1) follows by substituting 
  $\what{\kos{f}_!}=(\epsilon\kos{t}^!)\circ (\kos{f}_!(\what{\kos{f}^!})^{-1})
  \colon\kos{f}_!\kos{t}^{\pr!}\Rightarrow \kos{t}^!$,
  and (2) is a direct consequence of (1).
\end{proof}

We now show that 
a Galois $\kos{S}$-transformation
between Galois $\kos{S}$-morphisms
is precisely a Galois transformation
between the underlying Galois morphisms
that is compatible with the $\kos{S}$-equivariances
of the right adjoint parts.
In fact, we show a slightly more general statement.

Let $\kos{T}=(\CT,\tensor,\unit)$,
$\kos{T}^{\pr}=(\CT^{\pr},\tensor^{\pr},\unit^{\pr})$
be symmetric monoidal categories
and $\monbase{t}\colon\kos{S}\to \kos{T}$,
$\monbase{t}^{\pr}\colon\kos{S}\to \kos{T}^{\pr}$
be strong monoidal functors.
By the \emph{induced $\kos{S}$-action} $\bcts$ on $\kos{T}$
we mean the functor given by 
$X\bcts z\coloneqq X\tensor \monbase{t}(z)$,
and we write $\bcts^{\pr}$
for the induced $\kos{S}$-action on $\kos{T}^{\pr}$.
We introduce the following terminology.
\begin{itemize}
  \item 
  A \emph{colax $\kos{S}$-tensor functor}
  from $\monbase{t}^{\pr}$ to $\monbase{t}$ is a pair
  $(\mon{f},\what{\mon{f}})$
  consisting of a colax monoidal functor
  $\mon{f}\colon\kos{T}^{\pr}\to \kos{T}$
  and a comonoidal natural transformation
  $\what{\mon{f}}\colon\mon{f}\monbase{t}^{\pr}\Rightarrow\monbase{t}$.
  We omit $\what{\mon{f}}$ from the notation
  and write $\mon{f}\colon\monbase{t}^{\pr}\to \monbase{t}$.
  
  \item
  The \emph{$\kos{S}$-equivariance}
  of a colax $\kos{S}$-tensor functor
  $\mon{f}\colon \monbase{t}^{\pr}\to \monbase{t}$
  is the natural transformation
  $\vecar{\mon{f}}\colon\mon{f}(\slot \bcts^{\pr}\slot)
  \Rightarrow\mon{f}(\slot)\bcts\slot$
  between functors from $\CT^{\pr}\times \CS$ to $\CT$
  whose component at $X^{\pr}\in \CT^{\pr}$, $z\in \CS$ is
  \begin{equation*}
    \xymatrix@C=50pt{
      \mon{f}(X^{\pr}\tensor^{\pr} \monbase{t}^{\pr}(z))
      \ar[r]^-{(\mon{f}_{\tensor})_{X^{\pr},\monbase{t}^{\pr}(z)}}
      &\mon{f}(X^{\pr})\tensor \mon{f}\monbase{t}^{\pr}(z)
      \ar[r]^-{I_{\mon{f}(X^{\pr})}\tensor (\what{\mon{f}})_z}
      &\mon{f}(X^{\pr})\tensor \monbase{t}(z)
      .
    }
  \end{equation*}

  \item
  Let $\mon{f}$, $\mon{g}\colon\monbase{t}^{\pr}\to \monbase{t}$
  be colax $\kos{S}$-tensor functors.
  A \emph{colax $\kos{S}$-tensor natural transformation}
  $\vartheta\colon\mon{f}\Rightarrow\mon{g}$
  is a comonoidal natural transformation
  between the underlying colax monoidal functors
  $\mon{f}$, $\mon{g}\colon\kos{T}^{\pr}\to\kos{T}$
  that satisfies
  $\what{\mon{f}}=\what{\mon{g}}\circ (\vartheta\monbase{t}^{\pr})\colon
  \mon{f}\monbase{t}^{\pr}\Rightarrow\monbase{t}$.
  We say that $\vartheta$ is a
  \emph{colax $\kos{S}$-tensor natural isomorphism}
  if every component of $\vartheta$ is invertible.
\end{itemize}
Let $\kos{T}^{\ppr}$ be a symmetric monoidal category
and $\monbase{t}^{\ppr}\colon\kos{S}\to \kos{T}^{\ppr}$
be a strong monoidal functor,
with the induced $\kos{S}$-action $\bcts^{\ppr}$ on $\kos{T}^{\ppr}$.
The composite 
$\mon{f}\mon{f}^{\pr}\colon\monbase{t}^{\ppr}\to \monbase{t}$
of colax $\kos{S}$-tensor functors
$\mon{f}\colon\monbase{t}^{\pr}\to\monbase{t}$,
$\mon{f}^{\pr}\colon\monbase{t}^{\ppr}\to \monbase{t}^{\pr}$
consists of
$\mon{f}\mon{f}^{\pr}\colon\kos{T}^{\ppr}\to\kos{T}$
and
$\what{\mon{f}\mon{f}^{\pr}}
=\what{\mon{f}}\circ (\mon{f}\what{\mon{f}^{\pr}})
\colon\mon{f}\mon{f}^{\pr}\monbase{t}^{\ppr}
\Rightarrow \monbase{t}$,
and its $\kos{S}$-equivariance is
$(\vecar{\mon{f}\mon{f}^{\pr}})_{X^{\ppr},z}
=(\vecar{\mon{f}})_{\mon{f}^{\pr}(X^{\ppr}),z}
\circ \mon{f}((\vecar{\mon{f}^{\pr}})_{X^{\ppr},z})
\colon \mon{f}\mon{f}^{\pr}(X^{\ppr}\bcts^{\ppr} z)
\to \mon{f}\mon{f}^{\pr}(X^{\ppr})\bcts z$.
The identity colax $\kos{S}$-tensor functor of
$\monbase{t}\colon\kos{S}\to\kos{T}$,
denoted by $\vphantom{\id}\smash[t]{\id}_{\kos{T}\vphantom{!}}^{\vphantom{!}}$,
is the pair consisting of the identity functor of $\kos{T}$
and the identity natural transformation of $\monbase{t}$.

Let
$\kos{f}\colon\kos{t}^{\pr}\mor\kos{t}$
be a Galois $\kos{S}$-morphism
between Galois $\kos{S}$-prekosmoi
$\kos{T}\xmor{\kos{t}}\kos{S}$,
$\kos{T}^{\pr}\xmor{\kos{t}^{\pr}}\kos{S}$.
The right adjoint part $\kos{f}^!$ is a colax $\kos{S}$-tensor functor
from $\kos{t}^!$ to $\kos{t}^{\pr!}$
with $\what{\kos{f}^!}\colon\kos{f}^!\kos{t}^!\cong \kos{t}^{\pr!}$,
and the left adjoint part $\kos{f}_!$ is a colax $\kos{S}$-tensor functor
from $\kos{t}^{\pr!}$ to $\kos{t}^!$
with $\what{\kos{f}_!}\colon\kos{f}_!\kos{t}^{\pr!}\Rightarrow\kos{t}^!$.
One can check that 
the unit $\eta\colon\vphantom{\id}\smash[t]{\id}_{\kos{T}^{\pr}\vphantom{!}}^{\vphantom{!}}\Rightarrow\kos{f}^!\kos{f}_!$
and the counit $\epsilon\colon\kos{f}_!\kos{f}^!\Rightarrow\vphantom{\id}\smash[t]{\id}_{\kos{T}\vphantom{!}}^{\vphantom{!}}$
are colax $\kos{S}$-tensor natural transformations.
Let $\kos{g}\colon\kos{t}^{\pr}\mor\kos{t}$
also be a Galois $\kos{S}$-morphism.
A Galois $\kos{S}$-transformation
$\vartheta\colon\kos{f}\Rightarrow\kos{g}$
is precisely a colax $\kos{S}$-tensor natural transformation
$\vartheta\colon\kos{f}^!\Rightarrow\kos{g}^!$
between the right adjoint parts
$\kos{f}^!$, $\kos{g}^!\colon\kos{t}^!\to \kos{t}^{\pr!}$.

\begin{lemma}
  \label{lem-ActionbaseGal-colaxStnatequiv}
  Let $\kos{T}$, $\kos{T}^{\pr}$
  be symmetric monoidal categories
  and
  $\monbase{t}\colon\kos{S}\to \kos{T}$,
  $\monbase{t}^{\pr}\colon\kos{S}\to \kos{T}^{\pr}$
  be strong monoidal functors.
  Let $\bcts$ and $\bcts^{\pr}$ be the induced $\kos{S}$-actions
  on $\kos{T}$ and $\kos{T}^{\pr}$, respectively.
  Let $\mon{f}$, $\mon{g}\colon\monbase{t}^{\pr}\to \monbase{t}$
  be colax $\kos{S}$-tensor functors
  and $\vartheta\colon\mon{f}\Rightarrow\mon{g}$
  be a comonoidal natural transformation
  between the underlying colax monoidal functors
  $\mon{f}$, $\mon{g}\colon\kos{T}^{\pr}\to \kos{T}$.
  Then $\vartheta$ is a colax $\kos{S}$-tensor natural transformation
  if and only if
  the following diagram commutes for every pair
  $X^{\pr}\in \CT^{\pr}$, $z\in \CS$.
  \begin{equation}
    \label{eq-ActionbaseGal-colaxStnatequiv}
    \vcenter{\hbox{
      \xymatrix@R=15pt@C=35pt{
        \mon{f}(X^{\pr}\bcts^{\pr} z)
        \ar[r]^-{(\vecar{\mon{f}})_{X^{\pr},z}}
        \ar[d]_-{\vartheta_{X^{\pr}\bcts^{\pr} z}}
        &\mon{f}(X^{\pr})\bcts z
        \ar[d]^-{\vartheta_{X^{\pr}}\bcts I_z}
        \\
        \mon{g}(X^{\pr}\bcts^{\pr} z)
        \ar[r]^-{(\vecar{\mon{g}})_{X^{\pr},z}}
        &\mon{g}(X^{\pr})\bcts z
      }
    }}
  \end{equation}
\end{lemma}
\begin{proof}
  If $\vartheta$ is a colax $\kos{S}$-tensor natural transformation,
  the diagram \eqref{eq-ActionbaseGal-colaxStnatequiv}
  commutes by the definitions.
  Conversely, assume that
  \eqref{eq-ActionbaseGal-colaxStnatequiv} commutes
  for every pair $X^{\pr}$, $z$.
  One can check that
  $(\what{\mon{f}})_z
  \colon \mon{f}\monbase{t}^{\pr}(z)
  \to \monbase{t}(z)$
  is equal to
  $(\mon{f}_{\unit}\tensor I_{\monbase{t}(z)})
  \circ (\vecar{\mon{f}})_{\unit^{\pr},z}
  \colon \mon{f}(\unit^{\pr}\tensor^{\pr} \monbase{t}^{\pr}(z))
  \to \unit\tensor \monbase{t}(z)$
  up to coherence isomorphisms.
  Using this description for both $\what{\mon{f}}$ and $\what{\mon{g}}$
  and applying \eqref{eq-ActionbaseGal-colaxStnatequiv} at $X^{\pr}=\unit^{\pr}$,
  we deduce that
  $(\what{\mon{f}})_z=
  (\what{\mon{g}})_z
  \circ \vartheta_{\monbase{t}^{\pr}(z)}
  \colon \mon{f}\monbase{t}^{\pr}(z)
  \to \monbase{t}(z)$
  for every $z$.
  Thus $\vartheta$ is a colax $\kos{S}$-tensor natural transformation.
\end{proof}

Let $c$ be a cocommutative comonoid in $\kos{S}$.
Then $c\otimes \slot$ is a colax $\kos{S}$-tensor functor
from $\vphantom{\id}\smash[t]{\id}_{\kos{S}\vphantom{!}}^{\vphantom{!}}$ to $\vphantom{\id}\smash[t]{\id}_{\kos{S}\vphantom{!}}^{\vphantom{!}}$;
it is equipped with the comonoidal natural transformation
$\what{c\otimes}\colon c\otimes \slot\Rightarrow\vphantom{\id}\smash[t]{\id}_{\kos{S}\vphantom{!}}^{\vphantom{!}}$
whose component at $z\in \CS$ is 
$(\what{c\otimes})_z=e_c\otimes I_z$,
and its $\kos{S}$-equivariance is
$a_{c,x,z}\colon c\otimes (x\otimes z)\cong (c\otimes x)\otimes z$
for $x$, $z\in \CS$.
Recall $\kos{e}_{c\vphantom{!}}^{\vphantom{!}}\colon\kos{e}_{c\vphantom{!}}^{\vphantom{!}}\mor\vphantom{\id}\smash[t]{\id}_{\kos{S}\vphantom{!}}^{\vphantom{!}}$
in \eqref{eq-RepGalobj-basechange}.
The right adjoint part $\kos{e}_{c\vphantom{!}}^!$
is a colax $\kos{S}$-tensor functor
from $\vphantom{\id}\smash[t]{\id}_{\kos{S}\vphantom{!}}^{\vphantom{!}}$ to $\kos{e}_{c\vphantom{!}}^!$
with $\what{\kos{e}_{c\vphantom{!}}^!}\colon\kos{e}_{c\vphantom{!}}^!=\kos{e}_{c\vphantom{!}}^!$,
and the left adjoint part $\kos{e}_{c!}^{\vphantom{!}}$
is a colax $\kos{S}$-tensor functor
from $\kos{e}_{c\vphantom{!}}^!$ to $\vphantom{\id}\smash[t]{\id}_{\kos{S}\vphantom{!}}^{\vphantom{!}}$
with $\what{\vphantom{\kos{e}_{c!}}\smash[t]{\kos{e}_{c!}^{\vphantom{!}}}}\colon\kos{e}_{c!}^{\vphantom{!}}\kos{e}_{c\vphantom{!}}^!\Rightarrow\vphantom{\id}\smash[t]{\id}_{\kos{S}\vphantom{!}}^{\vphantom{!}}$
given by the adjunction counit.
We have $\kos{e}_{c!}^{\vphantom{!}}\kos{e}_{c\vphantom{!}}^!=c\otimes \slot\colon\vphantom{\id}\smash[t]{\id}_{\kos{S}\vphantom{!}}^{\vphantom{!}}\to \vphantom{\id}\smash[t]{\id}_{\kos{S}\vphantom{!}}^{\vphantom{!}}$
as colax $\kos{S}$-tensor functors.

Let $f\colon c^{\pr}\to c$
be a morphism in $\cat{Ens}(\kos{S})$.
The base change along $f$
is a colax $\kos{S}$-tensor functor
$\kos{f}^!\colon\kos{e}_{c\vphantom{!}}^!\to \kos{e}_{c^{\pr}}^!$
with $\what{\kos{f}^!}\colon\kos{f}^!\kos{e}_{c\vphantom{!}}^!=\kos{e}_{c^{\pr}}^!$.
We have a colax $\kos{S}$-tensor natural transformation
$c^{\pr}\otimes \slot\Rightarrow c\otimes \slot$
whose component at $z\in \CS$ is $f\otimes I_z$.

We give special attention to colax $\kos{S}$-tensor functors
from $\vphantom{\id}\smash[t]{\id}_{\kos{S}\vphantom{!}}^{\vphantom{!}}$ to $\vphantom{\id}\smash[t]{\id}_{\kos{S}\vphantom{!}}^{\vphantom{!}}$.
Let $\mon{f}\colon\vphantom{\id}\smash[t]{\id}_{\kos{S}\vphantom{!}}^{\vphantom{!}}\to \vphantom{\id}\smash[t]{\id}_{\kos{S}\vphantom{!}}^{\vphantom{!}}$
be a colax $\kos{S}$-tensor functor.
We denote by
\begin{equation}
  \label{eq-ActionbaseGal-hatardef}
  \hatar{\mon{f}}\colon\mon{f}\Rightarrow\mon{f}(s)\otimes \slot
\end{equation}
the natural transformation whose component at each $z\in \CS$ is
$(\hatar{\mon{f}})_z
=(\vecar{\mon{f}})_{s,z}\circ \mon{f}(\imath_z)\colon
\mon{f}(z)\to \mon{f}(s)\otimes z$.
One can check that
$\hatar{\mon{f}}$ is a colax $\kos{S}$-tensor natural transformation,
and that
$(\hatar{\mon{f}})_{s}
=\jmath_{\mon{f}(s)}\colon \mon{f}(s)\cong \mon{f}(s)\otimes s$.
Let $\mon{g}\colon\vphantom{\id}\smash[t]{\id}_{\kos{S}\vphantom{!}}^{\vphantom{!}}\to \vphantom{\id}\smash[t]{\id}_{\kos{S}\vphantom{!}}^{\vphantom{!}}$
also be a colax $\kos{S}$-tensor functor.
We have
$(\hatar{\mon{g}\mon{f}})_z=
(\vecar{\mon{g}})_{\mon{f}(s),z}\circ \mon{g}((\hatar{\mon{f}})_z)
\colon \mon{g}\mon{f}(z)\to \mon{g}\mon{f}(s)\otimes z$
for $z\in \CS$.
Let $\vartheta\colon\mon{f}\Rightarrow\mon{g}$
be a colax $\kos{S}$-tensor natural transformation.
By taking $\monbase{t}=\monbase{t}^{\pr}=\vphantom{\id}\smash[t]{\id}_{\kos{S}\vphantom{!}}^{\vphantom{!}}$,
$X^{\pr}=s$ in Lemma~\ref{lem-ActionbaseGal-colaxStnatequiv},
we obtain the following commutative diagram for every $z\in \CS$.
\begin{equation}\label{eq-ActionbaseGal-colaxStnathatar}
  \vcenter{\hbox{
    \xymatrix@R=15pt@C=35pt{
      \mon{f}(z)
      \ar[d]_-{\vartheta_z}
      \ar[r]^-{(\hatar{\mon{f}})_z}
      &\mon{f}(s)\otimes z
      \ar[d]^-{\vartheta_{s}\otimes I_z}
      \\%1
      \mon{g}(z)
      \ar[r]^-{(\hatar{\mon{g}})_z}
      &\mon{g}(s)\otimes z
    }
  }}
\end{equation}

We are now ready to prove Proposition~\ref{prop-Galcatoverpk-HomisIsom}.

\begin{proof}[Proof of Proposition~\ref{prop-Galcatoverpk-HomisIsom}]
  Let $\kos{T}\xmor{\kos{t}}\kos{S}$ be a Galois $\kos{S}$-prekosmos
  and let $\omega$, $\omega^{\pr}\colon\vphantom{\id}\smash[t]{\id}_{\kos{S}\vphantom{!}}^{\vphantom{!}}\mor\kos{t}$
  be Galois $\kos{S}$-morphisms.
  Put $p\coloneqq\omega^{\pr!}\omega_!(s)\in\cat{Ens}(\kos{S})$.
  Assume that $\omega$ is a fiber point of $\kos{t}$.
  The right adjoint part of $\omega^{\pr}$
  is a colax $\kos{S}$-tensor functor
  $\omega^{\pr!}\colon\kos{t}^!\to \vphantom{\id}\smash[t]{\id}_{\kos{S}\vphantom{!}}^{\vphantom{!}}$
  whose $\kos{S}$-equivariance
  is a natural isomorphism.
  The left adjoint part of $\omega$
  is a colax $\kos{S}$-tensor functor
  $\omega_!\colon\vphantom{\id}\smash[t]{\id}_{\kos{S}\vphantom{!}}^{\vphantom{!}}\to \kos{t}^!$
  whose $\kos{S}$-equivariance
  is a natural isomorphism by
  Lemma~\ref{lem-ActionbaseGal-projimplystrongequiv}~(2).
  Hence the composite
  $\omega^{\pr!}\omega_!\colon\vphantom{\id}\smash[t]{\id}_{\kos{S}\vphantom{!}}^{\vphantom{!}}\to \vphantom{\id}\smash[t]{\id}_{\kos{S}\vphantom{!}}^{\vphantom{!}}$
  is a colax $\kos{S}$-tensor functor whose $\kos{S}$-equivariance
  is a natural isomorphism.
  Thus, following the notation in \eqref{eq-ActionbaseGal-hatardef},
  we have a colax $\kos{S}$-tensor natural isomorphism
  \begin{equation}\label{eq-ActionbaseGal-hatarupperlowershriek}
    \hatar{\omega^{\pr!}\omega_!}
    \colon\omega^{\pr!}\omega_!\cong p\otimes\slot
  \end{equation}
  whose component at $z\in\CS$ is
  $(\hatar{\omega^{\pr!}\omega_!})_z=
  (\vecar{\omega^{\pr!}\omega_!})_{s,z}\circ 
  \omega^{\pr!}\omega_!(\imath_z)
  \colon \omega^{\pr!}\omega_!(z)\cong \omega^{\pr!}\omega_!(s)\otimes z$.
  This in particular proves (1).

  Let $c$ be a cocommutative comonoid in $\kos{S}$.
  A Galois $\kos{S}$-transformation
  $\vartheta\colon\omega\kos{e}_{c\vphantom{!}}^{\vphantom{!}}\Rightarrow\omega^{\pr}\!\kos{e}_{c\vphantom{!}}^{\vphantom{!}}$
  is precisely a colax $\kos{S}$-tensor natural transformation
  $\vartheta\colon\kos{e}_{c\vphantom{!}}^!\omega^!\Rightarrow \kos{e}_{c\vphantom{!}}^!\omega^{\pr!}$
  between the right adjoint parts.
  By taking mates, we see that such $\vartheta$ are in bijection with
  colax $\kos{S}$-tensor natural transformations
  $\tilde{\vartheta}\colon \kos{e}_{c!}^{\vphantom{!}}\kos{e}_{c\vphantom{!}}^!=c\otimes\slot\Rightarrow\omega^{\pr!}\omega_!$.
  For each $\vartheta$,
  the component of the corresponding $\tilde{\vartheta}$ at $z\in \CS$ is 
  $\vartheta_{\omega_!(z)}\circ (I_c\otimes \eta_z)
  \colon c\otimes z\to \omega^{\pr!}\omega_!(z)$,
  where $\vartheta_{\omega_!(z)}\colon c\otimes\omega^!\omega_!(z)\to \omega^{\pr!}\omega_!(z)$
  is the component of $\vartheta$ at $\omega_!(z)\in\CT$,
  regarded as a morphism in $\CS$.
  For each $\tilde{\vartheta}$,
  the component of the corresponding $\vartheta$ at $X\in \CT$,
  regarded as a morphism in $\CS$,
  is $\omega^{\pr!}(\epsilon_X)\circ \tilde{\vartheta}_{\omega^!(X)}
  \colon c\otimes \omega^!(X)\to \omega^{\pr!}(X)$.
  We illustrate $\vartheta$ and $\tilde{\vartheta}$ as follows.
  \begin{equation*}
    \begin{tikzcd}[row sep=0pt, column sep=15pt]
      \text{ }
      &\id_{\kos{S}}
      \ar[dl,"\kos{e}_{c}^!"']
      \ar[dd,2cell,"{\,\vartheta}"]
      &\text{ }
      &\text{ }
      &\text{ }
      &\id_{\kos{S}}
      \ar[dl,"\kos{e}_{c}^!"']
      \ar[dr,"\omega_!"]
      &\text{ }
      \\%1
      \kos{e}_{c}^!
      &\text{ }
      &\kos{t}^!
      \ar[ul,"\omega^!"']
      \ar[dl,"\omega^{\pr!}"]
      &\text{ }
      &\kos{e}_{c}^!
      \ar[dr,"\kos{e}_{c!}"']
      \ar[rr,2cell,shift right=0.5ex,"{\tilde{\vartheta}\,}"]
      &\text{ }
      &\kos{t}^!
      \ar[dl,"\omega^{\pr!}"]
      \\%2
      \text{ }
      &\id_{\kos{S}}
      \ar[ul,"\kos{e}_{c}^!"]
      &\text{ }
      &\text{ }
      &\text{ }
      &\id_{\kos{S}}
      &\text{ }
    \end{tikzcd}
  \end{equation*}
  On the other hand, colax $\kos{S}$-tensor natural transformations
  $\tilde{\vartheta}\colon c\otimes \slot\Rightarrow\omega^{\pr!}\omega_!$
  correspond bijectively to morphisms
  $g\colon c\to p$ in $\cat{Ens}(\kos{S})$;
  each $\tilde{\vartheta}$ corresponds to
  $\tilde{\vartheta}_{s}\circ \jmath_c\colon c\to p$,
  and each $g$ corresponds to the $\tilde{\vartheta}$
  whose component at $z\in\CS$ is
  $(\hatar{\omega^{\pr!}\omega_!})_z^{-1}\circ (g\otimes I_z)
  \colon c\otimes z\to \omega^{\pr!}\omega_!(z)$.
  One can check this by observing that each $\tilde{\vartheta}$
  satisfies the relation \eqref{eq-ActionbaseGal-colaxStnathatar}
  for every $z\in \CS$,
  which takes the following form.
  \begin{equation*}
    \vcenter{\hbox{
      \xymatrix@R=12pt@C=35pt{
        c\otimes z
        \ar[d]_-{\tilde{\vartheta}_z}
        \ar[r]^-{\jmath_c\otimes I_z}_-{\cong}
        &(c\otimes s)\otimes z
        \ar[d]^-{\tilde{\vartheta}_{s}\otimes I_z}
        \\%1
        \omega^{\pr!}\omega_!(z)
        \ar[r]^-{(\hatar{\omega^{\pr!}\omega_!})_z}_-{\cong}
        &\omega^{\pr!}\omega_!(s)\otimes z
      }
    }}
  \end{equation*}
  Thus we obtain a bijection between
  Galois $\kos{S}$-transformations
  $\vartheta\colon\omega\kos{e}_{c\vphantom{!}}^{\vphantom{!}}\Rightarrow\omega^{\pr}\!\kos{e}_{c\vphantom{!}}^{\vphantom{!}}$
  and morphisms $g\colon c\to p$ in $\cat{Ens}(\kos{S})$.
  For each $X\in\CT$, we define a morphism
  $\gamma_X\colon p\otimes \omega^!(X)\to \omega^{\pr!}(X)$
  in $\CS$ as
  \begin{equation}\label{eq-ActionbaseGal-gammaX}
    \xymatrix@C=12pt{
      \gamma_X\colon
      \omega^{\pr!}\omega_!(s)\otimes \omega^!(X)
      \ar[rrr]^-{(\hatar{\omega^{\pr!}\omega_!})_{\omega^!(X)}^{-1}}_-{\cong}
      &\text{ }
      &\text{ }
      &\omega^{\pr!}\omega_!\omega^!(X)
      \ar[rr]^-{\omega^{\pr!}(\epsilon_X)}
      &\text{ }
      &\omega^{\pr!}(X).
    }
  \end{equation}
  Each $\vartheta$ corresponds to 
  $g=\vartheta_{\omega_!(s)}\circ (I_c\otimes \eta_{s})\circ \jmath_c\colon c\to p$,
  where $\vartheta_{\omega_!(s)}\colon c\otimes\omega^!\omega_!(s)\to \omega^{\pr!}\omega_!(s)$
  is the component of $\vartheta$ at $\omega_!(s)\in \CT$,
  regarded as a morphism in $\CS$.
  Each $g\colon c\to p$ corresponds to the
  $\vartheta\colon\omega\kos{e}_{c\vphantom{!}}^{\vphantom{!}}\Rightarrow\omega^{\pr}\!\kos{e}_{c\vphantom{!}}^{\vphantom{!}}$
  whose component at $X\in \CT$,
  regarded as a morphism in $\CS$,
  is $\gamma_X\circ (g\otimes I_{\omega^!(X)})
  \colon c\otimes \omega^!(X)\to \omega^{\pr!}(X)$.
  The correspondence between
  $\vartheta$ and $g$ is natural in $c$
  and induces an isomorphism
  $\Hom_{\cat{Ens}(\kos{S})}(\slot,p)
  \cong\uHom_{\kos{S}}(\omega,\omega^{\pr})$
  of presheaves on $\cat{Ens}(\kos{S})$.
  We write $\xi\colon\omega\kos{e}_{p\vphantom{!}}^{\vphantom{!}}\Rightarrow\omega^{\pr}\!\kos{e}_{p\vphantom{!}}^{\vphantom{!}}$
  for the universal element
  corresponding to the identity morphism of $p$.
  The component of $\xi$ at $X\in \CT$,
  regarded as a morphism in $\CS$,
  is $\gamma_X$ in \eqref{eq-ActionbaseGal-gammaX}.
  This proves (2).
  
  The equality
  $\uHom_{\kos{S}}(\omega,\omega^{\pr})
  =\uIsom_{\kos{S}}(\omega,\omega^{\pr})$
  holds precisely when the universal element
  $\xi\colon\omega\kos{e}_{p\vphantom{!}}^{\vphantom{!}}\Rightarrow \omega^{\pr}\!\kos{e}_{p\vphantom{!}}^{\vphantom{!}}$
  is invertible;
  for if $\xi$ is invertible,
  since every
  $\vartheta\colon\omega\kos{e}_{c\vphantom{!}}^{\vphantom{!}}\Rightarrow \omega^{\pr}\!\kos{e}_{c\vphantom{!}}^{\vphantom{!}}$
  is the base change of $\xi$ along a unique morphism 
  $g\colon c\to p$ in $\cat{Ens}(\kos{S})$,
  the base change of $\xi^{-1}$ along the same $g$
  is the inverse of $\vartheta$.
  On the other hand,
  $\xi\colon\kos{e}_{p\vphantom{!}}^!\omega^!\Rightarrow \kos{e}_{p\vphantom{!}}^!\omega^{\pr!}$
  is invertible if and only if
  $\kos{e}_{p!}^{\vphantom{!}}\xi\colon p\otimes \omega^!(\slot)
  \Rightarrow p\otimes \omega^{\pr!}(\slot)$
  is invertible,
  since $\kos{e}_{p!}^{\vphantom{!}}\colon\CS_{p\vphantom{!}}^{\vphantom{!}}\to \CS$ is conservative.
  For each $X\in \CT$, we have the
  following explicit description of $\kos{e}_{p!}^{\vphantom{!}}(\xi_X)$.
  \begin{equation}\label{eq-ActionbaseGal-eplowershriekxi}
    \vcenter{\hbox{
      \xymatrix@R=15pt@C=6pt{
        \omega^{\pr!}\omega_!(s)\otimes \omega^!(X)
        \ar[d]_-{(\hatar{\omega^{\pr!}\omega_!})_{\omega^!(X)}^{-1}}^-{\cong}
        \ar[rrrrr]^-{\kos{e}_{p!}(\xi_X)}
        &\text{ }
        &\text{ }
        &\text{ }
        &\text{ }
        &\omega^{\pr!}\omega_!(s)\otimes \omega^{\pr!}(X)
        \\%1
        \omega^{\pr!}\omega_!\omega^!(X)
        \ar[rrr]^-{\omega^{\pr!}\omega_!(\imath_{\omega^!(X)})}_-{\cong}
        &\text{ }
        &\text{ }
        &\omega^{\pr!}\omega_!(s\otimes \omega^!(X))
        \ar[rr]^-{\omega^{\pr!}(\varphi_{s,X})}_-{\cong}
        &\text{ }
        &\omega^{\pr!}(\omega_!(s)\tensor X)
        \ar[u]_-{(\omega^{\pr!}_{\tensor})_{\omega_!(s),X}}^-{\cong}
      }
    }}
  \end{equation}
  Here $\varphi$ is the natural transformation
  associated to $\omega$ in \eqref{eq-Galprekos-projformula}.
  Since $\omega$ is locally connected,
  the morphism $\omega^{\pr!}(\varphi_{s,X})$
  in \eqref{eq-ActionbaseGal-eplowershriekxi}
  is invertible for every $X\in\CT$.
  Thus $\kos{e}_{p!}^{\vphantom{!}}\xi$ is invertible, and so is $\xi$.
  This proves (3).
  Throughout the proof, we have used only that $\omega$ is locally connected,
  which is equivalent to our assumption by Remark~\ref{rem-Galprekos-omegalowershriekconservative}.
\end{proof}

Let $\kos{t}$, $\omega$, $\omega^{\pr}$, and $p$
be as in Proposition~\ref{prop-Galcatoverpk-HomisIsom}
and recall the universal element
$\xi$ of the representation
$\Hom_{\cat{Ens}(\kos{S})}(\slot,p)
\cong\uHom_{\kos{S}}(\omega,\omega^{\pr})$.
We denote by
$\mathring{\gamma}_X\colon p\otimes \omega^{\pr!}(X)\to \omega^!(X)$
the component of the inverse $\xi^{-1}$ at $X\in \CT$,
regarded as a morphism in $\CS$;
it is explicitly given by
\begin{equation}
  \label{eq-ActionbaseGal-xiinverse}
  \xymatrix@C=45pt{
    \mathring{\gamma}_X\colon
    \omega^{\pr!}\omega_!(s)\otimes \omega^{\pr!}(X)
    \ar[r]^-{\kos{e}_{p!}(\xi_X^{-1})}_-{\cong}
    &\omega^{\pr!}\omega_!(s)\otimes \omega^!(X)
    \ar[r]^-{e_p\otimes I_{\omega^!(X)}}
    &\omega^!(X)
    ,
  }
\end{equation}
where $\kos{e}_{p!}^{\vphantom{!}}(\xi_X^{-1})=\kos{e}_{p!}^{\vphantom{!}}(\xi_X)^{-1}$
is the composite, in reverse order, of the inverses of the isomorphisms in 
\eqref{eq-ActionbaseGal-eplowershriekxi}.
Then for each $c\in\cat{Ens}(\kos{S})$
and each 
$\vartheta\colon\omega\kos{e}_{c\vphantom{!}}^{\vphantom{!}}\Rightarrow\omega^{\pr}\!\kos{e}_{c\vphantom{!}}^{\vphantom{!}}$,
the component of $\vartheta^{-1}$ at $X\in \CT$,
regarded as a morphism in $\CS$,
is given by
$\mathring{\gamma}_X\circ (g\otimes I_{\omega^{\pr!}(X)})
\colon c\otimes \omega^{\pr!}(X)\to \omega^!(X)$,
where $g=\vartheta_{\omega_!(s)}\circ (I_c\otimes \eta_{s})\circ \jmath_c\colon c\to p$
is the morphism in $\cat{Ens}(\kos{S})$ corresponding to $\vartheta$.

The following lemma, along with Proposition~\ref{prop-Galcatoverpk-HomisIsom},
will be used in the proof of Theorem~\ref{thm-Galcatoverpk-mainthm}.

\begin{lemma}
  \label{lem-ActionbaseGal-inverseandcomposition}
  Let $\kos{T}\xmor{\kos{t}}\kos{S}$ be a Galois $\kos{S}$-prekosmos
  and $\omega$, $\omega^{\pr}\colon\vphantom{\id}\smash[t]{\id}_{\kos{S}\vphantom{!}}^{\vphantom{!}}\mor\kos{t}$
  be 
  fiber points of $\kos{t}$.
  \begin{enumerate}
    \item 
    The isomorphism
    \begin{equation*}
      \uIsom_{\kos{S}}(\omega,\omega^{\pr})
      \cong \uIsom_{\kos{S}}(\omega^{\pr},\omega)
    \end{equation*}
    of presheaves
    on $\cat{Ens}(\kos{S})$,
    obtained by taking inverses,
    corresponds to the isomorphism
    $\omega^{\pr!}\omega_!(s)
    \cong \omega^!\omega^{\pr}_!(s)$
    between the representing objects,
    given by
    \begin{equation*}
      \xymatrix@C=11pt{
        \omega^{\pr!}\omega_!(s)
        \ar[rrr]^-{I_{\omega^{\pr!}\omega_!(s)}\otimes \eta^{\pr}_{s}}
        &\text{ }
        &\text{ }
        &\omega^{\pr!}\omega_!(s)\otimes \omega^{\pr!}\omega^{\pr}_!(s)
        \ar[rr]^-{\mathring{\gamma}_{\omega^{\pr}_!(s)}}
        &\text{ }
        &\omega^!\omega^{\pr}_!(s)
        .
      }
    \end{equation*}

    \item
    Let $\omega^{\ppr}\colon\vphantom{\id}\smash[t]{\id}_{\kos{S}\vphantom{!}}^{\vphantom{!}}\mor\kos{t}$
    be a Galois $\kos{S}$-morphism.
    The morphism 
    \begin{equation*}
      \uIsom_{\kos{S}}(\omega^{\pr},\omega^{\ppr})
      \times \uIsom_{\kos{S}}(\omega,\omega^{\pr})
      \to \uIsom_{\kos{S}}(\omega,\omega^{\ppr})
    \end{equation*}
    of presheaves
    on $\cat{Ens}(\kos{S})$,
    obtained by the vertical composition,
    corresponds to the morphism
    $\omega^{\ppr!}\omega^{\pr}_!(s)
    \otimes
    \omega^{\pr!}\omega_!(s)
    \to \omega^{\ppr!}\omega_!(s)$
    between the representing objects,
    given by
    \begin{equation*}
      \xymatrix@C=16pt{
        \omega^{\ppr!}\omega^{\pr}_!(s)
        \otimes
        \omega^{\pr!}\omega_!(s)
        \ar[rrr]^-{(\hatar{\omega^{\ppr!}\omega^{\pr}_!})_{\omega^{\pr!}\omega_!(s)}^{-1}}_-{\cong}
        &\text{ }
        &\text{ }
        &\omega^{\ppr!}\omega^{\pr}_!\omega^{\pr!}\omega_!(s)
        \ar[rr]^-{\omega^{\ppr!}(\epsilon^{\pr}_{\omega_!(s)})}
        &\text{ }
        &\omega^{\ppr!}\omega_!(s)
        .
      }
    \end{equation*}
  \end{enumerate}
\end{lemma}
\begin{proof}
  Put $p\coloneqq\omega^{\pr!}\omega_!(s)\in\cat{Ens}(\kos{S})$.
  Since $\omega$ is a fiber point of $\kos{t}$,
  Proposition~\ref{prop-Galcatoverpk-HomisIsom}
  implies that 
  $p$
  represents
  $\uHom_{\kos{S}}(\omega,\omega^{\pr})
  =\uIsom_{\kos{S}}(\omega,\omega^{\pr})$.
  Let $\xi\colon\omega\kos{e}_{p\vphantom{!}}^{\vphantom{!}}\cong \omega^{\pr}\!\kos{e}_{p\vphantom{!}}^{\vphantom{!}}$
  be the universal element of this representation.
  Since $\omega^{\pr}$ is a fiber point of $\kos{t}$ as well,
  by exchanging the roles of $\omega$ and $\omega^{\pr}$
  in Proposition~\ref{prop-Galcatoverpk-HomisIsom},
  we deduce that 
  $\omega^!\omega^{\pr}_!(s)$
  represents
  $\uHom_{\kos{S}}(\omega^{\pr},\omega)
  =\uIsom_{\kos{S}}(\omega^{\pr},\omega)$.
  The isomorphism $p\cong \omega^!\omega^{\pr}_!(s)$
  in (1) is the unique morphism in $\cat{Ens}(\kos{S})$
  that corresponds to
  $\xi^{-1}\colon\omega^{\pr}\!\kos{e}_{p\vphantom{!}}^{\vphantom{!}}\cong \omega\kos{e}_{p\vphantom{!}}^{\vphantom{!}}$,
  and is explicitly given by
  $\mathring{\gamma}_{\omega^{\pr}_!(s)}\circ (I_p\otimes \eta^{\pr}_{s})\circ \jmath_p
  \colon p\cong \omega^!\omega^{\pr}_!(s)$.
  Here $\mathring{\gamma}_{\omega^{\pr}_!(s)}$
  is the component of $\xi^{-1}$ at $\omega^{\pr}_!(s)\in\CT$,
  regarded as a morphism in $\CS$,
  given by substituting $X=\omega^{\pr}_!(s)$
  in \eqref{eq-ActionbaseGal-xiinverse}.
  This proves (1).
  
  Let $c$ be a cocommutative comonoid in $\kos{S}$
  and $\vartheta\colon\omega\kos{e}_{c\vphantom{!}}^{\vphantom{!}}\Rightarrow \omega^{\pr}\!\kos{e}_{c\vphantom{!}}^{\vphantom{!}}$,
  $\vartheta^{\pr}\colon\omega^{\pr}\!\kos{e}_{c\vphantom{!}}^{\vphantom{!}}\Rightarrow \omega^{\ppr}\!\kos{e}_{c\vphantom{!}}^{\vphantom{!}}$
  be Galois $\kos{S}$-transformations,
  with their vertical composite 
  $\vartheta^{\ppr}\coloneqq\vartheta^{\pr}\circ\vartheta\colon
  \omega\kos{e}_{c\vphantom{!}}^{\vphantom{!}}\Rightarrow \omega^{\ppr}\!\kos{e}_{c\vphantom{!}}^{\vphantom{!}}$.
  Note that $\vartheta$, $\vartheta^{\pr}$, and $\vartheta^{\ppr}$
  are all invertible by Proposition~\ref{prop-Galcatoverpk-HomisIsom}~(3).
  We denote by 
  $\tilde{\vartheta}\colon c\otimes\slot\Rightarrow \omega^{\pr!}\omega_!$,
  $\tilde{\vartheta}^{\pr}\colon c\otimes\slot\Rightarrow \omega^{\ppr!}\omega^{\pr}_!$,
  and
  $\tilde{\vartheta}^{\ppr}\colon c\otimes\slot\Rightarrow\omega^{\ppr!}\omega_!$
  the mates of $\vartheta$, $\vartheta^{\pr}$, and $\vartheta^{\ppr}$,
  respectively.
  One can check that the component of $\tilde{\vartheta}^{\ppr}$
  at $z\in\CS$ is given by
  \begin{equation}\label{eq-ActionbaseGal-inverseandcomposition}
    \xymatrix@C=35pt{
      \tilde{\vartheta}^{\ppr}_z\colon
      c\otimes z
      \ar[r]^-{\cp_c\otimes I_z}
      &c\otimes c\otimes z
      \ar[r]^-{(\tilde{\vartheta}^{\pr}\tilde{\vartheta})_z}
      &\omega^{\ppr!}\omega^{\pr}_!\omega^{\pr!}\omega_!(z)
      \ar[r]^-{\omega^{\ppr!}(\epsilon^{\pr}_{\omega_!(z)})}
      &\omega^{\ppr!}\omega_!(z)
      .
    }
  \end{equation}
  Let $g\colon c\to\omega^{\pr!}\omega_!(s)$,
  $g^{\pr}\colon c\to \omega^{\ppr!}\omega^{\pr}_!(s)$,
  and $g^{\ppr}\colon c\to \omega^{\ppr!}\omega_!(s)$
  be the morphisms in $\cat{Ens}(\kos{S})$
  corresponding to $\tilde{\vartheta}$, $\tilde{\vartheta}^{\pr}$, and $\tilde{\vartheta}^{\ppr}$,
  respectively.
  Using the relation \eqref{eq-ActionbaseGal-inverseandcomposition}
  for $z=s$,
  and the relation
  $\tilde{\vartheta}^{\pr}_z=
  (\hatar{\omega^{\ppr!}\omega^{\pr}_!})_z^{-1}\circ (g^{\pr}\otimes I_z)
  \colon c\otimes z\to \omega^{\ppr!}\omega^{\pr}_!(z)$
  for $z=\omega^{\pr!}\omega_!(s)$,
  we deduce that
  $g^{\ppr}=\omega^{\ppr!}(\epsilon^{\pr}_{\omega_!(s)})\circ (\hatar{\omega^{\ppr!}\omega^{\pr}_!})_{\omega^{\pr!}\omega_!(s)}^{-1}\circ (g^{\pr}\otimes g)\circ \cp_c$.
  This proves (2).
\end{proof}

\section{The main theorem}
\label{sec-GalMainthm}

We now prove Theorem~\ref{thm-Galcatoverpk-mainthm},
thereby completing the proof of Theorem~\ref{thm-Galcatoverpk-PrimaryObj}.

\begin{proof}[Proof of Theorem~\ref{thm-Galcatoverpk-mainthm}]
  Let $\kos{T}\xmor{\kos{t}}\kos{S}$
  be a Galois $\kos{S}$-prekosmos
  and $\omega\colon\vphantom{\id}\smash[t]{\id}_{\kos{S}\vphantom{!}}^{\vphantom{!}}\mor\kos{t}$
  be a surjective fiber point of $\kos{t}$.
  Put $\pi\coloneqq\omega^!\omega_!(s)\in\cat{Ens}(\kos{S})$.
  Proposition~\ref{prop-Galcatoverpk-HomisIsom}
  implies that the cocommutative comonoid $\pi$ in $\kos{S}$
  represents the presheaf
  $\uHom_{\kos{S}}(\omega,\omega)
  =\uIsom_{\kos{S}}(\omega,\omega)$
  on $\cat{Ens}(\kos{S})$.
  Let $\xi\colon\omega\kos{e}_{\pi\vphantom{!}}^{\vphantom{!}}\cong \omega\kos{e}_{\pi\vphantom{!}}^{\vphantom{!}}$
  be the universal element of this representation.
  By \eqref{eq-ActionbaseGal-gammaX} for $\omega^{\pr}=\omega$,
  the component of $\xi$ at $X\in \CT$,
  regarded as a morphism in $\CS$,
  is given by
  \begin{equation}\label{eq-GalMainthm-gamma}
    \xymatrix@C=12pt{
      \gamma_X\colon
      \omega^!\omega_!(s)\otimes \omega^!(X)
      \ar[rrr]^-{(\hatar{\omega^!\omega_!})_{\omega^!(X)}^{-1}}_-{\cong}
      &\text{ }
      &\text{ }
      &\omega^!\omega_!\omega^!(X)
      \ar[rr]^-{\omega^!(\epsilon_X)}
      &\text{ }
      &\omega^!(X)
      .
    }
  \end{equation}
  By \eqref{eq-ActionbaseGal-xiinverse} for $\omega^{\pr}=\omega$,
  the component of the inverse $\xi^{-1}$ at $X\in \CT$,
  regarded as a morphism in $\CS$,
  is given by
  \begin{equation}
    \label{eq-GalMainthm-mathringgamma}
    \xymatrix@C=45pt{
      \mathring{\gamma}_X\colon
      \omega^!\omega_!(s)\otimes \omega^!(X)
      \ar[r]^-{\kos{e}_{\pi!}(\xi_X^{-1})}_-{\cong}
      &\omega^!\omega_!(s)\otimes \omega^!(X)
      \ar[r]^-{e_{\pi}\otimes I_{\omega^!(X)}}
      &\omega^!(X)
      .
    }
  \end{equation}
  One can check that the identity Galois $\kos{S}$-transformation
  of $\omega$ corresponds to the morphism
  $\eta_{s}\colon s\to \omega^!\omega_!(s)$
  in $\cat{Ens}(\kos{S})$.
  Combining this with Lemma~\ref{lem-ActionbaseGal-inverseandcomposition}
  for $\omega^{\pr}=\omega^{\ppr}=\omega$,
  we deduce that $\pi$ equipped with
  \begin{equation}\label{eq-GalMainthm-pigroupstr}
    \begin{aligned}
      \pc_{\omega^!\omega_!(s)}
      \colon&\!\!
      \xymatrix@C=12pt{
        \omega^!\omega_!(s)\otimes\omega^!\omega_!(s)
        \ar[rrr]^-{(\hatar{\omega^!\omega_!})_{\omega^!\omega_!(s)}^{-1}}_-{\cong}
        &\text{ }
        &\text{ }
        &\omega^!\omega_!\omega^!\omega_!(s)
        \ar[rr]^-{\omega^!(\epsilon_{\omega_!(s)})}
        &\text{ }
        &\omega^!\omega_!(s)
      }
      \\%1
      u_{\omega^!\omega_!(s)}
      \colon&\!\!
      \xymatrix@C=25pt{
        s
        \ar[r]^-{\eta_{s}}
        &\omega^!\omega_!(s)
      }
      \\%2
      \varsigma_{\omega^!\omega_!(s)}
      \colon&\!\!
      \xymatrix@C=10pt{
        \omega^!\omega_!(s)
        \ar[rrr]^-{I_{\omega^!\omega_!(s)}\otimes \eta_s}
        &\text{ }
        &\text{ }
        &\omega^!\omega_!(s)\otimes \omega^!\omega_!(s)
        \ar[rr]^-{\mathring{\gamma}_{\omega_!(s)}}
        &\text{ }
        &\omega^!\omega_!(s)
      }
    \end{aligned}
  \end{equation}
  is a group object in $\kosEns(\kos{S})$
  that represents the presheaf of groups
  $\uAut_{\kos{S}}(\omega)$ on $\cat{Ens}(\kos{S})$,
  where the morphism $\mathring{\gamma}_{\omega_!(s)}$
  in the description of the antipode
  is given by substituting $X=\omega_!(s)$ in \eqref{eq-GalMainthm-mathringgamma}.
  The left adjoint part $\omega_!$ preserves reflexive coequalizers,
  and the right adjoint part $\omega^!$
  also preserves reflexive coequalizers,
  since $\omega$ is surjective.
  From $\hatar{\omega^!\omega_!}\colon \omega^!\omega_!\cong \pi\otimes\slot$
  we deduce that
  $\pi$ is a Galois object in $\kos{S}$.
  This proves (1).

  It remains to show (2).
  Put $\varpi\coloneqq\omega^!\omega_!$.
  We begin by showing that
  $\hatar{\varpi}\colon \varpi\cong \pi\otimes \slot$
  is an isomorphism of monads on $\CS$.
  Let $z\in\CS$.
  One can check that
  $(\hatar{\varpi})_z\circ \omega^!(\epsilon_{\omega_!(z)})
  =(\pc_{\pi}\otimes I_z)\circ (\hatar{\varpi}\hatar{\varpi})_z
  \colon \varpi\varpi(z)\to \pi\otimes z$,
  using \eqref{eq-ActionbaseGal-colaxStnathatar}
  for $\vartheta=\omega^!\epsilon\omega_!\colon\varpi\varpi\Rightarrow\varpi$,
  the relation
  $(\hatar{\varpi\varpi})_z
  =(\vecar{\varpi})_{\pi,z}\circ \varpi((\hatar{\varpi})_z)
  \colon \varpi\varpi(z)\cong \varpi(\pi)\otimes z$,
  and \eqref{eq-ActionbaseGal-colaxStnatequiv}
  for $\vartheta=\hatar{\varpi}\colon\varpi \cong \pi\otimes\slot$,
  $X^{\pr}=\pi$.
  We also have
  $u_{\pi}\otimes I_z=(\hatar{\varpi})_z\circ \eta_z
  \colon z\to \pi\otimes z$,
  as it is the relation \eqref{eq-ActionbaseGal-colaxStnathatar}
  for $\vartheta=\eta\colon\vphantom{\id}\smash[t]{\id}_{\kos{S}\vphantom{!}}^{\vphantom{!}}\Rightarrow\varpi$.
  Hence $\hatar{\varpi}\colon \varpi\cong \pi\otimes\slot$
  is an isomorphism of monads on $\CS$.

  The assignment
  $X\mapsto (\omega^!(X),\gamma_X)$
  defines a functor
  $\breve{\omega}^!\colon\CT\to\cat{Rep}_{\kos{S}}(\pi)$,
  since $\hatar{\varpi}\colon\varpi\cong \pi\otimes\slot$
  is an isomorphism of monads on $\CS$.
  The coherence isomorphisms
  $\omega^!_{\tensor}$, $\omega^!_{\unit}$ of $\omega^!$
  lift to isomorphisms
  $\breve{\omega}^!_{\tensor}$, $\breve{\omega}^!_{\unit}$
  in $\cat{Rep}_{\kos{S}}(\pi)$,
  since $\gamma\colon\pi\otimes \omega^!(\slot)\Rightarrow\omega^!(\slot)$
  is a comonoidal natural transformation;
  they are the coherence isomorphisms of the strong monoidal functor
  $\breve{\omega}^!\colon\kos{T}\to\kosRep_{\kos{S}}(\pi)$.
  The comonoidal natural isomorphism
  $\what{\omega^!}\colon\omega^!\kos{t}^!\cong \vphantom{\id}\smash[t]{\id}_{\kos{S}\vphantom{!}}^{\vphantom{!}}$
  lifts to a comonoidal natural isomorphism
  $\what{\breve{\omega}^!}\colon\breve{\omega}^!\kos{t}^!\cong \kos{t}_{\pi\vphantom{!}}^!$,
  since $\gamma\colon\pi\otimes \omega^!(\slot)\Rightarrow\omega^!(\slot)$
  is moreover a colax $\kos{S}$-tensor natural transformation;
  it is the structure of the colax $\kos{S}$-tensor functor
  $\breve{\omega}^!\colon\kos{t}^!\to \kos{t}_{\pi\vphantom{!}}^!$,
  omitted from the notation.
  We have an equality
  $\kos{u}_{\pi\vphantom{!}}^!\breve{\omega}^!=\omega^!\colon\kos{t}^!\to \vphantom{\id}\smash[t]{\id}_{\kos{S}\vphantom{!}}^{\vphantom{!}}$
  of colax $\kos{S}$-tensor functors.
  
  We present the rest of the proof in three steps.
  The arguments are essentially the same as
  those in the proof of the crude monadicity theorem
  \cite[Theorem~2.1]{Beck1968}.
  Again, we carry them out in full.
  Given $(x,\gamma_x)\in\cat{Rep}_{\kos{S}}(\pi)$,
  we write
  $\dot{\gamma}_x\coloneqq \gamma_x\circ(\hatar{\varpi})_x
  \colon \omega^!\omega_!(x)\to x$
  for short.

  First, we show that
  $\omega\colon\vphantom{\id}\smash[t]{\id}_{\kos{S}\vphantom{!}}^{\vphantom{!}}\mor\kos{t}$
  extends to a Galois $\kos{S}$-morphism
  $\breve{\omega}\colon\kos{t}_{\pi\vphantom{!}}^{\vphantom{!}}\mor\kos{t}$
  up to an invertible Galois $\kos{S}$-transformation
  $\breve{\omega}\kos{u}_{\pi\vphantom{!}}^{\vphantom{!}}\cong \omega$.
  For each $(x,\gamma_x)\in\cat{Rep}_{\kos{S}}(\pi)$,
  we fix a coequalizer
  $\omega_!(x)\twoheadrightarrow \breve{\omega}_!(x,\gamma_x)$
  in $\CT$ of the reflexive pair
  $\omega_!(\dot{\gamma}_x)$,
  $\epsilon_{\omega_!(x)}
  \colon\omega_!\omega^!\omega_!(x)\to \omega_!(x)$,
  with common section
  $\omega_!(\eta_x)$.
  Then $\breve{\omega}_!$ is left adjoint to $\breve{\omega}^!$.
  The component of the adjunction counit
  $\breve{\epsilon}\colon\breve{\omega}_!\breve{\omega}^!\Rightarrow\vphantom{\id}\smash[t]{\id}_{\CT\vphantom{!}}^{\vphantom{!}}$
  at $X\in \CT$ is the unique morphism
  $\breve{\omega}_!\breve{\omega}^!(X)\to X$
  in $\CT$ whose composite with
  $\omega_!\omega^!(X)
  \twoheadrightarrow \breve{\omega}_!\breve{\omega}^!(X)$ is $\epsilon_X$.
  The component of the adjunction unit
  $\breve{\eta}\colon \id_{\cat{Rep}_{\kos{S}}(\pi)}\Rightarrow \breve{\omega}^!\breve{\omega}_!$
  at $(x,\gamma_x)$,
  seen as a morphism in $\CS$,
  is the composite of
  $\eta_x\colon x\to \omega^!\omega_!(x)$
  followed by
  $\omega^!\omega_!(x)
  \twoheadrightarrow \omega^!\breve{\omega}_!(x,\gamma_x)$.
  Thus 
  $\breve{\omega}=(\breve{\omega}_!\dashv \breve{\omega}^!)\colon
  \kosRep_{\kos{S}}(\pi)\mor\kos{T}$
  is a Galois morphism.
  It is moreover a Galois $\kos{S}$-morphism
  $\breve{\omega}\colon\kos{t}_{\pi\vphantom{!}}^{\vphantom{!}}\mor\kos{t}$,
  where $\what{\breve{\omega}}\colon\kos{t}\breve{\omega}\cong \kos{t}_{\pi\vphantom{!}}^{\vphantom{!}}$
  is the comonoidal natural isomorphism
  $\what{\breve{\omega}^!}\colon\breve{\omega}^!\kos{t}^!\cong \kos{t}_{\pi\vphantom{!}}^!$
  between the right adjoint parts.
  The invertible Galois $\kos{S}$-transformation
  $\breve{\omega}\kos{u}_{\pi\vphantom{!}}^{\vphantom{!}}\cong \omega$
  is the equality
  $\kos{u}_{\pi\vphantom{!}}^!\breve{\omega}^!=\omega^!$
  between the right adjoint parts.

  Next, we show that
  $\breve{\omega}\colon\kos{t}_{\pi\vphantom{!}}^{\vphantom{!}}\mor\kos{t}$ is an embedding.
  Let $(x,\gamma_x)\in\cat{Rep}_{\kos{S}}(\pi)$.
  Since $\omega^!$ preserves reflexive coequalizers,
  $\omega^!\omega_!(x)\twoheadrightarrow \omega^!\breve{\omega}_!(x,\gamma_x)$
  is a coequalizer of the reflexive pair
  $\omega^!\omega_!(\dot{\gamma}_x)$,
  $\omega^!(\epsilon_{\omega_!(x)})
  \colon\omega^!\omega_!\omega^!\omega_!(x)\to \omega^!\omega_!(x)$.
  On the other hand, $\dot{\gamma}_x\colon
  \omega^!\omega_!(x)\to x$
  is a split coequalizer of the same pair
  with respect to the splitting
  $\eta_x$, $\eta_{\omega^!\omega_!(x)}$.
  Thus the underlying morphism in $\CS$
  of $\breve{\eta}_{(x,\gamma_x)}$ is an isomorphism,
  as it is the unique morphism between the two coequalizers.
  Hence $\breve{\eta}$ is a natural isomorphism,
  and $\breve{\omega}\colon\kos{t}_{\pi\vphantom{!}}^{\vphantom{!}}\mor\kos{t}$ is an embedding.

  Finally, we show that
  $\breve{\omega}\colon\kos{t}_{\pi\vphantom{!}}^{\vphantom{!}}\mor\kos{t}$
  is an equivalence.
  It suffices to show that
  $\breve{\epsilon}$ is a natural isomorphism as well.
  Let $X\in\CT$.
  Observe that $\dot{\gamma}_{\omega^!(X)}=\omega^!(\epsilon_X)$.
  Since $\omega^!$ preserves reflexive coequalizers,
  $\omega^!\omega_!\omega^!(X)
  \twoheadrightarrow \omega^!\breve{\omega}_!\breve{\omega}^!(X)$
  is a coequalizer of the reflexive pair
  $\omega^!\omega_!\omega^!(\epsilon_X)$,
  $\omega^!(\epsilon_{\omega_!\omega^!(X)})
  \colon\omega^!\omega_!\omega^!\omega_!\omega^!(X)\to \omega^!\omega_!\omega^!(X)$.
  On the other hand, $\omega^!(\epsilon_X)\colon\omega^!\omega_!\omega^!(X)\to \omega^!(X)$
  is a split coequalizer of the same pair
  with respect to the splitting
  $\eta_{\omega^!(X)}$,
  $\eta_{\omega^!\omega_!\omega^!(X)}$.
  Thus $\omega^!(\breve{\epsilon}_X)$ is an isomorphism,
  as it is the unique morphism between the two coequalizers.
  Since $\omega^!$ is conservative,
  $\breve{\epsilon}_X$ is an isomorphism.
  Hence $\breve{\epsilon}$ is a natural isomorphism,
  and $\breve{\omega}\colon\kos{t}_{\pi\vphantom{!}}^{\vphantom{!}}\mor\kos{t}$
  is an equivalence.
  This completes the proof of Theorem~\ref{thm-Galcatoverpk-mainthm}.
\end{proof}

\begin{remark}
  \label{rem-GalMainthm-reconstruction}
  Let $\pi$ be a Galois object in $\kos{S}$
  and $\kos{u}_{\pi\vphantom{!}}^{\vphantom{!}}\colon\vphantom{\id}\smash[t]{\id}_{\kos{S}\vphantom{!}}^{\vphantom{!}}\mor\kos{t}_{\pi\vphantom{!}}^{\vphantom{!}}$
  be the surjective fiber point of $\kos{t}_{\pi\vphantom{!}}^{\vphantom{!}}$
  in Lemma~\ref{lem-RepGalobj-Reppi}.
  We recover $\pi$ from $\kos{u}_{\pi\vphantom{!}}^{\vphantom{!}}$ as
  the representing object of $\uAut_{\kos{S}}(\kos{u}_{\pi\vphantom{!}}^{\vphantom{!}})$.
  Indeed, by Theorem~\ref{thm-Galcatoverpk-mainthm}~(1),
  $\uAut_{\kos{S}}(\kos{u}_{\pi\vphantom{!}}^{\vphantom{!}})$
  is represented by
  $\kos{u}_{\pi\vphantom{!}}^!\kos{u}_{\pi!}^{\vphantom{!}}(s)
  =\pi\otimes s$
  equipped with
  $\pc_{\pi\otimes s}$,
  $u_{\pi\otimes s}$,
  and $\varsigma_{\pi\otimes s}$
  in \eqref{eq-GalMainthm-pigroupstr}.
  One can check that 
  $\jmath_{\pi}^{-1}\colon \pi\otimes s\cong \pi$
  is an isomorphism of Galois objects in $\kos{S}$.
  Thus we have
  $\Hom_{\cat{Ens}(\kos{S})}(\slot,\pi)
  \cong \uAut_{\kos{S}}(\kos{u}_{\pi\vphantom{!}}^{\vphantom{!}})$;
  each morphism $g\colon c\to \pi$ in $\cat{Ens}(\kos{S})$
  is sent to the $\kos{u}_{\pi\vphantom{!}}^{\vphantom{!}}\kos{e}_{c\vphantom{!}}^{\vphantom{!}}\cong \kos{u}_{\pi\vphantom{!}}^{\vphantom{!}}\kos{e}_{c\vphantom{!}}^{\vphantom{!}}$
  whose component at $(x,\gamma_x)\in\cat{Rep}_{\kos{S}}(\pi)$ is
  $\upgamma_x(g)\colon\kos{e}_{c\vphantom{!}}^!(x)\cong \kos{e}_{c\vphantom{!}}^!(x)$
  in \eqref{eq-RepGalobj-actionofg}.
\end{remark}

\begin{remark}
  \label{rem-GalMainthm-flatfiberpt}
  Let $\kos{T}\xmor{\kos{t}}\kos{S}$ be a Galois $\kos{S}$-prekosmos
  and $\omega\colon\vphantom{\id}\smash[t]{\id}_{\kos{S}\vphantom{!}}^{\vphantom{!}}\mor\kos{t}$
  be a fiber point of $\kos{t}$.
  Assume that $\omega^!$ preserves reflexive coequalizers.
  Following the proof of Theorem~\ref{thm-Galcatoverpk-mainthm},
  we see that
  $\uAut_{\kos{S}}(\omega)$
  is represented by the Galois object
  $\pi\coloneqq\omega^!\omega_!(s)$ in $\kos{S}$,
  that $\omega$ extends to a Galois $\kos{S}$-morphism
  $\breve{\omega}\colon\kos{t}_{\pi\vphantom{!}}^{\vphantom{!}}\mor \kos{t}$
  up to an invertible Galois $\kos{S}$-transformation
  $\breve{\omega}\kos{u}_{\pi\vphantom{!}}^{\vphantom{!}}\cong \omega$,
  and that $\breve{\omega}$ is an embedding.
\end{remark}

We end this section with a remark
that revisits the proof of Theorem~\ref{thm-Galcatoverpk-mainthm}
under a weaker assumption
and relates, for each $X\in\CT$,
the invertibility of the component $\xi_X$ of the universal element $\xi$
to the dualizability of $\breve{\omega}^!(X)$ in $\kosRep_{\kos{S}}(\pi)$.

\begin{remark}
  \label{rem-GalMainthm-dualizableobj}
  Let $\kos{T}\xmor{\kos{t}}\kos{S}$
  be a Galois $\kos{S}$-prekosmos
  and $\omega\colon\vphantom{\id}\smash[t]{\id}_{\kos{S}\vphantom{!}}^{\vphantom{!}}\mor\kos{t}$
  be a Galois $\kos{S}$-morphism.
  Put $\pi\coloneqq\omega^!\omega_!(s)$.
  Assume that
  \begin{equation*}
    \omega^!(\varphi_{\slot,\kos{t}^!(\slot)})
    \colon\omega^!\omega_!(\slot\otimes \omega^!\kos{t}^!(\slot))
    \Rightarrow \omega^!(\omega_!(\slot)\tensor \kos{t}^!(\slot))
  \end{equation*}
  is a natural isomorphism,
  where $\varphi$ is the natural transformation
  associated to $\omega$ in \eqref{eq-Galprekos-projformula}.
  By Lemma~\ref{lem-ActionbaseGal-projimplystrongequiv}~(1)
  for $\kos{f}=\omega$,
  our assumption is equivalent to the condition that
  the $\kos{S}$-equivariance of the colax $\kos{S}$-tensor functor
  $\omega^!\omega_!\colon\vphantom{\id}\smash[t]{\id}_{\kos{S}\vphantom{!}}^{\vphantom{!}}\to \vphantom{\id}\smash[t]{\id}_{\kos{S}\vphantom{!}}^{\vphantom{!}}$
  is a natural isomorphism.
  Then Proposition~\ref{prop-Galcatoverpk-HomisIsom}~(1) and~(2)
  for $\omega^{\pr}=\omega$ still hold,
  and Lemma~\ref{lem-ActionbaseGal-inverseandcomposition}~(2)
  for $\omega^{\pr}=\omega^{\ppr}=\omega$
  also holds with $\uIsom$ replaced by $\uHom$.
  Hence $\pi$,
  equipped with $\pc_{\pi}$, $u_{\pi}$ in \eqref{eq-GalMainthm-pigroupstr},
  is a monoid object in $\kosEns(\kos{S})$
  that represents $\uEnd_{\kos{S}}(\omega)$.
  Let $\xi$ be the universal element of this representation,
  whose component at $X$, regarded as a morphism in $\CS$,
  is $\gamma_X$ in \eqref{eq-GalMainthm-gamma}.
  Using \eqref{eq-ActionbaseGal-eplowershriekxi} for $\omega^{\pr}=\omega$,
  we deduce that for each $X\in \CT$,
  the following are equivalent:
  \begin{itemize}
    \item
    $\omega^!(\varphi_{s,X})\colon
    \omega^!\omega_!(s\otimes \omega^!(X))
    \to \omega^!(\omega_!(s)\tensor X)$
    is an isomorphism in $\CS$.

    \item
    The component of the universal element $\xi$ at $X$
    is an isomorphism in $\CS_{\pi\vphantom{!}}^{\vphantom{!}}$.

    \item
    For each $c\in \cat{Ens}(\kos{S})$
    and each Galois $\kos{S}$-transformation
    $\vartheta\colon\omega\kos{e}_{c\vphantom{!}}^{\vphantom{!}}\Rightarrow \omega\kos{e}_{c\vphantom{!}}^{\vphantom{!}}$,
    the component of $\vartheta$ at $X$
    is an isomorphism in $\CS_{c\vphantom{!}}^{\vphantom{!}}$.
  \end{itemize}
  We similarly define a representation of $\pi$ in $\kos{S}$
  as an object in $\CS$ equipped with a functorial action of $\pi$,
  and identify it with a left $\pi$-module.
  We still have the symmetric monoidal category
  $\kosRep_{\kos{S}}(\pi)=(\cat{Rep}_{\kos{S}}(\pi),\tensor_{\pi},\unit_{\pi})$,
  the strong monoidal functor
  $\kos{u}_{\pi\vphantom{!}}^!\colon\kosRep_{\kos{S}}(\pi)\to \kos{S}$,
  the strong monoidal functor
  $\breve{\omega}^!\colon\kos{T}\to \kosRep_{\kos{S}}(\pi)$
  that sends each $X\in \CT$ to $(\omega^!(X),\gamma_X)$,
  and the equality
  $\kos{u}_{\pi\vphantom{!}}^!\breve{\omega}^!=\omega^!$.
  For each object $X$ in $\CT$, the following are equivalent:
  \begin{itemize}
    \item
    $\breve{\omega}^!(X)$ admits a dual in $\kosRep_{\kos{S}}(\pi)$.

    \item
    $\omega^!(X)$ admits a dual in $\kos{S}$
    and the component of the universal element $\xi$ at $X$
    is an isomorphism in $\CS_{\pi\vphantom{!}}^{\vphantom{!}}$.
  \end{itemize}
  Indeed, we write $x\coloneqq\omega^!(X)$ and $\gamma_x\coloneqq\gamma_X\colon \pi\otimes x\to x$
  so that $\breve{\omega}^!(X)=(x,\gamma_x)$.
  Whenever $\xi_X$ is an isomorphism in $\CS_{\pi\vphantom{!}}^{\vphantom{!}}$,
  we write $\mathring{\gamma}_x\colon \pi\otimes x\to x$
  for the inverse $\xi_X^{-1}$ regarded as a morphism in $\CS$.
  If $(x,\gamma_x)$ admits a dual $((x^{\vee},\gamma^{\vee}),\ev,\cv)$
  in $\kosRep_{\kos{S}}(\pi)$,
  then $(x^{\vee},\ev,\cv)$ is a dual of $x$ in $\kos{S}$,
  and one can check that the inverse $\xi_X^{-1}$ exists
  with $\mathring{\gamma}_x\colon \pi\otimes x\to x$
  given as follows.
  \begin{equation*}
    \xymatrix@R=12pt@C=12pt{
      \pi\otimes x
      \ar[d]_-{\mathring{\gamma}_x}
      \ar[rrrrr]^-{I_{\pi}\otimes \cv\otimes I_x}
      &\text{ }
      &\text{ }
      &\text{ }
      &\text{ }
      &\pi\otimes x\otimes x^{\vee}\otimes x
      \ar[d]^-{\sigma_{\pi,x}\otimes I_{x^{\vee}\otimes x}}_-{\cong}
      \\%1
      x
      &\text{ }
      &x\otimes x^{\vee}\otimes x
      \ar[ll]_-{I_x\otimes \ev}
      &\text{ }
      &\text{ }
      &x\otimes \pi\otimes x^{\vee}\otimes x
      \ar[lll]_-{I_x\otimes \gamma^{\vee}\otimes I_x}
      }
  \end{equation*}
  Conversely, assume that $x$ admits a dual $(x^{\vee},\ev,\cv)$ in $\kos{S}$
  and that $\xi_X$ is an isomorphism.
  Then one can check that
  $((x^{\vee},\gamma^{\vee}),\ev,\cv)$ is a dual of $(x,\gamma_x)$
  in $\kosRep_{\kos{S}}(\pi)$,
  with the left $\pi$-action $\gamma^{\vee}$ of $x^{\vee}$ given as follows.
  \begin{equation*}
    \xymatrix@R=12pt@C=12pt{
      \pi\otimes x^{\vee}
      \ar[d]_-{\gamma^{\vee}}
      \ar[rrrrr]^-{I_{\pi\otimes x^{\vee}}\otimes \cv}
      &\text{ }
      &\text{ }
      &\text{ }
      &\text{ }
      &\pi\otimes x^{\vee}\otimes x\otimes x^{\vee}
      \ar[d]^-{\sigma_{\pi,x^{\vee}}\otimes I_{x\otimes x^{\vee}}}_-{\cong}
      \\%1
      x^{\vee}
      &\text{ }
      &x^{\vee}\otimes x\otimes x^{\vee}
      \ar[ll]_-{\ev\otimes I_{x^{\vee}}}
      &\text{ }
      &\text{ }
      &x^{\vee}\otimes \pi\otimes x\otimes x^{\vee}
      \ar[lll]_-{I_{x^{\vee}}\otimes \mathring{\gamma}_x\otimes I_{x^{\vee}}}
    }
  \end{equation*}
\end{remark}

\section{Co-fpqc torsors under a Galois object}
\label{sec-GalTors}

Throughout this section,
$\pi$ is a Galois object in $\kos{S}$.
A \emph{co-fpqc right} $\pi$\emph{-torsor over $s$}
is a pair $(p,\lambda_p)$
consisting of a cocommutative comonoid $p$ in $\kos{S}$
that is a co-fpqc covering of $s$,
and a morphism
$\lambda_p\colon p\otimes \pi\to p$
in $\cat{Ens}(\kos{S})$
that is a right $\pi$-action,
such that the \emph{shear morphism}
$\tau_p\coloneqq(I_p\otimes \lambda_p)\circ(\cp_p\otimes I_{\pi})
\colon p\otimes \pi\to p\otimes p$
is an isomorphism.
We omit $\lambda_p$ from the notation and write $p$.
The \emph{division morphism} of $p$
is $\dvs_p\coloneqq(e_p\otimes I_{\pi})\circ\tau_p^{-1}\colon p \otimes p\to\pi$.

The pair consisting of $\pi$
and its product morphism $\pc_{\pi}$ in $\cat{Ens}(\kos{S})$
is a co-fpqc right $\pi$-torsor over $s$.
The shear morphism is
$\tau_{\pi}=(I_{\pi}\otimes \pc_{\pi})\circ (\cp_{\pi}\otimes I_{\pi})
\colon\pi\otimes \pi\cong \pi\otimes \pi$,
with its inverse
$\tau_{\pi}^{-1}=(I_{\pi}\otimes \pc_{\pi})\circ (I_{\pi}\otimes \varsigma_{\pi}\otimes I_{\pi})\circ (\cp_{\pi}\otimes I_{\pi})
\colon\pi\otimes \pi\cong \pi\otimes \pi$,
and the division morphism is
$\dvs_{\pi}=\pc_{\pi}\circ (\varsigma_{\pi}\otimes I_{\pi})
\colon\pi\otimes \pi\to \pi$.

A \emph{co-fpqc left} $\pi$\emph{-torsor over $s$}
is a pair consisting of a cocommutative comonoid $p$ in $\kos{S}$
that is a co-fpqc covering of $s$,
and a morphism 
$\varg_p\colon \pi\otimes p\to p$
in $\cat{Ens}(\kos{S})$
that is a left $\pi$-action,
such that
$(\varg_p\otimes I_p)\circ(I_{\pi}\otimes \cp_p)
\colon \pi\otimes p\to p\otimes p$
is an isomorphism.
Every co-fpqc right $\pi$-torsor $p$ over $s$
can be regarded as a co-fpqc left $\pi$-torsor over $s$,
and vice versa.
The two $\pi$-actions are related by
$\varg_p=
\lambda_p\circ \sigma_{\pi,p}\circ (\varsigma_{\pi}\otimes I_p)
\colon \pi\otimes p\to p$
and 
$\lambda_p=
\varg_p\circ \sigma_{p,\pi}\circ (I_p\otimes \varsigma_{\pi})
\colon p\otimes \pi\to p$.

The following lemma collects basic relations
among the morphisms introduced above.
These relations will be used freely in what follows;
we cite individual relations
only when we find it helpful to the reader.

\begin{lemma}
  \label{lem-GalTors-technicalcalculation}
  Let $p$ be a co-fpqc right $\pi$-torsor over $s$.
  The shear morphism
  $\tau_p\colon p\otimes \pi\cong p\otimes p$
  satisfies the following relations:
  \begin{enumerate}
    \item
    $(e_p\otimes I_p)\circ \tau_p=\lambda_p\colon p\otimes \pi\to p$.
    
    \item
    $(I_p\otimes e_p)\circ \tau_p=I_p\otimes e_{\pi}\colon p\otimes \pi\to p$.

    \item
    $(\cp_p\otimes I_p)\circ \tau_p
    =(I_p\otimes \tau_p)\circ (\cp_p\otimes I_{\pi})
    \colon p\otimes \pi\to p\otimes p\otimes p$.

    \item
    $\tau_p\circ (I_p\otimes \pc_{\pi})
    =(I_p\otimes \lambda_p)\circ (\tau_p\otimes I_{\pi})
    \colon p\otimes \pi\otimes \pi\to p\otimes p$.

    \item 
    $\tau_p\circ (I_p\otimes u_{\pi})=\cp_p\colon p\to p\otimes p$.
    
    \item 
    $(I_p\otimes \tau_p)\circ (\tau_p\otimes I_{\pi})
    =(I_{p\otimes p}\otimes \lambda_p)
    \circ
    (I_p\otimes \sigma_{p,p}\otimes I_{\pi})
    \circ
    (\cp_p\otimes I_{p\otimes \pi})
    \circ
    (\tau_p\otimes I_{\pi})
    \circ
    (I_p\otimes \tau_{\pi})
    \colon p\otimes \pi\otimes \pi\cong p\otimes p\otimes p$.
  \end{enumerate}
  The division morphism
  $\dvs_p\colon p \otimes p\to\pi$
  satisfies the following relations:
  \begin{enumerate}
    \item[(7)]
    $(I_p\otimes \dvs_p)\circ (\cp_p\otimes I_p)=\tau_p^{-1}
    \colon p\otimes p\cong p\otimes \pi$.
    
    \item[(8)]
    $\dvs_p\circ (I_p\otimes \lambda_p)=\pc_{\pi}\circ (\dvs_p\otimes I_{\pi})
    \colon p\otimes p\otimes \pi\to \pi$.
    
    \item[(9)]
    $\dvs_p\circ \cp_p=u_{\pi}\circ e_p\colon p\to \pi$.
    
    \item[(10)]
    $\varsigma_{\pi}\circ \dvs_p=\dvs_p\circ \sigma_{p,p}\colon p\otimes p\to \pi$.

    \item[(11)]
    $(\dvs_p\otimes I_p)\circ (I_p\otimes \cp_p)
    =(\varsigma_{\pi}\otimes I_p)
    \circ \sigma_{p,\pi}
    \circ \tau_p^{-1}
    \circ \sigma_{p,p}
    \colon p\otimes p\cong \pi\otimes p$.
  \end{enumerate}
  The induced left $\pi$-action
  $\varg_p=
  \lambda_p\circ \sigma_{\pi,p}\circ (\varsigma_{\pi}\otimes I_p)
  \colon \pi\otimes p\to p$
  satisfies the following relations:
  \begin{enumerate}
    \item[(12)]
    $(\varg_p\otimes I_p)\circ (I_{\pi}\otimes \cp_p)
    =\sigma_{p,p}
    \circ \tau_p\circ \sigma_{\pi,p} \circ (\varsigma_{\pi}\otimes I_p)
    \colon \pi\otimes p\cong p\otimes p$.

    \item[(13)]
    $\varg_p\circ (\dvs_p\otimes I_p)\circ (I_p\otimes \cp_p)
    =I_p\otimes e_p\colon p\otimes p\to p$.

    \item[(14)]
    $\dvs_p\circ (\varg_p\otimes I_p)
    =\pc_{\pi}\circ (I_{\pi}\otimes \dvs_p)\colon \pi\otimes p\otimes p\to \pi$.
  \end{enumerate}
\end{lemma}
\begin{proof}
  (1)--(5) are immediate.
  We obtain (6)
  using $\cp_p\circ\lambda_p
  =(\lambda_p\otimes \lambda_p)\circ \cp_{p\otimes \pi}
  \colon p\otimes \pi\to p\otimes p$.
  (7), (8), and (9) follow from 
  (3), (4), and (5), respectively.
  We obtain (10) by substituting
  $\varsigma_{\pi}=(e_{\pi}\otimes I_{\pi})\circ \tau_{\pi}^{-1}\circ (I_{\pi}\otimes u_{\pi})
  \colon \pi\cong \pi$
  and then using (6) in its equivalent form
  $(I_p\otimes \tau_{\pi}^{-1})\circ (\tau_p^{-1}\otimes I_{\pi})
  =(\tau_p^{-1}\otimes I_{\pi})
  \circ (I_p\otimes \tau_p^{-1})
  \circ (I_{p\otimes p}\otimes \lambda_p)
  \circ (I_p\otimes \sigma_{p,p}\otimes I_{\pi})
  \circ (\cp_p\otimes I_{p\otimes \pi})
  \colon p\otimes p\otimes \pi\cong p\otimes \pi\otimes \pi$.
  We obtain (11) by rewriting $\tau_p^{-1}$ using (7) and then using (10).
  (12) is immediate.
  We see directly from (11) and (12) that
  $(\dvs_p\otimes I_p)\circ (I_p\otimes \cp_p)\colon p\otimes p\cong \pi\otimes p$
  and
  $(\varg_p\otimes I_p)\circ (I_{\pi}\otimes \cp_p)\colon \pi\otimes p\cong p\otimes p$
  are inverse to each other.
  Thus we obtain (13).
  (14) is equivalent to
  $\dvs_p\circ (\varg_p\otimes I_p)\circ (I_{\pi}\otimes \tau_p)
  =\pc_{\pi}\circ (I_{\pi}\otimes e_p\otimes I_{\pi})
  \colon \pi\otimes p\otimes \pi\to \pi$,
  and one can check this using (8), (10), and (12).
\end{proof}

\begin{remark}
  \label{rem-GalTors-cofpqcexactseq}
  Let $p$ be a co-fpqc right $\pi$-torsor over $s$
  and $z\in\CS$.
  Then $e_p\otimes I_z\colon p\otimes z\to z$
  is a coequalizer in $\CS$ of the reflexive pair
  $I_p\otimes e_p\otimes I_z$, $e_p\otimes I_{p\otimes z}\colon 
  p\otimes p\otimes z\to p\otimes z$,
  since the diagram becomes a split coequalizer diagram
  after applying $p\otimes\slot$.
  By composing the pair with
  $\tau_p\otimes I_z\colon p\otimes \pi\otimes z\cong p\otimes p\otimes z$,
  we deduce that $e_p\otimes I_z$
  is a coequalizer in $\CS$ of the reflexive pair
  $\lambda_p\otimes I_z$, $I_p\otimes e_{\pi}\otimes I_z
  \colon p\otimes \pi\otimes z\to p\otimes z$.
  In particular, $e_p\colon p\to s$
  is compatible with the right $\pi$-actions
  when $s$ is equipped with the trivial action
  $I_s\otimes e_{\pi}\colon s\otimes \pi\to s$.
\end{remark}

We denote by $(\cat{Tors}\text{-}\pi)_{s}$
the category of co-fpqc right $\pi$-torsors over $s$,
whose morphisms are $\pi$-equivariant morphisms
between the underlying cocommutative comonoids in $\kos{S}$.
Note that $(\cat{Tors}\text{-}\pi)_{s}$ is a groupoid:
for each morphism
$f\colon q\to p$ in $(\cat{Tors}\text{-}\pi)_{s}$,
one can check that $I_q\otimes f\colon q\otimes q\to q\otimes p$
is inverse to
$\tau_q
\circ (I_q\otimes \dvs_p)
\circ (I_q\otimes f\otimes I_p)
\circ (\cp_q\otimes I_p)
\colon q\otimes p\to q\otimes q$,
and since $q\otimes \slot$ is conservative, $f$ is an isomorphism.
The category of co-fpqc left $\pi$-torsors over $s$
is also a groupoid and is isomorphic to $(\cat{Tors}\text{-}\pi)_{s}$.

Let $\pi^{\pr}$ also be a Galois object in $\kos{S}$.
A \emph{co-fpqc $(\pi^{\pr},\pi)$-bitorsor over $s$}
is a tuple $(p,\lambda_p,\varg^{\pr}_p)$
consisting of a cocommutative comonoid $p$ in $\kos{S}$
that is a co-fpqc covering of $s$,
a right $\pi$-action $\lambda_p$,
and a left $\pi^{\pr}$-action $\varg^{\pr}_p$, such that
(i) $p$ with $\lambda_p$
is a co-fpqc right $\pi$-torsor over $s$;
(ii) $p$ with $\varg^{\pr}_p$
is a co-fpqc left $\pi^{\pr}$-torsor over $s$;
and
(iii) the actions $\lambda_p$ and $\varg^{\pr}_p$ are compatible,
in the sense that
$\varg^{\pr}_p\circ (I_{\pi^{\pr}}\otimes \lambda_p)
=\lambda_p\circ (\varg^{\pr}_p\otimes I_{\pi})
\colon \pi^{\pr}\otimes p\otimes \pi\to p$.
We omit $\lambda_p$, $\varg^{\pr}_p$ from the notation
and write $p$.

\begin{proposition}\label{prop-GalTors-gaugegroup}
  Let $p$ be a co-fpqc right $\pi$-torsor over $s$.
  We have a Galois object $\tw{p}{\pi}$ in $\kos{S}$
  that acts on $p$ from the left
  and makes $p$ a co-fpqc $(\tw{p}{\pi},\pi)$-bitorsor over $s$.
\end{proposition}

We refer to $\tw{p}{\pi}$ as the \emph{gauge group} of $p$.

\begin{proof}[Proof of Proposition~\ref{prop-GalTors-gaugegroup}]
  Let $\varg_p$ be the induced left $\pi$-action on $p$.
  We fix a coequalizer $\mathrm{cq}_p\colon p\otimes p\twoheadrightarrow \tw{p}{\pi}$
  in $\cat{Ens}(\kos{S})$ of the reflexive pair
  $I_p\otimes \varg_p$, $\lambda_p\otimes I_p\colon p\otimes \pi\otimes p\to p\otimes p$,
  with common section $I_p\otimes u_{\pi}\otimes I_p$.
  We write $\tw{p}{\varg_p}\colon \tw{p}{\pi}\otimes p\to p$
  for the unique morphism in $\cat{Ens}(\kos{S})$ that satisfies 
  $\tw{p}{\varg_p}\circ (\mathrm{cq}_p\otimes I_p)
  =\lambda_p\circ (I_p\otimes \dvs_p)\colon p\otimes p\otimes p\to p$;
  one can check that $\tw{p}{\varg_p}$ exists
  using Lemma~\ref{lem-GalTors-technicalcalculation}~(14).
  Moreover, $\lambda_p$ and $\tw{p}{\varg_p}$ are compatible;
  that is, we have 
  $\tw{p}{\varg_p}\circ (I_{\tw{p}{\pi}}\otimes \lambda_p)
  =\lambda_p\circ (\tw{p}{\varg_p}\otimes I_{\pi})
  \colon \tw{p}{\pi}\otimes p\otimes \pi\to p$,
  as this can be checked after precomposing with
  the epimorphism
  $\mathrm{cq}_p\otimes I_{p\otimes \pi}$.
  
  Put
  $\Xi_p\coloneqq
  (I_p\otimes \mathrm{cq}_p)\circ (I_p\otimes \sigma_{p,p})\circ (\cp_p\otimes I_p)
  \colon p\otimes p\to p\otimes \tw{p}{\pi}$
  and
  $\Xi_p^{-1}\coloneqq
  (I_p\otimes \tw{p}{\varg_p})
  \circ (I_p\otimes \sigma_{p,\tw{p}{\pi}})
  \circ (\cp_p\otimes I_{\tw{p}{\pi}})
  \colon p\otimes \tw{p}{\pi}\to p\otimes p$.
  One can check that
  $\tw{p}{\varg_p}
  \circ \sigma_{p,\tw{p}{\pi}}
  \circ (I_p\otimes \mathrm{cq}_p)
  \circ (I_p\otimes \sigma_{p,p})
  \circ (\cp_p\otimes I_p)
  =e_p\otimes I_p
  \colon p\otimes p\to p$
  and that
  $\mathrm{cq}_p
  \circ \sigma_{p,p}
  \circ (I_p\otimes \tw{p}{\varg_p})
  \circ (I_p\otimes \sigma_{p,\tw{p}{\pi}})
  \circ (\cp_p\otimes I_{\tw{p}{\pi}})
  =e_p\otimes I_{\tw{p}{\pi}}
  \colon p\otimes \tw{p}{\pi}\to \tw{p}{\pi}$;
  the latter can be checked after precomposing
  with $I_p\otimes \mathrm{cq}_p$
  using Lemma~\ref{lem-GalTors-technicalcalculation}~(13).
  It thus follows that $\Xi_p$ and $\Xi_p^{-1}$ are inverse to each other.
  From this and the fact that $p$ is a co-fpqc covering of $s$,
  we deduce that $\tw{p}{\pi}$ is a co-fpqc covering of $s$.

  The product $\pc_{\tw{p}{\pi}}$ is the unique morphism in $\cat{Ens}(\kos{S})$
  that satisfies
  $\pc_{\tw{p}{\pi}}
  \circ (I_{\tw{p}{\pi}}\otimes \mathrm{cq}_p)
  =\mathrm{cq}_p
  \circ (\tw{p}{\varg_p}\otimes I_p)
  \colon \tw{p}{\pi}\otimes p\otimes p\to \tw{p}{\pi}$;
  one can check that $\pc_{\tw{p}{\pi}}$ exists
  using the fact that $\lambda_p$ and $\tw{p}{\varg_p}$ are compatible.
  We have
  $\tw{p}{\varg_p}
  \circ (\pc_{\tw{p}{\pi}}\otimes I_p)
  =\tw{p}{\varg_p}
  \circ (I_{\tw{p}{\pi}}\otimes \tw{p}{\varg_p})
  \colon 
  \tw{p}{\pi}\otimes \tw{p}{\pi}\otimes p\to p$,
  as this can be checked after precomposing with
  $I_{\tw{p}{\pi}}\otimes \mathrm{cq}_p\otimes I_p$.
  From this we deduce that $\pc_{\tw{p}{\pi}}$ is associative,
  as the relation can be checked after precomposing with
  $I_{\tw{p}{\pi}\otimes \tw{p}{\pi}}\otimes \mathrm{cq}_p$.

  The unit $u_{\tw{p}{\pi}}$ is the unique morphism in $\cat{Ens}(\kos{S})$
  that satisfies
  $u_{\tw{p}{\pi}}\circ e_p
  =\mathrm{cq}_p\circ \cp_p
  \colon p\to \tw{p}{\pi}$;
  one can check that $u_{\tw{p}{\pi}}$ exists
  using Lemma~\ref{lem-GalTors-technicalcalculation}~(1), (7), and~(13).
  We have
  $\tw{p}{\varg_p}\circ (u_{\tw{p}{\pi}}\otimes I_p)
  =I_p\colon p\to p$,
  as this can be checked after precomposing with $e_p\otimes I_p$.
  From this we deduce that $u_{\tw{p}{\pi}}$ is left unital,
  as the relation can be checked
  after precomposing with $\mathrm{cq}_p$.
  We also have the right unital relation of $u_{\tw{p}{\pi}}$,
  as the relation can be checked after precomposing with
  $(I_{\tw{p}{\pi}}\otimes e_p)\circ (\mathrm{cq}_p\otimes I_p)$
  using Lemma~\ref{lem-GalTors-technicalcalculation}~(13).

  We have shown that $\tw{p}{\pi}$ is a monoid object in $\kosEns(\kos{S})$,
  that it is a co-fpqc covering of $s$,
  that $\tw{p}{\varg_p}$ is a left $\tw{p}{\pi}$-action on $p$,
  and that the two actions $\lambda_p$ and $\tw{p}{\varg_p}$ 
  on $p$ are compatible.
  One can check that
  $(\tw{p}{\varg_p}\otimes I_p)\circ (I_{\tw{p}{\pi}}\otimes \cp_p)
  =\sigma_{p,p}
  \circ \Xi_p^{-1}
  \circ \sigma_{\tw{p}{\pi},p}
  \colon \tw{p}{\pi}\otimes p\cong p\otimes p$,
  from which we deduce that
  $(\tw{p}{\varg_p}\otimes I_p)\circ (I_{\tw{p}{\pi}}\otimes \cp_p)$
  is an isomorphism.
  
  It remains to show that $\tw{p}{\pi}$ has an antipode.
  The antipode of $\tw{p}{\pi}$ exists if and only if 
  $\nu\coloneqq (I_{\tw{p}{\pi}}\otimes \pc_{\tw{p}{\pi}})
  \circ(\cp_{\tw{p}{\pi}}\otimes I_{\tw{p}{\pi}})
  \colon \tw{p}{\pi}\otimes \tw{p}{\pi}\to \tw{p}{\pi}\otimes \tw{p}{\pi}$
  is an isomorphism,
  in which case the antipode is given by
  $\varsigma_{\tw{p}{\pi}}
  \coloneqq(e_{\tw{p}{\pi}}\otimes I_{\tw{p}{\pi}})\circ \nu^{-1}\circ (I_{\tw{p}{\pi}}\otimes u_{\tw{p}{\pi}})
  \colon \tw{p}{\pi}\to \tw{p}{\pi}$.
  Since $p\otimes \slot$ is conservative,
  the existence of $\varsigma_{\tw{p}{\pi}}$
  is also equivalent to the condition that
  $I_p\otimes \nu\colon 
  p\otimes \tw{p}{\pi}\otimes \tw{p}{\pi}\to p\otimes\tw{p}{\pi}\otimes \tw{p}{\pi}$
  is an isomorphism.

  Put
  $\mathring{\Upsilon}_p\coloneqq
  \pc_{\tw{p}{\pi}}
  \circ (\mathrm{cq}_p\otimes I_{\tw{p}{\pi}})
  \circ (\sigma_{p,p}\otimes I_{\tw{p}{\pi}})
  \colon p\otimes p\otimes \tw{p}{\pi}\to \tw{p}{\pi}$.
  One can check that
  $(I_p\otimes \cp_{\tw{p}{\pi}})
  \circ \Xi_p
  =(\Xi_p\otimes I_{\tw{p}{\pi}})
  \circ (I_{p\otimes p}\otimes \mathrm{cq}_p)
  \circ (I_{p\otimes p}\otimes \sigma_{p,p})
  \circ \cp_{p\otimes p}
  \colon p\otimes p
  \to p\otimes \tw{p}{\pi}\otimes \tw{p}{\pi}$,
  and from this we deduce that
  $(\Xi_p\otimes I_{\tw{p}{\pi}})
  \circ (I_{p\otimes p}\otimes \mathring{\Upsilon}_p)
  \circ (\cp_{p\otimes p}\otimes I_{\tw{p}{\pi}})
  =(I_p\otimes \nu)
  \circ (\Xi_p\otimes I_{\tw{p}{\pi}})
  \colon p\otimes p\otimes \tw{p}{\pi} 
  \to p\otimes \tw{p}{\pi}\otimes \tw{p}{\pi}$.
  Hence $\varsigma_{\tw{p}{\pi}}$ exists if and only if
  $(I_{p\otimes p}\otimes \mathring{\Upsilon}_p)
  \circ (\cp_{p\otimes p}\otimes I_{\tw{p}{\pi}})
  \colon p\otimes p\otimes \tw{p}{\pi}
  \to p\otimes p\otimes \tw{p}{\pi}$
  is an isomorphism.

  Put
  $\Upsilon_p\coloneqq
  \pc_{\tw{p}{\pi}}
  \circ (\mathrm{cq}_p\otimes I_{\tw{p}{\pi}})
  \colon p\otimes p\otimes \tw{p}{\pi}
  \to \tw{p}{\pi}$.
  One can check that
  $\pc_{\tw{p}{\pi}}
  \circ (\mathrm{cq}_p\otimes \mathrm{cq}_p)
  \circ (I_p\otimes \cp_p\otimes I_p)
  =\mathrm{cq}_p
  \circ (I_p\otimes e_p\otimes I_p)
  \colon p\otimes p\otimes p
  \to \tw{p}{\pi}$.
  From this, along with the associativity of $\pc_{\tw{p}{\pi}}$
  and the left unital relation of $u_{\tw{p}{\pi}}$,
  we deduce that
  $\Upsilon_p
  \circ (I_{p\otimes p}\otimes \mathring{\Upsilon}_p)
  \circ (\cp_{p\otimes p}\otimes I_{\tw{p}{\pi}})
  =\mathring{\Upsilon}_p
  \circ (I_{p\otimes p}\otimes \Upsilon_p)
  \circ (\cp_{p\otimes p}\otimes I_{\tw{p}{\pi}})
  =e_{p\otimes p}\otimes I_{\tw{p}{\pi}}
  \colon p\otimes p\otimes \tw{p}{\pi}
  \to \tw{p}{\pi}$.
  It thus follows that
  $(I_{p\otimes p}\otimes \Upsilon_p)
  \circ (\cp_{p\otimes p}\otimes I_{\tw{p}{\pi}})$
  and 
  $(I_{p\otimes p}\otimes \mathring{\Upsilon}_p)
  \circ (\cp_{p\otimes p}\otimes I_{\tw{p}{\pi}})$
  are inverse to each other.
  This shows that $\tw{p}{\pi}$ has an antipode.
  We conclude that
  $\tw{p}{\pi}$ is a Galois object in $\kos{S}$
  and that 
  $p$ is a co-fpqc $(\tw{p}{\pi},\pi)$-bitorsor over $s$.
\end{proof}

\section[Fiber points and co-fpqc torsors]{A perfect correspondence between fiber points and co-fpqc torsors}
\label{sec-pFibGalTors}

Recall from Corollary~\ref{cor-Galcatoverpk-pFibgroupoid}
that given any Galois $\kos{S}$-prekosmos
$\kos{T}\xmor{\kos{t}}\kos{S}$,
the category $\cat{pFib}_{\kos{S}}(\kos{t})$ of fiber points of $\kos{t}$
is a groupoid.
The goal of this section is to prove the following theorem.

\begin{theorem}
  \label{thm-pFibGalTors-perfcorres}
  Let $\kos{T}\xmor{\kos{t}}\kos{S}$ be a Galois $\kos{S}$-category.
  For each surjective fiber point $\omega$ of $\kos{t}$,
  we have an adjoint equivalence of groupoids
  $\cat{pFib}_{\kos{S}}(\kos{t})\simeq (\cat{Tors}\text{-}\pi)_{s}$,
  where $\pi\coloneqq\omega^!\omega_!(s)$
  is the Galois object in $\kos{S}$
  that represents $\uAut_{\kos{S}}(\omega)$
  in Theorem~\ref{thm-Galcatoverpk-mainthm}.
\end{theorem}

The proof of Theorem~\ref{thm-pFibGalTors-perfcorres}
is based on the following proposition.

\begin{proposition}
  \label{prop-pFibGalTors-twistedfib}
  Let $\pi$ be a Galois object in $\kos{S}$
  and $p$ be a co-fpqc right $\pi$-torsor over $s$.
  Recall from Lemma~\ref{lem-RepGalobj-Reppi}
  the surjective fiber point
  $\kos{u}_{\pi\vphantom{!}}^{\vphantom{!}}\colon \vphantom{\id}\smash[t]{\id}_{\kos{S}\vphantom{!}}^{\vphantom{!}}\mor\kos{t}_{\pi\vphantom{!}}^{\vphantom{!}}$
  of $\kos{t}_{\pi\vphantom{!}}^{\vphantom{!}}$.
  \begin{enumerate}
    \item 
    We have a fiber point
    $\tw{p}{\omega}\colon \vphantom{\id}\smash[t]{\id}_{\kos{S}\vphantom{!}}^{\vphantom{!}}\mor\kos{t}_{\pi\vphantom{!}}^{\vphantom{!}}$
    of $\kos{t}_{\pi\vphantom{!}}^{\vphantom{!}}$,
    whose right adjoint part $\tw{p}{\omega}^!$
    sends each $(x,\gamma_x)\in \cat{Rep}_{\kos{S}}(\pi)$
    to its $p$-twisted form in $\CS$.
    
    \item
    The presheaf
    $\uIsom_{\kos{S}}(\kos{u}_{\pi\vphantom{!}}^{\vphantom{!}},\tw{p}{\omega})$
    on $\cat{Ens}(\kos{S})$
    is represented by
    $p$.
  \end{enumerate}
  \begin{equation*}
    \begin{tikzcd}[row sep=5pt,column sep=4pt]
      \kos{S}
      \ar[rr,filled,"\tw{p}{\omega}"]
      \ar[ddr,equal]
      &\text{ }
      &\kosRep_{\kos{S}}(\pi)
      \ar[ddl,filled,"\kos{t}_{\pi}"]
      \\%1
      \text{ }
      &\text{ }
      &\text{ }
      \ar[ll,2cell,pos=0.58,shift left=0.3ex,"{\what{\tw{p}{\omega}}}"']
      \\%2
      \text{ }
      &\kos{S}
      &\text{ }
    \end{tikzcd}
    \qquad
    \Hom_{\cat{Ens}(\kos{S})}(\slot,p)
    \cong
    \uIsom_{\kos{S}}(\kos{u}_{\pi\vphantom{!}}^{\vphantom{!}},\tw{p}{\omega})
  \end{equation*}
\end{proposition}

We refer to $\tw{p}{\omega}\colon \vphantom{\id}\smash[t]{\id}_{\kos{S}\vphantom{!}}^{\vphantom{!}}\mor\kos{t}_{\pi\vphantom{!}}^{\vphantom{!}}$
as the \emph{$p$-twisted fiber point} of $\kos{t}_{\pi\vphantom{!}}^{\vphantom{!}}$.

\begin{proof}[Proof of Proposition~\ref{prop-pFibGalTors-twistedfib}]
  We write
  $X=(x,\gamma_x)$, $Y=(y,\gamma_y)$
  for objects in $\cat{Rep}_{\kos{S}}(\pi)$.
  For each $X$, we fix a coequalizer
  $\mathrm{cq}^p_X
  \colon p\otimes x
  \twoheadrightarrow \tw{p}{\omega}^!(X)$
  in $\CS$ of the reflexive pair
  $I_p\otimes \gamma_x$,
  $\lambda_p\otimes I_x
  \colon p\otimes \pi\otimes x
  \to p\otimes x$,
  with common section 
  $I_p\otimes u_{\pi}\otimes I_x$.
  The object $\tw{p}{\omega}^!(X)$ in $\CS$
  is what we mean by the $p$-twisted form of $X$ in (1).
  Then $\tw{p}{\omega}^!$ is a functor,
  sending each morphism $l\colon X\to Y$ in $\cat{Rep}_{\kos{S}}(\pi)$
  to the unique morphism
  $\tw{p}{\omega}^!(l)\colon \tw{p}{\omega}^!(X)\to \tw{p}{\omega}^!(Y)$
  in $\CS$ that satisfies
  $\tw{p}{\omega}^!(l)
  \circ \mathrm{cq}^p_X
  =\mathrm{cq}^p_Y
  \circ (I_p\otimes l)
  \colon p\otimes x
  \to \tw{p}{\omega}^!(Y)$.
  When we regard
  $\mathrm{cq}^p_X
  \colon p\otimes x
  =p\otimes \kos{u}_{\pi\vphantom{!}}^!(X)
  \twoheadrightarrow \tw{p}{\omega}^!(X)$
  as a morphism in the Kleisli category $\CS_{p\vphantom{!}}^{\vphantom{!}}$,
  we denote it by
  $\xi^p_X
  \colon \kos{e}_{p\vphantom{!}}^!\kos{u}_{\pi\vphantom{!}}^!(X)
  \to \kos{e}_{p\vphantom{!}}^!\tw{p}{\omega}^!(X)$.

  We first show that $\xi^p_X$ is an isomorphism in $\CS_{p\vphantom{!}}^{\vphantom{!}}$ for every $X$.
  It suffices to show that
  $\kos{e}_{p!}^{\vphantom{!}}(\xi^p_X)
  =(I_p\otimes \mathrm{cq}^p_X)
  \circ (\cp_p\otimes I_x)
  \colon p\otimes x
  \to p\otimes \tw{p}{\omega}^!(X)$
  is an isomorphism in $\CS$,
  since $\kos{e}_{p!}^{\vphantom{!}}\colon \CS_{p\vphantom{!}}^{\vphantom{!}}\to \CS$ is conservative.
  In the diagram
  \begin{equation*}
    \vcenter{\hbox{
      \xymatrix@R=20pt@C=50pt{
        p\otimes \pi\otimes \pi\otimes x
        \ar@<0.5ex>[r]^-{I_{p\otimes \pi}\otimes \gamma_x}
        \ar@<-0.5ex>[r]_-{I_p\otimes \pc_{\pi}\otimes I_x}
        \ar[d]_-{\tau_p\otimes I_{\pi\otimes x}}^-{\cong}
        &p\otimes \pi\otimes x
        \ar[r]^-{I_p\otimes \gamma_x}
        \ar[d]^-{\tau_p\otimes I_x}_-{\cong}
        &p\otimes x
        \ar@{.>}[d]^-{\exists!\, \kos{e}_{p!}(\xi^p_X)}_-{\cong}
        \\%1
        p\otimes p\otimes \pi\otimes x
        \ar@<0.5ex>[r]^-{I_{p\otimes p}\otimes \gamma_x}
        \ar@<-0.5ex>[r]_-{I_p\otimes \lambda_p\otimes I_x}
        &p\otimes p\otimes x
        \ar@{->>}[r]^-{I_p\otimes \mathrm{cq}^p_X}
        &p\otimes \tw{p}{\omega}^!(X)
      }
    }}
  \end{equation*}
  we see that
  (i) the bottom horizontal morphisms form a coequalizer diagram,
  since $p\otimes \slot$ preserves reflexive coequalizers;
  (ii) the top horizontal morphisms
  form a split coequalizer diagram
  with respect to the splitting
  $I_p\otimes u_{\pi}\otimes I_x$,
  $I_p\otimes u_{\pi}\otimes I_{\pi\otimes x}$;
  (iii) the middle and left vertical isomorphisms
  commute with the horizontal parallel pairs;
  and (iv)
  the unique isomorphism between the two coequalizers
  that makes the right square commute
  is given by
  $(I_p\otimes \mathrm{cq}^p_X)
  \circ (\tau_p\otimes I_x)
  \circ (I_p\otimes u_{\pi}\otimes I_x)
  =\kos{e}_{p!}^{\vphantom{!}}(\xi^p_X)$.
  Thus $\kos{e}_{p!}^{\vphantom{!}}(\xi^p_X)$ is an isomorphism in $\CS$,
  and hence $\xi^p_X$ is an isomorphism in $\CS_{p\vphantom{!}}^{\vphantom{!}}$.
  
  We next show that $\tw{p}{\omega}^!$ is a strong monoidal functor.
  Given $X$ and $Y$, recall the coherence morphism
  $(p\otimes)_{\tensor,x,y}\colon p\otimes x\otimes y\to p\otimes x\otimes p\otimes y$
  of the colax monoidal endofunctor
  $p\otimes \slot$ on $\kos{S}$.
  \begin{itemize}
    \item 
    Let $(\tw{p}{\omega}^!_{\tensor})_{X,Y}
    \colon \tw{p}{\omega}^!(X\tensor_{\pi} Y)
    \to \tw{p}{\omega}^!(X)\otimes \tw{p}{\omega}^!(Y)$
    be the unique morphism in $\CS$ that satisfies 
    $(\tw{p}{\omega}^!_{\tensor})_{X,Y}
    \circ \mathrm{cq}^p_{X\tensor_{\pi}Y}
    =(\mathrm{cq}^p_X\otimes \mathrm{cq}^p_Y)
    \circ (p\otimes)_{\tensor,x,y}$.
    One can check that $(\tw{p}{\omega}^!_{\tensor})_{X,Y}$ exists
    and that
    $(\kos{e}_{p!}^{\vphantom{!}}(\xi^p_X)\otimes I_{\tw{p}{\omega}^!(Y)})
    \circ (I_p\otimes \sigma_{\tw{p}{\omega}^!(Y),x})
    \circ (\kos{e}_{p!}^{\vphantom{!}}(\xi^p_Y)\otimes I_x)
    \circ (I_p\otimes \sigma_{x,y})
    = (I_p\otimes (\tw{p}{\omega}^!_{\tensor})_{X,Y})
    \circ \kos{e}_{p!}^{\vphantom{!}}(\xi^p_{X\tensor_{\pi}Y})
    \colon p\otimes x\otimes y
    \cong p\otimes \tw{p}{\omega}^!(X)\otimes \tw{p}{\omega}^!(Y)$.
    Hence $I_p\otimes (\tw{p}{\omega}^!_{\tensor})_{X,Y}$ is an isomorphism,
    and since $p\otimes \slot$ is conservative,
    we deduce that $(\tw{p}{\omega}^!_{\tensor})_{X,Y}$ is an isomorphism.
    
    \item
    Let 
    $\tw{p}{\omega}^!_{\unit}\colon \tw{p}{\omega}^!(\unit_{\pi})\to s$
    be the unique morphism in $\CS$ that satisfies
    $\tw{p}{\omega}^!_{\unit}
    \circ \mathrm{cq}^p_{\unit_{\pi}}
    =e_p\otimes I_{s}$.
    The morphism $\tw{p}{\omega}^!_{\unit}$ exists and is an isomorphism,
    since $e_p\colon p \to s$
    is a coequalizer of the reflexive pair
    $\lambda_p$, $I_p\otimes e_{\pi}
    \colon p\otimes \pi
    \to p$
    by Remark~\ref{rem-GalTors-cofpqcexactseq}.
  \end{itemize}
  One can check that
  $\tw{p}{\omega}^!_{\tensor}$, $\tw{p}{\omega}^!_{\unit}$
  are compatible with the coherence isomorphisms of
  $\kosRep_{\kos{S}}(\pi)$ and $\kos{S}$,
  using the compatibility relations that
  $(p\otimes)_{\tensor}$, $(p\otimes)_{\unit}$ satisfy.
  Thus $\tw{p}{\omega}^!$ equipped with
  $\tw{p}{\omega}^!_{\tensor}$, $\tw{p}{\omega}^!_{\unit}$
  is a strong monoidal functor.

  We then construct the Galois morphism
  $\tw{p}{\omega}=(\tw{p}{\omega}^!,\tw{p}{\omega}_!,\tw{p}{\eta},\tw{p}{\epsilon})
  \colon \kos{S}\mor \kosRep_{\kos{S}}(\pi)$.
  Let $\varg_p\colon \pi\otimes p\to p$
  be the induced left $\pi$-action on $p$.
  \begin{itemize}
    \item
    The left adjoint part $\tw{p}{\omega}_!$ sends each $z\in\CS$
    to the pair $(p\otimes z,\varg_p\otimes I_z)$.

    \item 
    The component $\tw{p}{\epsilon}_X\colon \tw{p}{\omega}_!\tw{p}{\omega}^!(X)\to X$
    of $\tw{p}{\epsilon}$ at each $X$
    is the unique morphism
    in $\cat{Rep}_{\kos{S}}(\pi)$
    whose underlying morphism
    $\tw{p}{\epsilon}_X\colon p\otimes \tw{p}{\omega}^!(X)\to x$
    in $\CS$ satisfies
    $\tw{p}{\epsilon}_X
    \circ (I_p\otimes \mathrm{cq}^p_X)
    =\gamma_x
    \circ (\dvs_p\otimes I_x)
    \colon p\otimes p\otimes x\to x$;
    the existence of the morphism $\tw{p}{\epsilon}_X$ in $\CS$
    follows from Lemma~\ref{lem-GalTors-technicalcalculation}~(8);
    $\tw{p}{\epsilon}_X$ is a morphism in $\cat{Rep}_{\kos{S}}(\pi)$,
    since the compatibility relation can be checked
    after precomposing with
    $I_{\pi\otimes p}\otimes \mathrm{cq}^p_X$
    using Lemma~\ref{lem-GalTors-technicalcalculation}~(14).

    \item
    The component $\tw{p}{\eta}_z\colon z\to \tw{p}{\omega}^!\tw{p}{\omega}_!(z)$
    of $\tw{p}{\eta}$ at each $z\in \CS$ is the unique morphism in $\CS$
    that satisfies 
    $\tw{p}{\eta}_z
    \circ (e_p\otimes I_z)
    =\mathrm{cq}^p_{\tw{p}{\omega}_!(z)}
    \circ (\cp_p\otimes I_z)
    \colon p\otimes z
    \to \tw{p}{\omega}^!\tw{p}{\omega}_!(z)$;
    one can check that $\tw{p}{\eta}_z$ exists
    using Lemma~\ref{lem-GalTors-technicalcalculation}~(1), (7), and~(13).

    \item
    We have one of the triangle identities
    $\tw{p}{\omega}^!(\tw{p}{\epsilon}_X)
    \circ \tw{p}{\eta}_{\tw{p}{\omega}^!(X)}
    =I_{\tw{p}{\omega}^!(X)}$,
    as this can be checked after precomposing with
    $(e_p\otimes I_{\tw{p}{\omega}^!(X)})
    \circ (I_p\otimes \mathrm{cq}^p_X)$.
    We also have the other triangle identity
    $\tw{p}{\epsilon}_{\tw{p}{\omega}_!(z)}\circ \tw{p}{\omega}_!(\tw{p}{\eta}_z)
    =I_{\tw{p}{\omega}_!(z)}$,
    as this can be checked, as an equality in $\CS$,
    after precomposing with
    $I_p\otimes e_p\otimes I_z$.
  \end{itemize}
  Thus $\tw{p}{\omega}\colon \kos{S}\mor\kosRep_{\kos{S}}(\pi)$
  is a Galois morphism.
  The left adjoint part $\tw{p}{\omega}_!$ is colax monoidal by doctrinal adjunction;
  we have $\kos{u}_{\pi\vphantom{!}}^!\tw{p}{\omega}_!=p\otimes \slot$
  as colax monoidal endofunctors on $\kos{S}$.
  
  We now lift $\tw{p}{\omega}$ to a Galois $\kos{S}$-morphism
  $\vphantom{\id}\smash[t]{\id}_{\kos{S}\vphantom{!}}^{\vphantom{!}}\mor\kos{t}_{\pi\vphantom{!}}^{\vphantom{!}}$.
  For each $z\in \CS$, 
  let
  $(\what{\tw{p}{\omega}^!})_z
  \colon \tw{p}{\omega}^!\kos{t}_{\pi\vphantom{!}}^!(z)
  \to z$
  be the unique morphism in $\CS$
  that satisfies
  $(\what{\tw{p}{\omega}^!})_z
  \circ \mathrm{cq}^p_{\kos{t}_{\pi\vphantom{!}}^!(z)}
  =e_p\otimes I_z$;
  the morphism $(\what{\tw{p}{\omega}^!})_z$ exists and is an isomorphism,
  since $e_p\otimes I_z$ is a coequalizer of
  the reflexive pair $\lambda_p\otimes I_z$, $I_p\otimes e_{\pi}\otimes I_z$
  by Remark~\ref{rem-GalTors-cofpqcexactseq}.
  One can check that $\what{\tw{p}{\omega}^!}
  \colon \tw{p}{\omega}^!\kos{t}_{\pi\vphantom{!}}^!\cong \vphantom{\id}\smash[t]{\id}_{\kos{S}\vphantom{!}}^{\vphantom{!}}$
  is a comonoidal natural isomorphism,
  using the fact that
  $\what{p\otimes}
  \colon p\otimes \slot\Rightarrow \vphantom{\id}\smash[t]{\id}_{\kos{S}\vphantom{!}}^{\vphantom{!}}$
  is a comonoidal natural transformation.
  Thus $\tw{p}{\omega}$
  is a Galois $\kos{S}$-morphism $\vphantom{\id}\smash[t]{\id}_{\kos{S}\vphantom{!}}^{\vphantom{!}}\mor\kos{t}_{\pi\vphantom{!}}^{\vphantom{!}}$,
  with $\what{\tw{p}{\omega}}\colon \kos{t}_{\pi\vphantom{!}}^{\vphantom{!}}\tw{p}{\omega}\cong \vphantom{\id}\smash[t]{\id}_{\kos{S}\vphantom{!}}^{\vphantom{!}}$
  given by $\what{\tw{p}{\omega}^!}\colon \tw{p}{\omega}^!\kos{t}_{\pi\vphantom{!}}^!\cong \vphantom{\id}\smash[t]{\id}_{\kos{S}\vphantom{!}}^{\vphantom{!}}$.
  We have $X\bcts z=X\tensor_{\pi}\kos{t}_{\pi\vphantom{!}}^!(z)=(x\otimes z,\gamma_x\otimes I_z)$.
  One can check that the $\kos{S}$-equivariance
  $(\vecar{\tw{p}{\omega}^!})_{X,z}
  \colon \tw{p}{\omega}^!(X\bcts z)
  \cong \tw{p}{\omega}^!(X)\otimes z$
  of $\tw{p}{\omega}^!$ is the unique morphism that satisfies
  $(\vecar{\tw{p}{\omega}^!})_{X,z}
  \circ \mathrm{cq}^p_{X\bcts z}
  =(\mathrm{cq}^p_X\otimes I_z)
  \circ a_{p,x,z}
  \colon p\otimes (x\otimes z)
  \to \tw{p}{\omega}^!(X)\otimes z$.

  Finally, we show that
  $\tw{p}{\omega}\colon \vphantom{\id}\smash[t]{\id}_{\kos{S}\vphantom{!}}^{\vphantom{!}}\mor\kos{t}_{\pi\vphantom{!}}^{\vphantom{!}}$
  is a fiber point of $\kos{t}_{\pi\vphantom{!}}^{\vphantom{!}}$.
  The left adjoint part $\tw{p}{\omega}_!$ is conservative
  by Remark~\ref{rem-Galprekos-omegalowershriekconservative}.
  It remains to show that 
  $\tw{p}{\omega}$ satisfies the projection formula.
  We write $\tw{p}{\varphi}$ for the natural transformation
  associated to $\tw{p}{\omega}$ in \eqref{eq-Galprekos-projformula}.
  Let $x\in \CS$ and $Z=(z,\gamma_z)\in\cat{Rep}_{\kos{S}}(\pi)$.
  The underlying morphism in $\CS$ of the component
  $\tw{p}{\varphi}_{x,Z}
  \colon \tw{p}{\omega}_!(x\otimes \tw{p}{\omega}^!(Z))
  \to \tw{p}{\omega}_!(x)\tensor_{\pi}Z$
  of $\tw{p}{\varphi}$ at $x$, $Z$ is given by
  $\kos{u}_{\pi\vphantom{!}}^!(\tw{p}{\varphi}_{x,Z})
  =(I_{p\otimes x}\otimes \tw{p}{\epsilon}_Z)
  \circ (p\otimes)_{\tensor,x,\tw{p}{\omega}^!(Z)}
  \colon p\otimes x\otimes \tw{p}{\omega}^!(Z)
  \to p\otimes x\otimes z$.
  One can check that 
  $\tw{p}{\epsilon}_Z
  \circ (I_p\otimes \mathrm{cq}^p_Z)
  \circ (\cp_p\otimes I_z)
  =e_p\otimes I_z
  \colon p\otimes z
  \to z$.
  We also have
  $\mathrm{cq}^p_Z
  \circ (I_p\otimes \tw{p}{\epsilon}_Z)
  \circ (\cp_p\otimes I_{\tw{p}{\omega}^!(Z)})
  =e_p\otimes I_{\tw{p}{\omega}^!(Z)}
  \colon p\otimes \tw{p}{\omega}^!(Z)
  \to \tw{p}{\omega}^!(Z)$,
  as this can be checked after precomposing with 
  $I_p\otimes \mathrm{cq}^p_Z$.
  From these two equalities we deduce that
  $\kos{u}_{\pi\vphantom{!}}^!(\tw{p}{\varphi}_{x,Z})$ is inverse to
  $(I_{p\otimes x}\otimes \mathrm{cq}^p_Z)
  \circ (p\otimes)_{\tensor,x,z}
  \colon p\otimes x\otimes z
  \to p\otimes x\otimes \tw{p}{\omega}^!(Z)$.
  Hence $\kos{u}_{\pi\vphantom{!}}^!\tw{p}{\varphi}$ is a natural isomorphism,
  and so is $\tw{p}{\varphi}$, since $\kos{u}_{\pi\vphantom{!}}^!$ is conservative.
  This proves (1).

  By Proposition~\ref{prop-Galcatoverpk-HomisIsom}
  for $\omega=\kos{u}_{\pi\vphantom{!}}^{\vphantom{!}}$, $\omega^{\pr}=\tw{p}{\omega}$,
  we deduce that
  $\tw{p}{\omega}^!\kos{u}_{\pi!}^{\vphantom{!}}(s)$
  represents
  $\uHom_{\kos{S}}(\kos{u}_{\pi\vphantom{!}}^{\vphantom{!}},\tw{p}{\omega})
  =\uIsom_{\kos{S}}(\kos{u}_{\pi\vphantom{!}}^{\vphantom{!}},\tw{p}{\omega})$.
  We have an isomorphism 
  $i_p
  \colon \tw{p}{\omega}^!\kos{u}_{\pi!}^{\vphantom{!}}(s)
  \cong p$
  in $\cat{Ens}(\kos{S})$,
  which is the unique one that satisfies
  $i_p
  \circ \mathrm{cq}^p_{\kos{u}_{\pi!}^{\vphantom{!}}(s)}
  =\lambda_p
  \circ (I_p\otimes \jmath_{\pi}^{-1})
  \colon p\otimes \pi\otimes s
  \to p$.
  Thus $p$ represents
  $\uIsom_{\kos{S}}(\kos{u}_{\pi\vphantom{!}}^{\vphantom{!}},\tw{p}{\omega})$.
  This proves (2).
  For later use,
  we show that the component of the universal element
  $\kos{u}_{\pi\vphantom{!}}^{\vphantom{!}}\kos{e}_{p\vphantom{!}}^{\vphantom{!}}\cong \tw{p}{\omega}\kos{e}_{p\vphantom{!}}^{\vphantom{!}}$
  of this representation
  at $X$ is given by
  $\xi^p_X
  \colon \kos{e}_{p\vphantom{!}}^!\kos{u}_{\pi\vphantom{!}}^!(X)
  \cong \kos{e}_{p\vphantom{!}}^!\tw{p}{\omega}^!(X)$.
  Indeed, one can check that for $z\in \CS$,
  the isomorphism
  $(\hatar{\tw{p}{\omega}^!\kos{u}_{\pi!}^{\vphantom{!}}})_z
  \colon \tw{p}{\omega}^!\kos{u}_{\pi!}^{\vphantom{!}}(z)
  \cong \tw{p}{\omega}^!\kos{u}_{\pi!}^{\vphantom{!}}(s)\otimes z$
  in \eqref{eq-ActionbaseGal-hatarupperlowershriek} satisfies
  $(\hatar{\tw{p}{\omega}^!\kos{u}_{\pi!}^{\vphantom{!}}})_z
  \circ \mathrm{cq}^p_{\kos{u}_{\pi!}^{\vphantom{!}}(z)}
  =(\mathrm{cq}^p_{\kos{u}_{\pi!}^{\vphantom{!}}(s)}\otimes I_z)
  \circ (I_{p\otimes \pi}\otimes \imath_z)
  \colon p\otimes \pi\otimes z
  \to \tw{p}{\omega}^!\kos{u}_{\pi!}^{\vphantom{!}}(s)\otimes z$.
  From this we deduce that the morphism
  $\gamma_X
  \colon \tw{p}{\omega}^!\kos{u}_{\pi!}^{\vphantom{!}}(s)\otimes x
  \to \tw{p}{\omega}^!(X)$
  in \eqref{eq-ActionbaseGal-gammaX} is equal to 
  $\mathrm{cq}^p_X
  \circ (i_p\otimes I_x)$,
  as this can be checked after precomposing with
  $\mathrm{cq}^p_{\kos{u}_{\pi!}^{\vphantom{!}}(s)}\otimes I_x$.
  Hence the isomorphisms
  $\xi^p_X
  \colon \kos{e}_{p\vphantom{!}}^!\kos{u}_{\pi\vphantom{!}}^!(X)
  \cong \kos{e}_{p\vphantom{!}}^!\tw{p}{\omega}^!(X)$
  in $\CS_{p\vphantom{!}}^{\vphantom{!}}$ constitute the universal element
  $\xi^p\colon \kos{u}_{\pi\vphantom{!}}^{\vphantom{!}}\kos{e}_{p\vphantom{!}}^{\vphantom{!}}\cong \tw{p}{\omega}\kos{e}_{p\vphantom{!}}^{\vphantom{!}}$.
\end{proof}

The following lemma,
along with Proposition~\ref{prop-pFibGalTors-twistedfib},
will be used in the proof of Theorem~\ref{thm-pFibGalTors-perfcorres}.

\begin{lemma}
  \label{lem-pFibGalTors-torsfibmor}
  Let $\pi$ be a Galois object in $\kos{S}$
  and $f\colon p\cong q$
  be an isomorphism of co-fpqc right $\pi$-torsors over $s$.
  We have an isomorphism
  $\tw{f}{\vartheta}\colon \tw{p}{\omega}\cong \tw{q}{\omega}$
  between the twisted fiber points
  $\tw{p}{\omega}$, $\tw{q}{\omega}\colon \vphantom{\id}\smash[t]{\id}_{\kos{S}\vphantom{!}}^{\vphantom{!}}\mor\kos{t}_{\pi\vphantom{!}}^{\vphantom{!}}$
  in Proposition~\ref{prop-pFibGalTors-twistedfib}.
\end{lemma}
\begin{proof}
  For each $X=(x,\gamma_x)\in \cat{Rep}_{\kos{S}}(\pi)$,
  let $\tw{f}{\vartheta}_X\colon \tw{p}{\omega}^!(X)\to \tw{q}{\omega}^!(X)$
  be the unique morphism in $\CS$ that satisfies
  $\tw{f}{\vartheta}_X
  \circ \mathrm{cq}^p_X
  =\mathrm{cq}^q_X
  \circ (f\otimes I_x)
  \colon p\otimes x
  \to \tw{q}{\omega}^!(X)$.
  One can check that the morphisms $\tw{f}{\vartheta}_X$ exist,
  are isomorphisms,
  and constitute an invertible Galois $\kos{S}$-transformation
  $\tw{f}{\vartheta}$ from $\tw{p}{\omega}$ to $\tw{q}{\omega}$.
\end{proof}

\begin{proof}[Proof of Theorem~\ref{thm-pFibGalTors-perfcorres}]
  Let $\kos{T}\xmor{\kos{t}}\kos{S}$ be a Galois $\kos{S}$-category
  and $\omega\colon \vphantom{\id}\smash[t]{\id}_{\kos{S}\vphantom{!}}^{\vphantom{!}}\mor\kos{t}$
  be a surjective fiber point of $\kos{t}$.
  Put $\pi\coloneqq\omega^!\omega_!(s)$
  and $\varpi\coloneqq \omega^!\omega_!$.
  By Theorem~\ref{thm-Galcatoverpk-mainthm},
  $\uAut_{\kos{S}}(\omega)$
  is represented by the Galois object
  $\pi$ in $\kos{S}$,
  and $\omega$ extends to $\breve{\omega}\colon \kos{t}_{\pi\vphantom{!}}^{\vphantom{!}}\simeq\kos{t}$
  up to an invertible Galois $\kos{S}$-transformation
  $\breve{\omega}\kos{u}_{\pi\vphantom{!}}^{\vphantom{!}}\cong \omega$.

  We define the functor
  $\cat{pFib}_{\kos{S}}(\kos{t})\to (\cat{Tors}\text{-}\pi)_{s}$
  as follows.
  Each fiber point $\omega^{\pr}$ of $\kos{t}$
  is sent to $\omega^{\pr!}\omega_!(s)$ equipped with
  $\lambda_{\omega^{\pr!}\omega_!(s)}\coloneqq
  \omega^{\pr!}(\epsilon_{\omega_!(s)})
  \circ (\hatar{\omega^{\pr!}\omega_!})^{-1}_{\pi}
  \colon \omega^{\pr!}\omega_!(s)\otimes \pi
  \to \omega^{\pr!}\omega_!(s)$.
  By Proposition~\ref{prop-Galcatoverpk-HomisIsom},
  $\omega^{\pr!}\omega_!(s)$
  represents 
  $\uHom_{\kos{S}}(\omega,\omega^{\pr})
  =\uIsom_{\kos{S}}(\omega,\omega^{\pr})$.
  Observe that $\lambda_{\omega^{\pr!}\omega_!(s)}$
  is the morphism in Lemma~\ref{lem-ActionbaseGal-inverseandcomposition}~(2)
  with $(\omega,\omega^{\pr},\omega^{\ppr})$
  replaced by $(\omega,\omega,\omega^{\pr})$.
  It thus follows that
  $\lambda_{\omega^{\pr!}\omega_!(s)}$
  is a right $\pi$-action,
  and that the shear morphism
  $\tau_{\omega^{\pr!}\omega_!(s)}$
  is an isomorphism.
  Moreover, $\omega^{\pr!}\omega_!(s)$
  is a co-fpqc covering of $s$;
  this follows from Lemma~\ref{lem-ActionbaseGal-inverseandcomposition}~(1)
  and from the conclusion in the proof
  of Lemma~\ref{lem-Galcatoverpk-pFiblocallyisom}
  that the representing object 
  $\omega^!\omega^{\pr}_!(s)$
  of $\uIsom_{\kos{S}}(\omega^{\pr},\omega)$
  is a co-fpqc covering of $s$.
  Thus $\omega^{\pr!}\omega_!(s)$
  equipped with $\lambda_{\omega^{\pr!}\omega_!(s)}$
  is indeed a co-fpqc right $\pi$-torsor over $s$.
  Each isomorphism $\vartheta\colon \omega^{\pr}\cong \omega^{\ppr}$
  between fiber points of $\kos{t}$
  is sent to the $\pi$-equivariant isomorphism
  $\vartheta_{\omega_!(s)}\colon 
  \omega^{\pr!}\omega_!(s)
  \cong 
  \omega^{\ppr!}\omega_!(s)$
  in $\cat{Ens}(\kos{S})$.
  
  We define the functor
  $(\cat{Tors}\text{-}\pi)_{s}\to \cat{pFib}_{\kos{S}}(\kos{t})$
  as follows.
  Each co-fpqc right $\pi$-torsor $p$ over $s$
  is sent to the fiber point
  $\breve{\omega}\tw{p}{\omega}$ of $\kos{t}$,
  where $\tw{p}{\omega}\colon \vphantom{\id}\smash[t]{\id}_{\kos{S}\vphantom{!}}^{\vphantom{!}}\mor\kos{t}_{\pi\vphantom{!}}^{\vphantom{!}}$
  is the $p$-twisted fiber point of $\kos{t}_{\pi\vphantom{!}}^{\vphantom{!}}$
  in Proposition~\ref{prop-pFibGalTors-twistedfib}.
  Each isomorphism $f\colon p\cong q$
  of co-fpqc right $\pi$-torsors over $s$
  is sent to the isomorphism
  $\breve{\omega}\tw{f}{\vartheta}
  \colon \breve{\omega}\tw{p}{\omega}\cong \breve{\omega}\tw{q}{\omega}$
  of fiber points of $\kos{t}$,
  where $\tw{f}{\vartheta}\colon \tw{p}{\omega}\cong \tw{q}{\omega}$
  is the isomorphism in Lemma~\ref{lem-pFibGalTors-torsfibmor}.

  We show that the functors defined above
  constitute an adjoint equivalence.

  Let $\omega^{\pr}$ be a fiber point of $\kos{t}$.
  Put
  $\varpi^{\pr}\coloneqq\omega^{\pr!}\omega_!$
  and $p\coloneqq \omega^{\pr!}\omega_!(s)$.
  We define an isomorphism 
  $\vartheta^{\omega^{\pr}}
  \colon \breve{\omega}\tw{p}{\omega}
  \cong \omega^{\pr}$
  of fiber points of $\kos{t}$ as follows.
  Recall that for $X\in \CT$,
  we have $\breve{\omega}^!(X)=(\omega^!(X),\gamma_X)$
  with $\gamma_X=\omega^!(\epsilon_X)\circ (\hatar{\varpi})_{\omega^!(X)}^{-1}
  \colon \pi\otimes \omega^!(X)
  \to \omega^!(X)$,
  and $\lambda_p
  =\omega^{\pr!}(\epsilon_{\omega_!(s)})
  \circ (\hatar{\varpi^{\pr}})^{-1}_{\pi}
  \colon p\otimes \pi
  \to p$.
  One can check that in the diagram
  \begin{equation*}
    \vcenter{\hbox{
      \xymatrix@R=20pt@C=15pt{
        p\otimes \pi\otimes \omega^!(X)
        \ar@<0.5ex>[rrr]^-{I_p\otimes \gamma_X}
        \ar@<-0.5ex>[rrr]_-{\lambda_p\otimes I_{\omega^!(X)}}
        \ar[d]_-{(\hatar{\varpi^{\pr}}\hatar{\varpi})_{\omega^!(X)}^{-1}}^-{\cong}
        &\text{ }
        &\text{ }
        &p\otimes \omega^!(X)
        \ar@{->>}[rr]^-{\mathrm{cq}^p_{\breve{\omega}^!(X)}}
        \ar[d]^-{(\hatar{\varpi^{\pr}})_{\omega^!(X)}^{-1}}_-{\cong}
        &\text{ }
        &\tw{p}{\omega}^!\breve{\omega}^!(X)
        \ar@{.>}[d]^-{\exists!\,\vartheta^{\omega^{\pr}}_X}
        \\%1
        \omega^{\pr!}\omega_!\omega^!\omega_!\omega^!(X)
        \ar@<0.5ex>[rrr]^-{\omega^{\pr!}\omega_!\omega^!(\epsilon_X)}
        \ar@<-0.5ex>[rrr]_-{\omega^{\pr!}(\epsilon_{\omega_!\omega^!(X)})}
        &\text{ }
        &\text{ }
        &\omega^{\pr!}\omega_!\omega^!(X)
        \ar[rr]^-{\omega^{\pr!}(\epsilon_X)}
        &\text{ }
        &\omega^{\pr!}(X)
      }
    }}
  \end{equation*}
  the vertical isomorphisms commute with the horizontal parallel pairs,
  using the relation
  $((\hatar{\varpi^{\pr}})_{\pi}\otimes I_{\omega^!(X)})
  \circ (\hatar{\varpi^{\pr}\varpi})_{\omega^!(X)}
  =(\hatar{\varpi^{\pr}}\hatar{\varpi})_{\omega^!(X)}
  \colon \varpi^{\pr}\varpi\omega^!(X)
  \cong p\otimes \pi\otimes \omega^!(X)$
  and \eqref{eq-ActionbaseGal-colaxStnathatar}
  for
  $\omega^{\pr!}\epsilon\omega_!
  \colon \varpi^{\pr}\varpi
  \Rightarrow\varpi^{\pr}$ with $z=\omega^!(X)$;
  thus the right vertical morphism
  $\vartheta^{\omega^{\pr}}_X$
  exists.
  One can also check that the morphisms
  $\vartheta^{\omega^{\pr}}_X$
  constitute a morphism
  $\vartheta^{\omega^{\pr}}
  \colon \breve{\omega}\tw{p}{\omega}
  \Rightarrow \omega^{\pr}$
  of fiber points of $\kos{t}$,
  using the fact that
  $\hatar{\varpi^{\pr}}
  \colon \varpi^{\pr}
  \cong p\otimes \slot$
  is a colax $\kos{S}$-tensor natural isomorphism.
  By Corollary~\ref{cor-Galcatoverpk-pFibgroupoid},
  $\vartheta^{\omega^{\pr}}$ is an isomorphism in $\cat{pFib}_{\kos{S}}(\kos{t})$.
  
  Moreover,
  $\vartheta^{\omega^{\pr}}
  \colon \breve{\omega}\tw{p}{\omega}
  \cong \omega^{\pr}$
  is natural in $\omega^{\pr}$.
  Indeed, let $\omega^{\ppr}$ also be a fiber point of $\kos{t}$
  and $\vartheta\colon \omega^{\pr}\cong\omega^{\ppr}$
  be an isomorphism in $\cat{pFib}_{\kos{S}}(\kos{t})$.
  Put $q\coloneqq\omega^{\ppr!}\omega_!(s)$
  and $f\coloneqq\vartheta_{\omega_!(s)}\colon p\cong q$.
  Then we have
  $\vartheta^{\omega^{\ppr}}_X
  \circ \tw{f}{\vartheta}_{\breve{\omega}^!(X)}
  =\vartheta_X
  \circ \vartheta^{\omega^{\pr}}_X
  \colon \tw{p}{\omega}^!\breve{\omega}^!(X)
  \cong \omega^{\ppr!}(X)$,
  as this can be checked after precomposing with
  $\mathrm{cq}^p_{\breve{\omega}^!(X)}$
  using \eqref{eq-ActionbaseGal-colaxStnathatar}
  for $\vartheta\omega_!\colon \omega^{\pr!}\omega_!\cong\omega^{\ppr!}\omega_!$
  with $z=\omega^!(X)$.

  Let $p$ be a co-fpqc right $\pi$-torsor over $s$.
  We write $f^p\colon \tw{p}{\omega}^!\breve{\omega}^!\omega_!(s)\to p$
  for the unique morphism in $\cat{Ens}(\kos{S})$ that satisfies
  $f^p
  \circ \mathrm{cq}^p_{\breve{\omega}^!\omega_!(s)}
  =\lambda_p
  \colon p\otimes \pi
  \to p$.
  One can check that $f^p$ exists,
  is an isomorphism,
  and is compatible with the right $\pi$-actions;
  the compatibility relation can be checked after precomposing with
  $\mathrm{cq}^p_{\breve{\omega}^!\omega_!(s)}\otimes I_{\pi}$
  using the relation
  $(\hatar{\tw{p}{\omega}^!\breve{\omega}^!\omega_!})_{\pi}
  \circ \mathrm{cq}^p_{\breve{\omega}^!\omega_!(\pi)}
  =(\mathrm{cq}^p_{\breve{\omega}^!\omega_!(s)}\otimes I_{\pi})
  \circ (I_p\otimes (\hatar{\varpi})_{\pi})
  \colon p\otimes \varpi(\pi)
  \to \tw{p}{\omega}^!\breve{\omega}^!\omega_!(s)\otimes \pi$.

  Moreover,
  $f^p\colon \tw{p}{\omega}^!\breve{\omega}^!\omega_!(s)\cong p$
  is natural in $p$.
  Indeed, let $q$ also be a co-fpqc right $\pi$-torsor over $s$
  and $f\colon p\cong q$
  be an isomorphism in $(\cat{Tors}\text{-}\pi)_{s}$.
  We have
  $f^q\circ \tw{f}{\vartheta}_{\breve{\omega}^!\omega_!(s)}
  =f\circ f^p
  \colon \tw{p}{\omega}^!\breve{\omega}^!\omega_!(s)
  \cong q$,
  as this can be checked after precomposing with
  $\mathrm{cq}^p_{\breve{\omega}^!\omega_!(s)}$.

  It remains to check the triangle identities.
  Let $\omega^{\pr}$ be a fiber point of $\kos{t}$
  and put $p\coloneqq\omega^{\pr!}\omega_!(s)$.
  We have one of the triangle identities
  $\vartheta^{\omega^{\pr}}_{\omega_!(s)}\circ (f^p)^{-1}
  =I_p\colon p\to p$,
  as this can be checked after precomposing with $\lambda_p$.

  Let $p$ be a co-fpqc right $\pi$-torsor over $s$.
  Put $\omega^{\pr}\coloneqq\breve{\omega}\tw{p}{\omega}$,
  $q\coloneqq \tw{p}{\omega}^!\breve{\omega}^!\omega_!(s)$,
  and $f\coloneqq (f^p)^{-1}\colon p\cong q$.
  Then we have the other triangle identity 
  $\vartheta^{\omega^{\pr}}\circ (\breve{\omega}\tw{f}{\vartheta})
  =I_{\omega^{\pr}}$;
  for each $X\in \CT$ we have
  $\vartheta^{\omega^{\pr}}_X
  \circ \tw{f}{\vartheta}_{\breve{\omega}^!(X)}
  =I_{\omega^{\pr!}(X)}$,
  as this can be checked after precomposing with
  $\mathrm{cq}^p_{\breve{\omega}^!(X)}
  \circ (\lambda_p\otimes I_{\omega^!(X)})$
  using the relation
  $(\hatar{\tw{p}{\omega}^!\breve{\omega}^!\omega_!})_{\omega^!(X)}
  \circ \mathrm{cq}^p_{\breve{\omega}^!\omega_!\omega^!(X)}
  =(\mathrm{cq}^p_{\breve{\omega}^!\omega_!(s)}\otimes I_{\omega^!(X)})
  \circ (I_p\otimes (\hatar{\varpi})_{\omega^!(X)})
  \colon p\otimes \varpi\omega^!(X)
  \to \tw{p}{\omega}^!\breve{\omega}^!\omega_!(s)\otimes \omega^!(X)$.
  This completes the proof of Theorem~\ref{thm-pFibGalTors-perfcorres}.
\end{proof}

We end this section with a lemma
that exhibits the relation between
the gauge group and the twisted fiber point
associated to a co-fpqc right torsor over $s$
under a Galois object in $\kos{S}$.

\begin{lemma}
  \label{lem-pFibGalTors-twppitwpomega}
  Let $\pi$ be a Galois object in $\kos{S}$
  and
  $p$ be a co-fpqc right $\pi$-torsor over $s$.
  Recall from Proposition~\ref{prop-GalTors-gaugegroup}
  the gauge group $\tw{p}{\pi}$ of $p$,
  and from Proposition~\ref{prop-pFibGalTors-twistedfib}
  the $p$-twisted fiber point
  $\tw{p}{\omega}$ of $\kos{t}_{\pi\vphantom{!}}^{\vphantom{!}}$.
  \begin{enumerate}
    \item
    The presheaf of groups
    $\uAut_{\kos{S}}(\tw{p}{\omega})$
    on $\cat{Ens}(\kos{S})$
    is represented by $\tw{p}{\pi}$.
    
    \item
    The fiber point $\tw{p}{\omega}$
    extends to an equivalence
    $\kos{t}_{\tw{p}{\pi}\vphantom{!}}^{\vphantom{!}}\simeq\kos{t}_{\pi\vphantom{!}}^{\vphantom{!}}$
    of Galois $\kos{S}$-prekosmoi
    up to an invertible Galois $\kos{S}$-transformation.
  \end{enumerate}
  \begin{equation*}
    \Hom_{\cat{Ens}(\kos{S})}(\slot,\tw{p}{\pi})
    \cong
    \uAut_{\kos{S}}(\tw{p}{\omega})
    \qquad\quad
    \vcenter{\hbox{
      \begin{tikzcd}[row sep=5pt,column sep=6pt]
        \text{ }
        &\id_{\kos{S}}
        \ar[ddr,filled,"{\tw{p}{\omega}}"]
        \ar[ddl,filled,"{\kos{u}_{\tw{p}{\pi}}}"']
        &\text{ }
        \\%1
        \text{ }
        \ar[rr,2cell,shift right=0.7ex,pos=0.5]
        &\text{ }
        &\text{ }
        \\%2
        \kos{t}_{\tw{p}{\pi}}
        \ar[rr,filled,"\simeq"]
        &\text{ }
        &\kos{t}_{\pi}
      \end{tikzcd}
    }}
  \end{equation*}
\end{lemma}
\begin{proof}
  By Corollary~\ref{cor-Galcatoverpk-pFibsurjective},
  the twisted fiber point
  $\tw{p}{\omega}$ of $\kos{t}_{\pi\vphantom{!}}^{\vphantom{!}}$ is surjective.
  By Theorem~\ref{thm-Galcatoverpk-mainthm}
  for $\kos{t}=\kos{t}_{\pi\vphantom{!}}^{\vphantom{!}}$, $\omega=\tw{p}{\omega}$,
  we deduce that $\uAut_{\kos{S}}(\tw{p}{\omega})$
  is represented by the Galois object
  $\tw{p}{\omega}^!\tw{p}{\omega}_!(s)$ in $\kos{S}$,
  and that $\tw{p}{\omega}$ extends to
  $\tw{p}{\breve{\omega}}\colon \kos{t}_{\tw{p}{\omega}^!\tw{p}{\omega}_!(s)}\simeq \kos{t}_{\pi\vphantom{!}}^{\vphantom{!}}$
  up to an invertible Galois $\kos{S}$-transformation
  $\tw{p}{\breve{\omega}}\kos{u}_{\tw{p}{\omega}^!\tw{p}{\omega}_!(s)}\cong \tw{p}{\omega}$.
  The product, unit, and antipode of $\tw{p}{\omega}^!\tw{p}{\omega}_!(s)$
  are given as in \eqref{eq-GalMainthm-pigroupstr}.
  Let
  $i_{\tw{p}{\pi}}\colon \tw{p}{\omega}^!\tw{p}{\omega}_!(s)\cong \tw{p}{\pi}$
  be the unique isomorphism in $\cat{Ens}(\kos{S})$ that satisfies
  $i_{\tw{p}{\pi}}
  \circ \mathrm{cq}^p_{\tw{p}{\omega}_!(s)}
  =\mathrm{cq}_p
  \circ (I_p\otimes \jmath_p^{-1})
  \colon p\otimes p\otimes s
  \to \tw{p}{\pi}$.
  Then the isomorphism
  $(\hatar{\tw{p}{\omega}^!\tw{p}{\omega}_!})_z
  \colon \tw{p}{\omega}^!\tw{p}{\omega}_!(z)
  \cong \tw{p}{\omega}^!\tw{p}{\omega}_!(s)\otimes z$
  in \eqref{eq-ActionbaseGal-hatarupperlowershriek}
  satisfies
  $(\hatar{\tw{p}{\omega}^!\tw{p}{\omega}_!})_z
  \circ \mathrm{cq}^p_{\tw{p}{\omega}_!(z)}
  =(i_{\tw{p}{\pi}}^{-1}\otimes I_z)
  \circ (\mathrm{cq}_p\otimes I_z)
  \colon p\otimes p\otimes z
  \to \tw{p}{\omega}^!\tw{p}{\omega}_!(s)\otimes z$,
  and using this one can check that
  $i_{\tw{p}{\pi}}
  \colon \tw{p}{\omega}^!\tw{p}{\omega}_!(s)
  \cong \tw{p}{\pi}$
  is an isomorphism of Galois objects in $\kos{S}$.
  This proves (1) and (2).
\end{proof}

\section{Fundamental morphisms over a Galois prekosmos}
\label{sec-fundmoroverGalpk}

In this section,
we define the $(2,1)$-category
$\mathsf{pGalCat}(\kos{S})$
whose $0$-cells are Galois $\kos{S}$-categories;
we begin by describing the $1$-cells.

\begin{definition}
  \label{def-fundmoroverGalpk-fundmor}
  Let
  $\kos{T}\xmor{\kos{t}}\kos{S}$,
  $\kos{T}^{\pr}\xmor{\kos{t}^{\pr}}\kos{S}$
  be Galois $\kos{S}$-categories.
  A \emph{fundamental $\kos{S}$-morphism},
  or a \emph{fundamental morphism over $\kos{S}$},
  from $\kos{t}^{\pr}$ to $\kos{t}$
  is a Galois $\kos{S}$-morphism
  $\kos{f}\colon \kos{t}^{\pr}\mor\kos{t}$
  for which there exists a surjective fiber point
  $\omega^{\pr}\colon \vphantom{\id}\smash[t]{\id}_{\kos{S}\vphantom{!}}^{\vphantom{!}}\mor\kos{t}^{\pr}$ of $\kos{t}^{\pr}$
  such that the composite $\kos{f}\omega^{\pr}\colon \vphantom{\id}\smash[t]{\id}_{\kos{S}\vphantom{!}}^{\vphantom{!}}\mor\kos{t}$
  is a surjective fiber point of $\kos{t}$.
  \begin{equation*}
    \begin{tikzcd}[row sep=4pt,column sep=7pt]
      \text{ }
      &\id_{\kos{S}}
      \ar[ddr,filled,"{\kos{f}\omega^{\pr}}"]
      \ar[ddl,filled,"{\omega^{\pr}}"']
      &\text{ }
      \\%1
      \text{ }
      &\text{ }
      &\text{ }
      \\%2
      \kos{t}^{\pr}
      \ar[rr,filled,"{\kos{f}}"]
      &\text{ }
      &\kos{t}
    \end{tikzcd}
  \end{equation*}
\end{definition}

\begin{lemma}
  \label{lem-fundmoroverGalpk-fundmor}
  Let 
  $f\colon \pi^{\pr}\to\pi$
  be a morphism of Galois objects in $\kos{S}$.
  We have a fundamental $\kos{S}$-morphism
  $\kos{f}\colon \kos{t}_{\pi^{\pr}\vphantom{!}}^{\vphantom{!}}\mor\kos{t}_{\pi\vphantom{!}}^{\vphantom{!}}$
  and an invertible Galois $\kos{S}$-transformation
  $\kos{f}\kos{u}_{\pi^{\pr}\vphantom{!}}^{\vphantom{!}}\cong \kos{u}_{\pi\vphantom{!}}^{\vphantom{!}}$.
  \begin{equation*}
    \begin{tikzcd}[row sep=5pt,column sep=6pt]
      \text{ }
      &\id_{\kos{S}}
      \ar[ddr,filled,"{\kos{u}_{\pi}}"]
      \ar[ddl,filled,"{\kos{u}_{\pi^{\pr}}}"']
      &\text{ }
      \\%1
      \text{ }
      \ar[rr,2cell,shift right=0.7ex,pos=0.5]
      &\text{ }
      &\text{ }
      \\%2
      \kos{t}_{\pi^{\pr}}
      \ar[rr,filled,"{\kos{f}}"]
      &\text{ }
      &\kos{t}_{\pi}
    \end{tikzcd}
  \end{equation*}
\end{lemma}
\begin{proof}
  We write $X=(x,\gamma_x)$
  for an object in $\cat{Rep}_{\kos{S}}(\pi)$
  and $Y=(y,\gamma^{\pr}_y)$
  for an object in $\cat{Rep}_{\kos{S}}(\pi^{\pr})$.
  The right adjoint part
  $\kos{f}^!\colon \kosRep_{\kos{S}}(\pi)\to \kosRep_{\kos{S}}(\pi^{\pr})$
  is the strong monoidal functor
  that sends each $X$ to
  $(x,\gamma_x\circ (f\otimes I_x))\in \cat{Rep}_{\kos{S}}(\pi^{\pr})$
  and whose coherence isomorphisms are the identity morphisms.
  For each $Y$, we fix a coequalizer
  $(\pi\otimes y,\pc_{\pi}\otimes I_y)\twoheadrightarrow \kos{f}_!(Y)$
  in $\cat{Rep}_{\kos{S}}(\pi)$
  of the reflexive pair
  $I_{\pi}\otimes \gamma^{\pr}_y$,
  $(\pc_{\pi}\otimes I_y)\circ (I_{\pi}\otimes f\otimes I_y)
  \colon (\pi\otimes \pi^{\pr}\otimes y, \pc_{\pi}\otimes I_{\pi^{\pr}\otimes y})
  \to (\pi\otimes y,\pc_{\pi}\otimes I_y)$,
  with common section
  $I_{\pi}\otimes u_{\pi^{\pr}}\otimes I_y$.
  Then $\kos{f}_!$ is the left adjoint part of $\kos{f}=(\kos{f}^!,\kos{f}_!,\eta^f,\epsilon^f)$.
  The component $\epsilon^f_X\colon \kos{f}_!\kos{f}^!(X)\to X$
  of $\epsilon^f$ at $X$ is the unique morphism in $\cat{Rep}_{\kos{S}}(\pi)$
  whose composite with 
  $(\pi\otimes x,\pc_{\pi}\otimes I_x)\twoheadrightarrow \kos{f}_!\kos{f}^!(X)$
  is $\gamma_x$.
  The component $\eta^f_Y\colon Y\to \kos{f}^!\kos{f}_!(Y)$
  of $\eta^f$ at $Y$,
  seen as a morphism in $\CS$,
  is the composite of
  $u_{\pi}\otimes I_y\colon y\to \pi\otimes y$
  followed by $\pi\otimes y\twoheadrightarrow \kos{u}_{\pi\vphantom{!}}^!\kos{f}_!(Y)$.
  The invertible Galois transformation
  $\what{\kos{f}}\colon \kos{t}_{\pi\vphantom{!}}^{\vphantom{!}}\kos{f}\cong \kos{t}_{\pi^{\pr}\vphantom{!}}^{\vphantom{!}}$
  is the equality
  $\what{\kos{f}^!}\colon \kos{f}^!\kos{t}_{\pi\vphantom{!}}^!=\kos{t}_{\pi^{\pr}}^!$
  between the right adjoint parts,
  and the invertible Galois $\kos{S}$-transformation
  $\kos{f}\kos{u}_{\pi^{\pr}\vphantom{!}}^{\vphantom{!}}\cong \kos{u}_{\pi\vphantom{!}}^{\vphantom{!}}$
  is the equality
  $\kos{u}_{\pi^{\pr}}^!\kos{f}^!=\kos{u}_{\pi\vphantom{!}}^!$
  between the right adjoint parts.
  Hence $\kos{f}\kos{u}_{\pi^{\pr}\vphantom{!}}^{\vphantom{!}}$
  is a surjective fiber point of $\kos{t}_{\pi\vphantom{!}}^{\vphantom{!}}$,
  and $\kos{f}\colon \kos{t}_{\pi^{\pr}\vphantom{!}}^{\vphantom{!}}\mor\kos{t}_{\pi\vphantom{!}}^{\vphantom{!}}$
  is a fundamental $\kos{S}$-morphism.
\end{proof}

\begin{proposition}
  \label{prop-fundmoroverGalpk-inducedtor}
  Let $f\colon \pi^{\pr}\to\pi$
  be a morphism of Galois objects in $\kos{S}$
  and $p^{\pr}$ be a co-fpqc right $\pi^{\pr}$-torsor over $s$.
  We have a co-fpqc right $\pi$-torsor
  $p$ over $s$
  and an invertible Galois $\kos{S}$-transformation
  $\kos{f}\tw{p^{\pr}}{\omega}\cong \tw{p}{\omega}$,
  where
  $\kos{f}\colon \kos{t}_{\pi^{\pr}\vphantom{!}}^{\vphantom{!}}\mor\kos{t}_{\pi\vphantom{!}}^{\vphantom{!}}$
  is the fundamental $\kos{S}$-morphism
  associated to $f$
  in Lemma~\ref{lem-fundmoroverGalpk-fundmor},
  and 
  $\tw{p}{\omega}$, $\tw{p^{\pr}}{\omega}$
  are the twisted fiber points
  associated to $p$, $p^{\pr}$, respectively,
  in Proposition~\ref{prop-pFibGalTors-twistedfib}.
  \begin{equation*}
    \begin{tikzcd}[row sep=5pt,column sep=6pt]
      \text{ }
      &\id_{\kos{S}}
      \ar[ddr,filled,"{\tw{p}{\omega}}"]
      \ar[ddl,filled,"{\tw{p^{\pr}}{\omega}}"']
      &\text{ }
      \\%1
      \text{ }
      \ar[rr,2cell,shift right=0.7ex,pos=0.5]
      &\text{ }
      &\text{ }
      \\%2
      \kos{t}_{\pi^{\pr}}
      \ar[rr,filled,"{\kos{f}}"]
      &\text{ }
      &\kos{t}_{\pi}
    \end{tikzcd}
  \end{equation*}
\end{proposition}

We refer to $p$ as the \emph{extension of structure group}
of $p^{\pr}$ along $f$.

\begin{proof}[Proof of Proposition~\ref{prop-fundmoroverGalpk-inducedtor}]
  Let 
  $\lambda^{\pr}_{p^{\pr}}
  \colon p^{\pr}\otimes \pi^{\pr}
  \to p^{\pr}$
  be the right $\pi^{\pr}$-action
  on $p^{\pr}$.
  We write
  $\pc^f
  \coloneqq \pc_{\pi}
  \circ (f\otimes I_{\pi})
  \colon \pi^{\pr}\otimes \pi
  \to \pi$
  for the left $\pi^{\pr}$-action on
  $\pi$ induced by $f$.
  We fix a coequalizer
  $\mathrm{cq}
  \colon p^{\pr}\otimes \pi
  \twoheadrightarrow p$
  in $\cat{Ens}(\kos{S})$ of the reflexive pair
  $\lambda^{\pr}_{p^{\pr}}\otimes I_{\pi}$,
  $I_{p^{\pr}}\otimes \pc^f
  \colon p^{\pr}\otimes \pi^{\pr}\otimes \pi
  \to p^{\pr}\otimes \pi$,
  with common section
  $I_{p^{\pr}}\otimes u_{\pi^{\pr}}\otimes I_{\pi}$.
  Put $\tau^f
  \coloneqq (I_{p^{\pr}}\otimes \mathrm{cq})
  \circ (\cp_{p^{\pr}}\otimes I_{\pi})
  \colon p^{\pr}\otimes \pi
  \to p^{\pr}\otimes p$.
  In the diagram
  \begin{equation*}
    \vcenter{\hbox{
      \xymatrix@R=23pt@C=60pt{
        p^{\pr}\otimes \pi^{\pr}\otimes \pi^{\pr}\otimes \pi
        \ar@<0.5ex>[r]^-{I_{p^{\pr}\otimes \pi^{\pr}}\otimes \pc^f}
        \ar@<-0.5ex>[r]_-{I_{p^{\pr}}\otimes \pc_{\pi^{\pr}}\otimes I_{\pi}}
        \ar[d]_-{\tau^{\pr}_{p^{\pr}}\otimes I_{\pi^{\pr}\otimes \pi}}^-{\cong}
        &p^{\pr}\otimes \pi^{\pr}\otimes \pi
        \ar[r]^-{I_{p^{\pr}}\otimes \pc^f}
        \ar[d]^-{\tau^{\pr}_{p^{\pr}}\otimes I_{\pi}}_-{\cong}
        &p^{\pr}\otimes \pi
        \ar@{.>}[d]^-{\exists!\,\tau^f}_-{\cong}
        \\%1
        p^{\pr}\otimes p^{\pr}\otimes \pi^{\pr}\otimes \pi
        \ar@<0.5ex>[r]^-{I_{p^{\pr}\otimes p^{\pr}}\otimes \pc^f}
        \ar@<-0.5ex>[r]_-{I_{p^{\pr}}\otimes \lambda^{\pr}_{p^{\pr}}\otimes I_{\pi}}
        &p^{\pr}\otimes p^{\pr}\otimes \pi
        \ar@{->>}[r]^-{I_{p^{\pr}}\otimes \mathrm{cq}}
        &p^{\pr}\otimes p
      }
    }}
  \end{equation*}
  we see that
  (i) the bottom horizontal morphisms form a coequalizer diagram,
  since $p^{\pr}\otimes \slot$ preserves reflexive coequalizers;
  (ii) the top horizontal morphisms
  form a split coequalizer diagram with respect to the splitting
  $I_{p^{\pr}}\otimes u_{\pi^{\pr}}\otimes I_{\pi}$,
  $I_{p^{\pr}}\otimes u_{\pi^{\pr}}\otimes I_{\pi^{\pr}\otimes \pi}$;
  (iii) the middle and left vertical isomorphisms
  commute with the horizontal parallel pairs;
  and
  (iv) the unique isomorphism between the two coequalizers
  that makes the right square commute is given by
  $(I_{p^{\pr}}\otimes \mathrm{cq})
  \circ (\tau^{\pr}_{p^{\pr}}\otimes I_{\pi})
  \circ (I_{p^{\pr}}\otimes u_{\pi^{\pr}}\otimes I_{\pi})
  =\tau^f$.
  Since $\tau^f$ is an isomorphism in $\cat{Ens}(\kos{S})$,
  we deduce that $p$ is a co-fpqc covering of $s$.

  Let $\lambda_p\colon p\otimes \pi\to p$
  be the unique morphism in $\cat{Ens}(\kos{S})$
  that satisfies
  $\lambda_p
  \circ (\mathrm{cq}\otimes I_{\pi})
  =\mathrm{cq}
  \circ (I_{p^{\pr}}\otimes \pc_{\pi})
  \colon p^{\pr}\otimes \pi\otimes \pi
  \to p$;
  one can check that $\lambda_p$
  exists and is a right $\pi$-action
  on $p$.
  Using
  $\cp_p
  \circ \mathrm{cq}
  =(\mathrm{cq}\otimes\mathrm{cq})
  \circ \cp_{p^{\pr}\otimes \pi}
  \colon p^{\pr}\otimes \pi
  \to p\otimes p$,
  we deduce that
  $(I_{p^{\pr}}\otimes \tau_p)
  \circ (\tau^f\otimes I_{\pi})
  =(I_{p^{\pr}}\otimes \sigma_{p,p})
  \circ (\tau^f\otimes I_p)
  \circ (I_{p^{\pr}}\otimes \sigma_{p,\pi})
  \circ (\tau^f\otimes I_{\pi})
  \circ (I_{p^{\pr}}\otimes \tau_{\pi})
  \colon p^{\pr}\otimes \pi\otimes \pi
  \cong p^{\pr}\otimes p\otimes p$.
  Hence $I_{p^{\pr}}\otimes \tau_p$ is an isomorphism,
  and since $p^{\pr}\otimes \slot$ is conservative,
  $\tau_p$ is an isomorphism.
  Thus $p$ equipped with $\lambda_p$
  is a co-fpqc right $\pi$-torsor over $s$.

  It remains to construct the invertible Galois $\kos{S}$-transformation
  $\kos{f}\tw{p^{\pr}}{\omega}\cong \tw{p}{\omega}$.
  Put
  $g\coloneqq \mathrm{cq}
  \circ (I_{p^{\pr}}\otimes u_{\pi})
  \colon p^{\pr}\to p$.
  For each $X=(x,\gamma_x)\in\cat{Rep}_{\kos{S}}(\pi)$,
  let $\vartheta_X
  \colon \tw{p^{\pr}}{\omega}^!\kos{f}^!(X)
  \to \tw{p}{\omega}^!(X)$
  be the unique morphism in $\CS$ that satisfies 
  $\vartheta_X
  \circ \mathrm{cq}^{p^{\pr}}_{\kos{f}^!(X)}
  =\mathrm{cq}^p_X
  \circ (g\otimes I_x)
  \colon p^{\pr}\otimes x
  \to \tw{p}{\omega}^!(X)$.
  One can check that $\vartheta_X$ exists
  and that 
  $\kos{g}^!(\xi^p_X)
  \colon \kos{e}_{p^{\pr}}^!\kos{u}_{\pi\vphantom{!}}^!(X)
  \cong \kos{e}_{p^{\pr}}^!\tw{p}{\omega}^!(X)$
  is equal to
  $\kos{e}_{p^{\pr}}^!(\vartheta_X)
  \circ \xi^{p^{\pr}}_{\kos{f}^!(X)}
  \colon \kos{e}_{p^{\pr}}^!\kos{u}_{\pi^{\pr}}^!\kos{f}^!(X)
  \to \kos{e}_{p^{\pr}}^!\tw{p}{\omega}^!(X)$
  under the equality
  $\kos{e}_{p^{\pr}}^!\kos{u}_{\pi\vphantom{!}}^!
  =\kos{e}_{p^{\pr}}^!\kos{u}_{\pi^{\pr}}^!\kos{f}^!$.
  Thus $\kos{e}_{p^{\pr}}^!(\vartheta_X)$ is an isomorphism in $\CS_{p^{\pr}\vphantom{!}}^{\vphantom{!}}$,
  and since $\kos{e}_{p^{\pr}}^!$ is conservative,
  $\vartheta_X$ is an isomorphism in $\CS$.
  One can also check that the isomorphisms
  $\vartheta_X$ constitute
  an invertible Galois $\kos{S}$-transformation
  $\vartheta\colon \kos{f}\tw{p^{\pr}}{\omega}\cong \tw{p}{\omega}$.
\end{proof}

\begin{proposition}
  \label{prop-fundmoroverGalpk-realization}
  Let $\kos{T}\xmor{\kos{t}}\kos{S}$,
  $\kos{T}^{\pr}\xmor{\kos{t}^{\pr}}\kos{S}$
  be Galois $\kos{S}$-categories
  and $\kos{f}\colon \kos{t}^{\pr}\mor\kos{t}$
  be a Galois $\kos{S}$-morphism.
  The following are equivalent:
  \begin{enumerate}
    \item 
    $\kos{f}\colon \kos{t}^{\pr}\mor\kos{t}$
    is a fundamental $\kos{S}$-morphism.

    \item
    There exist a surjective fiber point
    $\omega$ of $\kos{t}$,
    a surjective fiber point $\omega^{\pr}$ of $\kos{t}^{\pr}$,
    and an invertible Galois $\kos{S}$-transformation
    $\alpha\colon \kos{f}\omega^{\pr}\cong\omega$.
  \end{enumerate}
  In this case, 
  let $\omega$, $\omega^{\pr}$, and $\alpha$
  be as in \textup{(2)}.
  Put $\pi\coloneqq \omega^!\omega_!(s)$
  and
  $\pi^{\pr}\coloneqq \omega^{\pr!}\omega^{\pr}_!(s)$,
  and recall the equivalences
  $\breve{\omega}\colon \kos{t}_{\pi\vphantom{!}}^{\vphantom{!}}\simeq\kos{t}$,
  $\breve{\omega}^{\pr}\colon \kos{t}_{\pi^{\pr}\vphantom{!}}^{\vphantom{!}}\simeq\kos{t}^{\pr}$
  in Theorem~\ref{thm-Galcatoverpk-mainthm}.
  We have a morphism
  $\breve{f}\colon \pi^{\pr}\to\pi$
  of Galois objects in $\kos{S}$
  and an invertible Galois $\kos{S}$-transformation
  $\breve{\alpha}\colon \kos{f}\breve{\omega}^{\pr}\cong\breve{\omega}\breve{\kos{f}}$
  that satisfies the following equality of pasting diagrams,
  where we denote by
  $\breve{\kos{f}}\colon \kos{t}_{\pi^{\pr}\vphantom{!}}^{\vphantom{!}}\mor\kos{t}_{\pi\vphantom{!}}^{\vphantom{!}}$
  the fundamental $\kos{S}$-morphism
  associated to $\breve{f}$
  in Lemma~\ref{lem-fundmoroverGalpk-fundmor}.
  \begin{equation*}
    \begin{tikzcd}[row sep=5pt,column sep=10pt]
      \text{ }
      &\text{ }
      &\id_{\kos{S}}
      \ar[ddll,bend right=1.8pc,filled,"{\kos{u}_{\pi^{\pr}}}"']
      \ar[dd,filled,"{\kos{u}_{\pi}}"]
      \ar[ddddrr,bend left=2.5pc,filled,"{\omega}"]
      &\text{ }
      &\text{ }
      &\text{ }
      &\text{ }
      &\text{ }
      &\id_{\kos{S}}
      \ar[ddll,bend right=1.8pc,filled,"{\kos{u}_{\pi^{\pr}}}"']
      \ar[dddd,filled,"{\omega^{\pr}}" description]
      \ar[ddddrr,bend left=2.5pc,filled,"{\omega}"]
      &\text{ }
      &\text{ }
      \\%1
      \text{ }
      \ar[rr,2cell,shift right=0.5ex,pos=0.75]
      &\text{ }
      &\text{ }
      &\text{ }
      &\text{ }
      &\text{ }
      &\text{ }
      &\text{ }
      &\text{ }
      &\text{ }
      &\text{ }
      \\%2
      \kos{t}_{\pi^{\pr}}
      \ar[rr,filled,"\breve{\kos{f}}"]
      \ar[ddrr,filled,"\breve{\omega}^{\pr}"',"\simeq"]
      &\text{ }
      &\kos{t}_{\pi}
      \ar[ddrr,filled,"\breve{\omega}","\simeq"']
      \ar[rr,2cell,shift left=0.5ex,pos=0.4]
      &\text{ }
      &\text{ }
      &=
      &\kos{t}_{\pi^{\pr}}
      \ar[ddrr,filled,"\breve{\omega}^{\pr}"',"\simeq"]
      \ar[rr,2cell,shift left=1ex,pos=0.6]
      &\text{ }
      &\text{ }
      \ar[rr,2cell,shift right=1.5ex,pos=0.45,"{\!\alpha}"]
      &\text{ }
      &\text{ }
      \\%3
      \text{ }
      &\text{ }
      \ar[rr,2cell,shift right=1ex,"{\!\breve{\alpha}}"]
      &\text{ }
      &\text{ }
      &\text{ }
      &\text{ }
      &\text{ }
      &\text{ }
      &\text{ }
      &\text{ }
      &\text{ }
      \\%4
      \text{ }
      &\text{ }
      &\kos{t}^{\pr}
      \ar[rr,filled,"\kos{f}"']
      &\text{ }
      &\kos{t}
      &\text{ }
      &\text{ }
      &\text{ }
      &\kos{t}^{\pr}
      \ar[rr,filled,"\kos{f}"']
      &\text{ }
      &\kos{t}
    \end{tikzcd}
  \end{equation*}
\end{proposition}
\begin{proof}
  It is immediate that (1) and (2) are equivalent.
  Let $\omega$, $\omega^{\pr}$, and $\alpha$ be as in (2).
  We write $\alpha_!\colon \omega^{\pr}_!\Rightarrow\kos{f}^!\omega_!$
  for the mate of $\alpha^{-1}\colon \omega^!\cong \omega^{\pr!}\kos{f}^!$
  with respect to the adjunctions
  $\omega$ and $\omega^{\pr}$;
  it is the unique colax $\kos{S}$-tensor natural transformation
  that satisfies the following relations.
  \begin{equation}\label{eq-fundmoroverGalpk-matealpha}
    \vcenter{\hbox{
      \xymatrix@R=15pt@C=30pt{
        \id_{\kos{S}}
        \ar@2{->}[r]^-{\eta}
        \ar@2{->}[d]_-{\eta^{\pr}}
        &\omega^!\omega_!
        \ar@2{->}[d]^-{\alpha^{-1}\omega_!}_-{\cong}
        &\text{ }
        &\omega^{\pr}_!\omega^!
        \ar@2{->}[r]^-{\alpha_!\omega^!}
        \ar@2{->}[d]_-{\omega^{\pr}_!\alpha^{-1}}^-{\cong}
        &\kos{f}^!\omega_!\omega^!
        \ar@2{->}[d]^-{\kos{f}^!\epsilon}
        \\%1
        \omega^{\pr!}\omega^{\pr}_!
        \ar@2{->}[r]^-{\omega^{\pr!}\alpha_!}
        &\omega^{\pr!}\kos{f}^!\omega_!
        &\text{ }
        &\omega^{\pr}_!\omega^{\pr!}\kos{f}^!
        \ar@2{->}[r]^-{\epsilon^{\pr}\kos{f}^!}
        &\kos{f}^!
      }
    }}
  \end{equation}
  Put $\alpha^{\sharp}
  \coloneqq (\alpha\omega_!)
  \circ (\omega^{\pr!}\alpha_!)
  \colon \omega^{\pr!}\omega^{\pr}_!\Rightarrow\omega^!\omega_!$
  and define 
  $\breve{f}\colon \pi^{\pr}\to \pi$ as 
  $\alpha^{\sharp}_{s}\colon \omega^{\pr!}\omega^{\pr}_!(s)\to \omega^!\omega_!(s)$.
  By substituting $\vartheta=\alpha^{\sharp}$
  in \eqref{eq-ActionbaseGal-colaxStnathatar},
  we obtain the following commutative diagram for every $z\in \CS$.
  \begin{equation}
    \label{eq-fundmoroverGalpk-alphasharpequiv}
    \vcenter{\hbox{
      \xymatrix@R=15pt@C=40pt{
        \omega^{\pr!}\omega^{\pr}_!(z)
        \ar[r]^-{(\hatar{\omega^{\pr!}\omega^{\pr}_!})_z}_-{\cong}
        \ar[d]_-{\alpha^{\sharp}_z}
        &\omega^{\pr!}\omega^{\pr}_!(s)\otimes z
        \ar[d]^-{\breve{f}\otimes I_z}
        \\%1
        \omega^!\omega_!(z)
        \ar[r]^-{(\hatar{\omega^!\omega_!})_z}_-{\cong}
        &\omega^!\omega_!(s)\otimes z
      }
    }}
  \end{equation}
  One can check that
  $\alpha^{\sharp}\colon \omega^{\pr!}\omega^{\pr}_!\Rightarrow\omega^!\omega_!$
  is a morphism of monads,
  using \eqref{eq-fundmoroverGalpk-matealpha}.
  From this,
  along with \eqref{eq-fundmoroverGalpk-alphasharpequiv} for $z=\pi$,
  we deduce that $\breve{f}\colon \pi^{\pr}\to\pi$
  is a morphism of Galois objects in $\kos{S}$.
  It remains to show the existence of $\breve{\alpha}$.
  It suffices to show that for each $X\in \CT$,
  the isomorphism
  $\alpha_X
  \colon \omega^{\pr!}\kos{f}^!(X)
  \cong \omega^!(X)$
  in $\CS$ lifts to an isomorphism
  $\breve{\alpha}_X
  \colon \breve{\omega}^{\pr!}\kos{f}^!(X)
  \cong \breve{\kos{f}}^!\breve{\omega}^!(X)$
  in $\cat{Rep}_{\kos{S}}(\pi^{\pr})$.
  One can check this using
  \eqref{eq-fundmoroverGalpk-alphasharpequiv}
  for $z=\omega^!(X)$ and
  \eqref{eq-fundmoroverGalpk-matealpha}.
\end{proof}

\begin{lemma}
  \label{lem-fundmoroverGalpk-sendingptpt}
  Let $\kos{t}$, $\kos{t}^{\pr}$
  be Galois $\kos{S}$-categories
  and $\kos{f}\colon \kos{t}^{\pr}\mor\kos{t}$
  be a fundamental $\kos{S}$-morphism.
  The composite of $\kos{f}$
  with any fiber point of $\kos{t}^{\pr}$
  is a fiber point of $\kos{t}$.
\end{lemma}
\begin{proof}
  Let $\omega$, $\omega^{\pr}$, $\alpha$,
  $\pi$, $\pi^{\pr}$,
  $\breve{\omega}$,
  $\breve{\omega}^{\pr}$,
  $\breve{f}$, $\breve{\kos{f}}$,
  and $\breve{\alpha}$
  be as in Proposition~\ref{prop-fundmoroverGalpk-realization},
  and $\omega^{\ppr}\colon \vphantom{\id}\smash[t]{\id}_{\kos{S}\vphantom{!}}^{\vphantom{!}}\mor\kos{t}^{\pr}$
  be a fiber point of $\kos{t}^{\pr}$.
  Recall from the proof of Theorem~\ref{thm-pFibGalTors-perfcorres} that
  $p^{\pr}\coloneqq \omega^{\ppr!}\omega^{\pr}_!(s)$
  is a co-fpqc right $\pi^{\pr}$-torsor over $s$,
  and that $\breve{\omega}^{\pr}\tw{p^{\pr}}{\omega}\cong \omega^{\ppr}$.
  Let $p$
  be the extension of structure group
  of $p^{\pr}$ along $\breve{f}$
  with 
  $\breve{\kos{f}}\tw{p^{\pr}}{\omega}\cong \tw{p}{\omega}$
  as in Proposition~\ref{prop-fundmoroverGalpk-inducedtor}.
  From these isomorphisms, along with 
  $\breve{\alpha}
  \colon \kos{f}\breve{\omega}^{\pr}
  \cong \breve{\omega}\breve{\kos{f}}$,
  we deduce that
  $\kos{f}\omega^{\ppr}
  \cong \kos{f}\breve{\omega}^{\pr}\tw{p^{\pr}}{\omega}
  \cong \breve{\omega}\breve{\kos{f}}\tw{p^{\pr}}{\omega}
  \cong \breve{\omega}\tw{p}{\omega}$.
  Thus $\kos{f}\omega^{\ppr}$ is a fiber point of $\kos{t}$.
\end{proof}

\begin{lemma}
  \label{lem-fundmoroverGalpk-closedundercomp}
  Fundamental $\kos{S}$-morphisms
  between Galois $\kos{S}$-categories
  are closed under composition.
\end{lemma}
\begin{proof}
  Let
  $\kos{f}\colon \kos{t}^{\pr}\mor\kos{t}$,
  $\kos{f}^{\pr}\colon \kos{t}^{\ppr}\mor\kos{t}^{\pr}$
  be fundamental $\kos{S}$-morphisms
  between Galois $\kos{S}$-categories
  $\kos{t}$, $\kos{t}^{\pr}$, $\kos{t}^{\ppr}$.
  Choose a surjective fiber point $\omega^{\ppr}$ of $\kos{t}^{\ppr}$
  such that $\kos{f}^{\pr}\omega^{\ppr}$ is a surjective fiber point of $\kos{t}^{\pr}$.
  By Lemma~\ref{lem-fundmoroverGalpk-sendingptpt},
  the composite $\kos{f}\kos{f}^{\pr}\omega^{\ppr}$
  is a fiber point of $\kos{t}$,
  and it is surjective by Corollary~\ref{cor-Galcatoverpk-pFibsurjective}.
  Thus $\kos{f}\kos{f}^{\pr}\colon \kos{t}^{\ppr}\mor\kos{t}$
  is a fundamental $\kos{S}$-morphism.
\end{proof}

\begin{lemma}
  \label{lem-fundmoroverGalpk-HomisIsom}
  Let $\kos{t}$, $\kos{t}^{\pr}$
  be Galois $\kos{S}$-categories
  and
  $\kos{f}$, $\kos{g}\colon \kos{t}^{\pr}\mor\kos{t}$
  be fundamental $\kos{S}$-morphisms.
  Every Galois $\kos{S}$-transformation
  from $\kos{f}$ to $\kos{g}$ is invertible.
\end{lemma}
\begin{proof}
  Let $\vartheta\colon \kos{f}\Rightarrow\kos{g}$
  be a Galois $\kos{S}$-transformation
  and choose a surjective fiber point
  $\omega^{\pr}$ of $\kos{t}^{\pr}$.
  By Lemma~\ref{lem-fundmoroverGalpk-sendingptpt},
  both $\kos{f}\omega^{\pr}$ and $\kos{g}\omega^{\pr}$
  are fiber points of $\kos{t}$.
  Thus 
  $\vartheta\omega^{\pr}\colon \kos{f}\omega^{\pr}\Rightarrow \kos{g}\omega^{\pr}$
  is invertible by Corollary~\ref{cor-Galcatoverpk-pFibgroupoid},
  and since $\omega^{\pr}$ is surjective, $\vartheta$ is invertible.
\end{proof}

Note that $\vphantom{\id}\smash[t]{\id}_{\kos{S}\vphantom{!}}^{\vphantom{!}}$ is a Galois $\kos{S}$-category,
as the identity Galois $\kos{S}$-morphism of $\vphantom{\id}\smash[t]{\id}_{\kos{S}\vphantom{!}}^{\vphantom{!}}$
is a surjective fiber point of $\vphantom{\id}\smash[t]{\id}_{\kos{S}\vphantom{!}}^{\vphantom{!}}$.
Conversely, every fiber point $\omega$ of $\vphantom{\id}\smash[t]{\id}_{\kos{S}\vphantom{!}}^{\vphantom{!}}$
is isomorphic to the identity Galois $\kos{S}$-morphism
of $\vphantom{\id}\smash[t]{\id}_{\kos{S}\vphantom{!}}^{\vphantom{!}}$;
we have $\what{\omega}\colon \omega\cong \vphantom{\id}\smash[t]{\id}_{\kos{S}\vphantom{!}}^{\vphantom{!}}$.

\begin{lemma}
  \label{lem-fundmoroverGalpk-pointismor}
  Let $\kos{t}$ be a Galois $\kos{S}$-category.
  Fundamental $\kos{S}$-morphisms
  from $\vphantom{\id}\smash[t]{\id}_{\kos{S}\vphantom{!}}^{\vphantom{!}}$ to $\kos{t}$
  are precisely the surjective fiber points of $\kos{t}$.
\end{lemma}
\begin{proof}
  Let $\omega\colon \vphantom{\id}\smash[t]{\id}_{\kos{S}\vphantom{!}}^{\vphantom{!}}\mor\kos{t}$ be a Galois $\kos{S}$-morphism.
  It is immediate that if $\omega$ is a surjective fiber point of $\kos{t}$,
  then $\omega$ is a fundamental $\kos{S}$-morphism.
  Conversely, assume that $\omega$ is a fundamental $\kos{S}$-morphism.
  Then $\omega$ is isomorphic to a surjective fiber point of $\kos{t}$,
  and is thus a surjective fiber point of $\kos{t}$.
\end{proof}

\begin{definition}
  \label{def-fundmoroverGalpk-pGalCatS}
  We define the $(2,1)$-category
  $\mathsf{pGalCat}(\kos{S})$
  as follows.
  \begin{itemize}
    \item 
    A $0$-cell in $\mathsf{pGalCat}(\kos{S})$
    is a Galois $\kos{S}$-category.
    
    \item
    A $1$-cell in $\mathsf{pGalCat}(\kos{S})$
    is a fundamental $\kos{S}$-morphism.

    \item
    A $2$-cell in $\mathsf{pGalCat}(\kos{S})$
    is an invertible Galois $\kos{S}$-transformation.

    \item
    The rest of the $2$-categorical structure of
    $\mathsf{pGalCat}(\kos{S})$
    is given by that of $\mathtt{pGalKos}/\kos{S}$;
    the composition of $1$-cells is well defined by
    Lemma~\ref{lem-fundmoroverGalpk-closedundercomp}.
  \end{itemize}
\end{definition}

\section{A duality over a Galois prekosmos}
\label{sec-dualityoverGalpk}

In this section, we establish a duality
between Galois objects in $\kos{S}$
and Galois $\kos{S}$-categories
with a choice of a surjective fiber point;
the latter are defined as follows.

\begin{definition}
  \label{def-dualityoverGalpk-ptGalcat}
  A \emph{pointed Galois $\kos{S}$-category},
  also called a \emph{pointed Galois category over $\kos{S}$},
  is a pair $(\kos{t},\omega)$ 
  consisting of a Galois $\kos{S}$-category $\kos{t}$
  and a surjective fiber point $\omega$ of $\kos{t}$.
\end{definition}

Recall from Section~\ref{sec-Galobjinpk}
the $(2,1)$-category $\mathsf{pGalObj}(\kos{S})$
of Galois objects in $\kos{S}$,
defined in \eqref{eq-Galobjinpk-pGalObjS21cat}.
The goal of this section is to
define the $(2,1)$-category
$\mathsf{pGalCat}_{\mathsf{pt}}(\kos{S})$
whose $0$-cells are pointed Galois $\kos{S}$-categories,
and prove the following theorem.

\begin{theorem}
  \label{thm-dualityoverGalpk-equiv21cat}
  We have mutually quasi-inverse equivalences
  \begin{equation*}
    \xymatrix@C=35pt{
      \mathsf{pGalObj}(\kos{S})
      \ar@<1ex>[r]^-{\mathsf{Rep}_{\kos{S}}}
      \ar@{}[r]|-{\simeq}
      &\mathsf{pGalCat}_{\mathsf{pt}}(\kos{S})
      \ar@<1ex>[l]^-{\mathsf{Aut}_{\kos{S}}}
    }
  \end{equation*}
  of $(2,1)$-categories.
\end{theorem}

\begin{definition}
  \label{def-dualityoverGalpk-pGalCatptS}
  We define the $(2,1)$-category
  $\mathsf{pGalCat}_{\mathsf{pt}}(\kos{S})$
  as follows.
  \begin{itemize}
    \item 
    A $0$-cell in $\mathsf{pGalCat}_{\mathsf{pt}}(\kos{S})$
    is a pointed Galois $\kos{S}$-category.
 
    \item
    A $1$-cell in $\mathsf{pGalCat}_{\mathsf{pt}}(\kos{S})$
    is a \emph{pointed fundamental $\kos{S}$-morphism}.
    Let $(\kos{t},\omega)$, $(\kos{t}^{\pr},\omega^{\pr})$
    be
    $0$-cells.
    A pointed fundamental $\kos{S}$-morphism
    $(\kos{f},\alpha)\colon (\kos{t}^{\pr},\omega^{\pr})\mor (\kos{t},\omega)$,
    denoted by an arrow with a filled tip,
    is a pair consisting of a Galois $\kos{S}$-morphism
    $\kos{f}\colon \kos{t}^{\pr}\mor\kos{t}$
    and an invertible Galois $\kos{S}$-transformation
    $\alpha\colon \kos{f}\omega^{\pr}\cong \omega$.
    Such $\kos{f}$ is necessarily a fundamental $\kos{S}$-morphism
    by Proposition~\ref{prop-fundmoroverGalpk-realization}.
    We illustrate
    $(\kos{f},\alpha)\colon (\kos{t}^{\pr},\omega^{\pr})\mor (\kos{t},\omega)$
    as follows.
    \begin{equation*}
      \begin{tikzcd}[row sep=5pt,column sep=7pt]
        \text{ }
        &\id_{\kos{S}}
        \ar[ddr,filled,"{\omega}"]
        \ar[ddl,filled,"{\omega^{\pr}}"']
        &\text{ }
        \\%1
        \text{ }
        \ar[rr,2cell,shift right=0.8ex,pos=0.5,"{\alpha}"]
        &\text{ }
        &\text{ }
        \\%2
        \kos{t}^{\pr}
        \ar[rr,filled,"{\kos{f}}"]
        &\text{ }
        &\kos{t}
      \end{tikzcd}
    \end{equation*}

    \item
    A $2$-cell 
    $\vartheta\colon (\kos{f},\alpha)\cong(\kos{g},\beta)$
    between $1$-cells
    $(\kos{f},\alpha)$, $(\kos{g},\beta)\colon (\kos{t}^{\pr},\omega^{\pr})\mor(\kos{t},\omega)$
    in $\mathsf{pGalCat}_{\mathsf{pt}}(\kos{S})$
    is an invertible Galois $\kos{S}$-transformation
    $\vartheta\colon \kos{f}\cong \kos{g}$
    between the underlying fundamental $\kos{S}$-morphisms.

    \item
    The composite of $1$-cells
    $(\kos{f},\alpha)\colon (\kos{t}^{\pr},\omega^{\pr})\mor(\kos{t},\omega)$,
    $(\kos{f}^{\pr},\alpha^{\pr})\colon (\kos{t}^{\ppr},\omega^{\ppr})\mor(\kos{t}^{\pr},\omega^{\pr})$
    is the pair $(\kos{f}\kos{f}^{\pr},\alpha\circ (\kos{f}\alpha^{\pr}))
    \colon (\kos{t}^{\ppr},\omega^{\ppr})\mor(\kos{t},\omega)$.
    The identity $1$-cell of a $0$-cell 
    $(\kos{t},\omega)$
    is the pair consisting of the identity Galois $\kos{S}$-morphism
    of $\kos{t}$
    and the identity Galois $\kos{S}$-transformation of
    $\omega\colon \vphantom{\id}\smash[t]{\id}_{\kos{S}\vphantom{!}}^{\vphantom{!}}\mor\kos{t}$.
    Horizontal and vertical compositions of $2$-cells 
    are given by those of invertible Galois $\kos{S}$-transformations
    between the underlying fundamental $\kos{S}$-morphisms.
    Identity $2$-cells are given by the identity Galois $\kos{S}$-transformations
    of the underlying fundamental $\kos{S}$-morphisms.
  \end{itemize}
\end{definition}

\begin{remark}
  \label{rem-dualityoverGalpk-notcoslice}
  The $0$-cells and the $1$-cells
  of $\mathsf{pGalCat}_{\mathsf{pt}}(\kos{S})$
  are those of the coslice $(2,1)$-category
  of $\mathsf{pGalCat}(\kos{S})$ under $\vphantom{\id}\smash[t]{\id}_{\kos{S}\vphantom{!}}^{\vphantom{!}}$,
  by Lemma~\ref{lem-fundmoroverGalpk-pointismor},
  whereas the $2$-cells of $\mathsf{pGalCat}_{\mathsf{pt}}(\kos{S})$
  are inherited from $\mathsf{pGalCat}(\kos{S})$
  without any compatibility conditions.
\end{remark}

\begin{proposition}
  \label{prop-dualityoverGalpk-RepSpseudofunct}
  We have a pseudofunctor
  \begin{equation*}
    \mathsf{Rep}_{\kos{S}}
    \colon \mathsf{pGalObj}(\kos{S})
    \to \mathsf{pGalCat}_{\mathsf{pt}}(\kos{S})
  \end{equation*}
  that sends each Galois object $\pi$ in $\kos{S}$
  to the pointed Galois $\kos{S}$-category
  $(\kos{t}_{\pi\vphantom{!}}^{\vphantom{!}},\kos{u}_{\pi\vphantom{!}}^{\vphantom{!}})$.
\end{proposition}
\begin{proof}
  We describe how $\mathsf{Rep}_{\kos{S}}$ sends $1$-cells and $2$-cells.
  \begin{itemize}
    \item
    A $1$-cell $f\colon \pi^{\pr}\to\pi$
    is sent to the pair consisting of
    the fundamental $\kos{S}$-morphism
    $\kos{f}\colon \kos{t}_{\pi^{\pr}\vphantom{!}}^{\vphantom{!}}\mor \kos{t}_{\pi\vphantom{!}}^{\vphantom{!}}$
    and the invertible Galois $\kos{S}$-transformation
    $\kos{f}\kos{u}_{\pi^{\pr}\vphantom{!}}^{\vphantom{!}}\cong \kos{u}_{\pi\vphantom{!}}^{\vphantom{!}}$
    in Lemma~\ref{lem-fundmoroverGalpk-fundmor}.

    \item
    Let $f$, $g\colon \pi^{\pr}\to\pi$ be $1$-cells 
    in $\mathsf{pGalObj}(\kos{S})$
    and write $\kos{f}$, $\kos{g}\colon \kos{t}_{\pi^{\pr}\vphantom{!}}^{\vphantom{!}}\mor\kos{t}_{\pi\vphantom{!}}^{\vphantom{!}}$
    for the associated fundamental $\kos{S}$-morphisms.
    A $2$-cell $\theta\colon f\cong g$
    in $\mathsf{pGalObj}(\kos{S})$
    is sent to the $2$-cell
    $\vartheta\colon \kos{f}\cong \kos{g}$
    in $\mathsf{pGalCat}_{\mathsf{pt}}(\kos{S})$,
    defined as follows.
    The component $\vartheta_X\colon \kos{f}^!(X)\cong \kos{g}^!(X)$
    of $\vartheta$ at $X=(x,\gamma_x)\in\cat{Rep}_{\kos{S}}(\pi)$
    is the isomorphism in $\cat{Rep}_{\kos{S}}(\pi^{\pr})$
    whose underlying isomorphism in $\CS$
    is $\gamma_x\circ (\theta\otimes I_x)\colon x\cong x$.
    It is straightforward to check that
    $\vartheta$ is indeed an invertible Galois $\kos{S}$-transformation
    from $\kos{f}$ to $\kos{g}$.
  \end{itemize}
  It is straightforward to check that 
  $\mathsf{Rep}_{\kos{S}}$ preserves vertical composition of $2$-cells and identity $2$-cells;
  thus $\mathsf{Rep}_{\kos{S}}$ defines a functor
  on each hom-groupoid.
  The constraints of the pseudofunctor $\mathsf{Rep}_{\kos{S}}$
  are defined as follows.
  \begin{itemize}
    \item
    Let $f$, $f^{\pr}$
    be a composable pair of $1$-cells in $\mathsf{pGalObj}(\kos{S})$,
    with composite $f^{\ppr}\coloneqq ff^{\pr}$.
    We write
    $\kos{f}$, $\kos{f}^{\pr}$, and $\kos{f}^{\ppr}$
    for the fundamental $\kos{S}$-morphisms
    associated to $f$, $f^{\pr}$, and $f^{\ppr}$, respectively.
    The component of the functorial constraint of $\mathsf{Rep}_{\kos{S}}$
    at $f$, $f^{\pr}$
    is the invertible Galois $\kos{S}$-transformation
    $\kos{f}\kos{f}^{\pr}\cong \kos{f}^{\ppr}$
    given by the equality
    $\kos{f}^{\pr!}\kos{f}^!=\kos{f}^{\ppr!}$
    between the right adjoint parts.

    \item
    Let $\pi$ be a $0$-cell in $\mathsf{pGalObj}(\kos{S})$.
    We write 
    $\kos{I}_{\pi\vphantom{!}}^{\vphantom{!}}\colon \kos{t}_{\pi\vphantom{!}}^{\vphantom{!}}\mor\kos{t}_{\pi\vphantom{!}}^{\vphantom{!}}$
    for the fundamental $\kos{S}$-morphism
    associated to the identity morphism of $\pi$
    in Lemma~\ref{lem-fundmoroverGalpk-fundmor}.
    The component of the unity constraint of $\mathsf{Rep}_{\kos{S}}$
    at $\pi$
    is the invertible Galois $\kos{S}$-transformation
    from the identity Galois $\kos{S}$-morphism of $\kos{t}_{\pi\vphantom{!}}^{\vphantom{!}}$
    to $\kos{I}_{\pi\vphantom{!}}^{\vphantom{!}}$,
    given by the equality between the right adjoint parts.
  \end{itemize}
  One can check that $\mathsf{Rep}_{\kos{S}}$ preserves
  horizontal composition of $2$-cells
  up to the functorial constraint,
  or equivalently,
  that the functorial constraint
  $\kos{f}\kos{f}^{\pr}\cong \kos{f}^{\ppr}$
  is natural in $f$, $f^{\pr}$.
  It is immediate that the functorial and unity constraints of $\mathsf{Rep}_{\kos{S}}$
  satisfy the associativity and left/right unity axioms,
  as the axioms are just equalities between the right adjoint parts.
  This completes the construction of
  $\mathsf{Rep}_{\kos{S}}
  \colon \mathsf{pGalObj}(\kos{S})\to \mathsf{pGalCat}_{\mathsf{pt}}(\kos{S})$.
\end{proof}

\begin{proposition}
  \label{prop-dualityoverGalpk-AutS2funct}
  We have a $2$-functor
  \begin{equation*}
    \mathsf{Aut}_{\kos{S}}
    \colon \mathsf{pGalCat}_{\mathsf{pt}}(\kos{S})\to \mathsf{pGalObj}(\kos{S})
  \end{equation*}
  that sends each pointed Galois $\kos{S}$-category
  $(\kos{t},\omega)$
  to the Galois object $\omega^!\omega_!(s)$ in $\kos{S}$
  that represents $\uAut_{\kos{S}}(\omega)$ in Theorem~\ref{thm-Galcatoverpk-mainthm}.
\end{proposition}
\begin{proof}
  We describe how $\mathsf{Aut}_{\kos{S}}$ sends $1$-cells and $2$-cells.
  \begin{itemize}
    \item
    A $1$-cell
    $(\kos{f},\alpha)\colon (\kos{t}^{\pr},\omega^{\pr})\mor(\kos{t},\omega)$
    in $\mathsf{pGalCat}_{\mathsf{pt}}(\kos{S})$ is sent to
    the $1$-cell
    $\breve{f}
    \colon \omega^{\pr!}\omega^{\pr}_!(s)\to \omega^!\omega_!(s)$
    in $\mathsf{pGalObj}(\kos{S})$
    described in the proof of Proposition~\ref{prop-fundmoroverGalpk-realization}.

    \item
    A $2$-cell $\vartheta\colon (\kos{f},\alpha)\cong (\kos{g},\beta)$
    between $1$-cells
    $(\kos{f},\alpha)$, $(\kos{g},\beta)
    \colon (\kos{t}^{\pr},\omega^{\pr})\mor (\kos{t},\omega)$
    in $\mathsf{pGalCat}_{\mathsf{pt}}(\kos{S})$
    is sent to the $2$-cell
    $\breve{\theta}\colon \breve{f}\cong \breve{g}$
    in $\mathsf{pGalObj}(\kos{S})$,
    defined as follows.
    Put $\pi\coloneqq \omega^!\omega_!(s)$
    and 
    $\vartheta^{\sharp}
    \coloneqq (\beta\omega_!)
    \circ (\omega^{\pr!}\vartheta\omega_!)
    \circ (\alpha^{-1}\omega_!)
    \circ \eta
    \colon \vphantom{\id}\smash[t]{\id}_{\kos{S}\vphantom{!}}^{\vphantom{!}}
    \Rightarrow \omega^!\omega_!$.
    We define
    $\breve{\theta}\colon s\to \pi$
    as $\vartheta^{\sharp}_{s}$.
    By \eqref{eq-ActionbaseGal-colaxStnathatar}
    for the colax $\kos{S}$-tensor natural transformation $\vartheta^{\sharp}$,
    we deduce that the following diagram commutes
    for every $z\in \CS$.
    \begin{equation}
      \label{eq-dualityoverGalpk-varthetasharpequiv}
      \vcenter{\hbox{
        \xymatrix@R=12pt@C=35pt{
          z
          \ar[d]_-{\vartheta^{\sharp}_z}
          \ar[r]^-{\imath_z}_-{\cong}
          &s\otimes z
          \ar[d]^-{\breve{\theta}\otimes I_z}
          \\%1
          \omega^!\omega_!(z)
          \ar[r]^-{(\hatar{\omega^!\omega_!})_z}_-{\cong}
          &\omega^!\omega_!(s)\otimes z
        }
      }}
    \end{equation}
    Using \eqref{eq-fundmoroverGalpk-matealpha},
    one can check that 
    $\vartheta^{\sharp}
    =(\beta\omega_!)
    \circ (\omega^{\pr!}\vartheta\omega_!)
    \circ (\omega^{\pr!}\alpha_!)
    \circ \eta^{\pr}
    \colon \vphantom{\id}\smash[t]{\id}_{\kos{S}\vphantom{!}}^{\vphantom{!}}
    \Rightarrow \omega^!\omega_!$
    and that
    $(\omega^!\epsilon\omega_!)
    \circ (\vartheta^{\sharp}\alpha^{\sharp})
    =(\beta\omega_!)
    \circ (\omega^{\pr!}\vartheta\omega_!)
    \circ (\omega^{\pr!}\alpha_!)
    =(\omega^!\epsilon\omega_!)
    \circ (\beta^{\sharp}\vartheta^{\sharp})
    \colon \omega^{\pr!}\omega^{\pr}_!
    \Rightarrow \omega^!\omega_!$.
    From these equalities,
    along with \eqref{eq-fundmoroverGalpk-alphasharpequiv} for $z=\pi$
    and \eqref{eq-dualityoverGalpk-varthetasharpequiv} for $z=\pi$,
    we deduce that
    $\breve{\theta}$
    is indeed a $2$-isomorphism from $\breve{f}$ to $\breve{g}$.
  \end{itemize}
  It is straightforward to check that
  $\mathsf{Aut}_{\kos{S}}$
  preserves identity $1$-cells and identity $2$-cells.
  \begin{itemize}
    \item 
    Let $\vartheta\colon (\kos{f},\alpha)\cong (\kos{g},\beta)$
    and
    $\tilde{\vartheta}\colon (\kos{g},\beta)\cong (\tilde{\kos{g}},\tilde{\beta})$
    be $2$-cells in $\mathsf{pGalCat}_{\mathsf{pt}}(\kos{S})$,
    with vertical composite
    $\bar{\vartheta}\coloneqq \tilde{\vartheta}\circ \vartheta
    \colon (\kos{f},\alpha)\cong (\tilde{\kos{g}},\tilde{\beta})$.
    Put $\pi\coloneqq \omega^!\omega_!(s)$.
    One can check that
    $(\omega^!\epsilon\omega_!)
    \circ (\tilde{\vartheta}^{\sharp}\vartheta^{\sharp})
    =\bar{\vartheta}^{\sharp}
    \colon \vphantom{\id}\smash[t]{\id}_{\kos{S}\vphantom{!}}^{\vphantom{!}}
    \Rightarrow \omega^!\omega_!$,
    using \eqref{eq-fundmoroverGalpk-matealpha}.
    From this,
    along with \eqref{eq-dualityoverGalpk-varthetasharpequiv}
    for $z=\pi$,
    we deduce that $\breve{\tilde{\theta}}\star\breve{\theta}=\breve{\bar{\theta}}$.
    This shows that $\mathsf{Aut}_{\kos{S}}$
    preserves vertical composition of $2$-cells.
    
    \item 
    Let $(\kos{f},\alpha)\colon (\kos{t}^{\pr},\omega^{\pr})\mor(\kos{t},\omega)$
    and
    $(\kos{f}^{\pr},\alpha^{\pr})\colon (\kos{t}^{\ppr},\omega^{\ppr})\mor(\kos{t}^{\pr},\omega^{\pr})$
    be $1$-cells in $\mathsf{pGalCat}_{\mathsf{pt}}(\kos{S})$,
    with composite denoted by
    $(\kos{f}^{\ppr},\alpha^{\ppr})\colon (\kos{t}^{\ppr},\omega^{\ppr})\mor(\kos{t},\omega)$,
    where $\kos{f}^{\ppr}\coloneqq \kos{f}\kos{f}^{\pr}$
    and $\alpha^{\ppr}\coloneqq \alpha\circ (\kos{f}\alpha^{\pr})
    \colon \kos{f}\kos{f}^{\pr}\omega^{\ppr}\cong \omega$.
    One can check that 
    $\alpha^{\ppr}_!
    =(\kos{f}^{\pr!}\alpha_!)
    \circ \alpha^{\pr}_!
    \colon \omega^{\ppr}_!
    \Rightarrow \kos{f}^{\pr!}\kos{f}^!\omega_!$
    and that
    $\alpha^{\ppr\sharp}
    =\alpha^{\sharp}\circ \alpha^{\pr\sharp}
    \colon \omega^{\ppr!}\omega^{\ppr}_!
    \Rightarrow \omega^!\omega_!$,
    using \eqref{eq-fundmoroverGalpk-matealpha}.
    Thus we have $\breve{f}^{\ppr}=\breve{f}\circ \breve{f}^{\pr}
    \colon \omega^{\ppr!}\omega^{\ppr}_!(s)
    \to \omega^!\omega_!(s)$.
    This shows that $\mathsf{Aut}_{\kos{S}}$
    preserves composition of $1$-cells.
  
    \item
    Let $(\kos{f},\alpha)$,
    $(\kos{g},\beta)\colon (\kos{t}^{\pr},\omega^{\pr})\mor (\kos{t},\omega)$,
    $(\kos{f}^{\pr},\alpha^{\pr})$,
    $(\kos{g}^{\pr},\beta^{\pr})\colon (\kos{t}^{\ppr},\omega^{\ppr})\mor (\kos{t}^{\pr},\omega^{\pr})$
    be $1$-cells and
    $\vartheta\colon (\kos{f},\alpha)\cong (\kos{g},\beta)$,
    $\vartheta^{\pr}\colon (\kos{f}^{\pr},\alpha^{\pr})\cong (\kos{g}^{\pr},\beta^{\pr})$
    be $2$-cells in $\mathsf{pGalCat}_{\mathsf{pt}}(\kos{S})$.
    Write
    $(\kos{f}^{\ppr},\alpha^{\ppr})$,
    $(\kos{g}^{\ppr},\beta^{\ppr})\colon (\kos{t}^{\ppr},\omega^{\ppr})\mor(\kos{t},\omega)$
    for the composite $1$-cells
    and
    $\vartheta^{\ppr}\colon (\kos{f}^{\ppr},\alpha^{\ppr})\cong(\kos{g}^{\ppr},\beta^{\ppr})$
    for the horizontal composite of $\vartheta$ and $\vartheta^{\pr}$.
    One can check that
    $\vartheta^{\ppr\sharp}
    =(\omega^!\epsilon\omega_!)
    \circ (\beta^{\sharp}\omega^!\omega_!)
    \circ (\vartheta^{\pr\sharp}\vartheta^{\sharp})
    \colon \vphantom{\id}\smash[t]{\id}_{\kos{S}\vphantom{!}}^{\vphantom{!}}
    \Rightarrow \omega^!\omega_!$,
    using \eqref{eq-fundmoroverGalpk-matealpha}.
    From this,
    along with \eqref{eq-dualityoverGalpk-varthetasharpequiv}
    for $z=\pi\coloneqq \omega^!\omega_!(s)$,
    we deduce that
    $\breve{\theta}^{\ppr}
    =(\breve{g}\circ\breve{\theta}^{\pr})\star \breve{\theta}$.
    This shows that $\mathsf{Aut}_{\kos{S}}$
    preserves horizontal composition of $2$-cells.
  \end{itemize}
  This completes the construction of
  $\mathsf{Aut}_{\kos{S}}
  \colon \mathsf{pGalCat}_{\mathsf{pt}}(\kos{S})\to \mathsf{pGalObj}(\kos{S})$.
\end{proof}

\begin{proof}[Proof of Theorem~\ref{thm-dualityoverGalpk-equiv21cat}]
  We show that the pseudofunctor
  $\mathsf{Rep}_{\kos{S}}$
  in Proposition~\ref{prop-dualityoverGalpk-RepSpseudofunct}
  and the $2$-functor
  $\mathsf{Aut}_{\kos{S}}$
  in Proposition~\ref{prop-dualityoverGalpk-AutS2funct}
  are mutually quasi-inverse equivalences of $(2,1)$-categories;
  that is, the composites
  $\mathsf{Aut}_{\kos{S}}\circ \mathsf{Rep}_{\kos{S}}$
  and $\mathsf{Rep}_{\kos{S}}\circ \mathsf{Aut}_{\kos{S}}$
  are equivalent to the identity $2$-functors
  of $\mathsf{pGalObj}(\kos{S})$
  and $\mathsf{pGalCat}_{\mathsf{pt}}(\kos{S})$,
  respectively.
  
  The composite
  $\mathsf{Aut}_{\kos{S}}\circ \mathsf{Rep}_{\kos{S}}\colon 
  \mathsf{pGalObj}(\kos{S})\to \mathsf{pGalObj}(\kos{S})$
  sends each $0$-cell $\pi$
  to the $0$-cell $\pi\otimes s$,
  and we have an isomorphism
  $\jmath_{\pi}^{-1}\colon \pi\otimes s\cong \pi$
  of Galois objects in $\kos{S}$;
  see Remark~\ref{rem-GalMainthm-reconstruction}.
  It is straightforward to check that
  $\mathsf{Aut}_{\kos{S}}\circ \mathsf{Rep}_{\kos{S}}$
  is a $2$-functor, and that the isomorphisms $\jmath_{\pi}^{-1}$
  constitute a $2$-natural isomorphism
  from $\mathsf{Aut}_{\kos{S}}\circ \mathsf{Rep}_{\kos{S}}$
  to the identity $2$-functor of $\mathsf{pGalObj}(\kos{S})$.

  We now construct an equivalence
  from the composite
  $\mathsf{Rep}_{\kos{S}}\circ \mathsf{Aut}_{\kos{S}}$
  to the identity $2$-functor of $\mathsf{pGalCat}_{\mathsf{pt}}(\kos{S})$
  as follows.
  Such an equivalence is an invertible strong transformation
  in the sense of
  \cite[Definition~4.2.1, Corollary~4.4.13, and Definition~5.1.18]{Johnson2021};
  it is an equivalence internal to the bicategory of
  pseudofunctors from $\mathsf{pGalCat}_{\mathsf{pt}}(\kos{S})$ to itself,
  strong transformations, and modifications.
  We take the component $1$-cell of this transformation
  at each $0$-cell $(\kos{t},\omega)$ to be the $1$-cell
  $(\kos{t}_{\omega^!\omega_!(s)},\kos{u}_{\omega^!\omega_!(s)})
  \simeq
  (\kos{t},\omega)$
  consisting of the equivalence
  $\breve{\omega}\colon \kos{t}_{\omega^!\omega_!(s)}\simeq\kos{t}$
  and 
  the invertible Galois $\kos{S}$-transformation
  $\breve{\omega}\kos{u}_{\omega^!\omega_!(s)}\cong \omega$
  in Theorem~\ref{thm-Galcatoverpk-mainthm}.
  We take its component $2$-cell
  at each $1$-cell 
  $(\kos{f},\alpha)\colon (\kos{t}^{\pr},\omega^{\pr})\mor (\kos{t},\omega)$
  to be the invertible Galois $\kos{S}$-transformation
  $\breve{\alpha}\colon \kos{f}\breve{\omega}^{\pr}\cong \breve{\omega}\breve{\kos{f}}$
  in Proposition~\ref{prop-fundmoroverGalpk-realization}.
  It is straightforward to check,
  using \cite[Proposition~6.2.16]{Johnson2021}, that
  these data indeed define an invertible strong transformation.
  This completes the proof of Theorem~\ref{thm-dualityoverGalpk-equiv21cat}.
\end{proof}

%% file: PKGalCat.tex
\chapter{Prekosmic Galois Categories}
\label{chap-PKGalCat}

In this chapter, we gather the dualities
established in Section~\ref{sec-dualityoverGalpk},
one for each Galois prekosmos $\kos{S}$,
into a single duality.
We begin with the following definitions;
recall Definitions~\ref{def-Galprekos-Galcatoverpk},
\ref{def-Galobjinpk-Galobj},
and~\ref{def-dualityoverGalpk-ptGalcat}.

\begin{definition}
  \label{def-PKGalCat-PKGalCat}
  A \emph{prekosmic Galois category}
  is a pair $(\kos{S},\kos{t})$
  consisting of a Galois prekosmos $\kos{S}$
  and a Galois $\kos{S}$-category
  $\kos{T}\xmor{\kos{t}}\kos{S}$.
\end{definition}

\begin{definition}
  \label{def-PKGalCat-PKGalobj}
  A \emph{prekosmic Galois object}
  is a pair $(\kos{S},\pi)$
  consisting of a Galois prekosmos $\kos{S}$
  and a Galois object $\pi$ in $\kos{S}$.
\end{definition}

\begin{definition}
  \label{def-PKGalCat-ptPKGalCat}
  A \emph{pointed prekosmic Galois category}
  is a tuple $(\kos{S},\kos{t},\omega)$
  consisting of a Galois prekosmos $\kos{S}$,
  a Galois $\kos{S}$-category
  $\kos{T}\xmor{\kos{t}}\kos{S}$,
  and a surjective fiber point $\omega\colon \vphantom{\id}\smash[t]{\id}_{\kos{S}\vphantom{!}}^{\vphantom{!}}\mor\kos{t}$ of $\kos{t}$.
\end{definition}

The change of base along an \'{e}tale Galois morphism
$\kos{S}^{\pr}\mor\kos{S}$
between Galois prekosmoi
sends each Galois object in $\kos{S}$
to one in $\kos{S}^{\pr}$.
We write $\mathsf{pGalKos}_{\etsf}$
for the $(2,1)$-category 
of Galois prekosmoi,
\'{e}tale Galois morphisms,
and invertible Galois transformations.
Prekosmic Galois objects form a $(2,1)$-category
$\mathsf{pGalObj}$ over $\mathsf{pGalKos}_{\etsf}$,
and the strict fiber over each Galois prekosmos $\kos{S}$
is $2$-isomorphic to $\mathsf{pGalObj}(\kos{S})$
in \eqref{eq-Galobjinpk-pGalObjS21cat}.
Pointed prekosmic Galois categories
likewise form a $(2,1)$-category
$\mathsf{pGalCat}_{\mathsf{pt}}$
over $\mathsf{pGalKos}_{\etsf}$,
and the strict fiber over each Galois prekosmos $\kos{S}$
is $2$-isomorphic to $\mathsf{pGalCat}_{\mathsf{pt}}(\kos{S})$
in Definition~\ref{def-dualityoverGalpk-pGalCatptS}.
The following theorem gathers,
for each Galois prekosmos $\kos{S}$,
the mutually quasi-inverse equivalences
$\mathsf{Rep}_{\kos{S}}$ and $\mathsf{Aut}_{\kos{S}}$
in Theorem~\ref{thm-dualityoverGalpk-equiv21cat}.

\begin{theorem}
  \label{thm-PKGalCat-equiv21catoveret}
  We have mutually quasi-inverse equivalences
  \begin{equation*}
    \vcenter{\hbox{
      \xymatrix@R=25pt@C=0pt{
        \mathsf{pGalObj}
        \ar@<1ex>[rr]^-{\mathsf{Rep}}
        \ar[dr]_-{\mathsf{p}_{\mathsf{obj}}}
        \ar@{}[rr]|-{\simeq}
        &\text{ }
        &\mathsf{pGalCat}_{\mathsf{pt}}
        \ar@<1ex>[ll]^-{\mathsf{Aut}}
        \ar[dl]^-{\mathsf{p}_{\mathsf{ptcat}}}
        \\%1
        \text{ }
        &\mathsf{pGalKos}_{\etsf}
        &\text{ }
      }
    }}
  \end{equation*}
  of $(2,1)$-categories over $\mathsf{pGalKos}_{\etsf}$.
\end{theorem}

\begin{remark}
  \label{rem-PKGalCat-cleavagetransport}
  The $2$-functor $\mathsf{p}_{\mathsf{obj}}$
  is a split $2$-fibration,
  with the canonical cleavage described in Section~\ref{sec-PKGalobj}.
  The equivalences in Theorem~\ref{thm-PKGalCat-equiv21catoveret}
  transport this cleavage to $\mathsf{p}_{\mathsf{ptcat}}$,
  making it a cloven fibration of bicategories.
  We refer to \cite{Buckley2014} for the terminology.
\end{remark}

The following lemma ensures that
the $(2,1)$-category $\mathsf{pGalKos}_{\etsf}$ is well defined.

\begin{lemma}
  \label{lem-PKGalCat-etalecomposite}
  \'{E}tale Galois morphisms 
  between Galois prekosmoi are closed under composition.
\end{lemma}
\begin{proof}
  Let $\phi\colon \kos{S}^{\pr}\mor\kos{S}$,
  $\phi^{\pr}\colon \kos{S}^{\ppr}\mor\kos{S}^{\pr}$
  be \'{e}tale Galois morphisms between Galois prekosmoi,
  and $\varphi$, $\varphi^{\pr}$, and $\varphi^{\ppr}$
  be the natural transformations in \eqref{eq-Galprekos-projformula}
  associated to $\phi$, $\phi^{\pr}$, and $\phi^{\ppr}\coloneqq\phi\phi^{\pr}$, respectively.
  One can check that
  $\varphi^{\ppr}_{x^{\ppr},z}=
  \varphi_{\phi^{\pr}_!(x^{\ppr}),z}
  \circ \phi_!(\varphi^{\pr}_{x^{\ppr},\phi^!(z)})$
  for every $z\in \CS$, $x^{\ppr}\in\CS^{\ppr}$,
  and since $\varphi$ and $\varphi^{\pr}$ are invertible,
  we deduce that $\varphi^{\ppr}$ is invertible.
  The left adjoint part $\phi^{\ppr}_!=\phi_!\phi^{\pr}_!$
  is conservative, since $\phi_!$ and $\phi^{\pr}_!$ are.
  Hence $\phi^{\ppr}$ is \'{e}tale.
\end{proof}

\section{Prekosmic Galois objects}
\label{sec-PKGalobj}

Recall from Definition~\ref{def-PKGalCat-PKGalobj}
the definition of a prekosmic Galois object.
Let $(\kos{S},\pi)$,
$(\kos{S}^{\pr},\pi^{\pr})$
be prekosmic Galois objects.
A \emph{morphism}
$(\phi,f)\colon (\kos{S}^{\pr},\pi^{\pr})\to (\kos{S},\pi)$
is a pair consisting of an \'{e}tale Galois morphism
$\phi\colon \kos{S}^{\pr}\mor\kos{S}$
and a morphism
$f\colon \pi^{\pr}\to \phi^!(\pi)$
of Galois objects in $\kos{S}^{\pr}$;
see Remark~\ref{rem-PKGalobj-etalebasechange}.
Let 
$(\psi,g)\colon (\kos{S}^{\pr},\pi^{\pr})\to (\kos{S},\pi)$
be a morphism.
A \emph{$2$-isomorphism}
$(\chi,\theta)\colon (\phi,f)\cong (\psi,g)$
is a pair consisting of 
an invertible Galois transformation $\chi\colon \phi\cong \psi$
and a $2$-isomorphism
$\theta\colon f\cong \chi_{\pi}^{-1}g$
between morphisms
$f$, $\chi_{\pi}^{-1}g\colon \pi^{\pr}\to \phi^!(\pi)$.
Prekosmic Galois objects,
morphisms between them,
and $2$-isomorphisms between these morphisms
form a $(2,1)$-category,
which we denote by
\begin{equation}
  \label{eq-PKGalobj-21cat}
  \mathsf{pGalObj}.
\end{equation}
We describe the rest of the structure of $\mathsf{pGalObj}$ as follows.

Let 
$(\tilde{\psi},\tilde{g})\colon (\kos{S}^{\pr},\pi^{\pr})\to (\kos{S},\pi)$
be a morphism
and $(\tilde{\chi},\tilde{\theta})\colon (\psi,g)\cong (\tilde{\psi},\tilde{g})$
be a $2$-isomorphism.
The vertical composite of
$(\chi,\theta)$ and $(\tilde{\chi},\tilde{\theta})$
is the pair $(\tilde{\chi}\circ \chi,(\chi_{\pi}^{-1}\tilde{\theta})\star\theta)
\colon (\phi,f)\cong (\tilde{\psi},\tilde{g})$
consisting of
$\tilde{\chi}\circ\chi\colon \phi\cong \tilde{\psi}$
and
$(\chi_{\pi}^{-1}\tilde{\theta})\star\theta
\colon f
\cong \chi_{\pi}^{-1}\tilde{\chi}_{\pi}^{-1}\tilde{g}$.
The identity $2$-isomorphism of $(\phi,f)$
is the pair consisting of the identity Galois transformation of $\phi$
and the identity $2$-isomorphism of $f$.
Each $2$-isomorphism $(\chi,\theta)\colon (\phi,f)\cong (\psi,g)$
is indeed invertible with respect to the vertical composition;
its inverse is the pair
$(\chi^{-1},\chi_{\pi}\theta^{-1})\colon (\psi,g)\cong (\phi,f)$
consisting of
$\chi^{-1}\colon \psi\cong \phi$
and
$\chi_{\pi}\theta^{-1}\colon g\cong \chi_{\pi}f$,
where $\theta^{-1}\colon \chi_{\pi}^{-1}g\cong f$
is the inverse of $\theta$ with respect to $\star$.

Let $(\phi^{\pr},f^{\pr})$, $(\psi^{\pr},g^{\pr})
\colon (\kos{S}^{\ppr},\pi^{\ppr})\to (\kos{S}^{\pr},\pi^{\pr})$
be morphisms.
The composite of $(\phi,f)$ and $(\phi^{\pr},f^{\pr})$ is the pair
$(\phi\phi^{\pr},\phi^{\pr!}(f)f^{\pr})$
consisting of
$\phi\phi^{\pr}\colon \kos{S}^{\ppr}\mor\kos{S}$
and 
$\phi^{\pr!}(f)f^{\pr}\colon \pi^{\ppr}\to 
\phi^{\pr!}\phi^!(\pi)$.
The identity morphism of $(\kos{S},\pi)$
is the pair consisting of 
the identity Galois morphism $\vphantom{\id}\smash[t]{\id}_{\kos{S}\vphantom{!}}^{\vphantom{!}}$ of $\kos{S}$
and the identity morphism of $\pi$.
Let 
$(\chi^{\pr},\theta^{\pr})\colon (\phi^{\pr},f^{\pr})\cong (\psi^{\pr},g^{\pr})$
be a $2$-isomorphism.
The horizontal composite of 
$(\chi,\theta)$ and $(\chi^{\pr},\theta^{\pr})$
is the pair
$(\chi\chi^{\pr},\phi^{\pr!}(\theta)\theta^{\pr})\colon 
(\phi\phi^{\pr},\phi^{\pr!}(f)f^{\pr})\cong (\psi\psi^{\pr},\psi^{\pr!}(g)g^{\pr})$
consisting of
$\chi\chi^{\pr}\colon \phi\phi^{\pr}\cong \psi\psi^{\pr}$
and the horizontal composite $\phi^{\pr!}(\theta)\theta^{\pr}$ of
$\phi^{\pr!}(\theta)\colon \phi^{\pr!}(f)\cong \phi^{\pr!}(\chi_{\pi}^{-1})\phi^{\pr!}(g)$
and
$\theta^{\pr}\colon f^{\pr}\cong (\chi^{\pr}_{\pi^{\pr}})^{-1}g^{\pr}$.

\begin{remark}
  \label{rem-PKGalobj-etalebasechange}
  When we define a morphism 
  $(\phi,f)\colon (\kos{S}^{\pr},\pi^{\pr})\to (\kos{S},\pi)$
  of prekosmic Galois objects,
  we require $\phi$ to be \'{e}tale,
  for then $\phi^!(\pi)$ is
  a Galois object in $\kos{S}^{\pr}$.
  Indeed, the functor
  $\phi_!(\phi^!(\pi)\otimes^{\pr}\slot)$
  preserves reflexive coequalizers,
  as it is isomorphic to $\pi\otimes \phi_!(\slot)$
  by \eqref{eq-Galprekos-projformula},
  and since $\phi_!$ is conservative and preserves reflexive coequalizers,
  we deduce that the functor $\phi^!(\pi)\otimes^{\pr} \slot$
  preserves reflexive coequalizers.
  The argument still applies
  when we are merely given 
  a strong monoidal functor $\phi^!\colon \kos{S}\to\kos{S}^{\pr}$,
  a functor $\phi_!\colon \CS^{\pr}\to \CS$ that is conservative
  and preserves reflexive coequalizers,
  and a natural isomorphism
  $\varphi\colon
  \phi_!(\slot\otimes^{\pr}\phi^!(\slot))
  \cong \phi_!(\slot)\otimes \slot$.
  Thus such data $(\phi^!,\phi_!,\varphi)$,
  in which $\phi^!$ need not be right adjoint to $\phi_!$,
  could serve as $1$-cells of
  $\mathsf{pGalKos}_{\etsf}$,
  with $\mathsf{pGalObj}$ and $\mathsf{pGalCat}_{\mathsf{pt}}$
  enlarged accordingly.
  We do not pursue this extension in this paper.
\end{remark}

We denote by 
$\mathsf{p}_{\mathsf{obj}}\colon \mathsf{pGalObj}\to\mathsf{pGalKos}_{\etsf}$
the $2$-functor that sends each cell to its first component.
It is immediate that $\mathsf{pGalObj}$
is a split fibered $2$-category over $\mathsf{pGalKos}_{\etsf}$
with respect to $\mathsf{p}_{\mathsf{obj}}$,
whose cleavage is given as follows.
\begin{itemize}
  \item 
  For a $0$-cell $(\kos{S},\pi)$ in $\mathsf{pGalObj}$
  and a $1$-cell $\phi\colon \kos{S}^{\pr}\mor\kos{S}$ in $\mathsf{pGalKos}_{\etsf}$,
  the cartesian $1$-cell
  $(\kos{S}^{\pr},\phi^!(\pi))\to (\kos{S},\pi)$
  over $\phi$
  is the pair consisting of 
  $\phi$ and the identity morphism of $\phi^!(\pi)$.

  \item
  For a $1$-cell $(\psi,g)\colon (\kos{S}^{\pr},\pi^{\pr})\to (\kos{S},\pi)$
  in $\mathsf{pGalObj}$
  and a $2$-cell 
  $\chi\colon \phi\cong \psi$
  between $1$-cells $\phi$, $\psi\colon \kos{S}^{\pr}\mor\kos{S}$
  in $\mathsf{pGalKos}_{\etsf}$,
  the cartesian $2$-cell
  $(\phi,\chi_{\pi}^{-1}g)\cong (\psi,g)$
  over $\chi$
  is the pair consisting of 
  $\chi$ and the identity $2$-isomorphism
  of $\chi^{-1}_{\pi}g\colon \pi^{\pr}\to \phi^!(\pi)$.
\end{itemize}
The strict fiber of $\mathsf{p}_{\mathsf{obj}}$
over each Galois prekosmos $\kos{S}$
is $2$-isomorphic to $\mathsf{pGalObj}(\kos{S})$,
defined in \eqref{eq-Galobjinpk-pGalObjS21cat}.
We illustrate $\mathsf{p}_{\mathsf{obj}}$ as follows.
\begin{equation}
  \label{eq-PKGalobj-pobj21funct}
  \begin{tikzcd}[row sep=22pt,column sep=20pt]
    \mathsf{pGalObj}
    \ar[d,"{\mathsf{p}_{\mathsf{obj}}}"']
    &(\kos{S}^{\pr},\pi^{\pr})
    \ar[rr,shift left=0ex,bend left=1.5pc,"{(\phi,f)}" description]
    \ar[rr,shift right=0ex,bend right=1.5pc,"{(\psi,g)}" description]
    &\text{ }
    \ar[d,2cell,pos=-0.2,"{\,{(\chi,\theta)}}"]
    \ar[d,2cell,pos=1,"{\,\,\chi}"]
    &(\kos{S},\pi)
    \\%1
    {\mathsf{pGalKos}_{\etsf}}
    &\kos{S}^{\pr}
    \ar[rr,filled,shift left=0.3ex,bend left=1.5pc,"\phi" description]
    \ar[rr,filled,shift right=0.3ex,bend right=1.5pc,"\psi" description]
    &\text{ }
    &\kos{S}
  \end{tikzcd}
\end{equation}

\section{Change of base along a Galois morphism}
\label{sec-changebaseGal}

In this section, we introduce the notation for
changing the base Galois prekosmos,
which we use throughout the rest of this chapter.
Let $\phi\colon \kos{S}^{\pr}\mor\kos{S}$
be a Galois morphism between Galois prekosmoi
$\kos{S}$, $\kos{S}^{\pr}$.
We write
$(\slot)_{\phi}\colon \mathtt{pGalKos}/\kos{S}^{\pr}\to \mathtt{pGalKos}/\kos{S}$
for the change of base $2$-functor
induced by $\phi$
between the pseudoslice $2$-categories.
We explain the $2$-functor $(\slot)_{\phi}$ in detail.
\begin{itemize}
  \item 
  Let $\kos{T}^{\pr}\xmor{\kos{t}^{\pr}}\kos{S}^{\pr}$
  be a Galois $\kos{S}^{\pr}$-prekosmos.
  When we regard the composite 
  $\phi\kos{t}^{\pr}\colon \kos{T}^{\pr}\mor \kos{S}$
  as a Galois $\kos{S}$-prekosmos,
  we write it as $\kos{T}^{\pr}\xmor{\kos{t}^{\pr}_{\phi}}\kos{S}$,
  or simply as $\kos{t}^{\pr}_{\phi}$.

  \item
  In particular,
  when we regard $\phi$
  as a Galois $\kos{S}$-prekosmos,
  we write it as $\kos{S}^{\pr}\xmor{\phi}\kos{S}$, or simply as $\phi$.
  
  \item
  Let $\kos{f}^{\pr}\colon \tilde{\kos{t}}^{\pr}\mor \kos{t}^{\pr}$
  be a Galois $\kos{S}^{\pr}$-morphism
  between Galois $\kos{S}^{\pr}$-prekosmoi
  $\kos{t}^{\pr}$, $\tilde{\kos{t}}^{\pr}$;
  it consists of a Galois morphism
  $\kos{f}^{\pr}\colon \tilde{\kos{T}}^{\pr}\mor\kos{T}^{\pr}$
  and an invertible Galois transformation
  $\what{\kos{f}^{\pr}}\colon \kos{t}^{\pr}\kos{f}^{\pr}\cong \tilde{\kos{t}}^{\pr}$.
  We write
  $\kos{f}^{\pr}_{\phi}\colon \tilde{\kos{t}}^{\pr}_{\phi}\mor \kos{t}^{\pr}_{\phi}$
  for the Galois $\kos{S}$-morphism
  whose underlying Galois morphism is that of
  $\kos{f}^{\pr}\colon \tilde{\kos{t}}^{\pr}\mor \kos{t}^{\pr}$
  and whose invertible Galois transformation
  $\what{\kos{f}^{\pr}_{\phi}}\colon 
  \kos{t}^{\pr}_{\phi}\kos{f}^{\pr}_{\phi}\cong \tilde{\kos{t}}^{\pr}_{\phi}$
  is $\phi\what{\kos{f}^{\pr}}\colon 
  \phi\kos{t}^{\pr}\kos{f}^{\pr}\cong \phi\tilde{\kos{t}}^{\pr}$.
  We illustrate
  $\kos{f}^{\pr}_{\phi}\colon \tilde{\kos{t}}^{\pr}_{\phi}\mor\kos{t}^{\pr}_{\phi}$
  and its description as follows.
  \begin{equation*}
    \begin{tikzcd}[row sep=7pt,column sep=7pt]
      \tilde{\kos{T}}^{\pr}
      \ar[rr,filled,"\kos{f}^{\pr}_{\phi}"]
      \ar[ddr,filled,"\tilde{\kos{t}}^{\pr}_{\phi}"']
      &\text{ }
      &\kos{T}^{\pr}
      \ar[ddl,filled,"\kos{t}^{\pr}_{\phi}"]
      \\%1
      \text{ }
      &\text{ }
      &\text{ }
      \ar[ll,2cell,shift left=0.9ex,"{\what{\kos{f}^{\pr}_{\phi}}}"']
      \\%2
      \text{ }
      &\kos{S}
      &\text{ }
    \end{tikzcd}
    \,\,=\,\,
    \begin{tikzcd}[row sep=3pt,column sep=8pt]
      \tilde{\kos{T}}^{\pr}
      \ar[rr,filled,"\kos{f}^{\pr}"]
      \ar[ddr,filled,"\tilde{\kos{t}}^{\pr}"']
      &\text{ }
      &\kos{T}^{\pr}
      \ar[ddl,filled,"\kos{t}^{\pr}"]
      \\%1
      \text{ }
      &\text{ }
      &\text{ }
      \ar[ll,2cell,shift left=1ex,pos=0.5,"{\,\,\,\what{\kos{f}^{\pr}}}"']
      \\%2
      \text{ }
      &\kos{S}^{\pr}
      \ar[dd,filled,"\phi"]
      &\text{ }
      \\%3
      \text{ }
      &\text{ }
      &\text{ }
      \\%4
      \text{ }
      &\kos{S}
      &\text{ }
    \end{tikzcd}
  \end{equation*}

  \item
  Let $\kos{t}^{\pr}$,
  $\tilde{\kos{t}}^{\pr}$
  be Galois $\kos{S}^{\pr}$-prekosmoi,
  $\kos{f}^{\pr}$, $\kos{g}^{\pr}\colon \tilde{\kos{t}}^{\pr}\mor \kos{t}^{\pr}$
  be Galois $\kos{S}^{\pr}$-morphisms,
  and $\vartheta^{\pr}\colon \kos{f}^{\pr}\Rightarrow \kos{g}^{\pr}$
  be a Galois $\kos{S}^{\pr}$-transformation.
  Then $\vartheta^{\pr}$ is a Galois $\kos{S}$-transformation
  from $\kos{f}^{\pr}_{\phi}$ to $\kos{g}^{\pr}_{\phi}$,
  and we write it as
  $\vartheta^{\pr}_{\phi}\colon \kos{f}^{\pr}_{\phi}\Rightarrow\kos{g}^{\pr}_{\phi}$.
\end{itemize}
The pair $(\phi,\what{\phi})$
consisting of $\phi$
and the equality $\what{\phi}\colon \phi=\phi$
is a Galois $\kos{S}$-morphism from $\phi$ to $\vphantom{\id}\smash[t]{\id}_{\kos{S}\vphantom{!}}^{\vphantom{!}}$;
we omit $\what{\phi}$ from the notation and write
$\phi\colon \phi\mor\vphantom{\id}\smash[t]{\id}_{\kos{S}\vphantom{!}}^{\vphantom{!}}$.

Let 
$\psi\colon \kos{S}^{\pr}\mor\kos{S}$
be a Galois morphism
and $\chi\colon \phi\cong \psi$ be an invertible Galois transformation.
We explain the $2$-natural isomorphism
$(\slot)_{\chi}$
induced by $\chi$
from $(\slot)_{\phi}$ to $(\slot)_{\psi}$.
\begin{itemize}
  \item 
  Let $\kos{T}^{\pr}\xmor{\kos{t}^{\pr}}\kos{S}^{\pr}$
  be a Galois $\kos{S}^{\pr}$-prekosmos.
  We write $\kos{t}^{\pr}_{\chi}\colon \kos{t}^{\pr}_{\phi}\mor \kos{t}^{\pr}_{\psi}$
  for the Galois $\kos{S}$-morphism
  whose underlying Galois morphism is
  the identity Galois morphism $\vphantom{\id}\smash[t]{\id}_{\kos{T}^{\pr}\vphantom{!}}^{\vphantom{!}}$ of $\kos{T}^{\pr}$
  and whose invertible Galois transformation
  $\what{\kos{t}^{\pr}_{\chi}}$ is
  $\chi^{-1}\kos{t}^{\pr}\colon \psi\kos{t}^{\pr}\cong \phi\kos{t}^{\pr}$.
  We illustrate $\kos{t}^{\pr}_{\chi}\colon \kos{t}^{\pr}_{\phi}\mor\kos{t}^{\pr}_{\psi}$
  and its right adjoint part
  $\kos{t}^{\pr!}_{\chi}\colon \kos{t}^{\pr!}_{\psi}\to \kos{t}^{\pr!}_{\phi}$,
  as well as their descriptions as follows.
  \begin{equation*}
    \begin{tikzcd}[row sep=5pt,column sep=5pt]
      \kos{T}^{\pr}
      \ar[rr,filled,"\kos{t}^{\pr}_{\chi}"]
      \ar[ddr,filled,"\kos{t}^{\pr}_{\phi}"']
      &\text{ }
      &\kos{T}^{\pr}
      \ar[ddl,filled,"\kos{t}^{\pr}_{\psi}"]
      \\%1
      \text{ }
      &\text{ }
      &\text{ }
      \ar[ll,2cell,shift left=1ex,"{\!\!\what{\kos{t}^{\pr}_{\chi}}}"']
      \\%2
      \text{ }
      &\kos{S}
      &\text{ }
    \end{tikzcd}
    =\,\,
    \begin{tikzcd}[row sep=2pt,column sep=1pt]
      \kos{T}^{\pr}
      \ar[dd,filled,"\kos{t}^{\pr}"']
      \ar[rr,equal]
      &\text{ }
      &\kos{T}^{\pr}
      \ar[dd,filled,"\kos{t}^{\pr}"]
      \\%1
      \text{ }
      &\text{ }
      &\text{ }
      \\%2
      \kos{S}^{\pr}
      \ar[ddr,filled,"\phi"']
      \ar[rr,equal]
      &\text{ }
      &\kos{S}^{\pr}
      \ar[ddl,filled,"\psi"]
      \\%3
      \text{ }
      &\text{ }
      &\text{ }
      \ar[ll,2cell,shift left=0.6ex,"{\,\,\,\,\chi^{-1}}"']
      \\%4
      \text{ }
      &\kos{S}
      &\text{ }
    \end{tikzcd}
    \qquad
    \begin{tikzcd}[row sep=5pt,column sep=5pt]
      \kos{T}^{\pr}
      &\text{ }
      &\kos{T}^{\pr}
      \ar[ll,"\kos{t}^{\pr!}_{\chi}"']
      \\%1
      \text{ }
      &\text{ }
      &\text{ }
      \ar[ll,2cell,shift left=1ex,"{\!\what{\kos{t}^{\pr!}_{\chi}}}"']
      \\%2
      \text{ }
      &\kos{S}
      \ar[uul,"\kos{t}^{\pr!}_{\phi}"]
      \ar[uur,"\kos{t}^{\pr!}_{\psi}"']
      &\text{ }
    \end{tikzcd}
    =\,\,
    \begin{tikzcd}[row sep=2pt,column sep=1pt]
      \kos{T}^{\pr}
      \ar[rr,equal]
      &\text{ }
      &\kos{T}^{\pr}
      \\%1
      \text{ }
      &\text{ }
      &\text{ }
      \\%2
      \kos{S}^{\pr}
      \ar[uu,"\kos{t}^{\pr!}"]
      \ar[rr,equal]
      &\text{ }
      &\kos{S}^{\pr}
      \ar[uu,"\kos{t}^{\pr!}"']
      \\%3
      \text{ }
      &\text{ }
      &\text{ }
      \ar[ll,2cell,shift left=0.6ex,"{\,\,\,\,\chi^{-1}}"']
      \\%4
      \text{ }
      &\kos{S}
      \ar[uul,"\phi^!"]
      \ar[uur,"\psi^!"']
      &\text{ }
    \end{tikzcd}
  \end{equation*}

  \item
  In particular, when $\kos{t}^{\pr}=\vphantom{\id}\smash[t]{\id}_{\kos{S}^{\pr}\vphantom{!}}^{\vphantom{!}}$,
  we write $\id_{\chi}^{\pr}\colon \phi\mor\psi$
  for the Galois $\kos{S}$-morphism
  consisting of the identity Galois morphism $\vphantom{\id}\smash[t]{\id}_{\kos{S}^{\pr}\vphantom{!}}^{\vphantom{!}}$ of $\kos{S}^{\pr}$
  and the invertible Galois transformation
  $\chi^{-1}\colon \psi\cong\phi$.
  We illustrate $\id_{\chi}^{\pr}\colon \phi\mor\psi$ and
  its right adjoint part
  $\id_{\chi}^{\pr!}\colon \psi^!\to \phi^!$,
  as well as their descriptions as follows.
  \begin{equation*}
    \begin{tikzcd}[row sep=4pt,column sep=4pt]
      \kos{S}^{\pr}
      \ar[rr,filled,"\id_{\chi}^{\pr}"]
      \ar[ddr,filled,"\phi"']
      &\text{ }
      &\kos{S}^{\pr}
      \ar[ddl,filled,"\psi"]
      &\text{ }
      &\kos{S}^{\pr}
      \ar[ddr,filled,"\phi"']
      \ar[rr,equal]
      &\text{ }
      &\kos{S}^{\pr}
      \ar[ddl,filled,"\psi"]
      \\%1
      \text{ }
      &\text{ }
      &\text{ }
      \ar[ll,2cell,shift left=0.5ex,"{\what{\id_{\chi}^{\pr}}}"']
      &=
      &\text{ }
      &\text{ }
      &\text{ }
      \ar[ll,2cell,shift left=0ex,"{\,\,\,\,\chi^{-1}}"']
      \\%2
      \text{ }
      &\kos{S}
      &\text{ }
      &\text{ }
      &\text{ }
      &\kos{S}
      &\text{ }
    \end{tikzcd}
    \qquad
    \begin{tikzcd}[row sep=4pt,column sep=4pt]
      \kos{S}^{\pr}
      &\text{ }
      &\kos{S}^{\pr}
      \ar[ll,"\id_{\chi}^{\pr!}"']
      &\text{ }
      &\kos{S}^{\pr}
      \ar[rr,equal]
      &\text{ }
      &\kos{S}^{\pr}
      \\%1
      \text{ }
      &\text{ }
      &\text{ }
      \ar[ll,2cell,shift left=0.5ex,"{\what{\id_{\chi}^{\pr!}}}"']
      &=
      &\text{ }
      &\text{ }
      &\text{ }
      \ar[ll,2cell,shift left=0ex,"{\,\,\,\,\chi^{-1}}"']
      \\%2
      \text{ }
      &\kos{S}
      \ar[uul,"\phi^!"]
      \ar[uur,"\psi^!"']
      &\text{ }
      &\text{ }
      &\text{ }
      &\kos{S}
      \ar[uul,"\phi^!"]
      \ar[uur,"\psi^!"']
      &\text{ }
    \end{tikzcd}
  \end{equation*}
\end{itemize}
For each Galois $\kos{S}^{\pr}$-morphism
$\kos{f}^{\pr}\colon \tilde{\kos{t}}^{\pr}\mor\kos{t}^{\pr}$,
we have
$\kos{f}^{\pr}_{\psi}\tilde{\kos{t}}^{\pr}_{\chi}
=\kos{t}^{\pr}_{\chi}\kos{f}^{\pr}_{\phi}
\colon \tilde{\kos{t}}^{\pr}_{\phi}\mor\kos{t}^{\pr}_{\psi}$.

Recall from Section~\ref{sec-ActionbaseGal}
the definitions of a colax $\kos{S}$-tensor functor
and a colax $\kos{S}$-tensor natural transformation.
\begin{itemize}
  \item 
  The right adjoint part
  $\phi^!\colon \vphantom{\id}\smash[t]{\id}_{\kos{S}\vphantom{!}}^{\vphantom{!}}\to \phi^!$
  of $\phi\colon \phi\mor\vphantom{\id}\smash[t]{\id}_{\kos{S}\vphantom{!}}^{\vphantom{!}}$
  is the colax $\kos{S}$-tensor functor
  consisting of $\phi^!\colon \kos{S}\to \kos{S}^{\pr}$
  and $\what{\phi^!}\colon \phi^!=\phi^!$,
  and its $\kos{S}$-equivariance is
  $(\vecar{\phi^!})_{x,z}
  =(\phi^!_{\tensor})_{x,z}
  \colon \phi^!(x\otimes z)
  \cong \phi^!(x)\otimes^{\pr}\phi^!(z)$,
  where $x$, $z\in \CS$.
  
  \item
  Let $\kos{f}^{\pr}\colon \tilde{\kos{t}}^{\pr}\mor\kos{t}^{\pr}$
  be a Galois $\kos{S}^{\pr}$-morphism.
  We write
  $(\kos{f}^{\pr!})_{\phi}$
  and $(\kos{f}^{\pr}_!)_{\phi}$
  for the right and left adjoint parts
  of the Galois $\kos{S}$-morphism $\kos{f}^{\pr}_{\phi}$,
  both of which are colax $\kos{S}$-tensor functors;
  the notation indicates that they are obtained from 
  $\kos{f}^{\pr!}$ and $\kos{f}^{\pr}_!$
  by change of base along $\phi$.
  The relation
  $\kos{f}^{\pr}_{\psi}\tilde{\kos{t}}^{\pr}_{\chi}
  =\kos{t}^{\pr}_{\chi}\kos{f}^{\pr}_{\phi}$
  above is the equality
  $\tilde{\kos{t}}^{\pr!}_{\chi}(\kos{f}^{\pr!})_{\psi}
  =(\kos{f}^{\pr!})_{\phi}\kos{t}^{\pr!}_{\chi}
  \colon \kos{t}^{\pr!}_{\psi}\to \tilde{\kos{t}}^{\pr!}_{\phi}$
  between the right adjoint parts.
  
  \item
  Observe that
  $\chi\colon\phi^!\cong \id^{\pr!}_{\chi}\psi^!$
  is a colax $\kos{S}$-tensor natural isomorphism
  between colax $\kos{S}$-tensor functors
  from $\vphantom{\id}\smash[t]{\id}_{\kos{S}\vphantom{!}}^{\vphantom{!}}$ to $\phi^!$.

  \item
  Let $\mon{f}^{\pr}\colon \vphantom{\id}\smash[t]{\id}_{\kos{S}^{\pr}\vphantom{!}}^{\vphantom{!}}\to \vphantom{\id}\smash[t]{\id}_{\kos{S}^{\pr}\vphantom{!}}^{\vphantom{!}}$
  be a colax $\kos{S}^{\pr}$-tensor functor
  consisting of a colax monoidal functor
  $\mon{f}^{\pr}\colon \kos{S}^{\pr}\to \kos{S}^{\pr}$
  and a comonoidal natural transformation
  $\what{\mon{f}^{\pr}}\colon \mon{f}^{\pr}\Rightarrow\vphantom{\id}\smash[t]{\id}_{\kos{S}^{\pr}\vphantom{!}}^{\vphantom{!}}$.
  We write $\mon{f}^{\pr}_{\phi}\colon \phi^!\to \phi^!$
  for the colax $\kos{S}$-tensor functor
  whose underlying colax monoidal functor is that of 
  $\mon{f}^{\pr}\colon \vphantom{\id}\smash[t]{\id}_{\kos{S}^{\pr}\vphantom{!}}^{\vphantom{!}}\to \vphantom{\id}\smash[t]{\id}_{\kos{S}^{\pr}\vphantom{!}}^{\vphantom{!}}$
  and whose comonoidal natural transformation
  $\what{\mon{f}^{\pr}_{\phi}}\colon \mon{f}^{\pr}_{\phi}\phi^!\Rightarrow \phi^!$
  is $\what{\mon{f}^{\pr}}\phi^!\colon \mon{f}^{\pr}\phi^!\Rightarrow\phi^!$.
\end{itemize}
Let $\mon{f}^{\pr}\colon \vphantom{\id}\smash[t]{\id}_{\kos{S}^{\pr}\vphantom{!}}^{\vphantom{!}}\to \vphantom{\id}\smash[t]{\id}_{\kos{S}^{\pr}\vphantom{!}}^{\vphantom{!}}$
be a colax $\kos{S}^{\pr}$-tensor functor,
$\mon{g}\colon \vphantom{\id}\smash[t]{\id}_{\kos{S}\vphantom{!}}^{\vphantom{!}}\to \vphantom{\id}\smash[t]{\id}_{\kos{S}\vphantom{!}}^{\vphantom{!}}$
be a colax $\kos{S}$-tensor functor,
and $\vartheta\colon \mon{f}^{\pr}_{\phi}\phi^!\Rightarrow\phi^!\mon{g}$
be a colax $\kos{S}$-tensor natural transformation.
By taking
$\monbase{t}=\phi^!$,
$\monbase{t}^{\pr}=\vphantom{\id}\smash[t]{\id}_{\kos{S}\vphantom{!}}^{\vphantom{!}}$,
$X^{\pr}=s$ in Lemma~\ref{lem-ActionbaseGal-colaxStnatequiv},
we obtain the following commutative diagram for every $z\in \CS$.
\begin{equation}
  \label{eq-changebaseGal-colaxStnathatar}
  \vcenter{\hbox{
    \xymatrix@R=12pt@C=60pt{
      \mon{f}^{\pr}\phi^!(z)
      \ar[d]_-{\vartheta_z}
      \ar[r]^-{(\hatar{\mon{f}^{\pr}})_{\phi^!(z)}}
      &\mon{f}^{\pr}(s^{\pr})\otimes^{\pr} \phi^!(z)
      \ar[r]^-{\mon{f}^{\pr}(\phi^!_{\unit})^{-1}\otimes^{\pr} I_{\phi^!(z)}}_-{\cong}
      &\mon{f}^{\pr}\phi^!(s)\otimes^{\pr} \phi^!(z)
      \ar[d]^-{\vartheta_{s}\otimes^{\pr} I_{\phi^!(z)}}
      \\%1
      \phi^!\mon{g}(z)
      \ar[r]^-{\phi^!((\hatar{\mon{g}})_z)}
      &\phi^!(\mon{g}(s)\otimes z)
      \ar[r]^-{(\phi^!_{\tensor})_{\mon{g}(s),z}}_-{\cong}
      &\phi^!\mon{g}(s)\otimes^{\pr}\phi^!(z)
    }
  }}
\end{equation}
Observe that when $\phi=\vphantom{\id}\smash[t]{\id}_{\kos{S}\vphantom{!}}^{\vphantom{!}}$,
the commutative diagram \eqref{eq-changebaseGal-colaxStnathatar}
reduces to the commutative diagram \eqref{eq-ActionbaseGal-colaxStnathatar}.

\section{Prekosmic fundamental morphisms}
\label{sec-PKGalfund}

\begin{definition}
  \label{def-PKGalfund-pkfundmor}
  Let $(\kos{S},\kos{t})$, $(\kos{S}^{\pr},\kos{t}^{\pr})$
  be prekosmic Galois categories.
  A \emph{prekosmic fundamental morphism}
  $(\phi,\kos{f})\colon(\kos{S}^{\pr},\kos{t}^{\pr})\mor (\kos{S},\kos{t})$,
  denoted by an arrow with a filled tip,
  is a pair consisting of an \'{e}tale Galois morphism
  $\phi\colon \kos{S}^{\pr}\mor\kos{S}$
  and a Galois $\kos{S}$-morphism
  $\kos{f}\colon \kos{t}^{\pr}_{\phi}\mor \kos{t}$
  such that there exist
  a surjective fiber point $\omega\colon\vphantom{\id}\smash[t]{\id}_{\kos{S}\vphantom{!}}^{\vphantom{!}}\mor\kos{t}$ of $\kos{t}$,
  a surjective fiber point $\omega^{\pr}\colon\vphantom{\id}\smash[t]{\id}_{\kos{S}^{\pr}\vphantom{!}}^{\vphantom{!}}\mor\kos{t}^{\pr}$ of $\kos{t}^{\pr}$,
  and an invertible Galois $\kos{S}$-transformation
  $\kos{f}\omega^{\pr}_{\phi}\cong \omega\phi$.
  We illustrate 
  $(\phi,\kos{f})$
  and 
  $\kos{f}\omega^{\pr}_{\phi}\cong \omega\phi$
  as follows.
  \begin{equation*}
    \begin{tikzcd}[row sep=5pt,column sep=20pt]
      \kos{T}^{\pr}
      \ar[dd,filled,"\kos{t}^{\pr}"']
      \ar[r,filled,"\kos{f}"]
      &\kos{T}
      \ar[dd,filled,"\kos{t}"]
      &\text{ }
      &\phi
      \ar[dd,filled,"\omega^{\pr}_{\phi}"']
      \ar[r,filled,"\phi"]
      &\id_{\kos{S}}
      \ar[dd,filled,"\omega"]
      \\%1
      \text{ }
      &\text{ }
      \ar[l,2cell,shift left=0.8ex,"{\what{\kos{f}}}"']
      &\text{ }
      &\text{ }
      \ar[r,2cell,shift right=0.7ex,pos=0.5,"{\!\exists}"]
      &\text{ }
      \\%2
      \kos{S}^{\pr}
      \ar[r,filled,"{\phi}"]
      &\kos{S}
      &\text{ }
      &\kos{t}^{\pr}_{\phi}
      \ar[r,filled,"{\kos{f}}"]
      &\kos{t}
    \end{tikzcd}
  \end{equation*}
\end{definition}

\begin{remark}
  \label{rem-PKGalfund-fundmordef}
  By taking $\phi=\vphantom{\id}\smash[t]{\id}_{\kos{S}\vphantom{!}}^{\vphantom{!}}$ in Definition~\ref{def-PKGalfund-pkfundmor},
  we recover the definition of a fundamental $\kos{S}$-morphism,
  given in Definition~\ref{def-fundmoroverGalpk-fundmor}.
  Indeed, let $\kos{S}$ be a Galois prekosmos,
  $\kos{t}$, $\kos{t}^{\pr}$ be Galois $\kos{S}$-categories,
  and
  $\kos{f}\colon \kos{t}^{\pr}\mor\kos{t}$
  be a Galois $\kos{S}$-morphism.
  Then $(\vphantom{\id}\smash[t]{\id}_{\kos{S}\vphantom{!}}^{\vphantom{!}},\kos{f})\colon(\kos{S},\kos{t}^{\pr})\mor (\kos{S},\kos{t})$
  is a prekosmic fundamental morphism
  if and only if
  $\kos{f}\colon \kos{t}^{\pr}\mor\kos{t}$
  is a fundamental $\kos{S}$-morphism.
\end{remark}

\begin{lemma}
  \label{lem-PKGalfund-inducedPKfundmor}
  Let 
  $(\phi,f)\colon (\kos{S}^{\pr},\pi^{\pr})\to (\kos{S},\pi)$
  be a morphism of prekosmic Galois objects.
  We have a Galois $\kos{S}$-morphism
  $\kos{f}\colon (\kos{t}^{\pr}_{\pi^{\pr}})_{\phi}\mor\kos{t}_{\pi\vphantom{!}}^{\vphantom{!}}$
  and an invertible Galois $\kos{S}$-transformation
  $\kos{f}(\kos{u}^{\pr}_{\pi^{\pr}})_{\phi}\cong \kos{u}_{\pi\vphantom{!}}^{\vphantom{!}}\phi$,
  so that the pair $(\phi,\kos{f})\colon (\kos{S}^{\pr},\kos{t}^{\pr}_{\pi^{\pr}})\mor (\kos{S},\kos{t}_{\pi\vphantom{!}}^{\vphantom{!}})$
  is a prekosmic fundamental morphism.
  \begin{equation*}
    \begin{tikzcd}[row sep=5pt,column sep=15pt]
      \kosRep_{\kos{S}^{\pr}}(\pi^{\pr})
      \ar[dd,filled,"\kos{t}^{\pr}_{\pi^{\pr}}"']
      \ar[r,filled,"\kos{f}"]
      &\kosRep_{\kos{S}}(\pi)
      \ar[dd,filled,"\kos{t}_{\pi}"]
      &\text{ }
      &\phi
      \ar[dd,filled,"(\kos{u}^{\pr}_{\pi^{\pr}})_{\phi}"']
      \ar[r,filled,"\phi"]
      &\id_{\kos{S}}
      \ar[dd,filled,"\kos{u}_{\pi}"]
      \\%1
      \text{ }
      &\text{ }
      \ar[l,2cell,shift left=0.8ex,"{\what{\kos{f}}}"']
      &\text{ }
      &\text{ }
      \ar[r,2cell,shift right=0.7ex,pos=0.5]
      &\text{ }
      \\%2
      \kos{S}^{\pr}
      \ar[r,filled,"{\phi}"]
      &\kos{S}
      &\text{ }
      &(\kos{t}^{\pr}_{\pi^{\pr}})_{\phi}
      \ar[r,filled,"{\kos{f}}"]
      &\kos{t}_{\pi}
    \end{tikzcd}
  \end{equation*}
\end{lemma}
\begin{proof}
  We write 
  $X=(x,\gamma_x)$ for an object in $\cat{Rep}_{\kos{S}}(\pi)$
  and
  $Y^{\pr}=(y^{\pr},\gamma^{\pr}_{y^{\pr}})$
  for an object in $\cat{Rep}_{\kos{S}^{\pr}}(\pi^{\pr})$.
  Put
  $\kos{u}\coloneqq \kos{u}_{\pi\vphantom{!}}^{\vphantom{!}}$
  and
  $\kos{u}^{\pr}\coloneqq \kos{u}^{\pr}_{\pi^{\pr}}$,
  and write 
  $\kos{u}=(\kos{u}^!,\kos{u}_!,\eta,\epsilon)$,
  $\kos{u}^{\pr}=(\kos{u}^{\pr!},\kos{u}^{\pr}_!,\eta^{\pr},\epsilon^{\pr})$,
  and $\phi=(\phi^!,\phi_!,\eta^{\phi},\epsilon^{\phi})$.
  The right adjoint part
  $\kos{f}^!\colon \kosRep_{\kos{S}}(\pi)\to \kosRep_{\kos{S}^{\pr}}(\pi^{\pr})$
  is a strong monoidal functor,
  which sends each $X$ to the pair
  $\kos{f}^!(X)=(\phi^!(x),
  \phi^!(\gamma_x)
  \circ (\phi^!_{\tensor})^{-1}_{\pi,x}
  \circ (f\otimes^{\pr}I_{\phi^!(x)}))$,
  and whose coherence isomorphisms are
  given by those of $\phi^!\colon \kos{S}\to \kos{S}^{\pr}$
  so that we have an equality
  $\kos{u}^{\pr!}\kos{f}^!=\phi^!\kos{u}^!$
  of strong monoidal functors.
  We write
  $\zeta\colon \phi_!(\phi^!(\pi)\otimes^{\pr}\slot)\cong \pi\otimes \phi_!(\slot)$
  for the natural isomorphism
  whose component at $z^{\pr}\in \CS^{\pr}$ is
  $(\epsilon^{\phi}_{\pi}\otimes I_{\phi_!(z^{\pr})})
  \circ (\phi_{!\tensor})_{\phi^!(\pi),z^{\pr}}$;
  observe that $\zeta$ agrees with $\varphi^{\phi}_{\slot,\pi}
  \colon \phi_!(\slot\otimes^{\pr}\phi^!(\pi))\cong \phi_!(\slot)\otimes \pi$
  up to symmetric coherence isomorphisms,
  where $\varphi^{\phi}$
  is the natural isomorphism associated to $\phi$
  in \eqref{eq-Galprekos-projformula}.
  For each $Y^{\pr}$,
  we fix a coequalizer $\kos{f}_!(Y^{\pr})$ in $\cat{Rep}_{\kos{S}}(\pi)$
  of the reflexive pair
  $I_{\pi}\otimes \phi_!(\gamma^{\pr}_{y^{\pr}})$,
  $(\pc_{\pi}\otimes I_{\phi_!(y^{\pr})})
  \circ (I_{\pi}\otimes \zeta_{y^{\pr}})
  \circ (I_{\pi}\otimes \phi_!(f\otimes^{\pr}I_{y^{\pr}}))
  \colon
  (\pi\otimes \phi_!(\pi^{\pr}\otimes^{\pr}y^{\pr}),
  \pc_{\pi}\otimes I_{\phi_!(\pi^{\pr}\otimes^{\pr}y^{\pr})})
  \to
  (\pi\otimes \phi_!(y^{\pr}),
  \pc_{\pi}\otimes I_{\phi_!(y^{\pr})})$,
  with common section
  $I_{\pi}\otimes \phi_!(u^{\pr}_{\pi^{\pr}}\otimes^{\pr}I_{y^{\pr}})$.
  Then $\kos{f}_!$ is left adjoint to $\kos{f}^!$
  by the adjoint lifting theorem
  \cite[Theorem~4.5.6 and Exercise~4.8.6]{Borceux1994}.
  Hence $\kos{f}=(\kos{f}_!\dashv\kos{f}^!)$
  is a Galois morphism from
  $\kosRep_{\kos{S}^{\pr}}(\pi^{\pr})$ to $\kosRep_{\kos{S}}(\pi)$.
  The invertible Galois transformation
  $\what{\kos{f}}\colon \kos{t}_{\pi\vphantom{!}}^{\vphantom{!}}\kos{f}\cong \phi\kos{t}^{\pr}_{\pi^{\pr}}$
  is the equality
  $\kos{f}^!\kos{t}_{\pi\vphantom{!}}^!=\kos{t}^{\pr!}_{\pi^{\pr}}\phi^!$
  between the right adjoint parts,
  and the invertible Galois $\kos{S}$-transformation
  $\kos{f}\kos{u}^{\pr}_{\phi}\cong \kos{u}\phi$
  is the equality
  $(\kos{u}^{\pr!})_{\phi}\kos{f}^!=\phi^!\kos{u}^!$
  between the right adjoint parts
  as colax $\kos{S}$-tensor functors.
\end{proof}

\begin{proposition}
  \label{prop-PKGalfund-inducedPKfundmor}
  Let 
  $(\phi,\kos{f})\colon (\kos{S}^{\pr},\kos{t}^{\pr})\mor(\kos{S},\kos{t})$
  be a prekosmic fundamental morphism.
  Choose a surjective fiber point $\omega$ of $\kos{t}$,
  a surjective fiber point $\omega^{\pr}$ of $\kos{t}^{\pr}$,
  and an invertible Galois $\kos{S}$-transformation
  $\alpha\colon \kos{f}\omega^{\pr}_{\phi}\cong \omega\phi$.
  Put $\pi\coloneqq \omega^!\omega_!(s)$
  and
  $\pi^{\pr}\coloneqq \omega^{\pr!}\omega^{\pr}_!(s^{\pr})$,
  and recall the equivalences
  $\breve{\omega}\colon \kos{t}_{\pi\vphantom{!}}^{\vphantom{!}}\simeq\kos{t}$,
  $\breve{\omega}^{\pr}\colon \kos{t}^{\pr}_{\pi^{\pr}}\simeq\kos{t}^{\pr}$
  in Theorem~\ref{thm-Galcatoverpk-mainthm}.
  Then we have a morphism
  $(\phi,\breve{f})\colon (\kos{S}^{\pr},\pi^{\pr})\to (\kos{S},\pi)$
  of prekosmic Galois objects
  and an invertible Galois $\kos{S}$-transformation
  $\breve{\alpha}\colon \kos{f}\breve{\omega}^{\pr}_{\phi}\cong \breve{\omega}\breve{\kos{f}}$
  such that the following equality of pasting diagrams holds,
  where $\breve{\kos{f}}\colon (\kos{t}^{\pr}_{\pi^{\pr}})_{\phi}\mor\kos{t}_{\pi\vphantom{!}}^{\vphantom{!}}$
  is the Galois $\kos{S}$-morphism
  associated to $(\phi,\breve{f})$
  in Lemma~\ref{lem-PKGalfund-inducedPKfundmor}.
  \begin{equation*}
    \begin{tikzcd}[row sep=4pt,column sep=8pt]
      \phi
      \ar[rr,filled,"\phi"]
      \ar[dd,filled,"{(\kos{u}^{\pr}_{\pi^{\pr}})_{\phi}}"']
      &\text{ }
      &\id_{\kos{S}}
      \ar[dd,filled,"{\kos{u}_{\pi}}"]
      \ar[ddddrr,bend left=2.5pc,filled,"{\omega}"]
      &\text{ }
      &\text{ }
      &\text{ }
      &\phi
      \ar[rr,filled,"\phi"]
      \ar[dd,filled,"{(\kos{u}^{\pr}_{\pi^{\pr}})_{\phi}}"']
      \ar[ddddrr,filled,bend left=2.5pc,"{\omega^{\pr}_{\phi}}" description]
      &\text{ }
      &\id_{\kos{S}}
      \ar[ddddrr,bend left=2.5pc,filled,"{\omega}"]
      &\text{ }
      &\text{ }
      \\%1
      \text{ }
      \ar[rr,2cell,shift right=0.3ex,pos=0.55]
      &\text{ }
      &\text{ }
      &\text{ }
      &\text{ }
      &\text{ }
      &\text{ }
      &\text{ }
      &\text{ }
      &\text{ }
      &\text{ }
      \\%2
      (\kos{t}^{\pr}_{\pi^{\pr}})_{\phi}
      \ar[rr,filled,"\breve{\kos{f}}"]
      \ar[ddrr,filled,"\breve{\omega}^{\pr}_{\phi}"',"\simeq"]
      &\text{ }
      &\kos{t}_{\pi}
      \ar[ddrr,filled,"\breve{\omega}","\simeq"']
      \ar[rr,2cell,shift left=0ex,pos=0.4]
      &\text{ }
      &\text{ }
      &=
      &(\kos{t}^{\pr}_{\pi^{\pr}})_{\phi}
      \ar[ddrr,filled,"\breve{\omega}^{\pr}_{\phi}"',"\simeq"]
      \ar[rr,2cell,shift left=0ex,pos=0.3]
      &\text{ }
      &\text{ }
      \ar[rr,2cell,shift left=0ex,pos=0.25,"{\!\alpha}"]
      &\text{ }
      &\text{ }
      \\%3
      \text{ }
      &\text{ }
      \ar[rr,2cell,shift right=1ex,"{\!\breve{\alpha}}"]
      &\text{ }
      &\text{ }
      &\text{ }
      &\text{ }
      &\text{ }
      &\text{ }
      &\text{ }
      &\text{ }
      &\text{ }
      \\%4
      \text{ }
      &\text{ }
      &\kos{t}^{\pr}_{\phi}
      \ar[rr,filled,"\kos{f}"']
      &\text{ }
      &\kos{t}
      &\text{ }
      &\text{ }
      &\text{ }
      &\kos{t}^{\pr}_{\phi}
      \ar[rr,filled,"\kos{f}"']
      &\text{ }
      &\kos{t}
    \end{tikzcd}
  \end{equation*}
\end{proposition}
\begin{proof}
  We write
  $\alpha_!\colon 
  (\omega^{\pr}_!)_{\phi}\phi^!\Rightarrow \kos{f}^!\omega_!$
  for the mate of
  $\alpha^{-1}
  \colon \phi^!\omega^!
  \cong (\omega^{\pr!})_{\phi}\kos{f}^!$
  with respect to the adjunctions
  $\omega=(\omega_!\dashv\omega^!)$ and
  $\omega^{\pr}_{\phi}=((\omega^{\pr}_!)_{\phi}\dashv(\omega^{\pr!})_{\phi})$;
  it is the unique colax $\kos{S}$-tensor natural transformation
  that satisfies the following relations.
  \begin{equation}\label{eq-PKGalfund-matealpha}
    \vcenter{\hbox{
      \xymatrix@R=12pt@C=30pt{
        \phi^!
        \ar@2{->}[r]^-{\phi^!\eta}
        \ar@2{->}[d]_-{\eta^{\pr}\phi^!}
        &\phi^!\omega^!\omega_!
        \ar@2{->}[d]^-{\alpha^{-1}\omega_!}_-{\cong}
        &(\omega^{\pr}_!)_{\phi}\phi^!\omega^!
        \ar@2{->}[r]^-{\alpha_!\omega^!}
        \ar@2{->}[d]_-{(\omega^{\pr}_!)_{\phi}\alpha^{-1}}^-{\cong}
        &\kos{f}^!\omega_!\omega^!
        \ar@2{->}[d]^-{\kos{f}^!\epsilon}
        \\%1
        (\omega^{\pr!})_{\phi}(\omega^{\pr}_!)_{\phi}\phi^!
        \ar@2{->}[r]^-{(\omega^{\pr!})_{\phi}\alpha_!}
        &(\omega^{\pr!})_{\phi}\kos{f}^!\omega_!
        &(\omega^{\pr}_!)_{\phi}(\omega^{\pr!})_{\phi}\kos{f}^!
        \ar@2{->}[r]^-{\epsilon^{\pr}\kos{f}^!}
        &\kos{f}^!
      }
    }}
  \end{equation}
  Put
  $\alpha^\sharp
  \coloneqq (\alpha\omega_!)
  \circ ((\omega^{\pr!})_{\phi}\alpha_!)
  \colon (\omega^{\pr!})_{\phi}(\omega^{\pr}_!)_{\phi}\phi^!
  \Rightarrow \phi^!\omega^!\omega_!$
  and define 
  $\breve{f}\colon \pi^{\pr}\to \phi^!(\pi)$
  as
  $\alpha^{\sharp}_{s}\circ \omega^{\pr!}\omega^{\pr}_!(\phi^!_{\unit})^{-1}
  \colon
  \omega^{\pr!}\omega^{\pr}_!(s^{\pr})
  \to \phi^!\omega^!\omega_!(s)$.
  Put 
  $\varpi\coloneqq\omega^!\omega_!$,
  $\varpi^{\pr}\coloneqq\omega^{\pr!}\omega^{\pr}_!$,
  $\mu\coloneqq\omega^!\epsilon\omega_!\colon \varpi\varpi\Rightarrow\varpi$,
  and
  $\mu^{\pr}\coloneqq\omega^{\pr!}\epsilon^{\pr}\omega^{\pr}_!
  \colon \varpi^{\pr}\varpi^{\pr}\Rightarrow\varpi^{\pr}$.
  Observe that $\varpi^{\pr}_{\phi}=(\omega^{\pr!})_{\phi}(\omega^{\pr}_!)_{\phi}$
  as colax $\kos{S}$-tensor functors.
  By substituting 
  $\mon{f}^{\pr}=\varpi^{\pr}$,
  $\mon{g}=\varpi$,
  $\vartheta=\alpha^{\sharp}$
  in \eqref{eq-changebaseGal-colaxStnathatar},
  we obtain the following commutative diagram for every $z\in \CS$.
  \begin{equation}
    \label{eq-PKGalfund-alphasharpequiv}
    \vcenter{\hbox{
      \xymatrix@R=15pt@C=50pt{
        \omega^{\pr!}\omega^{\pr}_!\phi^!(z)
        \ar[rr]^-{(\hatar{\omega^{\pr!}\omega^{\pr}_!})_{\phi^!(z)}}_-{\cong}
        \ar[d]_-{\alpha^{\sharp}_z}
        &\text{ }
        &\omega^{\pr!}\omega^{\pr}_!(s^{\pr})\otimes^{\pr}\phi^!(z)
        \ar[d]^-{\breve{f}\otimes^{\pr} I_{\phi^!(z)}}
        \\%1
        \phi^!\omega^!\omega_!(z)
        \ar[r]^-{\phi^!((\hatar{\omega^!\omega_!})_z)}_-{\cong}
        &\phi^!(\omega^!\omega_!(s)\otimes z)
        \ar[r]^-{(\phi^!_{\tensor})_{\omega^!\omega_!(s),z}}_-{\cong}
        &\phi^!\omega^!\omega_!(s)\otimes^{\pr}\phi^!(z)
      }
    }}
  \end{equation}
  One can check that
  $(\phi^!\mu)
  \circ (\alpha^{\sharp}\varpi)
  \circ (\varpi^{\pr}\alpha^{\sharp})
  =\alpha^{\sharp}
  \circ (\mu^{\pr}\phi^!)
  \colon \varpi^{\pr}\varpi^{\pr}\phi^!
  \Rightarrow \phi^!\varpi$
  and that
  $\phi^!\eta
  =\alpha^{\sharp}
  \circ (\eta^{\pr}\phi^!)
  \colon \phi^!
  \Rightarrow \phi^!\varpi$
  as comonoidal natural transformations between colax monoidal functors,
  using \eqref{eq-PKGalfund-matealpha}.
  We also have
  $\pc^{\pr}_{\phi^!(\pi)}
  \circ (\breve{f}\otimes^{\pr}I_{\phi^!(\pi)})
  =\phi^!(\mu_{s})
  \circ \alpha^{\sharp}_{\pi}
  \circ (\hatar{\varpi^{\pr}})_{\phi^!(\pi)}^{-1}
  \colon \pi^{\pr}\otimes^{\pr}\phi^!(\pi)
  \to \phi^!(\pi)$,
  which follows from the relation
  \eqref{eq-PKGalfund-alphasharpequiv} for $z=\pi$.
  From these equalities, we deduce that
  $\breve{f}\colon \pi^{\pr}\to \phi^!(\pi)$
  is a morphism of Galois objects in $\kos{S}^{\pr}$.
  It remains to show the existence of $\breve{\alpha}$.
  It suffices to show that for each $X\in \CT$,
  the isomorphism
  $\alpha_X\colon \omega^{\pr!}\kos{f}^!(X)\cong \phi^!\omega^!(X)$
  in $\CS^{\pr}$ lifts to an isomorphism
  $\breve{\alpha}_X\colon \breve{\omega}^{\pr!}\kos{f}^!(X)\cong \breve{\kos{f}}^!\breve{\omega}^!(X)$
  in $\cat{Rep}_{\kos{S}^{\pr}}(\pi^{\pr})$.
  One can check this
  using \eqref{eq-PKGalfund-alphasharpequiv} for $z=\omega^!(X)$
  and \eqref{eq-PKGalfund-matealpha}.
\end{proof}

\section{The duality over Galois prekosmoi}
\label{sec-dualityPKGal}

In this section, we construct the $(2,1)$-category
$\mathsf{pGalCat}_{\mathsf{pt}}$
and prove Theorem~\ref{thm-PKGalCat-equiv21catoveret}.

\begin{definition}
  \label{def-dualityPKGal-pGalCatpt21cat}
  We define 
  the $(2,1)$-category
  $\mathsf{pGalCat}_{\mathsf{pt}}$
  as follows.
  \begin{itemize}
    \item 
    A $0$-cell in $\mathsf{pGalCat}_{\mathsf{pt}}$
    is a pointed prekosmic Galois category.

    \item
    A $1$-cell in $\mathsf{pGalCat}_{\mathsf{pt}}$
    is a \emph{pointed prekosmic fundamental morphism}.
    Let 
    $(\kos{S},\kos{t},\omega)$,
    $(\kos{S}^{\pr},\kos{t}^{\pr},\omega^{\pr})$
    be $0$-cells.
    A pointed prekosmic fundamental morphism
    $(\phi,\kos{f},\alpha)\colon (\kos{S}^{\pr},\kos{t}^{\pr},\omega^{\pr})\mor(\kos{S},\kos{t},\omega)$,
    denoted by an arrow with a filled tip,
    is a tuple consisting of
    an \'{e}tale Galois morphism $\phi\colon \kos{S}^{\pr}\mor\kos{S}$,
    a Galois $\kos{S}$-morphism
    $\kos{f}\colon \kos{t}^{\pr}_{\phi}\mor\kos{t}$,
    and an invertible Galois $\kos{S}$-transformation
    $\alpha\colon \kos{f}\omega^{\pr}_\phi\cong \omega\phi$.
    The pair $(\phi,\kos{f})$
    is a prekosmic fundamental morphism
    from $(\kos{S}^{\pr},\kos{t}^{\pr})$ to $(\kos{S},\kos{t})$.
    We illustrate $(\phi,\kos{f},\alpha)$ as follows.
    \begin{equation*}
      \begin{tikzcd}[row sep=5pt,column sep=20pt]
        \kos{T}^{\pr}
        \ar[dd,filled,"\kos{t}^{\pr}"']
        \ar[r,filled,"\kos{f}"]
        &\kos{T}
        \ar[dd,filled,"\kos{t}"]
        &\text{ }
        &\phi
        \ar[dd,filled,"\omega^{\pr}_{\phi}"']
        \ar[r,filled,"\phi"]
        &\id_{\kos{S}}
        \ar[dd,filled,"\omega"]
        \\%1
        \text{ }
        &\text{ }
        \ar[l,2cell,shift left=1ex,"{\what{\kos{f}}}"']
        &\text{ }
        &\text{ }
        \ar[r,2cell,shift right=0.7ex,pos=0.5,"{\!\alpha}"]
        &\text{ }
        \\%2
        \kos{S}^{\pr}
        \ar[r,filled,"{\phi}"]
        &\kos{S}
        &\text{ }
        &\kos{t}^{\pr}_{\phi}
        \ar[r,filled,"{\kos{f}}"]
        &\kos{t}
      \end{tikzcd}
    \end{equation*}
    
    \item
    A $2$-cell
    $(\chi,\vartheta)\colon (\phi,\kos{f},\alpha)\cong (\psi,\kos{g},\beta)$
    in $\mathsf{pGalCat}_{\mathsf{pt}}$
    between $1$-cells
    $(\phi,\kos{f},\alpha)$, $(\psi,\kos{g},\beta)
    \colon (\kos{S}^{\pr},\kos{t}^{\pr},\omega^{\pr})\mor (\kos{S},\kos{t},\omega)$
    is a pair consisting of an invertible Galois transformation
    $\chi\colon \phi\cong\psi$
    and an invertible Galois $\kos{S}$-transformation
    $\vartheta\colon \kos{f}\cong \kos{g}\kos{t}^{\pr}_{\chi}$
    between Galois $\kos{S}$-morphisms
    $\kos{f}$, $\kos{g}\kos{t}^{\pr}_{\chi}\colon \kos{t}^{\pr}_{\phi}\mor\kos{t}$.
  \end{itemize}
  We define the rest of the structure of $\mathsf{pGalCat}_{\mathsf{pt}}$ as follows.
  Let 
  $(\chi,\vartheta)\colon (\phi,\kos{f},\alpha)\cong (\psi,\kos{g},\beta)$
  be a $2$-cell
  between $1$-cells
  $(\phi,\kos{f},\alpha)$, $(\psi,\kos{g},\beta)
  \colon (\kos{S}^{\pr},\kos{t}^{\pr},\omega^{\pr})\mor (\kos{S},\kos{t},\omega)$.

  Let
  $(\tilde{\chi},\tilde{\vartheta})\colon (\psi,\kos{g},\beta)\cong (\tilde{\psi},\tilde{\kos{g}},\tilde{\beta})$
  be a $2$-cell.
  The vertical composite of
  $(\chi,\vartheta)$ and $(\tilde{\chi},\tilde{\vartheta})$
  is the pair
  $(\tilde{\chi}\circ\chi,(\tilde{\vartheta}\kos{t}^{\pr}_{\chi})
  \circ \vartheta)
  \colon (\phi,\kos{f},\alpha)\cong (\tilde{\psi},\tilde{\kos{g}},\tilde{\beta})$
  consisting of 
  $\tilde{\chi}\circ\chi\colon \phi\cong \tilde{\psi}$
  and
  $(\tilde{\vartheta}\kos{t}^{\pr}_{\chi})
  \circ \vartheta
  \colon \kos{f}\cong \tilde{\kos{g}}\kos{t}^{\pr}_{\tilde{\chi}}\kos{t}^{\pr}_{\chi}$.
  The identity $2$-cell of $(\phi,\kos{f},\alpha)$
  is the pair consisting of the identity Galois transformation of $\phi$
  and the identity Galois $\kos{S}$-transformation of $\kos{f}$.
  Each $2$-cell
  $(\chi,\vartheta)\colon (\phi,\kos{f},\alpha)\cong (\psi,\kos{g},\beta)$
  is indeed invertible with respect to the vertical composition;
  its inverse $(\chi^{-1},\vartheta^{-1}\kos{t}^{\pr}_{\chi^{-1}})\colon(\psi,\kos{g},\beta)\cong (\phi,\kos{f},\alpha)$
  is the pair consisting of 
  $\chi^{-1}\colon \psi\cong \phi$
  and 
  $\vartheta^{-1}\kos{t}^{\pr}_{\chi^{-1}}\colon \kos{g}\cong \kos{f}\kos{t}^{\pr}_{\chi^{-1}}$.

  Let 
  $(\phi^{\pr},\kos{f}^{\pr},\alpha^{\pr})$,
  $(\psi^{\pr},\kos{g}^{\pr},\beta^{\pr})
  \colon (\kos{S}^{\ppr},\kos{t}^{\ppr},\omega^{\ppr})\mor(\kos{S}^{\pr},\kos{t}^{\pr},\omega^{\pr})$
  be $1$-cells.
  The composite of
  $(\phi,\kos{f},\alpha)$
  and
  $(\phi^{\pr},\kos{f}^{\pr},\alpha^{\pr})$
  is the tuple 
  $(\phi^{\ppr},\kos{f}^{\ppr},\alpha^{\ppr})\colon (\kos{S}^{\ppr},\kos{t}^{\ppr},\omega^{\ppr})\mor (\kos{S},\kos{t},\omega)$
  consisting of 
  $\phi^{\ppr}\coloneqq \phi\phi^{\pr}\colon \kos{S}^{\ppr}\mor\kos{S}$,
  $\kos{f}^{\ppr}\coloneqq \kos{f}\kos{f}^{\pr}_{\phi}\colon \kos{t}^{\ppr}_{\phi^{\ppr}}\mor\kos{t}$,
  and 
  $\alpha^{\ppr}\coloneqq
  (\alpha\phi^{\pr}_{\phi})
  \circ (\kos{f}\alpha^{\pr}_{\phi})
  \colon \kos{f}^{\ppr}\omega^{\ppr}_{\phi^{\ppr}}\cong \omega\phi^{\ppr}$.
  We illustrate the composite 
  and its description as follows.
  \begin{equation*}
    \begin{tikzcd}[row sep=5pt,column sep=20pt]
      \phi^{\ppr}
      \ar[dd,filled,"\omega^{\ppr}_{\phi^{\ppr}}"']
      \ar[r,filled,"\phi^{\ppr}"]
      &\id_{\kos{S}}
      \ar[dd,filled,"\omega"]
      &\text{ }
      &\phi^{\pr}_{\phi}
      \ar[dd,filled,"(\omega^{\ppr}_{\phi^{\pr}})_{\phi}"']
      \ar[r,filled,"\phi^{\pr}_{\phi}"]
      &\phi
      \ar[dd,filled,"\omega^{\pr}_{\phi}"description]
      \ar[r,filled,"\phi"]
      &\id_{\kos{S}}
      \ar[dd,filled,"\omega"]
      \\%1
      \text{ }
      \ar[r,2cell,shift right=0ex,pos=0.5,"{\!\!\alpha^{\ppr}}"]
      &\text{ }
      &=
      &\text{ }
      \ar[r,2cell,shift right=0ex,pos=0.5,"{\!\!\alpha^{\pr}_{\phi}}"]
      &\text{ }
      \ar[r,2cell,shift right=0ex,pos=0.5,"{\!\!\alpha}"]
      &\text{ }
      \\%2
      \kos{t}^{\ppr}_{\phi^{\ppr}}
      \ar[r,filled,"{\kos{f}^{\ppr}}"]
      &\kos{t}
      &\text{ }
      &(\kos{t}^{\ppr}_{\phi^{\pr}})_{\phi}
      \ar[r,filled,"{\kos{f}^{\pr}_{\phi}}"]
      &\kos{t}^{\pr}_{\phi}
      \ar[r,filled,"{\kos{f}}"]
      &\kos{t}
    \end{tikzcd}
  \end{equation*}
  The identity $1$-cell of $(\kos{S},\kos{t},\omega)$ is 
  the tuple consisting of the identity Galois morphism $\vphantom{\id}\smash[t]{\id}_{\kos{S}\vphantom{!}}^{\vphantom{!}}$ of $\kos{S}$,
  the identity Galois $\kos{S}$-morphism of $\kos{t}$,
  and the identity Galois $\kos{S}$-transformation of $\omega\colon \vphantom{\id}\smash[t]{\id}_{\kos{S}\vphantom{!}}^{\vphantom{!}}\mor\kos{t}$.
  Let
  $(\chi^{\pr},\vartheta^{\pr})\colon (\phi^{\pr},\kos{f}^{\pr},\alpha^{\pr})\cong (\psi^{\pr},\kos{g}^{\pr},\beta^{\pr})$
  be a $2$-cell.
  We write $(\psi^{\ppr},\kos{g}^{\ppr},\beta^{\ppr})\colon (\kos{S}^{\ppr},\kos{t}^{\ppr},\omega^{\ppr})\mor (\kos{S},\kos{t},\omega)$
  for the composite of
  $(\psi,\kos{g},\beta)$
  and
  $(\psi^{\pr},\kos{g}^{\pr},\beta^{\pr})$.
  The horizontal composite
  $(\chi^{\ppr},\vartheta^{\ppr})
  \colon (\phi^{\ppr},\kos{f}^{\ppr},\alpha^{\ppr})\cong (\psi^{\ppr},\kos{g}^{\ppr},\beta^{\ppr})$
  of $(\chi,\vartheta)$ and $(\chi^{\pr},\vartheta^{\pr})$
  is the pair
  consisting of
  $\chi^{\ppr}\colon \phi^{\ppr}\cong \psi^{\ppr}$,
  given by
  $\chi\chi^{\pr}\colon \phi\phi^{\pr}\cong \psi\psi^{\pr}$,
  and 
  $\vartheta^{\ppr}\colon \kos{f}^{\ppr}\cong \kos{g}^{\ppr}\kos{t}^{\ppr}_{\chi^{\ppr}}$,
  given by
  $\vartheta\vartheta^{\pr}_{\phi}
  \colon \kos{f}\kos{f}^{\pr}_{\phi}
  \cong \kos{g}\kos{t}^{\pr}_{\chi}\kos{g}^{\pr}_{\phi}(\kos{t}^{\ppr}_{\chi^{\pr}})_{\phi}
  =\kos{g}\kos{g}^{\pr}_{\psi}(\kos{t}^{\ppr}_{\psi^{\pr}})_{\chi}(\kos{t}^{\ppr}_{\chi^{\pr}})_{\phi}
  =\kos{g}\kos{g}^{\pr}_{\psi}\kos{t}^{\ppr}_{\chi\chi^{\pr}}$.
\end{definition}

We denote by $\mathsf{p}_{\mathsf{ptcat}}
\colon \mathsf{pGalCat}_{\mathsf{pt}}\to\mathsf{pGalKos}_{\etsf}$
the $2$-functor that sends each cell to its first component.
The strict fiber of $\mathsf{p}_{\mathsf{ptcat}}$
over each Galois prekosmos $\kos{S}$
is $2$-isomorphic to $\mathsf{pGalCat}_{\mathsf{pt}}(\kos{S})$,
defined in Definition~\ref{def-dualityoverGalpk-pGalCatptS}.
We illustrate $\mathsf{p}_{\mathsf{ptcat}}$ as follows.
\begin{equation*}
  \begin{tikzcd}[row sep=22pt,column sep=20pt]
    \mathsf{pGalCat}_{\mathsf{pt}}
    \ar[d,"\mathsf{p}_{\mathsf{ptcat}}"']
    &(\kos{S}^{\pr},\kos{t}^{\pr},\omega^{\pr})
    \ar[rr,filled,shift left=0.3ex,bend left=1.4pc,"{(\phi,\kos{f},\alpha)}" description]
    \ar[rr,filled,shift right=0.3ex,bend right=1.4pc,"{(\psi,\kos{g},\beta)}" description]
    &\text{ }
    \ar[d,2cell,pos=-0.2,shift right=0.2ex,"{\,{(\chi,\vartheta)}}"]
    \ar[d,2cell,pos=1,shift right=0.2ex,"{\,\chi}"]
    &(\kos{S},\kos{t},\omega)
    \\%1
    {\mathsf{pGalKos}_{\etsf}}
    &\kos{S}^{\pr}
    \ar[rr,filled,shift left=0.3ex,bend left=1.4pc,"\phi" description]
    \ar[rr,filled,shift right=0.3ex,bend right=1.4pc,"\psi" description]
    &\text{ }
    &\kos{S}
  \end{tikzcd}
\end{equation*}

\begin{proposition}
  \label{prop-dualityPKGal-Rep}
  We have a pseudofunctor
  \begin{equation*}
    \mathsf{Rep}\colon \mathsf{pGalObj}\to \mathsf{pGalCat}_{\mathsf{pt}}
  \end{equation*}
  over $\mathsf{pGalKos}_{\etsf}$
  that sends each prekosmic Galois object $(\kos{S},\pi)$
  to the pointed prekosmic Galois category
  $(\kos{S},\kos{t}_{\pi\vphantom{!}}^{\vphantom{!}},\kos{u}_{\pi\vphantom{!}}^{\vphantom{!}})$.
\end{proposition}
\begin{proof}
  We describe how $\mathsf{Rep}$ sends
  $1$-cells and $2$-cells.
  \begin{itemize}
    \item
    A $1$-cell $(\phi,f)\colon (\kos{S}^{\pr},\pi^{\pr})\to(\kos{S},\pi)$
    is sent to the $1$-cell
    $(\kos{S}^{\pr},\kos{t}^{\pr}_{\pi^{\pr}},\kos{u}^{\pr}_{\pi^{\pr}})
    \mor(\kos{S},\kos{t}_{\pi\vphantom{!}}^{\vphantom{!}},\kos{u}_{\pi\vphantom{!}}^{\vphantom{!}})$
    consisting of
    $\phi\colon \kos{S}^{\pr}\mor\kos{S}$,
    the Galois $\kos{S}$-morphism
    $\kos{f}\colon (\kos{t}^{\pr}_{\pi^{\pr}})_{\phi}\mor \kos{t}_{\pi\vphantom{!}}^{\vphantom{!}}$,
    and the invertible Galois $\kos{S}$-transformation
    $\kos{f}(\kos{u}^{\pr}_{\pi^{\pr}})_{\phi}\cong \kos{u}_{\pi\vphantom{!}}^{\vphantom{!}}\phi$,
    described in the proof of Lemma~\ref{lem-PKGalfund-inducedPKfundmor}.

    \item
    Let
    $(\phi,f)$, $(\psi,g)\colon (\kos{S}^{\pr},\pi^{\pr})\to(\kos{S},\pi)$
    be $1$-cells in $\mathsf{pGalObj}$
    and denote by
    $\kos{f}\colon (\kos{t}^{\pr}_{\pi^{\pr}})_{\phi}\mor \kos{t}_{\pi\vphantom{!}}^{\vphantom{!}}$,
    $\kos{g}\colon (\kos{t}^{\pr}_{\pi^{\pr}})_{\psi}\mor \kos{t}_{\pi\vphantom{!}}^{\vphantom{!}}$
    the associated Galois $\kos{S}$-morphisms.
    A $2$-cell $(\chi,\theta)\colon (\phi,f)\cong (\psi,g)$
    in $\mathsf{pGalObj}$
    is sent to the $2$-cell
    in $\mathsf{pGalCat}_{\mathsf{pt}}$
    consisting of $\chi\colon \phi\cong \psi$
    and the invertible Galois $\kos{S}$-transformation
    $\vartheta\colon \kos{f}\cong \kos{g}(\kos{t}^{\pr}_{\pi^{\pr}})_{\chi}$,
    defined as follows.
    The component $\vartheta_X\colon \kos{f}^!(X)\cong \kos{g}^!(X)$
    of $\vartheta$ at $X=(x,\gamma_x)\in\cat{Rep}_{\kos{S}}(\pi)$
    is the isomorphism in $\cat{Rep}_{\kos{S}^{\pr}}(\pi^{\pr})$
    whose underlying isomorphism in $\CS^{\pr}$ is
    $\chi_x
    \circ \phi^!(\gamma_x)
    \circ (\phi^!_{\tensor})_{\pi,x}^{-1}
    \circ (\theta\otimes^{\pr} I_{\phi^!(x)})
    \colon \phi^!(x)
    \cong \psi^!(x)$.
    It is straightforward to check that
    $\vartheta$ is indeed an invertible Galois $\kos{S}$-transformation
    from $\kos{f}$ to $\kos{g}(\kos{t}^{\pr}_{\pi^{\pr}})_{\chi}$.
  \end{itemize}
  It is straightforward to check that $\mathsf{Rep}$
  preserves vertical composition of $2$-cells and identity $2$-cells;
  thus $\mathsf{Rep}$ defines a functor on each hom-groupoid.
  The constraints of the pseudofunctor $\mathsf{Rep}$
  are defined as follows.
  \begin{itemize}
    \item
    Let 
    $(\phi,f)$,
    $(\phi^{\pr},f^{\pr})$
    be a composable pair of $1$-cells in $\mathsf{pGalObj}$,
    with composite
    $(\phi^{\ppr},f^{\ppr})\coloneqq (\phi\phi^{\pr},\phi^{\pr!}(f)f^{\pr})$,
    where $\phi\colon \kos{S}^{\pr}\mor\kos{S}$
    and
    $\phi^{\pr}\colon\kos{S}^{\ppr}\mor\kos{S}^{\pr}$.
    We denote by
    $\kos{f}$, $\kos{f}^{\pr}$, and $\kos{f}^{\ppr}$
    the associated Galois morphisms
    over $\kos{S}$, $\kos{S}^{\pr}$, and $\kos{S}$, respectively.
    The component of the functorial constraint of $\mathsf{Rep}$
    at $(\phi,f)$, $(\phi^{\pr},f^{\pr})$
    is the invertible Galois $\kos{S}$-transformation
    $\kos{f}\kos{f}^{\pr}_{\phi}\cong\kos{f}^{\ppr}$
    given by the equality
    $(\kos{f}^{\pr!})_{\phi}\kos{f}^!=\kos{f}^{\ppr!}$
    between the right adjoint parts.

    \item
    Let $(\kos{S},\pi)$ be a $0$-cell in $\mathsf{pGalObj}$.
    We write 
    $\kos{I}_{\pi\vphantom{!}}^{\vphantom{!}}\colon \kos{t}_{\pi\vphantom{!}}^{\vphantom{!}}\mor\kos{t}_{\pi\vphantom{!}}^{\vphantom{!}}$
    for the Galois $\kos{S}$-morphism
    associated to the identity $1$-cell
    of $(\kos{S},\pi)$ in Lemma~\ref{lem-PKGalfund-inducedPKfundmor},
    which agrees with the $\kos{I}_{\pi\vphantom{!}}^{\vphantom{!}}$
    introduced in the proof of Proposition~\ref{prop-dualityoverGalpk-RepSpseudofunct}.
    The component
    of the unity constraint of $\mathsf{Rep}$
    at $(\kos{S},\pi)$
    is the $2$-cell
    consisting of the identity Galois transformation of
    $\vphantom{\id}\smash[t]{\id}_{\kos{S}\vphantom{!}}^{\vphantom{!}}$
    and the invertible Galois $\kos{S}$-transformation
    from the identity Galois $\kos{S}$-morphism of $\kos{t}_{\pi\vphantom{!}}^{\vphantom{!}}$
    to $\kos{I}_{\pi\vphantom{!}}^{\vphantom{!}}$,
    given by the equality between the right adjoint parts.
  \end{itemize}
  One can check that $\mathsf{Rep}$ preserves
  horizontal composition of $2$-cells
  up to the functorial constraint,
  or equivalently,
  that the functorial constraint
  $\kos{f}\kos{f}^{\pr}_{\phi}\cong \kos{f}^{\ppr}$
  is natural in $(\phi,f)$, $(\phi^{\pr},f^{\pr})$.
  It is immediate that 
  the functorial and unity constraints of $\mathsf{Rep}$
  satisfy the associativity and left/right unity axioms,
  as the axioms are just equalities between the right adjoint parts.
  This completes the construction of the pseudofunctor
  $\mathsf{Rep}\colon \mathsf{pGalObj}\to \mathsf{pGalCat}_{\mathsf{pt}}$
  over $\mathsf{pGalKos}_{\etsf}$.
  The strict fiber of $\mathsf{Rep}$ over
  each Galois prekosmos $\kos{S}$
  is $2$-isomorphic to the pseudofunctor $\mathsf{Rep}_{\kos{S}}$
  in Proposition~\ref{prop-dualityoverGalpk-RepSpseudofunct}.
\end{proof}

\begin{proposition}
  \label{prop-dualityPKGal-Aut}
  We have a $2$-functor
  \begin{equation*}
    \mathsf{Aut}\colon \mathsf{pGalCat}_{\mathsf{pt}}\to \mathsf{pGalObj}
  \end{equation*}
  over $\mathsf{pGalKos}_{\etsf}$
  that sends each pointed prekosmic Galois category
  $(\kos{S},\kos{t},\omega)$
  to the prekosmic Galois object $(\kos{S},\omega^!\omega_!(s))$,
  where $\omega^!\omega_!(s)$
  is the representing object of $\uAut_{\kos{S}}(\omega)$
  in Theorem~\ref{thm-Galcatoverpk-mainthm}.
\end{proposition}
\begin{proof}
  We describe how $\mathsf{Aut}$ sends $1$-cells and $2$-cells.
  \begin{itemize}
    \item
    A $1$-cell $(\phi,\kos{f},\alpha)\colon (\kos{S}^{\pr},\kos{t}^{\pr},\omega^{\pr})\mor(\kos{S},\kos{t},\omega)$
    in $\mathsf{pGalCat}_{\mathsf{pt}}$
    is sent to
    the $1$-cell
    $(\phi,\breve{f})\colon (\kos{S}^{\pr},\omega^{\pr!}\omega^{\pr}_!(s^{\pr}))
    \to (\kos{S},\omega^!\omega_!(s))$
    in $\mathsf{pGalObj}$
    described in the proof of Proposition~\ref{prop-PKGalfund-inducedPKfundmor}.

    \item
    A $2$-cell $(\chi,\vartheta)\colon (\phi,\kos{f},\alpha)\cong (\psi,\kos{g},\beta)$
    in $\mathsf{pGalCat}_{\mathsf{pt}}$
    between $1$-cells
    $(\phi,\kos{f},\alpha)$, $(\psi,\kos{g},\beta)
    \colon (\kos{S}^{\pr},\kos{t}^{\pr},\omega^{\pr})\mor (\kos{S},\kos{t},\omega)$
    is sent to the $2$-cell
    $(\chi,\breve{\theta})\colon (\phi,\breve{f})\cong (\psi,\breve{g})$
    in $\mathsf{pGalObj}$,
    defined as follows.
    Put
    $\grave{\vartheta}^{\sharp}
    \coloneqq (\id_{\chi}^{\pr!}\beta\omega_!)
    \circ ((\omega^{\pr!})_{\phi}\vartheta\omega_!)
    \circ (\alpha^{-1}\omega_!)
    \circ (\phi^!\eta)
    \colon \phi^!
    \Rightarrow \id_{\chi}^{\pr!}\psi^!\omega^!\omega_!$,
    $\vartheta^{\sharp}
    \coloneqq
    (\chi^{-1}\omega^!\omega_!)\circ \grave{\vartheta}^{\sharp}\circ \chi^{-1}
    \colon \id_{\chi}^{\pr!}\psi^!
    \Rightarrow \phi^!\omega^!\omega_!$,
    and $\pi\coloneqq \omega^!\omega_!(s)$.
    We define
    $\breve{\theta}\colon s^{\pr}\to \phi^!(\pi)$
    as
    $\vartheta^{\sharp}_{s}
    \circ (\psi^!_{\unit})^{-1}$.
    By \eqref{eq-ActionbaseGal-colaxStnatequiv}
    for the colax $\kos{S}$-tensor natural transformation $\vartheta^{\sharp}$
    with $\monbase{t}=\phi^!$, $\monbase{t}^{\pr}=\vphantom{\id}\smash[t]{\id}_{\kos{S}\vphantom{!}}^{\vphantom{!}}$, $X^{\pr}=s$,
    we deduce that the following diagram commutes
    for every $z\in \CS$.
    \begin{equation}
      \label{eq-dualityPKGal-thetasharpequiv}
      \vcenter{\hbox{
        \xymatrix@R=12pt@C=45pt{
          \psi^!(z)
          \ar[rr]^-{\chi^{-1}_z}_-{\cong}
          \ar[d]_-{\vartheta^{\sharp}_z}
          &\text{ }
          &\phi^!(z)
          \ar[d]^-{\breve{\theta}\otimes^{\pr}I_{\phi^!(z)}}
          \\%1
          \phi^!\omega^!\omega_!(z)
          \ar[r]^-{\phi^!((\hatar{\omega^!\omega_!})_z)}_-{\cong}
          &\phi^!(\omega^!\omega_!(s)\otimes z)
          \ar[r]^-{(\phi^!_{\tensor})_{\omega^!\omega_!(s),z}}_-{\cong}
          &\phi^!\omega^!\omega_!(s)\otimes^{\pr}\phi^!(z)
        }
      }}
    \end{equation}
    Put 
    $\varpi\coloneqq \omega^!\omega_!$,
    $\varpi^{\pr}\coloneqq \omega^{\pr!}\omega^{\pr}_!$,
    and $\mu\coloneqq \omega^!\epsilon\omega_!\colon \varpi\varpi\Rightarrow\varpi$.
    One can check that 
    $\grave{\vartheta}^{\sharp}
    =(\beta\omega_!)
    \circ (\omega^{\pr!}\vartheta\omega_!)
    \circ (\omega^{\pr!}\alpha_!)
    \circ (\eta^{\pr}\phi^!)
    \colon \phi^!\Rightarrow\psi^!\varpi$
    and that
    $(\psi^!\mu)
    \circ (\grave{\vartheta}^{\sharp}\varpi)
    \circ \alpha^{\sharp}
    =(\beta\omega_!)
    \circ (\omega^{\pr!}\vartheta\omega_!)
    \circ (\omega^{\pr!}\alpha_!)
    = (\psi^!\mu)
    \circ (\beta^{\sharp}\varpi)
    \circ (\varpi^{\pr}\grave{\vartheta}^{\sharp})
    \colon\varpi^{\pr}\phi^!\Rightarrow\psi^!\varpi$
    as comonoidal natural transformations between colax monoidal functors,
    using \eqref{eq-PKGalfund-matealpha}.
    From these equalities,
    along with the relations
    \eqref{eq-PKGalfund-alphasharpequiv} for $z=\pi$
    and
    \eqref{eq-dualityPKGal-thetasharpequiv} for $z=\pi$,
    we deduce that $\breve{\theta}$ is indeed a $2$-isomorphism
    from $\breve{f}$ to $\chi_{\pi}^{-1}\breve{g}$.
  \end{itemize}
  It is straightforward to check that
  $\mathsf{Aut}$ preserves identity $1$-cells and identity $2$-cells.
  \begin{itemize}
    \item 
    Let $(\chi,\vartheta)\colon (\phi,\kos{f},\alpha)\cong (\psi,\kos{g},\beta)$
    and $(\tilde{\chi},\tilde{\vartheta})\colon (\psi,\kos{g},\beta)\cong (\tilde{\psi},\tilde{\kos{g}},\tilde{\beta})$
    be $2$-cells in $\mathsf{pGalCat}_{\mathsf{pt}}$,
    with vertical composite
    $(\bar{\chi},\bar{\vartheta})
    \colon (\phi,\kos{f},\alpha)\cong (\tilde{\psi},\tilde{\kos{g}},\tilde{\beta})$.
    Put $\pi\coloneqq \omega^!\omega_!(s)$.
    One can check that
    $\grave{\bar{\vartheta}}^{\sharp}
    =(\tilde{\psi}^!\omega^!\epsilon\omega_!)
    \circ (\grave{\tilde{\vartheta}}^{\sharp}\omega^!\omega_!)
    \circ \grave{\vartheta}^{\sharp}
    \colon \phi^!\Rightarrow\tilde{\psi}^!\omega^!\omega_!$
    as comonoidal natural transformations
    between colax monoidal functors,
    using \eqref{eq-PKGalfund-matealpha}.
    From this,
    along with \eqref{eq-dualityPKGal-thetasharpequiv}
    for $z=\pi$,
    we deduce that 
    $(\chi_{\pi}^{-1}\circ \breve{\tilde{\theta}})\star \breve{\theta}
    =\breve{\bar{\theta}}
    \colon s^{\pr}\to \phi^!(\pi)$
    as morphisms in $\cat{Ens}(\kos{S}^{\pr})$.
    Hence $(\bar{\chi},\breve{\bar{\theta}})$
    is the vertical composite of $(\chi,\breve{\theta})$ and $(\tilde{\chi},\breve{\tilde{\theta}})$;
    this shows that $\mathsf{Aut}$ preserves vertical composition of $2$-cells.
    
    \item 
    Let 
    $(\phi^{\ppr},\kos{f}^{\ppr},\alpha^{\ppr})$
    be the composite of
    $(\phi,\kos{f},\alpha)\colon (\kos{S}^{\pr},\kos{t}^{\pr},\omega^{\pr})\mor(\kos{S},\kos{t},\omega)$
    and
    $(\phi^{\pr},\kos{f}^{\pr},\alpha^{\pr})\colon (\kos{S}^{\ppr},\kos{t}^{\ppr},\omega^{\ppr})\mor(\kos{S}^{\pr},\kos{t}^{\pr},\omega^{\pr})$
    in $\mathsf{pGalCat}_{\mathsf{pt}}$,
    with
    $\phi^{\ppr}\coloneqq \phi\phi^{\pr}$,
    $\kos{f}^{\ppr}\coloneqq \kos{f}\kos{f}^{\pr}_{\phi}$,
    and
    $\alpha^{\ppr}\coloneqq (\alpha\phi^{\pr}_{\phi})\circ (\kos{f}\alpha^{\pr}_{\phi})$.
    One can check that
    $\alpha^{\ppr}_!=(\kos{f}^{\pr!}\alpha_!)\circ (\alpha^{\pr}_!\phi^!)
    \colon \omega^{\ppr}_!\phi^{\ppr!}\Rightarrow\kos{f}^{\ppr!}\omega_!$
    and that
    $\alpha^{\ppr\sharp}
    =(\phi^{\pr!}\alpha^{\sharp})\circ (\alpha^{\pr\sharp}\phi^!)
    \colon \omega^{\ppr!}\omega^{\ppr}_!\phi^{\ppr!}
    \Rightarrow \phi^{\ppr!}\omega^!\omega_!$
    as comonoidal natural transformations between colax monoidal functors,
    using \eqref{eq-PKGalfund-matealpha}.
    Thus we have
    $\breve{f}^{\ppr}=\phi^{\pr!}(\breve{f})\circ \breve{f}^{\pr}
    \colon \omega^{\ppr!}\omega^{\ppr}_!(s^{\ppr})
    \to \phi^{\ppr!}\omega^!\omega_!(s)$.
    This shows that $\mathsf{Aut}$ preserves composition of $1$-cells.
    
    \item
    Let
    $(\chi,\vartheta)\colon (\phi,\kos{f},\alpha)\cong (\psi,\kos{g},\beta)$
    be a $2$-cell
    between $1$-cells from $(\kos{S}^{\pr},\kos{t}^{\pr},\omega^{\pr})$ to $(\kos{S},\kos{t},\omega)$,
    and $(\chi^{\pr},\vartheta^{\pr})\colon (\phi^{\pr},\kos{f}^{\pr},\alpha^{\pr})\cong (\psi^{\pr},\kos{g}^{\pr},\beta^{\pr})$
    be a $2$-cell
    between $1$-cells from
    $(\kos{S}^{\ppr},\kos{t}^{\ppr},\omega^{\ppr})$
    to $(\kos{S}^{\pr},\kos{t}^{\pr},\omega^{\pr})$.
    Write
    $(\phi^{\ppr},\kos{f}^{\ppr},\alpha^{\ppr})$,
    $(\psi^{\ppr},\kos{g}^{\ppr},\beta^{\ppr})$
    for the composite $1$-cells,
    and
    $(\chi^{\ppr},\vartheta^{\ppr})\colon (\phi^{\ppr},\kos{f}^{\ppr},\alpha^{\ppr})\cong(\psi^{\ppr},\kos{g}^{\ppr},\beta^{\ppr})$
    for the horizontal composite of $(\chi,\vartheta)$ and $(\chi^{\pr},\vartheta^{\pr})$.
    Put $\pi\coloneqq \omega^!\omega_!(s)$.
    One can check that
    $\grave{\vartheta}^{\ppr\sharp}
    =(\psi^{\ppr!}\omega^!\epsilon\omega_!)
    \circ (\psi^{\pr!}\beta^{\sharp}\omega^!\omega_!)
    \circ (\grave{\vartheta}^{\pr\sharp}\psi^!\omega^!\omega_!)
    \circ (\phi^{\pr!}\grave{\vartheta}^{\sharp})
    \colon \phi^{\ppr!}\Rightarrow\psi^{\ppr!}\omega^!\omega_!$
    as comonoidal natural transformations between colax monoidal functors,
    using \eqref{eq-PKGalfund-matealpha}.
    From this,
    along with \eqref{eq-dualityPKGal-thetasharpequiv} for $z=\pi$,
    we deduce that
    $\breve{\theta}^{\ppr}
    =(\phi^{\pr!}(\chi_{\pi}^{-1}\circ \breve{g})\circ \breve{\theta}^{\pr})
    \star \phi^{\pr!}(\breve{\theta})
    \colon s^{\ppr}\to \phi^{\ppr!}(\pi)$
    as morphisms in $\cat{Ens}(\kos{S}^{\ppr})$.
    Hence $(\chi^{\ppr},\breve{\theta}^{\ppr})$
    is the horizontal composite of
    $(\chi,\breve{\theta})$ and $(\chi^{\pr},\breve{\theta}^{\pr})$;
    this shows that $\mathsf{Aut}$ preserves horizontal composition of $2$-cells.
  \end{itemize}
  This completes the construction of the $2$-functor
  $\mathsf{Aut}\colon \mathsf{pGalCat}_{\mathsf{pt}}\to \mathsf{pGalObj}$
  over $\mathsf{pGalKos}_{\etsf}$.
  The strict fiber of $\mathsf{Aut}$ over each Galois prekosmos $\kos{S}$
  is $2$-isomorphic to the $2$-functor
  $\mathsf{Aut}_{\kos{S}}$ in Proposition~\ref{prop-dualityoverGalpk-AutS2funct}.
\end{proof}

\begin{proof}[Proof of Theorem~\ref{thm-PKGalCat-equiv21catoveret}]
  We show that the pseudofunctor
  $\mathsf{Rep}$
  in Proposition~\ref{prop-dualityPKGal-Rep}
  and the $2$-functor
  $\mathsf{Aut}$
  in Proposition~\ref{prop-dualityPKGal-Aut}
  are mutually quasi-inverse equivalences of $(2,1)$-categories
  over $\mathsf{pGalKos}_{\etsf}$;
  that is, the composites
  $\mathsf{Aut}\circ \mathsf{Rep}$
  and $\mathsf{Rep}\circ \mathsf{Aut}$
  are equivalent, over $\mathsf{pGalKos}_{\etsf}$,
  to the identity $2$-functors
  of $\mathsf{pGalObj}$
  and $\mathsf{pGalCat}_{\mathsf{pt}}$,
  respectively.
  
  The composite
  $\mathsf{Aut}\circ \mathsf{Rep}\colon 
  \mathsf{pGalObj}\to \mathsf{pGalObj}$
  sends each $0$-cell $(\kos{S},\pi)$
  to the $0$-cell $(\kos{S},\pi\otimes s)$,
  and $\jmath_{\pi}^{-1}\colon \pi\otimes s\cong \pi$
  is an isomorphism of Galois objects in $\kos{S}$;
  see Remark~\ref{rem-GalMainthm-reconstruction}.
  It is straightforward to check that
  $\mathsf{Aut}\circ \mathsf{Rep}$
  is a $2$-functor over $\mathsf{pGalKos}_{\etsf}$,
  and that the $1$-cells
  $(\vphantom{\id}\smash[t]{\id}_{\kos{S}\vphantom{!}}^{\vphantom{!}},\jmath_{\pi}^{-1})$
  constitute a $2$-natural isomorphism over $\mathsf{pGalKos}_{\etsf}$
  from $\mathsf{Aut}\circ \mathsf{Rep}$
  to the identity $2$-functor of $\mathsf{pGalObj}$.

  We now construct an equivalence
  over $\mathsf{pGalKos}_{\etsf}$
  from the composite $\mathsf{Rep}\circ \mathsf{Aut}$
  to the identity $2$-functor of $\mathsf{pGalCat}_{\mathsf{pt}}$
  as follows.
  Such an equivalence is an invertible strong transformation
  over $\mathsf{pGalKos}_{\etsf}$,
  with the same convention as in the proof of
  Theorem~\ref{thm-dualityoverGalpk-equiv21cat};
  it is an equivalence internal to the $2$-category of
  pseudofunctors from $\mathsf{pGalCat}_{\mathsf{pt}}$ to itself,
  strong transformations, and modifications,
  all over $\mathsf{pGalKos}_{\etsf}$.
  We take the component $1$-cell of this transformation
  at each $0$-cell
  $(\kos{S},\kos{t},\omega)$ 
  to be the equivalence
  $(\kos{S},\kos{t}_{\omega^!\omega_!(s)},\kos{u}_{\omega^!\omega_!(s)})
  \simeq
  (\kos{S},\kos{t},\omega)$
  in the strict fiber over $\kos{S}$,
  consisting of $\vphantom{\id}\smash[t]{\id}_{\kos{S}\vphantom{!}}^{\vphantom{!}}$,
  the equivalence
  $\breve{\omega}\colon \kos{t}_{\omega^!\omega_!(s)}\simeq\kos{t}$,
  and the invertible Galois $\kos{S}$-transformation
  $\breve{\omega}\kos{u}_{\omega^!\omega_!(s)}\cong \omega$;
  see Theorem~\ref{thm-Galcatoverpk-mainthm}.
  We take its component $2$-cell at each $1$-cell 
  $(\phi,\kos{f},\alpha)\colon (\kos{S}^{\pr},\kos{t}^{\pr},\omega^{\pr})\mor (\kos{S},\kos{t},\omega)$
  to be the $2$-cell
  consisting of the identity Galois transformation
  of $\phi$ and the invertible Galois $\kos{S}$-transformation
  $\breve{\alpha}\colon \kos{f}\breve{\omega}^{\pr}_{\phi}\cong \breve{\omega}\breve{\kos{f}}$
  in Proposition~\ref{prop-PKGalfund-inducedPKfundmor}.
  It is straightforward to check,
  using \cite[Proposition~6.2.16]{Johnson2021}, that
  these data indeed define an invertible strong transformation
  in the $2$-category of pseudofunctors
  over $\mathsf{pGalKos}_{\etsf}$.
  This completes the proof of Theorem~\ref{thm-PKGalCat-equiv21catoveret}.
\end{proof}

%% file: Galcovsp.tex
\chapter{Galois Theory of Covering Spaces}
\label{appdix-Galcovsp}

We denote by $\kosSet=(\cat{Set},\times,\pt)$
the Galois prekosmos of sets,
where $\times$ is the cartesian product
and $\pt$ is the singleton set.
In this appendix,
we reformulate the Galois theory of covering spaces
in the language of Galois categories over
the Galois prekosmos $\kosSet$.
A space means a topological space,
and a map means a continuous map.
By a \emph{uc-space} we mean a nonempty,
path-connected, locally path-connected,
semilocally simply connected space.
We denote by $\mathsf{Top}^{\mathsf{uc}}$
the $(2,1)$-category of uc-spaces,
maps, and homotopy classes of homotopies.
Here homotopies are taken up to homotopy relative to their end maps:
for maps $f$, $g\colon X^{\pr}\to X$,
two homotopies $H$, $\tilde{H}\colon f\simeq g$ are identified
if they are homotopic relative to $X^{\pr}\times \{0,1\}$.
The reformulation takes the following form.

\begin{theorem}
  \label{thm-Galcovsp-Cov}
  We have a pseudofunctor 
  \begin{equation*}
    \mathsf{Cov}
    \colon
    \mathsf{Top}^{\mathsf{uc}}
    \to 
    \mathsf{pGalCat}(\kosSet)
  \end{equation*}
  that sends each uc-space $X$
  to the Galois $\kosSet$-category
  of covering spaces of $X$.
\end{theorem}

From Theorem~\ref{thm-Galcovsp-Cov}
we obtain its pointed variant.
We denote by $\mathsf{Top}^{\mathsf{uc}}_{\mathsf{pt}}$
the $(2,1)$-category
whose $0$-cells are pointed uc-spaces,
whose $1$-cells are maps equipped with a homotopy class of paths
from the image of the base point to the base point,
and whose $2$-cells are homotopy classes of homotopies.
We denote by $\mathsf{Grp}$ the $(2,1)$-category 
of groups, morphisms of groups, and conjugations,
which is identified with the $(2,1)$-category
$\mathsf{pGalObj}(\kosSet)$
of Galois objects in the Galois prekosmos $\kosSet$.
The pointed variant takes the following form,
where 
$\pi_1\colon \mathsf{Top}^{\mathsf{uc}}_{\mathsf{pt}}\to\mathsf{Grp}$
is the $2$-functor that sends each pointed uc-space
to its fundamental group,
and 
$\mathsf{Rep}_{\kosSet}$ and $\mathsf{Aut}_{\kosSet}$
are the mutually quasi-inverse equivalences
in Theorem~\ref{thm-dualityoverGalpk-equiv21cat}
for $\kos{S}=\kosSet$.

\begin{theorem}
  \label{thm-Galcovsp-Covpt}
  We have a pseudofunctor 
  \begin{equation*}
    \mathsf{Cov}_{\mathsf{pt}}
    \colon
    \mathsf{Top}^{\mathsf{uc}}_{\mathsf{pt}}
    \to 
    \mathsf{pGalCat}_{\mathsf{pt}}(\kosSet)
  \end{equation*}
  that sends each pointed uc-space $(X,x)$
  to the Galois $\kosSet$-category
  of covering spaces of $X$
  equipped with the fiber point
  at $x$.
  Moreover, the following diagram of pseudofunctors
  commutes up to equivalence.
  \begin{equation}
    \label{eq-Galcovsp-Covpt}
    \vcenter{\hbox{
      \xymatrix@R=20pt@C=15pt{
        \text{ }
        &\mathsf{Top}^{\mathsf{uc}}_{\mathsf{pt}}
        \ar[dl]_-{\pi_1}
        \ar[dr]^-{\mathsf{Cov}_{\mathsf{pt}}}
        &\text{ }
        \\%1
        \mathsf{Grp}
        \ar@<1ex>[rr]^-{\mathsf{Rep}_{\kosSet}}
        \ar@{}[rr]|-{\simeq}
        &\text{ }
        &\mathsf{pGalCat}_{\mathsf{pt}}(\kosSet)
        \ar@<1ex>[ll]^-{\mathsf{Aut}_{\kosSet}}
      }
    }}
  \end{equation}
\end{theorem}

Let $X$ be a space and $x\in X$.
We denote by $\cat{Cov}(X)$
the category of coverings of $X$.
We allow the preimage of an evenly covered neighborhood to be empty,
so that the empty space is a covering of $X$.
For a covering $E\to X$,
we write $E|_U$ for the preimage of $U\subset X$
and $E_x$ for its fiber over $x$.
A morphism $\phi\colon E\to F$ in $\cat{Cov}(X)$
restricts to a map
$\phi|_U\colon E|_U\to F|_U$
over $U$
and to a map
$\phi_x\colon E_x\to F_x$
over $x$.
For a set $S$, regarded as a discrete space,
the projection $X\times S\to X$ is the trivial covering of $X$.
For a path $\gamma$ in $X$,
$[\gamma]$ denotes its homotopy class relative to its endpoints.
We write $\delta\gamma$
for the concatenation of paths $\delta$, $\gamma$
with $\gamma(1)=\delta(0)$,
which follows $\gamma$ first and then $\delta$.
For each point $y\in X$ that can be joined to $x$ by a path,
the set $\pi_1(X;x,y)$ of homotopy classes of paths from $x$ to $y$
is a right torsor under the fundamental group $\pi_1(X,x)$,
with the action given by concatenation.
In the proof of the following proposition,
we refer to \cite[section~1.3]{Hatcher2002}
for the homotopy lifting property,
the unique lifting property,
and the construction of the based universal covering.

\begin{proposition}
  \label{prop-Galcovsp-coveringGalsetcat}
  Let $X$ be a uc-space.
  \begin{enumerate}
    \item 
    Covering spaces of $X$ form a Galois prekosmos,
    denoted by $\kosCov(X)=(\cat{Cov}(X),\times_X,X)$.
    
    \item
    We have a Galois morphism
    $\kos{t}_{X\vphantom{!}}^{\vphantom{!}}\colon \kosCov(X)\mor\kosSet$
    whose right adjoint part $\kos{t}_{X\vphantom{!}}^!$
    sends each set $S$ to the trivial covering $X\times S\to X$.
    
    \item
    For each $x\in X$,
    we have a surjective fiber point
    $\omega_{x\vphantom{!}}^{\vphantom{!}}\colon \vphantom{\id}\smash[t]{\id}_{\kosSet\vphantom{!}}^{\vphantom{!}}\mor\kos{t}_{X\vphantom{!}}^{\vphantom{!}}$
    whose right adjoint part $\omega_{x\vphantom{!}}^!$
    sends each covering $E\to X$ to its fiber $E_x$ over $x$.
  \end{enumerate}
  In particular, $\kos{t}_{X\vphantom{!}}^{\vphantom{!}}$ is a Galois $\kosSet$-category.
  \begin{equation*}
    \begin{tikzcd}[row sep=5pt,column sep=2pt]
      \kosSet
      \ar[rr,filled,"\omega_{x}"]
      \ar[ddr,equal]
      &\text{ }
      &\kosCov(X)
      \ar[ddl,filled,"\kos{t}_{X}"]
      \\%1
      \text{ }
      &\text{ }
      &\text{ }
      \ar[ll,2cell,pos=0.58,shift left=0.3ex,"{\what{\omega_x}}"']
      \\%2
      \text{ }
      &\kosSet
      &\text{ }
    \end{tikzcd}
  \end{equation*}
\end{proposition}

We refer to $\omega_{x\vphantom{!}}^{\vphantom{!}}\colon \vphantom{\id}\smash[t]{\id}_{\kosSet\vphantom{!}}^{\vphantom{!}}\mor\kos{t}_{X\vphantom{!}}^{\vphantom{!}}$
as the \emph{fiber point} of $\kos{t}_{X\vphantom{!}}^{\vphantom{!}}$ \emph{at} $x\in X$.

\begin{proof}[Proof of Proposition~\ref{prop-Galcovsp-coveringGalsetcat}]
  Finite products and coequalizers of coverings of $X$
  are given by those as spaces over $X$,
  as the resulting spaces are again coverings of $X$.
  For finite products this is true for any $X$.
  Let $\phi$, $\psi\colon E\to F$
  be morphisms in $\cat{Cov}(X)$
  and $q\colon F\to Q$ be the coequalizer as spaces over $X$.
  Since $X$ is locally path-connected,
  $X$ can be covered by path-connected open subsets $U$
  that are evenly covered by both $E$ and $F$.
  For each such $U$, choose $u\in U$ and identify
  $E|_U=U\times E_u$ and $F|_U=U\times F_u$.
  Since $U$ is path-connected,
  we have $\phi|_U=I_U\times \phi_u$
  and $\psi|_U=I_U\times \psi_u$.
  Thus the saturation under $q$ of each sheet of $F|_U$
  is a union of sheets and hence is open in $F$.
  Since $q$ is a quotient map,
  the distinct images of these sheets form a disjoint open cover of $Q|_U$,
  each mapping homeomorphically onto $U$.
  Thus $Q|_U\cong U\times Q_u$,
  and $Q$ is a covering of $X$.
  This proves (1).

  The functor
  $\kos{t}_{X\vphantom{!}}^!
  \colon \cat{Set}\to \cat{Cov}(X)$,
  sending each $S$ to $X\times S\to X$,
  preserves finite products
  and is thus strong monoidal.
  It is right adjoint to the functor
  $\kos{t}_{X!}^{\vphantom{!}}
  \colon \cat{Cov}(X)\to \cat{Set}$
  that sends each covering $E$ of $X$ to
  the set $\pi_0(E)$ of its path components.
  Indeed, given $E$ and $S$,
  maps $E\to X\times S$ over $X$
  are in natural bijection with functions $\pi_0(E)\to S$,
  since the path components of $E$ are open.
  Thus
  $\kos{t}_{X\vphantom{!}}^{\vphantom{!}}
  =(\kos{t}_{X!}^{\vphantom{!}}\dashv
  \kos{t}_{X\vphantom{!}}^!)
  \colon \kosCov(X)\mor\kosSet$ is a Galois morphism.
  This proves (2).

  Fix $x\in X$
  and let $(\tilde{X}^x,\tilde{x})\to (X,x)$
  be the based universal covering.
  A point in $\tilde{X}^x$ is a homotopy class
  $[\gamma]$ of paths starting at $x$,
  and is sent to the endpoint $\gamma(1)$ in $X$.
  We have 
  $\tilde{X}^x_x=\pi_1(X,x)$,
  $\tilde{X}^x_y=\pi_1(X;x,y)$ for each $y\in X$,
  and
  $\tilde{x}=[c_x]$,
  where $c_x$ is the constant path at $x$.
  We write $T_{\delta}\colon E_y\to E_z$
  for the transport in a covering $E$
  along a path $\delta$ from $y$ to $z$ in $X$,
  sending each $e\in E_y$ to the endpoint of the lift starting at $e$.
  The homotopy lifting property and the unique lifting property imply
  that $T_{\delta}$ is bijective,
  is natural in $E$
  and depends only on $[\delta]$,
  and that $T_{\delta\gamma}=T_{\delta}\circ T_{\gamma}$
  for composable paths $\delta$, $\gamma$.
  The transport $\tilde{T}^x$ in $\tilde{X}^x$ is given by
  $\tilde{T}^x_{\delta}([\gamma])=[\delta\gamma]$.
  For each $e\in E_x$, the function
  $\tilde{X}^x\to E$
  sending $[\gamma]$ to $T_{\gamma}(e)$
  is a map.
  Indeed, on each sheet over a sufficiently small
  path-connected open subset $U$ evenly covered by both $\tilde{X}^x$ and $E$,
  it is the canonical homeomorphism onto a sheet of $E$ over $U$.
  By the unique lifting property,
  the map
  $[\gamma]\mapsto T_{\gamma}(e)$
  over $X$ is the unique one
  sending $\tilde{x}$ to $e$.
  By repeating this on each copy of $\tilde{X}^x$,
  we deduce that for each set $S$,
  every function $S\to E_x$ extends uniquely
  to a map $\tilde{X}^x\times S\to E$ over $X$.
  Thus we have the following natural bijection,
  given by evaluation on $\{\tilde{x}\}\times S$.
  \begin{equation}
    \label{eq-Galcovsp-univcov}
    \Hom_{\cat{Cov}(X)}(\tilde{X}^x\times S,E)
    \cong
    \Hom_{\cat{Set}}(S,E_x)
  \end{equation}

  The functor
  $\omega_{x\vphantom{!}}^!\colon\cat{Cov}(X)\to\cat{Set}$,
  sending each covering $E$ of $X$ to $E_x$,
  preserves finite products by (1) and is thus strong monoidal.
  By \eqref{eq-Galcovsp-univcov}
  it is right adjoint to
  the functor
  $\omega_{x!}^{\vphantom{!}}$
  that sends each set $S$
  to the covering $\tilde{X}^x\times S$ of $X$.
  The unit $\eta_S\colon S\to\pi_1(X,x)\times S$
  and counit $\epsilon_E\colon\tilde{X}^x\times E_x\to E$
  are given by $s\mapsto([c_x],s)$
  and $([\gamma],e)\mapsto T_{\gamma}(e)$, respectively.
  Hence $\omega_{x\vphantom{!}}^{\vphantom{!}}=(\omega_{x!}^{\vphantom{!}}\dashv\omega_{x\vphantom{!}}^!)\colon\kosSet\mor\kosCov(X)$
  is a Galois morphism.
  The natural transformation
  $\varphi$ associated to 
  $\omega_{x\vphantom{!}}^{\vphantom{!}}$
  in \eqref{eq-Galprekos-projformula} is given by
  $\varphi_{S,E}\colon
  \omega_{x!}^{\vphantom{!}}(S\times\omega_{x\vphantom{!}}^!(E))\to\omega_{x!}^{\vphantom{!}}(S)\times_X E$,
  $([\gamma],s,e)\mapsto(([\gamma],s),T_{\gamma}(e))$.
  Since $T_{\gamma}$ is bijective,
  $\varphi_{S,E}$ is a bijective morphism of coverings
  and hence a homeomorphism.
  Thus $\omega_{x\vphantom{!}}^{\vphantom{!}}$
  satisfies the projection formula.
  The invertible Galois transformation
  $\what{\vphantom{\omega_x}\smash[t]{\omega_{x\vphantom{!}}^{\vphantom{!}}}}\colon
  \kos{t}_{X\vphantom{!}}^{\vphantom{!}}\omega_{x\vphantom{!}}^{\vphantom{!}}\cong\vphantom{\id}\smash[t]{\id}_{\kosSet\vphantom{!}}^{\vphantom{!}}$
  is the comonoidal natural isomorphism
  $\what{\omega_{x\vphantom{!}}^!}\colon
  \omega_{x\vphantom{!}}^!\kos{t}_{X\vphantom{!}}^!\cong\vphantom{\id}\smash[t]{\id}_{\kosSet\vphantom{!}}^{\vphantom{!}}$
  between the right adjoint parts,
  whose components are given by $(X\times S)_x=\{x\}\times S\cong S$.
  Thus $\omega_{x\vphantom{!}}^{\vphantom{!}}\colon\vphantom{\id}\smash[t]{\id}_{\kosSet\vphantom{!}}^{\vphantom{!}}\mor\kos{t}_{X\vphantom{!}}^{\vphantom{!}}$
  is a locally connected Galois $\kosSet$-morphism,
  and by Remark~\ref{rem-Galprekos-omegalowershriekconservative},
  it is a fiber point of $\kos{t}_{X\vphantom{!}}^{\vphantom{!}}$.

  Let $\phi\colon E\to F$
  be a morphism in $\cat{Cov}(X)$.
  Assume that $\phi_x$ is bijective.
  Since $X$ is path-connected,
  every point $y$ in $X$ can be joined to $x$ by a path,
  and transport along this path
  commutes with the morphisms
  $\phi_x$, $\phi_y$ on the fibers.
  It follows that $\phi_y$ is bijective
  for every $y\in X$.
  Thus $\phi$ is a bijective morphism of coverings
  and hence a homeomorphism.
  This shows that 
  $\omega_{x\vphantom{!}}^!$ is conservative,
  and it also preserves coequalizers by (1).
  Hence the fiber point
  $\omega_{x\vphantom{!}}^{\vphantom{!}}$ is surjective.
  This proves (3).

  Since $X$ is nonempty,
  it follows from (3) that $\kos{t}_{X\vphantom{!}}^{\vphantom{!}}$
  is a Galois $\kosSet$-category.
\end{proof}

We explain how the concepts introduced in
Chapters~\ref{chap-Galprekos} and~\ref{chap-Galcatoverpk}
reduce to their classical counterparts when $\kos{S}=\kosSet$.
A cocommutative comonoid in $\kosSet$ is a set with its diagonal,
so that $\kosEns(\kosSet)=\kosSet$.
Since $\kosSet$ is cartesian closed,
Galois objects in $\kosSet$ are the groups,
and we identify $\mathsf{pGalObj}(\kosSet)$
in \eqref{eq-Galobjinpk-pGalObjS21cat}
with $\mathsf{Grp}$.
Let $\pi$ be a group and recall Lemma~\ref{lem-RepGalobj-Reppi}.
The Galois prekosmos
$\kosRep_{\kosSet}(\pi)$
is the cartesian monoidal category of left $\pi$-sets.
The right adjoint part
$\kos{t}_{\pi\vphantom{!}}^!$
of the Galois morphism
$\kos{t}_{\pi\vphantom{!}}^{\vphantom{!}}\colon \kosRep_{\kosSet}(\pi)\mor\kosSet$
sends each set to itself equipped with the trivial $\pi$-action.
The right adjoint part
$\kos{u}_{\pi\vphantom{!}}^!$
of the surjective fiber point
$\kos{u}_{\pi\vphantom{!}}^{\vphantom{!}}\colon 
\vphantom{\id}\smash[t]{\id}_{\kosSet\vphantom{!}}^{\vphantom{!}}\mor\kos{t}_{\pi\vphantom{!}}^{\vphantom{!}}$
sends each left $\pi$-set to its underlying set.
A set is a co-fpqc covering of $\pt$
if and only if it is nonempty.
Thus a co-fpqc right $\pi$-torsor over $\pt$
is a right $\pi$-torsor in the usual sense,
that is, a nonempty set with a free and transitive right $\pi$-action.

\begin{lemma}
  \label{lem-Galcovsp-comonoidalautomatic}
  Let $\kos{T}=(\CT,\tensor,\unit)$,
  $\kos{T}^{\pr}=(\CT^{\pr},\tensor^{\pr},\unit^{\pr})$
  be cartesian monoidal categories.
  \begin{enumerate}
    \item 
    Every functor $\mon{f}\colon \CT^{\pr}\to \CT$
    has a unique colax monoidal functor structure,
    whose coherence morphisms are
    the obstructions for $\mon{f}$
    to preserve binary products and terminal objects.
    Hence $\mon{f}$ equipped with this unique structure
    is strong monoidal
    if and only if $\mon{f}$ preserves finite products.

    \item
    Every natural transformation
    $\vartheta\colon \mon{f}\Rightarrow\mon{g}$
    between functors $\mon{f}$, $\mon{g}\colon \CT^{\pr}\to \CT$
    is comonoidal with respect to the unique colax monoidal functor structures
    on $\mon{f}$ and $\mon{g}$.
  \end{enumerate}
\end{lemma}
\begin{proof}
  We leave the verification to the reader.
\end{proof}

\begin{lemma}
  \label{lem-Galcovsp-recognisingAut}
  Let $\kos{T}\xmor{\kos{t}}\kosSet$,
  $\kos{T}^{\pr}\xmor{\kos{t}^{\pr}}\kosSet$
  be Galois $\kosSet$-categories,
  and $\kos{f}$, $\kos{g}\colon\kos{t}^{\pr}\mor\kos{t}$
  be Galois $\kosSet$-morphisms.
  Then every natural transformation
  $\kos{f}^!\Rightarrow\kos{g}^!$
  between the right adjoint parts
  defines a Galois $\kosSet$-transformation
  $\kos{f}\Rightarrow\kos{g}$.
\end{lemma}
\begin{proof}
  Let $\vartheta\colon\kos{f}^!\Rightarrow\kos{g}^!$
  be a natural transformation.
  By Theorem~\ref{thm-Galcatoverpk-PrimaryObj} for $\kos{S}=\kosSet$,
  we have equivalences
  $\kos{t}_{\pi\vphantom{!}}^{\vphantom{!}}\simeq\kos{t}$,
  $\kos{t}_{\pi^{\pr}\vphantom{!}}^{\vphantom{!}}\simeq\kos{t}^{\pr}$
  for some groups $\pi$, $\pi^{\pr}$
  so that $\kos{T}$, $\kos{T}^{\pr}$ are cartesian monoidal.
  Thus $\vartheta$ is a comonoidal natural transformation
  by Lemma~\ref{lem-Galcovsp-comonoidalautomatic}.
  By Corollary~\ref{cor-Galcatoverpk-connsurj},
  $\kos{t}^{\pr}$ is connected.
  Thus every natural endomorphism of
  the fully faithful functor $\kos{t}^{\pr!}$
  is uniquely written as $\kos{t}^{\pr!}\xi$
  for some natural endomorphism 
  $\xi$ of the identity functor of $\cat{Set}$.
  Since the identity functor of $\cat{Set}$
  has no natural endomorphisms other than the identity,
  the same is also true for $\kos{t}^{\pr!}$.
  Thus the composite
  $\what{\kos{g}^!}\circ(\vartheta\kos{t}^!)\circ(\what{\kos{f}^!})^{-1}
  \colon\kos{t}^{\pr!}\Rightarrow\kos{t}^{\pr!}$
  must be the identity,
  and $\vartheta$ defines a Galois $\kosSet$-transformation
  $\kos{f}\Rightarrow\kos{g}$.
\end{proof}

\begin{remark}
  Lemma~\ref{lem-Galcovsp-recognisingAut}
  still holds under the weaker assumption
  that $\kos{t}$, $\kos{t}^{\pr}$ are Galois $\kosSet$-prekosmoi,
  that $\kos{T}$, $\kos{T}^{\pr}$ are cartesian monoidal
  and that $\kos{t}^{\pr}$ is connected.
\end{remark}

Let $\kos{T}\xmor{\kos{t}}\kosSet$
be a Galois $\kosSet$-category
and let
$\omega$, $\omega^{\pr}$ be fiber points of $\kos{t}$.
By Proposition~\ref{prop-Galcatoverpk-HomisIsom},
the presheaves
$\uIsom_{\kosSet}(\omega,\omega^{\pr})$
and
$\uAut_{\kosSet}(\omega)$
are representable and hence are determined by
their evaluations at $\pt$,
as is 
true for
any representable presheaf of sets or groups on $\cat{Set}$.
By Lemma~\ref{lem-Galcovsp-recognisingAut}
for
$\kos{t}^{\pr}=
\vphantom{\id}\smash[t]{\id}_{\kosSet\vphantom{!}}^{\vphantom{!}}$,
these evaluations are the set 
$\Isom(\omega^!,\omega^{\pr!})$
of natural isomorphisms $\omega^!\cong\omega^{\pr!}$
and the group
$\Aut(\omega^!)$
of natural automorphisms of $\omega^!$.
In other words, these presheaves
associated with $\omega$, $\omega^{\pr}$
carry exactly the same information
as the corresponding set and group
associated with
$\omega_{\vphantom{!}}^!$, $\omega_{\vphantom{!}}^{\pr!}$
as ordinary functors.

By applying the results of Chapter~\ref{chap-Galcatoverpk},
we obtain the following lemmas.

\begin{lemma}
  \label{lem-Galcovsp-mainthmcovering}
  Let $X$ be a uc-space and $x\in X$.
  Recall from
  Proposition~\ref{prop-Galcovsp-coveringGalsetcat}
  the Galois $\kosSet$-category $\kos{t}_{X\vphantom{!}}^{\vphantom{!}}$ of coverings of $X$
  and its fiber point $\omega_{x\vphantom{!}}^{\vphantom{!}}$ at $x$.
  \begin{enumerate}
    \item
    The presheaf of groups
    $\uAut_{\kosSet}(\omega_{x\vphantom{!}}^{\vphantom{!}})$
    on $\cat{Set}$
    is represented by
    $\pi_1(X,x)$.
    
    \item
    The fiber point $\omega_{x\vphantom{!}}^{\vphantom{!}}$ extends to an equivalence
    $\breve{\omega}_{x\vphantom{!}}^{\vphantom{!}}\colon \kos{t}_{\pi_1(X,x)}\simeq \kos{t}_{X\vphantom{!}}^{\vphantom{!}}$    
    of Galois $\kosSet$-prekosmoi
    up to an invertible Galois $\kosSet$-transformation
    $\breve{\omega}_{x\vphantom{!}}^{\vphantom{!}}\kos{u}_{\pi_1(X,x)\vphantom{!}}^{\vphantom{!}}\cong\omega_{x\vphantom{!}}^{\vphantom{!}}$
    given by the equality
    $\kos{u}_{\pi_1(X,x)\vphantom{!}}^!\breve{\omega}_{x\vphantom{!}}^!=\omega_{x\vphantom{!}}^!$.
    The right adjoint part
    $\breve{\omega}_{x\vphantom{!}}^!$
    sends each covering of $X$ to its fiber over $x$
    equipped with the monodromy action.
    
    \item
    For each $y\in X$,
    let $\omega_{y\vphantom{!}}^{\vphantom{!}}$ be the fiber point at $y$.
    The presheaf
    $\uIsom_{\kosSet}(\omega_{x\vphantom{!}}^{\vphantom{!}},\omega_{y\vphantom{!}}^{\vphantom{!}})$
    on $\cat{Set}$
    is represented by the right $\pi_1(X,x)$-torsor $\pi_1(X;x,y)$,
    which corresponds to $\omega_{y\vphantom{!}}^{\vphantom{!}}$
    under the adjoint equivalence in Theorem~\ref{thm-pFibGalTors-perfcorres}
    for $\omega=\omega_{x\vphantom{!}}^{\vphantom{!}}$.
  \end{enumerate}
  \begin{equation*}
    \begin{aligned}
      \Hom_{\cat{Set}}(\slot,\pi_1(X,x))
      &\cong
      \uAut_{\kosSet}(\omega_{x\vphantom{!}}^{\vphantom{!}})
      \\%1
      \Hom_{\cat{Set}}(\slot,\pi_1(X;x,y))
      &\cong
      \uIsom_{\kosSet}(\omega_{x\vphantom{!}}^{\vphantom{!}},\omega_{y\vphantom{!}}^{\vphantom{!}})
    \end{aligned}
    \qquad\quad
    \begin{tikzcd}[row sep=5pt,column sep=2pt]
      \text{ }
      &\id_{\kosSet}
      \ar[ddr,filled,"{\omega_x}"]
      \ar[ddl,filled,"{\kos{u}_{\pi_1(X,x)}}"']
      &\text{ }
      \\%1
      \text{ }
      \ar[rr,2cell,shift right=0.7ex,pos=0.55]
      &\text{ }
      &\text{ }
      \\%2
      \kos{t}_{\pi_1(X,x)}
      \ar[rr,filled,"{\breve{\omega}_x}","\simeq"']
      &\text{ }
      &\kos{t}_X
    \end{tikzcd}
  \end{equation*}
\end{lemma}
\begin{proof}
  We have
  $\omega_{x\vphantom{!}}^!\omega_{x!}^{\vphantom{!}}(\pt)
  =\pi_1(X,x)\times\pt
  \cong\pi_1(X,x)$,
  and 
  $\omega_{y\vphantom{!}}^!\omega_{x!}^{\vphantom{!}}(\pt)
  =\pi_1(X;x,y)\times \pt
  \cong \pi_1(X;x,y)$
  for each $y\in X$.
  Up to these bijections,
  the natural isomorphism
  in \eqref{eq-ActionbaseGal-hatarupperlowershriek}
  for $\omega=\omega_{x\vphantom{!}}^{\vphantom{!}}$, $\omega^{\pr}=\omega_{y\vphantom{!}}^{\vphantom{!}}$
  agrees with the identity natural transformation of
  $\pi_1(X;x,y)\times\slot$,
  the group structure \eqref{eq-GalMainthm-pigroupstr}
  on $\omega_{x\vphantom{!}}^!\omega_{x!}^{\vphantom{!}}(\pt)$
  agrees with that of the fundamental group $\pi_1(X,x)$,
  the action \eqref{eq-GalMainthm-gamma}
  of $\omega_{x\vphantom{!}}^!\omega_{x!}^{\vphantom{!}}(\pt)$
  on $\omega_{x\vphantom{!}}^!(E)$ for each covering $E$ of $X$
  agrees with the monodromy action of $\pi_1(X,x)$ on the fiber $E_x$,
  and the right action
  of $\omega_{x\vphantom{!}}^!\omega_{x!}^{\vphantom{!}}(\pt)$
  on $\omega_{y\vphantom{!}}^!\omega_{x!}^{\vphantom{!}}(\pt)$
  in Lemma~\ref{lem-ActionbaseGal-inverseandcomposition}~(2)
  agrees with the right action of
  $\pi_1(X,x)$ on $\pi_1(X;x,y)$
  given by concatenation of paths.
  Thus (1) and (2) follow from 
  Theorem~\ref{thm-Galcatoverpk-mainthm},
  and (3) follows from
  Proposition~\ref{prop-Galcatoverpk-HomisIsom}
  and Theorem~\ref{thm-pFibGalTors-perfcorres}.
\end{proof}

\begin{lemma}
  \label{lem-Galcovsp-pullback}
  Let $f\colon X^{\pr}\to X$ be a map between uc-spaces.
  We have a fundamental $\kosSet$-morphism
  $\kos{f}\colon \kos{t}_{X^{\pr}\vphantom{!}}^{\vphantom{!}}\mor\kos{t}_{X\vphantom{!}}^{\vphantom{!}}$
  whose right adjoint part
  is the pullback along $f$.
\end{lemma}
\begin{proof}
  We write $E$ for a covering of $X$.
  The pullback $f^*\colon \cat{Cov}(X)\to \cat{Cov}(X^{\pr})$,
  sending each $E$ to $X^{\pr}\times_XE$,
  preserves finite products,
  and is thus strong monoidal.
  Choose a point $x^{\pr}$ in $X^{\pr}$
  and
  put $x\coloneqq f(x^{\pr})$.
  Recall $\omega_{x\vphantom{!}}^{\vphantom{!}}$,
  $\omega_{x^{\pr}\vphantom{!}}^{\vphantom{!}}$,
  as well as 
  $\breve{\omega}_{x\vphantom{!}}^{\vphantom{!}}$,
  $\breve{\omega}_{x^{\pr}\vphantom{!}}^{\vphantom{!}}$
  in Lemma~\ref{lem-Galcovsp-mainthmcovering}~(2).
  We have a natural isomorphism
  $\omega_{x^{\pr}\vphantom{!}}^!f^*
  \cong \omega_{x\vphantom{!}}^!$
  whose component at $E$ is
  $f^*(E)_{x^{\pr}}\cong E_x$.
  Write $h\colon \pi_1(X^{\pr},x^{\pr})\to\pi_1(X,x)$
  for the morphism of groups
  $[\gamma^{\pr}]\mapsto[f\circ \gamma^{\pr}]$,
  and $\kos{h}=(\cat{Ind}(h)\dashv\cat{Res}(h))$
  for the fundamental $\kosSet$-morphism
  associated to $h$
  in Lemma~\ref{lem-fundmoroverGalpk-fundmor}.
  The right adjoint part
  $\cat{Res}(h)\colon
  \cat{Rep}_{\kosSet}(\pi_1(X,x))
  \to
  \cat{Rep}_{\kosSet}(\pi_1(X^{\pr},x^{\pr}))$
  of $\kos{h}$ is the restriction along $h$.
  Then the natural isomorphism
  $\omega_{x^{\pr}\vphantom{!}}^!f^*
  \cong \omega_{x\vphantom{!}}^!$
  lifts to a natural isomorphism
  $\breve{\omega}_{x^{\pr}\vphantom{!}}^!f^*
  \cong
  \cat{Res}(h)\breve{\omega}_{x\vphantom{!}}^!$,
  since for each $E$,
  the transport in $f^*(E)$
  along a loop $\gamma^{\pr}$ at $x^{\pr}$
  is the transport in $E$ along $f\circ \gamma^{\pr}$
  under $f^*(E)_{x^{\pr}}\cong E_x$.
  Hence $f^*$ has a left adjoint,
  and we fix one, say $f_!$.
  Then 
  $\kos{f}=(f_!\dashv f^*)$
  is a Galois morphism.
  By Lemma~\ref{lem-Galcovsp-comonoidalautomatic},
  the natural isomorphism
  $\what{f^*}\colon f^*(X\times \slot)\cong X^{\pr}\times \slot$
  defines an invertible Galois transformation
  $\what{\kos{f}}\colon\kos{t}_{X\vphantom{!}}^{\vphantom{!}}\kos{f}\cong \kos{t}_{X^{\pr}\vphantom{!}}^{\vphantom{!}}$.
  By Lemma~\ref{lem-Galcovsp-recognisingAut},
  the natural isomorphism
  $\omega_{x^{\pr}\vphantom{!}}^!f^*
  \cong\omega_{x\vphantom{!}}^!$
  defines an invertible Galois $\kosSet$-transformation
  $\kos{f}\omega_{x^{\pr}\vphantom{!}}^{\vphantom{!}}
  \cong\omega_{x\vphantom{!}}^{\vphantom{!}}$.
  Thus $\kos{f}$ is a fundamental $\kosSet$-morphism.
\end{proof}

% We define the $(2,1)$-category
% $\mathsf{Top}^{\mathsf{uc}}$
% as follows.

\begin{proof}[Proof of Theorem~\ref{thm-Galcovsp-Cov}]
  The pseudofunctor $\mathsf{Cov}$
  sends each uc-space $X$
  to the Galois $\kosSet$-category
  $\kos{t}_{X\vphantom{!}}^{\vphantom{!}}$
  of coverings of $X$ in Proposition~\ref{prop-Galcovsp-coveringGalsetcat},
  and each map $f\colon X^{\pr}\to X$
  between uc-spaces
  to the fundamental $\kosSet$-morphism
  $\kos{f}=(f_!\dashv f^*)\colon
  \kos{t}_{X^{\pr}\vphantom{!}}^{\vphantom{!}}
  \mor\kos{t}_{X\vphantom{!}}^{\vphantom{!}}$
  in Lemma~\ref{lem-Galcovsp-pullback}.

  We describe how $\mathsf{Cov}$ sends $2$-cells.
  Let $f$, $g\colon X^{\pr}\to X$ be $1$-cells in $\mathsf{Top}^{\mathsf{uc}}$
  and $\kos{f}$, $\kos{g}\colon \kos{t}_{X^{\pr}\vphantom{!}}^{\vphantom{!}}\mor\kos{t}_{X\vphantom{!}}^{\vphantom{!}}$
  be the associated fundamental $\kosSet$-morphisms.
  A $2$-cell $[H]\colon f\cong g$ in $\mathsf{Top}^{\mathsf{uc}}$,
  represented by a homotopy $H\colon f\simeq g$,
  is sent to the invertible Galois $\kosSet$-transformation
  $\vartheta\colon\kos{f}\cong \kos{g}$
  defined as follows.
  The component of $\vartheta$ at each covering $E$ of $X$
  is the map $\vartheta_E\colon f^*(E)\to g^*(E)$ over $X^{\pr}$
  that sends $(x^{\pr},e)$ to $(x^{\pr},T_{H(x^{\pr},\slot)}(e))$.
  It is indeed a map,
  as its second component is obtained by lifting the homotopy
  $f^*(E)\times[0,1]\to X$,
  $((x^{\pr},e),t)\mapsto H(x^{\pr},t)$
  to $E$ with initial map $(x^{\pr},e)\mapsto e$,
  and then evaluating the lift at $t=1$.
  Note that $\vartheta_E$ depends only on $[H]$
  by homotopy invariance of transport,
  and its inverse is given by
  the transport along the reversed path
  $H(x^{\pr},\slot)^{-1}$.
  Since maps of coverings commute with transport,
  the components $\vartheta_E$
  define a natural isomorphism $\vartheta\colon f^*\cong g^*$.
  By Lemma~\ref{lem-Galcovsp-recognisingAut},
  $\vartheta$ is indeed an invertible Galois $\kosSet$-transformation
  from $\kos{f}$ to $\kos{g}$.

  The equality
  $T_{\delta\gamma}=T_\delta\circ T_{\gamma}$
  holds for composable paths $\delta$, $\gamma$,
  and the transport along a constant path is the identity.
  Hence $\mathsf{Cov}$ preserves vertical composition of $2$-cells
  and identity $2$-cells,
  and thus defines a functor on each hom-groupoid.
  The constraints of the pseudofunctor $\mathsf{Cov}$
  are defined as follows.
  \begin{itemize}
    \item 
    Let $f\colon X^{\pr}\to X$ and $f^{\pr}\colon X^{\ppr}\to X^{\pr}$
    be a composable pair of $1$-cells in $\mathsf{Top}^{\mathsf{uc}}$,
    with composite $f^{\ppr}\coloneqq ff^{\pr}$.
    We write $\kos{f}$, $\kos{f}^{\pr}$, and $\kos{f}^{\ppr}$
    for the fundamental $\kosSet$-morphisms
    associated to $f$, $f^{\pr}$, and $f^{\ppr}$, respectively.
    The component of the functorial constraint of $\mathsf{Cov}$
    at $f$, $f^{\pr}$
    is the invertible Galois $\kosSet$-transformation
    $\kos{f}\kos{f}^{\pr}\cong\kos{f}^{\ppr}$,
    given by
    $X^{\ppr}\times_{X^{\pr}}(X^{\pr}\times_XE)
    \cong X^{\ppr}\times_XE$
    for each covering $E$ of $X$.

    \item
    Let $X$ be a uc-space.
    We write $\kos{I}_X\colon \kos{t}_{X\vphantom{!}}^{\vphantom{!}}\mor\kos{t}_{X\vphantom{!}}^{\vphantom{!}}$
    for the fundamental $\kosSet$-morphism
    associated to $I_X\colon X\to X$
    in Lemma~\ref{lem-Galcovsp-pullback}.
    The component of the unity constraint of $\mathsf{Cov}$ at $X$
    is the invertible Galois $\kosSet$-transformation
    from the identity Galois $\kosSet$-morphism of $\kos{t}_{X\vphantom{!}}^{\vphantom{!}}$
    to $\kos{I}_X$,
    given by
    $E\cong X\times_XE$
    for each covering $E$ of $X$.
  \end{itemize}
  The above natural isomorphisms define invertible Galois
  $\kosSet$-transformations
  by Lemma~\ref{lem-Galcovsp-recognisingAut}.
  The functorial and unity constraints
  satisfy the associativity and left/right unity axioms,
  which follows from the universal property
  of iterated fiber products.

  It remains to show that $\mathsf{Cov}$ preserves
  horizontal composition of $2$-cells
  up to the functorial constraint.
  Let $f$, $g\colon X^{\pr}\to X$,
  $f^{\pr}$, $g^{\pr}\colon X^{\ppr}\to X^{\pr}$
  be maps between uc-spaces
  and
  $H\colon f\simeq g$,
  $H^{\pr}\colon f^{\pr}\simeq g^{\pr}$
  be homotopies.
  Put $f^{\ppr}\coloneqq ff^{\pr}$ and $g^{\ppr}\coloneqq gg^{\pr}$.
  The horizontal composite of $[H]$, $[H^{\pr}]$
  is represented by the homotopy
  $H^{\ppr}\colon f^{\ppr}\simeq g^{\ppr}$,
  $H^{\ppr}(x^{\ppr},t)\coloneqq H(H^{\pr}(x^{\ppr},t),t)$.
  It suffices to show that
  for each covering $E$ of $X$ and each $x^{\ppr}\in X^{\ppr}$,
  we have 
  $T_{g\circ H^{\pr}(x^{\ppr},\slot)}\circ
  T_{H(f^{\pr}(x^{\ppr}),\slot)}
  =T_{H^{\ppr}(x^{\ppr},\slot)}$.
  Consider the map $[0,1]^2\to X$ that sends 
  $(s,t)$ to $H(H^{\pr}(x^{\ppr},s),t)$.
  The image of the diagonal gives the transport on the right,
  whereas the image of the boundary path
  from $(0,0)$ through $(0,1)$ to $(1,1)$
  gives the composite of the transports on the left.
  The diagonal and this boundary path are homotopic
  in $[0,1]^2$ relative to their endpoints,
  so the equality follows from homotopy invariance of transport.
  This completes the construction of
  $\mathsf{Cov}\colon \mathsf{Top}^{\mathsf{uc}}\to\mathsf{pGalCat}(\kosSet)$.
\end{proof}

We define the $(2,1)$-category
$\mathsf{Top}^{\mathsf{uc}}_{\mathsf{pt}}$
as follows.
A $0$-cell is a pointed uc-space $(X,x)$.
A $1$-cell $(f,[\alpha])\colon (X^{\pr},x^{\pr})\to (X,x)$
is a pair consisting of a map $f\colon X^{\pr}\to X$
and the homotopy class of a path $\alpha$ in $X$
from $f(x^{\pr})$ to $x$.
The composite of $(f,[\alpha])$
and $(f^{\pr},[\alpha^{\pr}])\colon (X^{\ppr},x^{\ppr})\to (X^{\pr},x^{\pr})$
is $(ff^{\pr},[\alpha(f\circ\alpha^{\pr})])$,
and the identity $1$-cell of $(X,x)$ is $(I_X,[c_x])$.
A $2$-cell $[H]\colon (f,[\alpha])\cong (g,[\beta])$
is a $2$-cell $[H]\colon f\cong g$ in $\mathsf{Top}^{\mathsf{uc}}$.
Horizontal and vertical compositions of $2$-cells,
as well as identity $2$-cells,
are those of $\mathsf{Top}^{\mathsf{uc}}$.
Thus the $0$-cells and $1$-cells are those of the coslice
$(2,1)$-category of $\mathsf{Top}^{\mathsf{uc}}$ under the singleton space,
whereas the $2$-cells are inherited from $\mathsf{Top}^{\mathsf{uc}}$
without any compatibility condition;
see Remark~\ref{rem-dualityoverGalpk-notcoslice} for the analogy
with $\mathsf{pGalCat}_{\mathsf{pt}}(\kosSet)$.

We define the $2$-functor
$\pi_1\colon \mathsf{Top}^{\mathsf{uc}}_{\mathsf{pt}}\to\mathsf{Grp}$
as follows.
A $0$-cell $(X,x)$ is sent to $\pi_1(X,x)$.
A $1$-cell $(f,[\alpha])\colon (X^{\pr},x^{\pr})\to (X,x)$
is sent to the morphism of groups
$\pi_1(f,[\alpha])\colon \pi_1(X^{\pr},x^{\pr})\to\pi_1(X,x)$
given by
$[\gamma^{\pr}]\mapsto[\alpha(f\circ\gamma^{\pr})\alpha^{-1}]$
for each loop $\gamma^{\pr}$ at $x^{\pr}$.
A $2$-cell $[H]\colon (f,[\alpha])\cong (g,[\beta])$
is sent to the element
$[\beta H(x^{\pr},\slot)\alpha^{-1}]\in\pi_1(X,x)$.
It is straightforward to check that
the assignment of $\pi_1$ on each cell is well defined
and depends only on the indicated homotopy classes,
and that $\pi_1$ is indeed a $2$-functor.

\begin{proof}[Proof of Theorem~\ref{thm-Galcovsp-Covpt}]
  We describe the pseudofunctor $\mathsf{Cov}_{\mathsf{pt}}$.
  A pointed uc-space $(X,x)$ is sent to
  $(\kos{t}_{X\vphantom{!}}^{\vphantom{!}},\omega_{x\vphantom{!}}^{\vphantom{!}})$
  in Proposition~\ref{prop-Galcovsp-coveringGalsetcat}.
  A $1$-cell $(f,[\alpha])\colon (X^{\pr},x^{\pr})\to (X,x)$
  is sent to the pair
  consisting of the fundamental $\kosSet$-morphism
  $\kos{f}\colon \kos{t}_{X^{\pr}\vphantom{!}}^{\vphantom{!}}\mor\kos{t}_{X\vphantom{!}}^{\vphantom{!}}$
  associated to $f$ in Lemma~\ref{lem-Galcovsp-pullback},
  and the invertible Galois $\kosSet$-transformation
  $\check{T}_\alpha\colon\kos{f}\omega_{x^{\pr}\vphantom{!}}^{\vphantom{!}}\cong\omega_{x\vphantom{!}}^{\vphantom{!}}$
  whose component at each covering $E$ of $X$
  is the composite of
  $f^*(E)_{x^{\pr}}\cong E_{f(x^{\pr})}$
  followed by
  $T_{\alpha}\colon E_{f(x^{\pr})}\cong E_x$.
  The assignment on $1$-cells is well defined by
  naturality and homotopy invariance of transport,
  together with Lemma~\ref{lem-Galcovsp-recognisingAut}.
  A $2$-cell
  $[H]\colon (f,[\alpha])\cong (g,[\beta])$
  is sent to the $2$-cell
  $\vartheta\colon\kos{f}\cong\kos{g}$,
  the image of $[H]\colon f\cong g$ under $\mathsf{Cov}$.
  The functorial and unity constraints of $\mathsf{Cov}_{\mathsf{pt}}$
  are given by those of $\mathsf{Cov}$,
  and the pseudofunctor axioms for $\mathsf{Cov}_{\mathsf{pt}}$
  follow from those for $\mathsf{Cov}$.

  It remains to show the commutativity of \eqref{eq-Galcovsp-Covpt}
  up to equivalence.
  It suffices to show that
  the composite $\mathsf{Aut}_{\kosSet}\circ\mathsf{Cov}_{\mathsf{pt}}$
  is a $2$-functor,
  and is $2$-naturally isomorphic to $\pi_1$.
  We describe how the composite sends
  $0$-cells, $1$-cells, and $2$-cells.
  \begin{itemize}
    \item 
    A $0$-cell $(X,x)$ is sent to the group
    $\omega_{x\vphantom{!}}^!\omega_{x!}^{\vphantom{!}}(\pt)$.
    
    \item
    A $1$-cell $(f,[\alpha])\colon(X^{\pr},x^{\pr})\to(X,x)$
    is sent to the morphism of groups
    $(\check{T}_\alpha)^{\sharp}_{\pt}\colon
    \omega_{x^{\pr}\vphantom{!}}^!\omega_{x^{\pr}!}^{\vphantom{!}}(\pt)
    \to \omega_{x\vphantom{!}}^!\omega_{x!}^{\vphantom{!}}(\pt)$,
    where 
    $(\check{T}_\alpha)_!\colon
    \omega_{x^{\pr}!}^{\vphantom{!}}\Rightarrow f^*\omega_{x!}^{\vphantom{!}}$
    is the mate of
    $\check{T}_\alpha^{-1}\colon
    \omega_{x\vphantom{!}}^!
    \cong \omega_{x^{\pr}\vphantom{!}}^!f^*$
    with respect to the adjunctions
    $\omega_{x\vphantom{!}}^{\vphantom{!}}$
    and $\omega_{x^{\pr}\vphantom{!}}^{\vphantom{!}}$,
    as in \eqref{eq-fundmoroverGalpk-matealpha},
    and
    $(\check{T}_\alpha)^{\sharp}
    \coloneqq (\check{T}_\alpha\omega_{x!}^{\vphantom{!}})
    \circ (\omega_{x^{\pr}\vphantom{!}}^!(\check{T}_\alpha)_!)
    \colon \omega_{x^{\pr}\vphantom{!}}^!\omega_{x^{\pr}!}^{\vphantom{!}}
    \Rightarrow\omega_{x\vphantom{!}}^!\omega_{x!}^{\vphantom{!}}$
    is as in the proof of Proposition~\ref{prop-fundmoroverGalpk-realization}.
    
    \item
    A $2$-cell $[H]\colon(f,[\alpha])\cong(g,[\beta])$
    is sent to the $2$-cell represented by the map
    $\vartheta^{\sharp}_{\pt}
    \colon\pt
    \to\omega_{x\vphantom{!}}^!\omega_{x!}^{\vphantom{!}}(\pt)$,
    where
    $\vartheta\colon \kos{f}\cong \kos{g}$
    is the image of $[H]$ under $\mathsf{Cov}_{\mathsf{pt}}$,
    and
    $\vartheta^{\sharp}
    \coloneqq (\check{T}_\beta\omega_{x!}^{\vphantom{!}})
    \circ (\omega_{x^{\pr}\vphantom{!}}^!\vartheta
    \omega_{x!}^{\vphantom{!}})
    \circ (\check{T}_\alpha^{-1}\omega_{x!}^{\vphantom{!}})
    \circ \eta$
    is as in the proof of
    Proposition~\ref{prop-dualityoverGalpk-AutS2funct},
    with $\eta$ the unit of $\omega_{x\vphantom{!}}^{\vphantom{!}}$.
  \end{itemize}
  We show that
  the functorial and unity constraints of $\mathsf{Cov}_{\mathsf{pt}}$
  are sent to identity $2$-cells
  under $\mathsf{Aut}_{\kosSet}$.
  \begin{itemize}
    \item 
    Let
    $(f,[\alpha])\colon(X^{\pr},x^{\pr})\to(X,x)$,
    $(f^{\pr},[\alpha^{\pr}])\colon
    (X^{\ppr},x^{\ppr})\to(X^{\pr},x^{\pr})$
    be a composable pair of $1$-cells
    in $\mathsf{Top}^{\mathsf{uc}}_{\mathsf{pt}}$,
    with composite denoted by $(f^{\ppr},[\alpha^{\ppr}])$.
    The component
    $\vartheta\colon \kos{f}\kos{f}^{\pr}\cong\kos{f}^{\ppr}$
    of the functorial constraint at this pair is given by
    $X^{\ppr}\times_{X^{\pr}}(X^{\pr}\times_XE)
    \cong X^{\ppr}\times_XE$
    for each covering $E$ of $X$.
    On the fiber at $x^{\ppr}\in X^{\ppr}$,
    we have
    $\check{T}_{\alpha^{\ppr}}
    \circ (\omega_{x^{\ppr}\vphantom{!}}^!\vartheta)
    =\check{T}_\alpha\circ(\check{T}_{\alpha^{\pr}}f^*)
    \colon \omega_{x^{\ppr}\vphantom{!}}^!f^{\pr*}f^*\cong \omega_{x\vphantom{!}}^!$.
    Following the definition,
    we have
    $\vartheta^{\sharp}
    =(\check{T}_{\alpha^{\ppr}}\omega_{x!}^{\vphantom{!}})
    \circ(\omega_{x^{\ppr}\vphantom{!}}^!\vartheta
    \omega_{x!}^{\vphantom{!}})
    \circ((\check{T}_\alpha\circ(\check{T}_{\alpha^{\pr}}f^*))^{-1}
    \omega_{x!}^{\vphantom{!}})
    \circ\eta$,
    which is equal to $\eta$.
    Hence $\vartheta$ is sent to the identity $2$-cell
    under $\mathsf{Aut}_{\kosSet}$.

    \item
    Let $(X,x)$ be a pointed uc-space.
    The component $\vartheta$ of the unity constraint at $(X,x)$
    is given by $E\cong X\times_XE$
    for each covering $E$ of $X$.
    On the fiber at $x\in X$, we have
    $\check{T}_{c_x}\circ(\omega_{x\vphantom{!}}^!\vartheta)
    =I_{\omega_{x\vphantom{!}}^!}
    \colon\omega_{x\vphantom{!}}^!\cong\omega_{x\vphantom{!}}^!$.
    Following the definition,
    we have
    $\vartheta^{\sharp}
    =(\check{T}_{c_x}\omega_{x!}^{\vphantom{!}})
    \circ(\omega_{x\vphantom{!}}^!\vartheta
    \omega_{x!}^{\vphantom{!}})
    \circ\eta$,
    which is equal to $\eta$.
    Hence $\vartheta$ is sent to the identity $2$-cell
    under $\mathsf{Aut}_{\kosSet}$.
  \end{itemize}
  This shows that the composite
  $\mathsf{Aut}_{\kosSet}\circ\mathsf{Cov}_{\mathsf{pt}}$
  is a $2$-functor.

  Let $(X,x)$ be a pointed uc-space.
  Recall from Lemma~\ref{lem-Galcovsp-mainthmcovering}
  that the bijection
  $\omega_{x\vphantom{!}}^!\omega_{x!}^{\vphantom{!}}(\pt)
  =\pi_1(X,x)\times\pt
  \cong\pi_1(X,x)$
  is an isomorphism of groups.
  We show that these bijections constitute a $2$-natural isomorphism
  from $\mathsf{Aut}_{\kosSet}\circ\mathsf{Cov}_{\mathsf{pt}}$ to $\pi_1$.
  For the rest of the proof,
  we identify
  $\omega_{x!}^{\vphantom{!}}(\pt)$
  with the universal covering $\tilde{X}^x$,
  $\omega_{x\vphantom{!}}^!\omega_{x!}^{\vphantom{!}}(\pt)$
  with $\pi_1(X,x)$,
  and $\omega_{y\vphantom{!}}^!\omega_{x!}^{\vphantom{!}}(\pt)$
  with $\pi_1(X;x,y)$ for each $y\in X$,
  as well as the corresponding structures compared
  in Lemma~\ref{lem-Galcovsp-mainthmcovering}.
  It then suffices to verify the following.
  \begin{itemize}
    \item 
    Let 
    $(f,[\alpha])\colon(X^{\pr},x^{\pr})\to(X,x)$
    be a $1$-cell.
    We show that
    $(\check{T}_\alpha)^{\sharp}_{\pt}([\gamma^{\pr}])
    =[\alpha(f\circ\gamma^{\pr})\alpha^{-1}]$
    for each loop $\gamma^{\pr}$ at $x^{\pr}$.
    Indeed, the component of $(\check{T}_\alpha)_!$ at $\pt$
    is a map $((\check{T}_\alpha)_!)_{\pt}
    \colon (\tilde{X}^{\pr})^{x^{\pr}}\to f^*(\tilde{X}^x)$
    of coverings of $X^{\pr}$.
    It sends $[c_{x^{\pr}}]$ to $(x^{\pr},[\alpha^{-1}])$,
    by the left relation in \eqref{eq-fundmoroverGalpk-matealpha}
    and the description of the adjunction unit
    in Proposition~\ref{prop-Galcovsp-coveringGalsetcat}.
    Since $((\check{T}_\alpha)_!)_{\pt}$ is compatible with transport
    and the transport in $f^*(\tilde{X}^x)$ along $\gamma^{\pr}$
    projects to transport in $\tilde{X}^x$
    along $f\circ\gamma^{\pr}$,
    it sends each $[\gamma^{\pr}]$ to
    $(x^{\pr},\tilde{T}^x_{f\circ\gamma^{\pr}}([\alpha^{-1}]))
    =(x^{\pr},[(f\circ\gamma^{\pr})\alpha^{-1}])$.
    Hence the composite
    $(\check{T}_\alpha)^{\sharp}_{\pt}
    =(\check{T}_\alpha)_{\tilde{X}^x}
    \circ\omega_{x^{\pr}\vphantom{!}}^!
    (((\check{T}_\alpha)_!)_{\pt})$
    sends $[\gamma^{\pr}]$ to 
    $(x^{\pr},[(f\circ\gamma^{\pr})\alpha^{-1}])$
    under $\omega_{x^{\pr}\vphantom{!}}^!
    (((\check{T}_\alpha)_!)_{\pt})$,
    then to
    $\tilde{T}^x_{\alpha}[(f\circ\gamma^{\pr})\alpha^{-1}]
    =[\alpha(f\circ\gamma^{\pr})\alpha^{-1}]$
    under $(\check{T}_\alpha)_{\tilde{X}^x}$.
    
    \item
    Let $[H]\colon(f,[\alpha])\cong(g,[\beta])$
    be a $2$-cell between $1$-cells
    from $(X^{\pr},x^{\pr})$ to $(X,x)$,
    and write $\vartheta\colon\kos{f}\cong\kos{g}$
    for its image under $\mathsf{Cov}_{\mathsf{pt}}$.
    We show that the map
    $\vartheta^{\sharp}_{\pt}\colon\pt\to\pi_1(X,x)$
    has value $[\beta H(x^{\pr},\slot)\alpha^{-1}]$.
    Indeed, the composite
    $\vartheta^{\sharp}_{\pt}
    =(\check{T}_\beta)_{\tilde{X}^x}
    \circ
    \omega_{x^{\pr}\vphantom{!}}^!
    (\vartheta_{\tilde{X}^x})
    \circ(\check{T}_\alpha^{-1})_{\tilde{X}^x}
    \circ\eta_{\pt}$
    sends the unique element in $\pt$ to 
    $[c_x]$ under $\eta_{\pt}$,
    then to 
    $(x^{\pr},\tilde{T}^x_{\alpha^{-1}}([c_x]))
    =(x^{\pr},[\alpha^{-1}])$
    under $(\check{T}_\alpha^{-1})_{\tilde{X}^x}$,
    then to
    $(x^{\pr},\tilde{T}^x_{H(x^{\pr},\slot)}[\alpha^{-1}])
    =(x^{\pr},[H(x^{\pr},\slot)\alpha^{-1}])$
    under $\omega_{x^{\pr}\vphantom{!}}^!
    (\vartheta_{\tilde{X}^x})$,
    and then to
    $\tilde{T}^x_{\beta}([H(x^{\pr},\slot)\alpha^{-1}])
    =[\beta H(x^{\pr},\slot)\alpha^{-1}]$
    under $(\check{T}_\beta)_{\tilde{X}^x}$.
  \end{itemize}
  We thus obtain a $2$-natural isomorphism
  $\mathsf{Aut}_{\kosSet}\circ\mathsf{Cov}_{\mathsf{pt}}\cong \pi_1$.
  This completes the proof of Theorem~\ref{thm-Galcovsp-Covpt}.
\end{proof}
 

%% file: ltcatoversmc.tex
\chapter{Lax Tensor Categories over a Symmetric Monoidal Category}
\label{appdix-ltcatoversmc}

Throughout this appendix,
let $\kos{K}=(\CK,\otimes,\kappa)$
be a symmetric monoidal category.
We study symmetric monoidal categories
equipped with a lax monoidal functor from $\kos{K}$,
lax monoidal functors between these categories
that preserve the structures over $\kos{K}$
up to a monoidal natural transformation,
and monoidal natural transformations between these functors
that are compatible over $\kos{K}$.
Using these notions,
we provide an abstract characterization of 
the endofunctor on $\kos{K}$
given by tensoring with 
a commutative monoid,
a commutative bimonoid,
or a commutative Hopf monoid in $\kos{K}$.

We denote by 
$\mathtt{SMC}_{\mathtt{lax}}$
the $2$-category of symmetric monoidal categories,
lax monoidal functors,
and monoidal natural transformations.
We write 
$\kos{K}/\mathtt{SMC}_{\mathtt{lax}}$
for the lax coslice $2$-category of 
$\mathtt{SMC}_{\mathtt{lax}}$
under $\kos{K}$.
\begin{itemize}
  \item 
  A \emph{lax $\kos{K}$-tensor category},
  or a \emph{lax tensor category over $\kos{K}$},
  is a $0$-cell in $\kos{K}/\mathtt{SMC}_{\mathtt{lax}}$.
  It is a pair $(\kos{T},\monbase{t})$
  consisting of a symmetric monoidal category $\kos{T}$
  and a lax monoidal functor $\monbase{t}\colon\kos{K}\to\kos{T}$.
  We write the pair as $\kos{K}\xrightarrow{\monbase{t}}\kos{T}$,
  or simply as $\monbase{t}$.
  An example is the pair consisting of
  $\kos{K}$ and its identity functor $\id_{\kos{K}}$.

  \item
  A \emph{lax $\kos{K}$-tensor functor},
  or a \emph{lax tensor functor over $\kos{K}$},
  is a $1$-cell in $\kos{K}/\mathtt{SMC}_{\mathtt{lax}}$.
  Let $\kos{K}\xrightarrow{\monbase{t}}\kos{T}$,
  $\kos{K}\xrightarrow{\monbase{t}^{\pr}}\kos{T}^{\pr}$
  be lax $\kos{K}$-tensor categories.
  A lax $\kos{K}$-tensor functor
  $(\mon{f},\what{\mon{f}})\colon\monbase{t}^{\pr}\to\monbase{t}$
  is a pair consisting of
  a lax monoidal functor
  $\mon{f}\colon \kos{T}^{\pr}\to \kos{T}$
  and a monoidal natural transformation
  $\what{\mon{f}}\colon\monbase{t}\Rightarrow\mon{f}\monbase{t}^{\pr}$.
  We omit $\what{\mon{f}}$ from the notation
  and write $\mon{f}\colon\monbase{t}^{\pr}\to \monbase{t}$.

  \item
  A \emph{lax $\kos{K}$-tensor natural transformation},
  or a \emph{lax tensor natural transformation over $\kos{K}$},
  is a $2$-cell in $\kos{K}/\mathtt{SMC}_{\mathtt{lax}}$.
  Let $\mon{f}$, $\mon{g}\colon \monbase{t}^{\pr}\to \monbase{t}$
  be lax $\kos{K}$-tensor functors.
  A lax $\kos{K}$-tensor natural transformation
  $\vartheta\colon \mon{f}\Rightarrow\mon{g}$
  is a monoidal natural transformation
  between the underlying lax monoidal functors
  $\mon{f}$, $\mon{g}\colon \kos{T}^{\pr}\to \kos{T}$
  that satisfies
  $\what{\mon{g}}=(\vartheta\monbase{t}^{\pr})\circ \what{\mon{f}}
  \colon \monbase{t}\Rightarrow\mon{g}\monbase{t}^{\pr}$.
  A \emph{lax $\kos{K}$-tensor natural isomorphism},
  or a \emph{lax tensor natural isomorphism over $\kos{K}$},
  is an invertible $2$-cell in $\kos{K}/\mathtt{SMC}_{\mathtt{lax}}$.
\end{itemize}
Let $\mon{f}$, $\mon{g}\colon \monbase{t}^{\pr}\to \monbase{t}$
be lax $\kos{K}$-tensor functors.
We write
$(\kos{K}/\mathtt{SMC}_{\mathtt{lax}})[\monbase{t}^{\pr},\monbase{t}]$
for the category of lax $\kos{K}$-tensor functors
from $\monbase{t}^{\pr}$ to $\monbase{t}$.
We illustrate
$\mon{f}$, $\mon{g}$
and a lax $\kos{K}$-tensor natural transformation
$\vartheta\colon \mon{f}\Rightarrow\mon{g}$
as follows.
\begin{equation*}
  \begin{tikzcd}[row sep=5pt,column sep=5pt]
    \kos{T}^{\pr}
    \ar[rr,"\mon{f}"]
    &\text{ }
    &\kos{T}
    &\text{ }
    &\kos{T}^{\pr}
    \ar[rr,"\mon{g}"]
    &\text{ }
    &\kos{T}
    &\text{ }
    &\kos{T}^{\pr}
    \ar[rr,shift left=0.5ex,bend left=22pt,"\mon{g}" description]
    \ar[rr,shift right=0.5ex,bend right=22pt,"\mon{f}" description]
    &\text{ }
    &\kos{T}
    \\
    \text{ }
    &\text{ }
    &\text{ }
    \ar[ll,2cell,shift left=0.5ex,"{\,\,\what{\mon{f}}}"']
    &\text{ }
    &\text{ }
    &\text{ }
    &\text{ }
    \ar[ll,2cell,shift left=0.5ex,"{\,\what{\mon{g}}}"']
    &\text{ }
    &\text{ }
    &\text{ }
    \ar[u,2cell,pos=2.1,shift right=0.1ex,"{\vartheta}"']
    &\text{ }
    \\
    \text{ }
    &\kos{K}
    \ar[uul,"\monbase{t}^{\pr}"]
    \ar[uur,"\monbase{t}"']
    &\text{ }
    &\text{ }
    &\text{ }
    &\kos{K}
    \ar[uul,"\monbase{t}^{\pr}"]
    \ar[uur,"\monbase{t}"']
    &\text{ }
    &\text{ }
    &\text{ }
    &\kos{K}
    \ar[uul,"\monbase{t}^{\pr}"]
    \ar[uur,"\monbase{t}"']
    &\text{ }
  \end{tikzcd}
\end{equation*}

Let $\kos{K}\xrightarrow{\monbase{t}}\kos{T}$
be a lax $\kos{K}$-tensor category,
where $\kos{T}=(\CT,\tensor,\unit)$.
The \emph{induced $\kos{K}$-action} $\acts$ on $\kos{T}$
is the functor $\acts\colon \CK\times\CT\to\CT$
sending each pair $(z,X)$ to $z\acts X\coloneqq\monbase{t}(z)\tensor X$.
Let $\kos{K}\xrightarrow{\monbase{t}^{\pr}}\kos{T}^{\pr}$
be a lax $\kos{K}$-tensor category,
where $\kos{T}^{\pr}=(\CT^{\pr},\tensor^{\pr},\unit^{\pr})$
and $\acts^{\pr}$ is the induced $\kos{K}$-action on $\kos{T}^{\pr}$.
Let $\mon{f}\colon \monbase{t}^{\pr}\to \monbase{t}$
be a lax $\kos{K}$-tensor functor.
The \emph{$\kos{K}$-equivariance}
of $\mon{f}$ is the natural transformation
\begin{equation*}
  \cevar{\mon{f}}\colon\slot\acts \mon{f}(\slot)
  \Rightarrow\mon{f}(\slot \acts^{\pr}\slot)
\end{equation*}
between functors from $\CK\times \CT^{\pr}$ to $\CT$
whose component at $z\in \CK$, $X^{\pr}\in \CT^{\pr}$ is
\begin{equation*}
  \xymatrix@C=50pt{
    \monbase{t}(z)\tensor \mon{f}(X^{\pr})
    \ar[r]^-{(\what{\mon{f}})_z\tensor I_{\mon{f}(X^{\pr})}}
    &\mon{f}\monbase{t}^{\pr}(z)\tensor \mon{f}(X^{\pr})
    \ar[r]^-{(\mon{f}_{\tensor})_{\monbase{t}^{\pr}(z),X^{\pr}}}
    &\mon{f}(\monbase{t}^{\pr}(z)\tensor^{\pr} X^{\pr})
    .
  }
\end{equation*}
Let $\kos{K}\xrightarrow{\monbase{t}^{\ppr}} \kos{T}^{\ppr}$
be a lax $\kos{K}$-tensor category,
where $\kos{T}^{\ppr}=(\CT^{\ppr},\tensor^{\ppr},\unit^{\ppr})$
and $\acts^{\ppr}$ is the induced $\kos{K}$-action on $\kos{T}^{\ppr}$.
The composite 
$\mon{f}\mon{f}^{\pr}\colon\monbase{t}^{\ppr}\to \monbase{t}$
of lax $\kos{K}$-tensor functors
$\mon{f}\colon\monbase{t}^{\pr}\to\monbase{t}$,
$\mon{f}^{\pr}\colon\monbase{t}^{\ppr}\to \monbase{t}^{\pr}$
consists of
$\mon{f}\mon{f}^{\pr}\colon\kos{T}^{\ppr}\to\kos{T}$
and
$\what{\mon{f}\mon{f}^{\pr}}
=(\mon{f}\what{\mon{f}^{\pr}})\circ \what{\mon{f}}
\colon\monbase{t}\Rightarrow \mon{f}\mon{f}^{\pr}\monbase{t}^{\ppr}$.
The $\kos{K}$-equivariance of the composite is
$(\cevar{\mon{f}\mon{f}^{\pr}})_{z,X^{\ppr}}
=\mon{f}((\cevar{\mon{f}^{\pr}})_{z,X^{\ppr}})
\circ (\cevar{\mon{f}})_{z,\mon{f}^{\pr}(X^{\ppr})}
\colon z\acts\mon{f}\mon{f}^{\pr}(X^{\ppr})\to \mon{f}\mon{f}^{\pr}(z\acts^{\ppr} X^{\ppr})$
for $z\in \CK$, $X^{\ppr}\in\CT^{\ppr}$.
The identity lax $\kos{K}$-tensor functor of $\monbase{t}$,
denoted by $\id_{\kos{T}}$,
is the pair consisting of the identity functor of $\kos{T}$
and the identity natural transformation of $\monbase{t}$.

\begin{lemma}
  \label{lem-ltcatoversmc-laxKtnatequiv}
  Let
  $\kos{K}\xrightarrow{\monbase{t}} \kos{T}$,
  $\kos{K}\xrightarrow{\monbase{t}^{\pr}} \kos{T}^{\pr}$
  be lax $\kos{K}$-tensor categories,
  where
  $\kos{T}=(\CT,\tensor,\unit)$,
  $\kos{T}^{\pr}=(\CT^{\pr},\tensor^{\pr},\unit^{\pr})$.
  Let $\acts$ and $\acts^{\pr}$
  be the induced $\kos{K}$-actions
  on $\kos{T}$ and $\kos{T}^{\pr}$, respectively.
  Let 
  $\mon{f}$, $\mon{g}\colon \monbase{t}^{\pr}\to \monbase{t}$
  be lax $\kos{K}$-tensor functors
  and $\vartheta\colon\mon{f}\Rightarrow\mon{g}$
  be a monoidal natural transformation
  between the underlying lax monoidal functors
  $\mon{f}$, $\mon{g}\colon\kos{T}^{\pr}\to \kos{T}$.
  Then $\vartheta\colon\mon{f}\Rightarrow\mon{g}$
  is a lax $\kos{K}$-tensor natural transformation
  if and only if
  the following diagram commutes for every pair
  $z\in \CK$, $X^{\pr}\in \CT^{\pr}$.
  \begin{equation}
    \label{eq-ltcatoversmc-laxKtnatequiv}
    \vcenter{\hbox{
      \xymatrix@R=15pt@C=35pt{
        z\acts \mon{f}(X^{\pr})
        \ar[r]^-{(\cevar{\mon{f}})_{z,X^{\pr}}}
        \ar[d]_-{I_z\acts \vartheta_{X^{\pr}}}
        &\mon{f}(z\acts^{\pr} X^{\pr})
        \ar[d]^-{\vartheta_{z\acts^{\pr} X^{\pr}}}
        \\
        z\acts \mon{g}(X^{\pr})
        \ar[r]^-{(\cevar{\mon{g}})_{z,X^{\pr}}}
        &\mon{g}(z\acts^{\pr} X^{\pr})
      }
    }}
  \end{equation}
\end{lemma}
\begin{proof}
  If $\vartheta$ is a lax $\kos{K}$-tensor natural transformation,
  the diagram \eqref{eq-ltcatoversmc-laxKtnatequiv}
  commutes by the definitions.
  Conversely, assume that
  \eqref{eq-ltcatoversmc-laxKtnatequiv} commutes
  for every pair $z$, $X^{\pr}$.
  One can check that
  $(\what{\mon{f}})_z\colon\monbase{t}(z)\to \mon{f}\monbase{t}^{\pr}(z)$
  is equal to
  $(\cevar{\mon{f}})_{z,\unit^{\pr}}\circ
  (I_{\monbase{t}(z)}\tensor\mon{f}_{\unit})
  \colon\monbase{t}(z)\tensor \unit\to
  \mon{f}(\monbase{t}^{\pr}(z)\tensor^{\pr} \unit^{\pr})$
  up to coherence isomorphisms.
  Using this description for both $\what{\mon{f}}$ and $\what{\mon{g}}$
  and applying \eqref{eq-ltcatoversmc-laxKtnatequiv} at $X^{\pr}=\unit^{\pr}$,
  we deduce that
  $(\what{\mon{g}})_z=
  \vartheta_{\monbase{t}^{\pr}(z)}\circ (\what{\mon{f}})_z
  \colon\monbase{t}(z)\to\mon{g}\monbase{t}^{\pr}(z)$
  for every $z$.
  Thus $\vartheta$ is a lax $\kos{K}$-tensor natural transformation.
\end{proof}

Let $\kos{K}\xrightarrow{\monbase{t}}\kos{T}$
be a lax $\kos{K}$-tensor category,
where $\kos{T}=(\CT,\tensor,\unit)$.
We denote by
$\kosComm(\kos{T})=(\cat{Comm}(\kos{T}),\tensor,\unit)$
the cocartesian monoidal category
of commutative monoids in $\kos{T}$.
The product and unit morphisms of a commutative monoid $B$ in $\kos{T}$
are denoted by $\pc_B\colon B\tensor B\to B$ and $u_B\colon \unit\to B$.
Each commutative monoid $B$ in $\kos{T}$
determines a lax $\kos{K}$-tensor functor
$\slot\acts B\colon \id_{\kos{K}}\to \monbase{t}$
consisting of the lax monoidal functor
$\slot\acts B=\monbase{t}(\slot)\tensor B
\colon \kos{K}\to \kos{T}$
and the monoidal natural transformation
$\what{\acts B}\colon \monbase{t}\Rightarrow \slot\acts B$
given by
$(\what{\acts B})_z
\coloneqq I_{\monbase{t}(z)}\tensor u_B
\colon \monbase{t}(z)\to \monbase{t}(z)\tensor B$.
The $\kos{K}$-equivariance of $\slot\acts B$ is
$(\cevar{\acts B})_{z,x}
=(\monbase{t}_{\tensor})_{z,x}\tensor I_B
\colon \monbase{t}(z)\tensor \monbase{t}(x)\tensor B
\to \monbase{t}(z\otimes x)\tensor B$
for $z$, $x\in \CK$.
Each morphism $f\colon B\to B^{\pr}$ in $\cat{Comm}(\kos{T})$
determines a lax $\kos{K}$-tensor natural transformation
$\slot\acts f\colon \slot\acts B\Rightarrow\slot\acts B^{\pr}$.

\begin{proposition}
  \label{prop-ltcatoversmc-CommlaxKtadj}
  Let $\kos{K}\xrightarrow{\monbase{t}}\kos{T}$
  be a lax $\kos{K}$-tensor category.
  We have an adjunction
  \begin{equation*}
    \vcenter{\hbox{
      \xymatrix@C=25pt{
        \cat{Comm}(\kos{T})
        \ar@<1ex>[r]^-{\iota}
        \ar@{}[r]|-{\bot}
        &(\kos{K}/\mathtt{SMC}_{\mathtt{lax}})[\id_{\kos{K}},\monbase{t}]
        \ar@<1ex>[l]^-{\mathcal{R}}
      }
    }}
  \end{equation*}
  between
  the category of commutative monoids in $\kos{T}$
  and the category of lax $\kos{K}$-tensor functors
  from $\id_{\kos{K}}$ to $\monbase{t}$.
\end{proposition}
\begin{proof}
  The functor $\iota$ sends each commutative monoid $B$ in $\kos{T}$ to $\slot\acts B$,
  and each morphism $f\colon B\to B^{\pr}$
  in $\cat{Comm}(\kos{T})$ to
  $\slot\acts f\colon \slot\acts B\Rightarrow\slot\acts B^{\pr}$.
  The functor $\mathcal{R}$ sends each lax $\kos{K}$-tensor functor
  $\mon{f}\colon\id_{\kos{K}}\to \monbase{t}$
  to the commutative monoid $\mon{f}(\kappa)$ in $\kos{T}$,
  and each morphism
  $\vartheta\colon \mon{f}\Rightarrow\mon{g}$
  in $(\kos{K}/\mathtt{SMC}_{\mathtt{lax}})[\id_{\kos{K}},\monbase{t}]$
  to the morphism
  $\vartheta_{\kappa}\colon \mon{f}(\kappa)\to \mon{g}(\kappa)$
  in $\cat{Comm}(\kos{T})$.
  Here $\mon{f}(\kappa)$ has product
  $\pc_{\mon{f}(\kappa)}=\mon{f}(\jmath^{-1}_{\kappa})\circ (\mon{f}_{\tensor})_{\kappa,\kappa}$
  and unit $u_{\mon{f}(\kappa)}=\mon{f}_{\unit}$.
  The component of the adjunction unit at each
  $B\in \cat{Comm}(\kos{T})$ is
  $\monbase{t}_\unit\tensor I_B\colon
  B\to \monbase{t}(\kappa)\tensor B$.
  The component of the adjunction counit at each
  lax $\kos{K}$-tensor functor
  $\mon{f}\colon \id_{\kos{K}}\to \monbase{t}$
  is the lax $\kos{K}$-tensor natural transformation
  \begin{equation}
    \label{eq-ltcatoversmc-tahardef}
    \tahar{\mon{f}}\colon
    \slot\acts \mon{f}(\kappa)\Rightarrow\mon{f}
  \end{equation}
  where 
  $(\tahar{\mon{f}})_z\coloneqq
  \mon{f}(\jmath^{-1}_z)\circ (\cevar{\mon{f}})_{z,\kappa}
  \colon \monbase{t}(z)\tensor \mon{f}(\kappa)\to \mon{f}(z)$
  for $z\in \CK$.
  We have
  $(\tahar{\acts B})_z=
  (\monbase{t}(\jmath^{-1}_z)\circ (\monbase{t}_{\tensor})_{z,\kappa})\tensor I_B
  \colon z\acts (\kappa\acts B)
  \to z\acts B$
  for $B\in\cat{Comm}(\kos{T})$, $z\in \CK$.
  It is straightforward to check that these data
  constitute an adjunction $\iota\dashv \mathcal{R}$.
\end{proof}

A \emph{strong $\kos{K}$-tensor category},
or a \emph{strong tensor category over $\kos{K}$},
is a pair
$(\kos{T},\monbase{t})$
consisting of a symmetric monoidal category $\kos{T}$
and a strong monoidal functor $\monbase{t}$.
It is in particular a lax $\kos{K}$-tensor category,
and we write it as 
$\kos{K}\xrightarrow{\monbase{t}}\kos{T}$,
or simply as $\monbase{t}$.
For instance, $\id_{\kos{K}}$ is a strong $\kos{K}$-tensor category.

\begin{definition}
  \label{def-ltcatoversmc-coreflective}
  Let $\kos{K}\xrightarrow{\monbase{t}}\kos{T}$
  be a strong $\kos{K}$-tensor category
  and $\mon{f}\colon \id_{\kos{K}}\to \monbase{t}$
  be a lax $\kos{K}$-tensor functor.
  \begin{enumerate}
    \item 
    We refer to 
    $\tahar{\mon{f}}\colon \slot\acts \mon{f}(\kappa)\Rightarrow\mon{f}$
    in \eqref{eq-ltcatoversmc-tahardef}
    as the \emph{coreflection} of $\mon{f}$.

    \item
    We say that $\mon{f}$ is \emph{coreflective}
    if its coreflection $\tahar{\mon{f}}$ is invertible.
  \end{enumerate}
  We write
  $(\kos{K}/\mathtt{SMC}_{\mathtt{lax}})[\id_{\kos{K}},\monbase{t}]_{\mathrm{crfl}}$
  for the full subcategory of
  $(\kos{K}/\mathtt{SMC}_{\mathtt{lax}})[\id_{\kos{K}},\monbase{t}]$
  consisting of coreflective
  lax $\kos{K}$-tensor functors
  from $\id_{\kos{K}}$ to $\monbase{t}$.
\end{definition}

\begin{proposition}
  \label{prop-ltcatoversmc-CommstrongKtequiv}
  Let $\kos{K}\xrightarrow{\monbase{t}}\kos{T}$
  be a strong $\kos{K}$-tensor category.
  Recall the adjunction
  $\iota\dashv \mathcal{R}$
  in Proposition~\ref{prop-ltcatoversmc-CommlaxKtadj}.
  \begin{enumerate}
    \item 
    The left adjoint
    $\iota\colon
    \cat{Comm}(\kos{T})
    \to (\kos{K}/\mathtt{SMC}_{\mathtt{lax}})[\id_{\kos{K}},\monbase{t}]$
    is fully faithful.

    \item 
    The essential image of the left adjoint
    $\iota$ is 
    $(\kos{K}/\mathtt{SMC}_{\mathtt{lax}})[\id_{\kos{K}},\monbase{t}]_{\mathrm{crfl}}$,
    which is a coreflective full subcategory
    of $(\kos{K}/\mathtt{SMC}_{\mathtt{lax}})[\id_{\kos{K}},\monbase{t}]$.

    \item 
    The adjunction $\iota\dashv \mathcal{R}$
    restricts to an adjoint equivalence of categories
    \begin{equation*}
      \vcenter{\hbox{
        \xymatrix@C=25pt{
          \cat{Comm}(\kos{T})
          \ar@<1ex>[r]^-{\iota}
          \ar@{}[r]|-{\simeq}
          &(\kos{K}/\mathtt{SMC}_{\mathtt{lax}})[\id_{\kos{K}},\monbase{t}]_{\mathrm{crfl}}.
          \ar@<1ex>[l]^-{\mathcal{R}}
        }
      }}
    \end{equation*}
  \end{enumerate}
\end{proposition}
\begin{proof}
  Since $\monbase{t}_{\unit}\colon \unit\to \monbase{t}(\kappa)$
  is invertible,
  the adjunction unit of $\iota\dashv \mathcal{R}$ is a natural isomorphism.
  Thus $\iota$ is fully faithful.
  This proves (1).
  A lax $\kos{K}$-tensor functor
  $\mon{f}\colon \id_{\kos{K}}\to \monbase{t}$
  lies in the essential image of $\iota$
  if and only if the component $\tahar{\mon{f}}$ of the adjunction counit at $\mon{f}$
  is invertible.
  Hence the essential image of $\iota$ is the full subcategory
  of coreflective lax $\kos{K}$-tensor functors,
  and its inclusion functor has a right adjoint $\iota\mathcal{R}$.
  This proves (2).
  Then (3) follows, as the adjunction counit becomes invertible
  after restriction.
\end{proof}

\begin{lemma}
  \label{lem-ltcatoversmc-coreflectionKequiv}
  Let $\kos{K}\xrightarrow{\monbase{t}}\kos{T}$
  be a strong $\kos{K}$-tensor category
  and $\mon{f}\colon \id_{\kos{K}}\to\monbase{t}$
  be a lax $\kos{K}$-tensor functor.
  Then the component of its coreflection
  $\tahar{\mon{f}}$ at $\kappa$ is
  $(\tahar{\mon{f}})_{\kappa}
  =(\monbase{t}_{\unit})^{-1}\tensor I_{\mon{f}(\kappa)}
  \colon \monbase{t}(\kappa)\tensor \mon{f}(\kappa)
  \cong \mon{f}(\kappa)$.
  Moreover, the following are equivalent:
  \begin{enumerate}
    \item 
    $\mon{f}$ is coreflective.

    \item
    $\mon{f}$ is isomorphic to
    $\slot\acts B\colon \id_{\kos{K}}\to \monbase{t}$
    as lax $\kos{K}$-tensor functors
    for some commutative monoid $B$ in $\kos{T}$.

    \item
    The $\kos{K}$-equivariance
    $\cevar{\mon{f}}\colon \slot\acts \mon{f}(\slot)
    \Rightarrow \mon{f}(\slot\otimes \slot)$
    of $\mon{f}$ is invertible.
  \end{enumerate}
\end{lemma}
\begin{proof}
  We have
  $(\tahar{\mon{f}})_{\kappa}
  \circ (\monbase{t}_{\unit}\tensor I_{\mon{f}(\kappa)})
  =I_{\mon{f}(\kappa)}$,
  as it is one of the triangle identities for the adjunction
  $\iota\dashv\mathcal{R}$
  in Proposition~\ref{prop-ltcatoversmc-CommlaxKtadj}.
  Since $\monbase{t}_{\unit}$ is invertible,
  we deduce the desired description of $(\tahar{\mon{f}})_{\kappa}$.
  The equivalence of (1) and (2) follows from
  Proposition~\ref{prop-ltcatoversmc-CommstrongKtequiv}.
  The definition of $\tahar{\mon{f}}$ in \eqref{eq-ltcatoversmc-tahardef}
  shows that (3) implies (1).
  One can check that
  $(\tahar{\mon{f}})_{z\otimes x}
  \circ ((\monbase{t}_{\tensor})_{z,x}\tensor I_{\mon{f}(\kappa)})
  =(\cevar{\mon{f}})_{z,x}
  \circ (I_{\monbase{t}(z)}\tensor (\tahar{\mon{f}})_x)
  \colon \monbase{t}(z)\tensor \monbase{t}(x)\tensor \mon{f}(\kappa)
  \to \mon{f}(z\otimes x)$
  holds for $z$, $x\in\CK$.
  Since $(\monbase{t}_{\tensor})_{z,x}$ is invertible,
  we deduce that (1) implies (3).
\end{proof}

By taking $\monbase{t}=\id_{\kos{K}}$
in Proposition~\ref{prop-ltcatoversmc-CommstrongKtequiv},
we obtain an adjoint equivalence of categories
\begin{equation}
  \label{eq-ltcatoversmc-CommKcrfllaxKtfunct}
  \vcenter{\hbox{
    \xymatrix@C=25pt{
      \cat{Comm}(\kos{K})
      \ar@<1ex>[r]^-{\iota}
      \ar@{}[r]|-{\simeq}
      &(\kos{K}/\mathtt{SMC}_{\mathtt{lax}})[\id_{\kos{K}},\id_{\kos{K}}]_{\mathrm{crfl}}.
      \ar@<1ex>[l]^-{\mathcal{R}}
    }
  }}
\end{equation}
We now show that the adjoint equivalence
\eqref{eq-ltcatoversmc-CommKcrfllaxKtfunct}
lifts to an adjoint equivalence of monoidal categories.
Let $\mon{f}$, $\mon{g}\colon \id_{\kos{K}}\to \id_{\kos{K}}$
be coreflective lax $\kos{K}$-tensor functors.
By Lemma~\ref{lem-ltcatoversmc-coreflectionKequiv},
the $\kos{K}$-equivariances 
$\cevar{\mon{f}}$, $\cevar{\mon{g}}$
of $\mon{f}$, $\mon{g}$ are invertible.
Hence the $\kos{K}$-equivariance $\cevar{\mon{f}\mon{g}}$ of the composite
$\mon{f}\mon{g}\colon \id_{\kos{K}}\to \id_{\kos{K}}$
is invertible,
and thus the composite $\mon{f}\mon{g}$ is also coreflective.
We write
\begin{equation*}
  \kosEnd[\id_{\kos{K}}]_{\mathrm{crfl}}
\end{equation*}
for the monoidal category of
coreflective lax $\kos{K}$-tensor endofunctors
on $\id_{\kos{K}}$,
whose monoidal product is given by composition
and whose unit object is the identity lax $\kos{K}$-tensor functor of $\id_{\kos{K}}$.
We also write
$\kosEnd[\id_{\kos{K}}]_{\mathrm{crfl}}^{\mathrm{rev}}$
for the reversed monoidal category,
whose monoidal product is given by composition in reverse order.

\begin{proposition}
  \label{prop-ltcatoversmc-CommKcrfllaxKtmonoidal}
  The adjoint equivalence of categories
  \eqref{eq-ltcatoversmc-CommKcrfllaxKtfunct}
  lifts to an adjoint equivalence
  of monoidal categories
  \begin{equation*}
    \vcenter{\hbox{
      \xymatrix@C=25pt{
        \kosComm(\kos{K})
        \ar@<1ex>[r]^-{\iota}
        \ar@{}[r]|-{\simeq}
        &\kosEnd[\id_{\kos{K}}]_{\mathrm{crfl}}^{\mathrm{rev}}.
        \ar@<1ex>[l]^-{\mathcal{R}}
      }
    }}
  \end{equation*}
  As a consequence, the monoidal category
  $\kosEnd[\id_{\kos{K}}]_{\mathrm{crfl}}$
  is cocartesian.
\end{proposition}
\begin{proof}
  Let $b$, $b^{\pr}$ be commutative monoids in $\kos{K}$.
  It follows from the Mac Lane coherence theorem
  that the coherence isomorphisms
  $a_{\slot,b,b^{\pr}}^{-1}
  \colon \iota(b^{\pr})\iota(b)
  \cong \iota(b\otimes b^{\pr})$
  and
  $\jmath
  \colon \id_{\kos{K}}
  \cong \iota(\kappa)$
  are lax $\kos{K}$-tensor natural isomorphisms.
  With these as its monoidal coherence isomorphisms,
  the left adjoint $\iota$ is a strong monoidal functor
  from $\kosComm(\kos{K})$ to $\kosEnd[\id_{\kos{K}}]_{\mathrm{crfl}}^{\mathrm{rev}}$.
  By doctrinal adjunction
  and the fact that $\iota\dashv \mathcal{R}$ is an adjoint equivalence,
  the right adjoint $\mathcal{R}$
  has a unique strong monoidal functor structure
  such that $\iota\dashv \mathcal{R}$
  becomes an adjoint equivalence of monoidal categories.
  Let $\mon{f}$, $\mon{g}$ 
  be coreflective lax $\kos{K}$-tensor endofunctors on $\id_{\kos{K}}$.
  One can check that 
  the monoidal coherence isomorphisms of $\mathcal{R}$
  are the isomorphisms
  $(\tahar{\mon{g}})_{\mon{f}(\kappa)}
  \colon
  \mon{f}(\kappa)\otimes \mon{g}(\kappa)
  \cong \mon{g}\mon{f}(\kappa)$
  and
  $I_{\kappa}\colon \kappa\cong\kappa$
  in $\cat{Comm}(\kos{K})$.
  Since $\kosComm(\kos{K})$ is cocartesian,
  so is $\kosEnd[\id_{\kos{K}}]_{\mathrm{crfl}}^{\mathrm{rev}}$,
  and hence also $\kosEnd[\id_{\kos{K}}]_{\mathrm{crfl}}$.
\end{proof}

A \emph{lax $\kos{K}$-tensor monad}
is a monad internal to $\kos{K}/\mathtt{SMC}_{\mathtt{lax}}$.
Since the monoidal category
$\kosEnd[\id_{\kos{K}}]_{\mathrm{crfl}}$
is cocartesian,
(i) every coreflective lax $\kos{K}$-tensor endofunctor on $\id_{\kos{K}}$
has a unique structure of a lax $\kos{K}$-tensor monad on $\id_{\kos{K}}$;
and
(ii)
every lax $\kos{K}$-tensor natural transformation
between coreflective lax $\kos{K}$-tensor endofunctors
on $\id_{\kos{K}}$
is a morphism of lax $\kos{K}$-tensor monads on $\id_{\kos{K}}$.
Let $\mon{f}$ be a coreflective lax $\kos{K}$-tensor endofunctor on $\id_{\kos{K}}$.
One can check that
the product and unit lax $\kos{K}$-tensor natural transformations
$\mu^{\mon{f}}$, $\upsilon^{\mon{f}}$
of $\mon{f}$ are given by
$\mu^{\mon{f}}_z=
\mon{f}(\jmath^{-1}_z)
\circ (\mon{f}_{\tensor})_{z,\kappa}
\circ (\tahar{\mon{f}})_{\mon{f}(z)}^{-1}
\colon \mon{f}\mon{f}(z)\to\mon{f}(z)$
and
$\upsilon^{\mon{f}}_z=(\what{\mon{f}})_z\colon z\to \mon{f}(z)$
for $z\in\CK$.

We write $\cat{Aff}(\kos{K})$
for the opposite category of $\cat{Comm}(\kos{K})$,
and by an \emph{affine scheme in $\kos{K}$}
we mean an object in $\cat{Aff}(\kos{K})$.
The affine scheme in $\kos{K}$
corresponding to a commutative monoid $b$ in $\kos{K}$
is denoted by $\Spec(b)$.
The morphism in $\cat{Aff}(\kos{K})$
corresponding to a morphism $f\colon b\to b^{\pr}$ in $\cat{Comm}(\kos{K})$
is denoted by the same symbol $f\colon \Spec(b^{\pr})\to\Spec(b)$;
whenever confusion may arise,
we indicate the category in which the morphism lies,
in words or by specifying the source and target. 
We write
$\kosAff(\kos{K})=(\cat{Aff}(\kos{K}),\times,\Spec(\kappa))$
for the cartesian monoidal category
of affine schemes in $\kos{K}$,
the opposite monoidal category of $\kosComm(\kos{K})$,
with $\Spec(b)\times \Spec(b^{\pr})=\Spec(b\otimes b^{\pr})$.

An \emph{affine monoid scheme in $\kos{K}$}
is a monoid object in $\kosAff(\kos{K})$.
An affine scheme $\Spec(\pi)$ in $\kos{K}$
is an affine monoid scheme in $\kos{K}$
if and only if $\pi$ is a commutative bimonoid in $\kos{K}$.
In this case, we write
$\cp_{\pi}\colon \Spec(\pi)\times\Spec(\pi)\to \Spec(\pi)$
and
$e_{\pi}\colon \Spec(\kappa)\to \Spec(\pi)$
for the product and unit morphisms in $\cat{Aff}(\kos{K})$
of $\Spec(\pi)$;
the corresponding morphisms
$\cp_{\pi}\colon \pi\to \pi\otimes\pi$
and
$e_{\pi}\colon \pi\to \kappa$
in $\cat{Comm}(\kos{K})$
are the coproduct and counit of $\pi$.
We write $\cat{Mon}(\kosAff(\kos{K}))$
for the category of affine monoid schemes in $\kos{K}$,
which is the opposite of the category of comonoid objects in 
$\kosComm(\kos{K})$.

A \emph{lax $\kos{K}$-tensor comonad}
is a comonad internal to $\kos{K}/\mathtt{SMC}_{\mathtt{lax}}$.
We denote by $\delta^{\mon{f}}$ and $\varepsilon^{\mon{f}}$
the coproduct and counit lax $\kos{K}$-tensor natural transformations
of a lax $\kos{K}$-tensor comonad $\mon{f}$.
A lax $\kos{K}$-tensor comonad on $\id_{\kos{K}}$
is called \emph{coreflective}
if it is coreflective as a lax $\kos{K}$-tensor endofunctor
on $\id_{\kos{K}}$.
We write
$\cat{Comonad}[\id_{\kos{K}}]_{\mathrm{crfl}}$
for the category of coreflective lax $\kos{K}$-tensor comonads on $\id_{\kos{K}}$,
which is the category of comonoid objects in
$\kosEnd[\id_{\kos{K}}]_{\mathrm{crfl}}$.
By considering comonoid objects on both sides of
the adjoint equivalence of monoidal categories
in Proposition~\ref{prop-ltcatoversmc-CommKcrfllaxKtmonoidal},
we obtain the following corollary.

\begin{corollary}
  \label{cor-ltcatoversmc-MonAfflaxKtcomonadequiv}
  We have an adjoint equivalence of categories
  \begin{equation*}
    \vcenter{\hbox{
      \xymatrix@C=25pt{
        \cat{Mon}(\kosAff(\kos{K}))^{\op}
        \ar@<1ex>[r]^-{\iota}
        \ar@{}[r]|-{\simeq}
        &\cat{Comonad}[\id_{\kos{K}}]_{\mathrm{crfl}}.
        \ar@<1ex>[l]^-{\mathcal{R}}
      }
    }}
  \end{equation*}
\end{corollary}
\begin{proof}
  This is a direct consequence of
  Proposition~\ref{prop-ltcatoversmc-CommKcrfllaxKtmonoidal},
  since comonoid objects in
  $\kosEnd[\id_{\kos{K}}]_{\mathrm{crfl}}^{\mathrm{rev}}$
  are the same as comonoid objects in $\kosEnd[\id_{\kos{K}}]_{\mathrm{crfl}}$.
  The left adjoint $\iota$
  sends each affine monoid scheme
  $\Spec(\pi)$ in $\kos{K}$
  to the coreflective lax $\kos{K}$-tensor comonad
  $\iota(\Spec(\pi))=\slot\otimes \pi$ on $\id_{\kos{K}}$,
  with $\delta^{\slot\otimes \pi}$, $\varepsilon^{\slot\otimes\pi}$ given by
  $\delta^{\slot\otimes \pi}_z
  =I_z\otimes \cp_{\pi}
  \colon z\otimes \pi
  \to z\otimes \pi\otimes \pi$
  and
  $\varepsilon^{\slot\otimes \pi}_z
  =I_z\otimes e_{\pi}
  \colon z\otimes \pi
  \to z$
  for $z\in \CK$.
  The right adjoint $\mathcal{R}$ sends
  each coreflective lax $\kos{K}$-tensor comonad $\mon{f}$ on $\id_{\kos{K}}$
  to the affine monoid scheme
  $\mathcal{R}(\mon{f})
  =\Spec(\mon{f}(\kappa))$ in $\kos{K}$,
  where the coproduct and counit of $\mon{f}(\kappa)$
  are given by
  $\cp_{\mon{f}(\kappa)}
  =(\tahar{\mon{f}})_{\mon{f}(\kappa)}^{-1}
  \circ \delta^{\mon{f}}_{\kappa}
  \colon \mon{f}(\kappa)
  \to \mon{f}(\kappa)\otimes \mon{f}(\kappa)$
  and
  $e_{\mon{f}(\kappa)}
  =\varepsilon^{\mon{f}}_{\kappa}
  \colon\mon{f}(\kappa)\to \kappa$.
\end{proof}

An \emph{affine group scheme in $\kos{K}$}
is a group object in $\kosAff(\kos{K})$.
An affine monoid scheme $\Spec(\pi)$ in $\kos{K}$
is an affine group scheme in $\kos{K}$
if and only if $\pi$ is a commutative Hopf monoid in $\kos{K}$.
In this case, the inverse morphism
$\varsigma_{\pi}\colon \Spec(\pi)\to \Spec(\pi)$
in $\cat{Aff}(\kos{K})$
corresponds to the antipode
$\varsigma_{\pi}\colon \pi\to \pi$
in $\cat{Comm}(\kos{K})$
and is involutive.
The existence of $\varsigma_{\pi}$
is equivalent to the condition that
$\tau_{\pi}
\coloneqq (I_{\pi}\otimes \pc_{\pi})
\circ (\cp_{\pi}\otimes I_{\pi})
\colon \pi\otimes \pi\to \pi\otimes \pi$
is invertible.
Indeed, $\varsigma_{\pi}$ and $\tau_{\pi}^{-1}$
determine each other by the relations
$\tau_{\pi}^{-1}
=(I_{\pi}\otimes \pc_{\pi})
\circ (I_{\pi}\otimes \varsigma_{\pi}\otimes I_{\pi})
\circ (\cp_{\pi}\otimes I_{\pi})
\colon \pi\otimes \pi\cong\pi\otimes \pi$
and 
$\varsigma_{\pi}
=(e_{\pi}\otimes I_{\pi})
\circ \tau_{\pi}^{-1}
\circ (I_{\pi}\otimes u_{\pi})
\colon \pi\cong \pi$.
We write $\cat{Grp}(\kosAff(\kos{K}))$
for the category of affine group schemes in $\kos{K}$.

Let $\kos{K}\xrightarrow{\monbase{t}}\kos{T}$
be a lax $\kos{K}$-tensor category
and $\mon{f}$ be a lax $\kos{K}$-tensor comonad on $\monbase{t}$,
where $\kos{T}=(\CT,\tensor,\unit)$.
The \emph{fusion operator}
of $\mon{f}$
is the natural transformation
$\Sigma^{\mon{f}}\colon
\mon{f}(\slot)\tensor \mon{f}(\slot)
\Rightarrow \mon{f}(\mon{f}(\slot)\tensor\slot)$
between functors from $\CT\times \CT$ to $\CT$
whose component at $Z$, $X\in \CT$ is given by
\begin{equation*}
  \xymatrix@C=40pt{
    \mon{f}(Z)\tensor \mon{f}(X)
    \ar[r]^-{\delta^{\mon{f}}_Z\tensor I_{\mon{f}(X)}}
    &\mon{f}\mon{f}(Z)\tensor \mon{f}(X)
    \ar[r]^-{(\mon{f}_{\tensor})_{\mon{f}(Z),X}}
    &\mon{f}(\mon{f}(Z)\tensor X)
    .
  }
\end{equation*}
A \emph{lax $\kos{K}$-tensor Hopf comonad}
is a lax $\kos{K}$-tensor comonad
whose fusion operator is invertible.
We denote by
$\cat{HComonad}[\id_{\kos{K}}]_{\mathrm{crfl}}$
the category of coreflective lax $\kos{K}$-tensor Hopf comonads
on $\id_{\kos{K}}$,
which is a full subcategory of
$\cat{Comonad}[\id_{\kos{K}}]_{\mathrm{crfl}}$.
By restricting the adjoint equivalence of categories
in Corollary~\ref{cor-ltcatoversmc-MonAfflaxKtcomonadequiv},
we obtain the following proposition.

\begin{proposition}
  \label{prop-ltcatoversmc-GrpAfflaxKtHcomonadequiv}
  We have an adjoint equivalence of categories
  \begin{equation*}
    \vcenter{\hbox{
      \xymatrix@C=25pt{
        \cat{Grp}(\kosAff(\kos{K}))^{\op}
        \ar@<1ex>[r]^-{\iota}
        \ar@{}[r]|-{\simeq}
        &\cat{HComonad}[\id_{\kos{K}}]_{\mathrm{crfl}}.
        \ar@<1ex>[l]^-{\mathcal{R}}
      }
    }}
  \end{equation*}
\end{proposition}
\begin{proof}
  It suffices to show that the functors
  $\iota$ and $\mathcal{R}$
  in Corollary~\ref{cor-ltcatoversmc-MonAfflaxKtcomonadequiv}
  restrict to the full subcategories
  $\cat{Grp}(\kosAff(\kos{K}))^{\op}$
  and $\cat{HComonad}[\id_{\kos{K}}]_{\mathrm{crfl}}$.

  Let $\Spec(\pi)$ be an affine monoid scheme in $\kos{K}$.
  The component of the fusion operator of
  $\iota(\Spec(\pi))=\slot\otimes \pi$
  at $z$, $x\in \CK$ is
  $\Sigma^{\slot\otimes \pi}_{z,x}
  =(I_{z\otimes \pi\otimes x}\otimes \pc_{\pi})
  \circ (I_{z\otimes \pi}\otimes \sigma_{\pi,x}\otimes I_{\pi})
  \circ (I_z\otimes \cp_{\pi}\otimes I_{x\otimes \pi})
  \colon z\otimes \pi\otimes x\otimes \pi
  \to z\otimes \pi\otimes x\otimes \pi$.
  The fusion operator $\Sigma^{\slot\otimes \pi}$
  is invertible if and only if
  $\Spec(\pi)$ is an affine group scheme in $\kos{K}$.
  Indeed, if $\varsigma_{\pi}$ exists,
  then $\Sigma^{\slot\otimes \pi}_{z,x}$ is invertible
  with inverse
  $(\Sigma^{\slot\otimes \pi}_{z,x})^{-1}
  =(I_{z\otimes \pi\otimes x}\otimes \pc_{\pi})
  \circ (I_{z\otimes \pi}\otimes \sigma_{\pi,x}\otimes I_{\pi})
  \circ (I_{z\otimes \pi}\otimes \varsigma_{\pi}\otimes I_{x\otimes \pi})
  \circ (I_z\otimes \cp_{\pi}\otimes I_{x\otimes \pi})
  \colon z\otimes \pi\otimes x\otimes \pi
  \to z\otimes \pi\otimes x\otimes \pi$.
  Conversely, if $\Sigma^{\slot\otimes \pi}$ is invertible,
  then so is $\tau_{\pi}$,
  as it is equal to
  $\Sigma^{\slot\otimes \pi}_{\kappa,\kappa}$
  up to coherence isomorphisms,
  and hence $\varsigma_{\pi}$ exists.
  In particular, 
  $\iota$ restricts to a functor from
  $\cat{Grp}(\kosAff(\kos{K}))^{\op}$ 
  to $\cat{HComonad}[\id_{\kos{K}}]_{\mathrm{crfl}}$.

  Let $\mon{f}$ be a coreflective lax $\kos{K}$-tensor Hopf comonad on $\id_{\kos{K}}$,
  and put $\pi\coloneqq\mon{f}(\kappa)$.
  Then $\tahar{\mon{f}}\colon \slot\otimes \pi\cong \mon{f}$
  in \eqref{eq-ltcatoversmc-tahardef}
  is an isomorphism of lax $\kos{K}$-tensor comonads,
  as it is the component of the adjunction counit
  of $\iota\dashv\mathcal{R}$ at $\mon{f}$.
  Thus the fusion operator of $\slot\otimes \pi$ is invertible as well,
  and hence $\mathcal{R}(\mon{f})=\Spec(\pi)$
  is an affine group scheme in $\kos{K}$.
  This shows that $\mathcal{R}$ restricts to a functor from
  $\cat{HComonad}[\id_{\kos{K}}]_{\mathrm{crfl}}$
  to $\cat{Grp}(\kosAff(\kos{K}))^{\op}$.
  One can check that the antipode
  $\varsigma_{\pi}=(e_{\pi}\otimes I_{\pi})\circ \tau_{\pi}^{-1}\circ (I_{\pi}\otimes u_{\pi})$
  can also be written as
  $\varsigma_{\mon{f}(\kappa)}
  =(\varepsilon^{\mon{f}}_{\kappa}\otimes I_{\mon{f}(\kappa)})
  \circ (\Sigma^{\mon{f}}_{\kappa,\kappa})^{-1}
  \circ \mon{f}(\jmath_{\mon{f}(\kappa)})
  \circ (\what{\mon{f}})_{\mon{f}(\kappa)}
  \colon \mon{f}(\kappa)
  \to \mon{f}(\kappa)$.
\end{proof}
 

%% file: cltcatoversmc.tex
\chapter{Colax Tensor Categories over a Symmetric Monoidal Category}
\label{appdix-cltcatoversmc}

Throughout this appendix,
let $\kos{S}=(\CS,\otimes,s)$
be a symmetric monoidal category.
We study symmetric monoidal categories
equipped with a colax monoidal functor from $\kos{S}$,
colax monoidal functors between these categories
that preserve the structures over $\kos{S}$
up to a comonoidal natural transformation,
and comonoidal natural transformations between these functors
that are compatible over $\kos{S}$.
Using these notions,
we provide an abstract characterization of 
the endofunctor on $\kos{S}$
given by tensoring with 
a cocommutative comonoid,
a cocommutative bimonoid,
or a cocommutative Hopf monoid in $\kos{S}$.

We denote by 
$\mathtt{SMC}_{\mathtt{colax}}$
the $2$-category of symmetric monoidal categories,
colax monoidal functors,
and comonoidal natural transformations.
We write 
$\kos{S}/\mathtt{SMC}_{\mathtt{colax}}$
for the colax coslice $2$-category of 
$\mathtt{SMC}_{\mathtt{colax}}$
under $\kos{S}$.
\begin{itemize}
  \item 
  A \emph{colax $\kos{S}$-tensor category},
  or a \emph{colax tensor category over $\kos{S}$},
  is a $0$-cell in $\kos{S}/\mathtt{SMC}_{\mathtt{colax}}$.
  It is a pair $(\kos{T},\monbase{t})$
  consisting of a symmetric monoidal category $\kos{T}$
  and a colax monoidal functor $\monbase{t}\colon\kos{S}\to\kos{T}$.
  We write the pair as $\kos{S}\xrightarrow{\monbase{t}}\kos{T}$,
  or simply as $\monbase{t}$.
  An example is the pair consisting of
  $\kos{S}$ and its identity functor $\id_{\kos{S}}$.

  \item
  A \emph{colax $\kos{S}$-tensor functor},
  or a \emph{colax tensor functor over $\kos{S}$},
  is a $1$-cell in $\kos{S}/\mathtt{SMC}_{\mathtt{colax}}$.
  Let $\kos{S}\xrightarrow{\monbase{t}}\kos{T}$,
  $\kos{S}\xrightarrow{\monbase{t}^{\pr}}\kos{T}^{\pr}$
  be colax $\kos{S}$-tensor categories.
  A colax $\kos{S}$-tensor functor
  $(\mon{f},\what{\mon{f}})\colon\monbase{t}^{\pr}\to\monbase{t}$
  is a pair consisting of
  a colax monoidal functor
  $\mon{f}\colon \kos{T}^{\pr}\to \kos{T}$
  and a comonoidal natural transformation
  $\what{\mon{f}}
  \colon \mon{f}\monbase{t}^{\pr}
  \Rightarrow \monbase{t}$.
  We omit $\what{\mon{f}}$ from the notation
  and write $\mon{f}\colon\monbase{t}^{\pr}\to \monbase{t}$.

  \item
  A \emph{colax $\kos{S}$-tensor natural transformation},
  or a \emph{colax tensor natural transformation over $\kos{S}$},
  is a $2$-cell in $\kos{S}/\mathtt{SMC}_{\mathtt{colax}}$.
  Let $\mon{f}$, $\mon{g}\colon \monbase{t}^{\pr}\to \monbase{t}$
  be colax $\kos{S}$-tensor functors.
  A colax $\kos{S}$-tensor natural transformation
  $\vartheta\colon \mon{f}\Rightarrow\mon{g}$
  is a comonoidal natural transformation
  between the underlying colax monoidal functors
  $\mon{f}$, $\mon{g}\colon \kos{T}^{\pr}\to \kos{T}$
  that satisfies
  $\what{\mon{f}}
  =\what{\mon{g}}
  \circ (\vartheta\monbase{t}^{\pr})
  \colon \mon{f}\monbase{t}^{\pr}
  \Rightarrow \monbase{t}$.
  A \emph{colax $\kos{S}$-tensor natural isomorphism},
  or a \emph{colax tensor natural isomorphism over $\kos{S}$},
  is an invertible $2$-cell in $\kos{S}/\mathtt{SMC}_{\mathtt{colax}}$.
\end{itemize}
Let $\mon{f}$, $\mon{g}\colon \monbase{t}^{\pr}\to \monbase{t}$
be colax $\kos{S}$-tensor functors.
We write
$(\kos{S}/\mathtt{SMC}_{\mathtt{colax}})[\monbase{t}^{\pr},\monbase{t}]$
for the category of colax $\kos{S}$-tensor functors
from $\monbase{t}^{\pr}$ to $\monbase{t}$.
We illustrate
$\mon{f}$, $\mon{g}$
and a colax $\kos{S}$-tensor natural transformation
$\vartheta\colon \mon{f}\Rightarrow\mon{g}$
as follows.
\begin{equation*}
  \begin{tikzcd}[row sep=5pt,column sep=5pt]
    \kos{T}^{\pr}
    \ar[rr,"\mon{f}"]
    &\text{ }
    &\kos{T}
    &\text{ }
    &\kos{T}^{\pr}
    \ar[rr,"\mon{g}"]
    &\text{ }
    &\kos{T}
    &\text{ }
    &\kos{T}^{\pr}
    \ar[rr,shift left=0.5ex,bend left=22pt,"\mon{f}" description]
    \ar[rr,shift right=0.5ex,bend right=22pt,"\mon{g}" description]
    &\text{ }
    \ar[d,2cell,pos=-1.5,shift right=0ex,"{\,\vartheta}"]
    &\kos{T}
    \\
    \text{ }
    \ar[rr,2cell,shift right=0.5ex,"{\what{\mon{f}}}"]
    &\text{ }
    &\text{ }
    &\text{ }
    &\text{ }
    \ar[rr,2cell,shift right=0.5ex,"{\what{\mon{g}}}"]
    &\text{ }
    &\text{ }
    &\text{ }
    &\text{ }
    &\text{ }
    &\text{ }
    \\
    \text{ }
    &\kos{S}
    \ar[uul,"\monbase{t}^{\pr}"]
    \ar[uur,"\monbase{t}"']
    &\text{ }
    &\text{ }
    &\text{ }
    &\kos{S}
    \ar[uul,"\monbase{t}^{\pr}"]
    \ar[uur,"\monbase{t}"']
    &\text{ }
    &\text{ }
    &\text{ }
    &\kos{S}
    \ar[uul,"\monbase{t}^{\pr}"]
    \ar[uur,"\monbase{t}"']
    &\text{ }
  \end{tikzcd}
\end{equation*}

Let $\kos{S}\xrightarrow{\monbase{t}}\kos{T}$
be a colax $\kos{S}$-tensor category,
where $\kos{T}=(\CT,\tensor,\unit)$.
The \emph{induced $\kos{S}$-action} $\bcts$ on $\kos{T}$
is the functor $\bcts\colon \CT\times\CS\to\CT$
sending each pair $(X,z)$ to $X\bcts z\coloneqq X\tensor \monbase{t}(z)$.
Let $\kos{S}\xrightarrow{\monbase{t}^{\pr}}\kos{T}^{\pr}$
be a colax $\kos{S}$-tensor category,
where $\kos{T}^{\pr}=(\CT^{\pr},\tensor^{\pr},\unit^{\pr})$
and $\bcts^{\pr}$ is the induced $\kos{S}$-action on $\kos{T}^{\pr}$.
Let $\mon{f}\colon \monbase{t}^{\pr}\to \monbase{t}$
be a colax $\kos{S}$-tensor functor.
The \emph{$\kos{S}$-equivariance}
of $\mon{f}$ is the natural transformation
\begin{equation*}
  \vecar{\mon{f}}\colon
  \mon{f}(\slot \bcts^{\pr}\slot)
  \Rightarrow \mon{f}(\slot)\bcts \slot
\end{equation*}
between functors from $\CT^{\pr}\times \CS$ to $\CT$
whose component at $X^{\pr}\in \CT^{\pr}$, $z\in \CS$ is
\begin{equation*}
  \xymatrix@C=50pt{
    \mon{f}(X^{\pr}\tensor^{\pr} \monbase{t}^{\pr}(z))
    \ar[r]^-{(\mon{f}_{\tensor})_{X^{\pr},\monbase{t}^{\pr}(z)}}
    &\mon{f}(X^{\pr})\tensor \mon{f}\monbase{t}^{\pr}(z)
    \ar[r]^-{I_{\mon{f}(X^{\pr})}\tensor (\what{\mon{f}})_z}
    &\mon{f}(X^{\pr})\tensor \monbase{t}(z)
    .
  }
\end{equation*}
Let $\kos{S}\xrightarrow{\monbase{t}^{\ppr}} \kos{T}^{\ppr}$
be a colax $\kos{S}$-tensor category,
where $\kos{T}^{\ppr}=(\CT^{\ppr},\tensor^{\ppr},\unit^{\ppr})$
and $\bcts^{\ppr}$ is the induced $\kos{S}$-action on $\kos{T}^{\ppr}$.
The composite 
$\mon{f}\mon{f}^{\pr}\colon\monbase{t}^{\ppr}\to \monbase{t}$
of colax $\kos{S}$-tensor functors
$\mon{f}\colon\monbase{t}^{\pr}\to\monbase{t}$,
$\mon{f}^{\pr}\colon\monbase{t}^{\ppr}\to \monbase{t}^{\pr}$
consists of
$\mon{f}\mon{f}^{\pr}\colon\kos{T}^{\ppr}\to\kos{T}$
and
$\what{\mon{f}\mon{f}^{\pr}}
=\what{\mon{f}}
\circ (\mon{f}\what{\mon{f}^{\pr}})
\colon \mon{f}\mon{f}^{\pr}\monbase{t}^{\ppr}
\Rightarrow \monbase{t}$.
The $\kos{S}$-equivariance of the composite is given by
$(\vecar{\mon{f}\mon{f}^{\pr}})_{X^{\ppr},z}
=(\vecar{\mon{f}})_{\mon{f}^{\pr}(X^{\ppr}),z}
\circ \mon{f}((\vecar{\mon{f}^{\pr}})_{X^{\ppr},z})
\colon \mon{f}\mon{f}^{\pr}(X^{\ppr}\bcts^{\ppr} z)
\to \mon{f}\mon{f}^{\pr}(X^{\ppr})\bcts z$
for $X^{\ppr}\in\CT^{\ppr}$, $z\in \CS$.
The identity colax $\kos{S}$-tensor functor of $\monbase{t}$,
denoted by $\id_{\kos{T}}$,
is the pair consisting of the identity functor of $\kos{T}$
and the identity natural transformation of $\monbase{t}$.

\begin{lemma}
  \label{lem-cltcatoversmc-colaxKtnatequiv}
  Let
  $\kos{S}\xrightarrow{\monbase{t}} \kos{T}$,
  $\kos{S}\xrightarrow{\monbase{t}^{\pr}} \kos{T}^{\pr}$
  be colax $\kos{S}$-tensor categories,
  where
  $\kos{T}=(\CT,\tensor,\unit)$,
  $\kos{T}^{\pr}=(\CT^{\pr},\tensor^{\pr},\unit^{\pr})$.
  Let $\bcts$ and $\bcts^{\pr}$
  be the induced $\kos{S}$-actions
  on $\kos{T}$ and $\kos{T}^{\pr}$, respectively.
  Let 
  $\mon{f}$, $\mon{g}\colon \monbase{t}^{\pr}\to \monbase{t}$
  be colax $\kos{S}$-tensor functors
  and $\vartheta\colon\mon{f}\Rightarrow\mon{g}$
  be a comonoidal natural transformation
  between the underlying colax monoidal functors
  $\mon{f}$, $\mon{g}\colon\kos{T}^{\pr}\to \kos{T}$.
  Then $\vartheta\colon\mon{f}\Rightarrow\mon{g}$
  is a colax $\kos{S}$-tensor natural transformation
  if and only if
  the following diagram commutes for every pair
  $X^{\pr}\in \CT^{\pr}$, $z\in \CS$.
  \begin{equation}
    \label{eq-cltcatoversmc-colaxKtnatequiv}
    \vcenter{\hbox{
      \xymatrix@R=15pt@C=35pt{
        \mon{f}(X^{\pr}\bcts^{\pr} z)
        \ar[r]^-{(\vecar{\mon{f}})_{X^{\pr},z}}
        \ar[d]_-{\vartheta_{X^{\pr}\bcts^{\pr} z}}
        &\mon{f}(X^{\pr})\bcts z
        \ar[d]^-{\vartheta_{X^{\pr}}\bcts I_z}
        \\
        \mon{g}(X^{\pr}\bcts^{\pr} z)
        \ar[r]^-{(\vecar{\mon{g}})_{X^{\pr},z}}
        &\mon{g}(X^{\pr})\bcts z
      }
    }}
  \end{equation}
\end{lemma}
\begin{proof}
  If $\vartheta$ is a colax $\kos{S}$-tensor natural transformation,
  the diagram \eqref{eq-cltcatoversmc-colaxKtnatequiv}
  commutes by the definitions.
  Conversely, assume that
  \eqref{eq-cltcatoversmc-colaxKtnatequiv} commutes
  for every pair $X^{\pr}$, $z$.
  One can check that
  $(\what{\mon{f}})_z
  \colon \mon{f}\monbase{t}^{\pr}(z)
  \to \monbase{t}(z)$
  is equal to
  $(\mon{f}_{\unit}\tensor I_{\monbase{t}(z)})
  \circ (\vecar{\mon{f}})_{\unit^{\pr},z}
  \colon \mon{f}(\unit^{\pr}\tensor^{\pr} \monbase{t}^{\pr}(z))
  \to \unit\tensor \monbase{t}(z)$
  up to coherence isomorphisms.
  Using this description for both $\what{\mon{f}}$ and $\what{\mon{g}}$
  and applying \eqref{eq-cltcatoversmc-colaxKtnatequiv} at $X^{\pr}=\unit^{\pr}$,
  we deduce that
  $(\what{\mon{f}})_z=
  (\what{\mon{g}})_z
  \circ \vartheta_{\monbase{t}^{\pr}(z)}
  \colon \mon{f}\monbase{t}^{\pr}(z)
  \to \monbase{t}(z)$
  for every $z$.
  Thus $\vartheta$ is a colax $\kos{S}$-tensor natural transformation.
\end{proof}

Let $\kos{S}\xrightarrow{\monbase{t}}\kos{T}$
be a colax $\kos{S}$-tensor category,
where $\kos{T}=(\CT,\tensor,\unit)$.
We denote by
$\kosEns(\kos{T})=(\cat{Ens}(\kos{T}),\tensor,\unit)$
the cartesian monoidal category
of cocommutative comonoids in $\kos{T}$.
The coproduct and counit morphisms of
a cocommutative comonoid $C$ in $\kos{T}$
are denoted by $\cp_C\colon C\to C\tensor C$ and $e_C\colon C\to \unit$.
Each cocommutative comonoid $C$ in $\kos{T}$
determines a colax $\kos{S}$-tensor functor
$C\bcts\slot\colon \id_{\kos{S}}\to \monbase{t}$
consisting of the colax monoidal functor
$C\bcts \slot=C\tensor \monbase{t}(\slot)
\colon \kos{S}\to \kos{T}$
and the comonoidal natural transformation
$\what{C\bcts}
\colon C\bcts \slot
\Rightarrow \monbase{t}$
given by
$(\what{C\bcts})_z
\coloneqq e_C\tensor I_{\monbase{t}(z)}
\colon C\tensor \monbase{t}(z)
\to \monbase{t}(z)$.
The $\kos{S}$-equivariance of $C\bcts \slot$ is
$(\vecar{C\bcts})_{x,z}
=I_C\tensor (\monbase{t}_{\tensor})_{x,z}
\colon C\tensor \monbase{t}(x\otimes z)
\to C\tensor \monbase{t}(x)\tensor \monbase{t}(z)$
for $x$, $z\in \CS$.
Each morphism $f\colon C^{\pr}\to C$ in $\cat{Ens}(\kos{T})$
determines a colax $\kos{S}$-tensor natural transformation
$f\bcts \slot\colon C^{\pr}\bcts \slot\Rightarrow C\bcts \slot$.

\begin{proposition}
  \label{prop-cltcatoversmc-EnscolaxKtadj}
  Let $\kos{S}\xrightarrow{\monbase{t}}\kos{T}$
  be a colax $\kos{S}$-tensor category.
  We have an adjunction
  \begin{equation*}
    \vcenter{\hbox{
      \xymatrix@C=25pt{
        (\kos{S}/\mathtt{SMC}_{\mathtt{colax}})[\id_{\kos{S}},\monbase{t}]
        \ar@<1ex>[r]^-{\mathcal{L}}
        \ar@{}[r]|-{\bot}
        &\cat{Ens}(\kos{T})
        \ar@<1ex>[l]^-{\iota}
      }
    }}
  \end{equation*}
  between
  the category of cocommutative comonoids in $\kos{T}$
  and the category of colax $\kos{S}$-tensor functors
  from $\id_{\kos{S}}$ to $\monbase{t}$.
\end{proposition}
\begin{proof}
  The functor $\iota$ sends
  each cocommutative comonoid $C$ in $\kos{T}$ to $C\bcts \slot$,
  and each morphism $f\colon C^{\pr}\to C$
  in $\cat{Ens}(\kos{T})$ to
  $f\bcts \slot\colon C^{\pr}\bcts \slot\Rightarrow C\bcts \slot$.
  The functor $\mathcal{L}$ sends each colax $\kos{S}$-tensor functor
  $\mon{f}\colon\id_{\kos{S}}\to \monbase{t}$
  to the cocommutative comonoid $\mon{f}(s)$ in $\kos{T}$,
  and each morphism
  $\vartheta\colon \mon{f}\Rightarrow\mon{g}$
  in $(\kos{S}/\mathtt{SMC}_{\mathtt{colax}})[\id_{\kos{S}},\monbase{t}]$
  to the morphism
  $\vartheta_{s}\colon \mon{f}(s)\to \mon{g}(s)$
  in $\cat{Ens}(\kos{T})$.
  Here $\mon{f}(s)$ has coproduct
  $\cp_{\mon{f}(s)}
  =(\mon{f}_{\tensor})_{s,s}
  \circ \mon{f}(\imath_{s})$
  and counit $e_{\mon{f}(s)}=\mon{f}_{\unit}$.
  The component of the adjunction counit at each
  $C\in \cat{Ens}(\kos{T})$ is
  $I_C\tensor \monbase{t}_\unit
  \colon C\tensor \monbase{t}(s)
  \to C$.
  The component of the adjunction unit at each
  colax $\kos{S}$-tensor functor
  $\mon{f}\colon \id_{\kos{S}}\to \monbase{t}$
  is the colax $\kos{S}$-tensor natural transformation
  \begin{equation}
    \label{eq-cltcatoversmc-hatardef}
    \hatar{\mon{f}}\colon
    \mon{f}
    \Rightarrow \mon{f}(s)\bcts \slot
  \end{equation}
  where 
  $(\hatar{\mon{f}})_z
  \coloneqq
  (\vecar{\mon{f}})_{s,z}
  \circ \mon{f}(\imath_z)
  \colon \mon{f}(z)
  \to \mon{f}(s)\tensor \monbase{t}(z)$
  for $z\in \CS$.
  We have
  $(\hatar{C\bcts})_z=
  I_C\tensor ((\monbase{t}_{\tensor})_{s,z}\circ \monbase{t}(\imath_z))
  \colon C\bcts z
  \to (C\bcts s)\bcts z$
  for $C\in\cat{Ens}(\kos{T})$, $z\in \CS$.
  It is straightforward to check that these data
  constitute an adjunction $\mathcal{L}\dashv \iota$.
\end{proof}

A \emph{strong $\kos{S}$-tensor category},
or a \emph{strong tensor category over $\kos{S}$},
is a pair
$(\kos{T},\monbase{t})$
consisting of a symmetric monoidal category $\kos{T}$
and a strong monoidal functor $\monbase{t}$.
It is in particular a colax $\kos{S}$-tensor category,
and we write it as 
$\kos{S}\xrightarrow{\monbase{t}}\kos{T}$,
or simply as $\monbase{t}$.
For instance, $\id_{\kos{S}}$ is a strong $\kos{S}$-tensor category.

\begin{definition}
  \label{def-cltcatoversmc-reflective}
  Let $\kos{S}\xrightarrow{\monbase{t}}\kos{T}$
  be a strong $\kos{S}$-tensor category
  and $\mon{f}\colon \id_{\kos{S}}\to \monbase{t}$
  be a colax $\kos{S}$-tensor functor.
  \begin{enumerate}
    \item 
    We refer to 
    $\hatar{\mon{f}}\colon \mon{f}\Rightarrow\mon{f}(s)\bcts \slot$
    in \eqref{eq-cltcatoversmc-hatardef}
    as the \emph{reflection} of $\mon{f}$.

    \item
    We say that $\mon{f}$ is \emph{reflective}
    if its reflection $\hatar{\mon{f}}$ is invertible.
  \end{enumerate}
  We write
  $(\kos{S}/\mathtt{SMC}_{\mathtt{colax}})[\id_{\kos{S}},\monbase{t}]_{\mathrm{rfl}}$
  for the full subcategory of
  $(\kos{S}/\mathtt{SMC}_{\mathtt{colax}})[\id_{\kos{S}},\monbase{t}]$
  consisting of reflective
  colax $\kos{S}$-tensor functors
  from $\id_{\kos{S}}$ to $\monbase{t}$.
\end{definition}

\begin{proposition}
  \label{prop-cltcatoversmc-EnsstrongKtequiv}
  Let $\kos{S}\xrightarrow{\monbase{t}}\kos{T}$
  be a strong $\kos{S}$-tensor category.
  Recall the adjunction
  $\mathcal{L}\dashv \iota$
  in Proposition~\ref{prop-cltcatoversmc-EnscolaxKtadj}.
  \begin{enumerate}
    \item 
    The right adjoint
    $\iota\colon
    \cat{Ens}(\kos{T})
    \to (\kos{S}/\mathtt{SMC}_{\mathtt{colax}})[\id_{\kos{S}},\monbase{t}]$
    is fully faithful.

    \item 
    The essential image of the right adjoint
    $\iota$ is 
    $(\kos{S}/\mathtt{SMC}_{\mathtt{colax}})[\id_{\kos{S}},\monbase{t}]_{\mathrm{rfl}}$,
    which is a reflective full subcategory
    of $(\kos{S}/\mathtt{SMC}_{\mathtt{colax}})[\id_{\kos{S}},\monbase{t}]$.

    \item 
    The adjunction $\mathcal{L}\dashv \iota$
    restricts to an adjoint equivalence of categories
    \begin{equation*}
      \vcenter{\hbox{
        \xymatrix@C=25pt{
          (\kos{S}/\mathtt{SMC}_{\mathtt{colax}})[\id_{\kos{S}},\monbase{t}]_{\mathrm{rfl}}
          \ar@<1ex>[r]^-{\mathcal{L}}
          \ar@{}[r]|-{\simeq}
          &\cat{Ens}(\kos{T}).
          \ar@<1ex>[l]^-{\iota}
        }
      }}
    \end{equation*}
  \end{enumerate}
\end{proposition}
\begin{proof}
  Since $\monbase{t}_{\unit}\colon \monbase{t}(s)\to \unit$
  is invertible,
  the adjunction counit of $\mathcal{L}\dashv \iota$ is a natural isomorphism.
  Thus $\iota$ is fully faithful.
  This proves (1).
  A colax $\kos{S}$-tensor functor
  $\mon{f}\colon \id_{\kos{S}}\to \monbase{t}$
  lies in the essential image of $\iota$
  if and only if the component $\hatar{\mon{f}}$
  of the adjunction unit at $\mon{f}$
  is invertible.
  Hence the essential image of $\iota$ is the full subcategory
  of reflective colax $\kos{S}$-tensor functors,
  and its inclusion functor has a left adjoint $\iota\mathcal{L}$.
  This proves (2).
  Then (3) follows,
  as the adjunction unit becomes invertible
  after restriction.
\end{proof}

\begin{lemma}
  \label{lem-cltcatoversmc-reflectionKequiv}
  Let $\kos{S}\xrightarrow{\monbase{t}}\kos{T}$
  be a strong $\kos{S}$-tensor category
  and $\mon{f}\colon \id_{\kos{S}}\to\monbase{t}$
  be a colax $\kos{S}$-tensor functor.
  Then the component of its reflection
  $\hatar{\mon{f}}$ at $s$ is
  $(\hatar{\mon{f}})_{s}
  =I_{\mon{f}(s)}\tensor (\monbase{t}_{\unit})^{-1}
  \colon \mon{f}(s)
  \cong \mon{f}(s)\tensor \monbase{t}(s)$.
  Moreover, the following are equivalent:
  \begin{enumerate}
    \item 
    $\mon{f}$ is reflective.

    \item
    $\mon{f}$ is isomorphic to
    $C\bcts \slot\colon \id_{\kos{S}}\to \monbase{t}$
    as colax $\kos{S}$-tensor functors
    for some cocommutative comonoid $C$ in $\kos{T}$.

    \item
    The $\kos{S}$-equivariance
    $\vecar{\mon{f}}
    \colon \mon{f}(\slot\otimes \slot)
    \Rightarrow \mon{f}(\slot)\bcts \slot$
    of $\mon{f}$ is invertible.
  \end{enumerate}
\end{lemma}
\begin{proof}
  We have
  $(I_{\mon{f}(s)}\tensor \monbase{t}_{\unit})
  \circ (\hatar{\mon{f}})_{s}
  =I_{\mon{f}(s)}$,
  as it is one of the triangle identities for the adjunction
  $\mathcal{L}\dashv\iota$
  in Proposition~\ref{prop-cltcatoversmc-EnscolaxKtadj}.
  Since $\monbase{t}_{\unit}$ is invertible,
  we deduce the desired description of $(\hatar{\mon{f}})_{s}$.
  The equivalence of (1) and (2) follows from
  Proposition~\ref{prop-cltcatoversmc-EnsstrongKtequiv}.
  The definition of $\hatar{\mon{f}}$ in \eqref{eq-cltcatoversmc-hatardef}
  shows that (3) implies (1).
  One can check that
  $(I_{\mon{f}(s)}\tensor (\monbase{t}_{\tensor})_{x,z})
  \circ (\hatar{\mon{f}})_{x\otimes z}
  =((\hatar{\mon{f}})_x\tensor I_{\monbase{t}(z)})
  \circ (\vecar{\mon{f}})_{x,z}
  \colon \mon{f}(x\otimes z)
  \to \mon{f}(s)\tensor \monbase{t}(x)\tensor \monbase{t}(z)$
  holds for $x$, $z\in\CS$.
  Since $(\monbase{t}_{\tensor})_{x,z}$ is invertible,
  we deduce that (1) implies (3).
\end{proof}

By taking $\monbase{t}=\id_{\kos{S}}$
in Proposition~\ref{prop-cltcatoversmc-EnsstrongKtequiv},
we obtain an adjoint equivalence of categories
\begin{equation}
  \label{eq-cltcatoversmc-EnsSrflcolaxKtfunct}
  \vcenter{\hbox{
    \xymatrix@C=25pt{
      (\kos{S}/\mathtt{SMC}_{\mathtt{colax}})[\id_{\kos{S}},\id_{\kos{S}}]_{\mathrm{rfl}}
      \ar@<1ex>[r]^-{\mathcal{L}}
      \ar@{}[r]|-{\simeq}
      &\cat{Ens}(\kos{S}).
      \ar@<1ex>[l]^-{\iota}
    }
  }}
\end{equation}
We now show that the adjoint equivalence
\eqref{eq-cltcatoversmc-EnsSrflcolaxKtfunct}
lifts to an adjoint equivalence of monoidal categories.
Let $\mon{f}$, $\mon{g}\colon \id_{\kos{S}}\to \id_{\kos{S}}$
be reflective colax $\kos{S}$-tensor functors.
By Lemma~\ref{lem-cltcatoversmc-reflectionKequiv},
the $\kos{S}$-equivariances 
$\vecar{\mon{f}}$, $\vecar{\mon{g}}$
of $\mon{f}$, $\mon{g}$ are invertible.
Hence the $\kos{S}$-equivariance $\vecar{\mon{f}\mon{g}}$ of the composite
$\mon{f}\mon{g}\colon \id_{\kos{S}}\to \id_{\kos{S}}$
is invertible,
and thus the composite $\mon{f}\mon{g}$ is also reflective.
We write
\begin{equation*}
  \kosEnd[\id_{\kos{S}}]_{\mathrm{rfl}}
\end{equation*}
for the monoidal category of
reflective colax $\kos{S}$-tensor endofunctors
on $\id_{\kos{S}}$,
whose monoidal product is given by composition
and whose unit object is the identity colax $\kos{S}$-tensor functor of $\id_{\kos{S}}$.

\begin{proposition}
  \label{prop-cltcatoversmc-EnsSrflcolaxKtmonoidal}
  The adjoint equivalence of categories
  \eqref{eq-cltcatoversmc-EnsSrflcolaxKtfunct}
  lifts to an adjoint equivalence
  of monoidal categories
  \begin{equation*}
    \vcenter{\hbox{
      \xymatrix@C=25pt{
        \kosEnd[\id_{\kos{S}}]_{\mathrm{rfl}}
        \ar@<1ex>[r]^-{\mathcal{L}}
        \ar@{}[r]|-{\simeq}
        &\kosEns(\kos{S}).
        \ar@<1ex>[l]^-{\iota}
      }
    }}
  \end{equation*}
  As a consequence, the monoidal category
  $\kosEnd[\id_{\kos{S}}]_{\mathrm{rfl}}$
  is cartesian.
\end{proposition}
\begin{proof}
  Let $c$, $c^{\pr}$ be cocommutative comonoids in $\kos{S}$.
  It follows from the Mac Lane coherence theorem
  that the coherence isomorphisms
  $a_{c^{\pr},c,\slot}^{-1}
  \colon \iota(c^{\pr}\otimes c)
  \cong \iota(c^{\pr})\iota(c)$
  and
  $\imath^{-1}
  \colon \iota(s)
  \cong \id_{\kos{S}}$
  are colax $\kos{S}$-tensor natural isomorphisms.
  With these as its monoidal coherence isomorphisms,
  the right adjoint $\iota$ is a strong monoidal functor
  from $\kosEns(\kos{S})$ to $\kosEnd[\id_{\kos{S}}]_{\mathrm{rfl}}$.
  By doctrinal adjunction
  and the fact that $\mathcal{L}\dashv \iota$ is an adjoint equivalence,
  the left adjoint $\mathcal{L}$
  has a unique strong monoidal functor structure
  such that $\mathcal{L}\dashv \iota$
  becomes an adjoint equivalence of monoidal categories.
  Let $\mon{f}$, $\mon{g}$ 
  be reflective colax $\kos{S}$-tensor endofunctors on $\id_{\kos{S}}$.
  One can check that 
  the monoidal coherence isomorphisms of $\mathcal{L}$
  are the isomorphisms
  $(\hatar{\mon{g}})_{\mon{f}(s)}
  \colon
  \mon{g}\mon{f}(s)
  \cong \mon{g}(s)\otimes \mon{f}(s)$
  and
  $I_{s}\colon s\cong s$
  in $\cat{Ens}(\kos{S})$.
  Since $\kosEns(\kos{S})$ is cartesian,
  so is $\kosEnd[\id_{\kos{S}}]_{\mathrm{rfl}}$.
\end{proof}

A \emph{colax $\kos{S}$-tensor comonad}
is a comonad internal to $\kos{S}/\mathtt{SMC}_{\mathtt{colax}}$.
Since the monoidal category
$\kosEnd[\id_{\kos{S}}]_{\mathrm{rfl}}$
is cartesian,
(i) every reflective colax $\kos{S}$-tensor endofunctor on $\id_{\kos{S}}$
has a unique structure of a colax $\kos{S}$-tensor comonad on $\id_{\kos{S}}$;
and
(ii)
every colax $\kos{S}$-tensor natural transformation
between reflective colax $\kos{S}$-tensor endofunctors
on $\id_{\kos{S}}$
is a morphism of colax $\kos{S}$-tensor comonads on $\id_{\kos{S}}$.
Let $\mon{f}$ be a reflective colax $\kos{S}$-tensor endofunctor on $\id_{\kos{S}}$.
One can check that
the coproduct and counit colax $\kos{S}$-tensor natural transformations
$\delta^{\mon{f}}$, $\varepsilon^{\mon{f}}$
of $\mon{f}$ are given by
$\delta^{\mon{f}}_z=
(\hatar{\mon{f}})_{\mon{f}(z)}^{-1}
\circ (\mon{f}_{\tensor})_{s,z}
\circ \mon{f}(\imath_z)
\colon \mon{f}(z)\to \mon{f}\mon{f}(z)$
and
$\varepsilon^{\mon{f}}_z
=(\what{\mon{f}})_z\colon \mon{f}(z)\to z$
for $z\in\CS$.

A cocommutative comonoid $\pi$ in $\kos{S}$
is a monoid object in $\kosEns(\kos{S})$
if and only if $\pi$ is a cocommutative bimonoid in $\kos{S}$.
In this case, we write
$\pc_{\pi}\colon \pi\otimes \pi \to \pi$
and
$u_{\pi}\colon s\to \pi$
for the product and unit morphisms in $\cat{Ens}(\kos{S})$
of $\pi$.
We write $\cat{Mon}(\kosEns(\kos{S}))$
for the category of monoid objects in $\kosEns(\kos{S})$,
which is the category of cocommutative bimonoids in $\kos{S}$.

A \emph{colax $\kos{S}$-tensor monad}
is a monad internal to $\kos{S}/\mathtt{SMC}_{\mathtt{colax}}$.
The underlying colax monoidal monad is
called a Hopf monad in \cite{Moerdijk2002}.
We denote by $\mu^{\mon{f}}$ and $\upsilon^{\mon{f}}$
the product and unit colax $\kos{S}$-tensor natural transformations
of a colax $\kos{S}$-tensor monad $\mon{f}$.
A colax $\kos{S}$-tensor monad on $\id_{\kos{S}}$
is called \emph{reflective}
if it is reflective as a colax $\kos{S}$-tensor endofunctor
on $\id_{\kos{S}}$.
We write
$\cat{Monad}[\id_{\kos{S}}]_{\mathrm{rfl}}$
for the category of reflective colax $\kos{S}$-tensor monads on $\id_{\kos{S}}$,
which is the category of monoid objects in
$\kosEnd[\id_{\kos{S}}]_{\mathrm{rfl}}$.
By considering monoid objects on both sides of
the adjoint equivalence of monoidal categories
in Proposition~\ref{prop-cltcatoversmc-EnsSrflcolaxKtmonoidal},
we obtain the following corollary.

\begin{corollary}
  \label{cor-cltcatoversmc-MonEnscolaxKtmonadequiv}
  We have an adjoint equivalence of categories
  \begin{equation*}
    \vcenter{\hbox{
      \xymatrix@C=25pt{
        \cat{Monad}[\id_{\kos{S}}]_{\mathrm{rfl}}
        \ar@<1ex>[r]^-{\mathcal{L}}
        \ar@{}[r]|-{\simeq}
        &\cat{Mon}(\kosEns(\kos{S})).
        \ar@<1ex>[l]^-{\iota}
      }
    }}
  \end{equation*}
\end{corollary}
\begin{proof}
  This is a direct consequence of
  Proposition~\ref{prop-cltcatoversmc-EnsSrflcolaxKtmonoidal}.
  The right adjoint $\iota$
  sends each cocommutative bimonoid $\pi$ in $\kos{S}$
  to the reflective colax $\kos{S}$-tensor monad
  $\iota(\pi)=\pi\otimes \slot$ on $\id_{\kos{S}}$,
  with $\mu^{\pi\otimes \slot}$, $\upsilon^{\pi\otimes \slot}$ given by
  $\mu^{\pi\otimes \slot}_z
  =\pc_{\pi}\otimes I_z
  \colon \pi\otimes \pi\otimes z
  \to \pi\otimes z$
  and
  $\upsilon^{\pi\otimes \slot}_z
  =u_{\pi}\otimes I_z
  \colon z
  \to \pi\otimes z$
  for $z\in \CS$.
  The left adjoint $\mathcal{L}$ sends
  each reflective colax $\kos{S}$-tensor monad $\mon{f}$ on $\id_{\kos{S}}$
  to the cocommutative bimonoid
  $\mathcal{L}(\mon{f})
  =\mon{f}(s)$ in $\kos{S}$,
  where the product and unit of $\mon{f}(s)$
  are given by
  $\pc_{\mon{f}(s)}
  =\mu^{\mon{f}}_{s}
  \circ (\hatar{\mon{f}})_{\mon{f}(s)}^{-1}
  \colon \mon{f}(s)\otimes \mon{f}(s)
  \to \mon{f}(s)$
  and
  $u_{\mon{f}(s)}
  =\upsilon^{\mon{f}}_{s}
  \colon s\to \mon{f}(s)$.
\end{proof}

A cocommutative bimonoid $\pi$ in $\kos{S}$
is a group object in $\kosEns(\kos{S})$
if and only if it is a cocommutative Hopf monoid in $\kos{S}$.
In this case, the inverse morphism
$\varsigma_{\pi}\colon \pi\to \pi$
in $\cat{Ens}(\kos{S})$
is the antipode of $\pi$ and is involutive.
The existence of $\varsigma_{\pi}$
is equivalent to the condition that
$\tau_{\pi}
\coloneqq (I_{\pi}\otimes \pc_{\pi})
\circ (\cp_{\pi}\otimes I_{\pi})
\colon \pi\otimes \pi\to \pi\otimes \pi$
is invertible.
Indeed, $\varsigma_{\pi}$ and $\tau_{\pi}^{-1}$
determine each other by the relations
$\tau_{\pi}^{-1}
=(I_{\pi}\otimes \pc_{\pi})
\circ (I_{\pi}\otimes \varsigma_{\pi}\otimes I_{\pi})
\circ (\cp_{\pi}\otimes I_{\pi})
\colon \pi\otimes \pi\cong\pi\otimes \pi$
and 
$\varsigma_{\pi}
=(e_{\pi}\otimes I_{\pi})
\circ \tau_{\pi}^{-1}
\circ (I_{\pi}\otimes u_{\pi})
\colon \pi\cong \pi$.
We write $\cat{Grp}(\kosEns(\kos{S}))$
for the category of group objects in $\kosEns(\kos{S})$,
which is the category of cocommutative Hopf monoids in $\kos{S}$.

Let $\kos{S}\xrightarrow{\monbase{t}}\kos{T}$
be a colax $\kos{S}$-tensor category
and $\mon{f}$ be a colax $\kos{S}$-tensor monad on $\monbase{t}$,
where $\kos{T}=(\CT,\tensor,\unit)$.
The \emph{fusion operator}
of $\mon{f}$
is the natural transformation
$\Sigma^{\mon{f}}\colon
\mon{f}(\slot\tensor\mon{f}(\slot))
\Rightarrow \mon{f}(\slot)\tensor \mon{f}(\slot)$
between functors from $\CT\times \CT$ to $\CT$
whose component at $X$, $Z\in \CT$ is given by
\begin{equation*}
  \xymatrix@C=40pt{
    \mon{f}(X\tensor \mon{f}(Z))
    \ar[r]^-{(\mon{f}_{\tensor})_{X,\mon{f}(Z)}}
    &\mon{f}(X)\tensor \mon{f}\mon{f}(Z)
    \ar[r]^-{I_{\mon{f}(X)}\tensor \mu^{\mon{f}}_Z}
    &\mon{f}(X)\tensor \mon{f}(Z)
    .
  }
\end{equation*}
A \emph{colax $\kos{S}$-tensor Hopf monad}
is a colax $\kos{S}$-tensor monad
whose fusion operator is invertible.
The terminology
is taken from
\cite[sections~2.6 and~2.7]{BruguieresLackVirelizier2011},
where $\Sigma^{\mon{f}}$ is the left fusion operator
of the underlying bimonad.
We denote by
$\cat{HMonad}[\id_{\kos{S}}]_{\mathrm{rfl}}$
the category of reflective colax $\kos{S}$-tensor Hopf monads
on $\id_{\kos{S}}$,
which is a full subcategory of
$\cat{Monad}[\id_{\kos{S}}]_{\mathrm{rfl}}$.
By restricting the adjoint equivalence of categories
in Corollary~\ref{cor-cltcatoversmc-MonEnscolaxKtmonadequiv},
we obtain the following proposition.

\begin{proposition}
  \label{prop-cltcatoversmc-GrpEnscolaxKtHmonadequiv}
  We have an adjoint equivalence of categories
  \begin{equation*}
    \vcenter{\hbox{
      \xymatrix@C=25pt{
        \cat{HMonad}[\id_{\kos{S}}]_{\mathrm{rfl}}
        \ar@<1ex>[r]^-{\mathcal{L}}
        \ar@{}[r]|-{\simeq}
        &\cat{Grp}(\kosEns(\kos{S})).
        \ar@<1ex>[l]^-{\iota}
      }
    }}
  \end{equation*}
\end{proposition}
\begin{proof}
  It suffices to show that the functors
  $\iota$ and $\mathcal{L}$
  in Corollary~\ref{cor-cltcatoversmc-MonEnscolaxKtmonadequiv}
  restrict to the full subcategories
  $\cat{Grp}(\kosEns(\kos{S}))$
  and $\cat{HMonad}[\id_{\kos{S}}]_{\mathrm{rfl}}$.

  Let $\pi$ be a cocommutative bimonoid in $\kos{S}$.
  The component of the fusion operator of
  $\iota(\pi)=\pi\otimes \slot$
  at $x$, $z\in \CS$ is
  $\Sigma^{\pi\otimes \slot}_{x,z}
  =(I_{\pi\otimes x}\otimes \pc_{\pi}\otimes I_z)
  \circ (I_{\pi}\otimes \sigma_{\pi,x}\otimes I_{\pi\otimes z})
  \circ (\cp_{\pi}\otimes I_{x\otimes \pi\otimes z})
  \colon \pi\otimes x\otimes \pi\otimes z
  \to \pi\otimes x\otimes \pi\otimes z$.
  The fusion operator $\Sigma^{\pi\otimes \slot}$
  is invertible if and only if
  $\pi$ is a cocommutative Hopf monoid in $\kos{S}$.
  Indeed, if $\varsigma_{\pi}$ exists,
  then $\Sigma^{\pi\otimes \slot}_{x,z}$ is invertible
  with inverse
  $(\Sigma^{\pi\otimes \slot}_{x,z})^{-1}
  =(I_{\pi\otimes x}\otimes \pc_{\pi}\otimes I_z)
  \circ (I_{\pi\otimes x}\otimes \varsigma_{\pi}\otimes I_{\pi\otimes z})
  \circ (I_{\pi}\otimes \sigma_{\pi,x}\otimes I_{\pi\otimes z})
  \circ (\cp_{\pi}\otimes I_{x\otimes \pi\otimes z})
  \colon \pi\otimes x\otimes \pi\otimes z
  \to \pi\otimes x\otimes \pi\otimes z$.
  Conversely, if $\Sigma^{\pi\otimes \slot}$ is invertible,
  then so is $\tau_{\pi}$,
  as it is equal to
  $\Sigma^{\pi\otimes \slot}_{s,s}$
  up to coherence isomorphisms,
  and hence $\varsigma_{\pi}$ exists.
  In particular, 
  $\iota$ restricts to a functor from
  $\cat{Grp}(\kosEns(\kos{S}))$ 
  to $\cat{HMonad}[\id_{\kos{S}}]_{\mathrm{rfl}}$.

  Let $\mon{f}$ be a reflective colax $\kos{S}$-tensor Hopf monad on $\id_{\kos{S}}$,
  and put $\pi\coloneqq\mon{f}(s)$.
  Then $\hatar{\mon{f}}\colon \mon{f}\cong \pi\otimes \slot$
  in \eqref{eq-cltcatoversmc-hatardef}
  is an isomorphism of colax $\kos{S}$-tensor monads,
  as it is the component of the adjunction unit
  of $\mathcal{L}\dashv\iota$ at $\mon{f}$.
  Thus the fusion operator of $\pi\otimes \slot$ is invertible as well,
  and hence $\mathcal{L}(\mon{f})=\pi$
  is a cocommutative Hopf monoid in $\kos{S}$.
  This shows that $\mathcal{L}$ restricts to a functor from
  $\cat{HMonad}[\id_{\kos{S}}]_{\mathrm{rfl}}$
  to $\cat{Grp}(\kosEns(\kos{S}))$.
  One can check that the antipode
  $\varsigma_{\pi}=(e_{\pi}\otimes I_{\pi})\circ \tau_{\pi}^{-1}\circ (I_{\pi}\otimes u_{\pi})$
  can also be written as
  $\varsigma_{\mon{f}(s)}
  =(\what{\mon{f}})_{\mon{f}(s)}
  \circ \mon{f}(\imath_{\mon{f}(s)}^{-1})
  \circ (\Sigma^{\mon{f}}_{s,s})^{-1}
  \circ (I_{\mon{f}(s)}\otimes \upsilon^{\mon{f}}_{s})
  \colon \mon{f}(s)
  \to \mon{f}(s)$.
\end{proof}

%% file: PKGGC.bbl
\begin{thebibliography}{MLM94}

\bibitem[Bec68]{Beck1968}
Jonathan~Mock Beck, \emph{The tripleableness theorem}, 1968, Unpublished manuscript. Republished in Repr. Theory Appl. Categ., no.~31 (2025), pp.~1--11.

\bibitem[BLV11]{BruguieresLackVirelizier2011}
Alain Brugui{\`e}res, Steve Lack, and Alexis Virelizier, \emph{{Hopf} monads on monoidal categories}, Adv. Math. \textbf{227} (2011), no.~2, 745--800. \MR{2793022}

\bibitem[Bor94]{Borceux1994}
Francis Borceux, \emph{Handbook of categorical algebra. 2. {C}ategories and structures}, Encyclopedia Math. Appl., vol.~51, Cambridge University Press, Cambridge, 1994. \MR{1313497}

\bibitem[Bre94]{Breen1994}
Lawrence Breen, \emph{Tannakian categories}, Motives ({S}eattle, {WA}, 1991), {P}art 1 (Uwe Jannsen, Steven Kleiman, and Jean-Pierre Serre, eds.), Proc. Sympos. Pure Math., vol.~55, Amer. Math. Soc., Providence, RI, 1994, pp.~337--376. \MR{1265536}

\bibitem[Bre10]{Breen2010}
Lawrence Breen, \emph{Notes on 1- and 2-gerbes}, Towards higher categories (John~C. Baez and J.~Peter May, eds.), IMA Vol. Math. Appl., vol. 152, Springer, New York, 2010, pp.~193--235. \MR{2664623}

\bibitem[Buc14]{Buckley2014}
Mitchell Buckley, \emph{Fibred 2-categories and bicategories}, J. Pure Appl. Algebra \textbf{218} (2014), no.~6, 1034--1074. \MR{3153612}

\bibitem[Cad13]{Cadoret2013}
Anna Cadoret, \emph{Galois categories}, Arithmetic and geometry around {G}alois theory, Progr. Math., vol. 304, Birkh\"{a}user/Springer, Basel, 2013, pp.~171--246. \MR{3408165}

\bibitem[Del90]{Deligne1990}
P.~Deligne, \emph{Cat\'{e}gories tannakiennes}, The {G}rothendieck {F}estschrift, {V}ol. {II}, Progr. Math., vol.~87, Birkh\"{a}user Boston, Boston, MA, 1990, pp.~111--195. \MR{1106898}

\bibitem[DM82]{DeligneMilne1982}
P.~Deligne and J.~S. Milne, \emph{Tannakian categories}, Hodge cycles, motives, and {S}himura varieties, Lecture Notes in Math., vol. 900, Springer-Verlag, Berlin-New York, 1982, pp.~101--228; updated and corrected version dated February 11, 2022, available at \url{https://www.jmilne.org/math/xnotes/tc.html}.

\bibitem[Ems17]{Emsalem2017}
Michel Emsalem, \emph{Twisting by a torsor}, Analytic and algebraic geometry (Anilatmaja Aryasomayajula, Indranil Biswas, Archana~S. Morye, and A.~J. Parameswaran, eds.), Springer, Singapore, 2017, pp.~125--152. \MR{3728130}

\bibitem[Gir71]{Giraud1971}
Jean Giraud, \emph{Cohomologie non ab{\'e}lienne}, Grundlehren Math. Wiss., vol. 179, Springer-Verlag, Berlin-New York, 1971. \MR{344253}

\bibitem[Gro71]{Grothendieck1971}
A.~Grothendieck, \emph{Rev\^{e}tements \'{e}tales et groupe fondamental ({SGA} 1)}, Lecture Notes in Math., vol. 224, Springer, Berlin-New York, 1971, S\'{e}minaire de G\'{e}om\'{e}trie Alg\'{e}brique du Bois Marie 1960--61, with two papers by M. Raynaud; updated and annotated edition, Doc. Math. (Paris) 3, Soc. Math. France, Paris, 2003. \MR{354651}

\bibitem[Hat02]{Hatcher2002}
Allen Hatcher, \emph{Algebraic topology}, Cambridge University Press, Cambridge, 2002. \MR{1867354}

\bibitem[Joh02a]{Johnstone2002a}
Peter~T. Johnstone, \emph{Sketches of an elephant: a topos theory compendium. {V}ol. 1}, Oxford Logic Guides, vol.~43, Oxford University Press, Oxford, 2002. \MR{1953060}

\bibitem[Joh02b]{Johnstone2002b}
Peter~T. Johnstone, \emph{Sketches of an elephant: a topos theory compendium. {V}ol. 2}, Oxford Logic Guides, vol.~44, Oxford University Press, Oxford, 2002. \MR{2063092}

\bibitem[JY21]{Johnson2021}
Niles Johnson and Donald Yau, \emph{2-dimensional categories},
Oxford University Press, Oxford, 2021,
Prepublication version: \mbox{\url{https://arxiv.org/abs/2002.06055v3}}.
References to numbered items in this book are to the prepublication version. \MR{4261588}

\bibitem[Kel74]{Kelly1974}
G.~M. Kelly, \emph{Doctrinal adjunction}, Category {S}eminar ({P}roc. {S}em., {S}ydney, 1972/1973), Lecture Notes in Math., vol. 420, Springer, Berlin-New York, 1974, pp.~257--280. \MR{360749}

\bibitem[KS74]{KellyStreet1974}
G.~M. Kelly and Ross Street, \emph{Review of the elements of 2-categories}, Category {S}eminar ({P}roc. {S}em., {S}ydney, 1972/1973) (G.~M. Kelly, ed.), Lecture Notes in Math., vol. 420, Springer, Berlin, 1974, pp.~75--103. \MR{357542}

\bibitem[Mil25]{Milne2025}
James~S. Milne, \emph{Tannakian categories: origins and summary}, arXiv:2502.10945v2, 2025.

\bibitem[ML63]{MacLane1963}
Saunders Mac~Lane, \emph{Natural associativity and commutativity}, Rice Univ. Stud. \textbf{49} (1963), no.~4, 28--46. \MR{170925}

\bibitem[MLM94]{MacLaneMoerdijk1994}
Saunders Mac~Lane and Ieke Moerdijk, \emph{Sheaves in geometry and logic}, Universitext, Springer-Verlag, New York, 1994, A first introduction to topos theory. Corrected reprint of the 1992 edition. \MR{1300636}

\bibitem[Moe02]{Moerdijk2002}
Ieke Moerdijk, \emph{Monads on tensor categories}, J. Pure Appl. Algebra \textbf{168} (2002), no.~2--3, 189--208. \MR{1887157}

\bibitem[SR72]{SaavedraRivano1972}
Neantro Saavedra~Rivano, \emph{Cat\'{e}gories {T}annakiennes}, Lecture Notes in Math., vol. 265, Springer, Berlin-New York, 1972. \MR{338002}

\bibitem[Sza09]{Szamuely2009}
Tam{\'a}s Szamuely, \emph{{G}alois groups and fundamental groups}, Cambridge Stud. Adv. Math., vol. 117, Cambridge University Press, Cambridge, 2009. \MR{2548205}

\end{thebibliography}
